\documentclass{memo-l}
\usepackage{fontenc}
\usepackage{inputenc}
\usepackage{float}
\usepackage{mathrsfs}
\usepackage{minibox}
\usepackage{amsthm}
\usepackage{amsmath}
\usepackage{amssymb}
\usepackage{esint}
\usepackage{enumitem}
\usepackage{multirow}
\usepackage{graphicx}
\usepackage{color}
\usepackage{xcolor}
\usepackage[english]{babel}
\usepackage[title, titletoc]{appendix}
\definecolor{heavyblue}{cmyk}{1,1,0,0.25}
\usepackage{hyperref}

\makeatletter

\def\l@subsection{\@tocline{2}{0pt}{2pc}{5pc}{}}
\setlist{align=left,leftmargin=*,labelindent=1.5pc,labelsep=0.5pc, itemsep=5pt}
\setlist[description]{leftmargin=4.5cm,style=multiline}
\setenumerate[enumerate,1]{align=left,leftmargin=*,labelindent=1.5pc,labelsep=0pc}
\setenumerate[enumerate,2]{align=left,leftmargin=*,labelindent=0.5pc,labelsep=0pc}
\numberwithin{equation}{section}

\theoremstyle{plain}
\newtheorem{thm}{\protect\theoremname}[section]
  \theoremstyle{definition}
  
  \theoremstyle{definition}
  \newtheorem{defn}[thm]{\protect\definitionname}
  \theoremstyle{definition}
  \newtheorem{example}[thm]{\protect\examplename}
  \theoremstyle{remark}
  \newtheorem{rem}[thm]{\protect\remarkname}
  \theoremstyle{plain}
  \newtheorem{prop}[thm]{\protect\propositionname}
  \theoremstyle{plain}
  \newtheorem{cor}[thm]{\protect\corollaryname}
  \theoremstyle{plain}
  \newtheorem{lem}[thm]{\protect\lemmaname}
  \theoremstyle{definition}
  \newtheorem*{condition}{\protect\conditionname}
  \theoremstyle{definition}
  \newtheorem{notation}[thm]{\protect\notationname}
  \theoremstyle{plain}
  \newtheorem{fact}[thm]{\protect\factname}
  \theoremstyle{plain}
  
  \theoremstyle{plain}
  \newtheorem*{claim}{\protect\claimname}
  \theoremstyle{remark}
  \newtheorem{note}[thm]{\protect\notename}

\makeatother

  \providecommand{\assumptionname}{Assumption}
  \providecommand{\definitionname}{Definition}
  \providecommand{\examplename}{Example}
  \providecommand{\lemmaname}{Lemma}
  \providecommand{\propositionname}{Proposition}
  \providecommand{\remarkname}{Remark}
\providecommand{\corollaryname}{Corollary}
\providecommand{\theoremname}{Theorem}
\providecommand{\conditionname}{Condition}
\providecommand{\problemname}{Question}
\providecommand{\notationname}{Notation}
\providecommand{\factname}{Fact}
\providecommand{\claimname}{Property}
\providecommand{\notename}{Note}

\begin{document}

\frontmatter

\title{Replication and Its Application to Weak Convergence}

\author{Chi Dong}
\address{Department of Mathematical and Statistical Sciences, University of Alberta, Edmonton, AB T6G 2G1, Canada}
\email{cd3@ualberta.ca}
%\thanks{}

\author{Michael A. Kouritzin}
\address{Department of Mathematical and Statistical Sciences, University of Alberta, Edmonton, AB T6G 2G1, Canada}
\email{michaelk@ualberta.ca}
%\thanks{}

\date{September 12, 2019}



%\dedicatory{Dedication text (use \\[2pt] for line break if necessary)}

\maketitle

\tableofcontents

\begin{abstract}
Herein, a methodology is developed to replicate functions,
measures and stochastic processes onto a compact metric space.
Many results are easily established
for the replica objects and then transferred back 
to the original ones.
Two problems are solved within 
to demonstrate the method: 
(1) Finite-dimensional convergence for processes living on 
general topological spaces. 
(2) New tightness and relative compactness criteria are given
for the Skorokhod space $D(\mathbf{R}^{+};E)$ 
with $E$ being a general Tychonoff space.
The methods herein are also used 
in companion papers to establish the: 
(3) existence of, uniqueness of and convergence to 
martingale problem solutions,
(4) classical Fujisaki-Kallianpur-Kunita and Duncan-Mortensen-Zakai filtering equations 
and stationary filters,
(5) finite-dimensional convergence
to stationary signal-filter pairs,
(6) invariant measures of Markov processes, and
(7) Ray-Knight theory
all in general settings.
\end{abstract}

\mainmatter

\chapter*{Frequently Used Notations}

\itemsep2pt

$\;$

\subsection*{Set, numbers and mappings}

$\;$

$\;$
\begin{description}
\item [{\textmd{$\circeq$}}] ``Being defined by''.
\item [{\textmd{$\varnothing$}}] Empty set.
\item [{\textmd{$\mathbf{I}$}}] Non-empty index set.
\item [{\textmd{$\mathbf{N},\;\mathbf{N}_{0}$}}] Positive and nonnegative
integers respectively.
\item [{\textmd{$\mathbf{Q},\;\mathbf{Q}^{+}$}}] Rational and nonnegative
rational numbers respectively.
\item [{\textmd{$\mathbf{R},\;\mathbf{R}^{+}$}}] Real and nonnegative
real numbers respectively.
\item [{\textmd{$\mathbf{R}^{+}$}}] Nonnegative real numbers.
\item [{\textmd{$(\mathbf{R}^{k},\left|\cdot\right|)$}}] $k$-dimensional
($k\in\mathbf{N}$) Euclidean space with Euclidean norm $\left|\cdot\right|$.
\item [{\textmd{$\uparrow,\;\downarrow$}}] Non-decreasing and non-increasing
convergence of real numbers (including convergence to $\infty$ and
$-\infty$) respectively.
\item [{\textmd{$\subset,\;\supset$}}] Containment of sets including equalities.
\item [{\textmd{$\aleph(E)$}}] Cardinality of set $E$.
\item [{\textmd{$\mathscr{P}_{0}(E)$}}] All finite non-empty subsets of
$E$.
\item [{\textmd{$\circ$}}] Composition of mappings.
\item [{\textmd{$\mathbf{1}_{A}$}}] Indicator function of set $A$.
\item [{\textmd{$f|_{A}$}}] Restriction of mapping $f$ to subset $A$
of its domain.
\item [{\textmd{$\mathcal{D}|_{A}$}}] Restrictions of the members of mapping
family $\mathcal{D}$ to $A$ (see (\ref{eq:Map_Class_Restrict})).
\item [{\textmd{$\bigotimes\mathcal{D}$}}] Joint mapping of members of
$\mathcal{D}$ (see (\ref{eq:Prod_Map})).
\item [{\textmd{$\mathfrak{var}$}}] Variant of mapping (see Notation \ref{notation:Var}).
\end{description}
$\;$

\subsection*{Measurable and measure spaces}

$\;$

$\;$
\begin{description}
\item [{\textmd{$\mathscr{U}|_{A}$}}] Concentration of $\sigma$-algebra
$\mathscr{U}$ on $A$ (see (\ref{eq:Concentrated_Sigma_Algebra})).
\item [{\textmd{$\sigma(\mathcal{D})$}}] $\sigma$-algebra induced by
mapping family $\mathcal{D}$ (see (\ref{eq:Baire_Algebra})).
\item [{\textmd{$\delta_{x}$}}] Dirac measure at $x$.
\item [{\textmd{$\mathfrak{M}^{+}(E,\mathscr{U}),\;\mathfrak{P}(E,\mathscr{U})$}}] Non-trivial
finite and probability measures on measurable space $(E,\mathscr{U})$
respectively.
\item [{\textmd{$\mathscr{N}(\mu)$}}] Subsets of measure space $(E,\mathscr{U},\mu)$
with zero outer measure (see \S \ref{sub:Meas}).
\item [{\textmd{$\nu|^{E}$}}] Epansion of $\nu\in\mathfrak{M}^{+}(A,\mathscr{U}|_{A})$
onto superspace $E$ (see (\ref{eq:Measure_Extend})).
\item [{\textmd{$\mu|_{A}$}}] Concentration of $\mu\in\mathfrak{M}^{+}(E,\mathscr{U})$
on subset $A$ (see (\ref{eq:Measure_Restrict})).
\item [{\textmd{$\mu\circ f^{-1}$}}] Push-forward measure of $\mu$ by
mapping $f$ (see (\ref{eq:Measure_PushForward})).
\end{description}
$\;$

\subsection*{Topological spaces}

$\;$

$\;$
\begin{description}
\item [{\textmd{$\mathscr{O}(E),\;\mathscr{C}(E)$}}] Open and closed subsets
of topological space $E$ respectively.
\item [{\textmd{$\mathscr{K}(E),\;\mathscr{K}^{\mathbf{m}}(E)$}}] Compact
and metrizable compact subsets of $E$ respectively.
\item [{\textmd{$\mathscr{K}_{\sigma}(E),\;\mathscr{K}_{\sigma}^{\mathbf{m}}(E)$}}] $\sigma$-compact
and $\sigma$-metrizable compact subsets of $E$.
\item [{\textmd{$\mathscr{B}(E)$}}] Borel subsets of $E$.
\item [{\textmd{$\mathscr{O}_{E}(A)$}}] Subspace topology of $A$ induced
from $E$ (see (\ref{eq:O_E_(A)})).
\item [{\textmd{$\mathscr{B}_{E}(A)$}}] Subspace Borel $\sigma$-algebra
of $A$ induced from $E$ (see (\ref{eq:B_E_(A)})).
\item [{\textmd{$A^{\epsilon}$}}] $\epsilon$-envelope of subset $A$
of a metric space (see (\ref{eq:Epsilon_envelope})).
\item [{\textmd{$\mathscr{O}_{\mathcal{D}}(A)$}}] Topology induced by
mapping family $\mathcal{D}$ on $A$ (see (\ref{eq:O_D_(A)_Topo_Basis})).
\item [{\textmd{$\mathscr{B}_{\mathcal{D}}(A)$}}] Borel $\sigma$-algebra
induced by $\mathcal{D}$ on $A$, i.e. $\sigma[\mathscr{O}_{\mathcal{D}}(A)]$.
\end{description}
$\;$

\subsection*{Product spaces}

$\;$

$\;$
\begin{description}
\item [{\textmd{$\times$}}] Cartesian product.
\item [{\textmd{$\prod_{i\in\mathbf{I}}S_{i}$}}] Cartesian product of
$\{S_{i}\}_{i\in\mathbf{I}}$.
\item [{\textmd{$E^{\mathbf{I}}$}}] Cartesian power of $E$ with respect
to index set $\mathbf{I}$; in particular, $E^{\infty}$ abbreviates
$E^{\mathbf{I}}$ when $\aleph(\mathbf{I})=\aleph(\mathbf{N})$ and
$E^{d}$ abbreviates $E^{\mathbf{I}}$ when $\aleph(\mathbf{I})=d\in\mathbf{N}$.
\item [{\textmd{$\mathfrak{p}_{\mathbf{I}_{0}}$}}] Projection mapping
on $\prod_{i\in\mathbf{I}}S_{i}$ for non-empty sub-index-set $\mathbf{I}_{0}\subset\mathbf{I}$;
in particular, $\mathfrak{p}_{j}$ abbreviates $\mathfrak{p}_{\{j\}}$.
\item [{\textmd{$\otimes$}}] Product of $\sigma$-algebras (or topologies).
\item [{\textmd{$\bigotimes_{i\in\mathbf{I}}\mathscr{A}_{i}$}}] Product
$\sigma$-algebra (or product topology) on the Cartesian product $\prod_{i\in\mathbf{I}}S_{i}$
of measurable (or topological) spaces $\{(S_{i},\mathscr{A}_{i})\}_{i\in\mathbf{I}}$;
in particular, $\mathscr{A}^{\otimes\mathbf{I}}$ abbreviates $\bigotimes_{i\in\mathbf{I}}\mathscr{A}_{i}$
when all $\mathscr{A}_{i}=\mathscr{A}$, $\mathscr{A}^{\otimes\infty}$
abbreviates $\mathscr{A}^{\otimes\mathbf{I}}$ when $\aleph(\mathbf{I})=\aleph(\mathbf{N})$,
and $\mathscr{A}^{\otimes d}$ abbreviates $\mathscr{A}^{\otimes\mathbf{I}}$
when $\aleph(\mathbf{I})=d\in\mathbf{N}$.
\end{description}
$\;$

\subsection*{Spaces of general mappings}

$\;$

$\;$
\begin{description}
\item [{\textmd{$\varpi_{\mathbf{I}}(f)$}}] Associated path mapping of
mapping $f\in S^{E}$ (see (\ref{eq:Path_Mapping})); in particular,
$\varpi(f)$, $\varpi_{T}(f)$ and $\varpi_{a,b}(f)$ abbreviates
$\varpi_{\mathbf{R}^{+}}(f)$, $\varpi_{[0,T]}(f)$ and $\varpi_{[a,b]}(f)$
respectively.
\item [{\textmd{$\varpi_{\mathbf{I}}(\mathcal{D})$}}] Joint path mapping
of mapping family $\mathcal{D}\subset S^{E}$ (see (\ref{eq:Path_Mapping_Class})).
\item [{\textmd{$M(S;E)$}}] Measurable mappings from measurable space
$S$ to measurable space $E$.
\item [{\textmd{$C(S;E)$}}] Continuous mappings from topological space
$S$ to topological space $E$.
\item [{\textmd{$\mathbf{hom}(S;E),\;\mathbf{biso}(S;E)$}}] Homeomorphisms
and Borel isomorphisms between $S$ and $E$ respectively.
\item [{\textmd{$\mathbf{imb}(S;E)$}}] Imbeddings from $S$ to $E$.
\item [{\textmd{$\mathfrak{D}(\mathcal{L}),\;\mathfrak{R}(\mathcal{L})$}}] Domain
and range of single-valued linear operator $\mathcal{L}$ respectively
(see (\ref{eq:D(L)}) and (\ref{eq:R(L)})).
\item [{\textmd{$\mathcal{L}|_{\mathcal{D}}$}}] Restriction of $\mathcal{L}$
to $\mathcal{D}\subset\mathfrak{D}(\mathcal{L})$ (see (\ref{eq:L|D})).
\end{description}
$\;$

\subsection*{Skorokhod $\mathscr{J}_{1}$-spaces}

$\;$

$\;$
\begin{description}
\item [{\textmd{$w_{\mathfrak{r},\delta,T}^{\prime}(x)$}}] $\mathfrak{r}$-modulus
of continuity of $x\in E^{\mathbf{R}^{+}}$ (see (\ref{eq:w'})).
\item [{\textmd{$J(x)$}}] Set of left-jump times of $x\in E^{\mathbf{R}^{+}}$
(see (\ref{eq:J(x)})).
\item [{\textmd{$\mathbf{TC}(\mathbf{R}^{+}),\;\mathbf{TC}([a,b])$}}] All
time-changes on $\mathbf{R}^{+}$ and $[a,b]$ respectively (see \S\ref{sub:Meas_Cont_Cadlag_Map}).
\item [{\textmd{$\varrho^{\mathfrak{r}}(x,y)\;\varrho_{[a,b]}^{\mathfrak{r}}(x,y)$}}] Skorokhod
pseudometric for $x,y\in E^{\mathbf{R}^{+}}$ and $E^{[a,b]}$ induced
by $\mathfrak{r}$ respectively (see (\ref{eq:SkoMetric}), (\ref{eq:SkoMetric_[a,b]})
and (\ref{eq:r[a,b]})).
\item [{\textmd{$(D(\mathbf{R}^{+};E),\mathscr{J}(E))$}}] Skorokhod $\mathscr{J}_{1}$-space
of all c$\grave{\mbox{a}}$dl$\grave{\mbox{a}}$g members of $E^{\mathbf{R}^{+}}$
(see \S\ref{sub:Meas_Cont_Cadlag_Map}).
\item [{\textmd{$(D([a,b];E),\mathscr{J}_{a,b}(E))$}}] Skorokhod $\mathscr{J}_{1}$-space
of all c$\grave{\mbox{a}}$dl$\grave{\mbox{a}}$g members of $E^{[a,b]}$
(see \S\ref{sub:Meas_Cont_Cadlag_Map}).
\item [{\textmd{$J(\mu)$}}] Set of fixed left-jump times of a finite measure
$\mu$ on $D(\mathbf{R}^{+};E)$ equipped with the restriction of
$,\mathscr{B}(E)^{\otimes\mathbf{R}^{+}}$ (see (\ref{eq:J(Mu)})).
\end{description}
$\;$

\subsection*{Spaces of functions}

$\;$

$\;$
\begin{description}
\item [{\textmd{$f\vee g(x),\; f\wedge g(x)$}}] $\max\{f(x),g(x)\}$ and
$\min\{f(x),g(x)\}$ for $\{f,g\}\subset\mathbf{R}^{E}$ respectively.
\item [{\textmd{$f^{+}(x),\; f^{-}(x)$}}] $\max\{f(x),0\}$ and $\max\{-f(x),0\}$
respectively.
\item [{\textmd{$\mathfrak{ae}(\mathcal{D}),\;\mathfrak{ac}(\mathcal{D}),\;\mathfrak{mc}(\mathcal{D})$}}] Additive
expansion, additive closure and multiplicative closure of $\mathbf{R}^{k}$-valued
function family $\mathcal{D}$ respectively (see (\ref{eq:Add_Extension}),
(\ref{eq:Add_Closure}), (\ref{eq:Mult_Closure})).
\item [{\textmd{$\mathfrak{ag}_{\mathbf{Q}}(\mathcal{D}),\;\mathfrak{ag}(\mathcal{D})$}}] $\mathbf{Q}$-algebra
and algebra of $\mathcal{D}$ respectively (see (\ref{eq:ag_Q_(D)}),
(\ref{eq:ag(D)})).
\item [{\textmd{$\Pi^{\mathbf{I}}(\mathcal{D})$}}] Product functions of
$\mathbf{R}$-valued function family $\mathcal{D}$ on $E^{\mathbf{I}}$
for finite $\mathbf{I}$; in particular, $\Pi^{d}(\mathcal{D})$ abbreviates
$\Pi^{\mathbf{I}}(\mathcal{D})$ with $d\circeq\aleph(\mathbf{I})$
(see (\ref{eq:Pi^d})).
\item [{\textmd{$\overset{u}{\rightarrow}$}}] Uniform convergence of $\mathbf{R}^{k}$-valued
functions.
\item [{\textmd{$\Vert\cdot\Vert_{\infty}$}}] Supremum norm of bounded
functions.
\item [{\textmd{$\mathfrak{cl}(\mathcal{D}),\;\mathfrak{ca}(\mathcal{D})$}}] Closure
of and closed algebra generated by bounded $\mathbf{R}^{k}$-valued
function family $\mathcal{D}$ under $\Vert\cdot\Vert_{\infty}$ respectively
(see (\ref{eq:ca(D)})).
\item [{\textmd{$M_{b}(E;\mathbf{R}^{k})$}}] Banach space of all bounded
members of $M(E;\mathbf{R}^{k})$ equipped with $\Vert\cdot\Vert_{\infty}$.
\item [{\textmd{$C_{b}(E;\mathbf{R}^{k})$}}] Banach space over scalar
field $\mathbf{R}$ of all bounded members of $C(E;\mathbf{R}^{k})$
equipped with $\Vert\cdot\Vert_{\infty}$.
\item [{\textmd{$C_{c}(E;\mathbf{R}^{k})$}}] Compact-supported continuous
functions from $E$ to $\mathbf{R}^{k}$.
\item [{\textmd{$C_{0}(E;\mathbf{R}^{k})$}}] Subspace of all $f\in C(E;\mathbf{R}^{k})$
such that for any $\epsilon>0$, there exists a $K_{\epsilon}\in\mathscr{K}(E)$
satisfying $\Vert f|_{E\backslash K_{\epsilon}}\Vert_{\infty}<\epsilon$.
\end{description}
$\;$

\subsection*{Finite Borel measures}

$\;$

$\;$
\begin{description}
\item [{\textmd{$\mathfrak{be}(\mu)$}}] All Borel extension(s) of $\mu$
(if any) (see \S\ref{sec:Borel_Measure}).
\item [{\textmd{$\mathcal{M}^{+}(E),\;\mathcal{P}(E)$}}] Weak topological
spaces of all finite and probability Borel measures on topological
space $E$ respectively.
\item [{\textmd{$f^{*}$}}] Integral functional of $f\in M_{b}(E;\mathbf{R}^{k})$
(see (\ref{eq:f_star})).
\item [{\textmd{$\mathcal{D}^{*}$}}] Integral functionals of members of
function family $\mathcal{D}\subset M_{b}(E;\mathbf{R})$ (see (\ref{eq:D_Star})).
\item [{\textmd{$\Longrightarrow$}}] Weak convergence of finite Borel
measures.
\item [{\textmd{$\mbox{w-}\lim_{n\rightarrow\infty}\mu_{n}$}}] Weak limit
(see (\ref{eq:WL})).
\end{description}
$\;$

\subsection*{Random variables and stochastic processes}

$\;$

$\;$
\begin{description}
\item [{\textmd{$(\Omega,\mathscr{F},\mathbb{P})$}}] Probability space.
\item [{\textmd{$\mathbb{E}$}}] Expectation operator of $(\Omega,\mathscr{F},\mathbb{P})$.
\item [{\textmd{$\mathrm{pd}(X)$}}] Process distribution of stochastic
process $X$ (see \S\ref{sec:Proc}).
\item [{\textmd{$X_{\mathbf{T}_{0}}$}}] $\mathfrak{p}_{\mathbf{T}_{0}}\circ X$,
the section of stochastic process$X$ for $\mathbf{T}_{0}\in\mathscr{P}_{0}(\mathbf{R}^{+})$
(see \S\ref{sec:Proc}).
\item [{\textmd{$J(X)$}}] Set of fixed left-jump times of stochastic process
$X$.
\item [{\textmd{$\mathscr{F}_{t}^{X}$}}] Augmented natural filtration
of $X$ (see (\ref{eq:FX})).
\item [{$\xrightarrow{\mathrm{D}(\mathbf{T})}$}] Finite-dimensional convergence
along $\mathbf{T}\subset\mathbf{R}^{+}$ (see Definition \ref{notation:FLP}).
\item [{$\mathfrak{fl}_{\mathbf{T}}(\{X^{n}\}_{n\in\mathbf{N}})$}] finite-dimensional
limit of process sequence $\{X^{n}\}_{n\in\mathbf{N}}$ along $\mathbf{T}\subset\mathbf{R}^{+}$
(see Definition \ref{notation:FLP}).
\item [{$\mathfrak{flp}_{\mathbf{T}}(\{X^{i}\}_{i\in\mathbf{I}})$}] Family
of finite-dimensional limit points of process family $\{X^{i}\}_{i\in\mathbf{I}}$
along $\mathbf{T}\subset\mathbf{R}^{+}$ (see Definition \ref{notation:FLP}).
\end{description}
$\;$

\subsection*{Replication}

$\;$

$\;$
\begin{description}
\item [{$(E_{0},\mathcal{F};\widehat{E},\widehat{\mathcal{F}})$}] Replication
base over $E$ (see Definition \ref{def:Base}).
\item [{$\overline{f}$}] Abbreviation of $\mathfrak{var}(f;\widehat{E},E_{0},0)$
(see Notation \ref{notation:Rep_Fun}).
\item [{$\widehat{f}$}] Replica function of $f$ (see Definition \ref{def:RepFun}
and Notation \ref{notation:Rep_Fun}).
\item [{$\widetilde{f}$}] $f|_{E_{0}^{d}}$ for base $(E_{0},\mathcal{F};\widehat{E},\widehat{\mathcal{F}})$
and $\mathbf{R}^{k}$-valued function $f$ on $E^{d}$ (see Notation
\ref{notation:OP}).
\item [{$\widetilde{\mathcal{F}}$}] $\mathcal{F}|_{E_{0}}$ for base $(E_{0},\mathcal{F};\widehat{E},\widehat{\mathcal{F}})$
(see Notation \ref{notation:OP}).
\item [{$\mathcal{L}_{0},\;\mathcal{L}_{1}$}] $\mathcal{L}|_{\mathfrak{ag}(\mathcal{F})}$
and $\mathcal{L}\cap(\mathfrak{ca}(\mathcal{F})\times\mathfrak{ca}(\mathcal{F}))$
respectively given base $(E_{0},\mathcal{F};\widehat{E},\widehat{\mathcal{F}})$
and linear operator $\mathcal{L}$ (see Notation \ref{notation:OP}).
\item [{$\widetilde{\mathcal{L}}_{i}$}] $\{(\widetilde{f},\widetilde{g}):(f,g)\in\mathcal{L}_{i}\}$
for each $i\in\{0,1\}$ (see Notation \ref{notation:OP}).
\item [{$\widehat{\mathcal{L}}_{0},\;\widehat{\mathcal{L}}_{1}$}] core
and extended replicas of linear operator $\mathcal{L}$ respectively
(see Definition \ref{def:RepOP}).
\item [{$\overline{\mu}$}] Replica measure of $\mu$ (see Definition \ref{def:RepMeas}
and Notation \ref{notation:ReplicaMeas}).
\item [{$\mathfrak{rep}(X;E_{0},\mathcal{F})$}] replica processes of stochastic
process $X$ with respect to base $(E_{0},\mathcal{F};\widehat{E},\widehat{\mathcal{F}})$;
in particular, $\mathfrak{rep}_{\mathrm{m}}(X;E_{0},\mathcal{F})$,
$\mathfrak{rep}_{\mathrm{p}}(X;E_{0},\mathcal{F})$ and $\mathfrak{rep}_{\mathrm{c}}(X;E_{0},\mathcal{F})$
denotes measurable, progressive and c$\grave{\mbox{a}}$dl$\grave{\mbox{a}}$g
members of $\mathfrak{rep}(X;E_{0},\mathcal{F})$ respectively (see
Definition \ref{notation:RepProc}).
\end{description}

\chapter{\label{chap:Intro}Introduction}

Traditionally, researchers (e.g. \cite{K71}, \cite{M83}, \cite{EK86},
\cite{BK93b}, \cite{J97a}, \cite{BBK00} and \cite{KL08}) have
focused on stochastic processes living on ``good'' topological spaces
with some metrizable, separable, completely metrizable or (locally)
compact properties. However, there are many settings of interest that
violate these assumptions. For example, \cite{HS79} and \cite{M83}
considered probability measures on the Skorokhod $\mathscr{J}_{1}$-space
of tempered distributions. \cite{S78} considered the nonlinear filtering
problem for c$\grave{\mbox{a}}$dl$\grave{\mbox{a}}$g signals living
on Lusin spaces. \cite{MZ84} considered tightness in the space of
all c$\grave{\mbox{a}}$dl$\grave{\mbox{a}}$g functions from the
non-negative real numbers $\mathbf{R}^{+}$ to the real line $\mathbf{R}$
equipped with the pseudo-path topology, which was further discussed
by \cite{S85} and \cite{K91}. \cite{F88} constructed Markov branching
processes whose values are finite Borel measures on a Lusin space.
\cite{J97a} extended the Skorokhod Representation Theorem to tight
sequences of probability measures on non-metrizable spaces. \cite{J97b}
considered a sequentially defined topology on the space of all c$\grave{\mbox{a}}$dl$\grave{\mbox{a}}$g
functions from the compact interval $[0,T]$ to $\mathbf{R}$. \cite{DZ98}
considered the space of all Borel probability measures on a Polish
space equipped with the strong topology. \cite{J86} and \cite{K15}
considered probability measures on $D([a,b];E)$, the Skorokhod $\mathscr{J}_{1}$-space
of all c$\grave{\mbox{a}}$dl$\grave{\mbox{a}}$g mappings from a
compact interval $[a,b]$ to a Tychonoff space $E$. \cite{L95,L98},
\cite{FV10} and \cite{FH14} worked with non-separable Banach spaces
of rough paths equipped with homogeneous $p$-variation or $1/p$-H$\ddot{\mbox{o}}$lder
norms. None of these spaces are necessarily Polish nor compact. Some
are not even metrizable or separable.

\begin{figure}[H]
\begin{centering}
\includegraphics[scale=0.6]{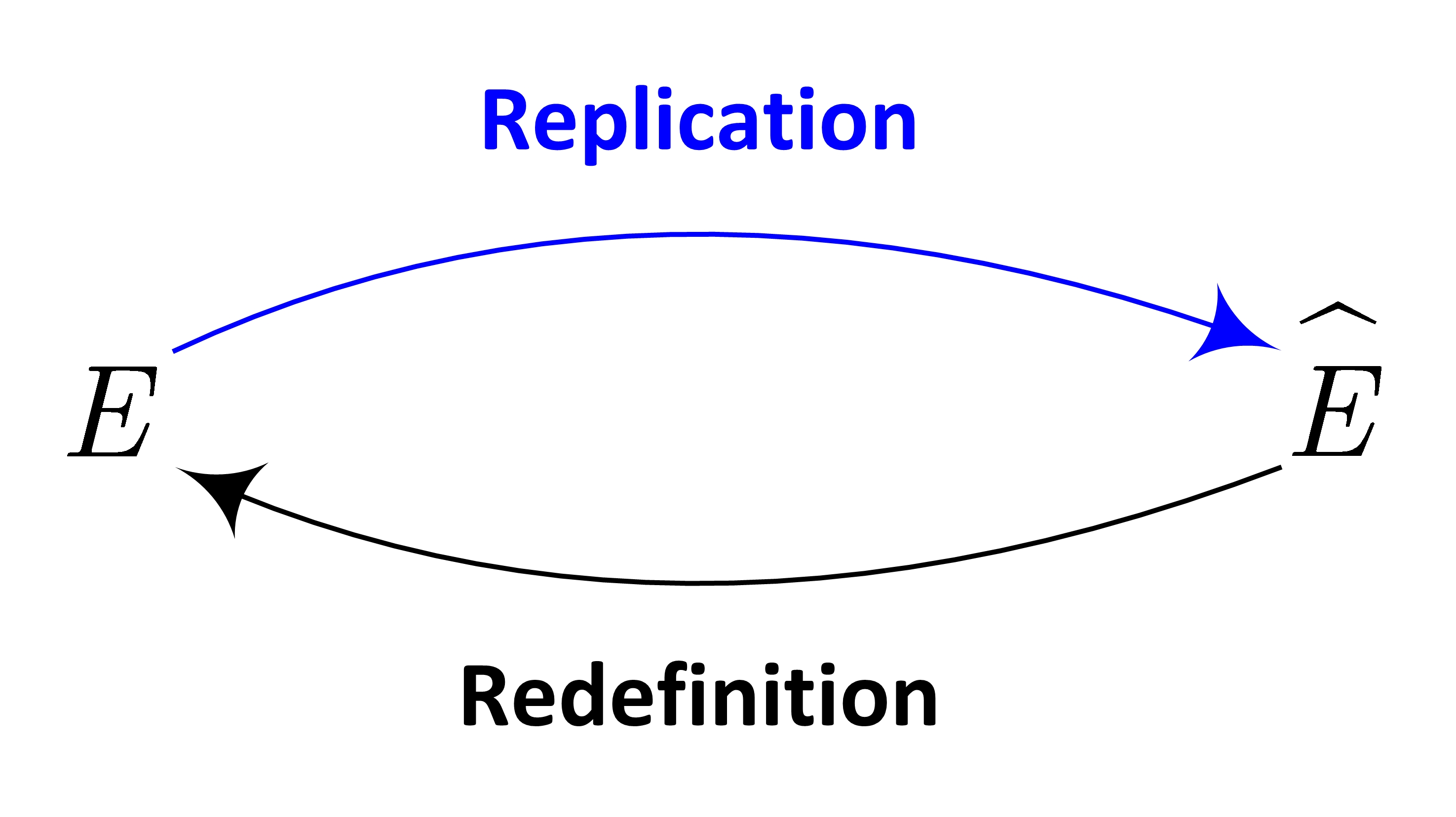}
\par\end{centering}

\caption{\label{fig:Main_Idea}\textit{The idea of replication}}
\end{figure}

The point of our work, as illustrated in Figure \ref{fig:Main_Idea}
above, is that Borel measurable functions, finite measures and stochastic
processes living on a \textit{general topological space} $E$ often
can be \textit{replicated} as \textit{replica functions},\textit{
measures} or\textit{ processes} living on some \textit{compact metric
space} $\widehat{E}$. Non-c$\grave{\mbox{a}}$dl$\grave{\mbox{a}}$g
processes can have c$\grave{\mbox{a}}$dl$\grave{\mbox{a}}$g replicas.
These replica objects are more easily analyzed on $\widehat{E}$ than
the original objects on $E$, and many results about the replica objects
are transferrable back to the original ones by proper (often indistinguishable)
redefinitions.

The idea of replication was motivated in part by \cite{EK86}, \cite{BK93b},
\cite{BK10} and \cite{K15} which exploit the use of imbedding and
compactification techniques in various aspects of probability theory.
One could extend results to various generalized settings in a case-by-case
manner. However, replication is probably an easier and more unified
approach of extending results from compact or Polish spaces to a large
category of exotic spaces simultaneously. This approach is believed
to have merit in areas such as weak convergence, martingale problems,
nonlinear filtering, large deviations, Markov processes etc., where
compactness or metric completeness can play a big role. Indeed, even
a Polish space can be improved by adding compactness.

The contributions herein are:

\begin{enumerate}
[label=\textbf{Theme \arabic*}, labelsep=0.5pc]

\item\label{enu:Theme1}The methodology of replication and its applications
to the Radon-Riesz Representation Theorem and Skorokhod Representation
Theorem (Chapter \ref{chap:Space_Change} - \ref{chap:Rep_Proc}).

\item\label{enu:Theme2}Finite-dimensional convergence of possibly
non-c$\grave{\mbox{a}}$dl$\grave{\mbox{a}}$g processes living on
general spaces (\S \ref{sec:RepProc_FC} and Chapter \ref{chap:FDDConv}).

\item\label{enu:Theme3}Tightness in Skorokhod $\mathscr{J}_{1}$-spaces,
the relationship between weak convergence on Skorokhod $\mathscr{J}_{1}$-spaces
and finite-dimensional convergence, and relative compactness in Skorokhod
$\mathscr{J}_{1}$-spaces (\S \ref{sec:RepProc_Path_Space} and Chapter
\ref{chap:Cadlag}).

\end{enumerate}

\ref{enu:Theme1} tells \textit{when} and \textit{how} one can perform
replication and serves as the theoretical foundation of all our developments.
The question \textit{what} replication can do is partially answered
by \ref{enu:Theme2} and \ref{enu:Theme3} herein. Further motivation
and illustration of the utility of replication can be found in the
companion works of this article (see \cite{DK20b,DK20c,DK21a}) and
other related papers (see e.g. \cite{KK20}).

\ref{enu:Theme2} grew out of the convergence of a stochastic evolution
system to its stationary distribution(s) or solution(s) over the long
term, which has been the central topic of many classical works on
both theory and application ends. For example, \cite[\S 10.2 and \S 10.4]{EK86}
considered the existence of stationary distributions for diffusion
approximations of the Wright-Fisher model. For a Fleming-Viot process
$X$, \cite{EK93} and \cite{DK99} considered the existence of a
stationary distribution $\mu$, and \cite{EK98} established the pointwise
ergodic theorem
\begin{equation}
\lim_{T\rightarrow\infty}\frac{1}{T}\int_{0}^{T}f(X_{t})dt=\int_{E}f(x)\mu(dx),\; a.s.\label{eq:LTB_exp-1}
\end{equation}
In nonlinear filtering, \cite[Theorem 4.1]{K71} considered the ``asymptotic
mean square filtering error''
\begin{equation}
\lim_{T\rightarrow\infty}\frac{1}{T}\int_{0}^{T}\mathbb{E}\left[\left(\pi_{t}^{\mu}(f)-f(X_{t}^{\mu})\right)^{2}\right]dt,\label{eq:LTB_exp-2}
\end{equation}
where $X^{\mu}$ is the signal with initial distribution $\mu$ and
$\pi^{\mu}(f)$ is the optimal filter for function $f$ of $X^{\mu}$.
\cite[(2.6)]{BK99} considered weak limit points of the ``pathwise
average error''
\begin{equation}
\frac{1}{T}\int_{0}^{T}\left(\pi_{t}^{\mu}(f)-f(X_{t}^{\mu})\right)^{2}dt\mbox{ as }T\uparrow\infty.\label{eq:LTB_exp-3}
\end{equation}
\cite[(1.2)]{B01} studied the ``$(\mu,\mu^{\prime})$-stability''
\begin{equation}
\lim_{T\rightarrow\infty}\frac{1}{T}\int_{0}^{T}\mathbb{E}\left[\left(\pi_{t}^{\mu}(f)-\pi_{t}^{\mu^{\prime}}(f)\right)^{2}\right]dt\label{eq:LTB_exp-4}
\end{equation}
of $\pi^{\mu}(f)$, where $\pi^{\mu^{\prime}}$ represents an approximate
filter with incorrect initiation $\mu^{\prime}$. In other areas,
\cite{CMP10} and \cite{CDP13} considered non-trivial stationary
solutions for the Lotka-Volterra model and those for perturbations
of the voter model. \cite{CG83} and \cite{C88} established (\ref{eq:LTB_exp-1})
for a basic voter process $X$ and an invariant measure $\mu$ of
$X$. All the works (and many others) above were based on separable
compact Hausdorff spaces, Polish spaces or compact metric spaces.
The following question, considering a weak and abstract form of long-time-average
limits like (\ref{eq:LTB_exp-1}), (\ref{eq:LTB_exp-2}), (\ref{eq:LTB_exp-3})
and (\ref{eq:LTB_exp-4}), still remains unanswered:

\begin{enumerate}
[label=\textbf{Q\arabic*}, labelsep=0.5pc]

\item\label{enu:Q_LTB}Let $E$ be a non-Polish, non-compact or even
non-metrizable space, and $X=\{X_{t}\}_{t\geq0}$ be an $E$-valued,
non-c$\grave{\mbox{a}}$dl$\grave{\mbox{a}}$g, measurable process.
Then, is there an $E$-valued stationary process $X^{\infty}$ such
that the long-time-averaged distributions
\begin{equation}
\frac{1}{T_{n}}\int_{0}^{T_{n}}\mathbb{P}\circ\left(X_{\tau+t_{1}},...,X_{\tau+t_{d}}\right)^{-1}d\tau\label{eq:LTB}
\end{equation}
converge weakly to the distribution of $(X_{t_{1}}^{\infty},...,X_{t_{d}}^{\infty})$
as $T_{n}\uparrow\infty$ for \textit{almost all} finite subsets $\{t_{1},...,t_{d}\}$
of $\mathbf{R}^{+}$?

\end{enumerate}

Weak convergence of finite-dimensional distributions of $E$-valued
c$\grave{\mbox{a}}$dl$\grave{\mbox{a}}$g processes is implied by
their weak convergence as $D(\mathbf{R}^{+};E)$-valued random variables
when $E$ is a separable metric space. Here, $D(\mathbf{R}^{+};E)$,
$\mathbf{N}$ and ``$\Rightarrow$'' denotes the Skorokhod $\mathscr{J}_{1}$-space
of all c$\grave{\mbox{a}}$dl$\grave{\mbox{a}}$g mappings from $\mathbf{R}^{+}$
to $E$, the positive integers and weak convergence respectively.
For processes $\{X^{n}\}_{n\in\mathbf{N}}$ and $X$ with paths in
$D(\mathbf{R}^{+};E)$,
\begin{equation}
X_{n}\Longrightarrow X\mbox{ as }n\uparrow\infty\mbox{ on }D(\mathbf{R}^{+};E)\label{eq:WC_on_Path_Space}
\end{equation}
has two implications:

\renewcommand{\labelenumi}{(\arabic{enumi}) }
\begin{enumerate}
\item $\{(X_{t_{1}}^{n},...,X_{t_{d}}^{n})\}_{n\in\mathbf{N}}$ converge
weakly to $(X_{t_{1}},...,X_{t_{d}})$ as $n\uparrow\infty$ for all
finite collection $\{t_{1},...,t_{d}\}$ in a dense subset of $\mathbf{R}^{+}$.
\item $\{X^{n}\}_{n\in\mathbf{N}}$ is relatively compact in $D(\mathbf{R}^{+};E)$.
\end{enumerate}

Often (2) requires strong or difficult-to-verify conditions in practice.
By contrast, (1) or a weaker form of it is believed to be establishable
for possibly non-c$\grave{\mbox{a}}$dl$\grave{\mbox{a}}$g processes
under much milder conditions than those for weak convergence on $D(\mathbf{R}^{+};E)$.
For instance, \cite{BK93b} discussed (1) with $X$ being a progressive
martingale problem solution and $\{X^{n}\}_{n\in\mathbf{N}}$ progressive
approximating processes, none of which is necessarily c$\grave{\mbox{a}}$dl$\grave{\mbox{a}}$g.
Their development was based on a Polish space $E$ and furthered that
of \cite[\S 4.8]{EK86}. Herein, we address the following general
questions on a more general $E$ without martingale problem setting:

\begin{enumerate}
[label=\textbf{Q\arabic*}, labelsep=0.5pc]
\setcounter{enumi}{1}

\item\label{enu:Q_FLP_Gen}When will a subsequence of $E$-valued
processes $\{X_{i}\}_{i\in\mathbf{I}}$ converge finite-dimensionally
to an $E$-valued process with general paths?

\item\label{enu:Q_FLP_Prog}When will a subsequence of $E$-valued
processes $\{X_{i}\}_{i\in\mathbf{I}}$ converge finite-dimensionally
to an $E$-valued progressive process?

\end{enumerate}

Either of these two questions may be answered in an individual way,
but replication helps to handle \ref{enu:Q_FLP_Gen}, \ref{enu:Q_FLP_Prog}
and the weak convergence of c$\grave{\mbox{a}}$dl$\grave{\mbox{a}}$g
processes on path spaces in one framework. We shall establish several
relatively mild and explicitly verifiable criteria for uniqueness
and existence of the limit processes in \ref{enu:Q_FLP_Gen} and \ref{enu:Q_FLP_Prog}
above. These criteria will be used to deduce (\ref{eq:LTB}) and answer
\ref{enu:Q_LTB} that motivates \ref{enu:Theme2}.

\ref{enu:Theme3} is concerned with two basic problems on a Tychonoff
space $E$:

\begin{enumerate}
[label=\textbf{Q\arabic*}, labelsep=0.5pc]
\setcounter{enumi}{3}

\item\label{enu:Q_Tight}When is a family of $E$-valued c$\grave{\mbox{a}}$dl$\grave{\mbox{a}}$g
processes bijectively indistinguishable from a tight family of $D(\mathbf{R}^{+};E)$-valued
random variables?

\item\label{enu:Q_RC}More generally, when is a family of $E$-valued
c$\grave{\mbox{a}}$dl$\grave{\mbox{a}}$g processes bijectively indistinguishable
from a relatively compact family of $D(\mathbf{R}^{+};E)$-valued
random variables?

\end{enumerate}

The main importance of tightness is that it implies relative compactness
for Borel probability measures on Hausdorff spaces. However, the verification
of tightness can be challenging. \cite{K75}, \cite{M83}, \cite{J86},
\cite{D91}, \cite{BK93b}, \cite{KX95} and \cite{K15}, to name
a few, all spent considerable efforts establishing tightness of c$\grave{\mbox{a}}$dl$\grave{\mbox{a}}$g
processes on exotic spaces. In particular, \cite{J86} developed systematic
tightness criteria for probability measures on both $D([0,1];E)$
and $D(\mathbf{R}^{+};E)$, which extended several results of \cite{M83}
and \cite[\S 3.7 - 3.9]{EK86} from the Polish to the possibly non-metrizable
Tychonoff case. \cite{K15} generalized the results of \cite{J86}
for $D([a,b];E)$ by loosening \cite[Theorem 3.1, (3.4)]{J86} to
the milder Weak Modulus of Continuity Condition (see \cite[\S 6]{K15}).
As a continuation of \cite{K15} on infinite time horizon, we answer
\ref{enu:Q_Tight} by establishing the equivalence among:
\begin{itemize}
\item Indistinguishability from a tight family of $D(\mathbf{R}^{+};E)$-valued
random variables;
\item Metrizable Compact Containment plus Weak Modulus of Continuity;
\item Metrizable Compact Containment plus Modulus of Continuity; and
\item Mild Pointwise Containment plus Modulus of Continuity for $\mathfrak{r}$
when $(E,\mathfrak{r})$ is a complete metric space.
\end{itemize}

As previously mentioned, weak convergnce on $D(\mathbf{R}^{+};E)$
is commonly thought to be composed of finite-dimensional convergence
along densely many times plus relative compactness in $D(\mathbf{R}^{+};E)$.
Herein, we give a more precise interpretation by showing that:
\begin{itemize}
\item Relative compactness in $D(\mathbf{R}^{+};E)$ implies the Modulus
of Continuity Condition for any Tychonoff space $E$.
\item When $E$ is a baseable space, weak convergence on $D(\mathbf{R}^{+};E)$
implies finite-dimensional convergence along densely many times. Baseable
spaces, defined in the sequel, are more general than metrizable and
separable spaces.
\item When $E$ is metrizable and separable, finite-dimensional convergence
along densely many times plus the Modulus of Continuity Condition
(weaker than relative compactness) are sufficient for weak convergence
on $D(\mathbf{R}^{+};E$).
\end{itemize}
Based on the results above, we answer \ref{enu:Q_RC} on metrizable
spaces by showing that:
\begin{itemize}
\item When $E$ is metrizable and separable, relative compactness in $D(\mathbf{R}^{+};E)$
is equivalent to the Modulus of Continuity Condition plus relative
compactness for finite-dimensional convergence to $E$-valued c$\grave{\mbox{a}}$dl$\grave{\mbox{a}}$g
processes along densely many times.
\item When $(E,\mathfrak{r})$ is a metric space, the combination of the
Modulus of Continuity Condition for $\mathfrak{r}$, Mild Pointwise
Containment Condition and relative compactness for finite-dimensional
convergence to $D(\mathbf{R}^{+};E)$-valued random variables along
densely many times is sufficient for relative compactness in $D(\mathbf{R}^{+};E)$.
\end{itemize}
Relative compactness can be weaker than tightness and require milder
conditions in non-Polish settings. The results of \ref{enu:Theme3}
demonstrate the Compact Containment Condition frequently used in verifying
tightness is superfluous for relative compactness in Skorokohod $\mathscr{J}_{1}$-spaces.
This is also why the second approach of \cite{BK93b} was well received.
While their work was restricted to a martingale problem setting, our
results are general.

The methods and results herein are also used in solving more concrete
problems.  For martingale problems in general settings:
\begin{itemize}
\item \cite[\S 3.1]{DK20b} establishes uniqueness of general or stationary
solution and existence of stationary solution under milder conditions
than \cite{BK93a}, especially removing the Polish space requirement
(see \cite[Remark 3.18]{DK20b} for comparison).
\item \cite{BK93b} gave a classical framework of establishing finite-dimensional
convergence of not-necessarily-c$\grave{\mbox{a}}$dl$\grave{\mbox{a}}$g
$\{X_{n}\}_{n\in\mathbf{N}}$ to not-necessarily-c$\grave{\mbox{a}}$dl$\grave{\mbox{a}}$g
solution $X$. \cite[\S 3.2.1]{DK20b} extended such convergence results
to non-Polish spaces and milder separability and tightness conditions
(see \cite[Remark 3.24]{DK20b} for comparison). As a byproduct, \cite[\S 3.3]{DK20b}
also give conditions for long-time typical behaviors of martingale
solutions to be stationary solutions.
\item Finite-dimensional convergence of c$\grave{\mbox{a}}$dl$\grave{\mbox{a}}$g
$\{X_{n}\}_{n\in\mathbf{N}}$ to c$\grave{\mbox{a}}$dl$\grave{\mbox{a}}$g
solution $X$ is known to imply their weak convergence on the path
space $D(\mathbf{R}^{+};E)$, which is generally stronger, when $E$
is a Polish space and $\{X_{n}\}_{n\in\mathbf{N}}$ is relatively
compact (see e.g. \cite[\S 4.8]{EK86} and \cite{BK93b}). \cite[\S 3.2.2]{DK20b}
shows that the completeness of $E$ can be exempted and relative compactness
can be reduced to weaker conditions for such implication (see \cite[Remark 3.37]{DK20b}
for comparison). When $E$ is a Tychonoff space, \cite{DK20b} establishes
similar implication by the tightness results of \ref{enu:Theme3}.
\item \cite{KK20} shows existence of solution to an infinite-dimensional
system of stochastic differential equations by constructing a class
of weighted empirical processes and establishing their weak convergence
to c$\grave{\mbox{a}}$dl$\grave{\mbox{a}}$g solutions to the associated
measure-valued martingale problem. They used one of the mild conditions
of modulus of continuity herein to show tightness, which is generally
weaker than the martingale-problem-type of conditions in \cite[\S 4.8]{EK86}
and \cite{BK93b}.
\end{itemize}

For the classical nonlinear filtering problem in general settings:
\begin{itemize}
\item The \textit{Fujisaki-Kallianpur-Kunita} (FKK) and \textit{Duncan-Mortensen-Zakai}
(DMZ) equations are known to be (infinite) equation systems that \textit{uniquely}%
\footnote{This uniqueness has various senses.%
} identify the (measure-valued) normalized and unnormalized filters
respectively (see e.g. \cite{S78}, \cite{KO88} and \cite{BBK95})
for Polish or compact Hausdorff state spaces. Assuming a more general
state space than those of \cite{S78}, \cite{KL08}, \cite{KO88}
and \cite{BBK95}, \cite[\S 2.3]{DK20c} uses the replication method
to establish uniqueness of general solution to the FKK and DMZ equations
systems. This uniqueness result generalizes that of \cite[\S V]{S78}
as it works for any usual filtration without imposing the c$\grave{\mbox{a}}$dl$\grave{\mbox{a}}$g
property of the solution (see \cite[Remarks 31 and 35]{DK20c}).
\item \cite[\S 2.4]{DK20c} extends the stationary filter construction of
\cite{K71} and \cite{BBK00} to non-Polish or non-compact spaces.
With the aid of replication, \cite{DK20c} reduces the problem in
a general setting to the case of \cite{K71} and \cite{BBK00} and
then transforms the known stationary filter back to the original setting.
\end{itemize}

The remainder of this manuscript is organized as follows. Chapter
\ref{chap:Pre} serves a preliminary collection of notations, terminologies
and facts for the three themes of this article. Moreover, this chapter
provides several examples of non-Polish settings of interest to probabilists
and convergence results that goes beyond the regular setting in e.g.
\cite[Chapter 4]{EK86} and \cite{BK93b}. \ref{enu:Theme1} occupies
four chapters: Chapter \ref{chap:Space_Change} develops the space
change method of replication. Chapter \ref{chap:RepFun} focuses on
the replication of function and linear operator. Chapter \ref{chap:RepMeas}
discusses weak convergence on general topological spaces by replication
of measure and weak convergence of replica measures. Chapter \ref{chap:Rep_Proc}
is devoted to the replication of stochastic process and the associated
convergence problems. Chapters \ref{chap:FDDConv} and \ref{chap:Cadlag}
correspond to \ref{enu:Theme2} and \ref{enu:Theme3} respectively.
We provide background content in Appendix \ref{chap:App1} and miscellaneous
results in Appendix \ref{chap:App2} for self-containment and referral
ease, especially for readers' convenience.

\chapter{\label{chap:Pre}Preliminaries}

\numberwithin{equation}{section}

This chapter contains preparatory materials for Chapter \ref{chap:Space_Change}
- \ref{chap:Cadlag}. \S \ref{sec:Basic} - \S \ref{sec:Convention}
introduce our general notations, terminologies and background. Further
background materials are provided in Appendix \ref{chap:App1}. \S
\ref{sec:Example_Space} motivates the replication approach by examples
of ``defective'' settings and results established by space change
in probability theory.

\section{\label{sec:Basic}Basic concepts}

\subsection{\label{sub:Set_Num}Numbers, sets and mappings}

``$\varnothing$'' denotes the empty set. $\mathbf{N}$ denotes
the positive integers and $\mathbf{N}_{0}\circeq\mathbf{N}\cup\{0\}$%
\footnote{``$\circeq$'' means ``being defined by''.%
}. $\mathbf{Q}$ denotes the rational numbers and $\mathbf{Q}^{+}\circeq\{q\in\mathbf{Q}:q\geq0\}$.
$\mathbf{R}$ denotes the real numbers and $\mathbf{R}^{+}\circeq\{x\in\mathbf{R}:x\geq0\}$.
``$\uparrow$'' and ``$\downarrow$'' denote the non-decreasing
and non-increasing convergence of real numbers (including convergence
to $\pm\infty$) respectively.

``$\subset$'' and ``$\supset$'' denote the containment of sets
\textit{including equalities}. Let $A\subset E$ be non-empty sets.
$\aleph(E)$ denotes the cardinality of $E$. $\mathscr{P}_{0}(E)$
denotes all finite \textit{non-empty} subsets of $E$. $A$ is a \textbf{cocountable}
subset of $E$ if $E\backslash A$ is a countable set. Empty, finite
and countably infinite sets are all considered as countable sets%
\footnote{So, ``countable'' is indifferent than ``at most countable''.%
}.

``$\times$'' denotes the Cartesian product of non-empty sets. Let
$\prod_{i\in\mathbf{I}}S_{i}$ denote the Cartesian product of non-empty
$\{S_{i}\}_{i\in\mathbf{I}}$. When $S_{i}=E$ for all $i\in\mathbf{I}$,
$\prod_{i\in\mathbf{I}}S_{i}$ is often denoted by $E^{\mathbf{I}}$,
or by $E^{\infty}$ if $\aleph(\mathbf{I})=\aleph(\mathbf{N})$, or
by $E^{d}$ if $\aleph(\mathbf{I})=d\in\mathbf{N}$. The projection
on $\prod_{i\in\mathbf{I}}S_{i}$ for non-empty sub-index-set $\mathbf{I}_{0}\subset\mathbf{I}$
is defined by
\begin{equation}
\begin{aligned}\mathfrak{p}_{\mathbf{I}_{0}}:\prod_{i\in\mathbf{I}}S_{i} & \longrightarrow\prod_{i\in\mathbf{I}_{0}}S_{i},\\
x\begin{aligned}\end{aligned}
 & \longmapsto\prod_{i\in\mathbf{I}_{0}}\left\{ x(i)\right\} .
\end{aligned}
\label{eq:Projection}
\end{equation}
In particular, $\mathfrak{p}_{j}\circeq\mathfrak{p}_{\{j\}}$ is called
the \textit{one-dimensional projection} on $\Pi_{i\in\mathbf{I}}S_{i}$
for $j\in\mathbf{I}$.

``$\circ$'' denotes the composition of mappings. $\mathbf{1}_{A}$
denotes the indicator function of $A$. For a mapping $f$ defined
on $E$, $f|_{A}$ denotes the restriction of $f$ to $A$. Given
$f_{i}:E\rightarrow S_{i}$ a mapping for each $i\in\mathbf{I}$ and
$\mathcal{D}\circeq\{f_{i}\}_{i\in\mathbf{I}}$, we define
\begin{equation}
\mathcal{D}|_{A}\circeq\left\{ f|_{A}:f\in\mathcal{D}\right\} \label{eq:Map_Class_Restrict}
\end{equation}
and
\begin{equation}
\begin{aligned}\bigotimes\mathcal{D}=\bigotimes_{i\in\mathbf{I}}f_{i}: & E\longrightarrow\prod_{i\in\mathbf{I}}S_{i},\\
 & x\longmapsto\prod_{i\in\mathbf{I}}\left\{ f_{i}(x)\right\} .
\end{aligned}
\label{eq:Prod_Map}
\end{equation}

\subsection{\label{sub:Meas}Measurable space and measure space}

Let $(E,\mathscr{U})$ and $(S,\mathscr{A})$ be measurable spaces
and $A\subset E$ be non-empty. The \textbf{concentration of $\mathscr{U}$
on} $A$ is
\begin{equation}
\mathscr{U}|_{A}\circeq\left\{ B\cap A:B\in\mathscr{U}\right\} ,\label{eq:Concentrated_Sigma_Algebra}
\end{equation}
which is apparently a $\sigma$-algebra on $A$. For a family of mappings
$\mathcal{D}$ from $E$ to $S$, the \textbf{$\sigma$-algebra induced
by} $\mathcal{D}$ is
\begin{equation}
\sigma(\mathcal{D})\circeq\sigma\left(\left\{ f^{-1}(B):B\in\mathscr{A},f\in\mathcal{D}\right\} \right).\label{eq:Baire_Algebra}
\end{equation}

$\delta_{x}$ denotes the \textit{Dirac measure at $x\in E$} so $\delta_{x}(B)=1$
when $B\in\mathscr{U}$ contains $x$.

$\mathfrak{M}^{+}(E,\mathscr{U})$ (resp.%
\footnote{``resp.'' abbreviates ``respectively''.%
} $\mathfrak{P}(E,\mathscr{U})$) denotes \textit{non-trivial} finite
measures (resp. probability measures) on $(E,\mathscr{U})$. Non-triviality
of a measure $\mu$ on $(E,\mathscr{U})$ means $\mu(E)$ is non-zero.
All measures are non-negative and countably additive.

Let $(E,\mathscr{U},\mu)$ be a measure space. $\mathscr{N}(\mu)$
denotes \textit{$\mu$-null subsets} of $E$, i.e. each member of
$\mathscr{N}(\mu)$ has zero \textit{outer measure induced by $\mu$}
(see \cite[p.89]{D02}). Their complements are called $\mu$\textbf{-conull}
sets. If $\mathscr{N}(\mu)\subset\mathscr{U}$, then $\mathscr{U}$
is called \textit{$\mu$-complete} and $(E,\mathscr{U},\mu)$ is called
\textit{complete}.

$A\subset E$ is a \textbf{support}%
\footnote{By our definition, a measure may have more than one supports.%
}\textbf{ of} $\mu$ (or $\mu$ is \textbf{supported on }$A$) if $E\backslash A\in\mathscr{N}(\mu)$.
The \textbf{expansion of }$\nu\in\mathfrak{M}^{+}(A,\mathscr{U}|_{A})$\textbf{
onto $E$} is defined by
\begin{equation}
\nu|^{E}(B)\circeq\nu(A\cap B),\;\forall B\in\mathscr{U}(E).\label{eq:Measure_Extend}
\end{equation}
When $A\in\mathscr{U}$, the \textbf{concentration of $\mu$ on }$A$
is defined by
\begin{equation}
\mu|_{A}(B)\circeq\mu(B),\;\forall B\in\mathscr{U}|_{A}\subset\mathscr{U}.\label{eq:Measure_Restrict}
\end{equation}
The following facts are well-known and we omit the proof for brevity.
\begin{fact}
\label{fact:Meas_Concen_Expan}Let $(E,\mathscr{U})$ be a measurable
space, $A\subset E$ be non-empty, $\mu\in\mathfrak{M}^{+}(E,\mathscr{U})$
and $\nu\in\mathfrak{M}^{+}(A,\mathscr{U}|_{A})$. Then:

\renewcommand{\labelenumi}{(\alph{enumi})}
\begin{enumerate}
\item If $A\in\mathscr{U}$, then (\ref{eq:Measure_Restrict}) well defines
$\mu|_{A}\in\mathfrak{M}^{+}(A,\mathscr{U}|_{A})$. If, in addition,
$\mu\in\mathfrak{P}(E,\mathscr{U})$, then $\mu|_{A}\in\mathfrak{P}(A,\mathscr{U}|_{A})$
precisely when $\mu(A)=1$.
\item (\ref{eq:Measure_Extend}) well defines $\nu|^{E}\in\mathfrak{M}^{+}(E,\mathscr{U})$
and $\nu|^{E}\in\mathfrak{P}(E,\mathscr{U})$ precisely when $\nu(A)=1$.
\item If $A\in\mathscr{U}$, then $\nu=(\nu|^{E})|_{A}$. If, in addition,
$\mu(E\backslash A)=0$, then $\mu=(\mu|_{A})|^{E}$.
\end{enumerate}
\end{fact}
For a measurable mapping $f:E\rightarrow S$, the \textit{push-forward
measure of} \textbf{$\mu$} \textit{by}\textbf{ }$f$ is $\mu\circ f^{-1}\in\mathfrak{M}^{+}(S,\mathscr{A})$
defined by
\begin{equation}
\mu\circ f^{-1}(B)\circeq\mu\left[f^{-1}(B)\right],\;\forall B\in\mathscr{A}.\label{eq:Measure_PushForward}
\end{equation}

Let $\mathscr{V}$ be another $\sigma$-algebra on $E$. If $\mathscr{V}\subset\mathscr{U}$
and $\nu$ is the restriction of $\mu$ to $\mathscr{V}$, then $\nu\in\mathfrak{M}^{+}(E,\mathscr{V})$
is called \textbf{the restriction of} $\mu$ \textbf{to} $\mathscr{V}$,
$\mu$\textbf{ an extension of} $\nu$\textbf{ to} \textbf{$\mathscr{U}$}
and $(E,\mathscr{U},\mu)$ \textbf{an extension of }$(E,\mathscr{V},\nu)$.
$(E,\mathscr{U},\mu)$ is the \textbf{completion of }$(E,\mathscr{V},\nu)$
if: (1) $(E,\mathscr{U},\mu)$ is a complete extension of $(E,\mathscr{V},\nu)$,
and (2) any complete extension $(E,\mathscr{U}^{\prime},\mu^{\prime})$
of $(E,\mathscr{V},\nu)$ is also an extension of $(E,\mathscr{U},\mu)$.

\subsection{\label{sub:Topo}Topological space}

Hereafter, we will often exclude the $\sigma$-algebra (resp. topology)
in a simplified notation of a measurable (resp. topological) space.

Let $E$ be a topological space. $\mathscr{O}(E)$, $\mathscr{C}(E)$,
$\mathscr{K}(E)$, $\mathscr{K}^{\mathbf{m}}(E)$, $\mathscr{K}_{\sigma}(E)$,
$\mathscr{K}_{\sigma}^{\mathbf{m}}(E)$ and $\mathscr{B}(E)\circeq\sigma[\mathscr{O}(E)]$
denote the families of open, closed, \textit{compact} (see p.\pageref{Compact_Set}),
\textit{metrizable }(see p.\pageref{Metrizable}) \textit{compact},
$\sigma$-compact, \textit{$\sigma$-metrizable compact} (i.e. countable
union of metrizble compact) and Borel subsets of $E$, respectively.
Also,
\begin{equation}
\mathscr{O}_{E}(A)\circeq\left\{ O\cap A:O\in\mathscr{O}(E)\right\} \label{eq:O_E_(A)}
\end{equation}
denotes the subspace topology of $A$ induced from $E$, and 
\begin{equation}
\mathscr{B}_{E}(A)\circeq\sigma\left[\mathscr{O}_{E}(A)\right]=\mathscr{B}(E)|_{A}\label{eq:B_E_(A)}
\end{equation}
denotes the subspace Borel $\sigma$-algebra of $A$ induced from
$E$ for non-empty $A\subset E$.

Let $\mathcal{R}$ be a family of \textit{pseudometrics} (see \cite[\S 2.1, Definition, p.26]{D02})
on $E$. The topology induced by $\mathcal{R}$ is generated by the
\textit{topological basis} (see \cite[\S 13, Definition, p.78]{M00})
\begin{equation}
\left\{ \bigcap_{\mathfrak{r}\in\mathcal{R}_{0}}\left\{ y\in E:\mathfrak{r}(x,y)<2^{-p}\right\} :x\in E,p\in\mathbf{N},\mathcal{R}_{0}\in\mathscr{P}_{0}(\mathcal{R})\right\} .\label{eq:Pseudometric_Topo}
\end{equation}
This topology is the metric topology of $(E,\mathfrak{r})$ when $\mathcal{R}$
is the singleton of a metric $\mathfrak{r}$ on $E$. When $(E,\mathfrak{r})$
is a metric space, we define
\begin{equation}
A^{\epsilon}\circeq\left\{ x\in E:\mathfrak{r}(x,y)<\epsilon\mbox{ for some }y\in A\right\} ,\;\forall\epsilon\in(0,\infty).\label{eq:Epsilon_envelope}
\end{equation}

Let $S$ be a topological space and $\mathcal{D}$ be a family of
mappings from $E$ to $S$. The topology generated by the topological
basis
\begin{equation}
\left\{ \bigcap_{f\in\mathcal{D}_{0}}\bigcap_{O\in\mathscr{U}_{0}}f^{-1}(O)\cap A:\mathscr{U}_{0}\in\mathscr{P}_{0}\left[\mathscr{O}(S)\right],\mathcal{D}_{0}\in\mathscr{P}_{0}(\mathcal{D})\right\} \label{eq:O_D_(A)_Topo_Basis}
\end{equation}
is called the \textbf{topology induced by $\mathcal{D}$ on }$A$
and is denoted by $\mathscr{O}_{\mathcal{D}}(A)$. The \textbf{Borel
$\sigma$-algebra induced by $\mathcal{D}$ on }$A$ refers to $\mathscr{B}_{\mathcal{D}}(A)\circeq\sigma[\mathscr{O}_{\mathcal{D}}(A)]$,
which may be different than $\sigma(\mathcal{D})|_{A}$ defined by
(\ref{eq:Concentrated_Sigma_Algebra}) and (\ref{eq:Baire_Algebra})
(see Fact \ref{fact:O_D_(A)_B_D_(A)}).

Let $\mathscr{U}$ be another topology on $E$. If $\mathscr{U}\subset\mathscr{O}(E)$,
then the topological space $(E,\mathscr{U})$ is called a \textbf{topological
coarsening} \textbf{of} $E$, or equivalently, $E$ is called a \textbf{topological
refinement} \textbf{of} $(E,\mathscr{U})$.

Hereafter, by ``$x_{n}\rightarrow x$ as $n\uparrow\infty$ in $E$''
we mean that: (1) $\{x_{n}\}_{n\in\mathbf{N}}$ and $x$ are members
of $E$, and (2) $\{x_{n}\}_{n\in\mathbf{N}}$ converges to $x$ with
respect to the topology of $E$.

\subsection{\label{sub:Morph}Morphisms}

Let $E$ and $S$ be topological spaces. A mapping $f:E\rightarrow S$
is a \textit{homeomorphism} between $E$ and $S$ if $f$ is bijective
and both $f$ and $f^{-1}$ are continuous. $E$ and $S$ are \textit{homeomorphic}
to and \textit{homeomorphs of} each other if a homeomorphism between
them exists. $f$ is an \textit{imbedding} from $E$ into $S$ if
$f$ is a homeomorphism between $E$ and $(f(E),\mathscr{O}_{S}(f(E)))$.
The less common notions of standard Borel subset and Borel isomorphism
are critical for our developments.
\begin{defn}
\label{def:SB}Let $E$ and $S$ be topological spaces%
\footnote{Standard Borel property can be defined for general measurable spaces.
Herein, we focus on the topological space case. There are other common
equivalent definitions to ours, which are discussed in \S \ref{sec:SB}.%
} and $A\subset E$ be non-empty.
\begin{itemize}
\item A mapping $f:E\rightarrow S$ is a \textbf{Borel isomorphism between
$E$ and $S$} if $f$ is bijective and both $f$ and $f^{-1}$ are
measurable with respect to $\mathscr{B}(E)$ and $\mathscr{B}(S)$.
\item $E$ and $S$ are \textbf{Borel isomorphic} \textbf{to }and \textbf{Borel
isomorphs of} \textbf{each other} if there exists a Borel isomorphism
between them.
\item $(A,\mathscr{O}_{E}(A))$ is a \textbf{Borel subspace of $E$} if
$A\in\mathscr{B}(E)$.
\item $E$ is a \textbf{standard Borel space} if $E$ is Borel isomorphic
to a Borel subspace of some \textit{Polish} (see p.\pageref{Polish})
space.
\item $A$ is a \textbf{standard Borel subset of $E$} if $(A,\mathscr{O}_{E}(A))$
is a standard Borel space. $\mathscr{B}^{\mathbf{s}}(E)$ denotes
the family of all standard Borel subsets of $E$.
\end{itemize}
\end{defn}
Standard Borel spaces inherit many nice properties of Borel $\sigma$-algebras
of Polish spaces. A brief review of standard Borel spaces/subsets
is provided in \S \ref{sec:SB}.

\subsection{\label{sub:Prod_Space}Product space}

``$\otimes$'' denotes product of $\sigma$-algebras. Given measurable
spaces $\{(S_{i},\mathscr{A}_{i})\}_{i\in\mathbf{I}}$, the product
$\sigma$-algebra of $\{\mathscr{A}_{i}\}_{i\in\mathbf{I}}$ on $\prod_{i\in\mathbf{I}}S_{i}$
refers to
\begin{equation}
\bigotimes_{i\in\mathbf{I}}\mathscr{A}_{i}\circeq\sigma\left(\{\mathfrak{p}_{i}\}_{i\in\mathbf{I}}\right).\label{eq:Prod_Sigma_Algebra}
\end{equation}
When $(S_{i},\mathscr{A}_{i})=(S,\mathscr{A})$ for all $i\in\mathbf{I}$,
$\bigotimes_{i\in\mathbf{I}}\mathscr{A}_{i}$ is often denoted by
$\mathscr{A}^{\otimes\mathbf{I}}$, or by $\mathscr{A}^{\otimes\infty}$
if $\aleph(\mathbf{I})=\aleph(\mathbf{N})$, or by $\mathscr{A}^{\otimes d}$
if $\aleph(\mathbf{I})=d\in\mathbf{N}$. The following facts are well-known.
\begin{fact}
\label{fact:Prod_Map_1}Let $(E,\mathscr{U})$ and $\{(S_{i},\mathscr{A}_{i})\}_{i\in\mathbf{I}}$
be measurable spaces, $\mathbf{I}_{0}\subset\mathbf{I}$ be non-empty,
$S\circeq\prod_{i\in\mathbf{I}}S_{i}$, $\mathscr{A}\circeq\bigotimes_{i\in\mathbf{I}}\mathscr{A}_{i}$,
$S_{\mathbf{I}_{0}}\circeq\prod_{i\in\mathbf{I}_{0}}S_{i}$, $\mathscr{A}_{\mathbf{I}_{0}}\circeq\bigotimes_{i\in\mathbf{I}_{0}}\mathscr{A}_{i}$
and $f:E\rightarrow S_{i}$ be a mapping for each $i\in\mathbf{I}$.
Then:

\renewcommand{\labelenumi}{(\alph{enumi})}
\begin{enumerate}
\item $\mathfrak{p}_{\mathbf{I}_{0}}$ is a measurable mapping from $(S,\mathscr{A})$
to $(S_{\mathbf{I}_{0}},\mathscr{A}_{\mathbf{I}_{0}})$.
\item $\bigotimes_{i\in\mathbf{I}}f_{i}:(E,\mathscr{U})\rightarrow(S,\mathscr{A})$
is measurable if and only if $f_{i}:(E,\mathscr{U})\rightarrow(S_{i},\mathscr{A}_{i})$
is measurable for all $i\in\mathbf{I}$.
\end{enumerate}
\end{fact}
``$\otimes$'' also denotes the product of topologies. Given topological
spaces $\{S_{i}\}_{i\in\mathbf{I}}$, the product topology of $\{\mathscr{O}(S_{i})\}_{i\in\mathbf{I}}$
on $\prod_{i\in\mathbf{I}}S_{i}$ refers to
\begin{equation}
\bigotimes_{i\in\mathbf{I}}\mathscr{O}(S_{i})\circeq\mathscr{O}_{\{\mathfrak{p}_{i}\}_{i\in\mathbf{I}}}\left(\prod_{i\in\mathbf{I}}S_{i}\right).\label{eq:Prod_Topology}
\end{equation}
When $S_{i}=E$ for all $i\in\mathbf{I}$, $\bigotimes_{i\in\mathbf{I}}\mathscr{O}(S_{i})$
is often denoted by $\mathscr{O}(E)^{\mathbf{I}}$, or by $\mathscr{O}(E)^{\infty}$
if $\aleph(\mathbf{I})=\aleph(\mathbf{N})$, or by $\mathscr{O}(E)^{d}$
if $\aleph(\mathbf{I})=d\in\mathbf{N}$. The following facts are well-known.
\begin{fact}
\label{fact:Prod_Map_2}Let $E$ and $\{(S_{i},\mathscr{A}_{i})\}_{i\in\mathbf{I}}$
be topological spaces, $\mathbf{I}_{0}\subset\mathbf{I}$ be non-empty,
$S\circeq\prod_{i\in\mathbf{I}}S_{i}$, $\mathscr{U}\circeq\bigotimes_{i\in\mathbf{I}}\mathscr{O}(E_{i})$,
$S_{\mathbf{I}_{0}}\circeq\prod_{i\in\mathbf{I}_{0}}S_{i}$, $\mathscr{U}_{\mathbf{I}_{0}}\circeq\bigotimes_{i\in\mathbf{I}_{0}}\mathscr{O}_{i}(E)$
and $f:E\rightarrow S_{i}$ be a mapping for each $i\in\mathbf{I}$.
Then:

\renewcommand{\labelenumi}{(\alph{enumi})}
\begin{enumerate}
\item $\mathfrak{p}_{\mathbf{I}_{0}}$ is a continuous mapping from $(S,\mathscr{U})$
to $(S_{\mathbf{I}_{0}},\mathscr{U}_{\mathbf{I}_{0}})$.
\item $\bigotimes_{i\in\mathbf{I}}f_{i}:E\rightarrow(S,\mathscr{U})$ is
continuous if and only if $f_{i}:E\rightarrow S_{i}$ is continuous
for all $i\in\mathbf{I}$.
\item $\bigotimes_{i\in\mathbf{I}}f_{i}:E\rightarrow(S,\mathscr{U})$ is
continuous at $x\in E$ (see \cite[p.104]{M00}) if and only if $f_{i}:E\rightarrow S_{i}$
is continuous at $x\in E$ for all $i\in\mathbf{I}$.
\end{enumerate}
\end{fact}
Standard discussions about product topological spaces can be found
in e.g. \cite[\S 15 and \S 19]{M00} and \cite[Vol.II, \S 6.4]{B07}.
Herein, we remind the readers of one basic but indispensable fact:
For general topological space $\{S_{i}\}_{i\in\mathbf{I}}$, the Borel
$\sigma$-algebra $\sigma[\bigotimes_{i\in\mathbf{I}}\mathscr{O}(S_{i})]$
generated by their product topology is likely to differ from $\bigotimes_{i\in\mathbf{I}}\mathscr{B}(S_{i})$,
the product of their individual Borel $\sigma$-algebras. Such difference
happens even in the two-dimensional case (see \cite[Vol.II, Example 6.4.3]{B07}).
This is why we use different notations for product $\sigma$-algebra
and product topologies. Avoidance of the difference above needs additional
countability of the product topology $\bigotimes_{i\in\mathbf{I}}\mathscr{O}(S_{i})$
(see Proposition \ref{prop:Prod_Space}).

Hereafter, $\mathbf{R}^{k}$ (with $k\in\mathbf{N}$) denotes the
$k$-dimensional Euclidean space equipped with the usual norm ``$\left|\cdot\right|$''
and the $k$-dimensional Lebesgue measure. \textit{Conull subsets}
of $\mathbf{R}^{k}$ are in the Lebesgue sense. $\left|\cdot\right|$
also denotes the norm metric on $\mathbf{R}^{k}$.

\section{\label{sec:Space_of_Map}Spaces of mappings}

\subsection{\label{sub:Map}General mappings}

Let $\mathbf{I}$, $E$ and $S$ be non-empty sets. The Cartesian
power $E^{\mathbf{I}}$ is the family of all mappings from $\mathbf{I}$
to $E$. When $\mathbf{I}$ has certain index meaning (e.g. time,
order), a member of $E^{\mathbf{I}}$ is often considered as a ``\textit{path
indexed by }$\mathbf{I}$''. So, we define the associated \textbf{path
mapping of} $f\in S^{E}$ by%
\footnote{$\varpi$ is the calligraphical form of the greek letter $\pi$.%
}
\begin{equation}
\varpi_{\mathbf{I}}(f)\circeq\bigotimes_{i\in\mathbf{I}}f\circ\mathfrak{p}_{i}\in(S^{\mathbf{I}})^{E^{\mathbf{I}}}.\label{eq:Path_Mapping}
\end{equation}
This mapping sends every $E$-valued path $x$ indexed by $\mathbf{I}$
to the $S$-valued path $f\circ x$ indexed by $\mathbf{I}$. We define
the associated \textbf{joint path mapping of} $\mathcal{D}\subset S^{E}$
by
\begin{equation}
\varpi_{\mathbf{I}}(\mathcal{D})\circeq\bigotimes\left\{ \varpi_{\mathbf{I}}(f)\in(S^{\mathbf{I}})^{E^{\mathbf{I}}}:f\in\mathcal{D}\right\} \in\left[(S^{\mathbf{I}})^{\mathcal{D}}\right]^{E^{\mathbf{I}}}.\label{eq:Path_Mapping_Class}
\end{equation}
For simplicity, $\varpi_{\mathbf{I}}(f)$ is often denoted by $\varpi(f)$
if $\mathbf{I}=\mathbf{R}^{+}$, or by $\varpi_{T}(f)$ if $\mathbf{I}=[0,T]$,
or by $\varpi_{a,b}(f)$ if $\mathbf{I}=[a,b]$. Similar notations
apply to $\varpi_{\mathbf{I}}(\mathcal{D})$.
\begin{rem}
\label{rem:Path_Mapping}$\varpi_{\mathbf{I}}(\mathcal{D})$ and $\varpi_{\mathbf{I}}(\bigotimes\mathcal{D})$
differ. The latter maps $E^{\mathbf{I}}$ to $(S^{\mathcal{D}})^{\mathbf{I}}$.
\end{rem}

Let $\delta,T\in(0,\infty)$, $[a,b]\subset\mathbf{R}^{+}$ and $\mathfrak{r}$
be a pseudometric on $E$. We define the $\mathfrak{r}$-modulus of
continuity
\begin{equation}
\begin{aligned}w_{\mathfrak{r},\delta,T}^{\prime}(x) & \circeq\inf\left\{ \max_{1\leq i\leq n}\sup_{s,t\in[t_{i-1},t_{i})}\mathfrak{r}\left(x(t),x(s)\right):0\leq t_{0}\right.\\
 & \left.<...<T<t_{n},\inf_{i\leq n}(t_{i}-t_{i-1})>\delta,n\in\mathbf{N}\right\} 
\end{aligned}
\label{eq:w'}
\end{equation}
for each $x\in E^{\mathbf{R}^{+}}$, define
\begin{equation}
\mathfrak{r}_{[a,b]}(x,y)\circeq\sup_{t\in[a,b]}1\wedge\mathfrak{r}\left(x(t),y(t)\right)\label{eq:r[a,b]}
\end{equation}
for each $x,y\in E^{[a,b]}$ or $E^{\mathbf{R}^{+}}$, and let
\begin{equation}
J(x)\circeq\left\{ t\in\mathbf{R}^{+}:x(t)\neq\lim_{s\rightarrow t-}x(s)\in E\right\} \label{eq:J(x)}
\end{equation}
denotes the\textbf{ set of left-jump times} \textbf{of} $x\in E^{\mathbf{R}^{+}}$.

\subsection{\label{sub:Meas_Cont_Cadlag_Map}Measurable, c$\grave{\mbox{a}}$dl$\grave{\mbox{a}}$g
and continuous mappings}

When $E$ and $S$ are measurable spaces, $M(S;E)$ denotes measurable
mappings from $S$ to $E$. When $E$ or $S$ is a topological space,
$M(S;E)$ abbreviates $M(S;E,\mathscr{B}(E))$ or $M(S,\mathscr{B}(S);E)$
respectively. When $E$ and $S$ are both topological spaces, $C(S;E)$,
$\mathbf{hom}(S;E)$, $\mathbf{imb}(S;E)$ and $\mathbf{biso}(S;E)$
denote continuous mappings, homeomorphisms, imbeddings and Borel isomorphisms
from $S$ to $E$, respectively.

$\mathbf{TC}(\mathbf{R}^{+})$ (resp. $\mathbf{TC}([a,b])$) denotes
the family of all \textit{time-changes} \textit{on} $\mathbf{R}^{+}$
(resp. $[a,b]\subset\mathbf{R}^{+}$). That is, each $\lambda\in\mathbf{TC}(\mathbf{R}^{+})$
(resp. $\lambda\in\mathbf{TC}([a,b])$) is a strictly increasing homeomorphism
from $\mathbf{R}^{+}$ (resp. $[a,b]$) to itself and satisfies
\begin{equation}
|||\lambda|||\circeq\sup_{t>s}\left|\ln\frac{\lambda(t)-\lambda(s)}{t-s}\right|<\infty.\label{eq:MonotoneNorm}
\end{equation}
Then, we define
\begin{equation}
\varrho_{[a,b]}^{\mathfrak{r}}(x,y)\circeq\inf_{\lambda\in\mathbf{TC}([a,b])}\left(|||\lambda|||\vee\mathfrak{r}_{[a,b]}(x\circ\lambda,y)\right)\label{eq:SkoMetric_[a,b]}
\end{equation}
for each $x,y\in E^{[a,b]}$, and define
\begin{equation}
\varrho^{\mathfrak{r}}(x,y)\circeq\inf_{\lambda\in\mathbf{TC}(\mathbf{R}^{+})}\left(|||\lambda|||\vee\int_{0}^{\infty}e^{-u}\mathfrak{r}_{[0,u]}(x\circ\lambda,y)du\right)\label{eq:SkoMetric}
\end{equation}
for each $x,y\in E^{\mathbf{R}^{+}}$.

When $E$ is a topological space, $x\in E^{\mathbf{R}^{+}}$ is \textit{c$\grave{\mbox{a}}$dl$\grave{\mbox{a}}$g}
(i.e. right-continuous and left-limited) if for every $t\in\mathbf{R}^{+}$,
there exists a unique $y^{t}\in E$ such that $x(u_{n})\rightarrow y^{t}$
as $n\uparrow\infty$ in $E$ for all $u_{n}\uparrow t$ and $x(v_{n})\rightarrow x(t)$
as $n\uparrow\infty$ in $E$ for all $v_{n}\downarrow t$. When $E$
is a \textit{Tychonoff space}%
\footnote{We use the terminologies ``Tychonoff space'' instead of ``completely
regular space'' since the latter sometimes is used in a non-Hausdorff
context.%
} (see p.\pageref{CR}), Proposition \ref{prop:CR} to follow shows
that the topology of $E$ is induced by a family $\mathcal{R}$ of
pseudometrics on $E$. Then, by $D(\mathbf{R}^{+};E)$ (resp. $D([a,b];E)$)
we denote the space of all c$\grave{\mbox{a}}$dl$\grave{\mbox{a}}$g
members of $E^{\mathbf{R}^{+}}$ (resp. $E^{[a,b]}$) equipped with
the Skorokhod $\mathscr{J}_{1}$-topology $\mathscr{J}(E)$ (resp.
$\mathscr{J}_{a,b}(E)$), that is, the topology induced by pseudometrics
$\{\varrho^{\mathfrak{r}}\}_{\mathfrak{r}\in\mathcal{R}}$ (resp.
$\{\varrho_{[a,b]}^{\mathfrak{r}}\}_{\mathfrak{r}\in\mathcal{R}}$).
$\mathscr{J}(E)$ and $\mathscr{J}_{a,b}(E)$ turn out to be independent
of the choice of the pseudometrics $\mathcal{R}$. More information
about Skorokhod $\mathscr{J}_{1}$-spaces can be found in \S \ref{sec:Sko}.

\subsection{\label{sub:Fun}$\mathbf{R}^{k}$-valued functions}

Let $E$ be a non-empty set and $\{f,g\}\subset\mathbf{R}^{E}$. We
define $f\vee g(x)\circeq\max\{f(x),g(x)\}$, $f\wedge g(x)\circeq\min\{f(x),g(x)\}$,
$f^{+}(x)\circeq\max\{f(x),0\}$ and $f^{-}(x)\circeq\max\{-f(x),0\}$
for all $x\in E$. A subset of $\mathbf{R}^{E}$ is a \textit{function
lattice} if it is closed under ``$\wedge$'' and ``$\vee$''.

Let $k\in\mathbf{N}$ and $\mathcal{D}\subset(\mathbf{R}^{k})^{E}$.
The \textbf{additive expansion} \textbf{of} $\mathcal{D}$ is defined
by
\begin{equation}
\mathfrak{ae}(\mathcal{D})\circeq\mathcal{D}\cup\left\{ f+g:f,g\in\mathcal{D}\right\} ,\label{eq:Add_Extension}
\end{equation}
and the \textbf{additive closure} \textbf{of} $\mathcal{D}$ is defined
by
\begin{equation}
\mathfrak{ac}(\mathcal{D})\circeq\left\{ \sum_{f\in\mathcal{D}_{0}}f:\mathcal{D}_{0}\in\mathscr{P}_{0}(\mathcal{D})\right\} .\label{eq:Add_Closure}
\end{equation}
When $k=1$, the \textbf{multiplicative closure} \textbf{of }$\mathcal{D}$
is defined by
\begin{equation}
\mathfrak{mc}(\mathcal{D})\circeq\left\{ \prod_{f\in\mathcal{D}_{0}}f:\mathcal{D}_{0}\in\mathscr{P}_{0}(\mathcal{D})\right\} ,\label{eq:Mult_Closure}
\end{equation}
the \textbf{$\mathbf{Q}$-algebra generated by} $\mathcal{D}$ is
defined by
\begin{equation}
\mathfrak{ag}_{\mathbf{Q}}(\mathcal{D})\circeq\mathfrak{ac}\left(\left\{ af:f\in\mathfrak{mc}(\mathcal{D}),a\in\mathbf{Q}\right\} \right),\label{eq:ag_Q_(D)}
\end{equation}
the \textbf{algebra generated by} $\mathcal{D}$ is defined by
\begin{equation}
\mathfrak{ag}(\mathcal{D})\circeq\mathfrak{ac}\left(\left\{ af:f\in\mathfrak{mc}(\mathcal{D}),a\in\mathbf{R}\right\} \right),\label{eq:ag(D)}
\end{equation}
and, for a \textit{finite} index set $\mathbf{I}$, we define%
\footnote{The function class in (\ref{eq:Pi^d}) is formed similarly to the
function class in \cite[\S 4.4, (4.15)]{EK86}.%
}
\begin{equation}
\Pi^{\mathbf{I}}(\mathcal{D})\circeq\left\{ g\in\mathbf{R}^{E^{\mathbf{I}}}:g=\prod_{j=1}^{i}f_{j}\circ\mathfrak{p}_{j},f_{j}\in\mathcal{D},1\leq i\leq\aleph(\mathbf{I})\right\} .\label{eq:Pi^d}
\end{equation}
Hereafter, $\Pi^{\mathbf{I}}(\mathcal{D})$ is often denoted by $\Pi^{d}(\mathcal{D})$
with $d\circeq\aleph(\mathbf{I})$. The enlargements of $\mathcal{D}$
above are often used to construct a rich but countable collection
of functions that includes $\mathcal{D}$. Some of their basic properties
are specified in \S \ref{sec:SP_Fun} and \S \ref{sec:Gen_Tech}
- \S \ref{sec:Comp_A1}.

``$\overset{u}{\rightarrow}$'' denotes uniform convergence of $\mathbf{R}^{k}$-valued
functions. When the members of $\mathcal{D}\subset(\mathbf{R}^{k})^{E}$
are bounded%
\footnote{$f$ is bounded if $\Vert f\Vert_{\infty}\in\mathbf{R}^{+}$.%
}, $\mathfrak{cl}(\mathcal{D})$ denotes the closure of $\mathcal{D}$
under the supremum norm $\Vert\cdot\Vert_{\infty}$ and, if $k=1$,
we define
\begin{equation}
\mathfrak{ca}(\mathcal{D})\circeq\mathfrak{cl}\left[\mathfrak{ag}(\mathcal{D})\right]=\mathfrak{cl}\left[\mathfrak{ag}_{\mathbf{Q}}(\mathcal{D})\right].\label{eq:ca(D)}
\end{equation}
The second equality above is immediate by the denseness of $\mathbf{Q}$
in $\mathbf{R}$ and properties of uniform convergence.

$M_{b}(E;\mathbf{R}^{k})$ (resp. $C_{b}(E;\mathbf{R}^{k})$) denotes
the Banach space over scalar field $\mathbf{R}$ of all bounded members
of $M(E;\mathbf{R}^{k})$ (resp. $C(E;\mathbf{R}^{k})$) equipped
with $\Vert\cdot\Vert_{\infty}$. $C_{c}(E;\mathbf{R}^{k})$ denotes
the subspace of all members of $C(E;\mathbf{R}^{k})$ that have \textit{compact
supports}, i.e. the closure of $E\backslash f^{-1}(\{0\})$ is compact.
$C_{0}(E;\mathbf{R}^{k})$ is the subspace of all $f\in C(E;\mathbf{R}^{k})$
such that given $\epsilon>0$, there exists a $K_{\epsilon}\in\mathscr{K}(E)$
satisfying $\Vert f|_{E\backslash K_{\epsilon}}\Vert_{\infty}<\epsilon$%
\footnote{This property is known as vanishing at infinity in the locally compact
case.%
}.

\subsection{\label{sub:SP}Functions and separation of points}

Let $E$ and $A\subset E$ be non-empty sets and $\mathcal{D}\subset\mathbf{R}^{E}$.
\label{SP}$\mathcal{D}$ \textbf{separates points on $A$} if $\bigotimes\mathcal{D}$
is injective, or equivalently, $f(x)=f(y)$ for all $f\in\mathcal{D}$
implies $x=y$ in $A$. Suppose $E$ is a topological space. Then,
\label{SSP}$\mathcal{D}$ \textbf{strongly separates points on $A$}
if $\mathscr{O}_{E}(A)\subset\mathscr{O}_{\mathcal{D}}(A)$. \label{DPC}$\mathcal{D}$
\textbf{determines point convergence on $A$} if $\bigotimes\mathcal{D}(x_{n})\rightarrow\bigotimes\mathcal{D}(x)$
as $n\uparrow\infty$%
\footnote{$\bigotimes\mathcal{D}(x_{n})\rightarrow\bigotimes\mathcal{D}(x)$
as $n\uparrow\infty$ is equivalent to $\lim_{n\rightarrow\infty}f(x_{n})=f(x)$
for all $f\in\mathcal{D}$ by Fact \ref{fact:Seq_Prod_Conv}.%
} in $(\mathbf{R}^{\mathcal{D}},\mathscr{O}(\mathbf{R})^{\mathcal{D}})$
implies $x_{n}\rightarrow x$ as $n\uparrow\infty$ in $(A,\mathscr{O}_{E}(A))$.
\begin{note}
\label{note:SP_Superclass}The point separability, strong point separability
or point convergence determining of $\mathcal{D}\subset\mathbf{R}^{E}$
on $A\subset E$ is inherited by any $\mathcal{D}^{\prime}\subset\mathbf{R}^{E}$
with $\mathcal{D}\subset\mathcal{D}^{\prime}$.
\end{note}
The following are several examples with these point-separation properties.
\begin{example}
\label{exp:SP}$\:$

\renewcommand{\labelenumi}{(\Roman{enumi}) }
\begin{enumerate}
\item Let $C([0,1];\mathbf{R})$ be alternatively equipped with the product
topology $\mathscr{O}(\mathbf{R})^{[0,1]}$. For each $x\in[0,1]$,
the one-dimensional projection $\mathfrak{p}_{x}$ is continuous on
$C([0,1];\mathbf{R})$ since the convergence under product topology
means pointwise convergence (see \cite[Theorem 46.1]{M00}). Note
that
\begin{equation}
\bigotimes_{x\in\mathbf{Q}\cap[0,1]}\mathfrak{p}_{x}(f)=\bigotimes_{x\in\mathbf{Q}\cap[0,1]}\mathfrak{p}_{x}(g)\label{eq:Projection_SP}
\end{equation}
implies $f=g$ by the denseness of $\mathbf{Q}$ in $\mathbf{R}$
and the continuity of $f$ and $g$. Hence, $\{\mathfrak{p}_{x}\}_{x\in\mathbf{Q}\cap[0,1]}$
is a countable collection of $\mathbf{R}$-valued continuous functions
that separates points on $C([0,1];\mathbf{R})$.
\item Let $(E,\mathfrak{r})$ be a metric space and define
\begin{equation}
g_{y,k}(x)\circeq\left[1-k\mathfrak{r}(x,y)\right]\vee0,\;\forall x,y\in E,k\in\mathbf{N}.\label{eq:g_y_k}
\end{equation}
\cite[(4)]{BK10} showed that $\{g_{y,k}\}_{y\in E,k\in\mathbf{N}}$
strongly separates points on $E$. It separates points and determines
point convergence on $E$ by Proposition \ref{prop:Fun_Sep_1} (a,
b) to follow.
\item Let $(E,\mathfrak{r})$ be a metric space. For each $x\in E$,
\begin{equation}
g_{x}(y)\circeq\mathfrak{r}(y,x),\;\forall y\in E\label{eq:g_x}
\end{equation}
is a Lipschitz function by triangular inequality and $g_{x}(y)>0$
for any $x\neq y$ in $E$. Hence, $\mathcal{D}\circeq\{g_{x}\}_{x\in E}$
separates points on $E$. For each $x\in E$, $g_{x}(x_{n})\rightarrow g_{x}(x)$
as $n\uparrow\infty$ in $\mathbf{R}$ implies
\begin{equation}
\lim_{n\rightarrow\infty}\mathfrak{r}(x_{n},x)=\lim_{n\rightarrow\infty}\left|g_{x}(x_{n})\right|=\lim_{n\rightarrow\infty}\left|g_{x}(x_{n})-g_{x}(x)\right|=0\label{eq:Check_Lipschitz_DPC}
\end{equation}
and so $x_{n}\rightarrow x$ as $n\uparrow\infty$ in $E$. Thus,
$\mathcal{D}$ determines point convergence and strongly separates
points on $E$ by (\ref{eq:Check_Lipschitz_DPC}) and Proposition
\ref{prop:Fun_Sep_1} (b) to follow. The family of all Lipschitz functions
on $E$ has the same point-separation properties as $\mathcal{D}$
by Note \ref{note:SP_Superclass}.
\item For each $n\in\mathbf{N}$,
\begin{equation}
f_{n}(x)\circeq\begin{cases}
1, & \mbox{if }x\in\left[0,\frac{1}{n}\right],\\
-\frac{n}{n^{2}-1}x+\frac{n^{2}}{n^{2}-1}, & \mbox{if }x\in\left(\frac{1}{n},n\right),\\
0, & \mbox{if }x\in[n,\infty)
\end{cases}\label{eq:Cc_SP_Example}
\end{equation}
defines a bounded continuous function on $\mathbf{R}^{+}$ which is
strictly decreasing on its compact support $[0,n]$. One immediately
observes that $\mathcal{D}\circeq\{f_{n}\}_{n\in\mathbf{N}}$ separates
points and determines point convergence on $\mathbf{R}^{+}$. $\mathcal{D}$
strongly separates points on $\mathbf{R}^{+}$ by Proposition \ref{prop:Fun_Sep_1}
(b) to follow. $C_{c}(\mathbf{R}^{+};\mathbf{R})$ has the same point-separation
properties as $\mathcal{D}$ by Note \ref{note:SP_Superclass}.
\end{enumerate}
\end{example}
\begin{note}
\label{note:C(E)_SP}$C(E;\mathbf{R})$ and $C_{b}(E;\mathbf{R})$
separate and strongly separate points on $E$ when $E$ is a Tychonoff
space (see Proposition \ref{prop:CR}). For more general $E$, however,
even $C(E;\mathbf{R})$ does not necessarily separate points on $E$
(see Example \ref{exp:P(E)_Hausdorff_2}).
\end{note}
For $\mathcal{D}=\{f_{j}\}_{j\in\mathbf{N}}\subset\mathbf{R}^{E}$
and $d\in\mathbf{N}$, the pseudometric
\begin{equation}
\rho_{\mathcal{D}}(x_{1},x_{2})\circeq\sum_{j=1}^{\infty}2^{-j+1}\left(\left|f_{j}(x_{1})-f_{j}(x_{2})\right|\wedge1\right),\;\forall x_{1},x_{2}\in E\label{eq:TF_Metric}
\end{equation}
on $E$ induces $\mathscr{O}_{\mathcal{D}}(E)$ (see Proposition \ref{prop:Fun_Sep_1}
(d)), and the pseudometric
\begin{equation}
\rho_{\mathcal{D}}^{d}(y_{1},y_{2})\circeq\max_{1\leq i\leq d}\rho_{\mathcal{D}}\left(\mathfrak{p}_{i}(y_{1}),\mathfrak{p}_{i}(y_{2})\right),\;\forall y_{1},y_{2}\in E^{d}\label{eq:TF_Metric^d}
\end{equation}
on $E^{d}$ induces $\mathscr{O}_{\mathcal{D}}(E)^{d}$%
\footnote{As mentioned in \S \ref{sub:Topo} and \S \ref{sub:Prod_Space},
$\mathscr{O}_{\mathcal{D}}(E)^{d}$ is the product topology on $E^{d}$
with $E$ being equipped with the topology $\mathscr{O}_{\mathcal{D}}(E)$.%
} (see Corollary \ref{cor:TF_Metric_FinDim}).

\subsection{\label{sub:OP}Linear Operator Review}

Let $(S,\Vert\cdot\Vert)$ be a Banach space over scalar field $\mathbf{R}$
(like $(C_{b}(E;\mathbf{R}),\Vert\cdot\Vert_{\infty})$). By a\textbf{
single-valued linear operator }$\mathcal{L}$ \textbf{on }$S$ we
refer to a linear subspace $\mathcal{L}\subset S\times S$ such that
for each $f\in S$, $\{g\in S:(f,g)\in\mathcal{L}\}$ is either $\varnothing$
or a singleton $\{\mathcal{L}f\}$. The domain $\mathfrak{D}(\mathcal{L})$
and range $\mathfrak{R}(\mathcal{L})$ of $\mathcal{L}$ are
\begin{equation}
\mathfrak{D}(\mathcal{L})\circeq\left\{ f\in S:\mathcal{L}\cap\left(\{f\}\times S\right)\neq\emptyset\right\} \label{eq:D(L)}
\end{equation}
and
\begin{equation}
\mathfrak{R}(\mathcal{L})\circeq\left\{ g\in S:\mathcal{L}\cap\left(S\times\{g\}\right)\neq\emptyset\right\} ,\label{eq:R(L)}
\end{equation}
which are well-known to be linear subspaces of $S$.

$\mathcal{L}$ is \textbf{closed} if $\mathcal{L}=\mathfrak{cl}(\mathcal{L})$,
where $\mathfrak{cl}(\mathcal{L})$ is the closure of $\mathcal{L}$
in $S\times S$. $\mathfrak{cl}(\mathcal{L})$ is single-valued if
$\mathcal{L}$ is. The \textbf{restriction of} \textbf{$\mathcal{L}$}
\textbf{to }$\mathcal{D}\subset\mathfrak{D}(\mathcal{L})$ is
\begin{equation}
\mathcal{L}|_{\mathcal{D}}\circeq\left\{ (f,\mathcal{L}f)\in\mathcal{L}:f\in\mathcal{D}\right\} .\label{eq:L|D}
\end{equation}
$\mathcal{L}$ is \textit{dissipative} if
\begin{equation}
\beta\Vert f\Vert\leq\Vert\beta f-\mathcal{L}f\Vert,\;\forall f\in\mathfrak{D}(\mathcal{L}),\beta\in(0,\infty).\label{eq:Dissipative}
\end{equation}
$\mathcal{L}$ satisfies \textit{positive maximum principle} if $\sup_{x\in E}f(x)=f(x_{0})\geq0$
implies $\mathcal{L}f(x_{0})\leq0$ for all $f\in\mathfrak{D}(\mathcal{L})$.
$\mathcal{L}$ is a\textbf{ strong} \textbf{generator} \textbf{on}
$S$ if $\mathfrak{cl}(\mathcal{L})$ is the \textit{infinitesimal
generator} (see \cite[p.231]{Y80}) of a \textit{strongly continuous
contraction semigroup} on $S$. When $E$ is a \textit{locally compact}
(see p.\pageref{Local_Compact}) \textit{Hausdorff} (see p.\pageref{Hausdorff})
space and $S=(C_{0}(E;\mathbf{R}),\Vert\cdot\Vert_{\infty})$, $\mathcal{L}$
is a \textbf{Feller generator on} $S$ if $\mathfrak{cl}(\mathcal{L})$
is the infinitesimal generator of a \textit{Feller semigroup} on $S$.
There are multiple Feller semigroup definitions in the literature.
We follow \cite[\S 4.2, p.166]{EK86} and define Feller semigroups
to be strongly continuous, positive, contraction semigroup with conservative
infinitesimal generators.

More details about operators on Banach spaces can be found in standard
texts like \cite[Chapter VIII and Chapter IX]{Y80} and \cite[Chapter 1]{EK86}.

\section{\label{sec:Borel_Measure}Spaces of finite Borel measures}

Let $E$ be a topological space. The members of $\mathfrak{M}^{+}(E,\mathscr{B}(E))$
and $\mathfrak{P}(E,\mathscr{B}(E))$ are called \textit{finite Borel
measures} and \textit{Borel probability measures} respectively.
\begin{defn}
\label{def:Borel_Ext}Let $E$ be a topological space and $\mathscr{U}$
be a sub-$\sigma$-algebra of $\mathscr{B}(E)$. Then, any extension
of $\mu\in\mathfrak{M}^{+}(E,\mathscr{U})$ to $\mathscr{B}(E)$ is
said to be a \textbf{Borel extension of }$\mu$.
\end{defn}
Hereafter, $\mathfrak{be}(\mu)$ denotes the family of all Borel extension(s)
of $\mu$ (if any). If $\mu^{\prime}$ is the unique member of $\mathfrak{be}(\mu)$,
then we specially denote $\mu^{\prime}=\mathfrak{be}(\mu)$. Any $\mu^{\prime}\in\mathfrak{be}(\mu)$
has the same total mass as $\mu$ since $E\in\mathscr{U}\subset\mathscr{B}(E)$.

$\mathcal{M}^{+}(E)$ denotes the space of all finite Borel measures
on $E$ equipped with the weak topology. $\mathcal{P}(E)$ denotes
the subspace of all probability members of $\mathcal{M}^{+}(E)$.
To be specific, the weak topology of $\mathcal{M}^{+}(E)$ is defined
by
\begin{equation}
\mathscr{O}\left[\mathcal{M}^{+}(E)\right]\circeq\mathscr{O}_{C_{b}(E;\mathbf{R})^{*}}\left[\mathcal{M}^{+}(E)\right],\label{eq:Weak_Topo_M(E)}
\end{equation}
where $f^{*}$ denotes the linear functional
\begin{equation}
\begin{aligned}f^{*}: & \mathcal{M}^{+}(E)\longrightarrow\mathbf{R},\\
 & \mu\longmapsto\int_{E}f(x)\mu(dx)
\end{aligned}
\label{eq:f_star}
\end{equation}
for each $f\in M_{b}(E;\mathbf{R})$, and
\begin{equation}
\mathcal{D}^{*}\circeq\left\{ f^{*}:f\in\mathcal{D}\right\} \label{eq:D_Star}
\end{equation}
for each $\mathcal{D}\subset M_{b}(E;\mathbf{R})$. Given $f\in M_{b}(E;\mathbf{R}^{k})$
(with $k\in\mathbf{N}$), we define $f^{*}\circeq\bigotimes_{i=1}^{k}(\mathfrak{p}_{i}\circ f)^{*}$.
\begin{rem}
\label{rem:Narrow_Topo}The weak topology of $\mathcal{M}^{+}(E)$,
sometimes called the ``narrow topology'', generally \textit{differs
from} the standard \textit{weak-$*$ topology induced by dual space}.
They become equal when $E$ is a locally compact Hausdorff space (see
\cite[Chapter II, \S 6.5 - 6.7]{M95}). Hereafter, we use ``$f^{*}$''
to denote the linear functional in (\ref{eq:f_star}). $\mathcal{D}^{*}$
is not necessarily a dual space of $\mathcal{D}$ even when $\mathcal{D}$
is a linear space.
\end{rem}

Weak convergence is one of the central interests of this work. As
specified in \S \ref{sub:Topo} for general topological spaces, the
statement
\begin{equation}
\mu_{n}\Longrightarrow\mu\mbox{ as }n\uparrow\infty\mbox{ in }\mathcal{M}^{+}(E)\label{eq:Mu_n_WC_Mu_M(E)}
\end{equation}
means that: (1) $\{\mu_{n}\}_{n\in\mathbf{N}}\cup\{\mu\}$ are members
of $\mathcal{M}^{+}(E)$, and (2) $\{\mu_{n}\}_{n\in\mathbf{N}}$
\textbf{converges to $\mu$ with respect to the weak topology of }$\mathcal{M}^{+}(E)$
(converge weakly to $\mu$ for short). Similar terminology and notation
apply to $\mathcal{P}(E)$.

$\mu\in\mathcal{M}^{+}(E)$ is a \textbf{weak limit point of $\Gamma\subset\mathcal{M}^{+}(E)$
}if there exist $\{\mu_{n}\}_{n\in\mathbf{N}}\subset\Gamma$ satisfying
(\ref{eq:Mu_n_WC_Mu_M(E)}). We denote the \textbf{weak limit }$\mu$\textbf{
of $\{\mu_{n}\}_{n\in\mathbf{N}}\subset\mathcal{M}^{+}(E)$} by
\begin{equation}
\mbox{w-}\lim_{n\rightarrow\infty}\mu_{n}=\mu\label{eq:WL}
\end{equation}
if (\ref{eq:Mu_n_WC_Mu_M(E)}) holds and $\mu$ is the \textit{unique}
weak limit point of $\{\mu_{n}\}_{n\in\mathbf{N}}$. \label{RC_Meas}$\Gamma\subset\mathcal{M}^{+}(E)$
is \textbf{relatively compact} if any infinite subset of $\Gamma$
has a weak limit point in $\mathcal{M}^{+}(E)$%
\footnote{This weak limit point need not belong to $\Gamma$.%
}.
\begin{rem}
\label{rem:WL}(\ref{eq:Mu_n_WC_Mu_M(E)}) does not necessarily imply
(\ref{eq:WL}) since $\mathcal{M}^{+}(E)$ in general is not guaranteed
to be a Hausdorff space so $\mu$ in (\ref{eq:Mu_n_WC_Mu_M(E)}) might
not be unique.
\end{rem}

\subsection{\label{sub:Weak_Topo_Priori}Weak topology in sequential view}

In this work, we consider weak convergence as the \textbf{topological
convergence} of the weak topology defined by (\ref{eq:Weak_Topo_M(E)}).
However, (\ref{eq:Mu_n_WC_Mu_M(E)}) is equivalent to the integral
test
\begin{equation}
\lim_{n\rightarrow\infty}\int_{E}f(x)\mu_{n}(dx)=\int_{E}f(x)\mu(dx),\;\forall f\in C_{b}(E;\mathbf{R}).\label{eq:Int_Test_WC}
\end{equation}
One often defines weak convergence by (\ref{eq:Int_Test_WC}) first
and then defines a topology sequentially by weak convergence, which
is called the \textit{topology of weak convergence} herein. We look
at the relationship between the weak topology and the topology of
weak convergence.

Every topology determines a sense of topological convergence. Conversely,
convergence is definable without any topological structure. Weak convergence
defined by (\ref{eq:Int_Test_WC}) is one example since it does not
formally involve any topology of $\mathcal{M}^{+}(E)$. Two other
examples are almost-sure convergence and \textit{bounded pointwise
convergence} (see \cite[\S 3.4]{EK86}). Suppose a sense of convergence
is given, which is called \textbf{convergence a priori}. Then, one
defines closedness of a set to be the containment of limits of all
convergent a priori sequences. A topology is generated by these ``sequentially
closed'' sets, which is called the \textbf{topology of convergence
a priori}.
\begin{rem}
\label{rem:Conv_Posteriori}The topology of convergence a priori has
its own topological convergence, which is called \textbf{convergnce
a posteriori}. Convergence a priori and convergence a posteriori need
not be the same in general (see \cite[\S 5]{J12} for further details).
Examples of convergence a priori that is stronger than convergence
a posteriori include almost-sure convergence, bounded pointwise convergence
and the $\mathcal{S}$-convergence introduced by \cite{J97b} (see
\S \ref{sub:S-Topo} for a short glance).
\end{rem}

When convergence a priori is the topological convergence of some given
topology, the topology of convergence a priori and the original topology
are not necessarily the same. One example is the weak topology and
the topology of weak convergence. A weak limit point is always a \label{Limit_Point}\textit{limit
point} (see \cite[p.97]{M00}) with respect to the weak topology (see
Proposition \ref{prop:WLP_RC_Metrizable} (b)). In other words, the
weak topology has no more closed sets and hence is coarser than the
topology of weak convergence. These two could be strictly different
as in the following example.
\begin{example}
\label{exp:Weak_Topo_Sequential}\cite[\S 21, Example 2]{M00} explained
that the product space $E\circeq(\mathbf{R}^{[0,1]},\mathscr{O}(\mathbf{R})^{[0,1]})$
is not \textit{first-countable} (see p.\pageref{First_Countable}).
More specifically, the constant function $0$ lies in the closure
of the set
\begin{equation}
A\circeq\left\{ x\in\mathbf{R}^{[0,1]}:x(i)=\begin{cases}
0, & \mbox{if }i\in\mathbf{I},\\
1, & \mbox{otherwise},
\end{cases}\mbox{ for some }\mathbf{I}\in\mathscr{P}_{0}([0,1])\right\} ,\label{eq:Weak_Topo_Sequential_1}
\end{equation}
but no sequence in $A$ converges to $0$. Let $\Gamma\circeq\{\delta_{x}:x\in A\}$,
$\Gamma_{1}$ be the closure of $\Gamma$ with respect to the weak
topology of $\mathcal{M}^{+}(E)$, and $\Gamma_{2}$ be the closure
of $\Gamma$ with resepct to the sequential topology induced by weak
convergence. We show $\Gamma_{1}\neq\Gamma_{2}$ by verifying $\delta_{0}\in\Gamma_{1}\backslash\Gamma_{2}$.
For any $\epsilon\in(0,\infty)$, $m\in\mathbf{N}$ and $\{f_{1},...,f_{m}\}\subset C_{b}(E;\mathbf{R})$,
there exists an $O_{\epsilon}\in\mathscr{O}(\mathbf{R})^{[0,1]}$
such that $0\in O_{\epsilon}$ and
\begin{equation}
\max_{1\leq i\leq m}\left|f_{i}^{*}(\delta_{x})-f_{i}^{*}(\delta_{0})\right|=\max_{1\leq i\leq m}\left|f_{i}(x)-f_{i}(0)\right|<\epsilon,\;\forall x\in O_{\epsilon}.\label{eq:Weak_Topo_Sequential_2}
\end{equation}
Since $0$ lies in the closure of $A$, there exists some $x_{\epsilon}\in(A\cap O_{\epsilon})\backslash\{0\}$
such that $\delta_{x_{\epsilon}}\in\Gamma$ and
\begin{equation}
\max_{1\leq i\leq m}\left|f_{i}^{*}(\delta_{x_{\epsilon}})-f_{i}^{*}(\delta_{0})\right|<\epsilon.\label{eq:Weak_Topo_Sequential_3}
\end{equation}
This proves $\delta_{0}$ is a limit point of $\Gamma$ with respect
to the weak topology, so $\delta_{0}\in\Gamma_{1}$. Meanwhile, $E$
is a Tychonoff space by Proposition \ref{prop:CR_Space} (c). According
to Lemma \ref{lem:CR_Dirac_Meas} (a, b), no sequence in $\Gamma$
may converge weakly to $\delta_{0}$, since no sequence in $A$ converges
to $0$. Thus, $\delta_{0}\notin\Gamma_{2}$.
\end{example}

Relative compactness of finite Borel meaures is defined in terms of
weak convergence (see (\ref{eq:Int_Test_WC})) not the weak topology.
Consequently, relative compactness of $\Gamma\subset\mathcal{M}^{+}(E)$
can be different than: (1) $\Gamma$ having a compact closure in $\mathcal{M}^{+}(E)$
with the weak topology, the usual interpretation of ``relative compactness''
(see p.\pageref{RC}) of $\Gamma$, or (2) $\Gamma$ having a \textit{limit
point compact} (see p.\pageref{LP_Compact}) closure with the weak
topology, or (3) $\Gamma$ having a \textit{sequentially compact}
(see p.\pageref{Seq_Compact}) closure with the weak topology.

When $E$ is a metrizable space, however, $\mathcal{M}^{+}(E)$ is
a metrizable space by Proposition \ref{prop:WLP_RC_Metrizable}),
the weak topology is the same as the topology of weak convergence
by Fact \ref{fact:First_Countable} (b), and the aforementioned ambiguities
about weak limit point and relative compactness will not occur.

\subsection{\label{sub:Sep_Meas}Separation of finite Borel measures by functions}

The measure-separation properties of $\mathcal{D}\subset M_{b}(E;\mathbf{R})$
(i.e. point-separation properties of $\mathcal{D}^{*}$) are vital
for studying weak convergence and $\mathcal{M}^{+}(E)$-valued or
$\mathcal{P}(E)$-valued processes (e.g. filters, measure-valued diffusions,
non-Markov branching particle systems). The following terminologies
are adapted from \cite[\S 3.4]{EK86}: \label{Separating}$\mathcal{D}\subset M_{b}(E;\mathbf{R})$
is \textbf{separating }or \label{CD}\textbf{convergence determining}
\textbf{on $E$} if $\mathcal{D}^{*}$ separates points or determines
point convergence on $\mathcal{M}^{+}(E)$ respectively.
\begin{note}
\label{note:CD_Sep}$C_{b}(E;\mathbf{R})^{*}$ by definition strongly
separates points and so determines point convergence on $\mathcal{M}^{+}(E)$
by Proposition \ref{prop:Fun_Sep_1} (b). Hence, $C_{b}(E;\mathbf{R})$
is convergence determining on $\mathcal{M}^{+}(E)$.
\end{note}

More details on the point-separation and other topological properties
of $\mathcal{M}^{+}(E)$ and $\mathcal{P}(E)$ can be found in \S
\ref{sec:Meas_Sep}, \cite[Part II]{T70} and \cite[Vol. II, Chapter 8]{B07}.

\subsection{\label{sub:Portmanteau}Portmanteau's Theorem}

One way of establishing (\ref{eq:Mu_n_WC_Mu_M(E)}) is establishing
$\lim_{n\rightarrow\infty}f^{*}(\mu_{n})=f^{*}(\mu)$ for all $f$
from a convergence determining collection. A useful alternative is
the Portamenteau's Theorem. This useful tool was commonly established
on metric spaces (see \cite[Lamma 2.2.2]{KX95}). \cite[p.XII and p.40 - 41]{T70}
gave the following partial generalization to Hausdorff spaces.
\begin{thm}
[\textbf{Portmanteau's Theorem}\textrm{, \cite[Theorem 8.1]{T70}}]\label{thm:Portamenteau}Let
$E$ be a Hausdorff space. Consider the following statements:

\renewcommand{\labelenumi}{(\alph{enumi})}
\begin{enumerate}
\item (\ref{eq:Mu_n_WC_Mu_M(E)}) holds.
\item $\limsup_{n\rightarrow\infty}\mu_{n}(F)\leq\mu(F)$ for all $F\in\mathscr{C}(E)$.
\item $\liminf_{n\rightarrow\infty}\mu_{n}(O)\geq\mu(O)$ for all $O\in\mathscr{O}(E)$.
\end{enumerate}
Then, (b) and (c) are equivalent and each of them implies (a). If,
in addition, $E$ is a Tychonoff space, then (a) - (c) are equivalent.
\end{thm}

\subsection{\label{sub:Tight}Tightness}

Tightness is often more easily verified than relative compactness.
Compact subsets are not necessarily Borel subsets in non-Hausdorff
spaces. At the same time, they can be in the domain of possibly non-Borel
measures (see \S \ref{sub:Metrizable_Compact}). So, we slightly
adjust the ordinary definition of and extend tightness to general
finite measures.
\begin{defn}
\label{def:Tight}Let $(E,\mathscr{U})$ be a measurable space, $S$
be a topological space and $\mathscr{A}$ be a $\sigma$-algebra on
$S$.
\begin{itemize}
\item When $S\subset E$, $\Gamma\subset\mathfrak{M}^{+}(E,\mathscr{U})$
is \textbf{tight in }$S$ (resp. $\mathbf{m}$\textbf{-tight in }$S$)
if for any $\epsilon\in(0,\infty)$, there exists a $K_{\epsilon}\in\mathscr{K}(S)$
(resp. $K_{\epsilon}\in\mathscr{K}^{\mathbf{m}}(S)$) such that $K_{\epsilon}\in\mathscr{U}$
and $\sup_{\mu\in\Gamma}\mu(E\backslash K_{\epsilon})\leq\epsilon$.
\item $\Gamma\subset\mathfrak{M}^{+}(S,\mathscr{A})$ is \textbf{tight in
}$A\subset S$ (resp. $\mathbf{m}$\textbf{-tight in }$A$) if $A$
is non-empty and $\Gamma$ is tight (resp. $\mathbf{m}$-tight) in
$(A,\mathscr{O}_{S}(A))$.
\item $\Gamma\subset\mathfrak{M}^{+}(S,\mathscr{A})$ is \textbf{tight }(resp.
$\mathbf{m}$\textbf{-tight}) if it is tight (resp. $\mathbf{m}$-tight)
in $S$.
\end{itemize}
\end{defn}
\begin{note}
\label{note:Singleton_Tight}Tightness of a measure $\mu$ refers
to that of the singleton $\{\mu\}$.\end{note}
\begin{rem}
\label{rem:Tight}$\mathbf{m}$-tightness is stronger than tightness,
and they are the same if every compact subset of the underlying space
is metrizable. We refer the readers to \S \ref{sub:Metrizable_Compact}
for specific discussion about metrizable compact subsets.
\end{rem}

The classical Ulam's Theorem (see \cite[Theorem 1.4]{B68}), showing
tightness of every finite set of Borel probability measures on a Polish
space $E$, has the following stronger form about $\mathbf{m}$-tightness.
\begin{thm}
[\textbf{Ulam's Theorem}\textrm{, \cite[Vol.II, Theorem 7.4.3]{B07}}]\label{thm:Ulam_m-Tight}If
$E$ is a Souslin space (see p.\pageref{Souslin}), especially if
$E$ is a Lusin (see p.\pageref{Lusin}) or Polish space, then any
finite subset of $\mathcal{M}^{+}(E)$ is $\mathbf{m}$-tight.
\end{thm}
The Prokhorov's Theorem is a fundamental result connecting relative
compactness and tightness of finite Borel measures. Part (a) below,
adapted from \cite[Vol.II, Theorem 8.6.2]{B07}, gives one direction
of Prokhorov's Theorem. Part (b), extending the other direction from
Polish to Hausdorff spaces, is adapted from \cite[Theorem 2.2.1]{KX95}.
\begin{thm}
[\textbf{Prokhorov's Theorem}]\label{thm:Prokhorov}$\,$

\renewcommand{\labelenumi}{(\alph{enumi})}
\begin{enumerate}
\item If $E$ is a Polish space, then relative compactness implies tightness
for any subset of $\mathcal{M}^{+}(E)$.
\item If $E$ is a Hausdorff space, then tightness implies relative compactness
for any subset of $\mathcal{P}(E)$.
\end{enumerate}
\end{thm}

\subsection{\label{sub:Sko_Meas}Finite Borel measures on $D(\mathbf{R}^{+};E)$}

When $E$ is a Tychonoff space, the Skorokhod $\mathscr{J}_{1}$-space
$D(\mathbf{R}^{+};E)$ always satisfies
\begin{equation}
\mathscr{B}\left[D(\mathbf{R}^{+};E)\right]=\sigma\left[\mathscr{J}(E)\right]\supset\left.\mathscr{B}(E)^{\otimes\mathbf{R}^{+}}\right|_{D(\mathbf{R}^{+};E)}\label{eq:Sko_Borel>Prod}
\end{equation}
and
\begin{equation}
\mathcal{M}^{+}\left[D(\mathbf{R}^{+};E)\right]\subset\mathfrak{M}^{+}\left(D(\mathbf{R}^{+};E),\left.\mathscr{B}(E)^{\otimes\mathbf{R}^{+}}\right|_{D(\mathbf{R}^{+};E)}\right).\label{eq:Sko_Meas}
\end{equation}
However, equality in (\ref{eq:Sko_Borel>Prod}) or (\ref{eq:Sko_Meas})
may not hold in general. The \textbf{set of fixed left-jump times
of} $\mu\in\mathfrak{M}^{+}(D(\mathbf{R}^{+};E),\mathscr{B}(E)^{\otimes\mathbf{R}^{+}}|_{D(\mathbf{R}^{+};E)})$
refers to
\begin{equation}
J(\mu)\circeq\left\{ t\in\mathbf{R}^{+}:\mu\left(\left\{ x\in D(\mathbf{R}^{+};E):t\in J(x)\right\} \right)>0\right\} \label{eq:J(Mu)}
\end{equation}
if it is well-defined. (\ref{eq:Sko_Borel>Prod}), (\ref{eq:Sko_Meas})
and (\ref{eq:J(Mu)}) are further discussed in \S \ref{sec:Sko}.

\section{\label{sec:RV}Random variable}

Let $(\Omega,\mathscr{F},\mathbb{P})$ be a probability space (i.e.
$\mathbb{P}\in\mathfrak{P}(\Omega,\mathscr{F})$), $(E,\mathscr{U})$
be a measurable space and $S$ be a topological space. Any $X\in M(\Omega,\mathscr{F};E,\mathscr{U})$
is said to be an $(E,\mathscr{U})$-valued random variable. $\mathbb{P}\circ X^{-1}\in\mathfrak{P}(E,\mathscr{U})$,
the push-forward measure of $\mathbb{P}$ by $X$, is called the \textit{distribution
of} $X$. $S$-valued random variables refer to $(S,\mathscr{B}(S))$-valued
random variables if not otherwise specified. Hereafter, we let $(\Omega,\mathscr{F},\mathbb{P};X)$
denote a random variable $X$ defined on probability space $(\Omega,\mathscr{F},\mathbb{P})$.

Any type of tightness in Definition \ref{def:Tight} is defined for
random variables by referring to the corresponding property of their
distributions. ``$X_{n}\Rightarrow X$ as $n\uparrow\infty$ on $S$''
means the distributions of $S$-valued random variables $\{X_{n}\}_{n\in\mathbf{N}}$
converge weakly%
\footnote{``converging weakly'' was defined in \S \ref{sec:Borel_Measure}
as converging with respect to the weak topology. %
} to that of $S$-valued random variable $X$ as $n\uparrow\infty$
in $\mathcal{P}(S)$. Similar interpretations apply to the statements
``$X$ is a the weak limit of $\{X_{n}\}_{n\in\mathbf{N}}$ on $S$'',
``$X$ is a weak limit point of $\{X_{i}\}_{i\in\mathbf{I}}$ on
$S$'' and ``$\{X_{i}\}_{i\in\mathbf{I}}$ is relatively compact
in $S$''.

\section{\label{sec:Proc}Stochastic process}

The stochastic processes treated in this work are indexed by time
horizon $\mathbf{R}^{+}$ and take values in topological spaces%
\footnote{Stochastic processes are definable on measurable spaces without any
topological structure.%
}. Throughout this section, we let $(\Omega,\mathscr{F},\mathbb{P})$
be a probability space, $E$ be a topological space and $X\in(E^{\mathbf{R}^{+}})^{\Omega}$.

\subsection{\label{sub:Proc_Def}Definition}

$X$ is an \textbf{$E$-valued} (\textbf{stochastic}) \textbf{process}
if $\mathscr{B}(E)^{\otimes\mathbf{R}^{+}}$ is a sub-$\sigma$-algebra
of
\begin{equation}
\mathscr{U}_{X}\circeq\left\{ B\subset E^{\mathbf{R}^{+}}:X^{-1}(B)\in\mathscr{F}\right\} ,\label{eq:Push_Forward_Algebra}
\end{equation}
or equivalently,
\begin{equation}
X\in M(\Omega,\mathscr{F};E^{\mathbf{R}^{+}},\mathscr{B}(E)^{\otimes\mathbf{R}^{+}}).\label{eq:Process_Def}
\end{equation}

\begin{rem}
\label{rem:Push_Forward_Algebra}The $\mathscr{U}_{X}$ in (\ref{eq:Push_Forward_Algebra})
is often called the ``push-forward $\sigma$-algebra of $X$''.
In any case, $X\in M(\Omega,\mathscr{F};E^{\mathbf{R}^{+}},\mathscr{U}_{X})$.
\end{rem}

Let $(\Omega,\mathscr{F},\mathbb{P};X)$ be an $E$-valued process%
\footnote{An $E$-valued process is an $(E^{\mathbf{R}^{+}},\mathscr{B}(E)^{\otimes\mathbf{R}^{+}})$-valued
random variable, hence it is consistent for $(\Omega,\mathscr{F},\mathbb{P};X)$
to denote an $E$-valued process $X$ defined on $(\Omega,\mathscr{F},\mathbb{P})$.%
}. For each $\omega\in\Omega$, $X(\omega)\in E^{\mathbf{R}^{+}}$
is called a (realization) \textbf{path} \textbf{of} $X$. The \textbf{process
distribution of} an $E$-valued process $X$ refers to the push-forward
measure of $\mathbb{P}$ by $X:(\Omega,\mathscr{F})\rightarrow(E^{\mathbf{R}^{+}},\mathscr{U}_{X})$
and is denoted by $\mathrm{pd}(X)\in\mathfrak{P}(E^{\mathbf{R}^{+}},\mathscr{U}_{X})$.
For each $t\in\mathbf{R}^{+}$, $X_{t}\circeq\mathfrak{p}_{t}\circ X$
denotes the (one-dimensional) \textbf{section of} $X$ \textbf{for}
$t$. For each $\mathbf{T}_{0}\in\mathscr{P}_{0}(\mathbf{R}^{+})$,
the \textbf{section of} $X$ \textbf{for} $\mathbf{T}_{0}$ refers
to the $E^{\mathbf{T}_{0}}$-valued mapping $X_{\mathbf{T}_{0}}\circeq\mathfrak{p}_{\mathbf{T}_{0}}\circ X$,
and the \textbf{finite-dimensional distribution of} \textbf{$X$ for
}$\mathbf{T}_{0}$ refers to $\mathrm{pd}(X)\circ\mathfrak{p}_{\mathbf{T}_{0}}^{-1}$.
From Fact \ref{fact:Prod_Map_1} we immediately observe that:
\begin{fact}
\label{fact:Proc_Basic_1}Let $(\Omega,\mathscr{F},\mathbb{P})$ be
a probability space and $E$ be a topological space. Then:

\renewcommand{\labelenumi}{(\alph{enumi})}
\begin{enumerate}
\item $X\in(E^{\mathbf{R}^{+}})^{\Omega}$ is an $E$-valued process if
and only if $\mathfrak{p}_{t}\circ X\in M(\Omega,\mathscr{F};E)$
for all $t\in\mathbf{R}^{+}$.
\item If $\zeta^{t}\in M(\Omega,\mathscr{F};E)$ for all $t\in\mathbf{R}^{+}$,
then
\begin{equation}
X(\omega)(t)\circeq\zeta^{t}(\omega),\;\forall t\in\mathbf{R}^{+},\omega\in\Omega\label{eq:Define_Proc_Marginal}
\end{equation}
well defines an $E$-valued process $X$ satisfying $\zeta^{t}=\mathfrak{p}_{t}\circ X$
for all $t\in\mathbf{R}^{+}$.
\item The section of an $E$-valued process $(\Omega,\mathscr{F},\mathbb{P};X)$
for each $\mathbf{T}_{0}\in\mathscr{P}_{0}(\mathbf{R}^{+})$ is a
member of $M(\Omega,\mathscr{F};E^{\mathbf{T}_{0}},\mathscr{B}(E)^{\otimes\mathbf{T}_{0}})$.
\item The finite-dimensional distribution of an $E$-valued process $X$
for each $\mathbf{T}_{0}\in\mathscr{P}_{0}(\mathbf{R}^{+})$ is the
distribution of $X_{\mathbf{T}_{0}}$ and belongs to $\mathfrak{P}(E^{\mathbf{T}_{0}},\mathscr{B}(E)^{\otimes\mathbf{T}_{0}})$.
In particular, the distribution of $X_{t}$ is a member of $\mathcal{P}(E)$
for all $t\in\mathbf{R}^{+}$.
\end{enumerate}
\end{fact}
\begin{rem}
\label{rem:FDD_Not_Borel}Given an $E$-valued process $X$ and a
general $\mathbf{T}_{0}\in\mathscr{P}_{0}(\mathbf{R}^{+})$, $X$
(resp. $X_{\mathbf{T}_{0}}$) need not be an $(E^{\mathbf{R}^{+}},\mathscr{B}(E^{\mathbf{R}^{+}}))$-valued
(resp. $(E^{\mathbf{T}_{0}},\mathscr{B}(E^{\mathbf{T}_{0}}))$-valued)
random variable, nor is the process distribution of $X$ (resp. the
finite-dimensional distribution of $X$ for $\mathbf{T}_{0}$) necessarily
a Borel measure. This is due to the possible difference between the
Borel $\sigma$-algebra generated by product topology and product
of Borel $\sigma$-algebras on each individual dimension, which was
mentioned in \S \ref{sub:Prod_Space}.
\end{rem}

Hereafter, an $E$-valued process $X$ is also denoted by $X=\{X_{t}\}_{t\geq0}$
or just by $\{X_{t}\}_{t\geq0}$; its section $X_{\mathbf{T}_{0}}$
for $\mathbf{T}_{0}=\{t_{1},...,t_{d}\}$ is also denoted by $(X_{t_{1}},...,X_{t_{d}})$.

Let $S$ be a topological space, $X$ be an $E$-valued process and
$f\in M(E;S)$. The process $\{f\circ X_{t}\}_{t\geq0}$ is the mapping
$\varpi(f)\circ X$ that sends every $\omega\in\Omega$ to the $S$-valued
path $\varpi(f)[X(\omega)]$. A popular notation of this process is
$f\circ X$ but we prefer $\varpi(f)\circ X$.

\subsection{\label{sub:Cadlag_Proc}C$\grave{\mbox{a}}$dl$\grave{\mbox{a}}$g
process}

$X\in M(\Omega,\mathscr{F};E^{\mathbf{R}^{+}})$ is an \textbf{$E$-valued
c$\grave{\mbox{a}}$dl$\grave{\mbox{a}}$g process} if
\begin{equation}
\left\{ \omega\in\Omega:X(\omega)\mbox{ is not a c}\grave{\mbox{a}}\mbox{dl}\grave{\mbox{a}}\mbox{g member of }E^{\mathbf{R}^{+}}\right\} \in\mathscr{F}\cap\mathscr{N}(\mathbb{P}).\label{eq:Cadlag_Proc}
\end{equation}
When $E$ is a Tychonoff space, the \textbf{path space of} an $E$-valued
c$\grave{\mbox{a}}$dl$\grave{\mbox{a}}$g process is thought to be
in $D(\mathbf{R}^{+};E)$. From (\ref{eq:Sko_Borel>Prod}) we immediately
have that:
\begin{fact}
\label{fact:Cadlag_Proc_Dist}Let $E$ be a Tychonoff space. Then:

\renewcommand{\labelenumi}{(\alph{enumi})}
\begin{enumerate}
\item Every $E$-valued c$\grave{\mbox{a}}$dl$\grave{\mbox{a}}$g process
$(\Omega,\mathscr{F},\mathbb{P};X)$ satisfies that
\begin{equation}
\Omega\backslash X^{-1}\left[D(\mathbf{R}^{+};E)\right]\in\mathscr{F}\cap\mathscr{N}(\mathbb{P}),\label{eq:Cadlag_2}
\end{equation}
\begin{equation}
X^{-1}(A)\in\mathscr{F},\;\forall A\in\left.\mathscr{B}(E)^{\otimes\mathbf{R}^{+}}\right|_{D(\mathbf{R}^{+};E)},\label{eq:Cadlag_3}
\end{equation}
and
\begin{equation}
\left.\mathrm{pd}(X)\right|_{D(\mathbf{R}^{+};E)}\in\mathfrak{P}\left(D(\mathbf{R}^{+};E),\left.\mathscr{B}(E)^{\otimes\mathbf{R}^{+}}\right|_{D(\mathbf{R}^{+};E)}\right).\label{eq:Cadlag_Proc_Dist}
\end{equation}

\item Every member of $M(\Omega,\mathscr{F};D(\mathbf{R}^{+};E))$%
\footnote{$M(\Omega,\mathscr{F};D(\mathbf{R}^{+};E))$ denotes the $D(\mathbf{R}^{+};E)$-valued
random variables defined on $(\Omega,\mathscr{F},\mathbb{P})$.%
} is an $E$-valued c$\grave{\mbox{a}}$dl$\grave{\mbox{a}}$g process
defined on $(\Omega,\mathscr{F},\mathbb{P})$.
\end{enumerate}
\end{fact}
\begin{rem}
\label{rem:Cadlag_Sko_RV}As $\mathscr{B}[D(\mathbf{R}^{+};E)]$ is
generically larger than $\mathscr{B}(E)^{\otimes\mathbf{R}^{+}}|_{D(\mathbf{R}^{+};E)}$,
an $E$-valued c$\grave{\mbox{a}}$dl$\grave{\mbox{a}}$g process
is not necessarily a $D(\mathbf{R}^{+};E)$-valued random variable.
More details about c$\grave{\mbox{a}}$dl$\grave{\mbox{a}}$g processes
are presented in \S \ref{sec:Cadlag}.
\end{rem}

\subsection{\label{sub:Proc_Terminology}Stochastic process terminologies}

Let $X$ and $Y$ be $E$-valued processes defined on $(\Omega,\mathscr{F},\mathbb{P})$.
The \textbf{set of fixed left-jump times of} $X$ refers to
\begin{equation}
J(X)\circeq\left\{ t\in\mathbf{R}^{+}:\mathbb{P}\left(\lim_{s\rightarrow t-}X_{s}\mbox{ exists and equals }X_{t}\right)<1\right\} \label{eq:J(X)}
\end{equation}
if it is well-defined. $X$ is a \textit{stationary process} if
\begin{equation}
\mathbb{P}\circ X_{\mathbf{T}_{0}}^{-1}=\mathbb{P}\circ X_{\mathbf{T}_{0}+c}^{-1},\;\forall\mathbf{T}_{0}\in\mathscr{P}_{0}(\mathbf{R}^{+}),c\in(0,\infty),\label{eq:Stationary}
\end{equation}
where
\begin{equation}
\mathbf{T}_{0}+c\circeq\left\{ t+c:t\in\mathbf{T}_{0}\right\} .\label{eq:Timeset_Translation}
\end{equation}
$X$ and $Y$ are (\textit{pathwisely}) \textit{indistinguishable}
if $\{X\neq Y\}\in\mathscr{N}(\mathbb{P})\cap\mathscr{F}$. $X$ and
$Y$ are \textit{modification}s \textit{of }each other if $\{X_{t}\neq Y_{t}\}\in\mathscr{N}(\mathbb{P})\cap\mathscr{F}$
for all $t\in\mathbf{R}^{+}$.

A \textit{filtration} (see \cite[p.453]{D02}) $\{\mathscr{G}_{t}\}_{t\geq0}$\textit{
on} $(\Omega,\mathscr{F},\mathbb{P})$ is \textit{$\mathbb{P}$-complete}
if $\mathscr{G}_{t}$ is $\mathbb{P}$-complete for all $t\geq0$.
We call $(\Omega,\mathscr{F},\{\mathscr{G}_{t}\}_{t\geq0},\mathbb{P})$
a \textit{stochastic basis} if both $\mathscr{F}$ and $\{\mathscr{G}_{t}\}_{t\geq0}$
are $\mathbb{P}$-complete. $X$ is $\mathscr{G}_{t}$-\textit{adapted}
if $\mathscr{F}_{t}^{X}\subset\mathscr{G}_{t}$ for all $t\geq0$,
where
\begin{equation}
\mathscr{F}_{t}^{X}\circeq\sigma\left[\sigma\left(\left\{ X_{u}:u\in[0,t]\right\} \right)\cup\mathscr{N}(\mathbb{P})\right],\;\forall t\geq0.\label{eq:FX}
\end{equation}
$\mathscr{F}^{X}\circeq\{\mathscr{F}_{t}^{X}\}_{t\geq0}$ is called
the \textit{augmented natural filtration of }$X$. Let $\xi(t,\omega)\circeq X_{t}(\omega)$
for each $\omega\in\Omega$ and $t\in\mathbf{R}^{+}$. Then, $X$
is a\textbf{ }\textit{measurable process} if
\begin{equation}
\xi\in M\left(\mathbf{R}^{+}\times\Omega,\mathscr{B}(\mathbf{R}^{+})\otimes\mathscr{F};E,\mathscr{B}(E)\right).\label{eq:Measurable_Proc}
\end{equation}
$X$ is a $\mathscr{G}_{t}$-\textit{progressive process} if
\begin{equation}
\xi|_{[0,t]\times\Omega}\in M\left([0,t]\times\Omega,\mathscr{B}([0,t])\otimes\mathscr{G}_{t};E,\mathscr{B}(E)\right),\;\forall t\in\mathbf{R}^{+}.\label{eq:Prog_Proc}
\end{equation}
$X$ is a \textbf{progressive process} if it is $\mathscr{F}_{t}^{X}$-progressive.

\section{\label{sec:Convention}Conventions}

The following conventions hold hereafter \textbf{if not otherwise
specified}:
\begin{itemize}
\item $\mathbf{I}$ is a non-empty index set.
\item Subsets are non-empty.
\item Measures are non-trivial.
\item Subsets of topological spaces are equipped with their subspace topologies.
\item Topological spaces are equipped with their Borel $\sigma$-algebras.
\item Any Cartesian product of topological spaces is equipped with the product
topology and, hence, is equipped with the Borel $\sigma$-algebra
generated by the product topology.
\item Linear spaces are over the scalar field $\mathbf{R}$.
\item Linear operators are single-valued.
\item $(\Omega,\mathscr{F},\mathbb{P})$, $\{(\Omega^{n},\mathscr{F}^{n},\mathbb{P}^{n})\}_{n\in\mathbf{N}_{0}}$
and $\{(\Omega^{i},\mathscr{F}^{i},\mathbb{P}^{i})\}_{i\in\mathbf{I}}$
are complete probability spaces with expectation operators $\mathbb{E}$,
$\{\mathbb{E}^{n}\}_{n\in\mathbf{N}_{0}}$ and $\{\mathbb{E}^{i}\}_{i\in\mathbf{I}}$,
respectively.
\end{itemize}

\section{\label{sec:Example_Space}Motivating examples}

The Introduction motivates the general use of replication. Herein,
we outline nine specific examples. \S \ref{sub:Pseudo_Path} - \S
\ref{sub:Rough_Path} contain brief reviews of the pseudo-path topology
of c$\grave{\mbox{a}}$dl$\grave{\mbox{a}}$g functions, $\mathcal{S}$-topology
of c$\grave{\mbox{a}}$dl$\grave{\mbox{a}}$g functions, strong topology
of Borel probability measures, strong dual of nuclear Frech$\acute{\mbox{e}}$t
space, Banach spaces of finite $p$-variation or $1/p$-H$\ddot{\mbox{o}}$lder
continuous paths, and Banach spaces of rough paths. These spaces fail
the traditional assumptions of metric \textit{completeness} (see \cite[\S 43, Definition, p.264]{M00}),
compactness, separability and/or metrizability. In \S \ref{sub:Kol},
we review a version of Kolmogorov's Extension Theorem for standard
Borel spaces to illustrate of boosting results by space change. Moreover,
\S \ref{sub:K15-Example} and \S \ref{sub:KK20} refer to two examples
of the convenience of our convergence results in Chapter \ref{chap:Cadlag}.

\subsection{\label{sub:Pseudo_Path}Pseudo-path topology}

In \cite{MZ84}, the pseudo-path topology (also known as ``Meyer-Zheng
topology'') was used to characterize tightness of c$\grave{\mbox{a}}$dl$\grave{\mbox{a}}$g
semimartingales with respect to the topology of convergence in measure.
This is an example of a non-Polish metrizable Lusin space.
\begin{example}
\label{exp:Pseudo-path}Let $D^{\mathrm{pp}}(\mathbf{R}^{+};\mathbf{R})$
denote the c$\grave{\mbox{a}}$dl$\grave{\mbox{a}}$g members of $\mathbf{R}^{\mathbf{R}^{+}}$
equipped with the pseudo-path topology%
\footnote{Pseudo-path topology can be defined similarly on the family of all
c$\grave{\mbox{a}}$dl$\grave{\mbox{a}}$g members of $E^{\mathbf{R}^{+}}$
when $E$ is a Polish space. In that case, $K$ will be a metrizable
compactification of $\mathbf{R}^{+}\times E$.%
}. This topology is induced by the mapping $\psi^{\mathrm{pp}}$ associating
each $x\in D^{\mathrm{pp}}(\mathbf{R}^{+};\mathbf{R})$ to its $\lambda^{\prime}$-almost
everywhere unique \textit{pseudo-path} $\psi^{\mathrm{pp}}(x)\in\mathcal{P}(K)$,
where
\begin{equation}
\lambda^{\prime}(A)\circeq\int_{A}e^{-t}\lambda(dt),\;\forall A\in\mathscr{B}(\mathbf{R}^{+}),\label{eq:lambda_prime}
\end{equation}
$\lambda$ is the Lebesgue measure on $\mathbf{R}^{+}$, $K\circeq[0,\infty]\times[-\infty,\infty]$,
and
\begin{equation}
\psi^{\mathrm{pp}}(x)(B)\circeq\lambda^{\prime}\left(\left\{ t\in\mathbf{R}^{+}:(t,x(t))\in B\right\} \right),\;\forall B\in\mathscr{B}(K).\label{eq:Pseudo-path}
\end{equation}
\cite[Theorem 2]{MZ84} showed that
\begin{equation}
\psi^{\mathrm{pp}}\in\mathbf{imb}\left(D^{\mathrm{pp}}(\mathbf{R}^{+};\mathbf{R});\mathcal{P}(K)\right)\label{eq:PseudoPath_Map}
\end{equation}
and
\begin{equation}
\psi^{\mathrm{pp}}\left[D^{\mathrm{pp}}(\mathbf{R}^{+};\mathbf{R})\right]\in\mathscr{B}\left[\mathcal{P}(K)\right].\label{eq:PseudoPath_Map_Image}
\end{equation}
$K$ is a Polish space by Proposition \ref{prop:Compact} (d). $\mathcal{P}(K)$
is a Polish space by Theorem \ref{thm:P(E)_Compact_Polish} (b). Hence,
$\psi^{\mathrm{pp}}[D^{\mathrm{pp}}(\mathbf{R}^{+};\mathbf{R})]$
is a metrizable Lusin space by Proposition \ref{prop:Borel_in_Polish}
(a, d), and so is its homeomorph $D^{\mathrm{pp}}(\mathbf{R}^{+};\mathbf{R})$.
However, \cite[p.355 - 356]{MZ84} pointed out that $\psi^{\mathrm{pp}}[D^{\mathrm{pp}}(\mathbf{R}^{+};\mathbf{R})]$
and $D^{\mathrm{pp}}(\mathbf{R}^{+};\mathbf{R})$ are not Polish spaces.
\end{example}

\subsection{\label{sub:S-Topo}$\mathcal{S}$-topology}

\cite{J97b} defined the $\mathcal{S}$-topology by introducing the
$\mathcal{S}$-convergence of c$\grave{\mbox{a}}$dl$\grave{\mbox{a}}$g
functions from $[0,T]\subset\mathbf{R}^{+}$ to $\mathbf{R}$, which
is related to the pseodo-path topology. The tightness conditions proposed
by \cite{S85} for the pseudo-path topology turns out to be superfluous
(see \cite{K91}) but serves as a tightness condition and motivation
for the $\mathcal{S}$-topology (see \cite{J97b} and \cite{J12}).
\begin{example}
\label{exp:S-Topo}We define the total variation of $x\in\mathbf{R}^{[0,T]}$
by
\begin{equation}
\Vert x\Vert_{1\mbox{-var},[0,T]}\circeq\left|x(0)\right|+\sup_{0\leq t_{0}<...<t_{n}\leq T,n\in\mathbf{N}}\sum_{i=1}^{n}\left|x(t_{i})-x(t_{i-1})\right|,\label{eq:1-Var}
\end{equation}
and put
\begin{equation}
\mathbf{V}\circeq\left\{ x\in\mathbf{R}^{[0,T]}:x\mbox{ is c}\grave{\mbox{a}}\mbox{dl}\grave{\mbox{a}}\mbox{g},\Vert x\Vert_{1\mbox{-var},[0,T]}<\infty\right\} .\label{eq:S-Topo_V}
\end{equation}
C$\grave{\mbox{a}}$dl$\grave{\mbox{a}}$g functions $\{x_{n}\}_{n\in\mathbf{N}}\subset\mathbf{R}^{[0,T]}$
$\mathcal{S}$-converge to c$\grave{\mbox{a}}$dl$\grave{\mbox{a}}$g
function $x_{0}\in\mathbf{R}^{[0,T]}$ if for any $\epsilon\in(0,\infty)$,
there exist $\{v_{n}^{\epsilon}\}_{n\in\mathbf{N}_{0}}\in\mathbf{V}$
such that
\begin{equation}
\sup_{n\in\mathbf{N}_{0}}\Vert x_{n}-v_{n}^{\epsilon}\Vert_{1\mbox{-var},[0,T]}<\epsilon\label{eq:Def_S-Topo_1}
\end{equation}
and
\begin{equation}
\lim_{n\rightarrow\infty}\int_{[0,T]}f(t)dv_{n}(t)=\int_{[0,T]}f(t)dv_{0}(t),\;\forall f\in C([0,T];\mathbf{R}).\label{eq:Def_S-Topo_2}
\end{equation}
Considering $\mathcal{S}$-convergence as convergence a priori, the
$\mathcal{S}$-topology is the topology of convergence a priori on
the c$\grave{\mbox{a}}$dl$\grave{\mbox{a}}$g members of $\mathbf{R}^{[0,T]}$.
This $\mathcal{S}$-topological space is a coarsening of the Skorokhod
$\mathscr{J}_{1}$-space $D([0,T];\mathbf{R})$ (see \cite[p.5]{J12})
that is neither necessarily a Tychonoff space, nor known to be a topological
vector space (see \cite{J97b} and \cite[p.5]{J12}). The corresponding
convergence a posteriori, called $\mathcal{S}*$-convergence, is different
than $\mathcal{S}$-convergence (see \cite[p.3 - 4]{J12}). Moreover,
$\mathcal{S}*$-convergence is also the topological convergence of
some coarserning of the $\mathcal{S}$-topology. The equality of these
two topologies is an open question (see \cite[p.4]{J12}).
\end{example}

\subsection{\label{sub:Strong_Topo}Strong topology of Borel probability measures}

\cite{DZ98} used the strong topology of Borel probability measures
on a Polish space in large deviations. This gives a non-metrizable
and non-\textit{separable} (see p.\pageref{separable}) Tychonoff
space.
\begin{example}
\label{exp:Strong-topo}Let $E$ be a Polish space and $\mathcal{P}_{S}$
be the space of all Borel probability measures on $E$ equipped with
the \textit{strong topology}
\begin{equation}
\mathscr{O}\left[\mathcal{P}_{S}(E)\right]\circeq\mathscr{O}_{M_{b}(E;\mathbf{R})^{*}}\left[\mathcal{P}(E)\right].\label{eq:Strong_Topology}
\end{equation}
Then, $\mathcal{P}_{S}(E)$ is a topological refinement of $\mathcal{P}(E)$,
\begin{equation}
M_{b}(E;\mathbf{R})^{*}\subset C_{b}\left(\mathcal{P}_{S}(E);\mathbf{R}\right),\label{eq:Mb_Star_SSP}
\end{equation}
and $M_{b}(E;\mathbf{R})^{*}$ strongly separates points on $\mathcal{P}_{S}(E)$.
Furthermore, from the fact
\begin{equation}
\left\{ (\mathbf{1}_{A})^{*}:A\in\mathscr{B}(E)\right\} \subset M_{b}(E;\mathbf{R})^{*}\label{eq:Mb_Star_SP}
\end{equation}
it follows that $M_{b}(E;\mathbf{R})^{*}$ separates points on $\mathcal{P}_{S}(E)$.
Hence, $\mathcal{P}_{S}(E)$ is a Tychonoff space by Proposition \ref{prop:CR}
(a, b). However, \cite[p.263]{DZ98} explained that $\mathcal{P}_{S}(E)$
is neither metrizable nor separable.
\end{example}

\subsection{\label{sub:Dual_Nuclear}Strong dual of nuclear Frech$\acute{\mbox{e}}$t
space}

\cite[\S 5.II]{J86} discussed tightness of probability measures on
the Skorokhod $\mathscr{J}_{1}$-space $D([0,1];E)$ with $E$ being
the strong dual of a general nuclear Frech$\acute{\mbox{e}}$t space.
This is an example of a possibly non-metrizable, Tychonoff topological
vector space.
\begin{example}
\label{exp:Dual_Nuclear}Let $E$ be the strong dual of some infinite-dimensional,
unnormable, nuclear Frech$\acute{\mbox{e}}$t space. $E$ is not metrizable
by \cite[\S 29.1, (7), p.394]{GK12}, nor is it necessarily separable.
However, $E$ is a nuclear space by \cite[\S IV.9.6, Theorem, p.172]{SW99},
the topology of $E$ is induced by a family of \textit{Hilbertian
semi-norms} (see \cite[Definition A.4]{SH12}), and $E$ is a Tychonoff
space by \cite[Theorem 2.1.1]{KX95}. Hence, $D([0,1];E)$ is a non-metrizable,
possibly non-separable, Tychonoff space by \cite[Proposition 1.6 ii) - iii)]{J86}
and Proposition \ref{prop:Sko_Basic_1} (e).
\end{example}

\subsection{\label{sub:Holder}Spaces of finite-variation or H$\ddot{\mbox{o}}$lder
continuous functions}

The spaces of $\mathbf{R}^{d}$-valued continuous functions with finite
$p$-variation or $\mathbf{R}^{d}$-valued $1/p$-H$\ddot{\mbox{o}}$lder
continuuous functions are frequently used in stochastic differential
equations driven by non-classical noises. They are examples of non-separable
Banach spaces.
\begin{example}
\label{exp:Holder}Let $d,N\in\mathbf{N}$, $p\in[1,\infty)$ and
$T\in(0,\infty)$. A path $x\in(\mathbf{R}^{d})^{[0,T]}$ has finite
$p$-variation if the homogeneous $p$-variation norm of $x$ defined
by
\begin{equation}
\Vert x\Vert_{p\mbox{-var},[0,T]}\circeq\left|x(0)\right|+\sup_{0\leq t_{0}<...<t_{n}\leq T,n\in\mathbf{N}}\left(\sum_{i=1}^{n}\left|x(t_{i})-x(t_{i-1})\right|^{p}\right)^{1/p},\label{eq:p-Var_Norm_Rd}
\end{equation}
or is $1/p$-H$\ddot{\mbox{o}}$lder continuous if the $1/p$-H$\ddot{\mbox{o}}$lder
norm of $x$ defined by
\begin{equation}
\Vert x\Vert_{\frac{1}{p}\mathrm{-H\ddot{\mathrm{o}}l},[0,T]}\circeq\sup_{0\leq s<t\leq T}\frac{\left|x(t)-x(s)\right|}{\left|t-s\right|^{\frac{1}{p}}}\label{eq:1/p-Holder_Norm_Rd}
\end{equation}
is finite respectively. The normed spaces
\begin{equation}
\left\{ x\in C\left([0,T];\mathbf{R}^{d}\right):\Vert x\Vert_{p\mathrm{-var},[0,T]}<\infty\right\} \label{eq:p-Var_Space_Rd}
\end{equation}
and
\begin{equation}
\left\{ x\in\left(\mathbf{R}^{d}\right)^{[0,T]}:\Vert x\Vert_{\frac{1}{p}\mathrm{-H\ddot{\mathrm{o}}l},[0,T]}<\infty\right\} \label{eq:1/p-Holder_Space_Rd}
\end{equation}
are non-separable Banach spaces (see \cite[Theorem 5.27]{FV10}).
\end{example}

\subsection{\label{sub:Rough_Path}Space of rough paths}

The rough path approach, initiated by the pioneering works of \cite{L95}
and \cite{L98}, is important to generalizing stochastic differential
equations like
\begin{equation}
dY_{t}=\alpha(t,Y_{t})dt+\sigma(t,Y_{t})dX_{t}\label{eq:SDE}
\end{equation}
to the case where the driving noise $X$ is not necessarily a semimartingale.
By this approach, the original noise $X$ is \textit{enhanced to a
random rough path $\mathfrak{X}$} (see \cite[\S 9.1]{FV10}) and
the Stratonovich solution of (\ref{eq:SDE}) is closely linked to
the solution of
\begin{equation}
dY_{t}=\alpha(t,Y_{t})dt+\sigma(t,Y_{t})d\mathfrak{X}_{t},\label{eq:SDE_RoughPath}
\end{equation}
where (\ref{eq:SDE_RoughPath}) is considered as \textit{rough differential
equations driven by the realization paths of $\mathfrak{X}$ as a
process} (see \cite[\S 10.3, \S 10.4, \S 17.1 and \S 17.2]{FV10}).
\cite{FV10} and \cite{FH14} considered the following spaces for
the paths of \textit{$\mathfrak{X}$}, which are also examples of
non-separable Banach spaces.
\begin{example}
\label{exp:Rough_Path}Let $d,N\in\mathbf{N}$, $p\in[1,\infty)$
and $T\in(0,\infty)$. A rough path is often considered as a mapping
from $[0,T]$ to $G^{N}(\mathbf{R}^{d})$, the \textit{free nilpotent
group of Step $N$ over $\mathbf{R}^{d}$} (see \cite[p.142-143]{FV10}).
As a Lie group, $G^{N}(\mathbf{R}^{d})$ is equipped with the usual
addition ``$+$'' of functions and the \textit{Carnot-Caratheodory
norm $\Vert\cdot\Vert_{\mathrm{cc}}$} (see \cite[Theorem 7.32]{FV10}).
Similar to $\mathbf{R}^{d}$-valued paths, a path $x\in G^{N}(\mathbf{R}^{d})^{[0,T]}$
has finite $p$-variation or is $1/p$-H$\ddot{\mbox{o}}$lder continuous
if the homogeneous $p$-variation $\mathrm{cc}$-norm of $x$ defined
by
\begin{equation}
\Vert x\Vert_{\mathrm{cc},p\mathrm{-var},[0,T]}\circeq\sup_{0\leq t_{0}<...<t_{n}\leq T,n\in\mathbf{N}}\left(\sum_{i=1}^{n}\left\Vert x(t_{i})-x(t_{i-1})\right\Vert _{\mathrm{cc}}^{p}\right)^{1/p}\label{eq:p-Var_cc_norm}
\end{equation}
or the homogeneous $1/p$-H$\ddot{\mbox{o}}$lder $\mathrm{cc}$-norm
of $x$ defined by
\begin{equation}
\Vert x\Vert_{\mathrm{cc},\frac{1}{p}\mathrm{-H\ddot{\mathrm{o}}l},[0,T]}\circeq\sup_{0\leq s<t\leq T}\frac{\left\Vert x(t)-x(s)\right\Vert _{\mathrm{cc}}}{\left|t-s\right|^{\frac{1}{p}}}\label{eq:1/p-Holder_cc_norm}
\end{equation}
is finite respectively. The random rough path $\mathfrak{X}$ in (\ref{eq:SDE_RoughPath})
may have paths in
\begin{equation}
\left\{ x\in C\left([0,T];G^{N}\left(\mathbf{R}^{d}\right)\right):\Vert x\Vert_{\mathrm{cc},p\mathrm{-var},[0,T]}<\infty\right\} \label{eq:p-Var_Space_Rough}
\end{equation}
or
\begin{equation}
\left\{ x\in G^{N}\left(\mathbf{R}^{d}\right)^{[0,T]}:\Vert x\Vert_{\mathrm{cc},\frac{1}{p}\mathrm{-H\ddot{\mathrm{o}}l},[0,T]}<\infty\right\} .\label{eq:1/p-Holder_Space_Rough}
\end{equation}
By \cite[Theorem 8.13]{FV10}, these normed spaces are non-separable
Banach spaces.
\end{example}

\subsection{\label{sub:Kol}Kolmogorov's Extension Theorem}

The Kolmogorov's Extension Theorem (see \cite[\S 15.6]{AB06}) is
a cornerstone of probability theory that depends purely on the relevant
$\sigma$-algebras. Hence, existence of Kolmogorov extension should
be transferrable from a ``nice'', typically Polish topological space
to ``defective'', typically standard Borel topological space which
are ``indifferent'' as measurable spaces.

Let $\{S_{i}\}_{i\in\mathbf{I}}$ be a family of standard Borel spaces,
\begin{equation}
(S,\mathscr{A})\circeq\left(\prod_{i\in\mathbf{I}}S_{i},\bigotimes_{i\in\mathbf{I}}\mathscr{B}(S_{i})\right),\label{eq:(S,A)_Prod_Meas_Space}
\end{equation}
and
\begin{equation}
\left(S_{\mathbf{I}_{0}},\mathscr{A}_{\mathbf{I}_{0}}\right)\circeq\left(\prod_{i\in\mathbf{I}_{0}}S_{i},\bigotimes_{i\in\mathbf{I}_{0}}\mathscr{B}(S_{i})\right),\;\forall\mathbf{I}_{0}\in\mathscr{P}_{0}(\mathbf{I}).\label{eq:(S_I0,A_I0)_Kol_Ext}
\end{equation}
For each $i\in\mathbf{I}$, Proposition \ref{prop:SB} (a, d) allows
us to change the topology of $S_{i}$ to a possibly different one
$\mathscr{U}_{i}$ such that $(S_{i},\mathscr{U}_{i})$ is a Polish
space and the Borel sets $\mathscr{B}(S_{i})=\mathscr{B}(S_{i},\mathscr{U}_{i})$
remain unchanged. So, any $\mu_{\mathbf{I}_{0}}\in\mathfrak{P}(S_{\mathbf{I}_{0}},\mathscr{A}_{\mathbf{I}_{0}})$
can be viewed as a probability measure on $(S_{\mathbf{I}_{0}},\bigotimes_{i\in\mathbf{I}_{0}}\mathscr{B}(S_{i},\mathscr{U}_{i}))$
for each $\mathbf{I}_{0}\subset\mathscr{P}_{0}(\mathbf{I})$, and
any Kolmogorov extension of $\{\mu_{\mathbf{I}_{0}}\}_{\mathbf{I}_{0}\in\mathscr{P}_{0}(\mathbf{I})}$
on $(S,\bigotimes_{i\in\mathbf{I}}\mathscr{B}(S_{i},\mathscr{U}_{i}))$
would be a desired Kolmogorov extension of them on $(S,\mathscr{A})$.
Therefore, the well-known version of Kolmogorov's Extension Theorem
for Polish spaces extends immediately to the standard Borel case.
\begin{thm}
[\textbf{Kolmogorov's Extension Theorem}\textrm{, \cite[Theorem 5.16]{K97}}]\label{thm:Kol_Ext}Let
$\{S_{i}\}_{i\in\mathbf{I}}$ be standard Borel spaces, $(S,\mathscr{A})$
be as in (\ref{eq:(S,A)_Prod_Meas_Space}), $\{(S_{\mathbf{I}_{0}},\mathscr{A}_{\mathbf{I}_{0}})\}_{\mathbf{I}_{0}\in\mathscr{P}_{0}(\mathbf{I})}$
be as in (\ref{eq:(S_I0,A_I0)_Kol_Ext}) and $\mu_{\mathbf{I}_{0}}\in\mathfrak{P}(S_{\mathbf{I}_{0}},\mathscr{A}_{\mathbf{I}_{0}})$
for each $\mathbf{I}_{0}\in\mathscr{P}_{0}(\mathbf{I})$. Suppose
in addition that for each $\mathbf{I}_{1},\mathbf{I}_{2}\in\mathscr{P}_{0}(\mathbf{I})$
with $\mathbf{I}_{1}\subset\mathbf{I}_{2}$, $\mu_{\mathbf{I}_{1}}$
is the push-forward measure of $\mu_{\mathbf{I}_{2}}$ by the projection
from $S_{\mathbf{I}_{2}}$ to $S_{\mathbf{I}_{1}}$. Then, there exists
a $\mu\in\mathfrak{P}(S,\mathscr{A})$ such that for each $\mathbf{I}_{0}\in\mathscr{P}_{0}(\mathbf{I})$,
$\mu_{\mathbf{I}_{0}}$ is the push-forward measure of $\mu$ by the
projection from $S$ to $S_{\mathbf{I}_{0}}$.
\end{thm}

\subsection{\label{sub:K15-Example}Approximation by Markov Empirical Processes}

A traditional way of constructing measure-valued processes is to show
convergence of weighted empirical processes to it. \cite{K15} confirmed
the general availability of such approximation on a Polish underlying
space. \cite{K15} employed similar compactification technique to
that of this paper and reduced tightness to the verification of Modulus
of Continuity Condition (see \S \ref{sub:Proc_Reg}). Unlike the
convergence results in \cite[\S 2]{BK93b}, \cite{K15} neither assumed
pointwise tightness condition, for either the limit or pre-limit processes,
nor required a martingale problem set-up. More importantly, an inspection
into the development of \cite[Theorem 2 and Theorem 3]{K15} shows
that the Polish space imposition was absent in verifying finite-dimensional
convergence or Modulus of Continuity Condition. Instead, metric completeness
was only used to compactify the space of measures and establish tightness.
By Corollary \ref{cor:Sko_FC_WC_Metrizable_Separable} to follow,
one only needs Modulus of Continuity Condition rather than tightness
for deriving weak convergence of c$\grave{\mbox{a}}$d$\grave{\mbox{a}}$g
processes from their finite-dimensional convergence. Therefore, the
following generalization of \cite[Theorem 3]{K15} is immediate from
Corollary \ref{cor:Sko_FC_WC_Metrizable_Separable} and Corollary
\ref{cor:Metrizable_Separable_P(E)}:
\begin{thm}
Let $E$ be a metrizable and separable space and $V=\{V_{t}\}_{t\geq0}$
be a process with $D(\mathbf{R}^{+};\mathcal{P}(E))$-valued paths.
Then, there exists $m_{N}\uparrow\infty$ and conditionally i.i.d.
$E$-valued c$\grave{\mbox{a}}$dl$\grave{\mbox{a}}$g Markov processes
$\{\Xi^{1,N},...,\Xi^{m_{N},N}\}_{N\in\mathbf{N}}$ such that the
empirical processes $V^{N}\circeq\frac{1}{m_{N}}\sum_{i=1}^{N}\delta_{\Xi^{i,N}}$
converge weakly to $V$ almost surely on $D(\mathbf{R}^{+};\mathcal{P}(E))$.
\end{thm}

\subsection{\label{sub:KK20}General Existence of Infinite Systems of Stochastic
Differential Equations}

Infinite-dimensional systems of stochastic differential equations
are important for modelling many vast and complicated real systems
and for analysis of stochastic partial differential equations. Several
authors, e.g. \cite{F87}, \cite{S99,S01} and \cite{ABGP06}, have
studied specific systems, like stochastic gradient equations or Ornstein-Uhlenbeck
perturbations, in detail. However, recently \cite{KK20} showed weak
existence for a relatively general system of the form:
\begin{equation}
\begin{aligned}X_{t}^{i} & =X_{0}^{i}+\int_{0}^{t}\sigma\left(X_{s}^{i},N_{s}\right)dB_{s}^{i}\\
 & +\int_{0}^{t}b\left(X_{s}^{i},N_{s}\right)ds+\int_{0}^{t}\alpha\left(X_{s}^{i},N_{s}\right)dW_{s},\; i\in\mathbf{Z},
\end{aligned}
\label{eq:KK20_SDE}
\end{equation}
where $\mathbf{Z}$ is the integers and $N_{t}\circeq\sum_{i\in\mathbf{Z}}\delta_{X_{t}^{i}}$
is the configuration at time $t\geq0$.

Let $\mathcal{M}^{\beta}(\mathbf{R}^{d})$ for $\beta>0$ denote the
family of all Borel measure $\mu$ on $\mathbf{R}^{d}$ satisfying
\begin{equation}
\int_{\mathbf{R}^{d}}(1+\left|x\right|^{2})^{-\beta}\mu(dx)<\infty.\label{eq:KK20_JC}
\end{equation}
Then, the coefficient functions $b(x,\nu)$, $\sigma(x,\nu)$ and
$\alpha(x,\nu)$ and initial configuration are regulated by the following
mild conditions:

\begin{enumerate}
[label=\textbf{(I)}, labelsep=0.5pc]

\item\label{enu:KK20_I}There exist $\beta,\delta\in(0,\infty)$
such that
\begin{equation}
\mathbb{E}\left[\left(\sum_{i\in\mathbf{Z}}\left(1+\left|X_{0}^{i}\right|^{2}\right)^{-\beta}\right)^{2+\delta}\right]<\infty.\label{eq:KK20_I}
\end{equation}

\end{enumerate}

\begin{enumerate}
[label=\textbf{(LG)}, labelsep=0.5pc]

\item\label{enu:KK20_LG}$b(x,\nu)$, $\sigma(x,\nu)$ and $\alpha(x,\nu)$
satisfy for $\beta$ in \ref{enu:KK20_I} that
\begin{equation}
\begin{aligned} & \sup_{x\in\mathbf{R}^{d},\nu\in\mathcal{M}^{\beta}(\mathbf{R}^{d})}\frac{\left|b(x,\nu)\right|}{1+\left|x\right|}<\infty,\\
 & \sup_{x\in\mathbf{R}^{d},\nu\in\mathcal{M}^{\beta}(\mathbf{R}^{d})}\frac{\left|\sigma(x,\nu)\right|}{1+\left|x\right|}<\infty,\\
 & \sup_{x\in\mathbf{R}^{d},\nu\in\mathcal{M}^{\beta}(\mathbf{R}^{d})}\frac{\left|\alpha(x,\nu)\right|}{1+\left|x\right|}<\infty.
\end{aligned}
\label{eq:KK20_LG}
\end{equation}

\end{enumerate}

\begin{enumerate}
[label=\textbf{(JC)}, labelsep=0.5pc]

\item\label{enu:KK20_JC}There is a $\overline{\beta}>\beta$ in
\ref{enu:KK20_I} such that $b$ and both $\sigma$, $\alpha$ are
continuous functions from $\mathbf{R}^{d}\times\mathcal{M}^{\overline{\beta}}(\mathbf{R}^{d})$
to $\mathbf{R}^{d}$ and to $\mathbf{R}^{d\times d}$ respectively.

\end{enumerate}

Their main result, establishing both existence of solutions to the
infinite system (\ref{eq:KK20_SDE}) and weak convergence of weighted-empirical-measures
into the set of all solutions, is stated below:
\begin{thm}
[\textrm{\cite[Theorem 2]{KK20}}]Suppose that (\ref{enu:KK20_LG},
\ref{enu:KK20_I}) hold. Then:

\renewcommand{\labelenumi}{(\alph{enumi})}
\begin{enumerate}
\item The weighted-empirical-measures
\begin{equation}
N_{t}^{n}\circeq\sum_{i\in\mathbf{Z}}\Phi_{t}^{i,n}\delta_{X_{t}^{i,n}},\; t\in\mathbf{R}^{+},n\in\mathbf{N}\label{eq:KK20_approx}
\end{equation}
constructed in \cite[p.6]{KK20} are bijectively indistinguishable
from a tight sequence of $D(\mathbf{R}^{+};\mathcal{M}^{+}(\mathbf{R}^{d}))$-valued
random variables, where $\Phi_{t}^{i,n}$ is a discretization of $(1+|X_{t}^{i,n}|^{2})^{-\overline{\beta}}$.
\item If, in addition, \ref{enu:KK20_JC} holds, then every weak limit point
of the tight sequence in (a) is a weak solution $N$ to (\ref{eq:KK20_SDE})
via density $(1+|x|^{2})^{\overline{\beta}}$.
\end{enumerate}
\end{thm}
$b$, $\sigma$ and $\alpha$ could be unbounded and non-Lipschitz
under the conditions above. The approximations $\{X^{i,n}\}_{n\in\mathbf{N}}$
of $X^{i}$ and the weight processes $\{\Phi^{i,n}\}_{n\in\mathbf{N}}$
of particle $\delta_{X^{i,n}}$ in (\ref{eq:KK20_approx}) are constructed
naturally in Euler-type form (see \cite[(2.11) - (2.14)]{KK20}) and
so computer implementable. \cite[\S 4 - 5]{KK20} showed their convergence
to a solution to the martingale problem form of (\ref{eq:KK20_SDE}).
Our Theorem \ref{thm:Sko_Tight} to follow reduced the desired tightness
to Weak Modulus of Continuity Condition (see \S  \ref{sub:Proc_Reg})
regarding a simple family of rational functions. \cite[Theorems 2.1 and Theorem 2.3]{BK93b}
did not apply as the well-posedness of the martingale problem was
not imposed and little was required on the coefficients.

\chapter{\label{chap:Space_Change}Space Change in Replication}

A topological space $E$ needs no enhancement when it is compact and
metrizable. Otherwise, a problem on $E$ might be simplified if it
is translated onto such a ``perfect'' space. Replication is a method
of space change for this purpose.

The current chapter discusses the space change aspect of replication.
\S \ref{sec:Base} introduces the notion of \textit{base} as our
core platform to implement space change and other goals of replication.
\S \ref{sec:Baseable_Space} and \S \ref{sec:Baseable_Subsets}
explain the existence and various properties of \textit{baseable spaces}
or \textit{baseable subsets} with which one can construct the desired
bases.

\section{\label{sec:Base}Base}

The goal of space change in replication is to create a compact metric
space $\widehat{E}$ related to the original space $E$. As illustrated
by the following figure, the most natural way is to establish a metrizable
\textit{compactification}%
\footnote{In this work, we always consider any compactification to be a Hausdorff
space.%
} (see p.\pageref{Compactification}) of $E$ itself or, more generally,
a Borel subset $E_{0}$ of $E$.

\begin{figure}[H]
\begin{centering}
\includegraphics[scale=0.8]{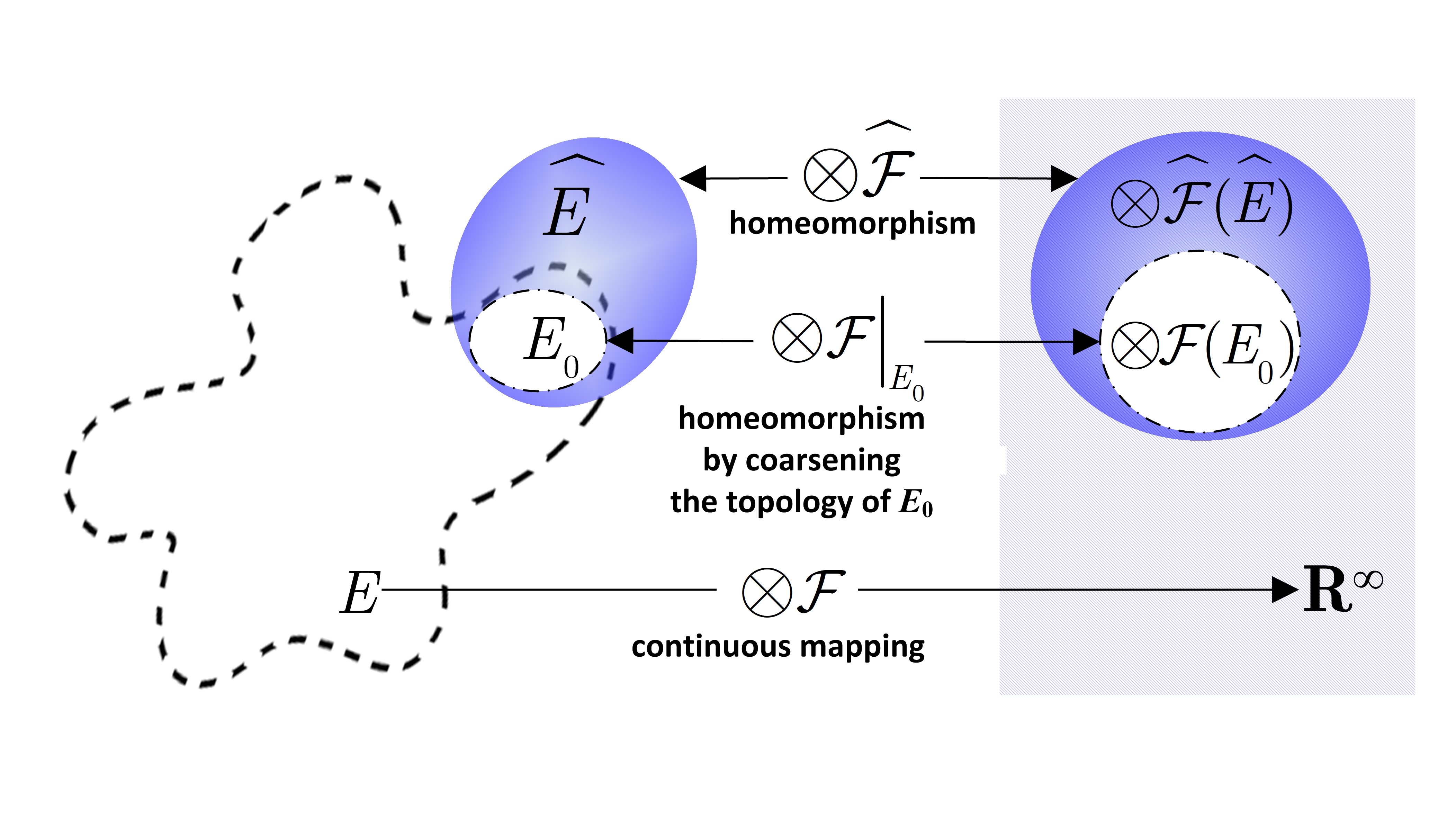}
\par\end{centering}

\caption{\textit{\label{fig:Space_Change}Space change in replication}}
\end{figure}

\subsection{\label{sub:Base_Def}Definition}

A base is a foundational notion of replication that concretizes the
space change idea mentioned above.
\begin{defn}
\label{def:Base}Let $E$ be a topological space. The quadruple $(E_{0},\mathcal{F};\widehat{E},\widehat{\mathcal{F}})$
is a \textbf{replication base over $E$} (a base over $E$ or a base
for short) if:
\begin{itemize}
\item $E_{0}$ is a non-empty Borel subset of $E$.
\item $\mathcal{F}\subset C_{b}(E;\mathbf{R})$ is countable and contains
the constant function $1$.
\item $\widehat{E}$ is a topological space containing $E_{0}$.
\item $\widehat{\mathcal{F}}\subset\mathbf{R}^{\widehat{E}}$ is a countable
collection, separates points on $\widehat{E}$ and satisfies%
\footnote{The notations ``$\mathcal{F}|_{E_{0}}$'' and ``$\bigotimes\mathcal{F}$''
were defined in \S \ref{sub:Meas}. ``$\mathscr{O}(\cdot)$'' denotes
the family of all open subsets.%
}
\begin{equation}
\bigotimes\mathcal{F}|_{E_{0}}=\bigotimes\widehat{\mathcal{F}}|_{E_{0}}\label{eq:F_Fhat_Coincide}
\end{equation}
and
\begin{equation}
\mathscr{O}(\widehat{E})=\mathscr{O}_{\widehat{\mathcal{F}}}(\widehat{E}).\label{eq:Fhat_SSP_on_Ehat}
\end{equation}

\item $\bigotimes\widehat{\mathcal{F}}(\widehat{E})$ is the closure of
$\bigotimes\mathcal{F}(E_{0})$ in $\mathbf{R}^{\infty}$.
\end{itemize}
\end{defn}
\begin{rem}
\label{rem:E_Ehat}In general, $E$ need not be a subset of $\widehat{E}$.
\end{rem}
The following lemma shows a base establishes the compactification
in Figure \ref{fig:Space_Change}.
\begin{lem}
\label{lem:Base}Let $E$ be a topological space, $(E_{0},\mathcal{F};\widehat{E},\widehat{\mathcal{F}})$
be a base over $E$ and $A\subset E_{0}$. Then:

\renewcommand{\labelenumi}{(\alph{enumi})}
\begin{enumerate}
\item $\widehat{\mathcal{F}}\subset C(\widehat{E};\mathbf{R})$ is countable,
contains the constant function $1$ and strongly separates points
on $\widehat{E}$. In particular,
\begin{equation}
\bigotimes\widehat{\mathcal{F}}\in\mathbf{imb}\left(\widehat{E};\mathbf{R}^{\infty}\right).\label{eq:Base_Imb}
\end{equation}

\item $\bigotimes\widehat{\mathcal{F}}(\widehat{E})$ is a compactification
of $\bigotimes\mathcal{F}(E_{0})$ and $\widehat{E}$ is a compactification
of $(E_{0},\mathscr{O}_{\mathcal{F}}(E_{0}))$.
\item $\widehat{E}$ is a Polish space and is completely metrized by $\rho_{\widehat{\mathcal{F}}}$%
\footnote{$\rho_{\widehat{\mathcal{F}}}$ is defined by (\ref{eq:TF_Metric})
with $\mathcal{D}=\widehat{\mathcal{F}}$.%
}.
\item $\bigotimes\mathcal{F}|_{A}\in\mathbf{imb}(A,\mathscr{O}_{\widehat{E}}(A);\mathbf{R}^{\infty})$.
Moreover, $(A,\mathscr{O}_{\widehat{E}}(A))$ is a metrizable and
separable topological coarsening of $(A,\mathscr{O}_{E}(A))$.
\item $\bigotimes\mathcal{F}\in C(E;\mathbf{R}^{\infty})$ is injective
on $A$. Moreover, $\mathcal{F}$ separates points on the Hausdorff
space $(A,\mathscr{O}_{E}(A))$.
\end{enumerate}
\end{lem}
\begin{proof}
(a) The members of $\widehat{\mathcal{F}}$ are continuous and $\widehat{\mathcal{F}}$
strongly separates points on $\widehat{E}$ by (\ref{eq:Fhat_SSP_on_Ehat}).
$\widehat{\mathcal{F}}$ is countable and contains $1$ by (\ref{eq:F_Fhat_Coincide})
and $1\in\mathcal{F}$. Moreover, (\ref{eq:Base_Imb}) follows by
Lemma \ref{lem:Compactification} (a, c) (with $E=\widehat{E}$ and
$\mathcal{D}=\widehat{\mathcal{F}}$).

(b) $\mathbf{R}^{\infty}$ is a Polish space by Proposition \ref{prop:Var_Polish}
(f). It follows by the Tychonoff Theorem (Proposition \ref{prop:Compact}
(b)) and Proposition \ref{prop:Compact} (a) that
\begin{equation}
K_{\mathcal{F}}\circeq\prod_{f\in\mathcal{F}}[-\Vert f\Vert_{\infty},\Vert f\Vert_{\infty}]\in\mathscr{K}(\mathbf{R}^{\infty})\subset\mathscr{C}(\mathbf{R}^{\infty}).\label{eq:Inf-Dim_Cube_TF}
\end{equation}
$\bigotimes\widehat{\mathcal{F}}(\widehat{E})$ is the closure of
$\bigotimes\mathcal{F}(E_{0})$ in $\mathbf{R}^{\infty}$ by definition.
So,
\begin{equation}
\bigotimes\widehat{\mathcal{F}}(\widehat{E})\in\mathscr{C}\left(K_{\mathcal{F}},\mathscr{O}_{\mathbf{R}^{\infty}}(K_{\mathcal{F}})\right)\subset\mathscr{K}(\mathbf{R}^{\infty})\subset\mathscr{C}(\mathbf{R}^{\infty})\label{eq:Prod(Fhat)(Ehat)_Compact_Rinf}
\end{equation}
and $\widehat{E}$ is compact by (\ref{eq:Inf-Dim_Cube_TF}), (\ref{eq:Base_Imb})
and Proposition \ref{prop:Compact} (a, e).

Moreover, it follows by (a) and (\ref{eq:F_Fhat_Coincide}) that
\begin{equation}
\mathscr{O}_{\widehat{E}}(E_{0})=\mathscr{O}_{\widehat{\mathcal{F}}}(E_{0})=\mathscr{O}_{\mathcal{F}}(E_{0}).\label{eq:Check_E0_Topo}
\end{equation}
$\bigotimes\mathcal{F}(E_{0})$ is dense in $\bigotimes\widehat{\mathcal{F}}(\widehat{E})$
by definition, so $E_{0}$ is a dense subset of $\widehat{E}$ by
(\ref{eq:Base_Imb}). $\widehat{E}$ is a Hausdorff space by (\ref{eq:Fhat_SSP_on_Ehat})
and Proposition \ref{prop:Fun_Sep_1} (c) (with $E=A=\widehat{E}$
and $\mathcal{D}=\widehat{\mathcal{F}}$). Hence, $\widehat{E}$ is
a compactification of $(E_{0},\mathscr{O}_{\mathcal{F}}(E_{0}))$.

(c) $\rho_{\widehat{\mathcal{F}}}$ metrizes $\widehat{E}$ by (a)
and Proposition \ref{prop:Fun_Sep_1} (d) (with $E=\widehat{E}$ and
$\mathcal{D}=\widehat{\mathcal{F}}$). $\bigotimes\widehat{\mathcal{F}}$
is an \textit{isometry} (see p.\pageref{Isometry}) between $(\widehat{E},\rho_{\widehat{\mathcal{F}}})$
and $(\bigotimes\widehat{\mathcal{F}}(\widehat{E}),d^{\infty})$,
where%
\footnote{$\mathfrak{p}_{n}$ denotes the one-dimensional projection on $\mathbf{R}^{\infty}$
for $n\in\mathbf{N}$.%
}
\begin{equation}
d^{\infty}(x,y)\circeq\sum_{n=1}^{\infty}2^{-n+1}\left(\left|\mathfrak{p}_{n}(x)-\mathfrak{p}_{n}(y)\right|\wedge1\right),\;\forall x,y\in\mathbf{R}^{\infty}\label{eq:R_Inf_Metric}
\end{equation}
completely metrizes $\mathbf{R}^{\infty}$ by Proposition \ref{prop:Metric_Prod}
(b) (with $(S_{i},\mathfrak{r}_{i})=\mathbf{R}$). $(\bigotimes\widehat{\mathcal{F}}(\widehat{E}),d^{\infty})$
is complete by its compactness and Proposition \ref{prop:Compact}
(c). Thus, $(\widehat{E},\rho_{\widehat{\mathcal{F}}})$ is a complete
metric space by Proposition \ref{prop:Completeness} (a).

(d) The first statement of (d) follows by (\ref{eq:Base_Imb}) and
(\ref{eq:F_Fhat_Coincide}). $(A,\mathscr{O}_{\widehat{E}}(A))$ is
metrizable and separable by (c) and Proposition \ref{prop:Var_Polish}
(c). Moreover, one finds by (\ref{eq:Check_E0_Topo}) and $\mathcal{F}\subset C(E;\mathbf{R})$
that
\begin{equation}
\mathscr{O}_{\widehat{E}}(A)=\mathscr{O}_{\mathcal{F}}(A)\subset\mathscr{O}_{E}(A).\label{eq:Fhat_Coarser_F_on_A}
\end{equation}

(e) The first statement of (e) follows by (d), $\mathcal{F}\subset C(E;\mathbf{R})$
and Fact \ref{fact:Prod_Map_2} (b). The second part follows by Proposition
\ref{prop:Fun_Sep_1} (e) (with $\mathcal{D}=\mathcal{F}$).\end{proof}

\begin{cor}
\label{cor:Base_Unique}Let $E$ be a topological space. If $\{(E_{0},\mathcal{F};\widehat{E}_{i},\widehat{\mathcal{F}}_{i})\}_{i=1,2}$
are bases over $E$, then $\widehat{E}_{1}$ and $\widehat{E}_{2}$
are isometric hence homeomorphic.
\end{cor}
\begin{proof}
Let $\mathcal{F}=\{f_{n}\}_{n\in\mathbf{N}}$. By (\ref{eq:F_Fhat_Coincide})
and Lemma \ref{lem:Base} (a) (with $\widehat{E}=\widehat{E}_{i}$
and $\widehat{\mathcal{F}}=\widehat{\mathcal{F}}_{i}$), $\widehat{\mathcal{F}}_{i}\subset C(\widehat{E}_{i};\mathbf{R})$
can be written as $\widehat{\mathcal{F}}_{i}=\{\widehat{f}_{n}^{i}\}_{n\in\mathbf{N}}$
for each $i\in\{1,2\}$ such that $\widehat{f}_{n}^{1}|_{E_{0}}=f_{n}|_{E_{0}}=\widehat{f}_{n}^{2}|_{E_{0}}$
for all $n\in\mathbf{N}$. Then, $(E_{0},\rho_{\widehat{\mathcal{F}}_{1}})$
and $(E_{0},\rho_{\widehat{\mathcal{F}}_{2}})$ are identical metric
spaces. Now, the corollary follows by Lemma \ref{lem:Base} (c) (with
$\widehat{E}=\widehat{E}_{i}$ and $\widehat{\mathcal{F}}=\widehat{\mathcal{F}}_{i}$)
and Proposition \ref{prop:Completeness} (a) (with $E=\widehat{E}_{1}$
and $S=\widehat{E}_{2}$).\end{proof}

\begin{rem}
\label{rem:Compactification}Compactification determined by \textit{extending}%
\footnote{Recall that mapping $g$ is a continuous extension of mapping $f$
if $g$ is continuous, the domain of $g$ contains that of $f$ and
$g=f$ restricted to $f$'s domain.%
} continuous functions was used in e.g. \cite{EK86}, \cite{BK93b},
\cite{BK10} and \cite{K15} to imbed stochastic processes into compact
metric spaces. The well-known \textit{Stone-$\check{\mbox{C}}$ech
compactification} (see p.\pageref{Stone-Cech}) exists for any compactifiable,
or equivalently, Tychonoff space and is determined by the continuous
extension of all bounded continuous function. The compactification
$\widehat{E}$ of $(E_{0},\mathscr{O}_{\mathcal{F}}(E_{0}))$ can
be thought of as a ``possibly smaller'' compactification which might
not extend all of $C_{b}(E_{0},\mathscr{O}_{\mathcal{F}}(E_{0});\mathbf{R})$.
\end{rem}

\begin{rem}
\label{rem:Change_Topo}As a cost of metrizability, $\widehat{E}$
\textit{compactifies} (see p.\pageref{Compactification}) $E_{0}$
with respect to a possibly \textit{coarser} topology than its natural
subspace topology induced from $E$. In fact, neither the original
space $E$ nor the subspace $(E_{0},\mathscr{O}_{E}(E_{0}))$ is necessarily
a Tychonoff space, hence need not have compactification.
\end{rem}

\begin{rem}
\label{rem:OnePointComp}\textit{One-point compactifications} (see
p.\pageref{One-Point}) exist for locally compact Hausdorff spaces
(see Proposition \ref{prop:One-Point}). We do not presume $E_{0}$
to be a locally compact subspace of $\widehat{E}$, and $\widehat{E}$
is not necessarily a one-point compactification. Nonetheless, even
Stone-$\check{\mbox{C}}$ech compactifications are sometimes one-point
compactifcations. Corollary \ref{prop:LC_Polish} to follow later
illustrates when the compactification establishing a base is of one-point
compactification type.
\end{rem}

The following theorem gives a partial converse of Lemma \ref{lem:Base}
(e) and answers the question \textit{when a base exists.}
\begin{thm}
\label{thm:Base}Let $E$ be a topological space, $E_{0}\in\mathscr{B}(E)$
and $\mathcal{F}\ni1$ be a countable subset of $C_{b}(E;\mathbf{R})$.
Then, there exists a base $(E_{0},\mathcal{F};\widehat{E},\widehat{\mathcal{F}})$
over $E$ if and only if $\mathcal{F}$ separates points on $E_{0}$.
\end{thm}
\begin{proof}
Necessity follows by Lemma \ref{lem:Base} (e). We prove sufficiency.
There exist a compactification $\widehat{E}$ of $(E_{0},\mathscr{O}_{\mathcal{F}}(E_{0}))$
and an extension $\varphi\in\mathbf{imb}(\widehat{E};\mathbf{R}^{\infty})$
of $\bigotimes\mathcal{F}|_{E_{0}}$ by Lemma \ref{lem:Compactification}
(a, b) (with $E=(E_{0},\mathscr{O}_{\mathcal{F}}(E_{0}))$ and $\mathcal{D}=\mathcal{F}|_{E_{0}}$).
Then, $(E_{0},\mathcal{F};\widehat{E},\widehat{\mathcal{F}})$ is
base over $E$ with $\widehat{\mathcal{F}}\circeq\{\mathfrak{p}_{n}\circ\varphi\}_{n\in\mathbf{N}}$.\end{proof}

\subsection{\label{sub:Base_Property}Properties}

The following four results specify basic finite-dimensional properties
of bases.
\begin{lem}
\label{lem:Base_Property}Let $E$ be a topological space, $(E_{0},\mathcal{F};\widehat{E},\widehat{\mathcal{F}})$
be a base over $E$, $d\in\mathbf{N}$ and $A\subset E_{0}^{d}$.
Then:

\renewcommand{\labelenumi}{(\alph{enumi})}
\begin{enumerate}
\item $(E_{0}^{d},\Pi^{d}(\mathcal{F});\widehat{E}^{d},\Pi^{d}(\widehat{\mathcal{F}}))$%
\footnote{The notation ``$\Pi^{d}(\cdot)$'' was defined in \S \ref{sub:Fun}.%
} is a base over $E^{d}$.
\item $\Pi^{d}(\widehat{\mathcal{F}})\subset C(\widehat{E}^{d};\mathbf{R})$
contains the constant function $1$, separates points and strongly
separates points on $\widehat{E}^{d}$. In particular,
\begin{equation}
\bigotimes\Pi^{d}(\widehat{\mathcal{F}})\in\mathbf{imb}\left(\widehat{E}^{d};\mathbf{R}^{\infty}\right).\label{eq:Base^d_Imb}
\end{equation}

\item $\widehat{E}^{d}$ is a compactification of $(E_{0}^{d},\mathscr{O}_{\mathcal{F}}(E_{0})^{d})$%
\footnote{As mentioned in \S \ref{sub:Topo} and \S \ref{sub:Prod_Space},
$\mathscr{O}_{\mathcal{F}}(E_{0})^{d}$ denotes the product topology
on $E_{0}^{d}$ where $E_{0}$ is equipped with the topology $\mathscr{O}_{\mathcal{F}}(E_{0})$.%
}. Moreover, $\widehat{E}^{d}$ is a Polish space and is completely
metrized by $\rho_{\widehat{\mathcal{F}}}^{d}$%
\footnote{$\rho_{\widehat{\mathcal{F}}}^{d}$ is defined by (\ref{eq:TF_Metric^d})
with $\mathcal{D}=\widehat{\mathcal{F}}$.%
}.
\item $\bigotimes\Pi^{d}(\mathcal{F})|_{A}\in\mathbf{imb}(A,\mathscr{O}_{\widehat{E}^{d}}(A);\mathbf{R}^{\infty})$.
Moreover, $(A,\mathscr{O}_{\widehat{E}^{d}}(A))$ is a metrizable
and separable coarsening of $(A,\mathscr{O}_{E^{d}}(A))$.
\item $\Pi^{d}(\mathcal{F}\backslash\{1\})$ separates points on the Hausdorff
space $(A,\mathscr{O}_{E^{d}}(A))$.
\end{enumerate}
\end{lem}
\begin{proof}
(a) We verify the four properties of Definition \ref{def:Base} in
four steps:

\textit{Step 1}. We have by $E_{0}\in\mathscr{B}(E)$ and Fact \ref{fact:Prod_Map_1}
(a) that
\begin{equation}
E_{0}^{d}=\bigcap_{i=1}^{d}\mathfrak{p}_{i}^{-1}(E_{0})\in\mathscr{B}(E)^{\otimes d},\;\forall d\in\mathbf{N}.\label{eq:E0d_Prod_Measurable_Ed}
\end{equation}

\textit{Step 2}. We have by $1\in\mathcal{F}\subset C_{b}(E;\mathbf{R})$,
Fact \ref{fact:Pi^d} (a) (with $\mathcal{D}=\mathcal{F}$) and Proposition
\ref{prop:Pi^d_SP} (a) (with $\mathcal{D}=\mathcal{F}$) that $\Pi^{d}(\mathcal{F})$
is a countable collection and
\begin{equation}
1\in\mathfrak{ca}\left(\Pi^{d}(\mathcal{F})\right)\subset C_{b}\left(E^{d},\mathscr{O}_{\mathcal{F}}(E)^{d};\mathbf{R}\right)\subset C_{b}(E^{d};\mathbf{R}).\label{eq:ca(Pi^d(F))_Cb}
\end{equation}

\textit{Step 3}. We have by (\ref{eq:F_Fhat_Coincide}) that
\begin{equation}
\left.\bigotimes\Pi^{d}(\mathcal{F})\right|_{E_{0}^{d}}=\left.\bigotimes\Pi^{d}(\widehat{\mathcal{F}})\right|_{E_{0}^{d}}.\label{eq:Pi^d(F)_Pi^d(Fhat)_Coincide}
\end{equation}
(\ref{eq:Fhat_SSP_on_Ehat}) and Proposition \ref{prop:Pi^d_SP} (a)
(with $E=\widehat{E}$ and $\mathcal{D}=\widehat{\mathcal{F}}$) imply
that
\begin{equation}
\Pi^{d}(\widehat{\mathcal{F}})\subset C(\widehat{E}^{d};\mathbf{R}).\label{eq:Pi^d(Fhat)_Cb}
\end{equation}
Then, $\Pi^{d}(\widehat{\mathcal{F}})$ separates points on $\widehat{E}^{d}$,
strongly separates points on $\widehat{E}^{d}$ and satisfies%
\footnote{In contrast, $\mathscr{O}(E^{d})$ is not necessarily the same as
$\mathscr{O}_{\mathcal{F}}(E)^{d}$.%
}
\begin{equation}
\mathscr{O}_{\Pi^{d}(\widehat{\mathcal{F}})}(\widehat{E}^{d})=\mathscr{O}(\widehat{E}^{d})=\mathscr{O}_{\widehat{\mathcal{F}}}(\widehat{E})^{d}\label{eq:Pi^d(Fhat)_SSP_Ehatd}
\end{equation}
by Lemma \ref{lem:Base} (a), Proposition \ref{prop:Pi^d_SP} (b)
(with $E=\widehat{E}$ and $\mathcal{D}=\widehat{\mathcal{F}}$),
(\ref{eq:Pi^d(Fhat)_Cb}) and Proposition \ref{prop:Fun_Sep_1} (b,
e) (with $E=A=\widehat{E}^{d}$ and $\mathcal{D}=\Pi^{d}(\widehat{\mathcal{F}})$).

\textit{Step 4}. $E_{0}^{d}$ is dense in $\widehat{E}^{d}$ by the
denseness of $E_{0}$ in $\widehat{E}$ and the product topology definition.
$\widehat{E}^{d}$ is compact by Lemma \ref{lem:Base} (b) and the
Tychonoff Theorem (see Proposition \ref{prop:Compact} (b)). $\bigotimes\Pi^{d}(\widehat{\mathcal{F}})\in C(\widehat{E}^{d};\mathbf{R}^{\infty})$
by (\ref{eq:Pi^d(Fhat)_Cb}) and Fact \ref{fact:Prod_Map_2} (b).
\begin{equation}
\bigotimes\Pi^{d}(\widehat{\mathcal{F}})(\widehat{E}^{d})\in\mathscr{K}(\mathbf{R}^{\infty})\subset\mathscr{C}(\mathbf{R}^{\infty})\label{eq:Pi^d(Fhat)(Ehatd)_Compact_Closed}
\end{equation}
by Proposition \ref{prop:Compact} (a, e). So, $\bigotimes\Pi^{d}(\widehat{\mathcal{F}})(\widehat{E}^{d})$
is the closure of $\bigotimes\Pi^{d}(\mathcal{F})(E_{0}^{d})$ in
$\mathbf{R}^{\infty}$ by (\ref{eq:Pi^d(F)_Pi^d(Fhat)_Coincide})
and \cite[Theorem 18.1 (a, b)]{M00}.

(b) follows by (a) and Lemma \ref{lem:Base} (a).

(c) The first statement follows by (a) and Lemma \ref{lem:Base} (b).
$\widehat{E}^{d}$ is a Polish space by (a) and Lemma \ref{lem:Base}
(c). Moreover, $\rho_{\widehat{\mathcal{F}}}^{d}$ metrizes $\widehat{E}^{d}$
by Lemma \ref{lem:Base} (c) and Proposition \ref{prop:Metric_Prod}
(a) (with $\mathbf{I}=\{1,...,d\}$ and $(E_{i},\mathfrak{r}_{i})=(\widehat{E},\rho_{\widehat{\mathcal{F}}})$).

(d) follows by (a) and Lemma \ref{lem:Base} (d).

(e) $\Pi^{d}(\mathcal{F}\backslash\{1\})$ separates points on $E_{0}^{d}$
by Lemma \ref{lem:Base} (e) and Proposition \ref{prop:Pi^d_SP} (b)
(with $\mathcal{D}=\mathcal{F}\backslash\{1\}$). The rest of (e)
follows by (a) and Lemma \ref{lem:Base} (e).\end{proof}

\begin{cor}
\label{cor:Base_Fun_Dense}Let $E$ be a topological space, $(E_{0},\mathcal{F};\widehat{E},\widehat{\mathcal{F}})$
be a base over $E$ and $d\in\mathbf{N}$. Then%
\footnote{The notation ``$\mathfrak{cl}(\cdot)$'', ``$\mathfrak{ag}_{\mathbf{Q}}(\cdot)$''
and ``$\mathfrak{ag}(\cdot)$'' were defined in \S \ref{sub:Fun}.%
},
\begin{equation}
C_{c}(\widehat{E}^{d};\mathbf{R})=C_{0}(\widehat{E}^{d};\mathbf{R})=C_{b}(\widehat{E}^{d};\mathbf{R})=C(\widehat{E}^{d};\mathbf{R})=\mathfrak{cl}\left[\mathfrak{ag}_{\mathbf{Q}}\left(\Pi^{d}(\widehat{\mathcal{F}})\right)\right].\label{eq:Q-Algebra_RepTF}
\end{equation}
and
\begin{equation}
\mathfrak{cl}\left[\mathfrak{ag}_{\mathbf{Q}}\left(\Pi^{d}\left(\mathcal{F}|_{E_{0}}\right)\right)\right]=\left.C(\widehat{E}^{d};\mathbf{R})\right|_{E_{0}^{d}}.\label{eq:ca(Pi^d(Ftilte))_C(Ehat)}
\end{equation}

\end{cor}
\begin{proof}
$\mathfrak{ag}(\Pi^{d}(\widehat{\mathcal{F}}))$ is uniformly dense%
\footnote{``Uniformly dense'' means dense with respect to the topology induced
by the supremum norm.%
} in $C(\widehat{E}^{d};\mathbf{R})$ by Lemma \ref{lem:Base_Property}
(b, c) and the Stone-Weierstrass Theorem (see \cite[Theorem 2.4.11]{D02}).
Thus, (\ref{eq:Q-Algebra_RepTF}) follows by Lemma \ref{lem:Base_Property}
(c), Fact \ref{fact:Cc_C0_Cb} (with $E=\widehat{E}$ and $k=1$)
and (\ref{eq:ca(D)}) (with $\mathcal{D}=\Pi^{d}(\widehat{\mathcal{F}})$).
(\ref{eq:ca(Pi^d(Ftilte))_C(Ehat)}) follows by (\ref{eq:Q-Algebra_RepTF}),
the denseness of $E_{0}^{d}$ in $\widehat{E}^{d}$ and properties
of uniform convergence.\end{proof}

\begin{cor}
\label{cor:Base_Sep_Meas}Let $E$ be a topological space, $(E_{0},\mathcal{F};\widehat{E},\widehat{\mathcal{F}})$
be a base over $E$, $d\in\mathbf{N}$, $A\subset\widehat{E}^{d}$
and $\mathcal{D}\circeq\mathfrak{mc}[\Pi^{d}(\widehat{\mathcal{F}}\backslash\{1\})]$%
\footnote{The notation ``$\mathfrak{mc}(\cdot)$'' was defined in \S \ref{sub:Fun}.%
}. Then:

\renewcommand{\labelenumi}{(\alph{enumi})}
\begin{enumerate}
\item $\mathcal{D}|_{A}^{*}$%
\footnote{The notation ``$\mathcal{D}^{*}$'' and the terminologies ``separating'',
``convergence determining'' were defined in \S \ref{sec:Borel_Measure}.%
} separates and strongly separates points on $\mathcal{P}(A,\mathscr{O}_{\widehat{E}^{d}}(A))$.
\item \textup{$\mathcal{D}|_{A}\cup\{1\}$ }(especially $\mathfrak{mc}[\Pi^{d}(\widehat{\mathcal{F}})]|_{A}$)
is separating and convergence determining on $(A,\mathscr{O}_{\widehat{E}^{d}}(A))$.
\item $\mathcal{M}^{+}(\widehat{E}^{d})$ and $\mathcal{P}(\widehat{E}^{d})$
are Polish spaces and, in particular, $\mathcal{P}(\widehat{E}^{d})$
is compact.
\end{enumerate}
\end{cor}
\begin{proof}
(a) follows by Lemma \ref{lem:Base} (a) and Lemma \ref{lem:Sep_Meas}
(b) (with $E=(A,\mathscr{O}_{\widehat{E}^{d}}(A))$).

(b) follows by (a), Fact \ref{fact:Sep_CD} and the fact $(\mathcal{D}\cup\{1\})\subset\mathfrak{mc}[\Pi^{d}(\widehat{\mathcal{F}})]$.

(c) follows by Lemma \ref{lem:Base_Property} (c) and Theorem \ref{thm:P(E)_Compact_Polish}
(with $E=\widehat{E}^{d}$).\end{proof}

\begin{fact}
\label{fact:Ehat_Borel_Prod}Let $E$ be a topological space, $(E_{0},\mathcal{F};\widehat{E},\widehat{\mathcal{F}})$
be a base over $E$, $d\in\mathbf{N}$ and $A\subset\widehat{E}^{d}$.
Then, $\mathscr{B}_{\widehat{E}^{d}}(A)=\mathscr{B}(\widehat{E})^{\otimes d}|_{A}$%
\footnote{The notations ``$\mathscr{B}(\widehat{E})^{\otimes d}$'' and ``$\left.\mathscr{B}(\widehat{E})^{\otimes d}\right|_{A}$''
were defined in \S \ref{sub:Prod_Space}.%
}.
\end{fact}
\begin{proof}
This fact follows immediately by Lemma \ref{lem:Base} (c) and Proposition
\ref{prop:Prod_Space} (d) (with $\mathbf{I}=\{1,...,d\}$ and $S_{i}=\widehat{E}$).\end{proof}

\begin{cor}
\label{cor:Base_Borel}Let $E$ be a topological space, $(E_{0},\mathcal{F};\widehat{E},\widehat{\mathcal{F}})$
be a base over $E$ and $d\in\mathbf{N}$. Then:

\renewcommand{\labelenumi}{(\alph{enumi})}
\begin{enumerate}
\item Any $A\subset E_{0}^{d}$ satisfies%
\footnote{The notation ``$\mathscr{B}_{\mathcal{F}}(E)$'' was defined in
\S \ref{sub:Topo}.%
}
\begin{equation}
\begin{aligned}\mathscr{B}_{\left(E^{d},\mathscr{O}_{\mathcal{F}}(E)^{d}\right)}(A) & =\left.\mathscr{B}(\widehat{E})^{\otimes d}\right|_{A}=\mathscr{B}_{\widehat{E}^{d}}(A)=\mathscr{B}_{\Pi^{d}(\mathcal{F})}(A)=\left.\mathscr{B}_{\mathcal{F}}(E)^{\otimes d}\right|_{A}\\
 & \subset\left.\mathscr{B}(E)^{\otimes d}\right|_{A}\subset\mathscr{B}_{E^{d}}(A).
\end{aligned}
\label{eq:Borel_Compare_1}
\end{equation}
In particular, $E_{0}^{d}$ satisfies
\begin{equation}
\mathscr{B}_{\widehat{E}^{d}}(E_{0}^{d})=\mathscr{B}_{\mathcal{F}}(E_{0})^{\otimes d}\subset\mathscr{B}_{E}(E_{0})^{\otimes d}\subset\mathscr{B}(E)^{\otimes d}\subset\mathscr{B}(E^{d})\label{eq:Borel_Compare_2}
\end{equation}

\item $\mathcal{F}$ satisfies
\begin{equation}
\begin{aligned}\mathfrak{ca}\left(\Pi^{d}(\mathcal{F})\right) & \subset M_{b}\left(E^{d},\mathscr{B}_{\mathcal{F}}(E)^{\otimes d};\mathbf{R}\right)\\
 & \subset M_{b}\left(E^{d},\mathscr{B}(E)^{\otimes d};\mathbf{R}\right)\subset M_{b}(E^{d};\mathbf{R}).
\end{aligned}
\label{eq:ca(Pi^d(F))_Mb}
\end{equation}

\end{enumerate}
\end{cor}
\begin{proof}
(a) We have that
\begin{equation}
\begin{aligned} & \mathscr{B}_{\left(E^{d},\mathscr{O}_{\mathcal{F}}(E)^{d}\right)}(A)=\mathscr{B}_{(E_{0}^{d},\mathscr{O}_{\mathcal{F}}(E_{0})^{d})}(A)\\
 & =\mathscr{B}_{(E_{0}^{d},\mathscr{O}_{\widehat{E}}(E_{0})^{d})}(A)=\mathscr{B}_{\widehat{E}^{d}}(A)=\left.\mathscr{B}(\widehat{E})^{\otimes d}\right|_{A}\\
 & =\left.\mathscr{B}_{\widehat{E}}(E_{0})^{\otimes d}\right|_{A}=\left.\mathscr{B}_{\mathcal{F}}(E_{0})^{\otimes d}\right|_{A}=\left.\mathscr{B}_{\mathcal{F}}(E)^{\otimes d}\right|_{A}
\end{aligned}
\label{eq:Check_Borel_Compare_1}
\end{equation}
by (\ref{eq:Check_E0_Topo}), Fact \ref{fact:Ehat_Borel_Prod} and
the fact $A\subset E_{0}^{d}\subset(E^{d}\cap\widehat{E}^{d})$. Then,
the first line of (\ref{eq:Borel_Compare_1}) follows by (\ref{eq:Check_Borel_Compare_1}),
(\ref{eq:Pi^d(F)_Pi^d(Fhat)_Coincide}) and Lemma \ref{lem:Base_Property}
(b). The second line of (\ref{eq:Borel_Compare_1}) follows by $\mathcal{F}\subset C(E;\mathbf{R})$
and Lemma \ref{prop:Prod_Space} (a). Now, (a) follows by (\ref{eq:E0d_Prod_Measurable_Ed}).

(b) follows by $\mathcal{F}\subset C_{b}(E;\mathbf{R})$, Proposition
\ref{prop:Pi^d_SP} (a) (with $\mathcal{D}=\mathcal{F}$) and (\ref{eq:Borel_Compare_2}).\end{proof}

$E_{0}$ is not necessarily a Borel subset of $\widehat{E}$ and $\widehat{E}$
does not endow $E_{0}$ with the same Borel sets as $E$ unless $E_{0}$
is a standard Borel subset of $E$.
\begin{lem}
\label{lem:SB_Base}Let $E$ be a topological space, $(E_{0},\mathcal{F};\widehat{E},\widehat{\mathcal{F}})$
be a base over $E$, $d\in\mathbf{N}$ and $A\subset E_{0}^{d}$.
Then:

\renewcommand{\labelenumi}{(\alph{enumi})}
\begin{enumerate}
\item $A\in\mathscr{B}^{\mathbf{s}}(E^{d})$%
\footnote{$\mathscr{B}^{\mathbf{s}}(E^{d})$ denotes the family of all standard
Borel subsets of $E^{d}$.%
} if and only if%
\footnote{$\mathscr{B}_{E^{d}}(A)=\mathscr{B}_{\widehat{E}^{d}}(A)$ plus $A\in\mathscr{B}(\widehat{E}^{d})$
is equivalent to $\mathscr{B}_{E^{d}}(A)=\mathscr{B}_{\widehat{E}^{d}}(A)\subset\mathscr{B}(\widehat{E}^{d})$.
Hereafter, we frequently use the latter notation.%
}
\begin{equation}
\mathscr{B}_{E^{d}}(A)=\mathscr{B}_{\widehat{E}^{d}}(A)\mbox{ and }A\in\mathscr{B}(\widehat{E}^{d}).\label{eq:A_SB_Borel_Equal_d}
\end{equation}

\item If $A\in\mathscr{B}^{\mathbf{s}}(E^{d})$, then
\begin{equation}
\mathscr{B}_{E^{d}}(A)=\left.\mathscr{B}(E)^{\otimes d}\right|_{A}\subset\mathscr{B}(E)^{\otimes d}\subset\mathscr{B}(E^{d}).\label{eq:A_SB_Borel_Prod_Equal}
\end{equation}

\item If $A=\bigcup_{n\in\mathbf{N}}A_{n}$ and $\{A_{n}\}_{n\in\mathbf{N}}\subset\mathscr{B}^{\mathbf{s}}(E^{d})$,
then $A\in\mathscr{B}^{\mathbf{s}}(E^{d})$.
\item If $A=\prod_{n=1}^{d}A_{n}$ and $\{A_{n}\}_{1\leq n\leq d}\subset\mathscr{B}^{\mathbf{s}}(E)$,
then
\begin{equation}
\mathscr{B}_{\widehat{E}^{d}}(A)=\left.\mathscr{B}(E)^{\otimes d}\right|_{A}\subset\left[\mathscr{B}(E)^{\otimes d}\cap\mathscr{B}(\widehat{E}^{d})\right].\label{eq:A_n_Prod_Borel_Ehatd}
\end{equation}

\end{enumerate}
\end{lem}
\begin{proof}
(a - Necessity) Let $f$ be the identity map on $A$, which is certainly
injective. $f\in M(A,\mathscr{B}_{E^{d}}(A);\widehat{E}^{d})$ by
(\ref{eq:Borel_Compare_1}). $A\in\mathscr{B}^{\mathbf{s}}(E^{d})$
and $A\subset E_{0}^{d}$ imply $A\in\mathscr{B}^{\mathbf{s}}(A,\mathscr{O}_{E^{d}}(A))$.
$\widehat{E}^{d}$ is a Polish space by Lemma \ref{lem:Base_Property}
(c) so $A\in\mathscr{B}^{\mathbf{s}}(\widehat{E}^{d})$ and $\mathscr{B}_{E^{d}}(A)=\mathscr{B}_{\widehat{E}^{d}}(A)$
by Proposition \ref{prop:SB_Map} (with $E=A$ and $S=\widehat{E}^{d}$).
Then, $A\in\mathscr{B}(\widehat{E}^{d})=\mathscr{B}^{\mathbf{s}}(\widehat{E}^{d})$
by Proposition \ref{prop:SB_Borel} (b) (with $E=\widehat{E}^{d}$).

(a - Sufficiency) follows by (\ref{eq:A_SB_Borel_Equal_d}) and Fact
\ref{fact:Polish_SB} (a).

(b) $A\in\mathscr{B}_{\widehat{E}^{d}}(E_{0}^{d})$ by the fact $A\subset E_{0}^{d}$
and (\ref{eq:A_SB_Borel_Equal_d}). $\mathscr{B}_{\widehat{E}^{d}}(A)\subset\mathscr{B}(E)^{\otimes d}\subset\mathscr{B}(E^{d})$
by (\ref{eq:Borel_Compare_2}). Now, (b) follows by (a) and (\ref{eq:Borel_Compare_1}).

(c) We find by (a) (with $A=A_{n}$) that
\begin{equation}
\mathscr{B}_{E^{d}}(A_{n})=\mathscr{B}_{\widehat{E}^{d}}(A_{n})\subset\mathscr{B}(\widehat{E}^{d}),\;\forall n\in\mathbf{N}.\label{eq:A_n_Borel_Equal_Ed}
\end{equation}
Then, $A$ satisfies (\ref{eq:A_SB_Borel_Equal_d}) by (\ref{eq:A_n_Borel_Equal_Ed}),
Fact \ref{fact:Union_Borel} (with $E=A$, $\mathscr{U}_{1}=\mathscr{B}_{E^{d}}(A)$
and $\mathscr{U}_{2}=\mathscr{B}_{\widehat{E}^{d}}(A)$) and (\ref{eq:Borel_Compare_1}).
Hence, $A\in\mathscr{B}^{\mathbf{s}}(E^{d})$ by (a).

(d) $A\subset E_{0}^{d}$ implies $A_{n}=\mathfrak{p}_{n}(A)\subset E_{0}$
for all $1\leq n\leq d$. We have that
\begin{equation}
\mathscr{B}_{E}(A_{n})=\mathscr{B}_{\widehat{E}}(A_{n})\subset\left[\mathscr{B}(E)\cap\mathscr{B}(\widehat{E})\right],\;\forall1\leq n\leq d\label{eq:A_n_Borel_Equal_E}
\end{equation}
by (a, b) (with $d=1$ and $A=A_{n}$). It then follows by (\ref{eq:A_n_Borel_Equal_E})
that
\begin{equation}
A=\bigcap_{n=1}^{d}\mathfrak{p}_{n}^{-1}(A_{n})\in\mathscr{B}(E)^{\otimes d}\cap\mathscr{B}(\widehat{E}^{d}),\label{eq:Check_A_n_Prod_Borel_Ehatd_1}
\end{equation}
and by Corollary \ref{cor:Base_Borel} (a) (with $d=1$ and $A=A_{n}$)
that
\begin{equation}
\mathscr{B}_{\widehat{E}^{d}}(A)=\left.\mathscr{B}(\widehat{E})^{\otimes d}\right|_{A}=\bigotimes_{n=1}^{d}\mathscr{B}_{\widehat{E}}(A_{n})=\bigotimes_{n=1}^{d}\mathscr{B}_{E}(A_{n})=\left.\mathscr{B}(E)^{\otimes d}\right|_{A}.\label{eq:Check_A_n_Prod_Borel_Ehatd_2}
\end{equation}
\end{proof}

\begin{cor}
\label{cor:Base_Compact}Let $E$ be a topological space, $(E_{0},\mathcal{F};\widehat{E},\widehat{\mathcal{F}})$
be a base over $E$, $d\in\mathbf{N}$ and $A\subset E_{0}^{d}$.
Then:

\renewcommand{\labelenumi}{(\alph{enumi})}
\begin{enumerate}
\item If $A\in\mathscr{K}(E_{0}^{d},\mathscr{O}_{E}(E_{0})^{d})$, then
\begin{equation}
\left(A,\mathscr{O}_{E^{d}}(A)\right)=\left(A,\rho_{\mathcal{F}}^{d}\right)=\left(A,\mathscr{O}_{\widehat{E}^{d}}(A)\right)\label{eq:Topo_Coincide_on_K_1}
\end{equation}
and\textup{}%
\footnote{The notations ``$\mathscr{C}(\cdot)$'', ``$\mathscr{K}(\cdot)$''
and ``$\mathscr{K}^{\mathbf{m}}(\cdot)$'', defined in \S \ref{sub:Topo},
denote the families of closed, compact and metrizable compact subsets,
respectively.%
}
\begin{equation}
A\in\mathscr{K}^{\mathbf{m}}(E^{d})\cap\mathscr{K}(\widehat{E}^{d})\cap\mathscr{C}(\widehat{E}^{d}).\label{eq:Topo_Coincide_on_K_2}
\end{equation}

\item If $A\in\mathscr{K}_{\sigma}(E_{0}^{d},\mathscr{O}_{E}(E_{0})^{d})$,
then
\begin{equation}
A\in\mathscr{B}^{\mathbf{s}}(E^{d})\cap\mathscr{K}_{\sigma}^{\mathbf{m}}(E^{d})\cap\mathscr{K}_{\sigma}(\widehat{E}^{d})\cap\mathscr{B}(E)^{\otimes d}\cap\mathscr{B}(\widehat{E}^{d}).\label{eq:E0_Sig_Compact_SB}
\end{equation}

\end{enumerate}
\end{cor}
\begin{proof}
(a) follows by Lemma \ref{lem:Base_Property} (b, c, e) and Fact \ref{fact:Compact_Topo_Coarsen}
(b) (with $E=E^{d}$, $\mathcal{D}=\Pi^{d}(\mathcal{F})$ and $K=A$).

(b) $A\in\mathscr{K}_{\sigma}^{\mathbf{m}}(E^{d})\cap\mathscr{K}_{\sigma}(\widehat{E}^{d})\cap\mathscr{B}(\widehat{E}^{d})$
by (\ref{eq:Topo_Coincide_on_K_2}). $\mathscr{K}(E_{0}^{d},\mathscr{O}_{E}(E_{0})^{d})\subset\mathscr{B}^{\mathbf{s}}(E^{d})$
by (a) and Lemma \ref{lem:SB_Base} (a). $A\in\mathscr{B}^{\mathbf{s}}(E^{d})$
by Lemma \ref{lem:SB_Base} (c). Moreover, $A\in\mathscr{B}(E)^{\otimes d}$
by Lemma \ref{lem:SB_Base} (b).\end{proof}

\begin{note}
\label{note:Base_Sigma_Compact_in_Domain}Given a base $(E_{0},\mathcal{F};\widehat{E},\widehat{\mathcal{F}})$
over $E$ and $d\in\mathbf{N}$, Corollary \ref{cor:Base_Compact}
(b) shows that $\mathscr{K}_{\sigma}(E_{0}^{d},\mathscr{O}_{E}(E_{0})^{d})$
lies in the domain of any $\mu\in\mathfrak{M}^{+}(E^{d},\mathscr{B}(E)^{\otimes d})$.\end{note}
\begin{cor}
\label{cor:SB_Base_Borel_Extension}Let $E$ be a topological space,
$(E_{0},\mathcal{F};\widehat{E},\widehat{\mathcal{F}})$ be a base
over $E$, $d\in\mathbf{N}$ and $\mu\in\mathfrak{M}^{+}(E^{d},\mathscr{B}(E)^{\otimes d})$.

\renewcommand{\labelenumi}{(\alph{enumi})}
\begin{enumerate}
\item If $\mu$ is supported on $A\subset E_{0}^{d}$%
\footnote{Support of measure was specified in \S \ref{sub:Meas}.%
} and $A\in\mathscr{B}^{\mathbf{s}}(E^{d})$, then $\mathfrak{be}(\mu)$%
\footnote{$\mathfrak{be}(\mu)$ as defined in \S \ref{sec:Borel_Measure} denotes
the family of all Borel extensions of $\mu$.%
} is a singleton.
\item If $\mu$ is tight in $(E_{0}^{d},\mathscr{O}_{E}(E_{0})^{d})$, then
there exists a $\mu^{\prime}=\mathfrak{be}(\mu)$ which is tight in
$(E_{0}^{d},\mathscr{O}_{E}(E_{0})^{d})$.
\end{enumerate}
\end{cor}
\begin{note}
\label{note:BExt}Any $\mu^{\prime}\in\mathfrak{be}(\mu)$ is \textit{not
an expansion} of $\mu$ to a superspace but rather an \textit{extension}
of $\mu$ to all Borel sets. Any support of $\mu$ is also that of
$\mu^{\prime}$.
\end{note}
\begin{proof}
[Proof of Corollary \ref{cor:SB_Base_Borel_Extension}](a) One finds
by Lemma \ref{lem:SB_Base} (b) that $A\in\mathscr{B}(E)^{\otimes d}$
and $\mathscr{B}_{E^{d}}(A)=\mathscr{B}(E)^{\otimes d}|_{A}$. Hence,
(a) follows by Lemma \ref{lem:Union_Borel_Prod_Equal} (c) (with $\mathbf{I}=\{1,...,d\}$,
$S_{i}=E$, $S=E^{d}$ and $\mathscr{A}=\mathscr{B}(E)^{\otimes d}$).

(b) $\mu$ is supported on some $A\in\mathscr{K}_{\sigma}(E_{0}^{d},\mathscr{O}_{E}(E_{0})^{d})$
by its tightness. $A\in\mathscr{B}^{\mathbf{s}}(E^{d})$ by Corollary
\ref{cor:Base_Compact} (b). Hence, (b) follows by (a).\end{proof}

The standard Borel property of $A\subset E_{0}^{d}$ also yields useful
properties of the weak topological space $\mathcal{M}^{+}(A,\mathscr{O}_{E^{d}}(A))$.
\begin{cor}
\label{cor:SB_Base_Sep_Meas}Let $E$ be a topological space, $(E_{0},\mathcal{F};\widehat{E},\widehat{\mathcal{F}})$
be a base over $E$, $d\in\mathbf{N}$, $A\in\mathscr{B}^{\mathbf{s}}(E^{d})$
with $A\subset E_{0}^{d}$ and $\mathcal{D}\circeq\mathfrak{mc}[\Pi^{d}(\mathcal{F}\backslash\{1\})]$.
Then:

\renewcommand{\labelenumi}{(\alph{enumi})}
\begin{enumerate}
\item $\mathcal{D}|_{A}^{*}$ separates points on $\mathcal{P}(A,\mathscr{O}_{E^{d}}(A))$.
Moreover, $\mathcal{D}|_{A}^{*}\cup\{1^{*}\}$ (especially $\mathfrak{mc}[\Pi^{d}(\mathcal{F})]|_{A}$)
is separating on $(A,\mathscr{O}_{E^{d}}(A))$.
\item $\mathcal{M}^{+}(A,\mathscr{O}_{E^{d}}(A))$ and $\mathcal{P}(A,\mathscr{O}_{E^{d}}(A))$
are Tychonoff spaces.
\end{enumerate}
\end{cor}
\begin{proof}
(a) $\mathcal{D}|_{A}^{*}$ separates points on $\mathcal{P}(A,\mathscr{O}_{\widehat{E}^{d}}(A))=\mathcal{P}(A,\mathscr{O}_{E^{d}}(A))$
by (\ref{eq:F_Fhat_Coincide}), Corollary \ref{cor:Base_Sep_Meas}
(a) and Lemma \ref{lem:SB_Base} (a). Now, (a) follows by Fact \ref{fact:Sep_CD}
(a) and the fact $(\mathcal{D}\cup\{1\})\subset\mathfrak{mc}[\Pi^{d}(\mathcal{F})]$.

(b) follows by (a) and Proposition \ref{prop:P(E)_CR} (a, c) (with
$E=(A,\mathscr{O}_{E^{d}}(A))$).\end{proof}

\section{\label{sec:Baseable_Space}Baseable space}

Theorem \ref{thm:Base} shows establishing a base is equivalent to
separating points of a subset of $E$ by countably many bounded continuous
functions. Boundedness is unecessary.
\begin{lem}
\label{lem:Baseable_Base}Let $E$ be a topological space, $E_{0}\in\mathscr{B}(E)$
and $\mathcal{D}\subset C(E;\mathbf{R})$ be countable and separate
points on $E$. Then, there exists a base $(E_{0},\mathcal{F};\widehat{E},\widehat{\mathcal{F}})$
over $E$ satisfying the following properties:

\renewcommand{\labelenumi}{(\alph{enumi})}
\begin{enumerate}
\item $\mathscr{O}_{\mathcal{D}}(E)\subset\mathscr{O}_{\mathcal{F}}(E)$.
\item $(\mathcal{D}\cap C_{b}(E;\mathbf{R}))\subset\mathcal{F}$.
\item $\mathcal{F}$ can be taken to equal $\mathcal{D}\cup\{1\}$ whenever
$\mathcal{D}\subset C_{b}(E;\mathbf{R})$.
\end{enumerate}
\end{lem}
\begin{proof}
$\{(f\wedge n)\vee(-n)\}_{n\in\mathbf{N},f\in\mathcal{D}\cup\{1\}}\subset C_{b}(E;\mathbf{R})$
so by Lemma \ref{lem:SP_Bounded} (with $\mathcal{G}=\mathcal{D}\cup\{1\}$
and $\mathcal{H}=C_{b}(E;\mathbf{R})$), there exists a countable
$\mathcal{F}\subset C_{b}(E;\mathbf{R})$ that satisfies (a) - (c).
$(E_{0},\mathscr{O}_{\mathcal{D}}(E_{0}))$ is a Hausdorff space by
Proposition \ref{prop:Fun_Sep_1} (c) (with $A=E_{0}$), hence $(E_{0},\mathscr{O}_{\mathcal{F}}(E_{0}))$
is also by (a) and Fact \ref{fact:Hausdorff_Refine}. $\mathcal{F}$
separates points by Proposition \ref{prop:Fun_Sep_1} (c) (with $A=E_{0}$
and $\mathcal{D}=\mathcal{F}$). The result follows by Theorem \ref{thm:Base}.\end{proof}

We use ``\textit{baseable}'' and ``\textit{baseability}'' to describe
the \textit{ability to create base}s.
\begin{defn}
\label{def:Baseable}Let $E$ be a topological space and $A\subset E$
be non-empty.
\begin{itemize}
\item $E$ is a \textbf{$\mathcal{D}$-baseable space} if $\mathcal{D}\subset C(E;\mathbf{R})$
has a countable subset separating points on $E$.
\item $A$ is a \textbf{$\mathcal{D}$-baseable subset of} $E$ if $A\in\mathscr{B}(E)$,
$\mathcal{D}\subset C(E;\mathbf{R})$ and $(A,\mathscr{O}_{E}(A))$
is a $\mathcal{D}|_{A}$-baseable space.
\item $E$ is a \textbf{baseable space} if $E$ is a $C(E;\mathbf{R})$-baseable
space.
\item $A$ is a \textbf{baseable subset} \textbf{of $E$} if $A$ is a $C(E;\mathbf{R})$-baseable
subset.
\end{itemize}
\end{defn}
\begin{rem}
\label{rem:Cont_Fun_Restrict}``being a baseable subset'' equals
``being a baseable subspace'' plus Borel measurability. Moreover,
the ``$\mathcal{D}|_{A}$-baseable space'' above is a proper statement
since
\begin{equation}
\mathcal{D}|_{A}\subset C(E;\mathbf{R})|_{A}\subset C(A,\mathscr{O}_{E}(A);\mathbf{R}).\label{eq:C(E;R)|A}
\end{equation}

\end{rem}

Baseable spaces or their analogues have appeared in many works such
as \cite[Chapter 3]{EK86}, \cite{J86}, \cite{KO88}, \cite{J97a},
\cite[Chapter 6]{B07}, \cite{BK10}, \cite{KS13} and \cite{K15}
etc. Herein, we are the first to characterize baseable topological
spaces. The next section will treat baseable subsets.

\subsection{\label{sub:Baseable_Space}Characterization}

Baseable spaces are a broad category of spaces.
\begin{thm}
\label{thm:Baseable_Space}Baseable spaces are precisely the topological
refinements of metrizable and separable spaces.
\end{thm}
The theorem above follows by two straightforward observations. First,
we note that baseable spaces sit between Hausdorff spaces and separable
metric spaces.
\begin{fact}
\label{fact:Baseable_Metrizable_Separable}The following statements
are true:

\renewcommand{\labelenumi}{(\alph{enumi})}
\begin{enumerate}
\item If $E$ is a baseable space, then $E$ is a Hausdorff space.
\item If $E$ is a metrizable and separable space, then $E$ is a $\mathcal{D}$-baseable
space for some countable $\mathcal{D}\subset C_{b}(E;\mathbf{R})$
that strongly separates points on $E$.
\end{enumerate}
\end{fact}
\begin{proof}
(a) follows by Proposition \ref{prop:Fun_Sep_1} (e) (with $A=E$).

(b) follows by Corollary \ref{cor:M_Compactification} (a, b).\end{proof}

\begin{cor}
\label{cor:Metrizable_Souslin_Baseable}Metrizable Souslin spaces,
metrizable Lusin spaces and metrizable standard Borel spaces are all
baseable spaces.
\end{cor}
\begin{proof}
Souslin and Lusin spaces are separable by Proposition \ref{prop:Var_Polish}
(d). Metrizable standard Borel spaces are Lusin spaces by Proposition
\ref{prop:Metrizable_SB_Lusin} (a, b). Now, the result follows by
Fact \ref{fact:Baseable_Metrizable_Separable} (b).\end{proof}

Refining the topology of a metrizable and separable space $E$ may
cause a function class $\mathcal{D}\subset C(E;\mathbf{R})$ to forfeit
strongly separation of points. However, it does not affect the $\mathcal{D}$-baseability
of $E$.
\begin{fact}
\label{fact:Baseable_Refine}If $E$ is $\mathcal{D}$-baseable, then
any topological refinement of $E$ is also.
\end{fact}
\begin{proof}
Note that $\mathcal{D}\subset C(E;\mathbf{R})\subset C(E,\mathscr{U};\mathbf{R})$
if $\mathscr{U}$ refines $\mathscr{O}(E)$.\end{proof}

\begin{proof}
[Proof of Theorem \ref{thm:Baseable_Space}](Necessity) If $E$ is
a baseable space and $\mathcal{D}\subset C(E;\mathbf{R})$ is countable
and separate points on $E$, then $(E,\mathscr{O}_{\mathcal{D}}(E))$
coarsens $E$ and is a metrizable and separable space by Proposition
\ref{prop:Fun_Sep_1} (d).

(Sufficiency) follows by Fact \ref{fact:Baseable_Metrizable_Separable}
(b) and Fact \ref{fact:Baseable_Refine}.\end{proof}

\subsection{\label{sub:Examples_Baseable_Space}Examples of baseable spaces}

The following figure illustrates the relationship of baseable spaces
and other major categories of topological spaces.

\begin{figure}[H]
\begin{centering}
\includegraphics[scale=0.85]{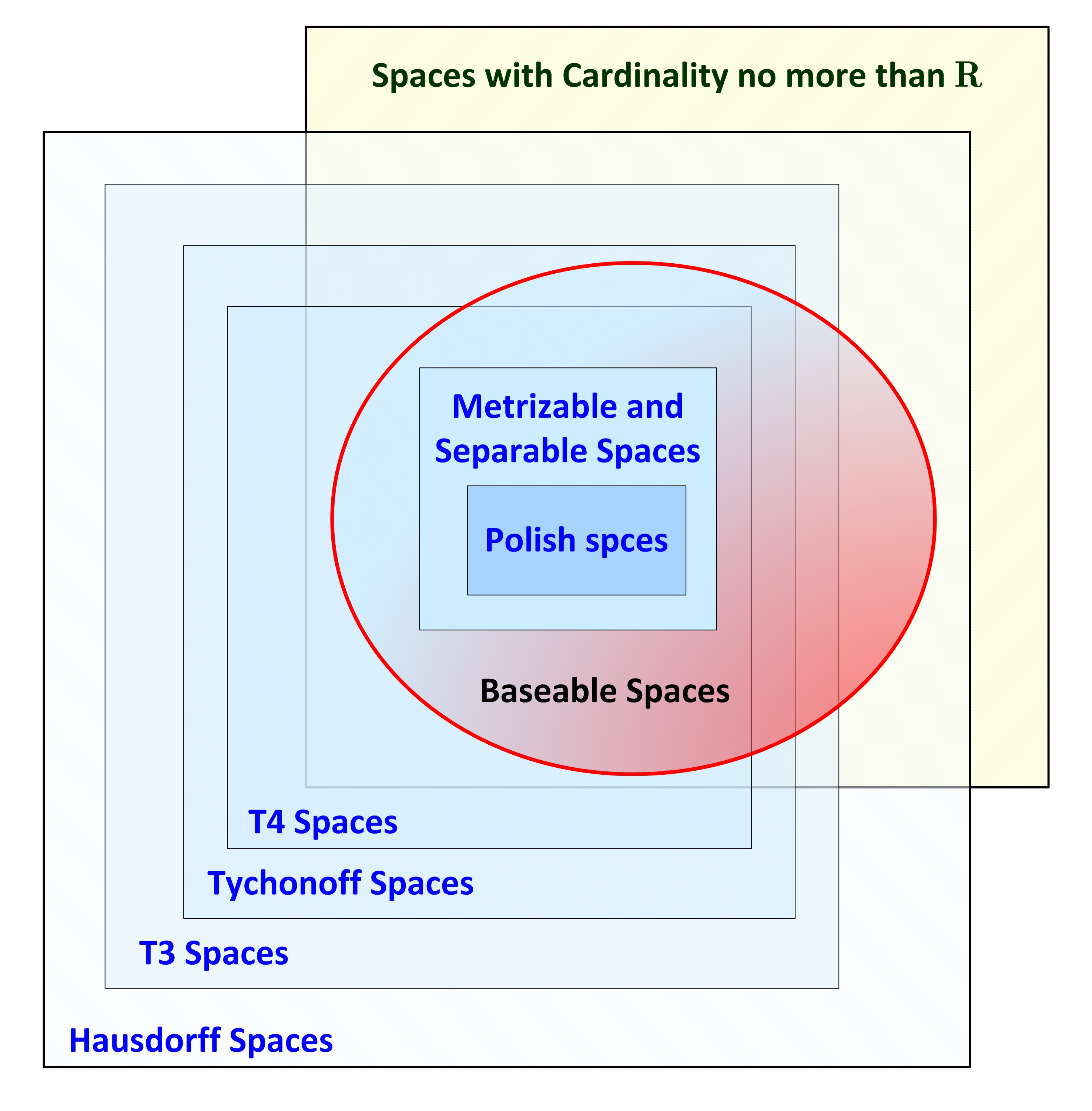}
\par\end{centering}

\caption{\textit{\label{fig:Baseable_Space}Comparison of baseable spaces and
other topological spaces}}
\end{figure}

Topological refinements of metrizable and separable spaces, i.e. so
baseable spaces range from Polish spaces to even non-\textit{T3}%
\footnote{Herein, we use the terminologies ``T3'' and ``T4'' instead of
``regular'' and ``normal'' since the latter sometimes are used
in a non-Hausdorff context.%
} (see p.\pageref{T3}) Hausdorff spaces as in the following examples.
\begin{example}
\label{exp:Baseable_Space}$\;$

\renewcommand{\labelenumi}{(\Roman{enumi}) }
\begin{enumerate}
\item \textit{Baseable, metrizable, non-Polish Lusin space}: Example \ref{exp:Pseudo-path}
mentioned that the pseudo-path topological space $D^{\mathrm{pp}}(\mathbf{R}^{+};\mathbf{R})$
is a metrizable but non-Polish Lusin space. Lusin spaces are separable
by Proposition \ref{prop:Var_Polish} (d). Hence, $D^{\mathrm{pp}}(\mathbf{R}^{+};\mathbf{R})$
is a baseable space.
\item \textit{Baseable, non-separable Banach space - 1}: Let $l^{\infty}$
be the space of all bounded $\mathbf{R}$-valued sequences equipped
with the supremum norm, i.e.
\begin{equation}
l^{\infty}\circeq\left\{ x\in\mathbf{R}^{\infty}:\Vert x\Vert_{\infty}=\sup_{n\in\mathbf{N}}\left|x_{n}\right|<\infty\right\} .\label{eq:l_inf}
\end{equation}
$(l^{\infty},\mathscr{O}_{\mathbf{R}^{\infty}}(l^{\infty}))$ is metrizable
and separable by Proposition \ref{prop:Var_Polish} (c, f). $(l^{\infty},\Vert\cdot\Vert_{\infty})$
is a Banach refinement of $(l^{\infty},\mathscr{O}_{\mathbf{R}^{\infty}}(l^{\infty}))$
by \cite[Theorem 43.5 and Theorem 20.4]{M00}, so $(l^{\infty},\Vert\cdot\Vert_{\infty})$
is a baseable space. However, $(l^{\infty},\Vert\cdot\Vert_{\infty})$
is non-separable by \cite[Example 6.6, p.83]{BN12}.
\item \textit{Baseable, non-separable Banach space - 2}: Consider the non-separable
Banach spaces mentioned in Examples \ref{exp:Holder} and \ref{exp:Rough_Path}.
 \cite[Corollary 7.50]{FV10} showed that $G^{N}(\mathbf{R}^{d})$,
the free nilpotent group of Step $N$ over $\mathbf{R}^{d}$ is a
separable Banach space with the Carnot-Caratheodory norm \textit{$\Vert\cdot\Vert_{\mathrm{cc}}$}.
$C([0,T];\mathbf{R}^{d})$ with $\Vert\cdot\Vert_{\infty}$ and $C([0,T];G^{N}(\mathbf{R}^{d}))$
equipped with the supremum $\mathrm{cc}$-norm
\begin{equation}
\Vert x\Vert_{\infty,\mathrm{cc}}\circeq\sup_{t\in[0,T]}\Vert x(t)\Vert_{\mathrm{cc}}\label{eq:CC_Sup_Norm}
\end{equation}
are Polish spaces by \cite[Theorem 2.4.3]{S98}. Then, the spaces
in (\ref{eq:p-Var_Space_Rd}) and (\ref{eq:1/p-Holder_Space_Rd})
equipped with $\Vert\cdot\Vert_{\infty}$, and those in (\ref{eq:p-Var_Space_Rough})
and (\ref{eq:1/p-Holder_Space_Rough}) equipped with $\Vert\cdot\Vert_{\infty,\mathrm{cc}}$
are all metrizable and separable spaces by Proposition \ref{prop:Var_Polish}
(c). It is known that the norms in (\ref{eq:p-Var_Norm_Rd}) and (\ref{eq:1/p-Holder_Norm_Rd})
both induce finer topologies than $\Vert\cdot\Vert_{\infty}$, while
the norms in (\ref{eq:p-Var_cc_norm}) and (\ref{eq:1/p-Holder_cc_norm})
both induce finer topologies than $\Vert\cdot\Vert_{\infty,\mathrm{cc}}$
(see \cite[p.262 and Remark 3.6]{BF13}). Thus, all the non-separable
Banach spaces in Example \ref{exp:Holder} and Example \ref{exp:Rough_Path}
are baseable spaces.
\item \textit{A baseable, non-second-countable }(see p.\pageref{Second_Countable})\textit{,
separable, Lindel$\ddot{\mbox{o}}$f }(see p.\pageref{Lindelof})\textit{
and non-metrizable T4 space}: The \textit{Sorgenfrey Line} $\mathbf{R}_{l}$
refers to the space of all real numbers equipped with the \textit{lower
limit topology}, which is generated by the topological basis
\begin{equation}
\left\{ [a,b):a,b\in\mathbf{R},a<b\right\} .\label{eq:Lower_Limit_Basis}
\end{equation}
$\mathbf{R}_{l}$ as a topological refinement of $\mathbf{R}$ is
a baseable space. $\mathbf{R}_{l}$ is separable, Lindel$\ddot{\mbox{o}}$f
but not second-countable by \cite[\S 30, Example 3]{M00}. So, $\mathbf{R}_{l}$
is non-metrizable by Proposition \ref{prop:Metrizable} (c). Nonetheless,
$\mathbf{R}_{l}$ is a \textit{T4 space} (see p.\pageref{T4}) by
\cite[\S 31, Example 2]{M00}.
\item \textit{A baseable, non-Lindel$\ddot{\mbox{o}}$f, separable and non-T4
Tychonoff space}: The \textit{Sorgenfrey Plane} $\mathbf{R}_{l}^{2}$
is a topological refinement of $\mathbf{R}^{2}$ hence baseable. Since
$\mathbf{R}_{l}$ is a separable Tychonoff space, $\mathbf{R}_{l}^{2}$
is also by Proposition \ref{prop:Countability} (c) and Proposition
\ref{prop:CR_Space} (c). However, $\mathbf{R}_{l}^{2}$ is neither
a Lindel$\ddot{\mbox{o}}$f space nor a T4 space by \cite[\S 30, Example 4 and \S 31, Example 3]{M00}.
\item \textit{A baseable, non-separable and non-metrizable Tychonoff space}:
When $E$ is a Polish space, $\mathcal{P}(E)$ is also by Theorem
\ref{thm:P(E)_Compact_Polish} (b). Example \ref{exp:Strong-topo}
explained that the strong topological space $\mathcal{P}_{S}(E)$
of all Borel probability measures on $E$ is a non-metrizable, non-separable,
Tychonoff refinement of $\mathcal{P}(E)$, so $\mathcal{P}_{S}(E)$
is a baseable space.
\item \textit{A baseable, second-countable and non-T3 space}: Let $\mathbf{R}_{K}$
denote the space of all real numbers equipped with the \textit{$K$-topology}
which is generated by the countable topological basis
\begin{equation}
\left\{ (a,b):a,b\in\mathbf{Q},a<b\right\} \cup\left\{ (a,b)\backslash\left\{ 1/n\right\} _{n\in\mathbf{N}}:a,b\in\mathbf{Q},a<b\right\} .\label{eq:K_Topo}
\end{equation}
So, $\mathbf{R}_{K}$ is a second-countable topological refinement
of $\mathbf{R}$ and hence is baseable. However, \cite[\S 31, Example 1]{M00}
explained that $\mathbf{R}_{K}$ is not a T3 space, nor is it a Tychonoff
space by Proposition \ref{prop:CR_Space} (a).
\end{enumerate}
\end{example}

A baseable space has no more points than $\mathbf{R}$ since its points
are distinguished by a countable function class.
\begin{fact}
\label{fact:Baseable_Cardinality}The cardinality of a baseable space
is no greater than $\aleph(\mathbf{R})$.
\end{fact}
\begin{proof}
The cardinality of a metrizable and separable space never exceeds
$\aleph(\mathbf{R}^{\infty})=\aleph(\mathbf{R})$ by Corollary \ref{cor:M_Compactification}
(a, c), nor can their topological refinements.\end{proof}

In general, Tychonoff spaces, metrizable spaces or separable spaces
are not necessarily ``small'' enough to be baseable spaces.
\begin{example}
\label{exp:Non_Baseable_Space}$\mathbf{R}^{[0,1]}$ equipped with
the product topology is a Tychonoff space by Proposition \ref{prop:CR_Space}
(c) and is separable by \cite[\S 30, Exercise 16 (a)]{M00}. $(\mathbf{R}^{[0,1]},\Vert\cdot\Vert_{\infty})$
is a Banach space by \cite[Theorem 43.5]{M00}. However, $\mathbf{R}^{[0,1]}$
is not baseable with any topology since its cardinality is strictly
greater than $\aleph(\mathbf{R})$.
\end{example}

It is worth mentioning that some of the baseable spaces in Example
\ref{exp:Baseable_Space} are also examples of non-Polish, non-separable
or non-metrizable standard Borel spaces.
\begin{example}
\label{exp:Baseable_Space_SB_1}Every metrizable Lusin space is a
standard Borel space by Proposition \ref{prop:Metrizable_SB_Lusin}.
Indeed, the pseudo-path topological space $D^{\mathrm{pp}}(\mathbf{R}^{+};\mathbf{R})$
is an example of a baseable, non-Polish, metrizable, separable, standard
Borel space by Example \ref{exp:Baseable_Space} (I).
\end{example}

By Proposition \ref{prop:SB} (a, d), a topological space is standard
Borel if its Borel $\sigma$-algebra can be generated by some Polish
topology. The following are examples:
\begin{example}
\label{exp:Baseable_Space_SB_2}$\,$

\renewcommand{\labelenumi}{(\Roman{enumi})   }
\begin{enumerate}
\item The Sorgenfrey Line $\mathbf{R}_{l}$ is a baseable, separable, non-metrizable,
standard Borel space. Example \ref{exp:Baseable_Space} (IV) notes
$\mathscr{B}(\mathbf{R}_{l})\supset\mathscr{B}(\mathbf{R})$. By definition,
\cite[Lemma 13.1]{M00} and Fact \ref{fact:Interval_Union}, any $O\in\mathscr{O}(\mathbf{R}_{l})$
satisfies
\begin{equation}
O=\bigcup_{i\in\mathbf{I}}[a_{i},b_{i})\in\mathscr{B}(\mathbf{R})\label{eq:Check_Rl_SB}
\end{equation}
for some $(\{a_{i}\}_{i\in\mathbf{I}}\cup\{b_{i}\}_{i\in\mathbf{I}})\subset\mathbf{R}$,
thus proving $\mathscr{B}(\mathbf{R}_{l})\subset\mathscr{B}(\mathbf{R})$.
\item For a Polish space $E$, the strong topological space $\mathcal{P}_{S}(E)$
has the same Borel $\sigma$-algebra as the Polish space $\mathcal{P}(E)$
by Lemma \ref{lem:Strong_Topo}. According to Example \ref{exp:Baseable_Space}
(VI), $\mathcal{P}_{S}(E)$ is a baseable, non-separable, non-metrizable,
standard Borel space.
\item The $K$-topological space $\mathbf{R}_{K}$ is a baseable, second-countable,
non-T3, standard Borel space. By Example \ref{exp:Baseable_Space}
(VII), $\mathscr{B}(\mathbf{R}_{K})\supset\mathscr{B}(\mathbf{R})$.
By the definition of $\mathbf{R}_{K}$ and \cite[Lemma 13.1]{M00},
any $O\in\mathscr{O}(\mathbf{R}_{K})$ satisfies
\begin{equation}
O=\left(\bigcup_{i\in\mathbf{I}}(a_{i},b_{i})\right)\backslash\{1/n\}_{n\in\mathbf{N}}\in\mathscr{B}(\mathbf{R})\label{eq:Check_K-Topo_SB}
\end{equation}
for some $(\{a_{i}\}_{i\in\mathbf{I}}\cup\{b_{i}\}_{i\in\mathbf{I}})\subset\mathbf{Q}$,
thus proving $\mathscr{B}(\mathbf{R}_{K})\subset\mathscr{B}(\mathbf{R})$.
\end{enumerate}
\end{example}
\begin{rem}
\label{rem:Rough_Path_Not_SB}The baseable but non-separable Banach
spaces in Example \ref{exp:Baseable_Space} (II, III) can not be standard
Borel, since Lemma \ref{lem:Separable_SB_Image} (b) shows metrizable
standard Borel spaces must be separable. The lack of the standard-Borel
property increases the complexity of random rough paths and their
distributions.
\end{rem}

\section{\label{sec:Baseable_Subsets}Baseable subset}

When a topological space is not necessarily baseable, its baseable
subsets are often used as ``blocks'' for building baseable subspaces.

\subsection{\label{sub:Baseable_Set_Property}Properties}

The following three facts describe baseable subsets.
\begin{fact}
\label{fact:Baseable_Subsets}Let $E$ be a topological space and
$A\subset E$. Consider the statements:

\renewcommand{\labelenumi}{(\alph{enumi})}
\begin{enumerate}
\item $A$ is a baseable subset of $E$.
\item There exits a base $(A,\mathcal{F};\widehat{E},\widehat{\mathcal{F}})$
over $E$.
\item $A$ is a $C_{b}(E;\mathbf{R})$-baseable subset of $E$.
\item $(A,\mathscr{O}_{E}(A))$ is a baseable space.
\item There exists a base $(A,\mathcal{F};\widehat{E},\widehat{\mathcal{F}})$
over $(A,\mathscr{O}_{E}(A))$.
\end{enumerate}
Then, (a) - (c) are equivalent, so are (d) and (e). Moreover, (a)
implies (d).
\end{fact}
\begin{proof}
((a) $\rightarrow$ (b)) follows by Lemma \ref{lem:Baseable_Base}
(with $E_{0}=A$). ((d) $\rightarrow$ (e)) follows by (a, b) (with
$E=(A,\mathscr{O}_{E}(A))$). The other parts are immediate by definition. \end{proof}

\begin{fact}
\label{fact:D-Baseable}Let $E$ be a topological space, $A\subset E$
and $\mathcal{D}\subset C(E;\mathbf{R})$. Consider the statements:

\renewcommand{\labelenumi}{(\alph{enumi})}
\begin{enumerate}
\item $(A,\mathscr{O}_{E}(A))$ is a $\mathcal{D}|_{A}$-baseable space.
\item $A$ is a $\mathcal{D}$-baseable subset of $E$.
\item $A$ is a $\mathcal{D}_{0}$-baseable subset of $E$ for some countable
$\mathcal{D}_{0}\subset\mathcal{D}$.
\item $A$ is a $\mathcal{D}^{\prime}$-baseable subset of $E$ for any
$\mathcal{D}^{\prime}\subset C(E;\mathbf{R})$ that includes $\mathcal{D}$.
\item Any $B\in\mathscr{B}_{E}(A)$ is a $\mathcal{D}$-baseable subset
of any topological refinement of $E$.
\end{enumerate}
Then, (b) - (e) are equivalent and any of them implies (a). Moreover,
if $A\in\mathscr{B}(E)$, then (a) - (e) are all equivalent.
\end{fact}
\begin{proof}
To show ((b) $\rightarrow$ (e)), we note that if $(E,\mathscr{U})$
is a topological refinement of $E$, then $A\in\mathscr{B}(E)\subset\sigma(\mathscr{U})$,
$\mathscr{B}_{E}(A)\subset\mathscr{B}_{(E,\mathscr{U})}(A)\subset\mathscr{B}(E,\mathscr{U})$
and $C(E;\mathbf{R})\subset C(E,\mathscr{U};\mathbf{R})$. The other
implications follow from definition.\end{proof}

\begin{fact}
\label{fact:Baseable_Space_Subset}Let $E$ be a topological space.
Then, the following statements are equivalent:

\renewcommand{\labelenumi}{(\alph{enumi})}
\begin{enumerate}
\item $E$ is a $\mathcal{D}$-baseable space (resp. baseable space).
\item Every $A\subset E$ is a $\mathcal{D}|_{A}$-baseable subspace (resp.
baseable subspace).
\item All members of $\mathscr{B}(E)$ are $\mathcal{D}$-baseable subsets
(resp. baseable subsets) of any topological refinement of $E$.
\end{enumerate}
\end{fact}
\begin{proof}
((a) $\rightarrow$ (b)) follows by (\ref{eq:C(E;R)|A}). ((b) $\rightarrow$
(c)) follows by Fact \ref{fact:D-Baseable} (a, e) (with $A=E$ and
$\mathcal{D}=\mathcal{D}$ or $C(E;\mathbf{R})$). ((c) $\rightarrow$
(a)) is automatic.\end{proof}

Next, we describe countable unions and products of baseable subsets.
\begin{fact}
\label{fact:D-Baseable_Union}Let $E$ be a topological space and
$A_{n}$ be a $\mathcal{D}_{n}$-baseable subset of $E$ for each
$n\in\mathbf{N}$. If $\{A_{n}\}_{n\in\mathbf{N}}$ is nested (i.e.
any $A_{n_{1}}$ and $A_{n_{2}}$ admit a common superset $A_{n_{3}}$),
then $\bigcup_{n\in\mathbf{N}}A_{n}$ is a $\bigcup_{n\in\mathbf{N}}\mathcal{D}_{n}$-baseable
subset of $E$.
\end{fact}
\begin{proof}
This result is immediate by Fact \ref{fact:SP_Nested_Union}.\end{proof}

\begin{prop}
\label{prop:D-Baseable_Prod}Let $A_{n}$ be a $\mathcal{D}_{n}$-baseable
subset of topological space $S_{n}$ for each $n\in\mathbf{N}$. Then:

\renewcommand{\labelenumi}{(\alph{enumi})}
\begin{enumerate}
\item $\prod_{n=1}^{m}A_{n}$ is a $\mathcal{D}^{m}$-baseable subset of
$\prod_{n=1}^{m}S_{n}$ with $\mathcal{D}^{m}\circeq\{f\circ\mathfrak{p}_{n}:1\leq n\leq m,f\in\mathcal{D}_{n}\}$
for all $m\in\mathbf{N}$.
\item $\prod_{n\in\mathbf{N}}A_{n}$ is a $\bigcup_{m\in\mathbf{N}}\mathcal{D}^{m}$-baseable
subset of $\prod_{n\in\mathbf{N}}S_{n}$.
\end{enumerate}
\end{prop}
\begin{proof}
(a) $\prod_{n=1}^{m}A_{n}=\bigcap_{n=1}^{m}\mathfrak{p}_{n}^{-1}(A_{n})\in\mathscr{B}(E)^{\otimes m}\subset\mathscr{B}(E^{m})$
by Proposition \ref{prop:Prod_Space} (a). Let $\{f_{n,k}\}_{k\in\mathbf{N}}\subset\mathcal{D}_{n}$
separate points on $A_{n}$ for each $n\in\mathbf{N}$ and $\mathcal{D}_{0}^{m}\circeq\{f_{n,k}\circ\mathfrak{p}_{n}\}_{1\leq n\leq m,k\in\mathbf{N}}$.
Then, $\bigotimes\mathcal{D}_{0}^{m}(x)=\bigotimes\mathcal{D}_{0}^{m}(y)$
in $\mathbf{R}^{\mathcal{D}_{0}^{m}}$ implies $x=\bigotimes_{n=1}^{m}\mathfrak{p}_{n}(x)=\bigotimes_{n=1}^{m}\mathfrak{p}_{n}(y)=y$
in $\prod_{n=1}^{m}S_{n}$. Thus, $\mathcal{D}_{0}^{m}$ is a countable
subset of $\mathcal{D}^{m}$ that separates points on $\prod_{n=1}^{m}A_{n}$.

(b) follows by an argument similar to (a).\end{proof}

\begin{cor}
\label{cor:D-Baseable_Ed}Let $E$ be a topological space, $\{A_{n}\}_{n\in\mathbf{N}}$
be $\mathcal{D}$-baseable subsets of $E$ and $d\in\mathbf{N}$.
Then, $\prod_{n=1}^{d}A_{n}$ is a $\Pi^{d}(\mathcal{D})$-baseable
subset of $E^{d}$.
\end{cor}
\begin{proof}
This corollary follows by Proposition \ref{prop:D-Baseable_Prod}
(a) (with $m=d$ and $\mathcal{D}_{n}=\mathcal{D}$), the definition
of $\Pi^{d}(\mathcal{D})$ and Fact \ref{fact:D-Baseable} (b, d).\end{proof}

\subsection{\label{sub:Test_Fun}Selection of point-separating functions}

When using $\mathcal{D}$-baseable subsets to construct a base, one
can include \textit{desired, up to countably many bounded} members
of $\mathcal{D}$ within the base. This is useful in applications
because one may include a desirable set of functions such as subdomains
of operators, observation functions of nonlinear filters and test
functions for measure-valued processes, etc.
\begin{lem}
\label{lem:Base_Construction}Let $E$ be a topological space, $E_{0}$
be a $\mathcal{D}$-baseable subset of $E$ and $\mathcal{D}_{0}\subset C_{b}(E;\mathbf{R})$
be countable. Then:

\renewcommand{\labelenumi}{(\alph{enumi})}
\begin{enumerate}
\item There exists a base $(E_{0},\mathcal{F};\widehat{E},\widehat{\mathcal{F}})$
over $E$ with $\mathcal{D}_{0}\subset\mathcal{F}$.
\item If $\mathcal{D}$ is countable, then the $\mathcal{F}$ in (a) can
be taken to contain $\mathcal{D}\cap C_{b}(E;\mathbf{R})$. If, in
addition, $\mathcal{D}\subset C_{b}(E;\mathbf{R})$, then $\mathcal{F}$
can be taken to equal $\mathcal{D}\cup\{1\}$.
\item If $\mathcal{D}_{0}\subset\mathcal{D}\subset C_{b}(E;\mathbf{R})$,
then the $\mathcal{F}$ in (a) can be taken within $\mathcal{D}\cup\{1\}$.
If, in addition, $\mathcal{D}$ is closed under addition or multiplication,
then $\mathcal{F}$ can be taken to have the same closedness.
\end{enumerate}
\end{lem}
\begin{proof}
Let $\mathcal{D}^{\prime}\subset\mathcal{D}$ be countable and separate
points on $E_{0}$. Then, (a) and (b) follow by Lemma \ref{lem:Baseable_Base}
(with $\mathcal{D}=\mathcal{D}^{\prime}\cup\mathcal{D}_{0}$).

For (c), we define%
\footnote{The notation ``$\mathfrak{ac}(\cdot)$'' was defined in \S \ref{sub:Fun}.%
} 
\begin{equation}
\mathcal{D}^{\prime\prime}\circeq\begin{cases}
\mathcal{D}_{0}\cup\mathcal{D}^{\prime}\cup\{1\}, & \mbox{in general},\\
\mathfrak{ac}\left(\mathcal{D}_{0}\cup\mathcal{D}^{\prime}\cup\{1\}\right), & \mbox{if }\mathcal{D}\ni1\mbox{ is closed under addition},\\
\mathfrak{mc}\left(\mathcal{D}_{0}\cup\mathcal{D}^{\prime}\cup\{1\}\right), & \mbox{if }\mathcal{D}\ni1\mbox{ is closed under multiplication},\\
\mathfrak{ac}\left(\mathfrak{mc}\left(\mathcal{D}_{0}\cup\mathcal{D}^{\prime}\cup\{1\}\right)\right), & \mbox{if }\mathcal{D}\ni1\mbox{ is closed under both}.
\end{cases}\label{eq:TF}
\end{equation}
In any case above, $\mathcal{D}^{\prime\prime}$ is a countable subset
of $\mathcal{D}\cup\{1\}$, contains $\{1\}\cup\mathcal{D}_{0}$ and
separates points on $E_{0}$. Now, (c) follows by (b) (with $\mathcal{D}=\mathcal{D}^{\prime\prime}$).\end{proof}

Often $\mathcal{D}\subset C(E;\mathbf{R})$ is uncountable but known
to separate points on $E$, and one desires a subset $A\subset E$
is $\mathcal{D}$-baseable. In other words, one desires reducing this
specific $\mathcal{D}$ to a countable subcollection separating points
on $A$. One sufficient condition for such reduction is the \textit{hereditary
Lindel$\ddot{\mbox{o}}$f property} (see p.\pageref{Hered_Lindelof})
of $A$.
\begin{prop}
\label{prop:Hered_Lindelof_Baseable}Let $E$ be a topological space,
$A\in\mathscr{B}(E)$ and $\mathcal{D}\subset C(E;\mathbf{R})$ separate
points on $A$. If $\{(x,x):x\in A\}$ is a Lindel$\ddot{\mbox{o}}$f
subspace of $E\times E$, especially if $A$ is a Souslin or second-countable
subspace of $E$, then $A$ is a $\mathcal{D}$-baseable subset of
$E$.
\end{prop}
\begin{proof}
If $A$ is a Souslin or second-countable subspace of $E$, then $A\times A$
is a hereditary Lindel$\ddot{\mbox{o}}$f subspace of $E\times E$
by Proposition \ref{prop:Var_Polish} (d, f) and Proposition \ref{prop:Countability}
(b, c), which implies $\{(x,x):x\in A\}$ is a Lindel$\ddot{\mbox{o}}$f
subspace of $E\times E$. Now, the result follows by Proposition \ref{prop:Fun_Sep_2}
(a) (with $E=(A,\mathscr{O}_{E}(A))$ and $\mathcal{D}=\mathcal{D}|_{A}$).\end{proof}

\begin{rem}
\label{rem:Base_Fun_Reduction}Separating points is usually weaker
than strongly separating points. Compared to Proposition \ref{prop:Fun_Sep_2}
(b), the selection result above uses hereditary Lindel$\ddot{\mbox{o}}$f
property, which is weaker than second-countability.
\end{rem}

\begin{rem}
\label{rem:Souslin_Hered_SP}When $E$ is a Tychonoff space, $C(E;\mathbf{R})$
separates points on $E$ by Proposition \ref{prop:CR} (a, b). So,
Proposition \ref{prop:Hered_Lindelof_Baseable} (with $\mathcal{D}=C(E;\mathbf{R})$)
slightly generalizes \cite[Vol. II, Theorem 6.7.7 (ii)]{B07}.
\end{rem}

Point-separating functions can be selected from a uniformly dense
collection.
\begin{prop}
\label{prop:Dense_Baseable}Let $E$ be a topological space and $A$
be a $\mathcal{D}$-baseable subset. If $\mathcal{D}_{0}\subset C_{b}(E;\mathbf{R})$
satisfies $\mathcal{D}\subset\mathfrak{cl}(\mathcal{D}_{0})$, then
$E$ is a $\mathcal{D}_{0}$-baseable subset.
\end{prop}
\begin{proof}
Let $\{f_{n}\}_{n\in\mathbf{N}}\subset\mathcal{D}$ separate points
on $A$. $\mathcal{D}$ and its superset $\mathfrak{cl}(\mathcal{D}_{0})$
both lie in the Banach space $(M_{b}(E;\mathbf{R}),\Vert\cdot\Vert_{\infty})$,
so there exist $\{f_{n,k}\}_{n,k\in\mathbf{N}}\subset\mathcal{D}_{0}$
such that $f_{n,k}\overset{u}{\rightarrow}f_{n}$ as $k\uparrow\infty$
for each $n\in\mathbf{N}$ by Fact \ref{fact:First_Countable}. Hence,
$\{f_{n}\}_{n\in\mathbf{N}}\subset\mathfrak{cl}(\{f_{n,k}\}_{n,k\in\mathbf{N}})$
and $\{f_{n,k}\}_{n,k\in\mathbf{N}}$ separates points on $A$ by
Corollary \ref{cor:SSP_Dense} (with $\mathcal{D}=\{f_{n}\}_{n\in\mathbf{N}}$
and $\mathcal{D}_{0}=\{f_{n,k}\}_{n,k\in\mathbf{N}}$).\end{proof}

\subsection{\label{sub:Baseable_SB}Baseable standard Borel subsets}

The standard Borel subsets and Borel subsets coincide for a baseable
standard Borel subspace. The following result generalizes its classical
version on metrizable spaces (see Proposition \ref{prop:SB_Borel}
(b)).
\begin{prop}
\label{prop:Baseable_SB_Borel}Let $E$ be a topological space and
$A\in\mathscr{B}^{\mathbf{s}}(E)$. Then:

\renewcommand{\labelenumi}{(\alph{enumi})}
\begin{enumerate}
\item If $(A,\mathscr{O}_{E}(A))$ is a baseable space, then $\mathscr{B}^{\mathbf{s}}(A,\mathscr{O}_{E}(A))=\mathscr{B}_{E}(A)\subset\mathscr{B}^{\mathbf{s}}(E)$.
\item If $E$ is a baseable space, then $\mathscr{B}^{\mathbf{s}}(E)\subset\mathscr{B}(E)$.
\item If $E$ is a baseable standard Borel space, then $\mathscr{B}(E)=\mathscr{B}^{\mathbf{s}}(E)$.
\end{enumerate}
\end{prop}
\begin{proof}
(a) $\mathscr{B}_{E}(A)\subset\mathscr{B}^{\mathbf{s}}(A,\mathscr{O}_{E}(A))\subset\mathscr{B}^{\mathbf{s}}(E)$
by Proposition \ref{prop:SB_Borel} (a). Now, let $B\in\mathscr{B}^{\mathbf{s}}(A,\mathscr{O}_{E}(A))$.
There exists a base $(A,\mathcal{F};\widehat{E},\widehat{\mathcal{F}})$
over $(A,\mathscr{O}_{E}(A))$ by Fact \ref{fact:Baseable_Subsets}
(d, e). Then, $B\in\mathscr{B}_{E}(A)$ by Lemma \ref{lem:SB_Base}
(b) (with $d=1$, $E=E_{0}=(A,\mathscr{O}_{E}(A))$ and $A=B$).

(b) There exists a base $(E,\mathcal{F};\widehat{E},\widehat{\mathcal{F}})$
over $E$ by Fact \ref{fact:Baseable_Subsets} (d, e) (with $A=E$).
Then, (b) follows by Lemma \ref{lem:SB_Base} (b) (with $d=1$ and
$E_{0}=E$).

(c) follows immediately by (a) (with $A=E$).\end{proof}

For $\mathcal{D}$-baseable standard Borel subsets, the function class
$\mathcal{D}$ not only separates their points but also determines
their subspace Borel $\sigma$-algebras.
\begin{prop}
\label{prop:Baseable_SB_B(D0)=00003DB(E)}Let $E$ be a topological
space and $A\in\mathscr{B}^{\mathbf{s}}(E)$. Then, $A$ is a $\mathcal{D}$-baseable
subset of $E$ if and only if $A\in\mathscr{B}(E)$ and there exists
a countable $\mathcal{D}_{0}\subset\mathcal{D}$ such that $\mathscr{B}_{E}(A)=\mathscr{B}_{\mathcal{D}_{0}}(A)=\sigma(\mathcal{D}_{0})|_{A}$%
\footnote{The $\sigma$-algebra $\sigma(\mathcal{D}_{0})$ was defined in \S
\ref{sub:Meas}. Recall that $\sigma(\mathcal{D}_{0})|_{A}$ is generally
smaller than $\mathscr{B}_{\mathcal{D}_{0}}(A)$.%
}.
\end{prop}
\begin{proof}
The $\mathcal{D}$-baseability of $A$ implies $A\in\mathscr{B}(E)$
and a countable $\mathcal{D}_{0}\subset\mathcal{D}\subset C(E;\mathbf{R})$
that separates points on $A$. Then, the result follows by Proposition
\ref{prop:SB_SP} (b, c) (with $E=(A,\mathscr{O}_{E}(A))$ and $\mathcal{D}=\mathcal{D}_{0}|_{A}$).\end{proof}

Baseability facilitates transferability of the standard Borel property.
\begin{prop}
\label{prop:Baseable_SB_Union}Let $E$ be a topological space, $\mathbf{I}$
be a countable set, $\{A_{i}\}_{i\in\mathbf{I}}\subset\mathscr{B}^{\mathbf{s}}(E)$
and $A=\bigcup_{i\in\mathbf{I}}A_{i}$. If $(A,\mathscr{O}_{E}(A))$
is a baseable space, then $A\in\mathscr{B}^{\mathbf{s}}(E)$.
\end{prop}
\begin{proof}
There exists a base $(A,\mathcal{F};\widehat{E},\widehat{\mathcal{F}})$
over $(A,\mathscr{O}_{E}(A))$ by Fact \ref{fact:Baseable_Subsets}
(d, e). Then, $A\in\mathscr{B}^{\mathbf{s}}(E)$ by Lemma \ref{lem:SB_Base}
(c) (with $E=E_{0}=(A,\mathscr{O}_{E}(A))$ and $d=1$).\end{proof}

Baseable standard Borel support often implies unique Borel extension.
\begin{prop}
\label{prop:Baseable_SB_Borel_Ext}Let $E$ be a topological space,
$A$ be a baseable subset of $E$, $d\in\mathbf{N}$ and $\mu\in\mathfrak{M}^{+}(E^{d},\mathscr{B}(E)^{\otimes d})$.
If $\mu$ is supported on $B\in\mathscr{B}^{\mathbf{s}}(E^{d})$ and
$B\subset A^{d}$, then $\mathfrak{be}(\mu)$ is a singleton.
\end{prop}
\begin{proof}
The baseability of $A$ implies $A\in\mathscr{B}(E)$ and $A^{d}\in\mathscr{B}(E)^{\otimes d}$.
One finds by Fact \ref{fact:Meas_Concen_Expan} (a) (with $\mathscr{U}=\mathscr{B}(E)^{\otimes d}$
and $A=A^{d}$) that $\mu|_{A^{d}}\in\mathfrak{M}^{+}(A^{d},\mathscr{B}(E)^{\otimes d}|_{A^{d}})$.
There exists a base $(A,\mathcal{F};\widehat{E},\widehat{\mathcal{F}})$
over $(A,\mathscr{O}_{E}(A))$ by Fact \ref{fact:Baseable_Subsets}
(d, e). $\mathfrak{be}(\mu|_{A^{d}})$ is a singleton by Corollary
\ref{cor:SB_Base_Borel_Extension} (a) (with $E=E_{0}=(A,\mathscr{O}_{E}(A))$,
$A=B$ and $\mu=\mu|_{A^{d}}$). Now, $\mathfrak{be}(\mu)$ is a singleton
by the fact $A^{d}\in\mathscr{B}(E)^{\otimes d}$ and Lemma \ref{lem:Union_Borel_Prod_Equal}
(b) (with $\mathbf{I}=\{1,...,d\}$, $S_{i}=E$, $S=E^{d}$, $\mathscr{A}=\mathscr{B}(E)^{\otimes d}$
and $A=A^{d}$).\end{proof}

The following three results relate the baseability of a standard Borel
space $E$ and that of $\mathcal{P}(E)$.
\begin{prop}
\label{prop:Baseable_SB_Meas_Sep}Let $E$ be a topological space,
$\mathcal{D}\subset C_{b}(E;\mathbf{R})$ and $A\in\mathscr{B}^{\mathbf{s}}(E)$.
If $(A,\mathscr{O}_{E}(A))$ is a $\mathcal{D}|_{A}$-baseable space,
then $\mathcal{P}(A,\mathscr{O}_{E}(A))$ is an $\mathfrak{mc}(\mathcal{D}|_{A})^{*}$-baseable
space.
\end{prop}
\begin{proof}
There exists a base $(A,\mathcal{F};\widehat{E},\widehat{\mathcal{F}})$
over $(A,\mathscr{O}_{E}(A))$ with $(\mathcal{F}\backslash\{1\})\subset\mathcal{D}|_{A}$
by Lemma \ref{lem:Base_Construction} (c) (with $E=E_{0}=(A,\mathscr{O}_{E}(A))$
and $\mathcal{D}=\mathcal{D}|_{A}$). Then, $\mathfrak{mc}(\mathcal{F}\backslash\{1\})^{*}$
is a countable subset of $\mathfrak{mc}(\mathcal{D}|_{A})^{*}$ by
Fact \ref{fact:ac_mc_Countable} and separates points on $\mathcal{P}(A,\mathscr{O}_{E}(A))$
by Corollary \ref{cor:SB_Base_Sep_Meas} (a) (with $E=(A,\mathscr{O}_{E}(A))$
and $d=1$).\end{proof}

\begin{cor}
\label{cor:E_Baseable_P(E)_Baseable}Let $E$ be a baseable standard
Borel space. Then, $\mathcal{M}^{+}(E)$ and $\mathcal{P}(E)$ are
baseable spaces.
\end{cor}
\begin{proof}
The result follows by Fact \ref{fact:Baseable_Subsets} (a, c) (with
$A=E$), Proposition \ref{prop:Baseable_SB_Meas_Sep} (with $\mathcal{D}=C_{b}(E;\mathbf{R})$
and $A=E$) and Fact \ref{fact:Sep_CD} (a) (with $\mathcal{D}=C_{b}(E;\mathbf{R})$).\end{proof}

\begin{prop}
\label{prop:P(E)_Baseable_E_Baseable}Let $E$ be a first-countable
space and $\{\{x\}:x\in E\}\subset\mathscr{B}(E)$%
\footnote{The property that singletones are Borel sets is milder than the Hausdorff
property or the \textit{T1 axiom} (see \cite[\S 17, p.99]{M00}).%
}. If $\mathcal{P}(E)$ is a baseable space, then $E$ is also.
\end{prop}
\begin{proof}
We suppose $\{g_{n}\}_{n\in\mathbf{N}}\subset C(\mathcal{P}(E);\mathbf{R})$
separates points on $\mathcal{P}(E)$ and define $f_{n}(x)\circeq g_{n}(\delta_{x})$%
\footnote{$\delta_{x}$ denotes the Dirac measure at $x$.%
} for all $x\in E$ and $n\in\mathbf{N}$. For distinct $x,y\in E$,
$\delta_{x}\neq\delta_{y}$ by Proposition \ref{prop:P(E)_Hausdorff}
(a) and so $\bigotimes_{n\in\mathbf{N}}f_{n}(x)=\bigotimes_{n\in\mathbf{N}}g_{n}(\delta_{x})\neq\bigotimes_{n\in\mathbf{N}}g_{n}(\delta_{y})=\bigotimes_{n\in\mathbf{N}}f_{n}(y)$.
Hence, $\{f_{n}\}_{n\in\mathbf{N}}\subset\mathbf{R}^{E}$ separates
points on $E$. We show each $f_{n}\in C(E;\mathbf{R})$. If $x_{k}\rightarrow x$
as $k\uparrow\infty$ in $E$, then $\delta_{x_{k}}\Rightarrow\delta_{x}$
as $k\uparrow\infty$ in $\mathcal{P}(E)$ by Fact \ref{fact:Dirac_WC}.
It follows by the continuity of $g_{n}$ that $\lim_{k\rightarrow\infty}f_{n}(x_{k})=\lim_{k\rightarrow\infty}g_{n}(\delta_{x_{k}})=g_{n}(\delta_{x})=f_{n}(x)$.
Now, the continuity of $f_{n}$ follows by the first-countability
of $E$ and \cite[Theorem 30.1 (b)]{M00}.\end{proof}

\subsection{\label{sub:Metrizable_Compact}Metrizable compact subsets}

Metrizable compact subsets form an essential class of hereditary Lindel$\ddot{\mbox{o}}$f,
standard Borel, baseable subsets for replication and weak convergence.
The following proposition gives several equivalent forms of metrizable
compact subsets.
\begin{prop}
\label{prop:MC}Let $E$ be a topological space, $K\in\mathscr{K}(E)$
and $\mathcal{D}\subset C(E;\mathbf{R})$. Consider the following
statements:

\renewcommand{\labelenumi}{(\alph{enumi})}
\begin{enumerate}
\item $K$ is a $\mathcal{D}$-baseable subset of $E$.
\item $K$ is a Souslin subspace of $E$.
\item $K$ is a Hausdorff subspace of $E$ and $\{(x,x):x\in K\}$ is a
Lindel$\ddot{\mbox{o}}$f subspace of $E\times E$.
\item $(K,\mathscr{O}_{E}(K))$ is a baseable space.
\item $K$ is a Hausdorff and second-countable subspace of $E$.
\item $K\in\mathscr{K}^{\mathbf{m}}(E)$.
\item $K\in\mathscr{B}^{\mathbf{s}}(E)$.
\end{enumerate}
Then, (b) - (f) are equivalent and implied by (a). (f) implies (g).
Morever, if $\mathcal{D}$ separates points on $E$, then (a) - (f)
are all equivalent.
\end{prop}
\begin{proof}
((b) $\rightarrow$ (c)) $K$ is a Hausdorff subspace by the definition
of Souslin spaces. $K\times K$ is a Souslin subspace of $E\times E$
by Proposition \ref{prop:Var_Polish} (f). Now, (c) follows by Proposition
\ref{prop:Var_Polish} (d).

((c) $\rightarrow$ (d)) $(K,\mathscr{O}_{E}(K))$ is Tychonoff by
Proposition \ref{prop:Compact} (d) and Proposition \ref{prop:CR_Space}
(a). Then, (d) follows by Proposition \ref{prop:CR} (a, b) (with
$E=(K,\mathscr{O}_{E}(K))$) and Proposition \ref{prop:Hered_Lindelof_Baseable}
(with $E=A=(K,\mathscr{O}_{E}(K))$ and $\mathcal{D}=C(K,\mathscr{O}_{E}(K);\mathbf{R})$).

((d) $\rightarrow$ (e, f)) $(K,\mathscr{O}_{E}(K))$ is Hausdorff
by Fact \ref{fact:Baseable_Metrizable_Separable} (a). Let $\mathcal{D}\subset C_{b}(K,\mathscr{O}_{E}(K);\mathbf{R})$
be countable and separate points on $(K,\mathscr{O}_{E}(K))$. $\mathcal{D}$
strongly separates points on $(K,\mathscr{O}_{E}(K))$ by Lemma \ref{lem:SP_on_Compact}.
Both (e) and (f) follow by Proposition \ref{prop:Fun_Sep_1} (d).

((e) $\rightarrow$ (c)) $K\times K$ is a second-countable subspace
of $E\times E$ by Proposition \ref{prop:Countability} (c). Then,
(c) follows by Proposition \ref{prop:Countability} (b).

((f) $\rightarrow$ (b, g)) follows by Proposition \ref{prop:Compact}
(d), Proposition \ref{prop:Var_Polish} (a) and Fact \ref{fact:Polish_SB}
(a).

((a) $\rightarrow$ (d)) follows by Fact \ref{fact:D-Baseable} (b,
d) (with $\mathcal{D}^{\prime}=C(E;\mathbf{R})$) and Fact \ref{fact:Baseable_Subsets}
(a, d) (with $A=K$).

Moreover, if $\mathcal{D}$ separates points on $E$, then $K\in\mathscr{B}(E)$
by Proposition \ref{prop:Fun_Sep_1} (e) (with $A=E$) and Proposition
\ref{prop:Compact} (a), and (c) implies (a) by Proposition \ref{prop:Hered_Lindelof_Baseable}
(with $A=K$).\end{proof}

Baseable spaces, Lusin spaces and Souslin spaces are not necessarily
metrizable, but all of them have metrizable compact subsets.
\begin{cor}
\label{cor:Baseable_MC}Let $E$ be a topological space and $K\in\mathscr{K}(E)$.
Then:

\renewcommand{\labelenumi}{(\alph{enumi})}
\begin{enumerate}
\item If $E$ is a baseable space, then $K$ is a metrizable standard Borel
subspace and is a baseable subset of $E$.
\item If $E$ is a Souslin or Lusin space, then $K$ is a metrizable, baseable,
standard Borel subspace of $E$. If, in addition, $C(E;\mathbf{R})$
separates points on $E$, then $K$ is a baseable subset of $E$.
\end{enumerate}
\end{cor}
\begin{proof}
(a) follows by Fact \ref{fact:Baseable_Space_Subset} and Proposition
\ref{prop:MC} (a, f, g) (with $\mathcal{D}=C(E;\mathbf{R})$). Next,
Proposition \ref{prop:Compact} (a), Proposition \ref{prop:Var_Polish}
(a, b) and Proposition \ref{prop:Separability} (c) imply that Lusin
(resp. Souslin) spaces are Souslin (resp. Hausdorff) spaces and compact
subsets of a Souslin space are closed, Souslin, Hausdorff subspaces.
Then, (b) follows by Proposition \ref{prop:MC} (a, b, d, f, g) (with
$\mathcal{D}=C(E;\mathbf{R})$).\end{proof}

The two results above indicated that many non-metrizable topological
spaces like those in Example \ref{exp:Baseable_Space} have metrizable
compact subsets. Still, having metrizable compact subsets is a strictly
milder property than baseability.
\begin{example}
\label{exp:Metrizable_Compact_non-Baseable}$([0,1]^{[0,1]},\Vert\cdot\Vert_{\infty})$
is a Banach space so it certainly has metrizable compact subsets.
However, it is not baseable as in Example \ref{exp:Non_Baseable_Space}.
\end{example}
Here are several constructive properties of metrizable compact subsets.
\begin{lem}
\label{lem:MC_Union}Let $E$ be a topological space, $m\in\mathbf{N}$
and $\{A_{i}\}_{1\leq i\leq m}\subset\mathscr{K}^{\mathbf{m}}(E)$.
If $A=\bigcup_{i=1}^{m}A_{i}$ is a Hausdorff subspace of $E$, then
$A\in\mathscr{K}^{\mathbf{m}}(E)$.
\end{lem}
\begin{proof}
$A\in\mathscr{K}(E)$ by Proposition \ref{prop:Compact} (b). $\{A_{i}\}_{1\leq i\leq m}$
are Souslin subspaces of $E$ by Proposition \ref{prop:MC} (b, f).
$A$ is a Souslin subspace of $E$ by Proposition \ref{prop:Var_Polish}
(g). Now, the result follows by Proposition \ref{prop:MC} (b, f).\end{proof}

\begin{lem}
\label{lem:MC_Prod}Let $\mathbf{I}$ be a countable index set, $\{S_{i}\}_{i\in\mathbf{I}}$
be topological spaces and $S\circeq\prod_{i\in\mathbf{I}}S_{i}$.
Then:

\renewcommand{\labelenumi}{(\alph{enumi})}
\begin{enumerate}
\item If $A_{i}\in\mathscr{K}^{\mathbf{m}}(S_{i})$ for all $i\in\mathbf{I}$,
then $\prod_{i\in\mathbf{I}}A_{i}\in\mathscr{K}^{\mathbf{m}}(S)$.
\item If $A\in\mathscr{K}^{\mathbf{m}}(S)$ and $\mathfrak{p}_{i}(A)$ is
a Hausdorff subspace of $S_{i}$ for some $i\in\mathbf{I}$, then
$\mathfrak{p}_{i}(A)\in\mathscr{K}^{\mathbf{m}}(S_{i})$.
\end{enumerate}
\end{lem}
\begin{proof}
(a) follows by Proposition \ref{prop:Compact} (b) and Proposition
\ref{prop:Metrizable_Prod}.

(b) $\{(x,x):x\in A\}$ is a Lindel$\ddot{\mbox{o}}$f subspace of
$S\times S$ by Proposition \ref{prop:MC} (c, f) (with $E=S$). It
follows that $\mathfrak{p}_{i}(A)\in\mathscr{K}(S_{i})$ and
\begin{equation}
\left\{ (y,y):y\in\mathfrak{p}_{i}(A)\right\} =\left\{ \left(\mathfrak{p}_{i}(x),\mathfrak{p}_{i}(x)\right):x\in A\right\} \label{eq:Check_MC_Proj}
\end{equation}
is a Lindel$\ddot{\mbox{o}}$f subspace of $S_{i}\times S_{i}$ by
Fact \ref{fact:Prod_Map_2} (a), Proposition \ref{prop:Compact} (e)
and Proposition \ref{prop:Countability} (d). Thus, $\mathfrak{p}_{i}(A)\in\mathscr{K}^{\mathbf{m}}(S_{i})$
by its Hausdorff proeprty and Proposition \ref{prop:MC} (c, f).\end{proof}

To handle $\mathbf{m}$-tightness of non-Borel measures, Definition
\ref{def:Tight} required a collection of metrizable compact sets
lying in their domains. The next lemma shows that if $E$ is a Hausdorff
space, then this requirement is automatically satisfied by the members
of $\mathfrak{M}^{+}(E^{d},\mathscr{B}(E)^{\otimes d})$%
\footnote{We remind the readers that: (i) the Borel $\sigma$-algebra of the
product topological space $E^{d}$ can be different than the product
$\sigma$-algebra $\mathscr{B}(E)^{\otimes d}$, and (ii) the notation
$\mathfrak{M}^{+}(\cdot)$ means the family of finite measures.%
}.
\begin{lem}
\label{lem:MC_Prod_Meas}Let $\mathbf{I}$ be a countable index set,
$\{S_{i}\}_{i\in\mathbf{I}}$ be topological spaces, $(S,\mathscr{A})$
be as in (\ref{eq:(S,A)_Prod_Meas_Space}) and $A\in\mathscr{K}^{\mathbf{m}}(S)$.
If $B_{i}\in\mathscr{B}(S_{i})$ is a Hausdorff subspace of $S_{i}$
and contains $\mathfrak{p}_{i}(A)$ for all $i\in\mathbf{I}$, then
$A\in\mathscr{A}$ and $\mathscr{B}_{S}(A)=\mathscr{A}|_{A}$.
\end{lem}
\begin{proof}
$\mathfrak{p}_{i}(A)\in\mathscr{K}^{\mathbf{m}}(S_{i})$ for all $i\in\mathbf{I}$
by Proposition \ref{prop:Separability} (c) and Lemma \ref{lem:MC_Prod}
(b). As $B_{i}\in\mathscr{B}(S_{i})$ for all $i\in\mathbf{I}$, we
have that
\begin{equation}
A\subset F\circeq\prod_{i\in\mathbf{I}}\mathfrak{p}_{i}(A)\in\bigotimes_{i\in\mathbf{I}}\mathscr{B}_{S_{i}}(B_{i})\subset\mathscr{A}\label{eq:Check_MC_Prod_Meas_1}
\end{equation}
by Corollary \ref{cor:Compact_Prod} (a) (with $S_{i}=B_{i}$). $F\in\mathscr{K}^{\mathbf{m}}(S)$
by Lemma \ref{lem:MC_Prod} (a). $F$ is a second-countable subspace
of $S$ by Proposition \ref{prop:MC} (e, f).
\begin{equation}
\mathscr{B}_{S}(F)=\mathscr{A}|_{F}\subset\mathscr{A}\label{eq:Check_MC_Prod_Meas_2}
\end{equation}
by Proposition \ref{prop:Prod_Space} (c) (with $S_{i}=\mathfrak{p}_{i}(A)$)
and (\ref{eq:Check_MC_Prod_Meas_1}). This implies $\mathscr{B}_{S}(A)=\mathscr{A}|_{A}$
since $A\subset F$. Moreover, $F$ is a Hausdorff subspace of $S$
by Proposition \ref{prop:Separability} (d). Hence, $A\in\mathscr{B}_{S}(F)\subset\mathscr{A}$
by Proposition \ref{prop:Compact} (a) and (\ref{eq:Check_MC_Prod_Meas_2}).\end{proof}

$\mathbf{m}$-tightness ensures a unique and tight Borel extension
on product space.
\begin{prop}
\label{prop:m-Tight_BExt}Let $\mathbf{I}$ be a countable index set,
$\{S_{i}\}_{i\in\mathbf{I}}$ be topological spaces, $(S,\mathscr{A})$
be as in (\ref{eq:(S,A)_Prod_Meas_Space}), $\Gamma\subset\mathfrak{M}^{+}(S,\mathscr{A})$
and $A\subset S$. Suppose in addition that $\mathfrak{p}_{i}(A)\in\mathscr{B}(S_{i})$
is a Hausdorff subspace of $S_{i}$ for all $i\in\mathbf{I}$. Then,
$\Gamma$ is $\mathbf{m}$-tight in $A$ if and only if $\{\mu^{\prime}=\mathfrak{be}(\mu)\}_{\mu\in\Gamma}$%
\footnote{The notation ``$\mu^{\prime}=\mathfrak{be}(\mu)$'' defined in \S
\ref{sec:Borel_Measure} means $\mu^{\prime}$ is the unique Borel
extension of $\mu$.%
} is $\mathbf{m}$-tight in $A$.
\end{prop}
\begin{proof}
We first show $\mathbf{m}$-tightness in $A$ implies the existence
of $\mu^{\prime}=\mathfrak{be}(\mu)$ for each $\mu\in\Gamma$. Given
such tightness, $\mu$ is supported on some $B\in\mathscr{K}_{\sigma}^{\mathbf{m}}(A,\mathscr{O}_{S}(A))$.
$\mathscr{B}_{S}(B)=\mathscr{A}|_{B}$ by Lemma \ref{lem:MC_Prod_Meas}
(with $B_{i}=\mathfrak{p}_{i}(A)$). Then, the unique existence of
$\mu^{\prime}$ follows by Lemma \ref{lem:Union_Borel_Prod_Equal}
(c) (with $A=B$). Now, $\mathscr{K}^{\mathbf{m}}(A,\mathscr{O}_{S}(A))\subset\mathscr{A}$
by Lemma \ref{lem:MC_Prod_Meas} (with $B_{i}=\mathfrak{p}_{i}(A)$).
So, the $\mathbf{m}$-tightness of $\Gamma$ and that of $\{\mu^{\prime}\}_{\mu\in\Gamma}$
(if any) are equivalent.\end{proof}

\subsection{\label{sub:Sigma_Metrizable_Compact}$\sigma$-metrizable compact
subsets}

$\sigma$-metrizable compact subsets inherit many nice properties
from its metrizable compact components.
\begin{prop}
\label{prop:Sigma_MC}Let $E$ be a topological space, $\{K_{n}\}_{n\in\mathbf{N}}\subset\mathscr{K}(E)$,
$A=\bigcup_{n\in\mathbf{N}}K_{n}$ and $\mathcal{D}\subset C(E;\mathbf{R})$.
Consider the following statements:

\renewcommand{\labelenumi}{(\alph{enumi})}
\begin{enumerate}
\item $(A,\mathscr{O}_{E}(A))$ is a Souslin space.
\item $\{K_{n}\}_{n\in\mathbf{N}}\subset\mathscr{K}^{\mathbf{m}}(E)$ (hence
$A\in\mathscr{K}_{\sigma}^{\mathbf{m}}(E)$).
\item $(K_{n},\mathscr{O}_{E}(K_{n}))$ is a baseable space for all $n\in\mathbf{N}$.
\item $(A,\mathscr{O}_{E}(A))$ is a baseable space.
\item $A$ is a $\mathcal{D}$-baseable subset of $E$.
\item $A\in\mathscr{B}^{\mathbf{s}}(E)$.
\end{enumerate}
Then, (a) - (e) are successively stronger. (b) implies (f). Moreover,
if $\mathcal{D}$ separates points on $E$, then (a) - (e) are all
equivalent.
\end{prop}
\begin{proof}
((b) $\rightarrow$ (a)) Each $(K_{n},\mathscr{O}_{E}(K_{n}))$ is
Souslin by Proposition \ref{prop:MC} (b, f). Hence, (a) follows by
Proposition \ref{prop:Var_Polish} (g).

((c) $\rightarrow$ (b)) follows by Proposition \ref{prop:MC} (d,
f) (with $K=K_{n}$). 

((d) $\rightarrow$ (c)) follows by Fact \ref{fact:Baseable_Space_Subset}
(a, b) (with $E=(A,\mathscr{O}_{E}(A))$ and $A=K_{n}$).

((e) $\rightarrow$ (d)) is automatic by definition.

((b) $\rightarrow$ (f)) follows by Proposition \ref{prop:MC} (f,
g) and Proposition \ref{prop:Baseable_SB_Union}.

When $\mathcal{D}$ separates points on $E$, $A\in\mathscr{B}(E)$
by Proposition \ref{prop:Fun_Sep_1} (e) (with $A=E$) and Proposition
\ref{prop:Compact} (a), and (a) implies (e) by Proposition \ref{prop:Hered_Lindelof_Baseable}.\end{proof}

Below are several constructive properties of $\sigma$-metrizable
compact subsets.
\begin{lem}
\label{lem:Sigma_MC_Union}Let $E$ be a topological space and $\{A_{n}\}_{n\in\mathbf{N}}\subset\mathscr{K}_{\sigma}^{\mathbf{m}}(E)$.
If $A=\bigcup_{n\in\mathbf{N}}A_{n}$ is a Hausdorff subspace of $E$,
then there exist $\{K_{q}\}_{q\in\mathbf{N}}\subset\mathscr{K}^{\mathbf{m}}(E)$
such that $A=\bigcup_{q\in\mathbf{N}}K_{q}$ and $K_{q}\subset K_{q+1}$
for all $q\in\mathbf{N}$.
\end{lem}
\begin{proof}
Let $A_{n}=\bigcup_{p\in\mathbf{N}}K_{p,n}$ with $\{K_{p,n}\}_{p\in\mathbf{N}}\subset\mathscr{K}^{\mathbf{m}}(E)$
for each $n\in\mathbf{N}$. Define $K_{q}\circeq\bigcup_{p=1}^{q}\bigcup_{n=1}^{q}K_{p,n}$
so $K_{q}\subset K_{q+1}$ for all $q\in\mathbf{N}$ and $A=\bigcup_{q\in\mathbf{N}}K_{q}$.
Each $(K_{q},\mathscr{O}_{E}(K_{q}))$ is Hausdorff by Proposition
\ref{prop:Separability} (c) and hence metrizable by Lemma \ref{lem:MC_Union}.\end{proof}

\begin{lem}
\label{lem:Sigma_MC_Prod}Let $\mathbf{I}$ be a countable index set,
$\{S_{i}\}_{i\in\mathbf{I}}$ be topological spaces and $S\circeq\prod_{i\in\mathbf{I}}S_{i}$.
Then:

\renewcommand{\labelenumi}{(\alph{enumi})}
\begin{enumerate}
\item If $\mathbf{I}$ is finite and $A_{i}\in\mathscr{K}_{\sigma}^{\mathbf{m}}(S_{i})$
for all $i\in\mathbf{I}$, then $\prod_{i\in\mathbf{I}}A_{i}\in\mathscr{K}_{\sigma}^{\mathbf{m}}(S)$.
\item If $A\in\mathscr{K}_{\sigma}^{\mathbf{m}}(S)$ and $\mathfrak{p}_{i}(A)$
is a Hausdorff subspace of $S_{i}$ for all $i\in\mathbf{I}$, then
$\mathfrak{p}_{i}(A)\in\mathscr{K}_{\sigma}^{\mathbf{m}}(S_{i})$
for all $i\in\mathbf{I}$.
\end{enumerate}
\end{lem}
\begin{proof}
(a) Without loss of generality, we suppose $\mathbf{I}=\{1,...,d\}$
and let $A_{i}=\bigcup_{p\in\mathbf{N}}K_{p,i}$ with $\{K_{p,i}\}_{p\in\mathbf{N}}\subset\mathscr{K}^{\mathbf{m}}(S_{i})$
for each $1\leq i\leq d$. We have that
\begin{equation}
\prod_{i=1}^{d}A_{i}\supset F_{p_{1},...,p_{d}}\circeq\prod_{i\in\mathbf{I}}K_{p_{i},i}\in\mathscr{K}^{\mathbf{m}}(S),\;\forall p_{1},...,p_{d}\in\mathbf{N}\label{eq:Check_Sigma_MC_Prod_1}
\end{equation}
by Lemma \ref{lem:MC_Prod} (a). For any $x\in\prod_{i=1}^{d}A_{i}$,
there exist $p_{1},...,p_{d}\in\mathbf{N}$ such that
\begin{equation}
\mathfrak{p}_{i}(x)\in K_{p_{i},i},\;\forall1\leq i\leq d.\label{eq:Check_Sigma_MC_Prod_2}
\end{equation}
It then follows by (\ref{eq:Check_Sigma_MC_Prod_1}) that
\begin{equation}
\prod_{i=1}^{d}A_{i}=\bigcup_{(p_{1},...,p_{d})\in\mathbf{N}^{d}}F_{p_{1},...,p_{d}}\in\mathscr{K}_{\sigma}^{\mathbf{m}}(S).\label{eq:Check_Sigma_MC_Prod_3}
\end{equation}

(b) Let $A=\bigcup_{p\in\mathbf{N}}K_{p}$ with $\{K_{p}\}_{p\in\mathbf{N}}\subset\mathscr{K}^{\mathbf{m}}(S)$
so $\mathfrak{p}_{i}(K_{p})\in\mathscr{K}^{\mathbf{m}}(S_{i})$ for
all $p\in\mathbf{N}$ and $i\in\mathbf{I}$ by the fact $\mathfrak{p}_{i}(K_{p})\subset\mathfrak{p}_{i}(A)$,
the Hausdorff property of $\mathfrak{p}_{i}(A)$, Proposition \ref{prop:Separability}
(c) and Lemma \ref{lem:MC_Prod} (b). Hence, $\mathfrak{p}_{i}(A)=\mathfrak{p}_{i}(\bigcup_{p\in\mathbf{N}}K_{p})=\bigcup_{p\in\mathbf{N}}\mathfrak{p}_{i}(K_{p})\in\mathscr{K}_{\sigma}^{\mathbf{m}}(S_{i})$
for all $i\in\mathbf{I}$.\end{proof}

\subsection{\label{sub:Sko_Baseable}Baseability about Skorokhod $\mathscr{J}_{1}$-space}

When $E$ is a Tychonoff space, the associated Skorokhod $\mathscr{J}_{1}$-space
$D(\mathbf{R}^{+};E)$ is also (see Proposition \ref{prop:Sko_Basic_1}
(c)). The following proposition shows that baseability of $E$ passes
to $D(\mathbf{R}^{+};E)$.
\begin{prop}
\label{prop:Sko_Baseable}Let $E$ be a Tychonoff space. Then:

\renewcommand{\labelenumi}{(\alph{enumi})}
\begin{enumerate}
\item If $E$ is a $\mathcal{D}$-baseable space with $\mathcal{D}\subset C_{b}(E;\mathbf{R})$,
then $D(\mathbf{R}^{+};E)$%
\footnote{The Skorokhod $\mathscr{J}_{1}$-space $D(\mathbf{R}^{+};E)$ was
defined in \S \ref{sub:Meas_Cont_Cadlag_Map}.%
} is a $\{\alpha_{t,n}^{f}:f\in\mathcal{D},t\in\mathbf{Q}^{+}\}$-baseable
space with
\begin{equation}
\alpha_{t,n}^{f}(x)\circeq n\int_{t}^{t+1/n}f(x(s))ds\label{eq:Sko_SP_Fun}
\end{equation}
for each $f\in\mathcal{D}$, $t\in\mathbf{Q}^{+}$ and $n\in\mathbf{N}$.
\item If $E$ is a baseable space, then $D(\mathbf{R}^{+};E)$ is also a
baseable space and $J(x)\subset(0,\infty)$%
\footnote{The notation ``$J(x)$'' was defined in \S \ref{sub:Map}.%
} is at most countable for all $x\in D(\mathbf{R}^{+};E)$.
\end{enumerate}
\end{prop}
\begin{proof}
(a) Without loss of generality, we suppose $\mathcal{D}$ is countable.
Then, (a) follows immediately by Proposition \ref{prop:Sko_Basic_1}
(b).

(b) There exists a countable $\mathcal{D}\subset C_{b}(E;\mathbf{R})$
separating points on $E$ by Fact \ref{fact:Baseable_Subsets} (a,
b) (with $A=E$) so $\varphi\circeq\bigotimes\mathcal{D}$ is injective.
$\varphi\in C(E;\mathbf{R}^{\mathcal{D}})$ by Fact \ref{fact:Prod_Map_2}
(b) and $D(\mathbf{R}^{+};E)$ is baseable by (a). $\mathbf{R}^{\mathcal{D}}$
is Polish by Proposition \ref{prop:Var_Polish} (f). Therefore,
\begin{equation}
J(x)=J\left[\varpi(\varphi)(x)\right],\;\forall x\in D(\mathbf{R}^{+};E)\label{eq:J(x)=00003DJ(phi(x))}
\end{equation}
by Proposition \ref{prop:Sko_Basic_1} (d) (with $S=E$, $E=\mathbf{R}^{\mathcal{D}}$
and $f=\varphi$) and $J[\varpi(\varphi)(x)]$ is a countable subset
of $(0,\infty)$ by \cite[\S 3.5, Lemma 5.1]{EK86}.\end{proof}

Metrizability of compact subsets of $E$ also passes to $D(\mathbf{R}^{+};E)$.
\begin{prop}
\label{prop:Sko_Compact}If $E$ is a Tychonoff space with $\mathscr{K}(E)=\mathscr{K}^{\mathbf{m}}(E)$,
then $\mathscr{K}(D(\mathbf{R}^{+};E))=\mathscr{K}^{\mathbf{m}}(D(\mathbf{R}^{+};E))$.
\end{prop}
\begin{proof}
Let $K\in\mathscr{K}(D(\mathbf{R}^{+};E))$. By Proposition \ref{prop:Sko_Compact_Containment},
there exist $\{A_{n}\}_{n\in\mathbf{N}}\subset\mathscr{K}(E)=\mathscr{K}^{\mathbf{m}}(E)$
such that $K\subset D(\mathbf{R}^{+};A)$%
\footnote{$D(\mathbf{R}^{+};\bigcup_{n\in\mathbf{N}}A_{n})$ is well-defined
by Corollary \ref{cor:Sko_Subspace}.%
} with $A\circeq\bigcup_{n\in\mathbf{N}}A_{n}$. $A$ is a baseable
space by Proposition \ref{prop:CR} (a, c) and Proposition \ref{prop:Sigma_MC}
(b, d) (with $K_{n}=A_{n}$ and $\mathcal{D}=C_{b}(E;\mathbf{R})$).
$D(\mathbf{R}^{+};A)$ is a baseable space by Proposition \ref{prop:Sko_Baseable}
(b) (with $E=A$). $K$ is a baseable subspace of $D(\mathbf{R}^{+};E)$
by Fact \ref{fact:Baseable_Space_Subset} (a, b) (with $E=A$) and
Corollary \ref{cor:Sko_Subspace}. Hence, $K$ is metrizable by Proposition
\ref{prop:MC} (d, f) (with $E=D(\mathbf{R}^{+};E)$).\end{proof}

The countability of an $E$-valued c$\grave{\mbox{a}}$dl$\grave{\mbox{a}}$g
process's fixed left-jump times is well-known when $E$ is metrizable
and separable. We extend this to baseable Tychonoff spaces.
\begin{prop}
\label{prop:J(Mu)_J(X)_Baseable}Let $E$ be a baseable Tychonoff
space. Then:

\renewcommand{\labelenumi}{(\alph{enumi})}
\begin{enumerate}
\item For any $\mu\in\mathfrak{M}^{+}(D(\mathbf{R}^{+};E),\mathscr{B}(E)^{\otimes\mathbf{R}^{+}}|_{D(\mathbf{R}^{+};E)})$%
\footnote{Recall that $\mathscr{B}(E)^{\otimes\mathbf{R}^{+}}|_{D(\mathbf{R}^{+};E)}$
is generally smaller than $\mathscr{B}(D(\mathbf{R}^{+};E))$. %
} (especially $\mu\in\mathcal{M}^{+}(D(\mathbf{R}^{+};E))$), $J(\mu)$%
\footnote{$J(\mu)$, the set of fixed left-jump times of $\mu$ was defined
in (\ref{eq:J(Mu)}).%
} is a well-defined countable subset of $(0,\infty)$.
\item For any $E$-valued c$\grave{\mbox{a}}$dl$\grave{\mbox{a}}$g process
$X$, $J(X)$%
\footnote{$J(X)$, the set of fixed left-jump times of $X$ was defined in (\ref{eq:J(X)}).%
} is a well-defined countable subset of $(0,\infty)$.
\end{enumerate}
In particular, (a) and (b) hold when $E$ is a metrizable and separable
space or a Polish space.
\end{prop}
\begin{proof}
(a) Metrizable and separable spaces are baseable spaces by Fact \ref{fact:Baseable_Metrizable_Separable}
(b). As $E$ is a baseable space, there exists a countable $\mathcal{D}\subset C_{b}(E;\mathbf{R})$
separating points on $E$ by Fact \ref{fact:Baseable_Subsets} (a,
b) (with $A=E$) so $\varphi\circeq\bigotimes\mathcal{D}$ is injective.
$\varphi\in C(E;\mathbf{R}^{\mathcal{D}})$ by Fact \ref{fact:Prod_Map_2}
(b).
\begin{equation}
\varpi(\varphi)\in C\left(D(\mathbf{R}^{+};E);D(\mathbf{R}^{+};\mathbf{R}^{\mathcal{D}})\right)\label{eq:Path_Mapping_(Prod(D))_Cont}
\end{equation}
by Proposition \ref{prop:Sko_Basic_1} (d) (with $S=E$, $E=\mathbf{R}^{\mathcal{D}}$
and $f=\varphi$), which implies (\ref{eq:J(x)=00003DJ(phi(x))}).
As $\varpi(\varphi)$ maps $D(\mathbf{R}^{+};E)$ into $D(\mathbf{R}^{+};\mathbf{R}^{\mathcal{D}})$,
we have that
\begin{equation}
\begin{aligned}\varpi(\varphi)\in & M\left(D(\mathbf{R}^{+};E),\mathscr{B}(E)^{\otimes\mathbf{R}^{+}}|_{D(\mathbf{R}^{+};E)};\right.\\
 & \left.D(\mathbf{R}^{+};\mathbf{R}^{\mathcal{D}}),\mathscr{B}(\mathbf{R}^{\mathcal{D}})^{\otimes\mathbf{R}^{+}}|_{D(\mathbf{R}^{+};\mathbf{R}^{\mathcal{D}})}\right)
\end{aligned}
\label{eq:Path_Mapping_(Prod(D))_Meas}
\end{equation}
by Fact \ref{fact:Path_Mapping} (b) (with $f=\varphi$).
\begin{equation}
\nu\circeq\mu\circ\varpi(\varphi)^{-1}\in\mathfrak{M}^{+}\left(D(\mathbf{R}^{+};\mathbf{R}^{\mathcal{D}}),\mathscr{B}(\mathbf{R}^{\mathcal{D}})^{\otimes\mathbf{R}^{+}}|_{D(\mathbf{R}^{+};\mathbf{R}^{\mathcal{D}})}\right)\label{eq:J(Path_Map_Dist)}
\end{equation}
by (\ref{eq:Path_Mapping_(Prod(D))_Meas}). $\mathbf{R}^{\mathcal{D}}$
is a Polish space, so (\ref{eq:J(x)=00003DJ(phi(x))}) implies 
\begin{equation}
\begin{aligned} & \mu\left(\left\{ x\in D(\mathbf{R}^{+};E):t\in J(x)\right\} \right)=\mu\left(\left\{ x\in D(\mathbf{R}^{+};E):t\in J\left[\varpi(\varphi)(x)\right]\right\} \right)\\
 & =\nu\left(\left\{ y\in D(\mathbf{R}^{+};\mathbf{R}^{\mathcal{D}}):t\in J(y)\right\} \right),\;\forall t\in\mathbf{R}^{+},
\end{aligned}
\label{eq:Check_J(Mu)_Countable}
\end{equation}
while the equalities in (\ref{eq:Check_J(Mu)_Countable}) as well
as $J(\mu)$ and $J(\nu)$ are well-defined by Fact \ref{fact:J(Mu)_Well_Defined}.
Hence, we have $J(\mu)=J(\nu)$ by (\ref{eq:Check_J(Mu)_Countable})
and this set is a countable subset of $(0,\infty)$ by \cite[\S 3.7, Lemma 7.7]{EK86}.

(b) follows immediately by Fact \ref{fact:Cadlag_Proc_Dist} (a),
the definition of $J(X)$ and (a) (with $\mu=\mathrm{pd}(X)|_{D(\mathbf{R}^{+};E)}$%
\footnote{$\mathrm{pd}(X)$ denotes the process distribution of $X$ and was
specified in \S \ref{sec:Proc}.%
}).\end{proof}

\chapter{\label{chap:RepFun}Replication of Function and Operator}

\chaptermark{Replica Function and Operator}

The previous chapter introduced space change through a base $(E_{0},\mathcal{F};\widehat{E},\widehat{\mathcal{F}})$
over topological space $E$ and the notions of baseable spaces and
subsets. Now, we discuss the replication of objects from $E$ onto
$\widehat{E}$ and the association of the original and replica objects.
\S \ref{sec:RepFun} introduces the replicas of continuous functions.
\S \ref{sec:RepOP} introduces the replicas of linear operators on
$C_{b}(E;\mathbf{R})$. The replica operators constructed in \S \ref{sub:RepOP_Base}
are strong generators of semigroups on $C(\widehat{E};\mathbf{R})$.

\section{\label{sec:RepFun}Replica function}

Given a base $(E_{0},\mathcal{F};\widehat{E},\widehat{\mathcal{F}})$
over topological space $E$, replicating a function $f\in M(E;\mathbf{R})$
onto $\widehat{E}$ basically means extending $f|_{E_{0}}$ onto $\widehat{E}$.
A naive approach preserves the values of $f$ on $E_{0}$ and assigns
$0$ on $\widehat{E}\backslash E_{0}$. We make the following general
notation for simplicity.
\begin{notation}
\label{notation:Var}Let $E$, $S_{1}$ and $S_{2}$ be non-empty
sets, $A$ be an arbitrary subset of $S_{1}\cap S_{2}$, $y_{0}\in E$
and $f\in E^{S_{1}}$. By $\mathfrak{var}(f;S_{2},A,y_{0})$%
\footnote{``$\mathfrak{var}$'' is ``var'' in fraktur font which stands
for ``variant''.%
} we denote the mapping 
\begin{equation}
\mathfrak{var}(f;S_{2},A,y_{0})(x)\circeq\begin{cases}
f(x), & \mbox{if }x\in A,\\
y_{0}, & \mbox{otherwise},
\end{cases}\;\forall x\in S_{2}\label{eq:(A,y0)_Ver}
\end{equation}
from $S_{2}$ to $E$.\end{notation}
\begin{note}
\label{note:Var}$\mathfrak{var}(f;S_{2},f^{-1}(\{y_{0}\}),y_{0})=\mathfrak{var}(f;S_{2},\varnothing,y_{0})$
maps all $x\in S_{2}$ to $y_{0}$.
\end{note}

$\mathfrak{var}(f;\widehat{E},E_{0},0)$ realizes the naive idea of
replicating $f$ above, but it may not preserve the Borel measurability
of $f$ if $E_{0}$ is not a standard Borel set. Noticing that (\ref{eq:F_Fhat_Coincide})
links members of $\mathcal{F}$ bijectively to a member of $\widehat{\mathcal{F}}\subset C(\widehat{E};\mathbf{R})$,
we can define the replica of a suitable continuous function on $E$
as a continuous function on $\widehat{E}$.
\begin{defn}
\label{def:RepFun}Let $E$ be a topological space, $(E_{0},\mathcal{F};\widehat{E},\widehat{\mathcal{F}})$
be a base over $E$ and $d,k\in\mathbf{N}$. The\textbf{ replica of
$f\in C(E^{d};\mathbf{R}^{k})$} with respect to $(E_{0},\mathcal{F};\widehat{E},\widehat{\mathcal{F}})$
(if any) refers to the continuous extension $\widehat{f}$ of $f|_{E_{0}^{d}}$
on $\widehat{E}^{d}$.\end{defn}
\begin{rem}
\label{rem:ContRep}Remark \ref{rem:Compactification} stated that
the compactification inducing $\widehat{E}$ does not necessarily
extend every member of $C_{b}(E_{0},\mathscr{O}_{\mathcal{F}}(E_{0});\mathbf{R})$
continuously onto $\widehat{E}$. So, a general $f\in C_{b}(E^{d};\mathbf{R}^{k})$
need not have a replica.\end{rem}
\begin{notation}
\label{notation:Rep_Fun}Let $(E_{0},\mathcal{F};\widehat{E},\widehat{\mathcal{F}})$
be a base over $E$ and $d,k\in\mathbf{N}$. Hereafter, we \textit{will
always} let $\overline{f}$ denote $\mathfrak{var}(f;\widehat{E},E_{0},0)$
for $f\in(\mathbf{R}^{k})^{E^{d}}$ and $\widehat{f}$ denote the
replica of $f\in C(E^{d};\mathbf{R}^{k})$ if no confusion is caused.
\end{notation}
Below are several basic properties of $\overline{f}$ and $\widehat{f}$.
\begin{prop}
\label{prop:RepFun_Basic}Let $E$ be a topological space, $(E_{0},\mathcal{F};\widehat{E},\widehat{\mathcal{F}})$
be a base over $E$ and $d,k\in\mathbf{N}$. Then:

\renewcommand{\labelenumi}{(\alph{enumi})}
\begin{enumerate}
\item If $f\in(\mathbf{R}^{k})^{E^{d}}$ is bounded, then $\overline{f}$
is also bounded.
\item If $f\in M(E^{d};\mathbf{R}^{k})$ and $f|_{E_{0}^{d}\backslash A}=0$
for some $A\in\mathscr{B}^{\mathbf{s}}(E^{d})$ with $A\subset E_{0}^{d}$,
then $\overline{f}\in M(\widehat{E}^{d};\mathbf{R}^{k})$. In particular,
this is true if $E_{0}^{d}\in\mathscr{B}^{\mathbf{s}}(E^{d})$.
\item The replica of $f\in C(E^{d};\mathbf{R}^{k})$ (if any) is unique.
\item If $f_{1},f_{2}\in C(E^{d};\mathbf{R}^{k})$ have replicas, then $a\widehat{f}_{1}+b\widehat{f}_{2}$
(resp. $\widehat{f}_{1}\widehat{f}_{2}$ when $k=1$) is the replica
of $af_{1}+bf_{2}$ for all $a,b\in\mathbf{R}$ (resp. $f_{1}f_{2}$).
\item $\widehat{\mathcal{F}}=\{\widehat{f}:f\in\mathcal{F}\}$, $\mathfrak{ag}(\widehat{\mathcal{F}})=\{\widehat{f}:f\in\mathfrak{ag}(\mathcal{F})\}$
and $\mathfrak{ag}[\Pi^{d}(\widehat{\mathcal{F}})]=\{\widehat{f}:f\in\mathfrak{ag}(\Pi^{d}(\mathcal{F}))\}$.
\item $f\in C(E^{d};\mathbf{R}^{k})$ admits a replica if and only if
\begin{equation}
\mathfrak{p}_{i}\circ f|_{E_{0}^{d}}\in\mathfrak{ca}\left[\Pi^{d}\left(\mathcal{F}|_{E_{0}}\right)\right],\;\forall1\leq i\leq k.\label{eq:ContRepFun_3}
\end{equation}
In particular, every $f\in\mathfrak{ca}[\Pi^{d}(\mathcal{F})]$ admits
a replica.
\end{enumerate}
\end{prop}
\begin{proof}
(a) The definition of $\overline{f}$ implies $\Vert\overline{f}\Vert_{\infty}\leq\Vert f\Vert_{\infty}$.

(b) $\overline{h}|_{A}=h|_{A}\in M(A,\mathscr{O}_{\widehat{E}^{d}}(A);\mathbf{R}^{k})$
by Lemma \ref{lem:SB_Base} (a). So, $\overline{h}=\overline{h}\mathbf{1}_{A}\in M(\widehat{E}^{d};\mathbf{R}^{k})$%
\footnote{$\mathbf{1}_{A}$ denotes the indicator function of $A$.%
} by Fact \ref{fact:Indicator_Modify} (with $E=\widehat{E}^{d}$,
$\mathscr{U}=\mathscr{B}(\widehat{E}^{d})$ and $f=\overline{h}$).

(c) follows since $E_{0}^{d}$ is dense in $\widehat{E}^{d}$, $f|_{E_{0}^{d}}=\widehat{f}|_{E_{0}^{d}}$
and $f$ and $\widehat{f}$ are continuous.

(d) follows by the fact that $(af_{1}+bf_{2})|_{E_{0}^{d}}=(a\widehat{f}_{1}+b\widehat{f}_{2})|_{E_{0}^{d}}$
and $f_{1}f_{2}|_{E_{0}^{d}}=\widehat{f}_{1}\widehat{f}_{2}|_{E_{0}^{d}}$,
that $a\widehat{f}_{1}+b\widehat{f}_{2}\in C(\widehat{E}^{d};\mathbf{R}^{k})$
and that $\widehat{f}_{1}\widehat{f}_{2}\in C(\widehat{E}^{d};\mathbf{R})$
when $k=1$.

(e) We let $\mathcal{F}=\{f_{n}\}_{n\in\mathbf{N}}$ and find $\widehat{\mathcal{F}}=\{\widehat{f}_{n}\}_{n\in\mathbf{N}}$
by (\ref{eq:F_Fhat_Coincide}) and Lemma \ref{lem:Base} (a). Let
$1\leq l\leq d$, $n_{1},...,n_{l}\in\mathbf{N}$, $f\circeq\prod_{i=1}^{l}f_{n_{i}}\circ\mathfrak{p}_{i}$
and $g\circeq\prod_{i=1}^{l}\widehat{f}_{n_{i}}\circ\mathfrak{p}_{i}$.
Then,
\begin{equation}
f|_{E_{0}^{d}}=\prod_{i=1}^{l}f_{n_{i}}|_{E_{0}}\circ\mathfrak{p}_{i}=\prod_{i=1}^{l}\widehat{f}_{n_{i}}|_{E_{0}}\circ\mathfrak{p}_{i}=g|_{E_{0}^{d}}.\label{eq:Check_ca(Pi^d(F))_Rep_1}
\end{equation}
and (\ref{eq:Pi^d(Fhat)_Cb}) imply $g=\widehat{f}$. The members
of $\mathfrak{ag}(\mathcal{F})$ (resp. $\mathfrak{ag}[\Pi^{d}(\mathcal{F})]$)
correspond bijectively to those of $\mathfrak{ag}(\widehat{\mathcal{F}})$
(resp. $\mathfrak{ag}[\Pi^{d}(\widehat{\mathcal{F}})]$). Hence, (e)
follows by (d).

(f - Necessity) If $\widehat{f}$ exists, then $\{\mathfrak{p}_{i}\circ\widehat{f}\}_{1\leq i\leq k}\subset C(\widehat{E}^{d};\mathbf{R})$
by Fact \ref{fact:Prod_Map_2} (a). Hence, (\ref{eq:ContRepFun_3})
follows by (\ref{eq:ca(Pi^d(Ftilte))_C(Ehat)}).

(f - Sufficiency) If (\ref{eq:ContRepFun_3}) holds, then $\{\widehat{\mathfrak{p}_{i}\circ f}\}_{1\leq i\leq k}$
exists by Corollary \ref{cor:Base_Fun_Dense}. Hence, $\bigotimes_{i=1}^{k}\widehat{\mathfrak{p}_{i}\circ f}=\widehat{f}$
by Fact \ref{fact:Prod_Map_2} (b) and the fact
\begin{equation}
f|_{E_{0}^{d}}=\bigotimes_{i=1}^{k}\mathfrak{p}_{i}\circ f|_{E_{0}^{d}}=\bigotimes_{i=1}^{k}\widehat{\mathfrak{p}_{i}\circ f}|_{E_{0}^{d}}.\label{eq:Check_ca(Pi^d(F))_Rep_2}
\end{equation}
\end{proof}

\begin{note}
\label{note:Rep_Fun}For the sake of brevity, hereafter we may use
the replica of $f\in\mathfrak{ca}(\Pi^{d}(\mathcal{F}))$ without
referring to Proposition \ref{prop:RepFun_Basic} (f) for its existence.
\end{note}
We next show a nice property of locally compact baseable spaces which
recovers \cite[Corollary 2.3.32]{S98}. This is also an example where
$\overline{f}$ and $\widehat{f}$ coincide.
\begin{prop}
\label{prop:LC_Polish}Let $E$ be a locally compact space and $\mathcal{D}\subset C_{0}(E;\mathbf{R})$%
\footnote{$C_{0}(E;\mathbf{R})$, the family of all $\mathbf{R}$-valued continuous
functions on $E$ that vanishes at infinity, was defined in \S \ref{sub:Fun}.%
}. Consider the following statements:

\renewcommand{\labelenumi}{(\alph{enumi})}
\begin{enumerate}
\item $E$ is a $\mathcal{D}$-baseable space.
\item There exists a base $(E,\mathcal{F};\widehat{E},\widehat{\mathcal{F}}\}$
over $E$ such that $\widehat{E}$ is a one-point compactification
of $E$, $(\mathcal{F}\backslash\{1\})\subset\mathcal{D}\subset C_{0}(E;\mathbf{R})\subset\mathfrak{ca}(\mathcal{F})$
and $\mathcal{F}$ strongly separates points on $E$.
\item $E$ is a Polish space.
\item $E$ is a metrizable and separable space.
\item $E$ is a $C_{0}(E;\mathbf{R})$-baseable space.
\end{enumerate}
Then, (a) - (e) are successively weaker. Moreover, (e) implies (a)
when $\mathcal{D}$ is uniformly dense in $C_{0}(E;\mathbf{R})$.
\end{prop}
\begin{proof}
((a) $\rightarrow$ (b)) By (a), there exists a countable $\mathcal{F}\subset(\mathcal{D}\cup\{1\})\subset C_{b}(E;\mathbf{R})$
that separates points on $E$. $E$ is a Hausdorff space by Fact \ref{fact:Baseable_Metrizable_Separable}
(a) and admits a one-point compactification $\widehat{E}$ by Proposition
\ref{prop:One-Point}. It follows by Lemma \ref{lem:SP_on_LC} (with
$\mathcal{D}=\mathcal{F}$) that $\mathcal{F}$ strongly separates
points on $E$ and
\begin{equation}
\widehat{\mathcal{F}}\circeq\left\{ \mathfrak{var}(f:\widehat{E},E,0):f\in\mathcal{F}\backslash\{1\}\right\} \cup\{1\}\subset C(\widehat{E};\mathbf{R})\label{eq:Fhat_LC}
\end{equation}
separates and strongly separates points on $\widehat{E}$. Hence,
$(E,\mathcal{F};\widehat{E},\widehat{\mathcal{F}})$ by definition
is a base over $E$. Moreover, we get $C_{0}(E;\mathbf{R})\subset\mathfrak{ca}(\mathcal{F})$
by Corollary \ref{cor:Base_Fun_Dense} (with $d=1$ and $E_{0}=E$)
and Fact \ref{fact:Cc_C0_Cb}.

((b) $\rightarrow$ (c)) $\mathscr{O}(E)=\mathscr{O}_{\mathcal{F}}(E)$
by (b), so $E$ is an open subspace of the Polish space $\widehat{E}$
by Proposition \ref{prop:Separability} (a) and Lemma \ref{lem:Base}
(b, c). Hence, (c) follows by Proposition \ref{prop:Var_Polish} (b).

((c) $\rightarrow$ (d)) follows by Proposition \ref{prop:Var_Polish}
(c).

((d) $\rightarrow$ (e)) $C_{0}(E;\mathbf{R})$ separates points on
$E$ by Proposition \ref{prop:Metrizable} (a) and Proposition \ref{prop:LC_CR}
(a, d). Then, (e) follows by Proposition \ref{prop:Metrizable} (c)
and Proposition \ref{prop:Hered_Lindelof_Baseable} (with $A=E$ and
$\mathcal{D}=C_{0}(E;\mathbf{R})$).

Moreover, if $C_{0}(E;\mathbf{R})\subset\mathfrak{cl}(\mathcal{D})$,
then (e) implies (a) by Proposition \ref{prop:Dense_Baseable} (with
$A=E$, $\mathcal{D}=C_{0}(E;\mathbf{R})$ and $\mathcal{D}_{0}=\mathcal{D}$).\end{proof}

\section{\label{sec:RepOP}Replica operator}

We now focus on replicating a linear operator $\mathcal{L}$ on $C_{b}(E;\mathbf{R})$
as a linear operator on $C(\widehat{E};\mathbf{R})$. Many general
concepts about linear operators used below were reviewed in \S \ref{sub:OP}
and, as noted in \S \ref{sec:Convention}, we always consider single-valued
operators.

\subsection{\label{sub:RepOP_Definition}Definition}

Replicating $\mathcal{L}$ from $C_{b}(E;\mathbf{R})$ onto $C(\widehat{E};\mathbf{R})$
means constructing a linear operator on $C(\widehat{E};\mathbf{R})$
whose domain and range are formed by the replicas of the member of
$\mathfrak{D}(\mathcal{L})$ and $\mathfrak{R}(\mathcal{L})$ respectively.
Below is a mild sufficient condition for its existence.
\begin{prop}
\label{prop:RepOP_Exist}Let $E$ be a topological space, $(E_{0},\mathcal{F};\widehat{E},\widehat{\mathcal{F}})$
be a base over $E$ and $\mathcal{L}$ be a linear operator on $C_{b}(E;\mathbf{R})$
such that
\begin{equation}
\begin{aligned} & \mathfrak{mc}(\mathcal{F})\subset\mathfrak{D}(\mathcal{L}),\\
 & \mathfrak{R}\left(\mathcal{L}|_{\mathfrak{mc}(\mathcal{F})}\right)\subset\mathfrak{ca}(\mathcal{F}).
\end{aligned}
\label{eq:Base_for_OP}
\end{equation}
Then:

\renewcommand{\labelenumi}{(\alph{enumi})}
\begin{enumerate}
\item There exists a unique linear operator $\widehat{\mathcal{L}}_{0}$
on $C(\widehat{E};\mathbf{R})$ such that
\begin{equation}
\mathfrak{D}(\widehat{\mathcal{L}}_{0})=\mathfrak{ag}(\widehat{\mathcal{F}})\label{eq:Lhat0_Domain}
\end{equation}
and
\begin{equation}
\widehat{\mathcal{L}}_{0}\widehat{f}=\widehat{\mathcal{L}f},\;\forall\widehat{f}\in\mathfrak{ag}(\widehat{\mathcal{F}}).\label{eq:Lhat0_Rep}
\end{equation}

\item There exists a unique linear operator $\widehat{\mathcal{L}}_{1}$
on $C(\widehat{E};\mathbf{R})$ such that
\begin{equation}
\mathfrak{D}(\widehat{\mathcal{L}}_{1})=\left\{ \widehat{f}:(f,\mathcal{L}f)\in\mathcal{L}\cap\left(\mathfrak{ca}(\mathcal{F})\times\mathfrak{ca}(\mathcal{F})\right)\right\} \label{eq:Lhat1_Domain}
\end{equation}
and
\begin{equation}
\widehat{\mathcal{L}}_{1}\widehat{f}=\widehat{\mathcal{L}f},\;\forall\widehat{f}\in\mathfrak{D}(\widehat{\mathcal{L}}_{1}).\label{eq:Lhat1_Rep}
\end{equation}

\end{enumerate}
\end{prop}
\begin{proof}
We have that
\begin{equation}
\mathfrak{ag}(\mathcal{F})=\mathfrak{ac}\left(\left\{ af:f\in\mathfrak{mc}(\mathcal{F}),a\in\mathbf{R}\right\} \right)\subset\mathfrak{D}(\mathcal{L})\label{eq:ag_(F)_in_D(L)}
\end{equation}
and
\begin{equation}
\begin{aligned}\mathfrak{R}\left(\mathcal{L}|_{\mathfrak{ag}(\mathcal{F})}\right) & =\left\{ \mathcal{L}h:h\in\mathfrak{ac}\left(\left\{ af:f\in\mathfrak{mc}(\mathcal{F}),a\in\mathbf{R}\right\} \right)\right\} \\
 & =\mathfrak{ac}\left(\left\{ ag:g\in\mathfrak{R}\left(\mathcal{L}|_{\mathfrak{mc}(\mathcal{F})}\right),a\in\mathbf{R}\right\} \right)\\
 & =\mathfrak{ag}\left[\mathfrak{R}\left(\mathcal{L}|_{\mathfrak{mc}(\mathcal{F})}\right)\right]\subset\mathfrak{ca}(\mathcal{F})
\end{aligned}
\label{eq:L(ag_F)_in_ca(F)}
\end{equation}
by (\ref{eq:Base_for_OP}) and the linear space properties of $\mathcal{L}$,
$\mathfrak{D}(\mathcal{L})$ and $\mathfrak{ca}(\mathcal{F})$. Thus,
\begin{equation}
\widehat{\mathcal{L}}_{0}\circeq\left\{ (\widehat{f},\widehat{\mathcal{L}f}):f\in\mathfrak{ag}(\mathcal{F})\right\} \label{eq:Lhat_0}
\end{equation}
and
\begin{equation}
\widehat{\mathcal{L}}_{1}\circeq\left\{ (\widehat{f},\widehat{g}):(f,g)\in\mathcal{L}\cap\left(\mathfrak{ca}(\mathcal{F})\times\mathfrak{ca}(\mathcal{F})\right)\right\} \label{eq:Lhat_1}
\end{equation}
satisfy (\ref{eq:ag_(F)_in_D(L)}), (\ref{eq:L(ag_F)_in_ca(F)}) and
Proposition \ref{prop:RepFun_Basic} (d, f) (with $d=k=1$).\end{proof}

The $\widehat{\mathcal{L}}_{0}$ and $\widehat{\mathcal{L}}_{1}$
above are defined as two (possibly) different replicas of $\mathcal{L}$.
\begin{defn}
\label{def:RepOP}Let $E$ be a topological space, $(E_{0},\mathcal{F};\widehat{E},\widehat{\mathcal{F}})$
be a base over $E$ and $\mathcal{L}$ be a linear operator on $C_{b}(E;\mathbf{R})$.
\begin{itemize}
\item $(E_{0},\mathcal{F};\widehat{E},\widehat{\mathcal{F}})$ is said to
be\textbf{ a base for $\mathcal{L}$} if (\ref{eq:Base_for_OP}) holds.
\item When $(E_{0},\mathcal{F};\widehat{E},\widehat{\mathcal{F}})$ is a
base for $\mathcal{L}$, the operator $\widehat{\mathcal{L}}_{0}$
in (\ref{eq:Lhat_0}) and the operator $\widehat{\mathcal{L}}_{1}$
in (\ref{eq:Lhat_1}) are called the \textbf{core replica} and the\textbf{
extended replica} \textbf{of} $\mathcal{L}$ respectively.
\end{itemize}
\end{defn}
Hereafter, we use the following notations for brevity if no confusion
is caused.
\begin{notation}
\label{notation:OP}Let $E$ be a topological space, $(E_{0},\mathcal{F};\widehat{E},\widehat{\mathcal{F}})$
be a base over $E$, $\mathcal{L}$ be a linear operator on $C_{b}(E;\mathbf{R})$
and $\beta\in\mathbf{R}$.
\begin{itemize}
\item We define $\widetilde{f}\circeq f|_{E_{0}^{d}}$ for each $d,k\in\mathbf{N}$
and $f\in M(E^{d};\mathbf{R}^{k})$. Moreover, $\widetilde{\mathcal{F}}\circeq\mathcal{F}|_{E_{0}}=\{\widetilde{f}:f\in\mathcal{F}\}$.
\item When $(E_{0},\mathcal{F};\widehat{E},\widehat{\mathcal{F}})$ is a
base for $\mathcal{L}$, we define
\begin{equation}
\begin{aligned} & \mathcal{L}_{0}\circeq\mathcal{L}|_{\mathfrak{ag}(\mathcal{F})},\\
 & \mathcal{L}_{1}\circeq\mathcal{L}\cap\left(\mathfrak{ca}(\mathcal{F})\times\mathfrak{ca}(\mathcal{F})\right),\\
 & \widetilde{\mathcal{L}}_{i}\circeq\left\{ (\widetilde{f},\widetilde{g}):(f,g)\in\mathcal{L}_{i}\right\} ,\;\forall i=0,1.
\end{aligned}
\label{eq:L0_L1}
\end{equation}

\item The operator $\beta-\mathcal{L}$ is defined by
\begin{equation}
(\beta-\mathcal{L})f\circeq\beta f-\mathcal{L}f,\;\forall f\in\mathfrak{D}(\mathcal{L}).\label{eq:Lambda-Ltilte0}
\end{equation}
Similar notations apply to $\mathcal{L}_{i}$, $\widetilde{\mathcal{L}}_{i}$
and $\widehat{\mathcal{L}}_{i}$ for each $i=0,1$.
\end{itemize}
\end{notation}
The domain and range of the operators above have the following properties.
\begin{prop}
\label{prop:OP_D_R}Let $E$ be a topological space, $\mathcal{L}$
be a linear operator on $C_{b}(E;\mathbf{R})$ and $(E_{0},\mathcal{F};\widehat{E},\widehat{\mathcal{F}})$
be a base over $E$ for $\mathcal{L}$. Then:

\renewcommand{\labelenumi}{(\alph{enumi})}
\begin{enumerate}
\item The linear operators $\mathcal{L}_{0}$ and $\mathcal{L}_{1}$ satisfy
$\mathcal{L}_{0}=\mathcal{L}_{1}|_{\mathfrak{ag}(\mathcal{F})}$ and
\begin{equation}
\mathfrak{R}(\mathcal{L}_{i})\subset\mathfrak{cl}\left(\mathfrak{D}(\mathcal{L}_{i})\right)=\mathfrak{ca}(\mathcal{F}),\;\forall i=0,1.\label{eq:L0_L1_DD}
\end{equation}

\item The linear operators $\widetilde{\mathcal{L}}_{0}$ and $\widetilde{\mathcal{L}}_{1}$
satisfy $\widetilde{\mathcal{L}}_{0}=\widetilde{\mathcal{L}}_{1}|_{\mathfrak{ag}(\widetilde{\mathcal{F}})}$
and
\begin{equation}
\mathfrak{R}(\widetilde{\mathcal{L}}_{i})\subset\mathfrak{cl}\left(\mathfrak{D}(\widetilde{\mathcal{L}}_{i})\right)=\mathfrak{ca}(\widetilde{\mathcal{F}}),\;\forall i=0,1.\label{eq:Ltilte0_Ltilte1_DD}
\end{equation}

\item The linear operators $\widehat{\mathcal{L}}_{0}$ and $\widehat{\mathcal{L}}_{1}$
satisfy $\widehat{\mathcal{L}}_{0}=\widehat{\mathcal{L}}_{1}|_{\mathfrak{ag}(\widehat{\mathcal{F}})}$
and
\begin{equation}
\mathfrak{R}(\widehat{\mathcal{L}}_{i})\subset\mathfrak{cl}\left(\mathfrak{D}(\widehat{\mathcal{L}}_{i})\right)=\mathfrak{ca}(\widehat{\mathcal{F}})=C(\widehat{E};\mathbf{R}),\;\forall i=0,1.\label{eq:Lhat0_Lhat1_DD}
\end{equation}

\item If $\mathcal{L}1=0$, then $\widehat{\mathcal{L}}_{0}1=0$ and $\widehat{\mathcal{L}}_{1}1=0$.
\end{enumerate}
\end{prop}
\begin{proof}
(a) The linearities of $\mathcal{L}_{0}$ and $\mathcal{L}_{1}$ follow
by that of $\mathcal{L}$. It follows by (\ref{eq:L0_L1}) and (\ref{eq:L(ag_F)_in_ca(F)})
that
\begin{equation}
\mathfrak{R}(\mathcal{L}_{0})=\mathfrak{R}\left(\mathcal{L}|_{\mathfrak{ag}(\mathcal{F})}\right)\subset\mathfrak{ca}(\mathcal{F})=\mathfrak{cl}\left(\mathfrak{ag}(\mathcal{F})\right)=\mathfrak{cl}\left(\mathfrak{D}(\mathcal{L}_{0})\right).\label{eq:L0_DD}
\end{equation}
It follows by (\ref{eq:L0_L1}) and (\ref{eq:L0_DD}) that
\begin{equation}
\mathfrak{R}(\mathcal{L}_{1})\subset\mathfrak{ca}(\mathcal{F})=\mathfrak{cl}\left(\mathfrak{D}(\mathcal{L}_{0})\right)\subset\mathfrak{cl}\left(\mathfrak{D}(\mathcal{L}_{1})\right)\subset\mathfrak{ca}(\mathcal{F}).\label{eq:L1_DD}
\end{equation}
Now, (a) follows by (\ref{eq:L0_DD}) and (\ref{eq:L1_DD}).

(b) The linearity of $\widetilde{\mathcal{L}}_{0}$ (resp. $\widetilde{\mathcal{L}}_{1}$)
follows by that of $\mathcal{L}_{0}$ (resp. $\mathcal{L}_{1}$).
It follows by (\ref{eq:L0_L1}), (\ref{eq:L(ag_F)_in_ca(F)}) and
properties of uniform convergence that 
\begin{equation}
\begin{aligned}\mathfrak{R}\left(\widetilde{\mathcal{L}}_{0}\right) & =\left\{ \widetilde{\mathcal{L}f}:f\in\mathfrak{ag}(\mathcal{F})\right\} =\left.\mathfrak{R}\left(\mathcal{L}|_{\mathfrak{ag}(\mathcal{F})}\right)\right|_{E_{0}}\\
 & \subset\left.\mathfrak{ca}(\mathcal{F})\right|_{E_{0}}\subset\mathfrak{ca}(\widetilde{\mathcal{F}})=\mathfrak{cl}\left(\mathfrak{ag}(\widetilde{\mathcal{F}})\right)=\mathfrak{cl}\left(\mathfrak{D}(\widetilde{\mathcal{L}}_{0})\right).
\end{aligned}
\label{eq:Ltilte0_DD}
\end{equation}
It follows by (\ref{eq:L0_L1}) and (\ref{eq:Ltilte0_DD}) that
\begin{equation}
\begin{aligned}\mathfrak{R}(\widetilde{\mathcal{L}}_{1}) & =\left\{ \widetilde{g}:g\in\mathfrak{R}(\mathcal{L}_{1})\right\} \subset\left.\mathfrak{ca}(\mathcal{F})\right|_{E_{0}}\\
 & \subset\mathfrak{ca}(\widetilde{\mathcal{F}})=\mathfrak{cl}\left(\mathfrak{D}(\widetilde{\mathcal{L}}_{0})\right)\subset\mathfrak{cl}\left(\mathfrak{D}(\widetilde{\mathcal{L}}_{1})\right)\subset\mathfrak{ca}(\widetilde{\mathcal{F}}).
\end{aligned}
\label{eq:Ltilte1_DD}
\end{equation}
Now, (b) follows by (\ref{eq:Ltilte0_DD}) and (\ref{eq:Ltilte1_DD}).

(c) follows by Proposition \ref{prop:RepOP_Exist} and Corollary \ref{cor:Base_Fun_Dense}
(with $d=1$).

(d) $1=\widehat{1}\in\widehat{\mathcal{F}}\subset\mathfrak{D}(\widehat{\mathcal{L}}_{0})$
and $\widehat{\mathcal{L}}_{0}1=\widehat{\mathcal{L}1}=\widehat{0}=0$
by the fact $1\in\mathcal{F}$, (\ref{eq:Lhat0_Domain}), (\ref{eq:Lhat0_Rep})
and the denseness of $E_{0}$ in $\widehat{E}$.\end{proof}

\subsection{\label{sec:Generator_RepOp}Markov-generator properties}

The replica operators $\widehat{\mathcal{L}}_{0}$ and $\widehat{\mathcal{L}}_{1}$
may inherit or refine properties of the original operator $\mathcal{L}$
by the compactness of $\widehat{E}$ and the association of the orignal
and replica functions. \cite{DK20b,DK20c,DK21a} use the following
Markov-generator-type properties of replica operators:
\begin{claim}
\label{Claim:OP_Property}$\,$

\begin{enumerate}
[label=\textbf{P\arabic*}, labelsep=0.5pc]

\item\label{enu:P1}$\widehat{\mathcal{L}}_{1}=\mathfrak{cl}(\widehat{\mathcal{L}}_{0})$.

\item\label{enu:P2}$\widehat{\mathcal{L}}_{0}$ satisfies the positive
maximum principle.

\item\label{enu:P3}$\widehat{\mathcal{L}}_{1}$ satisfies the positive
maximum principle.

\item\label{enu:P4}$\widehat{\mathcal{L}}_{0}$ is a strong generator
on $C(\widehat{E};\mathbf{R})$.

\item\label{enu:P5}$\widehat{\mathcal{L}}_{0}$ is a Feller generator
on $C(\widehat{E};\mathbf{R})$.

\end{enumerate}

\end{claim}
The following lemma gives a sufficient condition for \ref{enu:P1}
and explains why we call $\widehat{\mathcal{L}}_{0}$ and $\widehat{\mathcal{L}}_{1}$
the core and extended replica of $\mathcal{L}$.
\begin{lem}
\label{lem:OP_Dense}Let $E$ be a topological space, $\mathcal{L}$
be a linear operator on $C_{b}(E;\mathbf{R})$ and $(E_{0},\mathcal{F};\widehat{E},\widehat{\mathcal{F}})$
be a base over $E$ for $\mathcal{L}$. Then:

\renewcommand{\labelenumi}{(\alph{enumi})}
\begin{enumerate}
\item If \ref{enu:P2} (resp. \ref{enu:P3}) holds, then $\widehat{\mathcal{L}}_{0}$
(resp. $\widehat{\mathcal{L}}_{1}$) is dissipative. 
\item If \ref{enu:P4} holds, and $\widehat{\mathcal{L}}_{1}$ is dissipative
(especially \ref{enu:P3} holds), then \ref{enu:P1} holds.
\end{enumerate}
\end{lem}
\begin{proof}
$\widehat{E}$ is a compact Polish space by Lemma \ref{lem:Base}
(b, c). $C_{0}(\widehat{E};\mathbf{R})=C(\widehat{E};\mathbf{R})$
by (\ref{eq:Q-Algebra_RepTF}) (with $d=1$). Then, the result follows
by Proposition \ref{prop:OP_D_R} (c) and \cite[\S 4.2, Lemma 2.1 and \S 1.4, Proposition 4.1]{EK86}.\end{proof}

\begin{rem}
\label{rem:Core}$\widehat{\mathcal{L}}_{1}$ is a linear superspace
of $\widehat{\mathcal{L}}_{0}$, so we call $\widehat{\mathcal{L}}_{1}$
the extended replica. ``core replica'' comes from the fact that
$\mathfrak{D}(\widehat{\mathcal{L}}_{0})$ is a \textit{core} (see
\cite[\S 1.3, p.17]{EK86}) of $\widehat{\mathcal{L}}_{1}$ in the
setting of Lemma \ref{lem:OP_Dense} (b).
\end{rem}

The following lemma specifies when \ref{enu:P2} - \ref{enu:P5} hold.
\begin{lem}
\label{lem:RepOP_Generator}Let $E$ be a topological space, $\mathcal{L}$
be a linear operator on $C_{b}(E;\mathbf{R})$ and $(E_{0},\mathcal{F};\widehat{E},\widehat{\mathcal{F}})$
be a base over $E$ for $\mathcal{L}$. Then:

\renewcommand{\labelenumi}{(\alph{enumi})}
\begin{enumerate}
\item $\widehat{\mathcal{L}}_{0}$ (resp. $\widehat{\mathcal{L}}_{1}$)
is dissipative if and only if $\widetilde{\mathcal{L}}_{0}$ (resp.
$\widetilde{\mathcal{L}}_{1}$) is dissipative.
\item \ref{enu:P2} (resp. \ref{enu:P3}) holds if and only if for any $\epsilon\in(0,\infty)$
and $f\in\mathfrak{ag}(\mathcal{F})$ (resp. $f\in\mathfrak{D}(\mathcal{L}_{1})$),
there exists an $n_{\epsilon}^{f}\in\mathbf{N}$ such that%
\footnote{$f^{+}$ was defined in \S \ref{sub:Fun} and $\widetilde{f^{+}}$
was defined in Notation \ref{notation:OP}.%
}
\begin{equation}
\sup_{x\in E_{0}}\left[\mathcal{L}f(x)-n_{\epsilon}^{f}\left(\left\Vert \widetilde{f^{+}}\right\Vert _{\infty}-f(x)\right)\right]\leq\epsilon.\label{eq:Ltilte_PMP}
\end{equation}

\item \ref{enu:P4} holds if and only if (1) $\widetilde{\mathcal{L}}_{0}$
is dissipative, and (2) there exists a $\beta\in(0,\infty)$ such
that%
\footnote{The operator $\lambda-\widetilde{\mathcal{L}}$ was defined in Notation
\ref{notation:OP}.%
}
\begin{equation}
\widetilde{\mathcal{F}}\subset\mathfrak{ca}\left(\left\{ (\beta-\widetilde{\mathcal{L}})\widetilde{f}:f\in\mathfrak{mc}(\mathcal{F})\right\} \right).\label{eq:Ltilte_DR}
\end{equation}

\item \ref{enu:P5} holds if and only if (1) $\mathcal{L}1=0$, (2) there
exists a $\beta\in(0,\infty)$ such that (\ref{eq:Ltilte_DR}) holds,
and (3) for any $\epsilon\in(0,\infty)$ and $f\in\mathfrak{ag}(\mathcal{F})$,
there exists an $n_{\epsilon}^{f}\in\mathbf{N}$ such that (\ref{eq:Ltilte_PMP})
holds.
\end{enumerate}
\end{lem}
\begin{proof}
(a) We have by Fact \ref{fact:f+_Rep} (a) (with $d=k=1$) that
\begin{equation}
\Vert\widehat{f}\Vert_{\infty}=\Vert\widetilde{f}\Vert_{\infty},\;\forall f\in\mathfrak{ca}(\mathcal{F}).\label{eq:fhat_ftilte_Sup_Equal}
\end{equation}
Letting $g\circeq(\beta-\mathcal{L})f$, we have by (\ref{eq:Lhat1_Rep})
and Proposition \ref{prop:RepFun_Basic} (d) that
\begin{equation}
\left\Vert \beta\tilde{f}-\widetilde{\mathcal{L}f}\right\Vert _{\infty}=\left\Vert \widetilde{g}\right\Vert _{\infty}=\left\Vert \widehat{g}\right\Vert _{\infty}=\left\Vert \beta\widehat{f}-\widehat{\mathcal{L}f}\right\Vert _{\infty},\;\forall f\in\mathfrak{D}(\mathcal{L}_{1}),\beta\in(0,\infty).\label{eq:Resolvent_Sup_Equal}
\end{equation}
Now, (a) follows by (\ref{eq:fhat_ftilte_Sup_Equal}), (\ref{eq:Resolvent_Sup_Equal})
and the fact that $\widetilde{\mathcal{L}}_{0}\subset\widetilde{\mathcal{L}}_{1}$
and $\widehat{\mathcal{L}}_{0}\subset\widehat{\mathcal{L}}_{1}$.

(b - Sufficiency) We state the proof for $\widehat{\mathcal{L}}_{0}$.
$\widehat{\mathcal{L}}_{1}$ follows similarly. Let $\epsilon\in(0,\infty)$,
$f\in\mathfrak{ag}(\mathcal{F})$ and $n_{\epsilon}^{f}\in\mathbf{N}$
satisfy (\ref{eq:Ltilte_PMP}). We have by Fact \ref{fact:f+_Rep}
(with $d=k=1$) that 
\begin{equation}
\widehat{f^{+}}=\widehat{f}^{+}\mbox{ and }\Vert\widehat{f}^{+}\Vert_{\infty}=\Vert\widetilde{f^{+}}\Vert_{\infty}.\label{eq:f+hat=00003Dfhat+}
\end{equation}
Then, there exists an $x_{0}\in\widehat{E}$ such that
\begin{equation}
\left\Vert \widetilde{f^{+}}\right\Vert _{\infty}=\left\Vert \widehat{f}^{+}\right\Vert _{\infty}=\hat{f}(x_{0})\label{eq:fhat+_max}
\end{equation}
by (\ref{eq:f+hat=00003Dfhat+}), the compactness of $\widehat{E}$,
the continuity of $\widehat{f}^{+}$ and \cite[Theorem 27.4]{M00}.
$\widehat{E}$ is metrizable by Lemma \ref{lem:Base} (c), so there
exist $\{x_{k}\}_{k\in\mathbf{N}}\subset E_{0}$ such that
\begin{equation}
x_{k}\longrightarrow x_{0}\mbox{ as }k\uparrow\infty\mbox{ in }\widehat{E}\label{eq:xk_x0}
\end{equation}
by $E_{0}$'s denseness in $\widehat{E}$ and Fact \ref{fact:First_Countable}
(with $E=\widehat{E}$ and $A=E_{0}$). It follows that
\begin{equation}
\begin{aligned}\widehat{\mathcal{L}f}(x_{k}) & =\mathcal{L}f(x_{k})\leq n_{\epsilon}^{f}\left(\left\Vert \widetilde{f^{+}}\right\Vert _{\infty}-f(x_{k})\right)+\epsilon\\
 & =n_{\epsilon}^{f}\left(\left\Vert \widehat{f}^{+}\right\Vert _{\infty}-\hat{f}(x_{k})\right)+\epsilon\\
 & =n_{\epsilon}^{f}\left(\widehat{f}(x_{0})-\hat{f}(x_{k})\right)+\epsilon,\;\forall k\in\mathbf{N}
\end{aligned}
\label{eq:Check_Lhat_PMP_1}
\end{equation}
by (\ref{eq:Lhat0_Rep}), (\ref{eq:Ltilte_PMP}) and (\ref{eq:fhat+_max}).
Hence, we have that
\begin{equation}
\begin{aligned}\widehat{\mathcal{L}f}(x_{0}) & =\lim_{k\rightarrow\infty}\widehat{\mathcal{L}f}(x_{k})\\
 & \leq\lim_{k\rightarrow\infty}n_{\epsilon}^{f}\left(\widehat{f}(x_{0})-\hat{f}(x_{k})\right)+\epsilon\\
 & =n_{\epsilon}^{f}\left(\widehat{f}(x_{0})-\lim_{k\rightarrow\infty}\hat{f}(x_{k})\right)+\epsilon=\epsilon
\end{aligned}
\label{eq:Check_Lhat_PMP_2}
\end{equation}
by the continuities of $\widehat{f}$ and $\widehat{\mathcal{L}f}$,
(\ref{eq:xk_x0}), (\ref{eq:Check_Lhat_PMP_1}) and the independence
of $n_{\epsilon}^{f}$ and $\{x_{k}\}_{k\in\mathbf{N}}$. Letting
$\epsilon\downarrow0$ in (\ref{eq:Check_Lhat_PMP_2}), we get $\widehat{\mathcal{L}}\widehat{f}(x_{0})\leq0$.

(b - Necessity) Fix $\epsilon\in(0,\infty)$ and $f\in\mathfrak{ca}(\mathcal{F})$.
Then, $f$ satisfies (\ref{eq:f+hat=00003Dfhat+}) by Fact \ref{fact:f+_Rep}
(with $d=k=1$). $\widehat{\mathcal{L}}_{0}$ is dissipative by \ref{enu:P2}
and Lemma \ref{lem:OP_Dense} (a). By (\ref{eq:Lhat0_Rep}), (\ref{eq:fhat_ftilte_Sup_Equal}),
the compactness of $\widehat{E}$, \cite[\S 4.5, Lemmas 5.3]{EK86}
and \cite[\S 4.5, (5.15)]{EK86}, there exist an $n_{\epsilon}^{f}\in\mathbf{N}$
and a \textit{positive contraction} (see \cite[\S 1.1, p.6 and \S 4.2, p.165]{EK86})
$\mathcal{S}_{1/n_{\epsilon}^{f}}$ on $C(\widehat{E};\mathbf{R})$
such that
\begin{equation}
\begin{aligned}\mathcal{L}f(x) & =\widehat{\mathcal{L}f}(x)\leq n_{\epsilon}^{f}\left(\mathcal{S}_{1/n_{\epsilon}^{f}}\widehat{f}(x)-\widehat{f}(x)\right)+\epsilon\\
 & \leq n_{\epsilon}^{f}g(x)-\mathcal{S}_{1/n_{\epsilon}^{f}}g(x)+\epsilon\leq n_{\epsilon}^{f}\left(\left\Vert \widetilde{f^{+}}\right\Vert _{\infty}-f(x)\right)+\epsilon
\end{aligned}
\label{eq:Check_Lhat0_PMP_1}
\end{equation}
for all $x\in E_{0}\subset\widehat{E}$, where $g\circeq\Vert\widehat{f}^{+}\Vert_{\infty}-f\geq0$
satisfies $\mathcal{S}_{1/n_{\epsilon}^{f}}g(x)\geq0$ by the positiveness
of $\mathcal{S}_{1/n_{\epsilon}^{f}}$.

(c) It follows by (\ref{eq:Ltilte0_Ltilte1_DD}), (\ref{eq:Ltilte_DR}),
the linearity of $\widetilde{\mathcal{L}}_{0}$ and (\ref{eq:L0_L1})
that
\begin{equation}
\begin{aligned}\mathfrak{cl}\left[\mathfrak{D}(\widetilde{\mathcal{L}}_{0})\right] & =\mathfrak{ca}(\widetilde{\mathcal{F}})\\
 & \subset\mathfrak{ca}\left[\mathfrak{R}\left(\beta-\widetilde{\mathcal{L}}_{0}|_{\mathfrak{mc}(\widetilde{\mathcal{F}})}\right)\right]\\
 & \subset\mathfrak{ca}\left[\mathfrak{R}\left(\beta-\widetilde{\mathcal{L}}_{0}\right)\right]=\mathfrak{cl}\left[\mathfrak{ag}\left(\mathfrak{R}\left(\beta-\widetilde{\mathcal{L}}_{0}\right)\right)\right]\\
 & =\mathfrak{cl}\left[\left\{ \beta\widetilde{f}-\widetilde{\mathcal{L}f}:f\in\mathfrak{ag}\left(\mathfrak{D}(\mathcal{L}_{0})\right)\right\} \right]\\
 & =\mathfrak{cl}\left[\left\{ \beta\tilde{f}-\widetilde{\mathcal{L}f}:f\in\mathfrak{D}(\mathcal{L}_{0})\right\} \right]\\
 & =\mathfrak{cl}\left[\mathfrak{R}\left(\beta-\widetilde{\mathcal{L}}_{0}\right)\right]\subset\mathfrak{ca}(\widetilde{\mathcal{F}}),
\end{aligned}
\label{eq:Check_Ltilte0_DD_DR}
\end{equation}
thus proving the equivalence between (\ref{eq:Ltilte_DR}) and
\begin{equation}
\mathfrak{cl}\left[\mathfrak{D}(\widetilde{\mathcal{L}}_{0})\right]=\mathfrak{cl}\left[\mathfrak{R}\left(\beta-\widetilde{\mathcal{L}}_{0}\right)\right].\label{eq:Ltilte0_DD_DR}
\end{equation}
Next, we find by (\ref{eq:L0_L1}), (\ref{eq:Lhat0_Domain}), Proposition
\ref{prop:RepFun_Basic} (d, e) and (\ref{eq:Lhat0_Rep}) that
\begin{equation}
\mathfrak{D}(\widetilde{\mathcal{L}}_{0})=\mathfrak{ag}(\widetilde{\mathcal{F}})=\left.\mathfrak{ag}(\widehat{\mathcal{F}})\right|_{E_{0}}=\left.\mathfrak{D}(\widehat{\mathcal{L}}_{0})\right|_{E_{0}}\label{eq:Ltilte0_Domain_Extend_Ehat}
\end{equation}
and
\begin{equation}
\mathfrak{R}\left(\beta-\widetilde{\mathcal{L}}_{0}\right)=\left.\mathfrak{R}\left(\beta-\widehat{\mathcal{L}}_{0}\right)\right|_{E_{0}}.\label{eq:Ltilte0_Resolvent_Extend_Ehat}
\end{equation}
Then, (\ref{eq:Ltilte0_DD_DR}) is equivalent to
\begin{equation}
\mathfrak{cl}\left[\mathfrak{D}(\widehat{\mathcal{L}}_{0})\right]=C(\widehat{E};\mathbf{R})=\mathfrak{cl}\left[\mathfrak{R}\left(\beta-\widehat{\mathcal{L}}_{0}\right)\right]\label{eq:Lhat0_DD_DR}
\end{equation}
by (\ref{eq:Ltilte0_Domain_Extend_Ehat}), (\ref{eq:Ltilte0_Resolvent_Extend_Ehat}),
the denseness of $E_{0}$ in $\widehat{E}$, properties of uniform
convergence and (\ref{eq:Lhat0_Lhat1_DD}). So far, we have shown
the equivalence of (\ref{eq:Ltilte_DR}), (\ref{eq:Ltilte0_DD_DR})
and (\ref{eq:Lhat0_DD_DR}). Now, (c) follows by (a) and the Lumer-Phillips
Theorem (see \cite[\S IX.8]{Y80}).

(d) follows by (b), Proposition \ref{prop:OP_D_R} (d), \cite[\S 4.2, Theorem 2.2]{EK86}
and the equivalence between (\ref{eq:Ltilte_DR}) and (\ref{eq:Lhat0_DD_DR}).\end{proof}

\subsection{\label{sub:OP_Regularity}Conditions on operator}

We consider typical conditions on $\mathcal{L}$ under which: a) One
can construct bases in either of the following two forms, and b) The
associated replica operators $\widehat{\mathcal{L}}_{0}$ and $\widehat{\mathcal{L}}_{1}$
satisfy some of \ref{enu:P1} - \ref{enu:P5}.
\begin{claim}
\label{Claim:OP_Base}Suppose $\mathcal{D}_{0}\subset\mathfrak{D}(\mathcal{L})$
must be replicated. Then:

\begin{enumerate}
[label=\textbf{P\arabic*}, labelsep=0.5pc]
\setcounter{enumi}{5}

\item\label{enu:P6}There exists a base $(E,\mathcal{F};\widehat{E},\widehat{\mathcal{F}})$
over $E$ for $\mathcal{L}$ with $\mathcal{D}_{0}\subset\mathcal{F}=\mathfrak{ag}_{\mathbf{Q}}(\mathcal{F})$.

\item\label{enu:P7}There exists a base $(E_{0},\mathcal{F};\widehat{E},\widehat{\mathcal{F}})$
over $E$ for $\mathcal{L}$ such that $A\subset E_{0}$, $E_{0}\in\mathscr{K}_{\sigma}^{\mathbf{m}}(E)$
and $\mathcal{D}_{0}\subset\mathcal{F}=\mathfrak{ag}_{\mathbf{Q}}(\mathcal{F})$.

\end{enumerate}

\end{claim}
\begin{rem}
\label{rem:D0_OP}$\mathcal{D}_{0}$ also appeared in Lemma \ref{lem:Base_Construction}.
In martingale problems, $\mathcal{D}_{0}$ can be a rich collection
that approximates both the domain and range of $\mathcal{L}$ (see
e.g. \cite[\S 3.1]{DK20b}). In filtering, $\mathcal{D}_{0}$ may
approximate the given sensor function (see e.g. \cite{DK20c}). $A$
above stands for a desired subset of $E$ to be contained in $E_{0}$.
In martingale problems, the set $A$ could be a support of the given
initial distribution (see e.g. \cite[\S A.1]{DK20b}).
\end{rem}

Our regularity conditions consist of four types. The first type is
about the denseness of the domain and range of $\mathcal{L}$.
\begin{condition}
[\textbf{Denseness}]The operator $\mathcal{L}$ satisfies:

\begin{enumerate}
[label=\textbf{D\arabic*}, labelsep=0.5pc]

\item\label{enu:D1}$\mathfrak{R}(\mathcal{L})\subset\mathfrak{cl}(\mathfrak{D}(\mathcal{L}))$.

\item\label{enu:D2}$\mathfrak{D}(\mathcal{L})\subset\mathfrak{cl}(\mathfrak{R}(\beta-\mathcal{L}))$
for some $\beta\in(0,\infty)$.

\end{enumerate}

\end{condition}
\begin{rem}
\label{rem:D1_D2}\ref{enu:D1} and \ref{enu:D2} are true for any
strong generator $\mathcal{L}$ by the Lumer-Phillips Theorem.
\end{rem}

The second type is about the point-separability of the domain of $\mathcal{L}$.
\begin{condition}
[\textbf{Separability}]The operator $\mathcal{L}$ satisfies:

\begin{enumerate}
[label=\textbf{S\arabic*}, labelsep=0.5pc]

\item\label{enu:S1}$\mathfrak{D}(\mathcal{L})$ contains the constant
function $1$ and separates points on $E$.

\item\label{enu:S2}$\mathfrak{D}(\mathcal{L})$ contains the constant
function $1$ and $E$ is a $\mathfrak{D}(\mathcal{L})$-baseable
space.

\end{enumerate}

\end{condition}
\ref{enu:S2} is stronger than \ref{enu:S1}. They are mild requirements
about the richness of $\mathfrak{D}(\mathcal{L})$ by our discussions
in \S \ref{sec:Baseable_Space}, and the next example further illustrates
their generality.
\begin{example}
\label{exp:OP_Regularity_1}$\,$

\renewcommand{\labelenumi}{(\Roman{enumi}) }
\begin{enumerate}
\item For martingale problems and nonlinear filternig problems, a common
setting (see \cite{S78}, \cite{BBK95} and \cite{BBK00}) is that
$E$ is a metrizable Lusin (especially Polish) space and the domain
$\mathfrak{D}(\mathcal{L})$ of $\mathcal{L}$ contains $1$ and separates
points on $E$. Then, $E$ is a second-countable space by Proposition
\ref{prop:Var_Polish} (d) and Proposition \ref{prop:Metrizable}
(c). So, $\mathcal{L}$ satisfies \ref{enu:S2} by Proposition \ref{prop:Hered_Lindelof_Baseable}
(with $A=E$ and $\mathcal{D}=\mathfrak{D}(\mathcal{L})$).
\item Another classical setting for martingale problems is that $E$ is
a locally compact separable metric space with one-point compactification
$E\cup\{\Delta\}$, $\mathcal{L}$ is a linear operator on $C_{0}(E;\mathbf{R})$
and its domain $\mathfrak{D}(\mathcal{L})$ is uniformly dense in
$C_{0}(E;\mathbf{R})$ (see \cite[Chapter 4]{EK86} and \cite{KO88}).
In this case, one can simply extend $\mathcal{L}$ to a linear operator
$\mathcal{L}^{*}$ on $C_{b}(E;\mathbf{R})$ by defining $\mathfrak{D}(\mathcal{L}^{*})$
as the linear span of $\mathfrak{D}(\mathcal{L})\cup\{1\}$ and defining
\begin{equation}
\mathcal{L}^{*}(af+b)\circeq a\mathcal{L}f,\;\forall f\in\mathfrak{D}(\mathcal{L}),a,b\in\mathbf{R}.\label{eq:OP_Trivial_Extend}
\end{equation}
By Proposition \ref{prop:LC_Polish} (a, b, d, e) (with $\mathcal{D}=\mathfrak{D}(\mathcal{L})$),
there exists a countable $\mathcal{F}\subset C_{b}(E;\mathbf{R})$
such that $\mathcal{F}\backslash\{1\}\subset\mathfrak{D}(\mathcal{L})$
strongly separates points on $E$ and $C_{0}(E;\mathbf{R})\subset\mathfrak{ca}(\mathcal{F})$.
Thus, $\mathcal{L}^{*}$ satisfies \ref{enu:S2} by Proposition \ref{prop:Fun_Sep_1}
(a). Moreover, this $\mathcal{L}^{*}$ satisfies \ref{enu:D1}.
\item Suppose that $E$ is a (possibly non-metrizable) Tychonoff space and
$\mathcal{L}$ is a strong generator on $C_{b}(E;\mathbf{R})$. Without
loss of generality, one can consider $1\in\mathfrak{D}(\mathcal{L})$.
Otherwise, we extend $\mathcal{L}$ to $\mathcal{L}^{*}$ as in (II).
$\mathfrak{D}(\mathcal{L})$ is uniformly dense in $C_{b}(E;\mathbf{R})$
by the Lumer-Phillips Theorem, so \ref{enu:D1} holds. $C_{b}(E;\mathbf{R})$
separates points on $E$ by Proposition \ref{prop:CR} (a, c), so
\ref{enu:S1} holds by Corollary \ref{cor:SSP_Dense}.
\end{enumerate}
\end{example}
The third type of regularity conditions includes several analogues
of dissipativeness and positive maximum principle.
\begin{condition}
[\textbf{Generator}]The operator $\mathcal{L}$ satisfies:

\begin{enumerate}
[label=\textbf{G\arabic*}, labelsep=0.5pc]

\item\label{enu:G1}$\mathcal{L}$ is dissipative.

\item\label{enu:G2}For any $\epsilon,\beta\in(0,\infty)$ and $x\in E$,
there exist $\{K_{n,\epsilon}^{x,\beta}\}_{n\in\mathbf{N}}\subset\mathscr{K}^{\mathbf{m}}(E)$
(independent of $f$) such that each $f\in\mathfrak{D}(\mathcal{L})$
satisfies
\begin{equation}
\begin{aligned} & \beta\left|f(x)\right|-\left\Vert \left.(\beta f-\mathcal{L}f)\right|_{K_{n,\epsilon}^{x,\beta}}\right\Vert _{\infty}\\
 & \qquad\leq\left(\beta\Vert f\Vert_{\infty}+\Vert\mathcal{L}f\Vert_{\infty}+1\right)\epsilon,\;\forall f\in\mathfrak{D}(\mathcal{L})
\end{aligned}
\label{eq:G2}
\end{equation}
for some $n=n_{\epsilon}^{f,x,\beta}\in\mathbf{N}$.

\item\label{enu:G3}For any $f\in\mathfrak{D}(\mathcal{L})$,
\begin{equation}
\lim_{n\rightarrow\infty}\sup_{x\in E}\left[\mathcal{L}f(x)-n\left(\Vert f^{+}\Vert_{\infty}-f(x)\right)\right]\leq0.\label{eq:G3}
\end{equation}

\item\label{enu:G4}For any $\epsilon\in(0,\infty)$, $x\in E$ and
$f\in\mathfrak{D}(\mathcal{L})$, there exist $\{K_{n,\epsilon}^{x}\}_{n\in\mathbf{N}}\subset\mathscr{K}^{\mathbf{m}}(E)$
independent of $f$ and an $n_{\epsilon}^{f}\in\mathbf{N}$ independent
of $x$ such that
\begin{equation}
\mathcal{L}f(x)-n_{\epsilon}^{f}\left(\left\Vert \left.f^{+}\right|_{K_{n_{\epsilon}^{f},\epsilon}^{x}}\right\Vert _{\infty}-f(x)\right)\leq\epsilon.\label{eq:G4}
\end{equation}

\end{enumerate}

\end{condition}
The next example shows that \ref{enu:G2}, \ref{enu:G3} and \ref{enu:G4}
are not unnatural.
\begin{example}
\label{exp:OP_Regularity_2}Let $E$, $\mathcal{L}$ and $\mathcal{L}^{*}$
be as in Example \ref{exp:OP_Regularity_1} (II) and $\epsilon\in(0,\infty)$.

\renewcommand{\labelenumi}{(\Roman{enumi}) }
\begin{enumerate}
\item If $\mathcal{L}$ satisfies positive maximum principle, then $\mathcal{L}^{*}$
does also. Consequently, \ref{enu:G3} is satisfied by both $\mathcal{L}$
and $\mathcal{L}^{*}$ by an argument similar to the proof of Lemma
\ref{lem:RepOP_Generator} (b - ``only if'').
\item When $\mathcal{L}$ is a Feller generator, the Feller semigroup $\{\mathcal{S}_{t}\}_{t\geq0}$
generated by $\mathfrak{cl}(\mathcal{L})$ on $C_{0}(E;\mathbf{R})$
is often given by a \textit{transition function} (see \cite[\S 4.1, p.156]{EK86})
$\kappa:\mathbf{R}^{+}\times E\times\mathscr{B}(E)\rightarrow[0,1]$.
In the remainder of the example, we fix $x\in E$.

\begin{itemize}
\item Fix $\beta\in(0,\infty)$ and let $\mathbf{Q}=\{q_{n}\}_{n\in\mathbf{N}}$.
$E$ is a Polish space by Proposition \ref{prop:LC_Polish} (c, d).
So for each $n\in\mathbf{N}$, $\kappa(q_{n},x,\cdot)$ is tight in
$E$ by Ulam's Theorem (Theorem \ref{thm:Ulam_m-Tight}), i.e. there
exist $\{K_{n,\epsilon}^{x,\beta}\}_{n\in\mathbf{N}}\in\mathscr{K}(E)$
such that
\begin{equation}
\kappa\left(q_{n},x,E\backslash K_{n,\epsilon}^{x,\beta}\right)\leq\epsilon,\;\forall n\in\mathbf{N}.\label{eq:Check_G2_Example_1}
\end{equation}
One finds by \cite[\S 1.2, (2.1) and (2.6)]{EK86}, change of variable
and Jensen's Inequality that
\begin{equation}
\begin{aligned}\beta\left|f(x)\right| & =\left|\int_{0}^{\infty}\beta e^{-\beta t}\mathcal{S}_{t}(\beta-\mathcal{L})f(x)dt\right|\\
 & \leq\int_{0}^{1}\left|\mathcal{S}_{-\frac{\ln u}{\beta}}(\beta-\mathcal{L})f(x)\right|du,\;\forall f\in\mathfrak{D}(\mathcal{L}).
\end{aligned}
\label{eq:SemiGroup_Change_Variable}
\end{equation}
Then, there exist $\{t^{f,x,\beta}\}_{f\in\mathfrak{D}(\mathcal{L})}\subset(0,1)$
such that 
\begin{equation}
\beta\left|f(x)\right|\leq\left|\mathcal{S}_{t^{f,x,\beta}}(\beta-\mathcal{L})f(x)\right|,\;\forall f\in\mathfrak{D}(\mathcal{L})\label{eq:MeanValueThm}
\end{equation}
by (\ref{eq:SemiGroup_Change_Variable}), Mean-Value Theorem and Jensen's
Inequality. For each fixed $f\in\mathfrak{D}(\mathcal{L})$, there
exists an $n=n_{\epsilon}^{f,x,\beta}\in\mathbf{N}$ such that
\begin{equation}
\left\Vert \mathcal{S}_{q_{n}}(\beta-\mathcal{L})f-\mathcal{S}_{t^{f,x,\beta}}(\beta-\mathcal{L})f\right\Vert <\epsilon\label{eq:Check_G2_Example_2}
\end{equation}
 by the strong continuity of $\{\mathcal{S}_{t}\}_{t\geq0}$. From
(\ref{eq:MeanValueThm}), (\ref{eq:Check_G2_Example_1}) and (\ref{eq:Check_G2_Example_2})
it follows that
\begin{equation}
\begin{aligned}\beta\left|f(x)\right| & \leq\int_{E}\left|(\beta-\mathcal{L})f(y)\right|\kappa\left(q_{n},x,dy\right)+\epsilon\\
 & \leq\int_{K_{n,\epsilon}^{x,\beta}}\left|(\beta f-\mathcal{L}f)(y)\right|\kappa\left(q_{n},x,dy\right)\\
 & +\left(\Vert\beta f\Vert_{\infty}+\Vert\mathcal{L}f\Vert_{\infty}\right)\epsilon+\epsilon\\
 & \leq\left\Vert \left.(\beta f-\mathcal{L}f)\right|_{K_{n,\epsilon}^{x,\beta}}\right\Vert _{\infty}+\left(\Vert\beta f\Vert_{\infty}+\Vert\mathcal{L}f\Vert_{\infty}+1\right)\epsilon.
\end{aligned}
\label{eq:Check_G2_Example_3}
\end{equation}
Thus, $\mathcal{L}$ satisfies \ref{enu:G2} as $\{K_{n,\epsilon}^{x,\beta}\}_{n\in\mathbf{N}}$
does not involve any $f$.
\item For each $n\in\mathbf{N}$, the tightness of $\kappa(n^{-1},x,\cdot)$
implies a $K_{n,\epsilon}^{x}\in\mathscr{K}(E)$ satisfying
\begin{equation}
\kappa\left(n^{-1},x,E\backslash K_{n,\epsilon}^{x}\right)\leq\frac{\epsilon}{2n^{2}}.\label{eq:Check_G4_Example_1}
\end{equation}
Meanwhile, we fix $f\in\mathfrak{D}(\mathcal{L})$. $\mathfrak{cl}(\mathcal{L})$
is the infitesesimal generator of $\{\mathcal{S}_{t}\}_{t\geq0}$,
so there is a sufficiently large $n=n_{\epsilon}^{f}\in\mathbf{N}$
such that
\begin{equation}
n\geq\Vert f\Vert_{\infty}\label{eq:Choose_n_f_epsilon_large}
\end{equation}
and
\begin{equation}
\sup_{z\in E}\left|\mathcal{L}f(z)-n\left[\int_{E}f(y)\kappa\left(n^{-1},z,dy\right)-f(z)\right]\right|\leq\frac{\epsilon}{2}.\label{eq:Check_G4_Example_2}
\end{equation}
The sequence of compact sets $\{K_{n,\epsilon}^{x}\}_{n\in\mathbf{N}}$
is determined by $x$ and the transition function $\kappa$, which
is independent of $f$. The convergence rate $n_{0}=n_{\epsilon}^{f}$
is an intrinsic parameter of $f$ and is unrelated to $x$. From (\ref{eq:Check_G4_Example_1}),
(\ref{eq:Choose_n_f_epsilon_large}) and (\ref{eq:Check_G4_Example_2})
it follows that
\begin{equation}
\begin{aligned}\mathcal{L}f(x) & \leq n_{0}\left[\int_{K_{n_{0},\epsilon}^{x}}f^{+}(y)\kappa\left(n_{0}^{-1},x,dy\right)+\frac{\Vert f\Vert_{\infty}n_{0}^{2}\epsilon}{2}-f(x)\right]+\frac{\epsilon}{2}\\
 & \leq n_{0}\left(\left\Vert \left.f^{+}\right|_{K_{n_{0},\epsilon}^{x}}\right\Vert _{\infty}-f(x)\right)+\epsilon.
\end{aligned}
\label{eq:Check_G4_Example_3}
\end{equation}
Thus, $\mathcal{L}$ satisfies \ref{enu:G4}.
\end{itemize}
\end{enumerate}
\end{example}
The fourth condition type is a equivalent to assuming $\mathfrak{D}(\mathcal{L})$
is closed under multiplication since $\mathcal{L}$ and $\mathfrak{D}(\mathcal{L})$
are linear spaces.
\begin{condition}
[\textbf{DA}]\label{con:DA}The domain of the operator $\mathcal{L}$
is a subalgebra of $C_{b}(E;\mathbf{R})$.
\end{condition}

\subsection{\label{sub:RepOP_Base}Existence of Markov-generator-type replica
operator}

Now, we give four constructions of Markov-generator-type replica operators
under the aforementioned conditions. The first two assume \ref{enu:S2}
and construct bases satisfying \ref{enu:P6}.
\begin{lem}
\label{lem:OP_Base_1}Let $E$ be a topological space, $\mathcal{L}$
be a linear operator on $C_{b}(E;\mathbf{R})$ and $\mathcal{D}_{0}\subset\mathfrak{D}(\mathcal{L})$
be countable. Then, \ref{enu:S2}, \ref{enu:D1}, \ref{enu:G3} plus
\nameref{con:DA} imply \ref{enu:P6}, \ref{enu:P2} plus \ref{enu:P3}.
\end{lem}
\begin{proof}
We use induction to construct the $\mathcal{F}$ in \ref{enu:P6}.
By \ref{enu:S2}, there exists a countable $\mathcal{D}\subset\mathfrak{D}(\mathcal{L})$
that separates on $E$. For $k=0$,
\begin{equation}
\mathcal{F}_{0}\circeq\left(\mathcal{D}_{0}\cup\mathcal{D}\cup\{1\}\right)\subset\mathfrak{D}(\mathcal{L})\label{eq:Define_F0_for_OP}
\end{equation}
is countable, contains $1$ and separates points on $E$. For $k\in\mathbf{N}$,
we assume $\mathcal{F}_{0}\subset\mathcal{F}_{k-1}\subset\mathfrak{D}(\mathcal{L})$
and find by \nameref{con:DA} that
\begin{equation}
\mathfrak{ag}_{\mathbf{Q}}(\mathcal{F}_{k-1})\subset\mathfrak{D}(\mathcal{L}).\label{eq:mc(F_k-1)_in_D(L)}
\end{equation}
Then, we define
\begin{equation}
\mathcal{F}_{k}\circeq\bigcup_{f\in\mathfrak{ag}_{\mathbf{Q}}(\mathcal{F}_{k-1}),q\in\mathbf{N}}\left\{ f,g_{q}^{f,k}\right\} \subset\mathfrak{D}(\mathcal{L}),\label{eq:Define_F_k-1_for_OP_1}
\end{equation}
where each $g_{q}^{f,k}\in\mathfrak{D}(\mathcal{L})$ is chosen by
\ref{enu:D1} to satisfy 
\begin{equation}
\left\Vert g_{q}^{f,k}-\mathcal{L}f\right\Vert _{\infty}\leq2^{-q}.\label{eq:RangeApprox}
\end{equation}
It follows immediately that
\begin{equation}
\mathfrak{ag}_{\mathbf{Q}}(\mathcal{F}_{k-1})\subset\mathcal{F}_{k}\label{eq:Inductive_DD_1}
\end{equation}
and
\begin{equation}
\mathfrak{R}\left(\mathcal{L}|_{\mathfrak{ag}_{\mathbf{Q}}(\mathcal{F}_{k-1})}\right)\subset\mathfrak{cl}\left(\mathcal{F}_{k}\right).\label{eq:Inductive_DD_2}
\end{equation}
Based on the construction above%
\footnote{$\mathbf{N}_{0}=\{0\}\cup\mathbf{N}$ denotes the non-negative integers.%
},
\begin{equation}
\mathcal{F}\circeq\bigcup_{k\in\mathbf{N}_{0}}\mathcal{F}_{k}\label{eq:Define_F_for_OP}
\end{equation}
satisfies
\begin{equation}
\mathcal{D}_{0}\cup\{1\}\cup\mathcal{D}=\mathcal{F}_{0}\subset\mathcal{F}=\mathfrak{ag}_{\mathbf{Q}}(\mathcal{F})\subset\mathfrak{D}(\mathcal{L})\label{eq:Check_D0_in_F_and_SP}
\end{equation}
and (\ref{eq:Base_for_OP}). So, $\mathcal{F}$ separates points on
$E$ as $\mathcal{D}$ does. Now, \ref{enu:P6} follows by Lemma \ref{lem:Base_Construction}
(b) (with $E_{0}=E$ and $\mathcal{D}=\mathcal{F}$) and (\ref{eq:Check_D0_in_F_and_SP}).
\ref{enu:P2} and \ref{enu:P3} follow by \ref{enu:G3} and Lemma
\ref{lem:RepOP_Generator} (b) (with $E_{0}=E$, $\widetilde{\mathcal{L}}_{0}=\mathcal{L}_{0}$
and $\widetilde{\mathcal{L}}_{1}=\mathcal{L}_{1}$).\end{proof}

\begin{lem}
\label{lem:OP_Base_2}Let $E$ be a topological space, $\mathcal{L}$
be a linear operator on $C_{b}(E;\mathbf{R})$ such that \ref{enu:S2},
\ref{enu:D1}, \ref{enu:D2} and \nameref{con:DA} hold, and $\mathcal{D}_{0}\subset\mathfrak{D}(\mathcal{L})$
be countable. Then:

\renewcommand{\labelenumi}{(\alph{enumi})}
\begin{enumerate}
\item If \ref{enu:G1} holds, then \ref{enu:P6}, \ref{enu:P1} and \ref{enu:P4}
hold.
\item If $\mathcal{L}1=0$ and \ref{enu:G3} holds, then \ref{enu:P6},
\ref{enu:P1}, \ref{enu:P3} and \ref{enu:P5} hold.
\end{enumerate}
\end{lem}
\begin{proof}
(a) We choose $\beta\in(0,\infty)$ so \ref{enu:D2} is satisifed.
For each $f\in\mathfrak{D}(\mathcal{L})$, we choose each $g_{q}^{f,k}\in\mathfrak{D}(\mathcal{L})$
by \ref{enu:D1} so (\ref{eq:RangeApprox}) is satisfied and choose
each $h_{q}^{f,k}\in\mathfrak{D}(\mathcal{L})$ by\ref{enu:D2} so
\begin{equation}
\left\Vert (\beta-\mathcal{L})h_{q}^{f,k}-f\right\Vert _{\infty}<2^{-q}.\label{eq:DomainApprox}
\end{equation}
Now, we follow the proof of Lemma \ref{lem:OP_Base_1} to establish
the $\mathcal{F}=\bigcup_{k}\mathcal{F}_{k}$ in \ref{enu:P6}, where
\begin{equation}
\mathcal{F}_{k}\circeq\bigcup_{f\in\mathfrak{ag}_{\mathbf{Q}}(\mathcal{F}_{k-1}),q\in\mathbf{N}}\left\{ f,g_{q}^{f,k},h_{q}^{f,k}\right\} \subset\mathfrak{D}(\mathcal{L}).\label{eq:Define_F_k_for_OP2}
\end{equation}
Consequently, $\mathcal{F}_{k}$ defined in (\ref{eq:Define_F_k_for_OP2})
satisfies not only (\ref{eq:Inductive_DD_1}) and (\ref{eq:Inductive_DD_2})
but also
\begin{equation}
\mathfrak{ag}_{\mathbf{Q}}(\mathcal{F}_{k-1})\subset\mathfrak{cl}\left[\mathfrak{R}\left(\beta-\mathcal{L}|_{\mathcal{F}_{k}}\right)\right].\label{eq:Inductive_DR}
\end{equation}
$\mathcal{F}$ defined in (\ref{eq:Define_F_for_OP}) not only satisfies
(\ref{eq:Check_D0_in_F_and_SP}) and (\ref{eq:Base_for_OP}) but also
satisfies 
\begin{equation}
\mathcal{F}\subset\mathfrak{cl}\left[\mathfrak{R}\left(\beta-\mathcal{L}|_{\mathcal{F}}\right)\right]\subset\mathfrak{ca}\left[\mathfrak{R}\left(\beta-\mathcal{L}|_{\mathfrak{mc}(\mathcal{F})}\right)\right].\label{eq:Check_Ltilte_DR}
\end{equation}
Now, \ref{enu:P6} follows by Lemma \ref{lem:Base_Construction} (b)
(with $E_{0}=E$ and $\mathcal{D}=\mathcal{F}$) and (\ref{eq:Check_D0_in_F_and_SP}).
Both $\mathcal{L}_{0}$ and $\mathcal{L}_{1}$ are dissipative by
\ref{enu:G1}, so \ref{enu:P4} and the dissipativeness of $\widehat{\mathcal{L}}_{1}$
follow by (\ref{eq:Check_Ltilte_DR}) and Lemma \ref{lem:RepOP_Generator}
(a, c) (with $E_{0}=E$, $\widetilde{\mathcal{L}}_{0}=\mathcal{L}_{0}$
and $\widetilde{\mathcal{L}}_{1}=\mathcal{L}_{1}$). Moreover, \ref{enu:P1}
follows by Lemma \ref{lem:OP_Dense} (b).

(b) Let $\mathcal{F}$ be constructed and \ref{enu:P6} verified as
in (a). \ref{enu:P5} and \ref{enu:P3} follow by \ref{enu:G3}, (\ref{eq:Check_Ltilte_DR})
and Lemma \ref{lem:RepOP_Generator} (b, d) (with $E_{0}=E$, $\widetilde{\mathcal{L}}_{0}=\mathcal{L}_{0}$
and $\widetilde{\mathcal{L}}_{1}=\mathcal{L}_{1}$). \ref{enu:P1}
follows by Lemma \ref{lem:OP_Dense} (b).\end{proof}

The third method alternatively assumes \ref{enu:S1} and uses an arbitrary
$A\in\mathscr{K}_{\sigma}^{\mathbf{m}}(E)$ and metrizable compacts
provided by \ref{enu:G4} to construct a base satisfying \ref{enu:P7}.
\begin{lem}
\label{lem:OP_Base_3}Let $E$ be a topological space, $A\in\mathscr{K}_{\sigma}^{\mathbf{m}}(E)$,
$\mathcal{L}$ be a linear operator on $C_{b}(E;\mathbf{R})$ satisfy
\ref{enu:S1}, \ref{enu:G4} and \nameref{con:DA} hold, and $\mathcal{D}_{0}\subset\mathfrak{D}(\mathcal{L})$
be countable. Then:

\renewcommand{\labelenumi}{(\alph{enumi})}
\begin{enumerate}
\item If \ref{enu:D1} holds, then \ref{enu:P7}, \ref{enu:P2} and \ref{enu:P3}
hold.
\item If $\mathcal{L}1=0$, and if \ref{enu:D1} and \ref{enu:D2} hold,
then \ref{enu:P7}, \ref{enu:P1}, \ref{enu:P3} and \ref{enu:P5}
hold.
\end{enumerate}
\end{lem}
\begin{proof}
(a) We construct the $(E_{0},\mathcal{F})$ in \ref{enu:P7} by induction.
For $k=0$, we define
\begin{equation}
A_{0}\circeq A\in\mathscr{K}_{\sigma}^{\mathbf{m}}(E)\label{eq:Define_A_0_for_OP3}
\end{equation}
and
\begin{equation}
\mathcal{F}_{0}\circeq\left(\mathcal{D}_{0}\cup\{1\}\right)\subset\mathfrak{D}(\mathcal{L}).\label{eq:Define_F_0_for_OP3}
\end{equation}
For $k\in\mathbf{N}$, we assume $A_{k-1}\in\mathscr{K}_{\sigma}^{\mathbf{m}}(E)$
and $\mathcal{F}_{0}\subset\mathcal{F}_{k-1}\subset\mathfrak{D}(\mathcal{L})$.
$(A_{k-1},\mathscr{O}_{E}(A_{k-1}))$ is a separable space by Proposition
\ref{prop:Sigma_MC} (a, b) and Proposition \ref{prop:Var_Polish}
(d), so it has a countable dense subset $\{x_{j}^{k-1}\}_{j\in\mathbf{N}}$.
For each $i\in\mathbf{N}$ and $x\in A_{k-1}$, one finds by \ref{enu:G4}
a sequence of metrizable compact sets $\{K_{n,i}^{x}\}_{n\in\mathbf{N}}\subset\mathscr{K}^{\mathbf{m}}(E)$
such that each $f\in\mathfrak{D}(\mathcal{L})$ satisfies
\begin{equation}
\mathcal{L}f(x)-n\left(\left\Vert \left.f^{+}\right|_{K_{n,i}^{x}}\right\Vert _{\infty}-f(x)\right)\leq2^{-i}\label{eq:Compact_G4}
\end{equation}
for some sufficiently large $n=n_{i}^{f}\in\mathbf{N}$ independent
of $x$. Then, we redefine
\begin{equation}
A_{k}\circeq A_{k-1}\cup\bigcup_{i,j,n\in\mathbf{N}}K_{n,i}^{x_{j}^{k-1}}\in\mathscr{K}_{\sigma}^{\mathbf{m}}(E).\label{eq:Define_A_k_for_OP3}
\end{equation}

By \ref{enu:S1} and Proposition \ref{prop:Sigma_MC} (b, e) (with
$A=A_{k}$ and $\mathcal{D}=\mathfrak{D}(\mathcal{L})$), there exists
a countable $\mathcal{J}_{k}\subset\mathfrak{D}(\mathcal{L})$ that
separates points on $A_{k}$. \nameref{con:DA} implies
\begin{equation}
\mathfrak{ag}_{\mathbf{Q}}\left(\mathcal{F}_{k-1}\cup\mathcal{J}_{k}\right)\subset\mathfrak{D}(\mathcal{L}).\label{eq:mc(F_k-1_and_J_k)_in_D(L)}
\end{equation}
Then, we define
\begin{equation}
\mathcal{F}_{k}\circeq\bigcup_{f\in\mathfrak{ag}_{\mathbf{Q}}\left(\mathcal{F}_{k-1}\cup\mathcal{J}_{k}\right),q\in\mathbf{N}}\left\{ f,g_{q}^{f,k}\right\} \subset\mathfrak{D}(\mathcal{L}),\label{eq:Define_F_k_for_OP3_1}
\end{equation}
where each $g_{q}^{f,k}$ is chosen by \ref{enu:D1} to satisfy (\ref{eq:RangeApprox}).

By the construction above, $\{x_{j}^{k-1}\}_{j,k\in\mathbf{N}}$ is
a countable dense subset of
\begin{equation}
E_{0}\circeq\bigcup_{k\in\mathbf{N}_{0}}A_{k}\in\mathscr{K}_{\sigma}^{\mathbf{m}}(E)\label{eq:Define_E0_for_OP3}
\end{equation}
under the subspace topology $\mathscr{O}_{E}(E_{0})$. $E$ is a Hausdorff
space by \ref{enu:S1} and Proposition \ref{prop:Fun_Sep_1} (e) (with
$A=E$ and $\mathcal{D}=\mathfrak{D}(\mathcal{L})$), so $E_{0}\in\mathscr{B}(E)$
by Proposition \ref{prop:Compact} (a).

$\mathcal{F}$ defined by (\ref{eq:Define_F_for_OP}) satisfies (\ref{eq:Check_D0_in_F_and_SP})
and (\ref{eq:Base_for_OP}). $\mathcal{F}$ contains $\mathcal{J}_{k}$
and separates points on $A_{k}$ for all $k\in\mathbf{N}$. $\{A_{k}\}_{k\in\mathbf{N}}$
are nested%
\footnote{The terminology ``nested'' was explained in Fact \ref{fact:D-Baseable_Union}.%
} by (\ref{eq:Define_A_k_for_OP3}), so $\mathcal{F}$ separates points
on $E_{0}$ by Fact \ref{fact:SP_Nested_Union}. Fixing $i\in\mathbf{N}$
and $f\in\mathfrak{D}(\mathcal{L})$, we have by (\ref{eq:Compact_G4})
that
\begin{equation}
\begin{aligned} & \mathcal{L}f(x_{j}^{k-1})-n\left(\left\Vert \widetilde{f^{+}}\right\Vert _{\infty}-f(x_{j}^{k-1})\right)\\
 & \leq\mathcal{L}f((x_{j}^{k-1})-n\left(\left\Vert \left.f^{+}\right|_{K_{n,i}^{x_{j}^{k-1}}}\right\Vert _{\infty}-f(x_{j}^{k-1})\right)\leq2^{-i}
\end{aligned}
\label{eq:Check_PMP_Ltilte_D(L)_1}
\end{equation}
for all $j,k\in\mathbf{N}$ and a sufficiently large $n=n_{i}^{f}\in\mathbf{N}$
\textit{independent of any} $x_{j}^{k-1}$. Therefore, it follows
that
\begin{equation}
\mathcal{L}f(x)-n\left(\left\Vert \widetilde{f^{+}}\right\Vert _{\infty}-f(x)\right)\leq2^{-i},\;\forall x\in E_{0}\label{eq:Check_PMP_Ltilte_D(L)_2}
\end{equation}
by the denseness of $\{x_{j}^{k-1}\}_{j,k\in\mathbf{N}}$ in $(E_{0},\mathscr{O}_{E}(E_{0}))$
and the continuities of $f$ and $\mathcal{L}f$.

Now, \ref{enu:P7} follows by Lemma \ref{lem:Base_Construction} (b)
(with $\mathcal{D}=\mathcal{F}$). \ref{enu:P2} and \ref{enu:P3}
follow by (\ref{eq:Check_PMP_Ltilte_D(L)_2}) and Lemma \ref{lem:RepOP_Generator}
(b).

(b) We follow the proof of (a) to establish $(E_{0},\mathcal{F})$
except for reconstructing
\begin{equation}
\mathcal{F}_{k}\circeq\bigcup_{f\in\mathfrak{ag}_{\mathbf{Q}}(\mathcal{F}_{k-1}\cup\mathcal{J}_{k}),q\in\mathbf{N}}\left\{ f,g_{q}^{f,k},h_{q}^{f,k}\right\} \subset\mathfrak{D}(\mathcal{L}).\label{eq:Define_F_k_for_OP3_2}
\end{equation}
Here, we choose each $g_{q}^{f,k}\in\mathfrak{D}(\mathcal{L})$ by
\ref{enu:D1} to satisfy (\ref{eq:RangeApprox}). We find a constant
$\beta\in(0,\infty)$ and choose each $h_{q}^{f,k}\in\mathfrak{D}(\mathcal{L})$
by \ref{enu:D2} to satisfy (\ref{eq:DomainApprox}). Consequently,
$\mathcal{F}$ not only satisfies (\ref{eq:Check_D0_in_F_and_SP})
and (\ref{eq:Base_for_OP}) but also satisfies (\ref{eq:Check_Ltilte_DR}).

Now, \ref{enu:P7} follows by the same argument of (a). (\ref{eq:Ltilte_DR})
follows by (\ref{eq:Check_Ltilte_DR}) and properties of uniform convergence.
Hence, \ref{enu:P5} follows by (\ref{eq:Check_PMP_Ltilte_D(L)_2}),
Proposition \ref{prop:OP_D_R} (d) and Lemma \ref{lem:RepOP_Generator}
(b, d). Moreover, \ref{enu:P1} follows by Lemma \ref{lem:OP_Dense}
(b).\end{proof}

Our fourth construction is similar to the third with \ref{enu:G4}
replaced by \ref{enu:G2}.
\begin{lem}
\label{lem:OP_Base_4}Let $E$ be a topological space, $A\in\mathscr{K}_{\sigma}^{\mathbf{m}}(E)$,
$\mathcal{L}$ be a linear operator on $C_{b}(E;\mathbf{R})$ and
$\mathcal{D}_{0}\subset\mathfrak{D}(\mathcal{L})$ be countable. If
\ref{enu:S1}, \ref{enu:D1}, \ref{enu:D2}, \ref{enu:G2} and \nameref{con:DA}
hold, then \ref{enu:P7}, \ref{enu:P1} and \ref{enu:P4} hold.
\end{lem}
\begin{proof}
The proof of Lemma \ref{lem:OP_Base_3} (a) establishes the $(E_{0},\mathcal{F})$
in \ref{enu:P7} except $E_{0}=\bigcup_{k\in\mathbf{N}_{0}}A_{k}\in\mathscr{K}_{\sigma}^{\mathbf{m}}(E)$
is constructed as follows. For $k=0$, we still define $A_{0}$ as
in (\ref{eq:Define_A_0_for_OP3}). For each $i\in\mathbf{N}$, $x\in A_{k-1}$
and $\beta\in\mathbf{Q}^{+}$%
\footnote{$\mathbf{Q}^{+}$ denotes the non-negative rational numbers.%
}, one finds by \ref{enu:G2} a sequence of metrizable compact sets
$\{K_{n,i}^{x,\beta}\}_{n\in\mathbf{N}}\subset\mathscr{K}^{\mathbf{m}}(E)$
such that each $f\in\mathfrak{D}(\mathcal{L})$ satisfies 
\begin{equation}
\beta\left|f(x)\right|-\left\Vert (\beta f-\mathcal{L}f)|_{K_{n,i}^{x,\beta}}\right\Vert _{\infty}\leq\left(\beta\Vert f\Vert_{\infty}+\Vert\mathcal{L}f\Vert_{\infty}+1\right)2^{-i}\label{eq:Compact_G2}
\end{equation}
for some $n=n_{i}^{f,x,\beta}\in\mathbf{N}$.

Now, we take a countable dense subset $\{x_{j}^{k-1}\}_{j\in\mathbf{N}}$
of $(A_{k-1},\mathscr{O}_{E}(A_{k-1}))$ and define
\begin{equation}
A_{k}\circeq A_{k-1}\cup\left(\bigcup_{\beta\in\mathbf{Q}^{+}}\bigcup_{i,j,n\in\mathbf{N}}K_{n,i}^{x_{j}^{k-1},\beta}\right)\in\mathscr{K}_{\sigma}^{\mathbf{m}}(E).\label{eq:Define_A_k_for_OP4}
\end{equation}

By the reconstruction above, $(E_{0},\mathcal{F})$ has almost the
same properties as in the proof of Lemma \ref{lem:OP_Base_3} (b)
except for (\ref{eq:Compact_G4}). Instead, fixing $f\in\mathfrak{D}(\mathcal{L})$,
$\beta\in\mathbf{Q}^{+}$ and $i\in\mathbf{N}$, we have by (\ref{eq:Compact_G2})
and (\ref{eq:Define_A_k_for_OP4}) that 
\begin{equation}
\begin{aligned} & \beta\left|f(x_{j}^{k-1})\right|-\left\Vert \beta\tilde{f}-\widetilde{\mathcal{L}f}\right\Vert _{\infty}\\
 & \leq\beta\left|f(x_{j}^{k-1})\right|-\left\Vert (\beta f-\mathcal{L}f)|_{K_{n,i}^{x_{j}^{k-1},\beta}}\right\Vert _{\infty}\\
 & \leq\left(\beta\Vert f\Vert_{\infty}+\Vert\mathcal{L}f\Vert_{\infty}+1\right)2^{-i}
\end{aligned}
\label{eq:Check_DIS_Ltilte_D(L)_1}
\end{equation}
for all $j,k\in\mathbf{N}$ and some $n=n_{i}^{f,x_{j}^{k-1},\beta}\in\mathbf{N}$.
Therefore, it follows that
\begin{equation}
\beta\Vert\tilde{f}\Vert_{\infty}-\Vert\beta\tilde{f}-\widetilde{\mathcal{L}f}\Vert_{\infty}\leq\left(\beta\Vert f\Vert_{\infty}+\Vert\mathcal{L}f\Vert_{\infty}+1\right)2^{-i}\label{eq:Check_DIS_Ltilte_D(L)_2}
\end{equation}
by the denseness of $\{x_{j}^{k-1}\}_{j,k\in\mathbf{N}}$ in $(E_{0},\mathscr{O}_{E}(E_{0}))$
and the continuities of $f$. Letting $i\uparrow\infty$ in (\ref{eq:Check_DIS_Ltilte_D(L)_2}),
we obtain
\begin{equation}
\beta\Vert\tilde{f}\Vert_{\infty}\leq\Vert\beta\tilde{f}-\widetilde{\mathcal{L}f}\Vert_{\infty},\;\forall f\in\mathfrak{D}(\mathcal{L}).\label{eq:Check_DIS_Ltilte_D(L)_3}
\end{equation}
Next, we let $\beta\in(0,\infty)$, take $\{\beta_{m}\}_{m\in\mathbf{N}}\subset\mathbf{Q}^{+}\cap(0,\infty)$
with $\lim_{m\rightarrow\infty}\beta_{m}=\beta$ and find that
\begin{equation}
\begin{aligned}\beta\Vert\widetilde{f}\Vert_{\infty} & \leq\lim_{m\rightarrow\infty}\beta_{m}\Vert\widetilde{f}\Vert_{\infty}\\
 & \leq\left\Vert \beta\widetilde{f}-\widetilde{\mathcal{L}f}\right\Vert _{\infty}+\lim_{m\rightarrow\infty}\Vert\widetilde{f}\Vert_{\infty}\left|\beta_{m}-\beta\right|=\left\Vert \beta\widetilde{f}-\widetilde{\mathcal{L}f}\right\Vert _{\infty}
\end{aligned}
\label{eq:Check_DIS_Ltilte_D(L)_5}
\end{equation}
by (\ref{eq:Check_DIS_Ltilte_D(L)_3}) (with $\beta=\beta_{m}$),
thus proving the dissipativeness of $\widetilde{\mathcal{L}}_{0}$
and $\widetilde{\mathcal{L}}_{1}$.

Now, \ref{enu:P7} follows by Lemma \ref{lem:Base_Construction} (b)
(with $\mathcal{D}=\mathcal{F}$). (\ref{eq:Check_Ltilte_DR}) holds
as in the proof of Lemma \ref{lem:OP_Base_3} (a) and implies (\ref{eq:Ltilte_DR}).
Then, \ref{enu:P4} and the dissipativeness of $\widehat{\mathcal{L}}_{1}$
follow by Lemma \ref{lem:RepOP_Generator} (a, c). Moreover, \ref{enu:P1}
follows by Lemma \ref{lem:OP_Dense} (b).\end{proof}

Morever, it is worth noting that replica operators will exist under
much less regularity of $\mathcal{L}$ if no Markov-generator-type
property is required.
\begin{prop}
\label{prop:RepOP_Exist_Gen}Let $E$ be a topological space, $A\in\mathscr{K}_{\sigma}^{\mathbf{m}}(E)$,
$\mathcal{D}_{0}\subset\mathfrak{D}(\mathcal{L})$ be countable and
$\mathcal{L}$ be a linear operator on $C_{b}(E;\mathbf{R})$ satisfying
\ref{enu:D1} and \nameref{con:DA}. Then:

\renewcommand{\labelenumi}{(\alph{enumi})}
\begin{enumerate}
\item If \ref{enu:S2} holds, then \ref{enu:P6} holds.
\item If \ref{enu:S1} holds, then \ref{enu:P7} holds.
\end{enumerate}
\end{prop}
\begin{proof}
The construction of the desired bases for (a) and (b) are already
contained in the proofs of Lemma \ref{lem:OP_Base_1} and Lemma \ref{lem:OP_Base_3}
(a) respectively.\end{proof}

\chapter{\label{chap:RepMeas}Weak Convergence and Replication of Measure}

\chaptermark{Weak Convergence and Replica Measure}

Replication has an important role in weak convergence. Our baseable
space and baseable subset results are used in \S \ref{sec:WLP_Uni}
to establish mild conditions for the uniqueness, existence and consistency
of weak limit points on the finite-dimensional Cartesian power $E^{d}$
of a topological space $E$. \S \ref{sub:RepMeas_Def} introduces
the (Borel) replica of possibly non-Borel measure on $E^{d}$. \S
\ref{sub:RepMeas_WC} relates weak convergence of the replica measures
to that of the original ones, which will be a basic tool for our developments
in \ref{enu:Theme2} and \ref{enu:Theme3}. Using replica functions
and measures, we extend two fundamental theorems in probability theory
to non-classical settings. In \S \ref{sub:Riesz_Representation},
we establish a version of the Radon-Riesz Representation Theorem on
a non-locally-compact and even non-Tychonoff space. In \ref{sub:Sko_Rep},
we establish the Skorokhod Representation Theorem under mild conditions.

\section{\label{sec:WLP_Uni}Tightness and weak convergence}

Given a general topological space $E$, tightness and $\mathbf{m}$-tightness
unsurprisingly play key roles in establishing weak convergence on
$E^{d}$. Existence of a tight subsequence usually implies existence
of a weak limit point. Uniqueness requires slightly stronger tightness.
\begin{defn}
\label{def:Seq_Tight}Let $(E,\mathscr{U})$ be a measurable space,
$S$ be a topological space and $\mathscr{A}$ be a $\sigma$-algebra
on $S$.
\begin{itemize}
\item When $S\subset E$, $\Gamma\subset\mathfrak{M}^{+}(E,\mathscr{U})$
is \textbf{sequentially tight} (resp. \textbf{sequentially $\mathbf{m}$-tight})\textbf{
in }$S$ if: (1) $\Gamma$ is an infinite set, and (2) Any infinite
subset of $\Gamma$ admits a subsequence being tight (resp. $\mathbf{m}$-tight)
in $S$.
\item $\Gamma\subset\mathfrak{M}^{+}(S,\mathscr{A})$ is \textbf{sequentially
tight }(resp. \textbf{sequentially $\mathbf{m}$-tight}) \textbf{in}
$A\subset S$ if: (1) $\Gamma$ is an infinite set and $A$ is non-empty,
and (2) Any infinite subset of $\Gamma$ admits a subsequence being
tight (resp. $\mathbf{m}$-tight) in $(A,\mathscr{O}_{S}(A))$.
\item $\Gamma\subset\mathfrak{M}^{+}(S,\mathscr{A})$ is\textbf{ sequentially
tight} (resp. \textbf{$\mathbf{m}$-tight}) if $\Gamma$ is sequentially
tight (resp. $\mathbf{m}$-tight) in $S$.
\end{itemize}
\end{defn}
\begin{note}
\label{note:Seq_Tight_RV}Any type of sequential tightness in Definition
\ref{def:Seq_Tight} is defined for random variables as the corresponding
property of their distributions.
\end{note}
The following sequential version of Proposition \ref{prop:m-Tight_BExt}
relates sequential $\mathbf{m}$-tightness and unique existence of
Borel extension.
\begin{prop}
\label{prop:Seq_m-Tight_BExt}Let $\mathbf{I}$ be a countable index
set, $\{S_{i}\}_{i\in\mathbf{I}}$ be topological spaces, $(S,\mathscr{A})$
be as in (\ref{eq:(S,A)_Prod_Meas_Space}), $\Gamma\subset\mathfrak{M}^{+}(S,\mathscr{A})$
and $A\subset S$. Suppose in addition that $\mathfrak{p}_{i}(A)\in\mathscr{B}(S_{i})$
is a Hausdorff subspace of $S_{i}$ for all $i\in\mathbf{I}$. Then,
$\Gamma$ is sequentially $\mathbf{m}$-tight in $A$ if and only
if there exists a $\Gamma_{0}\in\mathscr{P}_{0}(\Gamma)$%
\footnote{$\mathscr{P}_{0}(\Gamma)$ is the family of all finite subsets of
$\Gamma$.%
} such that $\{\mu^{\prime}=\mathfrak{be}(\mu)\}_{\mu\in\Gamma\backslash\Gamma_{0}}$
is sequentially $\mathbf{m}$-tight in $A$.
\end{prop}
\begin{proof}
Suppose $\Gamma^{\prime}\subset\Gamma$ is an infinite set and $\aleph(\mathfrak{be}(\mu))\neq1$
for each $\mu\in\Gamma^{\prime}$. Given $\Gamma$'s sequential $\mathbf{m}$-tightness,
there exists an $\mathbf{m}$-tight subsequence $\{\mu_{n}\}_{n\in\mathbf{N}}\subset\Gamma^{\prime}$
and, by Proposition \ref{prop:m-Tight_BExt}, $\mathfrak{be}(\mu_{n})$
is a singleton for all $n\in\mathbf{N}$. Contradiction! Hence, there
exists $\{\mu^{\prime}=\mathfrak{be}(\mu)\}_{\mu\in\Gamma\backslash\Gamma_{0}}$
for some $\Gamma_{0}\in\mathscr{P}_{0}(\Gamma)$. Finally, the $\mathbf{m}$-tightness
of $\{\mu_{n}\}_{n\in\mathbf{N}}\backslash\Gamma_{0}$ in $A$ is
equivalent to that of $\{\mu^{\prime}:\mu\in\{\mu_{n}\}_{n\in\mathbf{N}}\backslash\Gamma_{0}\}$
(if any) by Proposition \ref{prop:m-Tight_BExt}.\end{proof}

We then give our conditions for Borel extensions of finite measures
on the product measurable space $(E^{d},\mathscr{B}(E)^{\otimes d})$
to have a unique weak limit point.
\begin{thm}
\label{thm:WLP_Uni}Let $E$ be a topological space, $d\in\mathbf{N}$,
$\Gamma\subset\mathfrak{M}^{+}(E^{d},\mathscr{B}(E)^{\otimes d})$,
$\mathcal{D}\subset C_{b}(E;\mathbf{R})$ and $\mathcal{G}\circeq\mathfrak{mc}[\Pi^{d}(\mathcal{D})]$.
Suppose that:

\renewcommand{\labelenumi}{(\roman{enumi})}
\begin{enumerate}
\item $\Gamma$ is sequentially $\mathbf{m}$-tight.
\item $\{\int_{E^{d}}f(x)\mu(dx)\}_{\mu\in\Gamma}$ has at most one%
\footnote{$\{\int_{E^{d}}f(x)\mu(dx)\}_{\mu\in\Gamma}$ lies in $[-\Vert f\Vert_{\infty}\mu(E^{d}),\Vert f\Vert_{\infty}\mu(E^{d})]$
and has at least one limit point in $\mathbf{R}$ by the Bolzano-Weierstrass
Theorem. So, we are effectively assuming it has a unique limit point.%
} limit point in $\mathbf{R}$ for all $f\in\mathcal{G}\cup\{1\}$.
\item $\mathcal{D}$ separates points on $E$.
\end{enumerate}
Then:

\renewcommand{\labelenumi}{(\alph{enumi})}
\begin{enumerate}
\item $\Gamma^{\prime}\circeq\{\mu^{\prime}=\mathfrak{be}(\mu)\}_{\mu\in\Gamma\backslash\Gamma_{0}}$
exists for some $\Gamma_{0}\in\mathscr{P}_{0}(\Gamma)$ and is sequentially
$\mathbf{m}$-tight.
\item $\Gamma^{\prime}$ has at most one weak limit point in $\mathcal{M}^{+}(E^{d})$.
\item If, in addition, $\{\mu(E^{d})\}_{\mu\in\Gamma}\subset[a,b]$ for
some $0<a<b$, then $\Gamma^{\prime}$ has a unique weak limit point
$\nu$ in $\mathcal{M}^{+}(E^{d})$. Moreover, $\nu$ is an $\mathbf{m}$-tight
measure with total mass%
\footnote{The notion of total mass was specified in \S \ref{sub:Meas}. The
notation ``w-$\lim_{n\rightarrow\infty}\mu_{n}=\nu$'' introduced
in \S \ref{sec:Borel_Measure} means that $\nu$ is the weak limit
of $\{\mu_{n}\}_{n\in\mathbf{N}}$. In other words, it means $\mu_{n}\Rightarrow\nu$
as $n\uparrow\infty$ and $\nu$ is the unique weak limit point of
$\{\mu_{n}\}_{n\in\mathbf{N}}$.%
} in $[a,b]$ and satisfies
\begin{equation}
\mathrm{w}\mbox{-}\lim_{n\rightarrow\infty}\mu_{n}^{\prime}=\nu,\;\forall\{\mu_{n}\}_{n\in\mathbf{N}}\subset\Gamma\backslash\Gamma_{0}.\label{eq:All_Seq_WC_Uni_WLP}
\end{equation}

\end{enumerate}
\end{thm}
\begin{note}
\label{note:Cb(E;R)_SP_Baseable}The condition (iii) above is true
for a wide subclass of Hausdorff spaces which need neither be Tychonoff
nor baseable.
\end{note}

\begin{note}
\label{note:Pi^d(D)_Mb_Cb}Any $\mathcal{D}\subset M_{b}(E;\mathbf{R})$
satisfies
\begin{equation}
\mathfrak{ca}\left[\Pi^{d}(\mathcal{D})\right]\subset M_{b}\left(E^{d},\mathscr{B}(E)^{\otimes d};\mathbf{R}\right)\label{eq:ca(Pi^d(D))_Mb}
\end{equation}
and any $\mathcal{D}\subset C_{b}(E;\mathbf{R})$ satisfies
\begin{equation}
\mathfrak{ca}\left[\Pi^{d}(\mathcal{D})\right]\subset C_{b}\left(E^{d},\mathscr{O}(E)^{d};\mathbf{R}\right)\label{eq:ca(Pi^d(D))_Cb}
\end{equation}
by Proposition \ref{prop:Pi^d_SP} (a) and properties of uniform convergence.
So, $\int_{E^{d}}f(x)\mu(dx)$ is well-defined for all $f\in\mathfrak{ca}[\Pi^{d}(M_{b}(E;\mathbf{R}))]$
and $\mu\in\mathfrak{M}^{+}(E^{d},\mathscr{B}(E)^{\otimes d})$.
\end{note}
Before proving Theorem \ref{thm:WLP_Uni}, we give a Portmanteau-type
lemma for compact sets.
\begin{lem}
\label{lem:Compact_Portmanteau}Let $C(E;\mathbf{R})$ separate points
on topological space $E$. Then:

\renewcommand{\labelenumi}{(\alph{enumi})}
\begin{enumerate}
\item (\ref{eq:Mu_n_WC_Mu_M(E)}) implies $\mu(K)\geq\limsup_{n\rightarrow\infty}\mu_{n}(K)$
for all $K\in\mathscr{K}(E)$.
\item If $\mu$ is a weak limit point of $\Gamma$ in $\mathcal{M}^{+}(E)$,
and if $\Gamma$ is tight (resp. $\mathbf{m}$-tight) in $A\subset E$,
then $\Gamma\cup\{\mu\}$ is tight (resp. $\mathbf{m}$-tight) in
$A$.
\end{enumerate}
\end{lem}
\begin{note}
\label{note:C(E;R)_Cb(E;R)_SP}$C(E;\mathbf{R})$ separating points
(resp. strongly separating points) on $E$ is equivalent to $C_{b}(E;\mathbf{R})$
separating points (resp. strongly separating points) on $E$ (see
Corollary \ref{cor:C(E;R)_Cb(E;R)_SP}).\end{note}
\begin{rem}
\label{rem:Mass_Confine}The classical Portmanteau's Theorem asserts
that the mass of a weakly convergent sequence does not escape any
\textit{closed} set. The mass may no longer be confined by a general
closed subset of the possibly non-Tychonoff space $E$ (see Theorem
\ref{thm:Portamenteau}). Nonetheless, Lemma \ref{lem:Compact_Portmanteau}
confirms that the subclass $\mathscr{K}(E)$ of $\mathscr{C}(E)$
still maintains this property.
\end{rem}
\begin{proof}
[Proof of Lemma \ref{lem:Compact_Portmanteau}](a) $(E,\mathscr{O}_{C(E;\mathbf{R})}(E))$
is a Tychonoff coarsening of $E$ by Proposition \ref{prop:CR} (a,
b). $K\in\mathscr{K}(E,\mathscr{O}_{C(E;\mathbf{R})}(E))\subset\mathscr{C}(E,\mathscr{O}_{C(E;\mathbf{R})}(E))$
by Fact \ref{fact:Compact_Topo_Coarsen} (b) (with $\mathcal{D}=C(E;\mathbf{R})$).
\begin{equation}
\mu_{n}\Longrightarrow\mu\mbox{ as }n\uparrow\infty\mbox{ in }\mathcal{M}^{+}\left(E,\mathscr{O}_{C(E;\mathbf{R})}(E)\right)\label{eq:WC_Coarsen_C(E;R)}
\end{equation}
by Fact \ref{fact:Weak_Topo_Coarsen} (b) (with $\mathscr{U}=\mathscr{O}_{C(E;\mathbf{R})}(E;\mathbf{R})$).
Now, (a) follows by Theorem \ref{thm:Portamenteau} (a, b) (with $E=(E,\mathscr{O}_{C(E;\mathbf{R})}(E))$).

(b) Let $\{\mu_{n}\}_{n\in\mathbf{N}}\subset\Gamma$ satisfy (\ref{eq:Mu_n_WC_Mu_M(E)}).
By tightness (resp. $\mathbf{m}$-tightness) of $\Gamma$ in $A$
and (a), there exist $\{K_{p}\}_{p\in\mathbf{N}}\subset\mathscr{K}(A,\mathscr{O}_{E}(A))$
(resp. $\mathscr{K}^{\mathbf{m}}(A,\mathscr{O}_{E}(A))$) such that
\begin{equation}
\mu(E\backslash K_{p})\leq\liminf_{n\rightarrow\infty}\mu_{n}(E\backslash K_{p})\leq\sup_{\mu\in\Gamma}\mu(E\backslash K_{p})\leq2^{-p},\;\forall p\in\mathbf{N}.\label{eq:Check_WC_Lim_Tight}
\end{equation}
\end{proof}

\begin{proof}
[Proof of Theorem \ref{thm:WLP_Uni}](a) $E$ is a Hausdorff space
by Proposition \ref{prop:Fun_Sep_1} (e) (with $A=E$). By Proposition
\ref{prop:Seq_m-Tight_BExt} (with $\mathbf{I}=\{1,...,d\}$, $S_{i}=E$
and $A=E^{d}$), there exists a $\Gamma_{0}\subset\mathscr{P}_{0}(\Gamma)$
such that (a) holds.

(b) Suppose $\{\mu_{i,n}^{\prime}:i=1,2,n\in\mathbf{N}\}\subset\Gamma^{\prime}$
satisfy
\begin{equation}
\mu_{i,n}^{\prime}\Longrightarrow\mu_{i}\mbox{ as }n\uparrow\infty\mbox{ in }\mathcal{M}^{+}(E^{d}),\;\forall i=1,2.\label{eq:Two_WC_Subseq}
\end{equation}
The sequential $\mathbf{m}$-tightness of $\Gamma^{\prime}$ implies
an $\mathbf{m}$-tight subsequence $\{\mu_{i,n_{k}}^{\prime}\}_{k\in\mathbf{N}}$
for each $i=1,2$. $\mathcal{G}$ separates points on $E^{d}$ by
Proposition \ref{prop:Pi^d_SP} (b), so does $C(E^{d};\mathbf{R})$
by (\ref{eq:ca(Pi^d(D))_Cb}). Then, $\{\mu_{i,n_{k}}^{\prime}\}_{k\in\mathbf{N},i=1,2}\cup\{\mu_{1}\}\cup\{\mu_{2}\}$
is $\mathbf{m}$-tight by Lemma \ref{lem:Compact_Portmanteau} (b)
(with $E=A=E^{d}$ and $\Gamma=\{\mu_{i,n_{k}}^{\prime}\}_{k\in\mathbf{N},i=1,2}$).
It follows that%
\footnote{The notation $f^{*}$ was defined in (\ref{eq:f_star}).%
}
\begin{equation}
\begin{aligned}\int_{E^{d}}f(x)\mu_{1}(dx) & =\lim_{n\rightarrow\infty}f^{*}(\mu_{1,n}^{\prime})\\
 & =\lim_{n\rightarrow\infty}f^{*}(\mu_{2,n}^{\prime})=\int_{E^{d}}f(x)\mu_{2}(dx)\;\forall f\in\mathcal{G}\cup\{1\}
\end{aligned}
\label{eq:Lim_Meas_SP_Func_Test}
\end{equation}
by (\ref{eq:Two_WC_Subseq}), (\ref{eq:ca(Pi^d(D))_Cb}) and the fact
that $\{\int_{E^{d}}f(x)\mu_{i,n}(dx)\}_{n\in\mathbf{N},i=1,2}$ has
at most one limit point in $\mathbf{R}$ for all $f\in\mathcal{G}\cup\{1\}$.
Now, $\mu_{1}=\mu_{2}$ by Lemma \ref{lem:Tight_Meas_Identical} (a)
(with $E=E^{d}$ and $\mathcal{D}=\mathcal{G}$).

(c) $E^{d}$ is a Hausdorff space by Proposition \ref{prop:Separability}
(d). So, $\Gamma^{\prime}$ has a unique weak limit point $\nu$ in
$\mathcal{M}^{+}(E^{d})$ with $\nu(E^{d})\in[a,b]$ by (b) and Lemma
\ref{lem:Seq_Prokhorov} (with $E=E^{d}$ and $\Gamma=\Gamma^{\prime}$).
$\nu$ is $\mathbf{m}$-tight by Lemma \ref{lem:Compact_Portmanteau}
(b) (with $E=A=E^{d}$ and $\Gamma=\Gamma^{\prime}$). Furthermore,
(\ref{eq:All_Seq_WC_Uni_WLP}) follows by the sequential $\mathbf{m}$-tightness
of $\Gamma^{\prime}$, the fact
\begin{equation}
\mu^{\prime}(E^{d})=\mu(E^{d})\in[a,b]\subset(0,\infty),\;\forall\mu\in\Gamma\backslash\Gamma_{0}\label{eq:BExt_Same_Total_Mass}
\end{equation}
and Corollary \ref{cor:Identify_Seq_WC} (with $E=A=E^{d}$, $\mu=\nu$,
$\Gamma=\Gamma^{\prime}$ and $\mu_{n}=\mu_{n}^{\prime}$).\end{proof}

For finite measures on infinite-dimensional Cartesian power of $E$,
Theorem \ref{thm:WLP_Uni} shows unique existence of weak limit point
for their finite-dimensional distributions under suitable conditions.
The following theorem considers the Kolmogorov Extension of these
finite-dimensional weak limit points.
\begin{thm}
\label{thm:WLP_Consistency}Let $E$ be a topological space, $\Gamma\subset\mathfrak{M}^{+}(E^{\mathbf{I}},\mathscr{B}(E)^{\otimes\mathbf{I}})$
and $\mathcal{D}\subset C_{b}(E;\mathbf{R})$. Suppose that:

\renewcommand{\labelenumi}{(\roman{enumi})}
\begin{enumerate}
\item $\{\mu\circ\mathfrak{p}_{\mathbf{I}_{0}}^{-1}\}_{\mu\in\Gamma}$ is
sequentially $\mathbf{m}$-tight for all $\mathbf{I}_{0}\in\mathscr{P}_{0}(\mathbf{I})$.
\item $\{\int_{E^{\mathbf{I}_{0}}}f(x)\mu\circ\mathfrak{p}_{\mathbf{I}_{0}}^{-1}(dx)\}_{\mu\in\Gamma}$
has at most one limit point in $\mathbf{R}$ for all $f\in\mathfrak{mc}[\Pi^{\mathbf{I}_{0}}(\mathcal{D})]\cup\{1\}$%
\footnote{The notation ``$\Pi^{\mathbf{I}_{0}}(\mathcal{D})$'' was defined
in \S \ref{sub:Fun}. $\Pi^{\mathbf{I}_{0}}(\mathcal{D})=\Pi^{d}(\mathcal{D})$
with $d\circeq\aleph(\mathbf{I}_{0})$.%
} for all $\mathbf{I}_{0}\in\mathscr{P}_{0}(\mathbf{I})$.
\item $\mathcal{D}$ separates points on $E$.
\item $\{\mu(E^{\mathbf{I}})\}_{\mu\in\Gamma}\subset[a,b]$ for some $0<a<b$.
\end{enumerate}
Then, there exist a unique $\mu^{\infty}\in\mathfrak{M}^{+}(E^{\mathbf{I}},\mathscr{B}(E)^{\otimes\mathbf{I}})$
such that:

\renewcommand{\labelenumi}{(\alph{enumi})}
\begin{enumerate}
\item The total mass of $\mu^{\infty}$ lies in $[a,b]$.
\item $\mu^{\infty}\circ\mathfrak{p}_{\mathbf{I}_{0}}^{-1}\in\mathcal{M}^{+}(E^{\mathbf{I}_{0}})$
is $\mathbf{m}$-tight for all $\mathbf{I}_{0}\in\mathscr{P}_{0}(\mathbf{I})$.
\item For each $\mathbf{I}_{0}\in\mathscr{P}_{0}(\mathbf{I})$, there is
some $\Gamma_{\mathbf{I}_{0}}^{0}\in\mathscr{P}_{0}(\Gamma)$ such
that $\mu^{\infty}\circ\mathfrak{p}_{\mathbf{I}_{0}}^{-1}$ is the
weak limit%
\footnote{The notion of weak limit was specified in \S \ref{sec:Borel_Measure}.%
} of any subsequence of $\{\mu_{\mathbf{I}_{0}}^{\prime}=\mathfrak{be}(\mu\circ\mathfrak{p}_{\mathbf{I}_{0}}^{-1})\}_{\mu\in\Gamma\backslash\Gamma_{\mathbf{I}_{0}}^{0}}$
in $\mathcal{M}^{+}(E^{\mathbf{I}_{0}})$.
\end{enumerate}
\end{thm}
\begin{proof}
[Proof of Theorem \ref{thm:WLP_Consistency}]We fix $\mathbf{I}_{1},\mathbf{I}_{2}\in\mathscr{P}_{0}(\mathbf{I})$
with $\mathbf{I}_{1}\subset\mathbf{I}_{2}$, let $\mathfrak{p}_{\mathbf{I}_{j}}$
denote the projection from $E^{\mathbf{I}}$ to $E^{\mathbf{I}_{j}}$
for each $j=1,2$, use $\widetilde{\mathfrak{p}}$ to specially denote
the projection from $E^{\mathbf{I}_{2}}$ to $E^{\mathbf{I}_{1}}$
and observe $\mu\circ\mathfrak{p}_{\mathbf{I}_{j}}^{-1}(E^{\mathbf{I}_{j}})=\mu(E^{\mathbf{I}})\in[a,b]$
for all $\mu\in\Gamma$ and $j=1,2$.

By Theorem \ref{thm:WLP_Uni} (a, c) (with $d=\aleph(\mathbf{I}_{j})$
and $\Gamma=\{\mu\circ\mathfrak{p}_{\mathbf{I}_{j}}^{-1}\}_{\mu\in\Gamma}$),
there exist $\mu_{\mathbf{I}_{j}}^{\infty}\in\mathcal{M}^{+}(E^{\mathbf{I}_{j}})$
and $\Gamma_{\mathbf{I}_{j}}^{0}\in\mathscr{P}_{0}(\Gamma)$ for each
$j=1,2$ such that $\mu_{\mathbf{I}_{j}}^{\infty}$ is an $\mathbf{m}$-tight
measure with total mass in $[a,b]$ and is the weak limit of any subsequence
of
\begin{equation}
\Gamma_{\mathbf{I}_{j}}^{\prime}\circeq\left\{ \mu_{\mathbf{I}_{j}}^{\prime}=\mathfrak{be}(\mu\circ\mathfrak{p}_{\mathbf{I}_{j}}^{-1})\right\} _{\mu\in\Gamma\backslash\Gamma_{\mathbf{I}_{j}}^{0}}\subset\mathcal{M}^{+}(E^{\mathbf{I}_{j}}).\label{eq:Gamma_I_i_Prime}
\end{equation}

Suppose that $\{\mu_{n}\}_{n\in\mathbf{N}}\subset\Gamma\backslash(\Gamma_{\mathbf{I}_{1}}^{0}\cup\Gamma_{\mathbf{I}_{2}}^{0})$
satisfies
\begin{equation}
\mbox{w-}\lim_{n\rightarrow\infty}\mu_{n,\mathbf{I}_{j}}^{\prime}=\mu_{\mathbf{I}_{j}}^{\infty},\;\forall j=1,2,\label{eq:E^I1_E^I2_WL}
\end{equation}
where $\mu_{n,\mathbf{I}_{j}}^{\prime}=\mathfrak{be}(\mu_{n}\circ\mathfrak{p}_{\mathbf{I}_{j}}^{-1})\in\Gamma_{\mathbf{I}_{j}}^{\prime}$
for each $n\in\mathbf{N}$ and $j=1,2$. $\widetilde{\mathfrak{p}}\in C(E^{\mathbf{I}_{2}};E^{\mathbf{I}_{1}})$
by Fact \ref{fact:Prod_Map_2} (a). It follows that
\begin{equation}
\mu_{n,\mathbf{I}_{1}}^{\prime}=\mu_{n}\circ\mathfrak{p}_{\mathbf{I}_{2}}^{-1}\circ\widetilde{\mathfrak{p}}^{-1}\Longrightarrow\mu_{\mathbf{I}_{2}}^{\infty}\circ\widetilde{\mathfrak{p}}^{-1}\mbox{ as }n\uparrow\infty\mbox{ in }\mathcal{M}^{+}(E^{\mathbf{I}_{1}})\label{eq:E^I2_WC_Proj}
\end{equation}
by (\ref{eq:E^I1_E^I2_WL}) and the Continuous Mapping Theorem (Theorem
\ref{thm:ContMapTh} (a)). So, (\ref{eq:E^I1_E^I2_WL}) and (\ref{eq:E^I2_WC_Proj})
imply
\begin{equation}
\mu_{\mathbf{I}_{1}}^{\infty}=\mu_{\mathbf{I}_{2}}^{\infty}\circ\mathfrak{p}_{\mathbf{I}_{1}}^{-1}.\label{eq:Consistency}
\end{equation}

By the argument above, $\Gamma$ uniquely determines $\{\mu_{\mathbf{I}_{0}}^{\infty}\in\mathcal{M}^{+}(E^{\mathbf{I}_{0}})\}_{\mathbf{I}_{0}\in\mathscr{P}_{0}(\mathbf{I})}$
such that: (1) $\{\mu_{\mathbf{I}_{0}}^{\infty}\}_{\mathbf{I}_{0}\in\mathscr{P}_{0}(\mathbf{I})}$
satisfies the Kolmogorov consistency and admits a common total mass
$c\in[a,b]$, and (2) each $\mu_{\mathbf{I}_{0}}^{\infty}$ is $\mathbf{m}$-tight
and is the weak limit of any subsequence of $\{\mu_{\mathbf{I}_{0}}^{\prime}\}_{\mu\in\Gamma\backslash\Gamma_{\mathbf{I}_{0}}^{0}}$
for some $\Gamma_{\mathbf{I}_{0}}^{0}\in\mathscr{P}_{0}(\Gamma)$.
$\{E^{\mathbf{I}_{0}}\}_{\mathbf{I}_{0}\in\mathscr{P}_{0}(\mathbf{I})}$
are all Hausdorff spaces by (\ref{eq:ca(Pi^d(D))_Cb}) (with $d=\aleph(\mathbf{I}_{0})$)
and Proposition \ref{prop:Fun_Sep_1} (e) (with $E=A=E^{\mathbf{I}_{0}}$
and $\mathcal{D}=C(E^{\mathbf{I}_{0}};\mathbf{R})$). Now, the unique
existence of $\mu^{\infty}\in\mathfrak{M}^{+}(E^{\mathbf{I}},\mathscr{B}(E)^{\otimes\mathbf{I}})$
satisfying $\mu^{\infty}(E^{\mathbf{I}})=c$ and
\begin{equation}
\mu^{\infty}\circ\mathfrak{p}_{\mathbf{I}_{0}}^{-1}=\mu_{\mathbf{I}_{0}}^{\infty},\;\forall\mathbf{I}_{0}\in\mathscr{P}_{0}(\mathbf{I})\label{eq:Kol_Extension}
\end{equation}
follows by a suitable version of the Kolmogorov's Extension Theorem
(see \cite[Corollary 15.28]{AB06}).\end{proof}

\section{\label{sec:RepMeas}Replica measure}

\subsection{\label{sub:RepMeas_Def}Definition}

Given a base $(E_{0},\mathcal{F};\widehat{E},\widehat{\mathcal{F}})$
over $E$ and $d\in\mathbf{N}$, replicating a finite measure $\mu$
from the product measurable space $(E^{d},\mathscr{B}(E)^{\otimes d})$
onto $\widehat{E}^{d}$ means expanding the concentrated measure $\mu|_{E_{0}^{d}}$
to a Borel replica measure on $\widehat{E}^{d}$.
\begin{defn}
\label{def:RepMeas}Let $E$ be a topological space and $d\in\mathbf{N}$.
The \textbf{replica of $\mu\in\mathfrak{M}^{+}(E^{d},\mathscr{B}(E)^{\otimes d})$}
with respect to a base $(E_{0},\mathcal{F};\widehat{E},\widehat{\mathcal{F}})$
over $E$ is defined by%
\footnote{The concentrated measure ``$\mu|_{A}$'' and expanded measure ``$\nu|^{E}$''
were defined in \S \ref{sub:Meas}.%
}
\begin{equation}
\overline{\mu}\circeq\left.\left(\left.\mu\right|_{E_{0}}\right)\right|^{\widehat{E}^{d}}.\label{eq:RepMeas}
\end{equation}

\end{defn}
The following fact justifies our definition of replica measure.
\begin{fact}
\label{fact:Base_Meas_Basic}Let $E$ be a topological space, $(E_{0},\mathcal{F};\widehat{E},\widehat{\mathcal{F}})$
be a base over $E$ and $d\in\mathbf{N}$. Then, any $\mu\in\mathfrak{M}^{+}(E^{d},\mathscr{B}(E)^{\otimes d})$
satisfies
\begin{equation}
\begin{aligned}\mu|_{E_{0}^{d}} & \in\mathfrak{M}^{+}\left(E_{0}^{d},\left.\mathscr{B}(E)^{\otimes d}\right|_{E_{0}^{d}}\right)=\mathfrak{M}^{+}\left(E_{0}^{d},\mathscr{B}_{E}(E_{0})^{\otimes d}\right)\\
 & \subset\mathfrak{M}^{+}\left(E_{0}^{d},\mathscr{B}_{\mathcal{F}}(E_{0})^{\otimes d}\right)=\mathcal{M}^{+}\left(E_{0}^{d},\mathscr{O}_{\widehat{E}}(E_{0})^{d}\right).
\end{aligned}
\label{eq:mu_E0d}
\end{equation}
Moreover,
\begin{equation}
\mathcal{M}^{+}\left(E_{0}^{d},\mathscr{O}_{E}(E_{0})^{d}\right)\subset\mathfrak{M}^{+}\left(E_{0}^{d},\mathscr{B}_{E}(E_{0})^{\otimes d}\right).\label{eq:E0d_Borel_Meas}
\end{equation}

\end{fact}
\begin{proof}
The first line of (\ref{eq:mu_E0d}) follows by (\ref{eq:E0d_Prod_Measurable_Ed})
and Fact \ref{fact:Meas_Concen_Expan} (a) (with $E=E^{d}$, $\mathscr{U}=\mathscr{B}(E)^{\otimes d}$
and $A=E_{0}^{d}$). The second line of (\ref{eq:mu_E0d}) and (\ref{eq:E0d_Borel_Meas})
follow by (\ref{eq:Borel_Compare_2}).\end{proof}

\begin{notation}
\label{notation:ReplicaMeas}If no confusion is caused, we will let
$\overline{\mu}$ denote the replica of $\mu$ with respect to the
underlying base and will not make special mention.
\end{notation}
Below are several basic properties of replica measure.
\begin{prop}
\label{prop:RepMeas_Basic}Let $E$ be a topological space, $(E_{0},\mathcal{F};\widehat{E},\widehat{\mathcal{F}})$
be a base over $E$, $d\in\mathbf{N}$ and $\mu\in\mathfrak{M}^{+}(E^{d},\mathscr{B}(E)^{\otimes d})$.
Then:

\renewcommand{\labelenumi}{(\alph{enumi})}
\begin{enumerate}
\item (\ref{eq:RepMeas}) well defines $\overline{\mu}\in\mathcal{M}^{+}(\widehat{E}^{d})$.
Moreover,
\begin{equation}
\overline{\mu}(A)=\mu\left(A\cap E_{0}^{d}\right),\;\forall A\in\mathscr{B}(\widehat{E}^{d})\label{eq:Mu_MuBar}
\end{equation}

\item $\overline{\mu}\in\mathcal{P}(\widehat{E}^{d})$ if and only if $\mu(E_{0}^{d})=1$.
\item Any $\nu\in\mathfrak{be}(\mu)$ satisfies $\overline{\nu}=\overline{\mu}$.
\item If $f\in M_{b}(E^{d};\mathbf{R})$ satisfies $\overline{f}\in M_{b}(\widehat{E}^{d};\mathbf{R})$%
\footnote{$\overline{f}$ was defined in Notation \ref{notation:Rep_Fun}.%
}, then
\begin{equation}
\begin{aligned}\int_{E^{d}}f(x)\mathbf{1}_{E_{0}^{d}}(x)\mu(dx) & =\int_{\left(E_{0}^{d},\mathscr{B}_{E}(E_{0})^{\otimes d}\right)}f|_{E_{0}^{d}}(x)\mu|_{E_{0}^{d}}(dx)\\
 & =\int_{\left(E_{0}^{d},\mathscr{B}_{\widehat{E}^{d}}(E_{0}^{d})\right)}f|_{E_{0}^{d}}(x)\mu|_{E_{0}^{d}}(dx)=\int_{\widehat{E}^{d}}\overline{f}(x)\overline{\mu}(dx).
\end{aligned}
\label{eq:f_fbar_RepMeas_Same_Int}
\end{equation}

\item If $f\in C(E^{d};\mathbf{R})$ has a replica $\widehat{f}$, then
\begin{equation}
\begin{aligned}\int_{E^{d}}f(x)\mathbf{1}_{E_{0}^{d}}(x)\mu(dx) & =\int_{\left(E_{0}^{d},\mathscr{B}_{E}(E_{0})^{\otimes d}\right)}f|_{E_{0}^{d}}(x)\mu|_{E_{0}^{d}}(dx)\\
 & =\int_{\left(E_{0}^{d},\mathscr{B}_{\widehat{E}^{d}}(E_{0}^{d})\right)}f|_{E_{0}^{d}}(x)\mu|_{E_{0}^{d}}(dx)=\int_{\widehat{E}^{d}}\widehat{f}(x)\overline{\mu}(dx).
\end{aligned}
\label{eq:f_fhat_RepMeas_Same_Int}
\end{equation}
In particular, (\ref{eq:f_fhat_RepMeas_Same_Int}) is true for all
$f\in\mathfrak{ca}[\Pi^{d}(\mathcal{F})]$.
\end{enumerate}
\end{prop}
\begin{proof}
(a) follows by (\ref{eq:RepMeas}), (\ref{eq:mu_E0d}) and Fact \ref{fact:Meas_Concen_Expan}
(b)  (with $E=\widehat{E}^{d}$, $A=E_{0}^{d}$ and $\nu=\mu|_{E_{0}^{d}}$)
and (b) is immediate by (\ref{eq:Mu_MuBar}).

(c) $\nu\in\mathfrak{M}^{+}(E_{0}^{d},\mathscr{B}(E)^{\otimes d})$
and $\mu|_{E_{0}^{d}}=\nu|_{E_{0}^{d}}$ as members of $\mathfrak{M}^{+}(E_{0}^{d},\mathscr{B}(E)^{\otimes d}|_{E_{0}^{d}})$.
Thus $\overline{\mu}=\overline{\nu}$ by definition of replica measure.

(d) We have by (\ref{eq:Borel_Compare_2}) that
\begin{equation}
f|_{E_{0}^{d}}=\overline{f}|_{E_{0}^{d}}\in M_{b}\left(E_{0}^{d},\mathscr{B}_{\widehat{E}^{d}}(E_{0}^{d});\mathbf{R}\right)\subset M_{b}\left(E_{0}^{d},\mathscr{B}_{E}(E_{0})^{\otimes d};\mathbf{R}\right).\label{eq:Check_RepMeas_Same_Int_1}
\end{equation}
By (\ref{eq:Check_RepMeas_Same_Int_1}), (\ref{eq:E0d_Prod_Measurable_Ed})
and Fact \ref{fact:Indicator_Modify} (with $E=E^{d}$, $\mathscr{U}=\mathscr{B}(E)^{\otimes d}$
and $A=E_{0}^{d}$),
\begin{equation}
f\mathbf{1}_{E_{0}^{d}}\in M_{b}\left(E^{d},\mathscr{B}(E)^{\otimes d};\mathbf{R}\right).\label{eq:f*1_E0d_Prod_Measurable}
\end{equation}
Now, (\ref{eq:f_fbar_RepMeas_Same_Int}) is well-defined and follows
by (\ref{eq:Check_RepMeas_Same_Int_1}) and (a).

(e) We have by (\ref{eq:Borel_Compare_2}) that
\begin{equation}
\begin{aligned} & f|_{E_{0}^{d}}=\widehat{f}|_{E_{0}^{d}}\in C_{b}\left(E_{0}^{d},\mathscr{O}_{\widehat{E}^{d}}(E_{0}^{d});\mathbf{R}\right)\\
 & \subset M_{b}\left(E_{0}^{d},\mathscr{O}_{\widehat{E}^{d}}(E_{0}^{d});\mathbf{R}\right)\subset M_{b}\left(E_{0}^{d},\mathscr{B}_{E}(E_{0})^{\otimes d};\mathbf{R}\right).
\end{aligned}
\label{eq:Check_RepMeas_Same_Int_2}
\end{equation}
So, (\ref{eq:f_fhat_RepMeas_Same_Int}) is well-defined and follows
by (\ref{eq:Check_RepMeas_Same_Int_2}) and (a).\end{proof}

The next proposition gives a sufficient condition for a Borel measure
on $\widehat{E}^{d}$ to be the replica of some Borel measure on $E^{d}$.
\begin{prop}
\label{prop:Redefine_FDD}Let $E$ be a topological space, $(E_{0},\mathcal{F};\widehat{E},\widehat{\mathcal{F}})$
be a base over $E$ and $d\in\mathbf{N}$. If $\nu\in\mathcal{M}^{+}(\widehat{E}^{d})$
is supported on $A\subset E_{0}^{d}$ and $A\in\mathscr{B}^{\mathbf{s}}(E^{d})$,
then
\begin{equation}
\mu\circeq\left.(\nu|_{A})\right|^{E^{d}}\in\mathcal{M}^{+}(E^{d})\label{eq:Redefine_FDD}
\end{equation}
satisfies $\mu|_{A}=\nu|_{A}\in\mathcal{M}^{+}(A,\mathscr{O}_{E^{d}}(A))$
and $\nu=\overline{\mu}$.
\end{prop}
\begin{proof}
We have by $A\in\mathscr{B}^{\mathbf{s}}(E^{d})$ and Lemma \ref{lem:SB_Base}
(a, b) that
\begin{equation}
\mathscr{B}_{E^{d}}(A)=\mathscr{B}_{\widehat{E}^{d}}(A)\subset\mathscr{B}(E^{d})\cap\mathscr{B}(\widehat{E}^{d}).\label{eq:SB_Support}
\end{equation}
By Fact \ref{fact:Meas_Concen_Expan} (a) (with $\mu=\nu$ and $(E,\mathscr{U})=(\widehat{E}^{d},\mathscr{B}(\widehat{E}^{d})$)
and (\ref{eq:SB_Support}),
\begin{equation}
\nu|_{A}\in\mathcal{M}^{+}\left(A,\mathscr{O}_{E^{d}}(A)\right)=\mathcal{M}^{+}\left(A,\mathscr{O}_{\widehat{E}^{d}}(A)\right).\label{eq:Nu|A_is_Borel_Meas_Ed}
\end{equation}
Then, $\mu\in\mathcal{M}^{+}(E^{d})$ by (\ref{eq:Nu|A_is_Borel_Meas_Ed})
and Fact \ref{fact:Meas_Concen_Expan} (b) (with $E=E^{d}$, $\mathscr{U}=\mathscr{B}(E^{d})$
and $\nu=\nu|_{A}$). It follows by Fact \ref{fact:Meas_Concen_Expan}
(c) (with $E=E^{d}$, $\mathscr{U}=\mathscr{B}(E^{d})$ and $\nu=\nu|_{A}$)
and (\ref{eq:Nu|A_is_Borel_Meas_Ed}) that
\begin{equation}
\mu|_{A}=\nu|_{A}\in\mathcal{M}^{+}\left(A,\mathscr{O}_{\widehat{E}^{d}}(A)\right).\label{eq:FDD_Redef_FDD_Same_on_Support_Ed_Ehatd}
\end{equation}
It follows by the fact $\nu(\widehat{E}^{d}\backslash A)=0$, (\ref{eq:FDD_Redef_FDD_Same_on_Support_Ed_Ehatd})
and Fact \ref{fact:Meas_Concen_Expan} (c) (with $E=\widehat{E}^{d}$,
$\mathscr{U}=\mathscr{B}(\widehat{E}^{d})$ and $\mu=\nu$) that
\begin{equation}
\nu=\left.\left(\nu|_{A}\right)\right|^{\widehat{E}^{d}}=\left.\left(\mu|_{A}\right)\right|^{\widehat{E}^{d}}=\overline{\mu}.\label{eq:Check_Redef_Rep}
\end{equation}
\end{proof}

\subsection{\label{sub:RepMeas_WC}Weak convergence of replica measures}

We now consider the association of weak convergence of Borel extensions
on $E^{d}$ and that of replica measures on $\widehat{E}^{d}$. The
next proposition discusses the direction from $E^{d}$ to $\widehat{E}^{d}$.
\begin{prop}
\textup{\label{prop:WC_Push_Forward}}Let $E$ be a topological space,
$(E_{0},\mathcal{F};\widehat{E},\widehat{\mathcal{F}})$ be a base
over $E$, $d\in\mathbf{N}$, $\mathcal{G}\circeq\mathfrak{mc}[\Pi^{d}(\mathcal{F}\backslash\{1\})]$
and $\{\mu_{n}\}_{n\in\mathbf{N}}\cup\{\mu\}\subset\mathfrak{M}^{+}(E^{d},\mathscr{B}(E)^{\otimes d})$.
Consider the following statements:

\renewcommand{\labelenumi}{(\alph{enumi})}
\begin{enumerate}
\item The replica measures $\{\overline{\mu}_{n}\}_{n\in\mathbf{N}}$ and
$\overline{\mu}$ satisfy
\begin{equation}
\mathrm{w}\mbox{-}\lim_{n\rightarrow\infty}\overline{\mu}_{n}=\overline{\mu}\mbox{ in }\mathcal{M}^{+}(\widehat{E}^{d}).\label{eq:MuBar_n_WC_MuBar_M(Ehatd)}
\end{equation}

\item The concentrated measures $\{\mu_{n}|_{E_{0}^{d}}\}_{n\in\mathbf{N}}$
and $\mu|_{E_{0}^{d}}$ satisfy
\begin{equation}
\mathrm{w}\mbox{-}\lim_{n\rightarrow\infty}\mu_{n}|_{E_{0}^{d}}=\mu|_{E_{0}^{d}}\mbox{ in }\mathcal{M}^{+}\left(E_{0}^{d},\mathscr{O}_{\widehat{E}}(E_{0})^{d}\right).\label{eq:Mu_n_Rep_WC_Mu_M(E0d)}
\end{equation}

\item The original measures $\{\mu_{n}\}_{n\in\mathbf{N}}$ and $\mu$ satisfy%
\footnote{The integrals in (\ref{eq:Int_Test_WC_on_E0d}) are well-defined by
the proof of Proposition \ref{prop:RepMeas_Basic} (d, e). Those in
(\ref{eq:Int_Test_WC_on_Ed}) are well-defined by Note \ref{note:Pi^d(D)_Mb_Cb}
(with $\mathcal{D}=\mathcal{F}$).%
}
\begin{equation}
\lim_{n\rightarrow\infty}\int_{E^{d}}f(x)\mathbf{1}_{E_{0}^{d}}(x)\mu_{n}(dx)=\int_{E^{d}}f(x)\mathbf{1}_{E_{0}^{d}}(x)\mu(dx),\;\forall f\in\mathcal{G}\cup\{1\}.\label{eq:Int_Test_WC_on_E0d}
\end{equation}

\item $E_{0}^{d}$ is a common support of $\{\mu_{n}\}_{n\in\mathbf{N}}\cup\{\mu\}$.
Moreover,
\begin{equation}
\lim_{n\rightarrow\infty}\int_{E^{d}}f(x)\mu_{n}(dx)=\int_{E^{d}}f(x)\mu(dx),\;\forall f\in\mathcal{G}\cup\{1\}.\label{eq:Int_Test_WC_on_Ed}
\end{equation}

\item $E_{0}^{d}$ is a common support of $\{\mu_{n}\}_{n\in\mathbf{N}}\cup\{\mu\}$.
Moreover, there exist $\{\mu_{n}^{\prime}\in\mathfrak{be}(\mu_{n})\}_{n\in\mathbf{N}}$
and $\mu^{\prime}\in\mathfrak{be}(\mu)$ such that
\begin{equation}
\mu_{n}^{\prime}\Longrightarrow\mu^{\prime}\mbox{ as }n\uparrow\infty\mbox{ in }\mathcal{M}^{+}(E^{d}).\label{eq:Mu_n_BExt_WC_Mu_BExt_M(Ed)}
\end{equation}

\item $E_{0}^{d}$ is a common support of $\{\mu_{n}\}_{n\in\mathbf{N}}\cup\{\mu\}$.
Moreover, there exist $\{\mu_{n}^{\prime}\in\mathfrak{be}(\mu_{n})\}_{n\in\mathbf{N}}$
and $\mu^{\prime}\in\mathfrak{be}(\mu)$ such that
\begin{equation}
\mu_{n}^{\prime}|_{E_{0}^{d}}\Longrightarrow\mu^{\prime}|_{E_{0}^{d}}\mbox{ as }n\uparrow\infty\mbox{ in }\mathcal{M}^{+}\left(E_{0}^{d},\mathscr{O}_{E}(E_{0})^{d}\right).\label{eq:Mu_n_BExt_WC_Mu_BExt_M(E0d)}
\end{equation}

\end{enumerate}
Then, (a) - (c) are equivalent. (c) - (f) are successively stronger.
Moreover, (e) and (f) are equivalent when $E^{d}$ is a Tychonoff
space.
\end{prop}
\begin{proof}
((a) $\rightarrow$ (b)) $\widehat{E}^{d}$ is a Tychonoff space by
Lemma \ref{lem:Base_Property} (c) and Proposition \ref{prop:CR_Space}
(a). $(E_{0}^{d},\mathscr{O}_{\widehat{E}}(E_{0})^{d})$ is a metrizable
and separable space by Lemma \ref{lem:Base_Property} (d) (with $A=E_{0}^{d}$)
and is a Tychonoff space by Proposition \ref{prop:CR_Space} (a).
Now, (b) follows by (\ref{eq:MuBar_n_WC_MuBar_M(Ehatd)}), (\ref{eq:RepMeas}),
Lemma \ref{lem:WC_Expansion} (with $E=\widehat{E}^{d}$, $A=E_{0}^{d}$,
$\nu_{n}=\mu_{n}|_{E_{0}^{d}}$ and $\nu=\mu|_{E_{0}^{d}}$) and the
Tychonoff property of $(E_{0}^{d},\mathscr{O}_{\widehat{E}}(E_{0})^{d}$.

((b) $\rightarrow$ (c)) We have by Lemma \ref{lem:Base_Property}
(b, d) (with $A=E_{0}^{d}$) that $\mathcal{G}|_{E_{0}^{d}}\subset C_{b}(E_{0}^{d},\mathscr{O}_{\widehat{E}}(E_{0})^{d};\mathbf{R})$.
Hence, we have by (\ref{eq:Mu_n_Rep_WC_Mu_M(E0d)}) and (\ref{eq:mu_E0d})
that
\begin{equation}
\begin{aligned} & \lim_{n\rightarrow\infty}\int_{E^{d}}f(x)\mathbf{1}_{E_{0}^{d}}(x)\mu_{n}(dx)\\
 & =\lim_{n\rightarrow\infty}\int_{\left(E_{0}^{d},\mathscr{O}_{\widehat{E}}(E_{0})^{d}\right)}f|_{E_{0}^{d}}(x)\mu_{n}|_{E_{0}^{d}}(dx)\\
 & =\int_{\left(E_{0}^{d},\mathscr{O}_{\widehat{E}}(E_{0})^{d}\right)}f|_{E_{0}^{d}}(x)\mu|_{E_{0}^{d}}(dx)\\
 & =\int_{E^{d}}f(x)\mathbf{1}_{E_{0}^{d}}(x)\mu(dx),\;\forall f\in\mathcal{G}\cup\{1\}.
\end{aligned}
\label{eq:Check_Int_Test_WC_on_E0d}
\end{equation}

((c) $\rightarrow$ (a)) It follows by Proposition \ref{prop:RepMeas_Basic}
(e) (with $\mu=\mu_{n}$ and $\mu$) that
\begin{equation}
\begin{aligned}\lim_{n\rightarrow\infty}\widehat{f}^{*}(\overline{\mu}_{n}) & =\lim_{n\rightarrow\infty}\int_{E^{d}}f(x)\mathbf{1}_{E_{0}^{d}}(x)\mu_{n}(dx)\\
 & =\int_{E^{d}}f(x)\mathbf{1}_{E_{0}^{d}}(x)\mu(dx)=\widehat{f}^{*}(\overline{\mu}),\;\forall f\in\mathcal{G}\cup\{1\}.
\end{aligned}
\label{eq:Int_Test_WC_on_Ehatd}
\end{equation}
$\mathcal{M}^{+}(\widehat{E}^{d})$ is a metrizable space by Corollary
\ref{cor:Base_Sep_Meas} (c). Now, (a) follows by (\ref{eq:Int_Test_WC_on_Ehatd}),
Corollary \ref{cor:Base_Sep_Meas} (b) (with $A=\widehat{E}^{d}$)
and the Hausdorff property of $\mathcal{M}^{+}(\widehat{E}^{d})$.

((d) $\rightarrow$ (c)) Note that if $E_{0}^{d}$ is a common support
of $\{\mu_{n}\}_{n\in\mathbf{N}}\cup\{\mu\}$, then
\begin{equation}
\begin{aligned} & \int_{E^{d}}f(x)\mathbf{1}_{E_{0}^{d}}(x)\mu_{n}(dx)-\int_{E^{d}}f(x)\mu_{n}(dx)\\
 & =\int_{E^{d}}f(x)\mathbf{1}_{E_{0}^{d}}(x)\mu(dx)-\int_{E^{d}}f(x)\mu(dx)=0,\;\forall n\in\mathbf{N}.
\end{aligned}
\label{eq:Check_Same_Int_E0d_1}
\end{equation}

((e) $\rightarrow$ (d)) follows by $\mathcal{F}\subset C_{b}(E;\mathbf{R})$
and Fact \ref{fact:BExt_Same_Int} (with $\mu=\mu_{n}$ or $\mu$).

In both (e) and (f), $E_{0}^{d}$ is a common support of $\{\mu_{n}^{\prime}\}_{n\in\mathbf{N}}\cup\{\mu^{\prime}\}$
and so $\mu_{n}^{\prime}=(\mu_{n}^{\prime}|_{E_{0}^{d}})|^{E^{d}}$
for all $n\in\mathbf{N}$ and $\mu^{\prime}=(\mu^{\prime}|_{E_{0}^{d}})|^{E^{d}}$
by Fact \ref{fact:Meas_Concen_Expan} (c) (with $E=E^{d}$, $\mathscr{U}=\mathscr{B}(E^{d})$,
$A=E_{0}^{d}$ and $\mu=\mu_{n}^{\prime}$ or $\mu^{\prime}$). It
then follows by Lemma \ref{lem:WC_Expansion} (with $E=E^{d}$, $A=E_{0}^{d}$,
$\mu_{n}=\mu_{n}^{\prime}|_{E_{0}^{d}}$ and $\mu=\mu^{\prime}|_{E_{0}^{d}}$)
that (f) implies (e) in general, and (e) implies (f) when $E^{d}$
is a Tychonoff space.\end{proof}

The following corollary specializes Proposition \ref{prop:WC_Push_Forward}
to probability measures.
\begin{cor}
\textup{\label{cor:WC_Push_Forward_P(E)}}Let $E$ be a topological
space, $(E_{0},\mathcal{F};\widehat{E},\widehat{\mathcal{F}})$ be
a base over $E$, $d\in\mathbf{N}$, $\mathcal{G}\circeq\mathfrak{mc}[\Pi^{d}(\mathcal{F}\backslash\{1\})]$
and $\{\mu_{n}\}_{n\in\mathbf{N}}\cup\{\mu\}\subset\mathfrak{P}(E^{d},\mathscr{B}(E)^{\otimes d})$.
Then, the following statements are equivalent:

\renewcommand{\labelenumi}{(\alph{enumi})}
\begin{enumerate}
\item The replica measures $\{\overline{\mu}_{n}\}_{n\in\mathbf{N}}$ and
$\overline{\mu}$ satisfy
\begin{equation}
\overline{\mu}_{n}\Longrightarrow\overline{\mu}\mbox{ as }n\uparrow\infty\mbox{ in }\mathcal{P}(\widehat{E}^{d}).\label{eq:MuBar_n_WC_MuBar_P(Ehatd)}
\end{equation}

\item The concentrated measures $\{\mu_{n}|_{E_{0}^{d}}\}_{n\in\mathbf{N}}$
and $\mu|_{E_{0}^{d}}$ satisfy
\begin{equation}
\mu_{n}|_{E_{0}^{d}}\Longrightarrow\mu|_{E_{0}^{d}}\mbox{ as }n\uparrow\infty\mbox{ in }\mathcal{P}\left(E_{0}^{d},\mathscr{O}_{\widehat{E}}(E_{0})^{d}\right).\label{eq:Mu_n_Rep_WC_Mu_P(E0d)}
\end{equation}

\item $E_{0}^{d}$ is a common support of $\{\mu_{n}\}_{n\in\mathbf{N}}\cup\{\mu\}$
and (\ref{eq:Int_Test_WC_on_Ed}) holds.
\end{enumerate}
\end{cor}
\begin{proof}
The result follows by Proposition \ref{prop:RepMeas_Basic} (b) and
Proposition \ref{prop:WC_Push_Forward} (a, b, d).\end{proof}

One can use sequential tightness to transform a weak limit point of
replica measures back into that of Borel extensions of the original
measures based upon the following two results.
\begin{prop}
\label{prop:WC_Pull_Back_FinDim}Let $E$ be a topological space,
$(E_{0},\mathcal{F};\widehat{E},\widehat{\mathcal{F}})$ be a base
over $E$, $d\in\mathbf{N}$, $\mathcal{G}\circeq\mathfrak{mc}[\Pi^{d}(\mathcal{F}\backslash\{1\})]$
and $\{\mu_{n}\}_{n\in\mathbf{N}}\subset\mathfrak{M}^{+}(E^{d},\mathscr{B}(E)^{\otimes d})$.
Suppose that:

\renewcommand{\labelenumi}{(\roman{enumi})}
\begin{enumerate}
\item $\{\mu_{n}\}_{n\in\mathbf{N}}$ is sequentially tight in $E_{0}^{d}$.
\item $\{\int_{E^{d}}f(x)\mu_{n}(dx)\}_{n\in\mathbf{N}}$ is convergent
in $\mathbf{R}$ for all $f\in\mathcal{G}\cup\{1\}$.
\item $\{\mu_{n}(E^{d})\}_{n\in\mathbf{N}}\subset[a,b]$ for some $0<a<b$.
\end{enumerate}
Then, there exist $\mu\in\mathcal{M}^{+}(E^{d})$ and $N\in\mathbf{N}$
such that:

\renewcommand{\labelenumi}{(\alph{enumi})}
\begin{enumerate}
\item $\mu$ is $\mathbf{m}$-tight in $E_{0}^{d}$ and $\{\mu_{n}^{\prime}=\mathfrak{be}(\mu_{n})\}_{n>N}$
exists.
\item $\{\mu_{n}^{\prime}\}_{n>N}$ satisfies
\begin{equation}
\mathrm{w}\mbox{-}\lim_{n\rightarrow\infty}\mu_{n}^{\prime}|_{E_{0}^{d}}=\mu|_{E_{0}^{d}}\mbox{ in }\mathcal{M}^{+}\left(E_{0}^{d},\mathscr{O}_{E}(E_{0})^{d}\right)\label{eq:Mu_n_BExt_WLim_Mu_M(E0d)}
\end{equation}
and
\begin{equation}
\mu_{n}^{\prime}\Longrightarrow\mu\mbox{ as }n\uparrow\infty\mbox{ in }\mathcal{M}^{+}(E^{d}).\label{eq:Mu_n_BExt_WC_Mu_M(Ed)}
\end{equation}

\end{enumerate}
\end{prop}
\begin{proof}
$\{\mu_{n}\}_{n\in\mathbf{N}}$ is sequentially $\mathbf{m}$-tight
in $(E_{0}^{d},\mathscr{O}_{E}(E_{0})^{d})$ by Corollary \ref{cor:Base_Compact}
(a). There exists an $N_{1}\in\mathbf{N}$ such that $\{\mu_{n}\}_{n>N_{1}}$
are all supported on $E_{0}^{d}$ by Fact \ref{fact:Seq_Tight_Support}
(with $(E,\mathscr{U})=(E^{d},\mathscr{B}(E)^{\otimes d})$, $A=E_{0}^{d}$
and $\Gamma=\{\mu_{n}\}_{n\in\mathbf{N}}$). There exist $\nu\in\mathcal{M}^{+}(E_{0}^{d},\mathscr{O}_{E}(E_{0})^{d})$
and $N_{2}\in\mathbf{N}$ such that $\nu$ is $\mathbf{m}$-tight
in $(E_{0}^{d},\mathscr{O}_{E}(E_{0})^{d})$, $\{\nu_{n}^{\prime}=\mathfrak{be}(\mu_{n}|_{E_{0}^{d}})\}_{n>N_{2}}$
exists and
\begin{equation}
\mathrm{w}\mbox{-}\lim_{n\rightarrow\infty}\nu_{n}^{\prime}=\nu\mbox{ in }\mathcal{M}^{+}\left(E_{0}^{d},\mathscr{O}_{E}(E_{0})^{d}\right)\label{eq:Check_Mu_n_BExt_WLim_Mu_M(E0d)}
\end{equation}
by Lemma \ref{lem:Base_Property} (e) (with $A=E_{0}^{d}$) and Theorem
\ref{thm:WLP_Uni} (a, c) (with $E=(E_{0},\mathscr{O}_{E}(E_{0}))$,
$\Gamma=\{\mu_{n}|_{E_{0}^{d}}\}_{n\in\mathbf{N}}$ and $\mathcal{D}=\mathcal{F}\backslash\{1\}$).

$\mu\circeq\nu|^{E^{d}}$ satisfies $\mu|_{E_{0}^{d}}=\nu$ by (\ref{eq:E0d_Prod_Measurable_Ed})
and Fact \ref{fact:Meas_Concen_Expan} (c) (with $(E,\mathscr{U})=(E^{d},\mathscr{B}(E)^{\otimes d}$
and $A=E_{0}^{d}$).
\begin{equation}
\nu_{n}^{\prime}|^{E^{d}}=\mathfrak{be}\left[\left.(\mu_{n}|_{E_{0}^{d}})\right|^{E^{d}}\right]=\mathfrak{be}(\mu_{n})=\mu_{n}^{\prime},\;\forall n>N\circeq N_{1}\vee N_{2}\label{eq:Mu_n_BExt=00003DNu_n}
\end{equation}
by (\ref{eq:E0d_Prod_Measurable_Ed}), Fact \ref{fact:Meas_Concen_Expan}
(c) (with $(E,\mathscr{U})=(E^{d},\mathscr{B}(E)^{\otimes d})$, $A=E_{0}^{d}$
and $\nu=\nu_{n}$) and Lemma \ref{lem:Union_Borel_Prod_Equal} (b)
(with $\mathbf{I}=\{1,...,d\}$, $S_{i}=E$, $A=E_{0}^{d}$, $\mu=\mu_{n}$
and $\mathfrak{be}(\mu|_{A})=\nu_{n}$).

Hence, (\ref{eq:Mu_n_BExt_WLim_Mu_M(E0d)}) follows by (\ref{eq:Check_Mu_n_BExt_WLim_Mu_M(E0d)})
and Fact \ref{fact:Meas_Concen_Expan} (c) (with $E=E^{d}$, $\mathcal{U}=\mathcal{B}(E)^{\otimes d}$
and $A=E_{0}^{d}$). (\ref{eq:Mu_n_BExt_WC_Mu_M(Ed)}) follows by
(\ref{eq:Mu_n_BExt_WLim_Mu_M(E0d)}) and Lemma \ref{lem:WC_Expansion}
(with $E=E^{d}$, $A=E_{0}^{d}$, $\mu_{n}=\nu_{n}$ and $\mu=\nu$).\end{proof}

\begin{cor}
\label{cor:WC_Pull_Back_FinDim}Let $(E_{0},\mathcal{F};\widehat{E},\widehat{\mathcal{F}})$
be a base over topological space $E$ and $d\in\mathbf{N}$. If $\{\mu_{n}\}_{n\in\mathbf{N}}\subset\mathfrak{M}^{+}(E^{d},\mathscr{B}(E)^{\otimes d})$
is sequentially tight in $(E_{0}^{d},\mathscr{O}_{E}(E_{0})^{d})$,
and if their replicas satisfy
\begin{equation}
\overline{\mu}_{n}\Longrightarrow\nu\mbox{ as }n\uparrow\infty\mbox{ in }\mathcal{M}^{+}(\widehat{E}^{d}),\label{eq:MuBar_n_WC_Nu_M(Ehatd)}
\end{equation}
then there exist $\mu\in\mathcal{M}^{+}(E^{d})$ and $N\in\mathbf{N}$
satisfying Proposition \ref{prop:WC_Pull_Back_FinDim} (a, b) and,
in particular, $\nu=\overline{\mu}$.
\end{cor}
\begin{proof}
There exists an $N_{1}\in\mathbf{N}$ such that $\{\mu_{n}\}_{n>N_{1}}$
are all supported on $E_{0}^{d}$ by Fact \ref{fact:Seq_Tight_Support}
(with $(E,\mathscr{U})=(E^{d},\mathscr{B}(E)^{\otimes d})$, $A=E_{0}^{d}$
and $\Gamma=\{\mu_{n}\}_{n\in\mathbf{N}}$). It follows by (\ref{eq:MuBar_n_WC_Nu_M(Ehatd)})
and Proposition \ref{prop:RepMeas_Basic} (e) (with $\mu=\mu_{n}$)
that
\begin{equation}
\lim_{n\rightarrow\infty}\int_{E^{d}}f(x)\mu_{n}(dx)=\lim_{n\rightarrow\infty}\widehat{f}^{*}(\overline{\mu}_{n})=\widehat{f}^{*}(\nu),\;\forall f\in\mathfrak{mc}\left[\Pi^{d}(\mathcal{F})\right].\label{eq:Check_Int_Test_Lim}
\end{equation}
$\nu\in\mathcal{M}^{+}(\widehat{E}^{d})$ means $\nu(\widehat{E}^{d})>0$
by our convention in \S \ref{sub:Meas}. It then follows by the fact
$1\in\Pi^{d}(\mathcal{F})$ and (\ref{eq:Check_Int_Test_Lim}) that
\begin{equation}
\mu_{n}(E^{d})\in\left(\frac{\nu(\widehat{E}^{d})}{2},\frac{3\nu(\widehat{E}^{d})}{2}\right)\subset(0,\infty),\;\forall n>N_{2}\label{eq:Check_Total_Mass_Not_Vanish}
\end{equation}
for some $N_{2}\in\mathbf{N}\cap(N_{1},\infty)$. Now, we obtain the
desired $\mu$ and $N$ by (\ref{eq:Check_Int_Test_Lim}), (\ref{eq:Check_Total_Mass_Not_Vanish})
and Proposition \ref{prop:WC_Pull_Back_FinDim} (with $n=N_{2}+n$,
$a=\nu(\widehat{E})/2$ and $b=3a$). $\{\mu_{n}^{\prime}=\mathfrak{be}(\mu_{n})\}_{n\in\mathbf{N}}$
satisfies (\ref{eq:Mu_n_BExt_WC_Mu_M(Ed)}), so we have%
\footnote{$\overline{\mu}_{n}^{\prime}$ denote the replica of $\mu_{n}^{\prime}$.%
}
\begin{equation}
\overline{\mu}=\mathrm{w}\mbox{-}\lim_{n\rightarrow\infty}\overline{\mu}_{n}^{\prime}=\mathrm{w}\mbox{-}\lim_{n\rightarrow\infty}\overline{\mu}_{n}=\nu\label{eq:Check_Lim_RepMeas}
\end{equation}
by Proposition \ref{prop:RepMeas_Basic} (c) (with $\mu=\mu_{n}$
and $\nu=\mu_{n}^{\prime}$) and Proposition \ref{prop:WC_Push_Forward}
(a, e).\end{proof}

\section{\label{sec:Riesz_Sko}Generalization of two fundamental results}

\subsection{\label{sub:Riesz_Representation}Integral representation of linear
functional}

The celebrated Riesz-Radon Representation Theorem was established
for \textit{positive} \textit{linear functionals on }$C_{0}(E;\mathbf{R})$\textit{}%
\footnote{Positiveness of a functional on $C_{0}(E;\mathbf{R})$ means it maps
non-negative functions into $\mathbf{R}^{+}$.%
} with $E$ being a locally compact Hausdorff space. This result is
now extended to baseable spaces by replication, avoiding the local
compactness assumption which is violated by many infinite-dimensional
spaces. As mentioned in \S \ref{sub:Examples_Baseable_Space}, baseable
spaces need not be locally compact nor Tychonoff.
\begin{thm}
\label{thm:Riesz}Let $E$ be a $C_{c}(E;\mathbf{R})$-baseable space,
$\varphi$ be a linear functional on $C_{c}(E;\mathbf{R})$ and
\begin{equation}
\mathcal{B}\circeq\left\{ g\in C_{c}(E;\mathbf{R}):0<\Vert g\Vert_{\infty}\leq1\right\} .\label{eq:Unit_Ball}
\end{equation}
Then, the following statements are equivalent:

\renewcommand{\labelenumi}{(\alph{enumi})}
\begin{enumerate}
\item There exists a positive linear functional $\Lambda$ on $C_{b}(E;\mathbf{R})$
such that
\begin{equation}
\varphi(g)\leq\Lambda(g),\;\forall g\in C_{c}(E;\mathbf{R})\label{eq:Original_Functional_Dominated}
\end{equation}
and
\begin{equation}
\lambda_{0}\circeq\sup_{g\in\mathcal{B}}\varphi(g)=\Lambda(1)<\infty.\label{eq:Define_Total_Mass}
\end{equation}

\item There exists unqiue $\mathbf{m}$-tight $\mu\in\mathcal{M}^{+}(E)$
such that
\begin{equation}
\varphi(g)=g^{*}(\mu),\;\forall g\in C_{c}(E;\mathbf{R})\label{eq:Riesz_Representation}
\end{equation}
and
\begin{equation}
\mu(E)=\sup_{g\in\mathcal{B}}g^{*}(\mu).\label{eq:Approx_Total_Mass}
\end{equation}

\end{enumerate}
\end{thm}
\begin{rem}
\label{rem:Reisz_Sup_Well_Defined}In the theorem above, $E$ is Hausdorff
by Fact \ref{fact:Baseable_Metrizable_Separable} (a). $C_{c}(E;\mathbf{R})$
is a possibly \textit{non-unit}%
\footnote{``non-unit'' means excluding the constant function $1$.%
} subalgebra of $C_{b}(E;\mathbf{R})$ and is a function lattice%
\footnote{The terminology ``function lattice'' was specified in \S \ref{sub:Fun}.%
} by Proposition \ref{prop:Cc_Lattice} (a). $C_{c}(E;\mathbf{R})\neq\{0\}$
since $C_{c}(E;\mathbf{R})$ separates points on $E$, so $\mathcal{B}\neq\varnothing$
and the supremum in (\ref{eq:Define_Total_Mass}) is well-defined.
\end{rem}
\begin{proof}
[Proof of Theorem \ref{thm:Riesz}]((a) $\rightarrow$ (b)) We divide
our proof into six steps.

\textit{Step 1: Extend $\varphi$ to a positive linear functional
on $C_{b}(E;\mathbf{R})$}. $\varphi(g)\le\Lambda(g)\leq0$ for all
non-positive $g\in C_{c}(E;\mathbf{R})$ by (\ref{eq:Original_Functional_Dominated})
and the positiveness of $\Lambda$, so $\varphi$ is also a positive
linear functional. Then, there exists a positive linear functional
$\Phi$ on $C_{b}(E;\mathbf{R})$ satisfying
\begin{equation}
\varphi=\Phi|_{C_{c}(E;\mathbf{R})}\label{eq:Extend_Original_Functional}
\end{equation}
and
\begin{equation}
\Phi(g)\leq\Lambda(g),\;\forall g\in C_{b}(E;\mathbf{R})\label{eq:Extend_Functional_Dominated}
\end{equation}
by a suitable version of the Hahn-Banach Theorem (see \cite[Theorem 8.31]{AB06}).
In particular,
\begin{equation}
\lambda_{0}=\Lambda(1)\geq\Phi(1)\geq\sup_{g\in C_{b}(E;\mathbf{R})\backslash\{0\}}\Phi\left(\frac{g}{\Vert g\Vert_{\infty}}\right)\geq\sup_{g\in\mathcal{B}}\Phi(g)=\lambda_{0}\label{eq:Extend_Functional_Bounded}
\end{equation}
by (\ref{eq:Define_Total_Mass}), (\ref{eq:Extend_Functional_Dominated}),
the positiveness of $\Phi$ and the fact
\begin{equation}
g\leq\Vert g\Vert_{\infty}\le1,\;\forall g\in\mathcal{B}.\label{eq:g_in_B_<=00003D1}
\end{equation}

\textit{Step 2: Construct a suitable base}. Letting
\begin{equation}
f_{g,a,b}\circeq ag+b,\;\forall g\in C_{c}(E;\mathbf{R}),a,b\in\mathbf{R},\label{eq:f_g_a_b}
\end{equation}
we have that
\begin{equation}
\mathcal{D}\circeq\mathfrak{ag}\left(C_{c}(E;\mathbf{R})\cup\{1\}\right)=\left\{ f_{g,a,b}:g\in C_{c}(E;\mathbf{R}),a,b\in\mathbf{R}\right\} .\label{eq:D_f_g_a_b}
\end{equation}
$E$ is a $\mathcal{D}$-baseable space by Fact \ref{fact:D-Baseable}
(d) (with $A=E$, $\mathcal{D}=C_{c}(E;\mathbf{R})$ and $\mathcal{D}^{\prime}=\mathcal{D}$).
There exist $\{g_{p}\}_{p\in\mathbf{N}}\subset\mathcal{B}$ satisfying
\begin{equation}
\Phi(1)=\lambda_{0}=\lim_{p\rightarrow\infty}\varphi(g_{p})=\lim_{p\rightarrow\infty}\Phi(g_{p})\label{eq:Total_Mass_Approx_Seq}
\end{equation}
by (\ref{eq:Extend_Functional_Bounded}), (\ref{eq:Define_Total_Mass})
and (\ref{eq:Extend_Original_Functional}). We then find a base $(E,\mathcal{F};\widehat{E},\widehat{\mathcal{F}})$
over $E$ satisfying
\begin{equation}
\{g_{p}\}_{p\in\mathbf{N}}\subset(\mathcal{F}\cap\mathcal{B})\subset(\mathcal{F}\backslash\{1\})\subset C_{c}(E;\mathbf{R})\subset\mathcal{D}\label{eq:gp_in_F_B}
\end{equation}
by Lemma \ref{lem:Base_Construction} (c) (with $E_{0}=E$ and $\mathcal{D}_{0}=\{g_{p}\}_{p\in\mathbf{N}}$).

\textit{Step 3: Construct a replica positive linear functional $\widehat{\Phi}$
on $C(\widehat{E};\mathbf{R})$} satisfying
\begin{equation}
\widehat{\Phi}(h)=\Phi(h|_{E}),\;\forall h\in C(\widehat{E};\mathbf{R}).\label{eq:Rep_Functional}
\end{equation}
$E_{0}=E$ here, so%
\footnote{We noted in Notation \ref{notation:Rep_Fun} that $\overline{g}\circeq\mathfrak{var}(g;\widehat{E},E_{0},0)$.%
}
\begin{equation}
\widehat{f}_{g,a,b}=a\widehat{g}+b=a\overline{g}+b,\;\forall f_{g,a,b}\in\mathcal{D}\label{eq:Rep_Cc_Fun}
\end{equation}
by Proposition \ref{prop:RepFun_Basic} (d) (with $d=k=1$ and $E_{0}=E$),
Lemma \ref{lem:Base} (c) and Fact \ref{fact:Cc_Ext} (b) (with $E=\widehat{E}$,
$A=E$ and $f=\widehat{g}$).
\begin{equation}
\mathfrak{ag}(\widehat{\mathcal{F}})=\left\{ \widehat{f}_{g,a,b}:f_{g,a,b}\in\mathfrak{ag}(\mathcal{F})\right\} \label{eq:Rep_Functional_Subdomain}
\end{equation}
is a linear subspace of $C(\widehat{E};\mathbf{R})$ on which
\begin{equation}
\widehat{\Phi}(\widehat{f}_{g,a,b})\circeq a\varphi(g)+b\lambda_{0},\;\forall\widehat{f}_{g,a,b}\in\mathfrak{ag}(\widehat{\mathcal{F}})\label{eq:Define_Rep_Functional}
\end{equation}
defines a positive linear functional. Moreover, 
\begin{equation}
\widehat{\Phi}(\widehat{f}_{g,a,b})=\Phi(f_{g,a,b})\leq\lambda_{0}\left\Vert f_{g,a,b}\right\Vert _{\infty}=\lambda_{0}\left\Vert \widehat{f}_{g,a,b}\right\Vert _{\infty},\;\forall f_{g,a,b}\in\mathfrak{ag}(\mathcal{F})\label{eq:Check_Rep_Functional_Bounded}
\end{equation}
by (\ref{eq:Extend_Original_Functional}), the first equality of (\ref{eq:Total_Mass_Approx_Seq})
and Fact \ref{fact:f+_Rep} (a) (with $d=k=1$, $E_{0}=E$ and $f=f_{g,a,b}$).
Hence, $\widehat{\Phi}$ extends linearly onto $C(\widehat{E};\mathbf{R})$
and satisfies
\begin{equation}
\lambda_{0}=\Phi(1)=\widehat{\Phi}(1)=\sup_{h\in C(\widehat{E};\mathbf{R})}\frac{\widehat{\Phi}(h)}{\Vert h\Vert_{\infty}}<\infty\label{eq:Rep_Functional_Bounded}
\end{equation}
by (\ref{eq:Define_Rep_Functional}), (\ref{eq:Check_Rep_Functional_Bounded})
(with $a=0$ and $b=1$) and the classical Hahn-Banach Theorem (see
\cite[Theorem 6.1.4]{D02}).

$\mathfrak{ag}(\widehat{\mathcal{F}})$ is uniformly dense in $C(\widehat{E};\mathbf{R})$
by Corollary \ref{cor:Base_Fun_Dense} (with $d=1$ and $E_{0}=E$).
For each fixed $h\in C(\widehat{E};\mathbf{R})$, there exist $\{f_{n}\}\subset\mathfrak{ag}(\mathcal{F})$
such that
\begin{equation}
\lim_{n\rightarrow\infty}\Vert h|_{E}-f_{n}\Vert_{\infty}=\lim_{n\rightarrow\infty}\Vert h|_{E}-\widehat{f}_{n}|_{E}\Vert_{\infty}\leq\lim_{n\rightarrow\infty}\Vert h-\widehat{f}_{n}\Vert_{\infty}=0.\label{eq:Check_Rep_Functional_Positive_1}
\end{equation}
$\Phi$ and $\widehat{\Phi}$ are continuous functionals by (\ref{eq:Extend_Functional_Bounded}),
(\ref{eq:Rep_Functional_Bounded}) and \cite[Theorem 6.1.2]{D02}.
So, (\ref{eq:Check_Rep_Functional_Bounded}) and (\ref{eq:Check_Rep_Functional_Positive_1})
imply
\begin{equation}
\widehat{\Phi}(h)=\lim_{n\rightarrow\infty}\widehat{\Phi}(\widehat{f}_{n})=\lim_{n\rightarrow\infty}\Phi(f_{n})=\Phi(h|_{E}).\label{eq:Check_Rep_Functional_Positive_2}
\end{equation}
(\ref{eq:Check_Rep_Functional_Positive_2}) verifies (\ref{eq:Rep_Functional})
and implies $\widehat{\Phi}(h)\geq0$ for all non-negative $h\in C(\widehat{E};\mathbf{R})$
since $\Phi$ is positive, thus proving the positivenss of $\widehat{\Phi}$.

\textit{Step 4: Establish integral representation of the replica functional}.
Since $\widehat{E}$ is a compact Polish space, we apply the classical
Riesz Representation Theorem (see \cite[Theorem 2.1.5]{KX95}) to
$\widehat{\Phi}$ and obtain a $\nu\in\mathcal{M}^{+}(\widehat{E})$
satisfying
\begin{equation}
\widehat{\Phi}(h)=h^{*}(\nu),\;\forall h\in C(\widehat{E};\mathbf{R}).\label{eq:Rep_Riesz_Representation_1}
\end{equation}
It follows by (\ref{eq:Rep_Functional}) and (\ref{eq:Rep_Riesz_Representation_1})
that
\begin{equation}
\Phi(h|_{E})=h^{*}(\nu),\;\forall h\in C(\widehat{E};\mathbf{R}).\label{eq:Rep_Riesz_Representation_2}
\end{equation}
Moreover, it follows by (\ref{eq:Rep_Riesz_Representation_2}), (\ref{eq:Total_Mass_Approx_Seq})
and (\ref{eq:Rep_Riesz_Representation_1}) that
\begin{equation}
\nu(\widehat{E})=\Phi(1)=\lambda_{0}=\lim_{p\rightarrow\infty}\Phi(g_{p})=\lim_{p\rightarrow\infty}\widehat{\Phi}(\widehat{g}_{p})=\lim_{p\rightarrow\infty}\widehat{g}_{p}^{*}(\nu).\label{eq:Nu(Ehat)=00003DTotal_Mass}
\end{equation}

\textit{Step 5: Establish the desired measure $\mu$}. We define 
\begin{equation}
\mathcal{A}\circeq\left\{ g\in C_{c}(E;\mathbf{R}):f_{g,a,b}\in\mathfrak{ag}_{\mathbf{Q}}(\mathcal{F})\mbox{ for some }a,b\in\mathbf{Q}\mbox{ with }ab\neq0\right\} ,\label{eq:Riesz_Support_Fun}
\end{equation}
let $K_{g}\in\mathscr{K}(E)$ denote the closure of $E\backslash g^{-1}(\{0\})$
in $E$ for each $g\in\mathcal{A}$, and have by Corollary \ref{cor:Base_Compact}
(a) (with $d=1$ and $E_{0}=E$) that
\begin{equation}
\{K_{g}\}_{g\in\mathcal{A}}\subset\mathscr{K}(\widehat{E})\subset\mathscr{B}(\widehat{E}).\label{eq:Riesz_Support_Component}
\end{equation}
$\mathfrak{ag}_{\mathbf{Q}}(\mathcal{F})$ is a countable collection
by Fact \ref{fact:ac_mc_Countable}, so $\mathcal{A}$ is also countable
and
\begin{equation}
A\circeq\bigcup_{g\in\mathcal{A}}K_{g}\in\mathscr{K}_{\sigma}^{\mathbf{m}}(E)\cap\mathscr{B}(\widehat{E})\label{eq:Riesz_Support}
\end{equation}
by Corollary \ref{cor:Base_Compact} (b) (with $d=1$ and $E_{0}=E$).
We have $\{g_{p}\}_{p\in\mathbf{N}}\subset\mathcal{A}$ and
\begin{equation}
\widehat{g}_{p}=\overline{g}_{p}=\widehat{g}_{p}\mathbf{1}_{K_{g_{p}}}\leq1,\;\forall p\in\mathbf{N}\label{eq:Check_Riesz_Meas_Concentration_1}
\end{equation}
by (\ref{eq:gp_in_F_B}), (\ref{eq:Rep_Cc_Fun}) (with $g=g_{p}$,
$a=1$ and $b=0$) and (\ref{eq:g_in_B_<=00003D1}).
\begin{equation}
\begin{aligned}\nu(\widehat{E}) & \geq\nu(A)\ge\lim_{p\rightarrow\infty}\nu(K_{g_{p}})\\
 & \geq\lim_{p\rightarrow\infty}\left(\widehat{g}_{p}\mathbf{1}_{K_{g_{p}}}\right)^{*}(\nu)=\lim_{p\rightarrow\infty}\widehat{g}_{p}^{*}(\nu)=\nu(\widehat{E})
\end{aligned}
\label{eq:Check_Riesz_Meas_Concentration_2}
\end{equation}
by (\ref{eq:Riesz_Support_Component}), (\ref{eq:Riesz_Support}),
(\ref{eq:Check_Riesz_Meas_Concentration_1}), (\ref{eq:Rep_Riesz_Representation_1})
and (\ref{eq:Nu(Ehat)=00003DTotal_Mass}). Hence, we have by (\ref{eq:Check_Riesz_Meas_Concentration_2})
and Proposition \ref{prop:Redefine_FDD} that $\nu$ is the replica
of $\mu\circeq(\nu|_{A})|^{E}$ and
\begin{equation}
\mu(E)=\mu(A)=\nu(A)=\nu(\widehat{E})=\lambda_{0}.\label{eq:Riesz_Mass_Concentrate_A}
\end{equation}
Moreover, the $\mathbf{m}$-tightness of $\mu$ follows by (\ref{eq:Riesz_Mass_Concentrate_A})
and the fact $A\in\mathscr{K}_{\sigma}^{\mathbf{m}}(E)$.

\textit{Step 6: Redefine the integral representation of the replica
functional as that of $\varphi$}. First, we find that
\begin{equation}
\varphi(g)=\Phi(g)=\widehat{g}^{*}(\nu)=g^{*}(\mu),\;\forall g\in C_{c}(E;\mathbf{R})\label{eq:Check_Riesz_Representation}
\end{equation}
by (\ref{eq:Extend_Original_Functional}), (\ref{eq:Rep_Cc_Fun})
(with $a=1$ and $b=0$), (\ref{eq:Rep_Riesz_Representation_2}),
the fact $\nu=\overline{\mu}$ and Proposition \ref{prop:RepMeas_Basic}
(e) (with $d=1$ and $E_{0}=E$), thus proving (\ref{eq:Riesz_Representation}).
Secondly, we have
\begin{equation}
\mu(E)=\lambda_{0}=\lim_{p\rightarrow\infty}\varphi(g_{p})=\lim_{p\rightarrow\infty}g_{p}^{*}(\mu)\label{eq:Check_Approx_Total_Mass}
\end{equation}
by (\ref{eq:Riesz_Mass_Concentrate_A}), (\ref{eq:Total_Mass_Approx_Seq})
and (\ref{eq:Check_Riesz_Representation}) (with $g=g_{p}$). Then,
(\ref{eq:Approx_Total_Mass}) follows by (\ref{eq:Check_Approx_Total_Mass})
and (\ref{eq:g_in_B_<=00003D1}).

((b) $\rightarrow$ (a)) $\Lambda(g)\circeq g^{*}(\mu)$ for each
$g\in C_{b}(E;\mathbf{R})$ defines a positive linear functional with
$\Lambda|_{C_{c}(E;\mathbf{R})}=\varphi$. (\ref{eq:Define_Total_Mass})
follows by (\ref{eq:Approx_Total_Mass}).\end{proof}

\begin{cor}
\label{cor:Riesz}Let $E$ be a $C_{c}(E;\mathbf{R})$-baseable space,
$\varphi$ be a linear functional on $C_{0}(E;\mathbf{R})$ and $\mathcal{B}$
be as in (\ref{eq:Unit_Ball}). Then, the following statements are
equivalent:

\renewcommand{\labelenumi}{(\alph{enumi})}
\begin{enumerate}
\item $\varphi$ is continuous and there exists a positive linear functional
$\Lambda$ on $C_{b}(E;\mathbf{R})$ satisfying (\ref{eq:Original_Functional_Dominated})
and (\ref{eq:Define_Total_Mass}).
\item There exists an $\mathbf{m}$-tight $\mu\in\mathcal{M}^{+}(E)$ satisfying
(\ref{eq:Approx_Total_Mass}) and
\begin{equation}
\varphi(g)=g^{*}(\mu),\;\forall g\in C_{0}(E;\mathbf{R}).\label{eq:Riesz_Representation_C0}
\end{equation}

\end{enumerate}
\end{cor}
\begin{proof}
((a) $\rightarrow$ (b)) There exists an $\mathbf{m}$-tight $\mu\in\mathcal{M}^{+}(E)$
satisfying (\ref{eq:Approx_Total_Mass}) and (\ref{eq:Riesz_Representation})
by Theorem \ref{thm:Riesz}. $C_{c}(E;\mathbf{R})$ is uniformly dense
in $C_{0}(E;\mathbf{R})$ by Fact \ref{fact:Baseable_Metrizable_Separable}
(a) and Proposition \ref{prop:Cc_Lattice} (b). Hence, (\ref{eq:Riesz_Representation_C0})
follows by (\ref{eq:Riesz_Representation}), the continuity of $\varphi$
and the Dominated Convergence Theorem.

((b) $\rightarrow$ (a)) The functional defined by $\Lambda(g)\circeq g^{*}(\mu)$
for each $g\in C_{b}(E;\mathbf{R})$ satisfies $\Lambda|_{C_{0}(E;\mathbf{R})}=\varphi$,
has linearity and is continuous by the Dominated Convergence Theorem.
Moreover, (\ref{eq:Define_Total_Mass}) is immediate by (\ref{eq:Approx_Total_Mass}).\end{proof}

\subsection{\label{sub:Sko_Rep}Almost-sure representation of weak convergence}

We now generalize the Skorokhod Representation Theorem in \cite{J97a}.
Commonly, the Skorokhod Representation Theorem is established on separable
metric spaces. \cite[Theorem 2]{J97a} extended this result to sequences
of tight probability measures on baseable spaces%
\footnote{While \cite{J97a} did not use the term ``baseable'', he did assume
point-separability by countably many continuous functions.%
}. $\mathbf{m}$-tightness is equivalent to tightness in a baseable
space $E$ by Corollary \ref{cor:Baseable_MC} (a). Hence, the conditions
of the following theorem are strictly milder than those in \cite{J97a}.
\begin{thm}
\label{thm:Sko_Rep}Let $E$ be a topological space, $C(E;\mathbf{R})$
separate points on $E$,
\begin{equation}
\mu_{n}\Longrightarrow\mu_{0}\mbox{ as }n\uparrow\infty\mbox{ in }\mathcal{P}(E),\label{eq:Mu_n_WC_Mu0_P(E)}
\end{equation}
and $\{\mu_{n}\}_{n\in\mathbf{N}}$ be $\mathbf{m}$-tight. Then,
there exist $E$-valued random variables $\{\xi_{n}\}_{n\in\mathbf{N}_{0}}$
defined on the same probability space such that $\xi_{n}$ has distribution
$\mu_{n}$ for all $n\in\mathbf{N}_{0}$ and $\{\xi_{n}\}_{n\in\mathbf{N}}$
converges to $\xi_{0}$ as $n\uparrow\infty$ almost surely.
\end{thm}
\begin{proof}
$\{\mu_{n}\}_{n\in\mathbf{N}_{0}}$ is $\mathbf{m}$-tight by Lemma
\ref{lem:Compact_Portmanteau} (b) (with $\Gamma=\{\mu_{n}\}_{n\in\mathbf{N}}$).
There exists a base $(E_{0},\mathcal{F};\widehat{E},\widehat{\mathcal{F}})$
such that $\{\mu_{n}\}_{n\in\mathbf{N}_{0}}$ is tight in
\begin{equation}
E_{0}\in\mathscr{K}_{\sigma}^{\mathbf{m}}(E)\cap\mathscr{B}^{\mathbf{s}}(E)\cap\mathscr{B}(E)\cap\mathscr{B}(\widehat{E})\label{eq:Sko_Rep_E0}
\end{equation}
by Lemma \ref{lem:m-Tight_Base} (with $\Gamma_{i}=\{\mu_{n}\}_{n\in\mathbf{N}_{0}}$
and $\mathcal{D}=C(E;\mathbf{R})$) and Corollary \ref{cor:Base_Compact}
(b).
\begin{equation}
\inf_{n\in\mathbf{N}_{0}}\mu_{n}(E_{0})=\inf_{n\in\mathbf{N}_{0}}\overline{\mu}_{n}(E_{0})=1\label{eq:E0_Common_Support}
\end{equation}
by the tightness of $\{\mu_{n}\}_{n\in\mathbf{N}_{0}}$ in $E_{0}$
and Proposition \ref{prop:RepMeas_Basic} (a). Furthermore,
\begin{equation}
\overline{\mu}_{n}\Longrightarrow\overline{\mu}_{0}\mbox{ as }k\uparrow\infty\mbox{ in }\mathcal{P}(\widehat{E})\label{eq:MuBar_nk_WC_Mu_Bar_P(Ehat)}
\end{equation}
by (\ref{eq:Mu_n_WC_Mu0_P(E)}) and Proposition \ref{prop:WC_Push_Forward}
(a, e) (with $d=1$, $\mu_{n}^{\prime}=\mu_{n}$ and $\mu^{\prime}=\mu_{0}$).

$\widehat{E}$ is a Polish space by Lemma \ref{lem:Base} (c), so
the classical Skorokhod Representation Theorem (see \cite[Theorem 11.7.2]{D02})
is applicable to $\{\overline{\mu}_{n}\}_{n\in\mathbf{N}_{0}}$, yielding
random variables $\{\overline{\xi}_{n}\}_{n\in\mathbf{N}_{0}}\subset M(\Omega,\mathscr{F},\mathbb{P};\widehat{E})$
that satisfy $\mathbb{P}\circ\overline{\xi}_{n}^{-1}=\overline{\mu}_{n}$
for all $n\in\mathbf{N}_{0}$ and
\begin{equation}
\mathbb{P}\left(\overline{\xi}_{n}\longrightarrow\overline{\xi}_{0}\mbox{ as }n\uparrow\infty\right)=1.\label{eq:Rep_Sko_Rep}
\end{equation}

Singletons are Borel sets in $\widehat{E}$ by Proposition \ref{prop:Separability}
(a, b). Hence, there exist
\begin{equation}
\xi_{n}\in M\left(\Omega,\mathscr{F};E_{0},\mathscr{B}_{\widehat{E}}(E_{0})\right),\;\forall n\in\mathbf{N}_{0}\label{eq:Check_Sko_Rep_Redef_RV_Measurable_1}
\end{equation}
such that
\begin{equation}
\inf_{n\in\mathbf{N}_{0}}\mathbb{P}\left(\xi_{n}=\overline{\xi}_{n}\right)\geq\inf_{n\in\mathbf{N}_{0}}\mathbb{P}\left(\overline{\xi}_{n}\in E_{0}\right)=\inf_{n\in\mathbf{N}_{0}}\overline{\mu}_{n}(E_{0})=1\label{eq:Sko_Rep_Redef_RV_Rep}
\end{equation}
by (\ref{eq:E0_Common_Support}) and Fact \ref{fact:var(f)} (b) (with
$(S,\mathscr{A})=(\Omega,\mathscr{F})$, $(E,\mathscr{U})=(\widehat{E},\mathscr{B}(\widehat{E}))$,
$A=E_{0}$ and $f=\overline{\xi}_{n}$).

Now, we have by (\ref{eq:Rep_Sko_Rep}) and (\ref{eq:Sko_Rep_Redef_RV_Rep})
that
\begin{equation}
\mathbb{P}\left(\xi_{n}\longrightarrow\xi_{0}\mbox{ as }n\uparrow\infty\right)\geq\mathbb{P}\left(\overline{\xi}_{n}\longrightarrow\overline{\xi}_{0}\mbox{ as }n\uparrow\infty\right)=1.\label{eq:Check_Sko_Rep}
\end{equation}
(\ref{eq:Sko_Rep_E0}) gives $E_{0}\in\mathscr{B}^{\mathbf{s}}(E)$.
It then follows by Lemma \ref{lem:SB_Base} (a) (with $d=1$ and $A=E_{0}$)
and (\ref{eq:Check_Sko_Rep_Redef_RV_Measurable_1}) that
\begin{equation}
\xi_{n}\in M\left(\Omega,\mathscr{F};E_{0},\mathscr{B}_{E}(E_{0})\right)\subset M\left(\Omega,\mathscr{F};E\right),\;\forall n\in\mathbf{N}_{0}.\label{eq:Check_Sko_Rep_Redef_RV_Measurable_2}
\end{equation}
\end{proof}

\begin{rem}
\label{rem:Generalize_J97}$([0,1]^{[0,1]},\Vert\cdot\Vert_{\infty})$
mentioned in Example \ref{exp:Metrizable_Compact_non-Baseable} is
non-baseable. Compact subsets of this normed space are automatically
metrizable. This space is Tychonoff by Proposition \ref{prop:CR_Space}
(a) and so its points are separated by $C([0,1]^{[0,1]},\Vert\cdot\Vert_{\infty};\mathbf{R})$
by Proposition \ref{prop:CR} (a, b). Theorem \ref{thm:Sko_Rep} applies
in this case whereas \cite[Theorem 2]{J97a} does not.
\end{rem}

\chapter{\label{chap:Rep_Proc}Replication of Stochastic Process}

\chaptermark{Replica Process}

This chapter is devoted to the replication of $E$-valued stochastic
process via a base $(E_{0},\mathcal{F};\widehat{E},\widehat{\mathcal{F}})$
over $E$. \S \ref{sec:RepProc} introduces and discusses the basic
properties of replica process. \S \ref{sec:Cadlag_RepProc} focuses
on the special case of c$\grave{\mbox{a}}$dl$\grave{\mbox{a}}$g
replica. Properties like tightness and relative compactness are simple
to verify or even automatic on the compact Polish space $\widehat{E}$.
\S \ref{sec:RepProc_FC} associates the finite-dimensional convergence
of general processes to that of their general replicas. \S \ref{sec:RepProc_Path_Space}
discusses tightness and weak convergence of c$\grave{\mbox{a}}$dl$\grave{\mbox{a}}$g
replicas as path-space-valued random variables. Finally, \S \ref{sec:Base_Proc}
considers when a family of processes can be contained in a baseable
set to perform the desired replication. If necessary, the readers
are referred to \S \ref{sec:Proc} where we specify our terminologies
and notations about stochastic processes.

\section{\label{sec:RepProc}Introduction to replica process}

\subsection{\label{sub:RepProc}Definition}

Given a base $(E_{0},\mathcal{F};\widehat{E},\widehat{\mathcal{F}})$
over topological space $E$, a replica of $E$-valued process $X$
is a related process that takes values in the compact Polish space
$\widehat{E}$. Since $X$ and its replicas may live in different
spaces, they are related by the mappings $\bigotimes\mathcal{F}$
and $\bigotimes\widehat{\mathcal{F}}$ rather than their own values.
\begin{defn}
\label{def:RepProc}Let $E$ be a topological space and $(\Omega,\mathscr{F},\mathbb{P};X)$%
\footnote{``$(\Omega,\mathscr{F},\mathbb{P};X)$'' as defined in \S \ref{sec:RV}
means an $E$-valued random variable or process $X$ defined on probability
space $(\Omega,\mathscr{F},\mathbb{P})$. We imposed in \S \ref{sec:Convention}
that the probability space $(\Omega,\mathscr{F},\mathbb{P})$ is complete.%
} be an $E$-valued process. With respect to a base $(E_{0},\mathcal{F};\widehat{E},\widehat{\mathcal{F}})$
over $E$, an $\widehat{E}$-valued process $(\Omega,\mathscr{F},\mathbb{P};\widehat{X})$
is said to be a \textbf{replica of $X$} if
\begin{equation}
\mathbb{P}\left(\bigotimes\widehat{\mathcal{F}}\circ\widehat{X}_{t}=\bigotimes\mathcal{F}\circ X_{t}\right)\geq\mathbb{P}\left(\bigotimes\mathcal{F}\circ X_{t}\in\bigotimes\widehat{\mathcal{F}}(\widehat{E})\right),\;\forall t\in\mathbf{R}^{+}.\label{eq:Define_RepProc}
\end{equation}
\end{defn}
\begin{note}
\label{note:Mutiple_Replica}An $E$-valued process $X$ may have
multiple replicas, divisible into equivalence classes by indistinguishability%
\footnote{The terminology ``indistinguishability'' was explained in \S \ref{sec:Proc}.%
}.
\end{note}
We make the following notations for simplicity.
\begin{notation}
\label{notation:RepProc}Let $E$ be a topological space, $(E_{0},\mathcal{F};\widehat{E},\widehat{\mathcal{F}})$
be a base over $E$ and $X$ be an $E$-valued process.
\begin{itemize}
\item $\mathfrak{rep}(X;E_{0},\mathcal{F})$%
\footnote{``$\mathfrak{rep}$'' is ``rep'' in fraktur font which stands
for ``replica''.%
} denotes the family of all equivalence classes of $X$'s replicas
with respect to $(E_{0},\mathcal{F};\widehat{E},\widehat{\mathcal{F}})$
under the equivalence relation of indistinguishability.
\item $\mathfrak{rep}_{\mathrm{m}}(X;E_{0},\mathcal{F})$, $\mathfrak{rep}_{\mathrm{p}}(X;E_{0},\mathcal{F})$
and $\mathfrak{rep}_{\mathrm{c}}(X;E_{0},\mathcal{F})$ denote the
measurable%
\footnote{The notions of measurable process and progressive process were reviewed
in \S \ref{sec:Proc}.%
}, progressive and c$\grave{\mbox{a}}$dl$\grave{\mbox{a}}$g members
of $\mathfrak{rep}(X;E_{0},\mathcal{F})$, respectively.
\item $\widehat{X}=\mathfrak{rep}(X;E_{0},\mathcal{F})$ means $\mathfrak{rep}(X;E_{0},\mathcal{F})$
is the single equivalence class $\{\widehat{X}\}$. Similar notations
apply to the above-mentioned subfamilities of $\mathfrak{rep}(X;E_{0},\mathcal{F})$.
\end{itemize}
\end{notation}
\begin{rem}
\label{rem:RepBase}The notation ``$\mathfrak{rep}(X;E_{0},\mathcal{F})$''
merely specifies the first two components $(E_{0},\mathcal{F})$ of
the base $(E_{0},\mathcal{F};\widehat{E},\widehat{\mathcal{F}})$
since Corollary \ref{cor:Base_Unique} and Theorem \ref{thm:Base}
showed that this base is totally determined by $(E_{0},\mathcal{F})$.\end{rem}
\begin{note}
\label{note:Ehat_Valued_Proc_FDD}$\,$
\begin{itemize}
\item $\mathbf{R}$, $\mathbf{R}^{\infty}$, $\widehat{E}$ and $\bigotimes\widehat{\mathcal{F}}(\widehat{E})$
(as a subspace of $\mathbf{R}^{\infty}$) are Polish spaces by Proposition
\ref{prop:Var_Polish} (f), Lemma \ref{lem:Base} (c) and (\ref{eq:Base_Imb}).
\item $D(\mathbf{R}^{+};\mathbf{R})$, $D(\mathbf{R}^{+};\mathbf{R}^{\infty})$,
$D(\mathbf{R}^{+};\widehat{E})$ and $D(\mathbf{R}^{+};\bigotimes\widehat{\mathcal{F}}(\widehat{E}))$
are well-defined Polish spaces by Proposition \ref{prop:Sko_Basic_2}
(d).
\item $\mathscr{B}(\widehat{E}^{d})=\mathscr{B}(\widehat{E})^{\otimes d}$
for all $d\in\mathbf{N}$ by Fact \ref{fact:Ehat_Borel_Prod}, so
finite-dimensional distributions of any $\widehat{E}$-valued process
(especially any replica process) are all Borel probability measures%
\footnote{While the finite-dimensional distributions of any Polish-space-valued
process are Borel, this is not necessarily true for general processes
(see \S \ref{sec:Proc}).%
}.
\end{itemize}
\end{note}
The following proposition justifies the general existence of replica
proecsses.
\begin{prop}
\label{prop:RepProc_Exist}Let $E$ be a topological space, $(E_{0},\mathcal{F};\widehat{E},\widehat{\mathcal{F}})$
be a base over $E$ and $(\Omega,\mathscr{F},\mathbb{P};X)$ be an
$E$-valued process. Then:

\renewcommand{\labelenumi}{(\alph{enumi})}
\begin{enumerate}
\item $\mathfrak{rep}(X;E_{0},\mathcal{F})$ is non-empty.
\item If $X$ is a measurable process, then $\mathfrak{rep}_{\mathrm{m}}(X;E_{0},\mathcal{F})$
is non-empty.
\end{enumerate}
\end{prop}
\begin{proof}
(a) We fix $x_{0}\in E_{0}$, let $\varphi$ be the identity mapping
on $\mathbf{R}^{\infty}$ and define%
\footnote{``$\mathfrak{var}(\cdot)$'' was introduced in Notation \ref{notation:Var}.%
}
\begin{equation}
\varphi_{x_{0}}\circeq\mathfrak{var}\left(\varphi;\mathbf{R}^{\infty},\bigotimes\widehat{\mathcal{F}}(\widehat{E}),\bigotimes\mathcal{F}(x_{0})\right).\label{eq:Define_RepProc_Gen_Map}
\end{equation}
$\bigotimes\widehat{\mathcal{F}}(\widehat{E})\in\mathscr{B}(\mathbf{R}^{\infty})$
by (\ref{eq:Prod(Fhat)(Ehat)_Compact_Rinf}). $\mathbf{R}^{\infty}$
is a Polish space, so $\{\bigotimes\mathcal{F}(x_{0})\}\in\mathscr{B}(\mathbf{R}^{\infty})$
by Proposition \ref{prop:Separability} (a, b). We then have that
\begin{equation}
\varphi_{x_{0}}(y)=y,\;\forall y\in\bigotimes\widehat{\mathcal{F}}(\widehat{E})\label{eq:phi_x0_Identity_Prod(Fhat)(Ehat)}
\end{equation}
and
\begin{equation}
\varphi_{x_{0}}\in M\left(\mathbf{R}^{\infty};\bigotimes\widehat{\mathcal{F}}(\widehat{E})\right)\label{eq:phi_x0_Measurable}
\end{equation}
by Fact  \ref{fact:var(f)} (b) (with $(S,\mathscr{A})=(E,\mathscr{U})=(\mathbf{R}^{\infty},\mathscr{B}(\mathbf{R}^{\infty}))$,
$A=\bigotimes\widehat{\mathcal{F}}(\widehat{E})$, $f=\varphi$ and
$y_{0}=\bigotimes\mathcal{F}(x_{0})$). It follows that
\begin{equation}
\left(\bigotimes\widehat{\mathcal{F}}\right)^{-1}\circ\varphi_{x_{0}}\circ\bigotimes\mathcal{F}\in M(E;\widehat{E})\label{eq:General_RepProc_Map}
\end{equation}
by (\ref{eq:Base_Imb}), (\ref{eq:phi_x0_Measurable}) and Lemma \ref{lem:Base}
(e). Hence,
\begin{equation}
\widehat{X}\circeq\varpi\left[\left(\bigotimes\widehat{\mathcal{F}}\right)^{-1}\circ\varphi_{x_{0}}\circ\bigotimes\mathcal{F}\right]\circ X\label{eq:Define_General_RepProc}
\end{equation}
well defines an $\widehat{E}$-valued process $(\Omega,\mathscr{F},\mathbb{P};\widehat{X})$
by Fact \ref{fact:Proc_Path_Mapping} (a) (with $S=\widehat{E}$ and
$f=(\bigotimes\widehat{\mathcal{F}})^{-1}\circ\varphi_{x_{0}}\circ\bigotimes\mathcal{F}$).
It follows by (\ref{eq:phi_x0_Identity_Prod(Fhat)(Ehat)}) and (\ref{eq:Define_General_RepProc})
that%
\footnote{By (\ref{eq:Define_RepProc_Gen_Map}), $\varphi_{x_{0}}\circ\bigotimes\mathcal{F}\circ X_{t}\in\bigotimes\widehat{\mathcal{F}}(\widehat{E})$
might not imply $\bigotimes\mathcal{F}\circ X_{t}\in\bigotimes\widehat{\mathcal{F}}(\widehat{E})$
when $\bigotimes\mathcal{F}(x_{0})\in\bigotimes\widehat{\mathcal{F}}(\widehat{E})$. %
}
\begin{equation}
\begin{aligned} & \mathbb{P}\left(\bigotimes\mathcal{F}\circ X_{t}\in\bigotimes\widehat{\mathcal{F}}(\widehat{E})\right)\\
 & \leq\mathbb{P}\left(\bigotimes\widehat{\mathcal{F}}\circ\left(\bigotimes\widehat{\mathcal{F}}\right)^{-1}\circ\varphi_{x_{0}}\circ\bigotimes\mathcal{F}\circ X_{t}\in\bigotimes\widehat{\mathcal{F}}(\widehat{E})\right)\\
 & \leq\mathbb{P}\left(\bigotimes\widehat{\mathcal{F}}\circ\widehat{X}_{t}=\bigotimes\mathcal{F}\circ X_{t}\right),\;\forall t\in\mathbf{R}^{+},
\end{aligned}
\label{eq:Check_RepProc}
\end{equation}
thus proving $\widehat{X}\in\mathfrak{rep}(X;E_{0},\mathcal{F})$
by (\ref{eq:Define_RepProc}).

(b) Let $\widehat{X}$ be as above and define $\xi(t,\omega)\circeq X_{t}(\omega)$
and $\widehat{\xi}(t,\omega)\circeq\widehat{X}_{t}(\omega)$ for each
$(t,\omega)\in\mathbf{R}^{+}\times\Omega$. If $X$ is a measurable
process, then
\begin{equation}
\widehat{\xi}=\left(\bigotimes\widehat{\mathcal{F}}\right)^{-1}\circ\varphi_{x_{0}}\circ\bigotimes\mathcal{F}\circ\xi\in M\left(\mathbf{R}^{+}\times\Omega,\mathscr{B}(\mathbf{R}^{+})\otimes\mathscr{F};\widehat{E},\mathscr{B}(\widehat{E})\right)\label{eq:Check_RepProc_Measurable}
\end{equation}
by (\ref{eq:General_RepProc_Map}), thus proving $\widehat{X}\in\mathfrak{rep}_{\mathrm{m}}(X;E_{0},\mathcal{F})$.\end{proof}

\subsection{\label{sub:RepProc_Association}Association with the original process}

The next two propositions expose the connection between the original
and replica processes using properties of the base.
\begin{prop}
\label{prop:RepProc_T-Base}Let $E$ be a topological space, $(E_{0},\mathcal{F};\widehat{E},\widehat{\mathcal{F}})$
be a base over $E$, $\mathbf{T}\subset\mathbf{R}^{+}$, $(\Omega,\mathscr{F},\mathbb{P};X)$
be an $E$-valued process satisfying
\begin{equation}
\inf_{t\in\mathbf{T}}\mathbb{P}\left(X_{t}\in E_{0}\right)=1\label{eq:T-Base}
\end{equation}
and $\widehat{X},\widehat{X}^{1},\widehat{X}^{2}\in\mathfrak{rep}(X;E_{0},\mathcal{F})$.
Then:

\renewcommand{\labelenumi}{(\alph{enumi})}
\begin{enumerate}
\item $\widehat{X}$ satisfies
\begin{equation}
\inf_{t\in\mathbf{T}}\mathbb{P}\left(X_{t}=\widehat{X}_{t}\in E_{0}\right)=1.\label{eq:RepProc_(T,E0)-Mod}
\end{equation}
Moreover, $\widehat{X}^{1}$ and $\widehat{X}^{2}$ satisfy
\begin{equation}
\inf_{t\in\mathbf{T}}\mathbb{P}\left(\widehat{X}_{t}^{1}=\widehat{X}_{t}^{2}\in E_{0}\right)=1.\label{eq:RepProc_(T,E0)-Mod_Uni}
\end{equation}

\item $\mathbb{P}\circ\widehat{X}_{\mathbf{T}_{0}}^{-1}$ is the replica
measure of $\mathbb{P}\circ X_{\mathbf{T}_{0}}^{-1}$ for all $\mathbf{T}_{0}\in\mathscr{P}_{0}(\mathbf{T})$%
\footnote{$\widehat{X}_{\mathbf{T}_{0}}$, the section of $\widehat{X}$ for
$\mathbf{T}_{0}$ was defined in \S \ref{sec:Proc}.%
}.
\end{enumerate}
\end{prop}
\begin{proof}
(a) (\ref{eq:RepProc_(T,E0)-Mod}) follows by (\ref{eq:Define_RepProc})
and Fact \ref{fact:TF_(T,E0)-Mod} (with $Y=\widehat{X}$). (\ref{eq:RepProc_(T,E0)-Mod_Uni})
is immediate by (\ref{eq:RepProc_(T,E0)-Mod}) (with $\widehat{X}=\widehat{X}^{1}$
or $\widehat{X}^{2}$).

(b) Note \ref{note:Ehat_Valued_Proc_FDD} implies that $\mathbb{P}\circ\widehat{X}_{\mathbf{T}_{0}}^{-1}\in\mathcal{P}(\widehat{E}^{\mathbf{T}_{0}})$.
Moreover, we have by (a) that
\begin{equation}
\mathbb{P}\left(\widehat{X}_{\mathbf{T}_{0}}\in A\right)=\mathbb{P}\left(X_{\mathbf{T}_{0}}\in A\cap E_{0}^{\mathbf{T}_{0}}\right),\;\forall A\in\mathscr{B}(\widehat{E}^{\mathbf{T}_{0}})\label{eq:Check_Rep_Proc_Rep_Meas}
\end{equation}
and
\begin{equation}
\mathbb{P}\left(X_{\mathbf{T}_{0}}=\widehat{X}_{\mathbf{T}_{0}}\in E_{0}^{\mathbf{T}_{0}}\right)=1.\label{eq:RepProc_FDD_T_Mod}
\end{equation}
\end{proof}

\begin{prop}
\label{prop:RepProc_FR-Base}Let $E$ be a topological space, $(E_{0},\mathcal{F};\widehat{E},\widehat{\mathcal{F}})$
be a base over $E$, $\mathbf{T}\subset\mathbf{R}^{+}$ and $(\Omega,\mathscr{F},\mathbb{P};X)$
be an $E$-valued process satisfying
\begin{equation}
\inf_{t\in\mathbf{T}}\mathbb{P}\left(\bigotimes\mathcal{F}\circ X_{t}\in\bigotimes\widehat{\mathcal{F}}(\widehat{E})\right)=1.\label{eq:FR-Base}
\end{equation}
Then:

\renewcommand{\labelenumi}{(\alph{enumi})}
\begin{enumerate}
\item Any $\widehat{X}\in\mathfrak{rep}(X;E_{0},\mathcal{F})$ satisfies
\begin{equation}
\mathbb{P}\left(f\circ X_{\mathbf{T}_{0}}=\widehat{f}\circ\widehat{X}_{\mathbf{T}_{0}}\right)=1,\;\forall f\in\mathfrak{ca}\left[\Pi^{\mathbf{T}_{0}}(\mathcal{F})\right],\mathbf{T}_{0}\in\mathscr{P}_{0}(\mathbf{T})\label{eq:RepProc_TF_T-Mod}
\end{equation}
and
\begin{equation}
\mathbb{E}\left[f\circ X_{\mathbf{T}_{0}}\right]=\mathbb{E}\left[\widehat{f}\circ\widehat{X}_{\mathbf{T}_{0}}\right],\;\forall f\in\mathfrak{ca}\left[\Pi^{\mathbf{T}_{0}}(\mathcal{F})\right],\mathbf{T}_{0}\in\mathscr{P}_{0}(\mathbf{T}).\label{eq:RepProc_Int_Rep}
\end{equation}

\item If $\mathbf{T}\subset\mathbf{R}^{+}$ is dense, then $\mathfrak{rep}_{\mathrm{c}}(X;E_{0},\mathcal{F})$
is at most a singleton.
\end{enumerate}
\end{prop}
\begin{proof}
(a) Any $\widehat{X}\in\mathfrak{rep}(X;E_{0},\mathcal{F})$ satisfies
\begin{equation}
\begin{aligned} & \inf_{t\in\mathbf{T}}\mathbb{P}\left(\bigotimes\widehat{\mathcal{F}}\circ\widehat{X}_{t}=\bigotimes\mathcal{F}\circ X_{t}\right)\\
 & \geq\inf_{t\in\mathbf{T}}\mathbb{P}\left(\bigotimes\mathcal{F}\circ X_{t}\in\bigotimes\widehat{\mathcal{F}}(\widehat{E})\right)=1
\end{aligned}
\label{eq:Check_RepProc_TF_T-Mod}
\end{equation}
by (\ref{eq:Define_RepProc}) and (\ref{eq:FR-Base}). Now, (\ref{eq:RepProc_TF_T-Mod})
follows by (\ref{eq:Check_RepProc_TF_T-Mod}) and properties of uniform
convergence. (\ref{eq:RepProc_Int_Rep}) is immediate by (\ref{eq:RepProc_TF_T-Mod}).

(b) $\mathbf{T}$ must have a countable subset $\mathbf{T}_{0}$ being
dense in $\mathbf{R}^{+}$. $\varpi(\bigotimes\widehat{\mathcal{F}})$%
\footnote{The notations ``$\varpi(f)$'' and ``$\varpi(\bigotimes\widehat{\mathcal{F}})$''
were defined in \S \ref{sub:Map}.%
} is injective on $D(\mathbf{R}^{+};\widehat{E})$ by Lemma \ref{lem:Base}
(a) and Fact \ref{fact:Path_Mapping_Injective} (with $E=A=\widehat{E}$
and $\mathcal{D}=\widehat{\mathcal{F}}$). Given any $\widehat{X}^{1},\widehat{X}^{2}\in\mathfrak{rep}_{\mathrm{c}}(X;E_{0},\mathcal{F})$,
$\{\varpi(\bigotimes\widehat{\mathcal{F}})\circ\widehat{X}\}_{i=1,2}$
are c$\grave{\mbox{a}}$dl$\grave{\mbox{a}}$g processes by (\ref{eq:Base_Imb})
and Fact \ref{fact:Cadlag_Proc} (a) (with $E=\widehat{E}$, $S=\mathbf{R}^{\infty}$
and $f=\bigotimes\widehat{\mathcal{F}}$), and
\begin{equation}
\begin{aligned}\mathbb{P}\left(\widehat{X}^{1}=\widehat{X}^{2}\right) & =\mathbb{P}\left(\varpi(\bigotimes\widehat{\mathcal{F}})\circ\widehat{X}^{1}=\varpi(\bigotimes\widehat{\mathcal{F}})\circ\widehat{X}^{2}\right)\\
 & \geq\mathbb{P}\left(\bigotimes\widehat{\mathcal{F}}\circ\widehat{X}_{t}^{1}=\bigotimes\widehat{\mathcal{F}}\circ\widehat{X}_{t}^{2},\forall t\in\mathbf{T}_{0}\right)=1
\end{aligned}
\label{eq:Check_Cad_Rep_Ind}
\end{equation}
by (a) (with $\widehat{X}=\widehat{X}^{i}$), Proposition \ref{prop:Proc_Basic_2}
(g) and the injectiveness of $\varpi(\bigotimes\widehat{\mathcal{F}})$.\end{proof}

The following consequence of (\ref{eq:F_Fhat_Coincide}) is apparent
but indispensable.
\begin{fact}
\label{fact:T-Base_FR-Base}(\ref{eq:T-Base}) is stronger than (\ref{eq:FR-Base}).
\end{fact}

\subsection{\label{sub:P(E)_RV_Rep}Application to replicating measure-valued
processes}
\begin{lem}
\label{lem:M(E)_RV_Rep}Let $E$ be a topological space, $(E_{0},\mathcal{F};\widehat{E},\widehat{\mathcal{F}})$
be a base over $E$, $\varphi\circeq\bigotimes\mathfrak{mc}(\mathcal{F})^{*}$%
\footnote{The notation ``$\mathfrak{mc}(\mathcal{F})^{*}$'' was specified
in \S \ref{sec:Borel_Measure}.%
}, $\widehat{\varphi}\circeq\bigotimes\mathfrak{mc}(\widehat{\mathcal{F}})^{*}$,
$S_{0}\in\mathscr{B}(\mathbf{R}^{\infty})$ be contained in $\widehat{\varphi}[\mathcal{M}^{+}(\widehat{E})]$,
$y_{0}\in S_{0}$, $\phi$ be the identity mapping on $\mathbf{R}^{\infty}$
and $X\in M(\Omega,\mathscr{F};\mathcal{M}^{+}(E))$%
\footnote{$X\in M(\Omega,\mathscr{F};\mathcal{M}^{+}(E))$ means $X$ is a finite
Borel random measure on $E$. %
} satisfy
\begin{equation}
\mathbb{P}\left(\varphi\circ X\in S_{0}\right)=1.\label{eq:FR-Base_M(E)}
\end{equation}
Then:

\renewcommand{\labelenumi}{(\alph{enumi})}
\begin{enumerate}
\item $\Psi\circeq\widehat{\varphi}^{-1}\circ\mathfrak{var}(\phi;\mathbf{R}^{\infty},S_{0},y_{0})$
satisfies
\begin{equation}
\Psi\in M\left[\mathbf{R}^{\infty};\widehat{\varphi}^{-1}(S_{0}),\mathscr{O}_{\mathcal{M}^{+}(\widehat{E})}\left(\widehat{\varphi}^{-1}(S_{0})\right)\right].\label{eq:M(E)_RV_Rep_Map}
\end{equation}

\item $Y\circeq\Psi\circ\varphi\circ X\in M(\Omega,\mathscr{F};\mathcal{M}^{+}(\widehat{E}))$
satisfies
\begin{equation}
\mathbb{P}\left(f^{*}\circ X=\widehat{f}^{*}\circ Y\right)=1,\;\forall f\in\mathfrak{ca}(\mathcal{F}).\label{eq:M(E)_RV_TF_Rep}
\end{equation}

\item If
\begin{equation}
\left\{ \omega\in\Omega:X(\omega)(E\backslash E_{0})>0\right\} \in\mathscr{N}(\mathbb{P}),\label{eq:Almostsure_Support_E0}
\end{equation}
then $Y(\omega)$ equals the replica (measure) of $X(\omega)$ for
$\mathbb{P}$ almost all $\omega\in\Omega$.
\item If $A\in\mathscr{B}^{\mathbf{s}}(E)$ satisfies $A\subset E_{0}$
and
\begin{equation}
\left\{ \omega\in\Omega:X(\omega)(E\backslash A)>0\right\} \in\mathscr{N}(\mathbb{P}),\label{eq:Almostsure_Support_A}
\end{equation}
then $(hf)^{*}\circ X$%
\footnote{The $k$-dimensional integration function $(hf)^{*}$ was defined
in \S \ref{sec:Borel_Measure}.%
} and $(\overline{h\mathbf{1}_{A}}\widehat{f})^{*}\circ Y$%
\footnote{$\overline{h\mathbf{1}_{A}}$ denotes the function $\mathfrak{var}(h1_{A};E,A,0)$.%
} belong to $M(\Omega,\mathscr{F};\mathbf{R}^{k})$ and satisfy
\begin{equation}
\mathbb{P}\left((hf)^{*}\circ X=(\overline{h\mathbf{1}_{A}}\widehat{f})^{*}\circ Y\right)=1\label{eq:M(E)_RV_TF_bar_Rep}
\end{equation}
for all $f\in\mathfrak{ca}(\mathcal{F})$, $h\in M_{b}(E;\mathbf{R}^{k})$
and $k\in\mathbf{N}$.
\end{enumerate}
\end{lem}
\begin{rem}
\label{rem:fstar_May_Not_Measurable}Every $f\in C_{b}(E;\mathbf{R})$
satisfies $f^{*}\in C(\mathcal{M}^{+}(E);\mathbf{R})$ by the definition
of weak topology and so $f^{*}\circ X\in M(\Omega,\mathscr{F};\mathbf{R})$.
For $f\in M_{b}(E;\mathbf{R})$, however, $f^{*}$ does not necessarily
belong to $M(\mathcal{M}^{+}(E);\mathbf{R})$ in general, nor is $f^{*}\circ X$
always a random variable (see e.g. Example \ref{exp:P(E)_Hausdorff_1}).
\end{rem}
\begin{proof}
[Proof of Lemma \ref{lem:M(E)_RV_Rep}](a) $\mathfrak{mc}(\mathcal{F})^{*}$
(resp. $\mathfrak{mc}(\widehat{\mathcal{F}})^{*}$) is a countable
subset of $C(\mathcal{M}^{+}(E);\mathbf{R})$ (resp. $C(\mathcal{M}^{+}(\widehat{E});\mathbf{R})$)
by Fact \ref{fact:ac_mc_Countable} (with $E=E$ or $\widehat{E}$,
$\mathcal{D}=\mathcal{F}$ or $\widehat{\mathcal{F}}$ and $d=k=1$),
Definition \ref{def:Base} and Lemma \ref{lem:Base} (a). Then, we
have by Fact \ref{fact:Prod_Map_2} (b) that
\begin{equation}
\varphi\in C\left(\mathcal{M}^{+}(E);\mathbf{R}^{\infty}\right).\label{eq:Check_M(E)_Rep_1}
\end{equation}
At the same time, we have that
\begin{equation}
\widehat{\varphi}\in\mathbf{imb}\left(\mathcal{M}^{+}(\widehat{E});\mathbf{R}^{\infty}\right)\label{eq:Check_M(E)_Rep_2}
\end{equation}
by Corollary \ref{cor:Base_Sep_Meas} (b) (with $d=1$ and $A=\widehat{E}$)
and Lemma \ref{lem:Imb_SSP} (b) (with $E=\mathcal{M}^{+}(\widehat{E})$,
$S=\mathbf{R}^{\infty}$ and $\mathcal{D}=\mathfrak{mc}(\widehat{\mathcal{F}})^{*}$).
Furthermore,
\begin{equation}
\mathfrak{var}\left(\phi;\mathbf{R}^{\infty},S_{0},y_{0}\right)\in M\left(\mathbf{R}^{\infty};S_{0}\right)\label{eq:Check_M(E)_Rep_3}
\end{equation}
by the fact $S_{0}\in\mathscr{B}(\mathbf{R}^{\infty})$ and Fact  \ref{fact:var(f)}
(b) (with $(S,\mathscr{A})=(E,\mathscr{U})=(\mathbf{R}^{\infty},\mathscr{B}(\mathbf{R}^{\infty}))$,
$A=S_{0}$ and $f=\phi$). Now, (a) follows by (\ref{eq:Check_M(E)_Rep_2})
and (\ref{eq:Check_M(E)_Rep_3}).

(b) $Y\in M(\Omega,\mathscr{F};\mathcal{M}^{+}(\widehat{E}))$ by
(\ref{eq:Check_M(E)_Rep_1}) and (a). $\Psi|_{S_{0}}=\widehat{\varphi}^{-1}|_{S_{0}}$,
hence
\begin{equation}
\mathbb{P}\left(\varphi\circ X=\widehat{\varphi}\circ Y\in S_{0}\right)\geq\mathbb{P}\left(\varphi\circ X\in S_{0}\right)=1\label{eq:Check_M(E)_Rep_4}
\end{equation}
by (\ref{eq:FR-Base_M(E)}), which implies
\begin{equation}
\mathbb{P}\left(g^{*}\circ X=\widehat{g}^{*}\circ Y,\forall g\in\mathfrak{ag}(\mathcal{F})\right)=1\label{eq:Check_M(E)_Rep_5}
\end{equation}
by linearity of integral. Fixing $f\in\mathfrak{ca}(\mathcal{F})$
and $g\in\mathfrak{ag}(\mathcal{F})$, we find that
\begin{equation}
\begin{aligned}\left|f^{*}\circ X-\widehat{f}^{*}\circ Y\right|(\omega) & \leq\Vert f-g\Vert_{\infty}+\Vert\widehat{f}-\widehat{g}\Vert_{\infty}+\left|g^{*}\circ X-\widehat{g}^{*}\circ Y\right|(\omega)\\
 & \leq2\Vert f-g\Vert_{\infty}+\left|g^{*}\circ X-\widehat{g}^{*}\circ Y\right|(\omega),\;\forall\omega\in\Omega
\end{aligned}
\label{eq:Check_M(E)_Rep_6}
\end{equation}
by (\ref{eq:Check_M(E)_Rep_5}), Proposition \ref{prop:RepFun_Basic}
(d) (with $a=1$ and $b=-1$) and Fact \ref{fact:f+_Rep} (a) (with
$f=f-g$). Now, (b) follows by (\ref{eq:Check_M(E)_Rep_5}), (\ref{eq:Check_M(E)_Rep_6})
and (\ref{eq:ca(D)}) (with $\mathcal{D}=\widehat{\mathcal{F}}$).

(c) The replica $\nu^{\omega}$ of $X(\omega)$ exists for each fix
$\omega\in\Omega$ by Proposition \ref{prop:RepMeas_Basic} (a) (with
$d=1$, $\mu=X(\omega)$ and $\overline{\mu}=\nu^{\omega}$).
\begin{equation}
\begin{aligned} & \Omega\backslash\left\{ \omega\in\Omega:\varphi\circ X(\omega)=\widehat{\varphi}(\nu^{\omega})\in\widehat{\varphi}\left[\mathcal{M}^{+}(\widehat{E})\right]\right\} \\
 & =\left\{ \omega\in\Omega:X(\omega)(E\backslash E_{0})>0\right\} \in\mathscr{N}(\mathbb{P})
\end{aligned}
\label{eq:Check_M(E)_Rep_7}
\end{equation}
by the countability of $\mathfrak{mc}(\mathcal{F})$, Proposition
\ref{prop:RepMeas_Basic} (e) (with $d=1$, $\mu=X(\omega)$ and $\overline{\mu}=\nu^{\omega}$)
and (\ref{eq:Almostsure_Support_E0}). Since $S_{0}$ is contained
in the closure of $\widehat{\varphi}^{*}[\mathcal{M}^{+}(\widehat{E})]$,
it follows by (\ref{eq:Check_M(E)_Rep_4}), (\ref{eq:Check_M(E)_Rep_7})
and the completeness of $(\Omega,\mathcal{F},\mathbb{P})$ that
\begin{equation}
\mathbb{P}\left(\left\{ \omega\in\Omega:\widehat{\varphi}\circ Y(\omega)=\widehat{\varphi}(\nu^{\omega})\right\} \right)=1.\label{eq:Check_M(E)_Rep_8}
\end{equation}
Hence, (c) follows by (\ref{eq:Check_M(E)_Rep_2}) and (\ref{eq:Check_M(E)_Rep_8}).

(d) We fix $f\in\mathfrak{ca}(\mathcal{F})$ and $h\in M_{b}(E;\mathbf{R}^{k})$,
get $\{\overline{h\mathbf{1}_{A}},\overline{hf\mathbf{1}_{A}}\}\subset M_{b}(\widehat{E};\mathbf{R}^{k})$
from Proposition \ref{prop:RepFun_Basic} (b) (with $d=1$ and $f=h\mathbf{1}_{A}$
or $hf\mathbf{1}_{A}$), and find
\begin{equation}
\begin{aligned} & \left\{ \omega\in\Omega:(hf)^{*}\circ X(\omega)=(hf\mathbf{1}_{A})^{*}\circ X(\omega)=(\overline{hf\mathbf{1}_{A}})^{*}(\nu^{\omega})=(\overline{h\mathbf{1}_{A}}\widehat{f})^{*}(\nu^{\omega})\right\} \\
 & \supset\left\{ \omega\in\Omega:X(\omega)(E\backslash A)=0\right\} 
\end{aligned}
\label{eq:Check_M(E)_Rep_9}
\end{equation}
by Proposition \ref{prop:RepMeas_Basic} (d) (with $d=1$), the fact
$\overline{hf\mathbf{1}_{A}}|_{E_{0}}=\overline{h\mathbf{1}_{A}}\widehat{f}|_{E_{0}}$
and the definition of $\nu^{\omega}$. Thus, (\ref{eq:M(E)_RV_TF_bar_Rep})
follows by (\ref{eq:Check_M(E)_Rep_9}) and (c). Moreover, we have
$(\overline{h\mathbf{1}_{A}}\widehat{f})^{*}\circ Y\in M(\Omega,\mathscr{F};\mathbf{R}^{k})$
by the fact $\overline{h\mathbf{1}_{A}}\widehat{f}\in M_{b}(\widehat{E};\mathbf{R}^{k})$,
Lemma \ref{lem:Base} (c) and Proposition \ref{prop:M(E)_RV_Proc_TF}
(ii) (with $E=\widehat{E}$, $f=\overline{h\mathbf{1}_{A}}\widehat{f}$
and $\xi=Y$). Hence, $(hf)^{*}\circ X\in M(\Omega,\mathscr{F};\mathbf{R}^{k})$
by Lemma \ref{lem:var(X)} (a) (with $E=S=\mathbf{R}^{k}$, $\mathscr{U}=\mathscr{B}(\mathbf{R}^{k})$,
$X=(hf)^{*}\circ X$ and $Z=(\overline{h\mathbf{1}_{A}}\widehat{f})^{*}\circ Y$).\end{proof}

\begin{prop}
\label{prop:M(E)_Proc_Rep}Let $E$ be a topological space, $(E_{0},\mathcal{F};\widehat{E},\widehat{\mathcal{F}})$
be a base over $E$, and $(\Omega,\mathscr{F},\mathbb{P};X)$ be an
$\mathcal{M}^{+}(E)$-valued process satisfying
\begin{equation}
\left\{ \omega\in\Omega:X_{t}(\omega)(E\backslash E_{0})>0\right\} \in\mathscr{N}(\mathbb{P}),\;\forall t\in\mathbf{R}^{+}.\label{eq:Almostsure_Support_Proc}
\end{equation}
Then, there exists an $\mathcal{M}^{+}(\widehat{E})$-valued $\mathscr{F}_{t}^{X}$-adapted
process $(\Omega,\mathscr{F},\mathbb{P};Y)$ satisfying the following
properties:

\renewcommand{\labelenumi}{(\alph{enumi})}
\begin{enumerate}
\item $\varpi(\widehat{f}^{*})\circ Y$ is a modification of $\varpi(f^{*})\circ X$
for all $f\in\mathfrak{ca}(\mathcal{F})$.
\item If $X$ satisfies
\begin{equation}
\inf_{t\in\mathbf{R}^{+}}\mathbb{P}\left(X_{t}\in\mathcal{P}(E)\right)=1,\label{eq:R+-Base_P(E)}
\end{equation}
then $Y$ can be a $\mathcal{P}(\widehat{E})$-valued process.
\item For each $t\in\mathbf{R}^{+}$, there exists an $\Omega_{t}\in\mathscr{N}(\mathbb{P})$
such that $Y_{t}(\omega)$ equals the replica (measure) of $X_{t}(\omega)$
for all $\omega\in\Omega\backslash\Omega_{t}$.
\item If $A\in\mathscr{B}^{\mathbf{s}}(E_{0})$ satisfies
\begin{equation}
\left\{ \omega\in\Omega:X_{t}(\omega)(E\backslash A)>0\right\} \in\mathscr{N}(\mathbb{P}),\;\forall t\in\mathbf{R}^{+},\label{eq:Almostsure_Support_Proc-1}
\end{equation}
then $\varpi((\overline{h\mathbf{1}_{A}}\widehat{f})^{*})\circ Y$
is a modification of $\varpi((hf)^{*})\circ X$ for all $f\in\mathfrak{ca}(\mathcal{F})$,
$h\in M_{b}(E;\mathbf{R}^{k})$ and $k\in\mathbf{N}$.
\end{enumerate}
\end{prop}
\begin{rem}
\label{rem:M(E)_Proc_Mod}An implication of the statements ``$\varpi(\widehat{f}^{*})\circ Y$
is a modification of $\varpi(f^{*})\circ X$'' in (a) and ``$\varpi((\overline{h\mathbf{1}_{A}}\widehat{f})^{*})\circ Y$
is a modification of $\varpi((hf)^{*})\circ X$'' in (d) is that
$\varpi(f^{*})\circ X$, $\varpi(\widehat{f}^{*})\circ Y$, $\varpi((\overline{h\mathbf{1}_{A}}\widehat{f})^{*})\circ Y$
and $\varpi((hf)^{*})\circ X$ are processes. For (a), we know $f\in\mathfrak{ca}(\mathcal{F})\subset C_{b}(E;\mathbf{R})$
and $\widehat{f}\in C(\widehat{E};\mathbf{R})$, so $\varpi(f^{*})\circ X$
and $\varpi(\widehat{f}^{*})\circ Y$ are processes by Proposition
\ref{prop:M(E)_RV_Proc_TF} (a). For (d), $\varpi((\overline{h\mathbf{1}_{A}}\widehat{f})^{*})\circ Y$
and $\varpi((hf)^{*})\circ X$ are processes by Lemma \ref{lem:M(E)_RV_Rep}
(d) (with $X=X_{t}$ and $Y=Y_{t}$).
\end{rem}
\begin{proof}
We set $\varphi$, $\widehat{\varphi}$ $y_{0}$ and $\Psi$ as in
Lemma \ref{lem:M(E)_RV_Rep}. If (\ref{eq:R+-Base_P(E)}) holds, we
let $S_{0}=\widehat{\varphi}[\mathcal{P}(\widehat{E})]$ or else we
let $S_{0}=\widehat{\varphi}[\mathcal{M}^{+}(\widehat{E})]$. Recall
that $\widehat{\varphi}$ satisfies (\ref{eq:Check_M(E)_Rep_2}).
$\mathcal{M}^{+}(\widehat{E})$, $\mathcal{P}(\widehat{E})$ and $\mathbf{R}^{\infty}$
are Polish spaces by Corollary \ref{cor:Base_Sep_Meas} (c) (with
$d=1$) and Note \ref{note:Ehat_Valued_Proc_FDD}. Consequently, $S_{0}\in\mathscr{B}(\mathbf{R}^{\infty})$
in both cases by Proposition \ref{prop:SB_Map} (with $E=A=\mathcal{M}^{+}(\widehat{E})$
or $\mathcal{P}(\widehat{E})$, $f=\widehat{\varphi}$ and $S=\mathbf{R}^{\infty}$)
and Proposition \ref{prop:SB_Borel} (b) (with $E=\mathbf{R}^{\infty}$).
Hence, $Y\circeq\varpi(\Psi\circ\varphi)\circ X$ is the desired process
by (\ref{eq:Check_M(E)_Rep_1}), Lemma \ref{lem:M(E)_RV_Rep} (with
$X=\mathfrak{p}_{t}\circ X$ and $Y=\mathfrak{p}_{t}\circ Y$) and
Fact \ref{fact:Proc_Path_Mapping} (a) (with $E=\mathcal{M}^{+}(E)$,
$f=\Psi\circ\varphi$ and $\mathscr{G}_{t}=\mathscr{F}_{t}^{X}$).\end{proof}

\section{\label{sec:RepProc_FC}Finite-dimensional convergence about replica
process}

\subsection{\label{sub:Def_FC}Definitions}

Given a general space $E$, the finite-dimensional convergence of
$E$-valued processes is about the Borel extensions of their possibly
non-Borel finite-dimensional distributions.
\begin{defn}
\label{def:FC}Let $E$ be a topological space and $\{(\Omega^{i},\mathscr{F}^{i},\mathbb{P}^{i};X^{i})\}_{i\in\mathbf{I}}$%
\footnote{$\{(\Omega^{i},\mathscr{F}^{i},\mathbb{P}^{i})\}_{i\in\mathbf{I}}$
and $\{(\Omega^{n},\mathscr{F}^{n},\mathbb{P}^{n})\}_{n\in\mathbf{N}_{0}}$
were assumed in \S \ref{sec:Convention} to be complete probability
spaces. Completeness of measure space was specified in \ref{sub:Meas}.%
}, $\{(\Omega^{n},\mathscr{F}^{n},\mathbb{P}^{n};X^{n})\}_{n\in\mathbf{N}}$
and $(\Omega,\mathscr{F},\mathbb{P};X)$ be $E$-valued processes.
\begin{itemize}
\item $\{X^{n}\}_{n\in\mathbf{N}}$ \textbf{converges finite-dimensionally
to }$X$\textbf{ along $\mathbf{T}$} if: (1) $\mathbf{T}\subset\mathbf{R}^{+}$
is non-empty, (2) For each $\mathbf{T}_{0}\in\mathscr{P}_{0}(\mathbf{T})$,
there exist $N_{\mathbf{T}_{0}}\in\mathbf{N}$, $\{\mu_{n}\in\mathfrak{be}(\mathbb{P}^{n}\circ(X_{\mathbf{T}_{0}}^{n})^{-1})\}_{n>N_{\mathbf{T}_{0}}}$
and $\mu\in\mathfrak{be}(\mathbb{P}\circ X_{\mathbf{T}_{0}}^{-1})$
such that $\mu_{n}\Rightarrow\mu$ as $n\uparrow\infty$ in $\mathcal{P}(E^{\mathbf{T}_{0}})$.
\item $X$ is a \textbf{finite-dimensional limit point of $\{X^{i}\}_{i\in\mathbf{I}}$
along $\mathbf{T}$} if there exists a subsequence of $\{X^{i}\}_{i\in\mathbf{I}}$
converging finite-dimensionally to $X$ along $\mathbf{T}$.
\item Two finite-dimensional limit points of $\{X^{i}\}_{i\in\mathbf{I}}$
along $\mathbf{T}$ are equivalent if their finite-dimensional distributions
for any $\mathbf{T}_{0}\in\mathscr{P}_{0}(\mathbf{T})$ are identical.
\item $X$ is the \textbf{finite-dimensional limit of $\{X^{n}\}_{n\in\mathbf{N}}$
along $\mathbf{T}$} if $\{X^{n}\}_{n\in\mathbf{N}}$ converges finite-dimensionally
to $X$ along $\mathbf{T}$ and any finite-dimensional limit point
of \textbf{$\{X^{n}\}_{n\in\mathbf{N}}$} is equivalent to $X$.
\item $\{X^{i}\}_{i\in\mathbf{I}}$ is \textbf{finite-dimensionally convergent
along }$\mathbf{T}$ \textbf{under }$\mathcal{D}$ ($(\mathbf{T},\mathcal{D})$-FDC%
\footnote{Hereafter, ``$(\mathbf{T},\mathcal{D})$-FDC'' and ``$(\mathbf{T},\mathcal{D})$-AS''
also stand for ``$(\mathbf{T},\mathcal{D})$-finite-dimensional convergence''
and ``$(\mathbf{T},\mathcal{D})$-asymptotical stationarity''.%
} for short) if: (1) $\mathbf{I}$ is infinite, (2) $\mathbf{T}\subset\mathbf{R}^{+}$
and $\mathcal{D}\subset M_{b}(E;\mathbf{R})$ are non-empty, and (3)
$\{\mathbb{E}^{i}[f\circ X_{\mathbf{T}_{0}}^{i}]\}_{i\in\mathbf{I}}$%
\footnote{$\mathbb{E}^{i}$ denotes the expectation operator of $(\Omega^{i},\mathscr{F}^{i},\mathbb{P}^{i})$.%
} has a unique%
\footnote{In the definitions of $(\mathbf{T},\mathcal{D})$-FDC and $(\mathbf{T},\mathcal{D})$-AS,
$\{\mathbb{E}^{i}[f\circ X_{\mathbf{T}_{0}}^{i}]\}_{i\in\mathbf{I}}$
and $\{\mathbb{E}^{i}[f\circ X_{\mathbf{T}_{0}}^{i}-f\circ X_{\mathbf{T}_{0}+c}^{i}]\}_{i\in\mathbf{I}}$
both lie in $[-2\Vert f\Vert_{\infty},2\Vert f\Vert_{\infty}]$. Each
of them has at least one limit point in $\mathbf{R}$ by the Bolzano-Weierstrass
Theorem, so it is enough to assume ``at most one limit point''.%
} limit point in $\mathbf{R}$ for all $f\in\mathfrak{mc}[\Pi^{\mathbf{T}_{0}}(\mathcal{D})]$
and $\mathbf{T}_{0}\in\mathscr{P}_{0}(\mathbf{T})$.
\item $\{X^{i}\}_{i\in\mathbf{I}}$ is \textbf{asymptotically stationary
along }$\mathbf{T}$ \textbf{under} $\mathcal{D}$ ($(\mathbf{T},\mathcal{D})$-AS
for short) if: (1) $\mathbf{I}$ is infinite, (2) $\mathbf{T}\subset\mathbf{R}^{+}$
and $\mathcal{D}\subset M_{b}(E;\mathbf{R})$ are non-empty, and (3)
The unique limit point of $\{\mathbb{E}^{i}[f\circ X_{\mathbf{T}_{0}}^{i}-f\circ X_{\mathbf{T}_{0}+c}^{i}]\}_{i\in\mathbf{I}}$%
\footnote{The notation ``$\mathbf{T}_{0}+c$'' was defined in (\ref{eq:Timeset_Translation}).%
} in $\mathbf{R}$ is $0$ for all $c\in(0,\infty)$, $f\in\mathfrak{mc}[\Pi^{\mathbf{T}_{0}}(\mathcal{D})]$
and $\mathbf{T}_{0}\in\mathscr{P}_{0}(\mathbf{T})$.
\end{itemize}
\end{defn}
\begin{note}
\label{note:Proc_Int_Test_Integrability}Given $E$-valued process
$X$, the expectation $\mathbb{E}[f\circ X_{\mathbf{T}_{0}}]$ is
well-defined for any $f\in\mathfrak{ca}[\Pi^{d}(M_{b}(E;\mathbf{R}))]$
by Fact \ref{fact:Proc_Basic_1} (d) and Note \ref{note:Pi^d(D)_Mb_Cb}
(with $\mu=\mathbb{P}\circ X_{\mathbf{T}_{0}}^{-1}$).
\end{note}

\begin{note}
\label{note:FLP_Equi}Let $X$ and $Y$ be $E$-valued processes and
$\mathbf{T}\subset\mathbf{R}^{+}$. We relate $X\sim Y$ if $X_{\mathbf{T}_{0}}$
and $Y_{\mathbf{T}_{0}}$ have the same distribution for all $\mathbf{T}_{0}\in\mathscr{P}_{0}(\mathbf{T})$.
This ``$\sim$'', which Definition \ref{def:FC} uses to define
equivalence of finite-dimensional limit points, is an equivalence
relation among $E$-valued stochastic processes.
\end{note}
We make the following notations for simplicity.
\begin{notation}
\label{notation:FLP}Let $X$, $\{X^{i}\}_{i\in\mathbf{I}}$ and $\{X^{n}\}_{n\in\mathbf{N}}$
be $E$-valued processes.

$\{X^{n}\}_{n\in\mathbf{N}}$ converging finite-dimensionally to $X$
along $\mathbf{T}$ is denoted by
\begin{equation}
X^{n}\xrightarrow{\quad\mathrm{D}(\mathbf{T})\quad}X\mbox{ as }n\uparrow\infty.\label{eq:FC_along_T}
\end{equation}

\begin{itemize}
\item By $X=\mathfrak{fl}_{\mathbf{T}}(\{X^{n}\}_{n\in\mathbf{N}})$%
\footnote{``$\mathfrak{fl}$'' and ``$\mathfrak{flp}$'' are ``fl'' and
``flp'' in fraktur font which stand for ``finite-dimensional limit''
and ``finite-dimensional limit point'' respectively. Members of
$\mathfrak{fl}_{\mathbf{T}}(\cdot)$ or $\mathfrak{flp}_{\mathbf{T}}(\cdot)$
are processes with time horizon $\mathbf{R}^{+}$ no matter $\mathbf{T}=\mathbf{R}^{+}$
or not. %
} we mean $X$ is the finite-dimensional limit of $\{X^{n}\}_{n\in\mathbf{N}}$
along $\mathbf{T}$.
\item By $\mathfrak{flp}_{\mathbf{T}}(\{X^{i}\}_{i\in\mathbf{I}})$ we denote
the family of all equivalence classes (see Note \ref{note:FLP_Equi})
of finite-dimensional limit points of $\{X^{i}\}_{i\in\mathbf{I}}$
along $\mathbf{T}$.
\item By $X=\mathfrak{flp}_{\mathbf{T}}(\{X^{i}\}_{i\in\mathbf{I}})$ we
mean $X$ is the unique member of $\mathfrak{flp}_{\mathbf{T}}(\{X^{i}\}_{i\in\mathbf{I}})$.
\end{itemize}
\end{notation}
\begin{rem}
\label{rem:FL}In general, $X=\mathfrak{fl}_{\mathbf{T}}(\{X^{n}\}_{n\in\mathbf{N}})$
is stronger than (\ref{eq:FC_along_T}) because: (1) Each of the finite-dimensional
distributions of $\{X^{n}\}_{n\in\mathbf{N}}$ may have multiple Borel
extensions. (2) $\mathcal{P}(E^{\mathbf{T}_{0}})$ is not necessarily
a Hausdorff space and a weakly convergent sequence may have multiple
limits.
\end{rem}
The following is a straightforward property of finite-dimensional
convergence.
\begin{fact}
\label{fact:FC_FDC}Let $E$ be a topological space, $\mathbf{T}\subset\mathbf{R}^{+}$
and $\{(\Omega^{n},\mathscr{F}^{n},\mathbb{P}^{n};X^{n})\}_{n\in\mathbf{N}}$
and $(\Omega,\mathscr{F},\mathbb{P};X)$ be $E$-valued processes.
If (\ref{eq:FC_along_T}) holds, then%
\footnote{$\mathbb{E}^{n}$ denotes the expectation operator of $(\Omega^{n},\mathscr{F}^{n},\mathbb{P}^{n})$.%
}
\begin{equation}
\lim_{n\rightarrow\infty}\mathbb{E}^{n}\left[f\circ X_{\mathbf{T}_{0}}^{n}\right]=\mathbb{E}\left[f\circ X_{\mathbf{T}_{0}}\right]\label{eq:Exp_Test}
\end{equation}
for all $f\in\mathfrak{mc}[\Pi^{\mathbf{T}_{0}}(C_{b}(E;\mathbf{R}))]$
and $\mathbf{T}_{0}\in\mathscr{P}_{0}(\mathbf{T})$. As a consequence,
$\{X^{n}\}_{n\in\mathbf{N}}$ is $(\mathbf{T},C_{b}(E;\mathbf{R}))$-FDC.
\end{fact}
\begin{proof}
Fixing $f\in\mathfrak{mc}[\Pi^{\mathbf{T}_{0}}(C_{b}(E;\mathbf{R}))]$
and $\mathbf{T}_{0}\in\mathscr{P}_{0}(\mathbf{T})$, it follows by
Fact \ref{fact:BExt_Same_Int} (with $d=\aleph(\mathbf{T}_{0})$ and
$X=X_{\mathbf{T}_{0}}^{n}$ or $X_{\mathbf{T}_{0}}$), (\ref{eq:ca(Pi^d(D))_Cb})
(with $d=\aleph(\mathbf{T}_{0})$ and $\mathcal{D}=C_{b}(E;\mathbf{R})$)
and (\ref{eq:FC_along_T}) that
\begin{equation}
\lim_{n\rightarrow\infty}\mathbb{E}^{n}\left[f\circ X_{\mathbf{T}_{0}}^{n}\right]=\lim_{n\rightarrow\infty}f^{*}(\mu_{n})=f^{*}(\mu)=\mathbb{E}\left[f\circ X_{\mathbf{T}_{0}}\right]\label{eq:Check_Proc_FC_Int_Test}
\end{equation}
for some $\{\mu_{n}\in\mathfrak{be}(\mathbb{P}\circ(X_{\mathbf{T}_{0}}^{n})^{-1}\}_{n>N_{\mathbf{T}_{0}}}$
with $N_{\mathbf{T}_{0}}\in\mathbf{N}$ and $\mu\in\mathfrak{be}(\mathbb{P}\circ X_{\mathbf{T}_{0}}^{-1})$.\end{proof}

\subsection{\label{sub:TransFC}Transformation of finite-dimensional convergence}

The following theorem is our main tool for transforming a finite-dimensional
limit point of replica processes back into that of original processes.
\begin{thm}
\label{thm:TransFC_1}Let $E$ be a topological space, $\{(\Omega^{n},\mathscr{F}^{n},\mathbb{P}^{n};X^{n})\}_{n\in\mathbf{N}}$
be $E$-valued processes, $(E_{0},\mathcal{F};\widehat{E},\widehat{\mathcal{F}})$
be a base over $E$, $x_{0}\in E_{0}$, $\mathbf{T}\subset\mathbf{R}^{+}$%
\footnote{Note \ref{note:Var} mentioned that $\mathfrak{var}(Y_{t};\Omega,Y_{t}^{-1}(\{x_{0}\}),x_{0})$
is the constant mapping that sends every $\omega\in\Omega$ to $x_{0}$.
We do not use $x_{0}$ to denote this mapping for clarity.%
}, and
\begin{equation}
X_{t}\circeq\begin{cases}
\mathfrak{var}\left(Y_{t};\Omega,Y_{t}^{-1}(E_{0}),x_{0}\right), & \mbox{if }t\in\mathbf{T},\\
\mathfrak{var}\left(Y_{t};\Omega,Y_{t}^{-1}(\{x_{0}\}),x_{0}\right), & \mbox{if }t\in\mathbf{R}^{+}\backslash\mathbf{T},
\end{cases}\label{eq:Trans_FC_Lim_Proc}
\end{equation}
where $(\Omega,\mathscr{F},\mathbb{P};Y)$ is an $\widehat{E}$-valued
process. Suppose that:

\renewcommand{\labelenumi}{(\roman{enumi})}
\begin{enumerate}
\item $\{X_{t}^{n}\}_{n\in\mathbf{N}}$ is sequentially tight%
\footnote{Sequential tightness of random variabels was specified in Note \ref{note:Seq_Tight_RV}.%
} in $E_{0}$ for all $t\in\mathbf{T}$.
\item $\{X^{n}\}_{n\in\mathbf{N}}$ and $Y$ satisfy
\begin{equation}
\lim_{n\rightarrow\infty}\mathbb{E}^{n}\left[f\circ X_{\mathbf{T}_{0}}^{n}\right]=\mathbb{E}\left[\widehat{f}\circ Y_{\mathbf{T}_{0}}\right]\label{eq:Exp_Test_Y}
\end{equation}
for all $f\in\mathfrak{mc}[\Pi^{\mathbf{T}_{0}}(\mathcal{F}\backslash\{1\})]$
and $\mathbf{T}_{0}\in\mathscr{P}_{0}(\mathbf{T})$.
\end{enumerate}
Then:

\renewcommand{\labelenumi}{(\alph{enumi})}
\begin{enumerate}
\item $X\circeq\{X_{t}\}_{t\geq0}$ is an $E$-valued process with paths
in $E_{0}^{\mathbf{R}^{+}}$%
\footnote{``with paths in $E_{0}^{\mathbf{R}^{+}}$'' means all paths of the
process lying in $E_{0}^{\mathbf{R}^{+}}$. Of course, an $E$-valued
process with paths in $E_{0}^{\mathbf{R}^{+}}$ is equivalent to an
$(E_{0},\mathscr{O}_{E}(E_{0}))$-valued process.%
} and for each $\mathbf{T}_{0}\in\mathscr{P}_{0}(\mathbf{T})$, there
exist%
\footnote{The notation ``$\mu_{\mathbf{T}_{0},n}^{\prime}=\mathfrak{be}(\mathbb{P}\circ(X_{\mathbf{T}_{0}}^{n})^{-1})$''
as defined in \S \ref{sub:Topo} means $\mu_{\mathbf{T}_{0},n}^{\prime}$
is the unique Borel extension of $\mathbb{P}\circ(X_{\mathbf{T}_{0}}^{n})^{-1}$.%
}
\begin{equation}
\mu_{\mathbf{T}_{0},n}^{\prime}=\mathfrak{be}\left(\mathbb{P}\circ\left(X_{\mathbf{T}_{0}}^{n}\right)^{-1}\right)\in\mathcal{P}(E^{\mathbf{T}_{0}}),\;\forall n>N_{\mathbf{T}_{0}}\label{eq:Proc_FDD_BExt_Uni}
\end{equation}
for some $N_{\mathbf{T}_{0}}\in\mathbf{N}$ and
\begin{equation}
\mu_{\mathbf{T}_{0}}=\mathfrak{be}\left(\mathbb{P}\circ X_{\mathbf{T}_{0}}^{-1}\right)\in\mathcal{P}(E^{\mathbf{T}_{0}})\label{eq:Lim_Proc_FDD_Borel}
\end{equation}
such that%
\footnote{The notation ``$\mathrm{w}$-$\lim$'', introduced in \S \ref{sec:Borel_Measure},
means weak limit of a sequence of finite Borel measures.%
}
\begin{equation}
\mathrm{w}\mbox{-}\lim_{n\rightarrow\infty}\left.\mu_{\mathbf{T}_{0},n}^{\prime}\right|_{E_{0}^{\mathbf{T}_{0}}}=\left.\mu_{\mathbf{T}_{0}}\right|_{E_{0}^{\mathbf{T}_{0}}}\mbox{ in }\mathcal{P}\left(E_{0}^{\mathbf{T}_{0}},\mathscr{O}_{E}(E_{0})^{\mathbf{T}_{0}}\right).\label{eq:FC_along_T_E0}
\end{equation}
Moreover, $X_{\mathbf{T}_{0}}$ is $\mathbf{m}$-tight in $E_{0}^{\mathbf{T}_{0}}$
for all $\mathbf{T}_{0}\in\mathscr{P}_{0}(\mathbf{T})$,
\begin{equation}
\inf_{t\in\mathbf{T}}\mathbb{P}\left(X_{t}=Y_{t}\in E_{0}\right)=1\label{eq:Lim_Proc_(T,E0)-Mod}
\end{equation}
and (\ref{eq:FC_along_T}) holds.
\item If $\mathbf{T}=\mathbf{R}^{+}$ and $\{X^{n}\}_{n\in\mathbf{N}}$
is $(\mathbf{R}^{+},\mathcal{F}\backslash\{1\})$-AS, then $X$ and
$Y$ are both stationary processes%
\footnote{The notion of stationary process was specified in \S \ref{sec:Proc}%
}.
\item If $C_{b}(E;\mathbf{R})$ separates points on $E$%
\footnote{Note \ref{note:C(E;R)_Cb(E;R)_SP} showed that $C_{b}(E;\mathbf{R})$
separating points is equivalent to $C(E;\mathbf{R})$ separating points
on $E$.%
}, then $X=\mathfrak{fl}_{\mathbf{T}}(\{X_{n}\}_{n\in\mathbf{N}})$.
\end{enumerate}
\end{thm}
\begin{rem}
\label{rem:TransFC_Seq_Tight}The condition (i) above will ensure
the following two facts for $\mathbf{T}_{0}\in\mathscr{P}_{0}(\mathbf{T})$
with finite exception of $n\in\mathbf{N}$: (1) $\mathbb{P}^{n}\circ(X_{\mathbf{T}_{0}}^{n})^{-1}$
admits a unique Borel extension, and (2) $\mathbb{P}^{n}\circ(\widehat{X}_{\mathbf{T}_{0}}^{n})^{-1}$
is the replica measure of $\mathbb{P}^{n}\circ(X_{\mathbf{T}_{0}}^{n})^{-1}$.
Hence, transforming finite-dimensional convergence from the replicas
$\{\widehat{X}^{n}\}_{n\in\mathbf{N}}$ to the original processes
$\{X^{n}\}_{n\in\mathbf{N}}$ comes down to transforming weak convergence
from replica measures to Borel extensions of original measures.
\end{rem}
\begin{proof}
[Proof of Theorem \ref{thm:TransFC_1}](a) We divide the proof of
(a) into four steps.

\textit{Step 1: Constructing $\mu_{\mathbf{T}_{0}}$ for each} $\mathbf{T}_{0}\in\mathscr{P}_{0}(\mathbf{T})$.
We know from (i) and Fact \ref{fact:Proc_FDD_BExt} (with $\mathbf{I}=\mathbf{N}$)
that $\{\mu_{\mathbf{T}_{0},n}\circeq\mathbb{P}^{n}\circ(X_{\mathbf{T}_{0}}^{n})^{-1}\}_{n\in\mathbf{N}}$
is sequentially tight in $E_{0}^{\mathbf{T}_{0}}$ and there exists
an $N_{\mathbf{T}_{0}}\in\mathbf{N}$ such that
\begin{equation}
\inf_{n>N_{\mathbf{T}_{0}},t\in\mathbf{T}_{0}}\mathbb{P}^{n}\left(X_{t}^{n}\in E_{0}\right)=1\label{eq:TransFC_Common_T-Base}
\end{equation}
and the $\{\mu_{\mathbf{T}_{0},n}^{\prime}=\mathfrak{be}(\mu_{\mathbf{T}_{0},n}^{\prime})\}_{n>N_{\mathbf{T}_{0}}}$
in (\ref{eq:Proc_FDD_BExt_Uni}) exist. We then have
\begin{equation}
\inf_{n>N_{\mathbf{T}_{0}},t\in\mathbf{T}_{0}}\mathbb{P}^{n}\left(\bigotimes\mathcal{F}\circ X_{t}^{n}\in\bigotimes\widehat{\mathcal{F}}(\widehat{E})\right)=1\label{eq:TransFC_Common_FR-Base}
\end{equation}
by (\ref{eq:TransFC_Common_T-Base}) and Fact \ref{fact:T-Base_FR-Base}
(with $X=X^{n}$). From (\ref{eq:TransFC_Common_FR-Base}), the condition
(ii) above and Lemma \ref{lem:RepProc_Int_Test} (c, e) it follows
that
\begin{equation}
\widehat{X}^{n}\xrightarrow{\quad\mathrm{D}(\mathbf{T})\quad}Y\mbox{ as }n\uparrow\infty.\label{eq:RepProc_FC_along_T_Y}
\end{equation}
From (\ref{eq:RepProc_FC_along_T_Y}) and Proposition \ref{prop:RepProc_T-Base}
(b) (with $X=X^{n}$ and $\mathbf{T}=\mathbf{T}_{0}$) it follows
that%
\footnote{$\overline{\mu}_{\mathbf{T}_{0},n}$, as specified in Notation \ref{notation:ReplicaMeas},
denotes the replica measure of $\mu_{\mathbf{T}_{0},n}$.%
}
\begin{equation}
\overline{\mu}_{\mathbf{T}_{0},n}=\mathbb{P}^{n}\circ\left(\widehat{X}_{\mathbf{T}_{0}}^{n}\right)^{-1}\Longrightarrow\nu_{\mathbf{T}_{0}}\circeq\mathbb{P}\circ Y_{\mathbf{T}_{0}}^{-1}\mbox{ as }n\uparrow\infty\mbox{ in }\mathcal{P}(\widehat{E}^{\mathbf{T}_{0}}).\label{eq:RepProc_FDD_WC_Lim}
\end{equation}
Now, by Corollary \ref{cor:WC_Pull_Back_FinDim} (with $d\circeq\aleph(\mathbf{T}_{0})$,
$\mu_{n}=\mu_{\mathbf{T}_{0},n+N_{\mathbf{T}_{0}}}$, $\mu_{n}^{\prime}=\mu_{\mathbf{T}_{0},n+N_{\mathbf{T}_{0}}}^{\prime}$
and $\nu=\nu_{\mathbf{T}_{0}}$) and Fact \ref{fact:P(E)_Closed_M(E)}
(with $E=E^{\mathbf{T}_{0}}$ or $(E_{0}^{\mathbf{T}_{0}},\mathscr{O}_{E}(E_{0})^{\mathbf{T}_{0}})$),
there exists a $\mu_{\mathbf{T}_{0}}\in\mathcal{P}(E^{\mathbf{T}_{0}})$
such that (\ref{eq:FC_along_T_E0}) holds, $\mu_{\mathbf{T}_{0}}$
is tight in $E_{0}^{\mathbf{T}_{0}}$, $\nu_{\mathbf{T}_{0}}=\overline{\mu}_{\mathbf{T}_{0}}$
and
\begin{equation}
\mu_{\mathbf{T}_{0},n}^{\prime}\Longrightarrow\mu_{\mathbf{T}_{0}}\mbox{ as }n\uparrow\infty\mbox{ in }\mathcal{P}\left(E^{\mathbf{T}_{0}}\right).\label{eq:Check_FC}
\end{equation}

\textit{Step 2: Verify $X=\{X_{t}\}_{t\geq0}$ defined by (\ref{eq:Trans_FC_Lim_Proc})
is a process and satisfies }(\ref{eq:Lim_Proc_(T,E0)-Mod}). For each
$t\in\mathbf{T}$, we let $\mu_{\{t\}}\in\mathcal{P}(E)$ be the measure
constructed in Step 1 with $\mathbf{T}_{0}=\{t\}$. By our argument
above, each $\mu_{\{t\}}$ is tight in $(E_{0},\mathscr{O}_{E}(E_{0}))$
and so is supported on some $S_{t}\in\mathscr{K}_{\sigma}(E_{0},\mathscr{O}_{E}(E_{0}))$.
$S_{t}\in\mathscr{B}(\widehat{E})$ and $\mathscr{B}_{E}(S_{t})=\mathscr{B}_{\widehat{E}}(S_{t})$
by Corollary \ref{cor:Base_Compact} (b) and Lemma \ref{lem:SB_Base}
(a). Let $\nu_{\{t\}}=\mathbb{P}\circ Y_{t}^{-1}=\overline{\mu}_{\{t\}}$
be defined as in (\ref{eq:RepProc_FDD_WC_Lim}) with $\mathbf{T}_{0}=\{t\}$.
It follows by Proposition \ref{prop:RepMeas_Basic} (a) that
\begin{equation}
\mathbb{P}\left(Y_{t}\in S_{t}\right)=\nu_{\{t\}}(S_{t})=\mu_{\{t\}}(S_{t})=1,\;\forall t\in\mathbf{T}.\label{eq:Lim_FDD_Tight_E0}
\end{equation}
Hence,
\begin{equation}
X_{t}\in M\left(\Omega,\mathscr{F};E_{0},\mathscr{B}_{E}(E_{0})\right),\;\forall t\in\mathbf{T}\label{eq:Lim_Proc_at_T}
\end{equation}
satisfy (\ref{eq:Lim_Proc_(T,E0)-Mod}) by Lemma \ref{lem:var(X)}
(b, c) (with $(E,\mathscr{U})=(\widehat{E},\mathscr{B}(\widehat{E}))$,
$S_{0}=S_{t}$, $(S,\mathscr{U}^{\prime})=(E_{0},\mathscr{B}_{E}(E_{0}))$,
$X=Y_{t}$ and $Y=X_{t}$). Furthermore, we have
\begin{equation}
\{x_{0}\}\in\mathscr{B}(\widehat{E})\cap\mathscr{B}\left(E_{0},\mathscr{O}_{E}(E_{0})\right)\cap\mathscr{B}(E),\label{eq:{x0}_Base_Borel}
\end{equation}
by Lemma \ref{lem:Base} (c, e), the fact $E_{0}\in\mathscr{B}(E)$
and Proposition \ref{prop:Separability} (a), which implies
\begin{equation}
X_{t}\in M\left(\Omega,\mathscr{F};E_{0},\mathscr{O}_{E}(E_{0})\right),\;\forall t\in\mathbf{R}^{+}\backslash\mathbf{T}.\label{eq:Lim_Proc_Outside_T}
\end{equation}
by (\ref{eq:Trans_FC_Lim_Proc}). Now, $X$ is an $(E_{0},\mathscr{O}_{E}(E_{0}))$-valued
process by Fact \ref{fact:Proc_Basic_1} (b).

\textit{Step 3: Verify the $\mathbf{m}$-tightness of $X_{\mathbf{T}_{0}}$
in} $E_{0}^{\mathbf{T}_{0}}$\textit{ and (\ref{eq:Lim_Proc_FDD_Borel})
for each $\mathbf{T}_{0}\in\mathscr{P}_{0}(\mathbf{T})$}. Letting
$\{S_{t}\}_{t\in\mathbf{T}}$ be as in Step 2, we have that
\begin{equation}
S_{\mathbf{T}_{0}}\circeq\prod_{t\in\mathbf{T}_{0}}S_{t}\in\mathscr{K}_{\sigma}\left(E_{0}^{\mathbf{T}_{0}},\mathscr{O}_{E}(E_{0})^{\mathbf{T}_{0}}\right)\label{eq:Lim_Proc_FDD_Support_1}
\end{equation}
by Corollary \ref{cor:Sigma_Compact_Prod} (b) (with $\mathbf{I}=\mathbf{T}_{0}$
and $S_{i}=(E_{0},\mathscr{O}_{E}(E_{0}))$). It follows that
\begin{equation}
S_{\mathbf{T}_{0}}\in\mathscr{B}(\widehat{E}^{\mathbf{T}_{0}})\cap\mathscr{B}(E)^{\otimes\mathbf{T}_{0}}\cap\mathscr{K}_{\sigma}^{\mathbf{m}}(E^{\mathbf{T}_{0}})\label{eq:Lim_Proc_FDD_Support_2}
\end{equation}
and
\begin{equation}
\mathscr{B}_{E^{\mathbf{T}_{0}}}(S_{\mathbf{T}_{0}})=\mathscr{B}_{\widehat{E}^{\mathbf{T}_{0}}}(S_{\mathbf{T}_{0}})\label{eq:Lim_Proc_FDD_Support_3}
\end{equation}
by Corollary \ref{cor:Base_Compact} (b) and Lemma \ref{lem:SB_Base}
(a). Now,
\begin{equation}
\nu_{\mathbf{T}_{0}}\left(S_{\mathbf{T}_{0}}\right)=\mathbb{P}\left(Y_{\mathbf{T}_{0}}\in S_{\mathbf{T}_{0}}\right)=1\label{eq:Check_Lim_Proc_FDD_Tight}
\end{equation}
by (\ref{eq:Lim_Proc_FDD_Support_2}) and (\ref{eq:Lim_FDD_Tight_E0}).
Moreover,
\begin{equation}
\nu_{\mathbf{T}_{0}}\left(A\cap S_{\mathbf{T}_{0}}\right)=\mu_{\mathbf{T}_{0}}\left(A\cap S_{\mathbf{T}_{0}}\right)=\mu_{\mathbf{T}_{0}}(A),\;\forall A\in\mathscr{B}(E^{\mathbf{T}_{0}})\label{eq:Lim_Proc_FDD_RepMeas}
\end{equation}
by (\ref{eq:Lim_Proc_FDD_Support_2}), (\ref{eq:Lim_Proc_FDD_Support_3}),
the fact $\nu_{\mathbf{T}_{0}}=\overline{\mu}_{\mathbf{T}_{0}}$,
(\ref{eq:Check_Lim_Proc_FDD_Tight}) and Proposition \ref{prop:RepMeas_Basic}
(a) (with $\mu=\mu_{\mathbf{T}_{0}}$). We established (\ref{eq:Lim_Proc_(T,E0)-Mod})
in Step 2. It follows that
\begin{equation}
\begin{aligned}\mathbb{P}\left(X_{\mathbf{T}_{0}}\in A\right) & =\mathbb{P}\left(Y_{\mathbf{T}_{0}}\in A\right)\\
 & =\nu_{\mathbf{T}_{0}}\left(A\cap S_{\mathbf{T}_{0}}\right)=\mu_{\mathbf{T}_{0}}(A),\;\forall A\in\mathscr{B}(E)^{\bigotimes\mathbf{T}_{0}}
\end{aligned}
\label{eq:Check_Lim_Proc_FDD_BExt}
\end{equation}
by (\ref{eq:Lim_Proc_(T,E0)-Mod}), (\ref{eq:Check_Lim_Proc_FDD_Tight})
and (\ref{eq:Lim_Proc_FDD_RepMeas}). Thus, $X_{\mathbf{T}_{0}}$
is $\mathbf{m}$-tight in $E_{0}^{\mathbf{T}_{0}}$ by (\ref{eq:Lim_Proc_FDD_Support_1}),
(\ref{eq:Lim_Proc_FDD_Support_2}), (\ref{eq:Check_Lim_Proc_FDD_Tight})
and (\ref{eq:Check_Lim_Proc_FDD_BExt}). (\ref{eq:Lim_Proc_FDD_Borel})
follows by (\ref{eq:Check_Lim_Proc_FDD_BExt}) and Proposition \ref{prop:m-Tight_BExt}
(with $\mathbf{I}=\mathbf{T}_{0}$, $S_{i}=E$, $A=E_{0}^{\mathbf{T}_{0}}$
and $\Gamma=\{\mathbb{P}\circ X_{\mathbf{T}_{0}}^{-1}\}$).

\textit{Step 4: Verify (\ref{eq:FC_along_T})}. The desired convergence
follows from (\ref{eq:Check_FC}) and (\ref{eq:Lim_Proc_FDD_Borel}).

(b) We have by (a) (with $\mathbf{T}=\mathbf{R}^{+}$) that
\begin{equation}
\inf_{t\in\mathbf{R}^{+}}\mathbb{P}\left(X_{t}=Y_{t}\in E_{0}\right)=1\label{eq:Lim_Proc_(R+,E0)_Mod}
\end{equation}
and
\begin{equation}
X^{n}\xrightarrow{\quad\mathrm{D}(\mathbf{R}^{+})\quad}X\mbox{ as }n\uparrow\infty.\label{eq:FC_along_R+}
\end{equation}
It then follows that

\begin{equation}
\begin{aligned}\mathbb{E}\left[\widehat{f}\circ Y_{\mathbf{T}_{0}}-\widehat{f}\circ Y_{\mathbf{T}_{0}+c}\right] & =\mathbb{E}\left[f\circ X_{\mathbf{T}_{0}}-f\circ X_{\mathbf{T}_{0}+c}\right]\\
 & =\lim_{n\rightarrow\infty}\mathbb{E}^{n}\left[f\circ X_{\mathbf{T}_{0}}^{n}-f\circ X_{\mathbf{T}_{0}+c}^{n}\right]=0
\end{aligned}
\label{eq:Check_Lim_Proc_R+_Mod_Sta}
\end{equation}
for all $c>0$, $\mathbf{T}_{0}\in\mathscr{P}_{0}(\mathbf{R}^{+})$
and $f\in\mathfrak{mc}[\Pi^{\mathbf{T}_{0}}(\mathcal{F}\backslash\{1\})]$
by (\ref{eq:FC_along_R+}), Fact \ref{fact:FDC_AS} (b) (with $\mathbf{T}=\mathbf{R}^{+}$),
Fact \ref{fact:FC_FDC} (with $\mathbf{T}=\mathbf{R}^{+}$), (\ref{eq:Lim_Proc_(R+,E0)_Mod})
and Lemma \ref{lem:Proc_Rep} (a) (with $\mathbf{T}=\mathbf{R}^{+}$).
Hence, the stationarity of $Y$ follows by Corollary \ref{cor:Base_Sep_Meas}
(a)%
\footnote{$\mathbb{P}\circ\widehat{X}_{\mathbf{T}_{0}}^{-1}$ and each $\mathbb{P}^{n}\circ(\widehat{X}_{\mathbf{T}_{0}}^{n})^{-1}$
are Borel probability measures as mentioned in Note \ref{note:Ehat_Valued_Proc_FDD}.%
} (with $d=\aleph(\mathbf{T}_{0})$ and $A=\widehat{E}^{d}$). $\mathbb{P}\circ X_{0}^{-1}$
is $\mathbf{m}$-tight in $E_{0}$ by (a), so there exists an $A\in\mathscr{K}_{\sigma}(E_{0},\mathscr{O}_{E}(E_{0}))$
such that
\begin{equation}
\inf_{t\in\mathbf{R}^{+}}\mathbb{P}\left(X_{t}=Y_{t}\in A\right)=1,\label{eq:Lim_Proc_(R+,A)-Mod}
\end{equation}
by (\ref{eq:Lim_Proc_(R+,E0)_Mod}) and the stationarity of $Y$.
Now, $X$ is stationary by Lemma \ref{lem:Proc_Rep} (e).

(c) We fix $\mathbf{T}_{0}\in\mathscr{P}_{0}(\mathbf{T})$ and let
each $\mu_{\mathbf{T}_{0},n}$, $\mu_{\mathbf{T}_{0},n}^{\prime}$
and $\mu_{\mathbf{T}_{0}}$ be as in (a). It follows by (\ref{eq:Check_FC}),
(\ref{eq:ca(Pi^d(D))_Cb}) (with $\mathcal{D}=C_{b}(E;\mathbf{R})$)
and Fact \ref{fact:BExt_Same_Int} (with $d=\aleph(\mathbf{T}_{0})$,
$\mu=\mu_{\mathbf{T}_{0},n}$ and $\nu_{1}=\mu_{\mathbf{T}_{0},n}^{\prime}$)
that
\begin{equation}
\lim_{n\rightarrow\infty}\int_{E^{\mathbf{T}_{0}}}f(x)\mu_{\mathbf{T}_{0},n}(dx)=\lim_{n\rightarrow\infty}f^{*}\left(\mu_{\mathbf{T}_{0},n}^{\prime}\right)=f^{*}(\mu_{\mathbf{T}_{0}})\label{eq:Check_Int_Test_All_Cb(E;R)}
\end{equation}
for all $f\in\mathfrak{mc}[\Pi^{\mathbf{T}_{0}}(C_{b}(E;\mathbf{R}))]$.
Hence, (\ref{eq:Check_FC}) implies
\begin{equation}
\mathrm{w}\mbox{-}\lim_{n\rightarrow\infty}\mu_{\mathbf{T}_{0},n}^{\prime}=\mu_{\mathbf{T}_{0}}\label{eq:Check_FL}
\end{equation}
by Theorem \ref{thm:WLP_Uni} (a, b) (with $d=\aleph(\mathbf{T}_{0})$,
$\Gamma=\{\mu_{\mathbf{T}_{0},n}\}_{n\in\mathbf{N}}$ and $\mathcal{D}=C_{b}(E;\mathbf{R})$).
Now, (c) follows by (\ref{eq:Proc_FDD_BExt_Uni}), (\ref{eq:Check_FL})
and Fact \ref{fact:FLP_Uni} (with $\mathbf{I}=\mathbf{N}$).\end{proof}

\begin{rem}
\label{rem:RepProc_FL}As mentioned in Note \ref{note:Ehat_Valued_Proc_FDD},
Polish-space-valued processes (especially replica processes) have
Borel finite-dimensional distributions and their finite-dimensional
convergence refers exactly to the weak convergence of their finite-dimensional
distributions. Moreover, $\mathcal{P}(\widehat{E})$ is a Polish space
by Corollary \ref{cor:Base_Sep_Meas} (c). Hence, (\ref{eq:RepProc_FC_along_T_Y})
is equivalent to $Y=\mathfrak{fl}_{\mathbf{T}}(\{\widehat{X}^{n}\}_{n\in\mathbf{N}})$.
\end{rem}

The next corollary leverages Theorem \ref{thm:TransFC_1} to establish
finite-dimensional convergence to a given limit process.
\begin{cor}
\label{cor:TransFC_2}Let $E$ be a topological space, $\{(\Omega^{n},\mathscr{F}^{n},\mathbb{P}^{n};X^{n})\}_{n\in\mathbf{N}}$
be $E$-valued processes, $(E_{0},\mathcal{F};\widehat{E},\widehat{\mathcal{F}})$
be a base over $E$ and $\mathbf{T}\subset\mathbf{R}^{+}$. Suppose
that:

\renewcommand{\labelenumi}{(\roman{enumi})}
\begin{enumerate}
\item $\{X_{t}^{n}\}_{n\in\mathbf{N}}$ is sequentially tight in $E_{0}$
for all $t\in\mathbf{T}$.
\item (\ref{eq:Exp_Test}) holds for all $f\in\mathfrak{mc}[\Pi^{\mathbf{T}_{0}}(\mathcal{F}\backslash\{1\})]$
and $\mathbf{T}_{0}\in\mathscr{P}_{0}(\mathbf{T})$.
\item $E$-valued process $(\Omega,\mathscr{F},\mathbb{P};X)$ satisfies
(\ref{eq:FR-Base}) (or (\ref{eq:T-Base})).
\end{enumerate}
Then:

\renewcommand{\labelenumi}{(\alph{enumi})}
\begin{enumerate}
\item (\ref{eq:Proc_FDD_BExt_Uni}), (\ref{eq:Lim_Proc_FDD_Borel}), (\ref{eq:FC_along_T_E0})
and (\ref{eq:FC_along_T}) hold for some $\{N_{\mathbf{T}_{0}}\}_{\mathbf{T}_{0}\in\mathscr{P}_{0}(\mathbf{T})}\subset\mathbf{N}$.
Moreover, $X_{\mathbf{T}_{0}}$ is $\mathbf{m}$-tight for all $\mathbf{T}_{0}\in\mathscr{P}_{0}(\mathbf{T})$.
\item If $\mathbf{T}=\mathbf{R}^{+}$ and $\{X^{n}\}_{n\in\mathbf{N}}$
is $(\mathbf{R}^{+},\mathcal{F}\backslash\{1\})$-AS, then $X$ is
stationary.
\item If $C_{b}(E;\mathbf{R})$ separates points on $E$, then $X=\mathfrak{fl}_{\mathbf{T}}(\{X_{n}\}_{n\in\mathbf{N}})$.
\end{enumerate}
\end{cor}
\begin{proof}
Letting $Y=\widehat{X}$, we obtain (\ref{eq:Exp_Test_Y}) for all
$f\in\mathfrak{mc}[\Pi^{\mathbf{T}_{0}}(\mathcal{F}\backslash\{1\})]$
and $\mathbf{T}_{0}\in\mathscr{P}_{0}(\mathbf{T})$ by the conditions
(ii) and (iii) above, Fact \ref{fact:T-Base_FR-Base} and Proposition
\ref{prop:RepProc_FR-Base} (a). Now, the result follows by Theorem
\ref{thm:TransFC_1} immediately.\end{proof}

\section{\label{sec:Cadlag_RepProc}C$\grave{\mbox{a}}$dl$\grave{\mbox{a}}$g
replica}

Properties like relative compactness, tightness, well-posedness of
martingale problems and convergence of nonlinear filters are often
easily established for c$\grave{\mbox{a}}$dl$\grave{\mbox{a}}$g
replicas. As Polish-space-valued c$\grave{\mbox{a}}$dl$\grave{\mbox{a}}$g
processes, they have the following nice measurability.
\begin{fact}
\label{fact:Cadlag_RepProc}Let $E$ be a topological space, $(E_{0},\mathcal{F};\widehat{E},\widehat{\mathcal{F}})$
be a base over $E$ and $(\Omega,\mathscr{F},\mathbb{P};X)$ be an
$E$-valued process. Then, $\mathfrak{rep}_{\mathrm{c}}(X;E_{0},\mathcal{F})\subset M(\Omega,\mathscr{F};D(\mathbf{R}^{+};\widehat{E}))$%
\footnote{$M(\Omega,\mathscr{F};D(\mathbf{R}^{+};\widehat{E}))$ denotes the
family of all $D(\mathbf{R}^{+};\widehat{E})$-valued random variables
defined on measurable space $(\Omega,\mathscr{F})$.%
} and $\mathfrak{rep}_{\mathrm{c}}(X;E_{0},\mathcal{F})\subset\mathfrak{rep}_{\mathrm{p}}(X;E_{0},\mathcal{F})\subset\mathfrak{rep}_{\mathrm{m}}(X;E_{0},\mathcal{F})$%
\footnote{$\mathfrak{rep}_{\mathrm{m}}(X;E_{0},\mathcal{F})$, $\mathfrak{rep}_{\mathrm{p}}(X;E_{0},\mathcal{F})$
and $\mathfrak{rep}_{\mathrm{c}}(X;E_{0},\mathcal{F})$, introduced
in Notation \ref{notation:RepProc}, are all equivalence classes of
measurable, progressive and c$\grave{\mbox{a}}$dl$\grave{\mbox{a}}$g
replicas of $X$ with respect to $(E_{0},\mathcal{F};\widehat{E},\widehat{\mathcal{F}})$.%
}.
\end{fact}
\begin{proof}
The first statement follows by Fact \ref{fact:Cadlag_Sko_RV} (b)
(with $E=\widehat{E}$ and $X=\widehat{X}\in\mathfrak{rep}_{\mathrm{c}}(X;E_{0},\mathcal{F})$)
and Fact \ref{fact:Sko_RV_Cadlag} (a) (with $E=\widehat{E}$). The
second statement follows by Proposition \ref{prop:Proc_Basic_2} (a,
c).\end{proof}

Due to the topological difference of $E$ and $\widehat{E}$, non-c$\grave{\mbox{a}}$dl$\grave{\mbox{a}}$g
$E$-valued processes can have c$\grave{\mbox{a}}$dl$\grave{\mbox{a}}$g
replicas if they are almost c$\grave{\mbox{a}}$dl$\grave{\mbox{a}}$g
on the functions in $\mathcal{F}$.
\begin{defn}
\label{def:Weakly_Cadlag}Let $E$ be a topological space and $(\Omega,\mathscr{F},\mathbb{P};X)$
be an $E$-valued process. $X$ is said to be \textbf{weakly c$\grave{\mbox{a}}$dl$\grave{\mbox{a}}$g
along $\mathbf{T}$ under $\mathcal{D}$} (\textbf{$(\mathbf{T},\mathcal{D})$}-c$\grave{\mbox{a}}$dl$\grave{\mbox{a}}$g
for short) if: (1) $\mathbf{T}\subset\mathbf{R}^{+}$ and $\mathcal{D}\subset M(E;\mathbf{R})$
are non-empty, and (2) There exist $\mathbf{R}$-valued c$\grave{\mbox{a}}$dl$\grave{\mbox{a}}$g
processes $\{(\Omega,\mathscr{F},\mathbb{P};\zeta^{f})\}_{f\in\mathcal{D}}$%
\footnote{Please be reminded that all processes in this article are indexed
by $\mathbf{R}^{+}$.%
} such that
\begin{equation}
\inf_{f\in\mathcal{D},t\in\mathbf{T}}\mathbb{P}\left(f\circ X_{t}=\zeta_{t}^{f}\right)=1.\label{eq:(T,D)_Weakly_Cadlag}
\end{equation}
\end{defn}
\begin{note}
\label{note:Weakly_Cadlag}$X$ is $(\mathbf{R}^{+},\mathcal{D})$-c$\grave{\mbox{a}}$dl$\grave{\mbox{a}}$g
if $\{\varpi(f)\circ X\}_{f\in\mathcal{D}}$ are all c$\grave{\mbox{a}}$dl$\grave{\mbox{a}}$g
processes, especially if $X$ is c$\grave{\mbox{a}}$dl$\grave{\mbox{a}}$g
and $\mathcal{D}\subset C(E;\mathbf{R})$ by Fact \ref{fact:Cadlag_Proc}
(a) (with $S=\mathbf{R}$). Apparently, $(\mathbf{R}^{+},\mathcal{D})$-c$\grave{\mbox{a}}$dl$\grave{\mbox{a}}$g
property is transitive between modifications.\end{note}
\begin{rem}
\label{rem:Weakly_Cadlag}A special case of an $(\mathbf{R}^{+},\mathcal{D}$)-c$\grave{\mbox{a}}$dl$\grave{\mbox{a}}$g
$X$ is when $\varpi(f)\circ X$ is c$\grave{\mbox{a}}$dl$\grave{\mbox{a}}$g
for all $f\in\mathcal{D}$.
\end{rem}

Here are three typical sufficient conditions for unique existence
of c$\grave{\mbox{a}}$dl$\grave{\mbox{a}}$g replica.
\begin{prop}
\label{prop:FR}Let $E$ be a topological space, $(E_{0},\mathcal{F};\widehat{E},\widehat{\mathcal{F}})$
be a base over $E$, $(\Omega,\mathscr{F},\mathbb{P};X)$ be an $E$-valued
process and $\mathbf{T}\subset\mathbf{R}^{+}$ be dense. Then:

\renewcommand{\labelenumi}{(\alph{enumi})}
\begin{enumerate}
\item If $X$ is $(\mathbf{R}^{+},\mathcal{F}$)-c$\grave{\mbox{a}}$dl$\grave{\mbox{a}}$g
and (\ref{eq:FR-Base}) holds, then $\widehat{X}=\mathfrak{rep}_{\mathrm{c}}(X;E_{0},\mathcal{F})$%
\footnote{We specified in Notation \ref{notation:RepProc} that ``$\widehat{X}=\mathfrak{rep}_{\mathrm{c}}(X;E_{0},\mathcal{F})$''
means $\widehat{X}$ is the unique c$\grave{\mbox{a}}$dl$\grave{\mbox{a}}$g
replica of $X$ up to indistinguishability. ``Unique up to indistinguishability''
means any two processes with the relevant property is indistinguishable.%
} exists and $\varpi(\bigotimes\widehat{\mathcal{F}})\circ\widehat{X}$
(resp. $\varpi(\widehat{f})\circ\widehat{X}$) is the unique c$\grave{\mbox{a}}$dl$\grave{\mbox{a}}$g
modification%
\footnote{The terminology ``modification'' was explained in \S \ref{sec:Proc}.%
} of $\varpi(\bigotimes\mathcal{F})\circ X$ (resp. $\varpi(f)\circ X$
for each $f\in\mathcal{F}$) up to indistinguishability.
\item If $\{\varpi(f)\circ X\}_{f\in\mathcal{F}}$ are all c$\grave{\mbox{a}}$dl$\grave{\mbox{a}}$g
and (\ref{eq:FR-Base}) holds, then $\widehat{X}=\mathfrak{rep}_{\mathrm{c}}(X;E_{0},\mathcal{F})$
exists and%
\footnote{Please be reminded that $\widehat{f}$ denotes the continuous replica
of $f$.%
}
\begin{equation}
\mathbb{P}\left(\varpi(f)\circ X=\varpi(\widehat{f})\circ\widehat{X},\forall f\in\mathfrak{ca}(\mathcal{F})\right)=1.\label{eq:ca(F)_Ind_FR}
\end{equation}

\item If $X$ is c$\grave{\mbox{a}}$dl$\grave{\mbox{a}}$g and satisfies
(\ref{eq:T-Base}), then $\widehat{X}=\mathfrak{rep}_{\mathrm{c}}(X;E_{0},\mathcal{F})$
exists and%
\footnote{The replica of functions was discussed in \S \ref{sec:RepFun}.%
}
\begin{equation}
\mathbb{P}\left(\varpi(f)\circ X=\varpi(\widehat{f})\circ\widehat{X},\forall f\in C(E;\mathbf{R})\mbox{ that has a replica }\widehat{f}\right)=1.\label{eq:All_ContRep_Ind}
\end{equation}

\end{enumerate}
\end{prop}
\begin{rem}
\label{rem:TFInd}The \textit{functional} indistinguishability of
$X$ and $\widehat{X}$ in (\ref{eq:ca(F)_Ind_FR}) is a valuable
property of c$\grave{\mbox{a}}$dl$\grave{\mbox{a}}$g replica. Corollary
\ref{cor:Base_Fun_Dense} showed $C(\widehat{E};\mathbf{R})=\mathfrak{ca}(\widehat{\mathcal{F}})=\{\widehat{f}:f\in\mathfrak{ca}(\mathcal{F})\}$,
so (\ref{eq:ca(F)_Ind_FR}) allows many properties to be transferred
between $X$ and $\widehat{X}$.
\end{rem}

Our construction of $\widehat{X}$ is based on the following technical
lemma.
\begin{lem}
\label{lem:FR_Construct}Let $E$ be a topological space, $(E_{0},\mathcal{F};\widehat{E},\widehat{\mathcal{F}})$
be a base over $E$, $(\Omega,\mathscr{F},\mathbb{P};X)$ be an $E$-valued
process satisfying (\ref{eq:FR-Base}) for some dense $\mathbf{T}\subset\mathbf{R}^{+}$
and $\mathbf{T}\subset\mathbf{S}\subset\mathbf{R}^{+}$. Then, the
following statements are equivalent:

\renewcommand{\labelenumi}{(\alph{enumi})}
\begin{enumerate}
\item $X$ is $(\mathbf{S},\mathcal{F})$-c$\grave{\mbox{a}}$dl$\grave{\mbox{a}}$g.
\item There exists an $\mathbf{R}^{\infty}$-valued c$\grave{\mbox{a}}$dl$\grave{\mbox{a}}$g
process $(\Omega,\mathscr{F},\mathbb{P};\zeta)$ such that
\begin{equation}
\inf_{t\in\mathbf{S}}\mathbb{P}\left(\bigotimes\mathcal{F}\circ X_{t}=\zeta_{t}\right)=1.\label{eq:Phi(X)_S_Mod}
\end{equation}

\item There exists an $\widehat{X}\in M(\Omega,\mathscr{F};D(\mathbf{R}^{+};\widehat{E}))$
such that
\begin{equation}
\inf_{t\in\mathbf{S}}\mathbb{P}\left(\bigotimes\mathcal{F}\circ X_{t}=\bigotimes\widehat{\mathcal{F}}\circ\widehat{X}_{t}\right)=1.\label{eq:Phi(X)_Phihat(Xhat)_S-Mod}
\end{equation}

\end{enumerate}
\end{lem}
\begin{proof}
((a) $\rightarrow$ (b)) is immediate by Fact \ref{fact:Weakly_Cadlag}
(with $\mathcal{D}=\mathcal{F}$ and $\mathbf{T}=\mathbf{S}$).

((b) $\rightarrow$ (c)) Let $\mathbf{T}_{0}\subset\mathbf{T}$ be
countable and dense in $\mathbf{R}^{+}$. $\bigotimes\widehat{\mathcal{F}}(\widehat{E})$
is a closed subspace of $\mathbf{R}^{\infty}$ by (\ref{eq:Prod(Fhat)(Ehat)_Compact_Rinf}).
\begin{equation}
\begin{aligned} & \mathbb{P}\left[\zeta\in D\left(\mathbf{R}^{+};\bigotimes\widehat{\mathcal{F}}(\widehat{E})\right)\right]\\
 & \geq\mathbb{P}\left(\zeta_{t}=\bigotimes\mathcal{F}\circ X_{t}\in\bigotimes\widehat{\mathcal{F}}(\widehat{E}),\forall t\in\mathbf{T}_{0}\right)=1
\end{aligned}
\label{eq:All_Path_in_Phihat(Ehat)}
\end{equation}
by (\ref{eq:Phi(X)_S_Mod}), $\mathbf{T}_{0}\subset\mathbf{S}$, (\ref{eq:FR-Base}),
the c$\grave{\mbox{a}}$dl$\grave{\mbox{a}}$g property of $\zeta$
and the closedness of $\bigotimes\widehat{\mathcal{F}}(\widehat{E})$.
Then, there exists a
\begin{equation}
\zeta^{\prime}\in M\left[\Omega,\mathscr{F};D\left(\mathbf{R}^{+};\bigotimes\widehat{\mathcal{F}}(\widehat{E})\right)\right]\label{eq:Construct_FR_1}
\end{equation}
satisfying
\begin{equation}
\mathbb{P}\left(\zeta=\zeta^{\prime}\right)=1\label{eq:Construct_FR_2}
\end{equation}
by (\ref{eq:All_Path_in_Phihat(Ehat)}), Proposition \ref{prop:Sko_Basic_2}
(b) (with $E=\mathbf{R}^{\infty}$) and Lemma \ref{lem:Path_Space_RV}
(b) (with $E=\mathbf{R}^{\infty}$, $E_{0}=\bigotimes\widehat{\mathcal{F}}(\widehat{E})$,
$S_{0}=D(\mathbf{R}^{+};\bigotimes\widehat{\mathcal{F}}(\widehat{E}))$
and $X=\zeta$).

It follows by by (\ref{eq:Base_Imb}) and Proposition \ref{prop:Sko_Basic_1}
(d) (with $S=\bigotimes\widehat{\mathcal{F}}(\widehat{E})$, $E=\widehat{E}$
and $f=(\bigotimes\widehat{\mathcal{F}})^{-1}$) that
\begin{equation}
\varpi\left[\left(\bigotimes\widehat{\mathcal{F}}\right)^{-1}\right]\in C\left[D\left(\mathbf{R}^{+};\bigotimes\widehat{\mathcal{F}}(\widehat{E})\right);D(\mathbf{R}^{+};\widehat{E})\right].\label{eq:Construct_FR_3}
\end{equation}
It follows by (\ref{eq:Construct_FR_1}) and (\ref{eq:Construct_FR_3})
that
\begin{equation}
\widehat{X}\circeq\varpi\left[\left(\bigotimes\widehat{\mathcal{F}}\right)^{-1}\right]\circ\zeta^{\prime}\in M\left(\Omega,\mathscr{F};D(\mathbf{R}^{+};\widehat{E})\right).\label{eq:Construct_FR_4}
\end{equation}
It follows by (\ref{eq:Construct_FR_2}), (\ref{eq:Construct_FR_4})
and (\ref{eq:Base_Imb}) that
\begin{equation}
\mathbb{P}\left(\zeta=\varpi\left(\bigotimes\widehat{\mathcal{F}}\right)\circ\varpi\left[\left(\bigotimes\widehat{\mathcal{F}}\right)^{-1}\right]\circ\zeta^{\prime}=\varpi\left(\bigotimes\widehat{\mathcal{F}}\right)\circ\widehat{X}\right)=1.\label{eq:Check_Phihat(Xhat)_Ind_Phi(X)_Mod}
\end{equation}
Now, (\ref{eq:Phi(X)_Phihat(Xhat)_S-Mod}) follows by (\ref{eq:Phi(X)_S_Mod})
and (\ref{eq:Check_Phihat(Xhat)_Ind_Phi(X)_Mod}).

((c) $\rightarrow$ (a)) is automatic.\end{proof}

\begin{proof}
[Proof of Proposition \ref{prop:FR}](a) By Lemma \ref{lem:FR_Construct}
(a - c) (with $\mathbf{S}=\mathbf{R}^{+}$), there exists an $\widehat{X}\in M(\Omega,\mathscr{F};D(\mathbf{R}^{+};\widehat{E}))$
such that
\begin{equation}
\inf_{t\in\mathbf{R}^{+}}\mathbb{P}\left(\bigotimes\mathcal{F}\circ X_{t}=\bigotimes\widehat{\mathcal{F}}\circ\widehat{X}_{t}\in\bigotimes\widehat{\mathcal{F}}(\widehat{E})\right)=1.\label{eq:Check_FR}
\end{equation}
$\varpi(\bigotimes\widehat{\mathcal{F}})\circ\widehat{X}$ (resp.
$\varpi(\widehat{f})\circ\widehat{X}$) is a c$\grave{\mbox{a}}$dl$\grave{\mbox{a}}$g
modification of $\varpi(\bigotimes\mathcal{F})\circ X$ (resp. $\varpi(f)\circ X$
for each $f\in\mathcal{F}$) by (\ref{eq:Base_Imb}), the fact $\widehat{\mathcal{F}}\subset C(\widehat{E};\mathbf{R})$,
(\ref{eq:Check_FR}) and Fact \ref{fact:Cadlag_Proc} (a) (with $E=\widehat{E}$,
$X=\widehat{X}$ and $f=\widehat{f}$ or $\bigotimes\widehat{\mathcal{F}}$).
Now, (a) follows by Proposition \ref{prop:RepProc_FR-Base} (b) and
Proposition \ref{prop:Proc_Basic_2} (h).

(b) Given any $\widehat{X}\in\mathfrak{rep}_{\mathrm{c}}(X;E_{0},\mathcal{F})$,
one finds that
\begin{equation}
\begin{aligned} & \left\{ \omega\in\Omega:\varpi(f)\circ X(\omega)=\varpi(\widehat{f})\circ\widehat{X}(\omega),\forall f\in\mathfrak{ca}(\mathcal{F})\right\} \\
 & =\left\{ \omega\in\Omega:\varpi\left(\bigotimes\mathcal{F}\right)\circ X(\omega)=\varpi\left(\bigotimes\widehat{\mathcal{F}}\right)\circ\widehat{X}(\omega)\right\} 
\end{aligned}
\label{eq:Check_ca(F)_Ind_Rep}
\end{equation}
by properties of uniform convergence. Then, (b) follows by (\ref{eq:Check_ca(F)_Ind_Rep})
and (a).

(c) Let $\mathbf{T}_{0}\subset\mathbf{T}$ be countable and dense
in $\mathbf{R}^{+}$. Given a c$\grave{\mbox{a}}$dl$\grave{\mbox{a}}$g
$X$, $\widehat{X}\in\mathfrak{rep}_{\mathrm{c}}(X;E_{0},\mathcal{F})$
exists by (a), and $\varpi(f)\circ X$ and $\varpi(\widehat{f})\circ\widehat{X}$
are c$\grave{\mbox{a}}$dl$\grave{\mbox{a}}$g process for all $f\in C(E;\mathbf{R})$
by Fact \ref{fact:Cadlag_Proc} (a).
\begin{equation}
\begin{aligned} & \left\{ \omega\in\Omega:\varpi(f)\circ X(\omega)=\varpi(\widehat{f})\circ\widehat{X}(\omega),\forall f\in C(E;\mathbf{R})\mbox{ having a replica }\widehat{f}\right\} \\
 & \supset\left\{ \omega\in\Omega:X_{t}(\omega)=\widehat{X}_{t}(\omega)\in E_{0},\forall t\in\mathbf{T}_{0}\right\} 
\end{aligned}
\label{eq:Check_All_ContRep_Ind}
\end{equation}
by the fact $f|_{E_{0}}=\widehat{f}|_{E_{0}}$ and Proposition \ref{prop:Proc_Basic_2}
(g). Now, (c) follows by Fact \ref{fact:T-Base_FR-Base}, (b), Proposition
\ref{prop:RepProc_T-Base} (a) and (\ref{eq:Check_All_ContRep_Ind}).\end{proof}

\begin{rem}
\label{rem:NonCad_Cad_Rep}The c$\grave{\mbox{a}}$dl$\grave{\mbox{a}}$g
property of $\varpi(\bigotimes\mathcal{F})\circ X(\omega)$ in $(\mathbf{R}^{\infty})^{\mathbf{R}^{+}}$
does not guarantee that of $X(\omega)$ in $E^{\mathbf{R}^{+}}$ since
$\bigotimes\mathcal{F}$ is not necessarily an imbedding on $E$.
\end{rem}

If $E_{0}$ is large enough for $X$ to almost surely live in $E_{0}^{\mathbf{R}^{+}}$,
then the following result shows one can modify merely a $\mathbb{P}$-negligible
amount of paths of $X$ and obtain an indistinguishable replica of
$X$.
\begin{prop}
\label{prop:PR}Let $E$ be a topological space, $E_{0}\in\mathscr{B}(E)$,
$S_{0}\subset E_{0}^{\mathbf{R}^{+}}$ and $(\Omega,\mathscr{F},\mathbb{P};X)$
be an $E$-valued process satisfying
\begin{equation}
\mathbb{P}\left(X\in S_{0}\right)=1.\label{eq:Path-Base}
\end{equation}
Then, there exists an $\widehat{X}\in S_{0}^{\Omega}$ satisfying
the following properties:

\renewcommand{\labelenumi}{(\alph{enumi})}
\begin{enumerate}
\item $\widehat{X}$ is an $(E_{0},\mathscr{O}_{E}(E_{0}))$-valued process
and $\mathbb{P}(X=\widehat{X}\in S_{0})=1$.
\item $\widehat{X}\in\mathfrak{rep}(X;E_{0},\mathcal{F})$ for any base
$(E_{0},\mathcal{F};\widehat{E},\widehat{\mathcal{F}})$ over $E$.
\item If every element of $S_{0}$ is a c$\grave{\mbox{a}}$dl$\grave{\mbox{a}}$g
member of $(E_{0}^{\mathbf{R}^{+}},\mathscr{O}_{E}(E_{0})^{\mathbf{R}^{+}})$%
\footnote{$(E_{0},\mathscr{O}_{E}(E_{0}))$ need not be a Tychonoff space, so
we avoid the notation $D(\mathbf{R}^{+};E_{0},\mathscr{O}_{E}(E_{0}))$.%
}, then $\widehat{X}=\mathfrak{rep}_{\mathrm{c}}(X;E_{0},\mathcal{F})$
for any base $(E_{0},\mathcal{F};\widehat{E},\widehat{\mathcal{F}})$
over $E$.
\end{enumerate}
\end{prop}
\begin{proof}
We fix $y_{0}\in S_{0}$ and define $\widehat{X}\circeq\mathfrak{var}(X;\Omega,X^{-1}(S_{0}),y_{0})$%
\footnote{$\mathfrak{var}(\cdot)$ was introduced in Notation \ref{notation:Var}.%
}. Then, (a) follows by Lemma \ref{lem:var(X)} (b, c) (with $(E,\mathscr{U})=(E^{\mathbf{R}^{+}},\mathscr{B}(E)^{\otimes\mathbf{R}^{+}})$,
$S=S_{0}$, $\mathscr{U}^{\prime}=\mathscr{U}|_{S_{0}}$ and $Y=\widehat{X}$).
Given any base $(E_{0},\mathcal{F};\widehat{E},\widehat{\mathcal{F}})$,
$(E_{0},\mathscr{O}_{\widehat{E}}(E_{0}))$ is coarser than $(E_{0},\mathscr{O}_{E}(E_{0}))$
by Lemma \ref{lem:Base} (d) and so $\widehat{X}$ is an $\widehat{E}$-valued
process. We have by (\ref{eq:F_Fhat_Coincide}) that
\begin{equation}
\begin{aligned} & \inf_{t\in\mathbf{R}^{+}}\mathbb{P}\left(\bigotimes\mathcal{F}\circ X_{t}=\bigotimes\widehat{\mathcal{F}}\circ\widehat{X}_{t}\in\bigotimes\widehat{\mathcal{F}}(\widehat{E})\right)\\
 & \geq\mathbb{P}\left(X_{t}=\widehat{X}_{t}\in E_{0},\forall t\in\mathbf{R}^{+}\right)\geq\mathbb{P}\left(X=\widehat{X}\in S_{0}\right)=1,
\end{aligned}
\label{eq:Check_PR}
\end{equation}
thus proving (b). The c$\grave{\mbox{a}}$dl$\grave{\mbox{a}}$g property
of $\widehat{X}(\omega):\mathbf{R}^{+}\rightarrow(E_{0},\mathscr{O}_{E}(E_{0}))$
implies $\widehat{X}(\omega)\in D(\mathbf{R}^{+};\widehat{E})$ for
all $\omega\in\Omega$%
\footnote{This statement is stronger than being an $\widehat{E}$-valued c$\grave{\mbox{a}}$dl$\grave{\mbox{a}}$g
process.%
} by Fact \ref{fact:Cadlag_Path} (b) (with $E=(E_{0},\mathscr{O}_{E}(E_{0}))$,
$S=(E_{0},\mathscr{O}_{\widehat{E}}(E_{0}))$ and $f$ being the identity
mapping on $E_{0}$). Hence, (c) follows by (b), (\ref{eq:Check_PR})
and Proposition \ref{prop:RepProc_FR-Base} (b) (with $\mathbf{T}=\mathbf{R}^{+}$).\end{proof}

\begin{rem}
\label{rem:E0_Size}In general, many paths of a c$\grave{\mbox{a}}$dl$\grave{\mbox{a}}$g
replica could live outside $E_{0}^{\mathbf{R}^{+}}$. Such a c$\grave{\mbox{a}}$dl$\grave{\mbox{a}}$g
replica is not necessarily an $E$-valued process, nor is it (necessarily)
indistinguishable from $X$. Specifically, if $X$ is c$\grave{\mbox{a}}$dl$\grave{\mbox{a}}$g
and satisfies
\begin{equation}
\inf_{t\in\mathbf{R}^{+}}\mathbb{P}(X_{t}\in E_{0})=1,\label{eq:R+-Base}
\end{equation}
then $\widehat{X}=\mathfrak{rep}_{\mathrm{c}}(X;E_{0},\mathcal{F})$
satisfies
\begin{equation}
\inf_{t\in\mathbf{R}^{+}}\mathbb{P}\left(X_{t}=\widehat{X}_{t}\in E_{0}\right)=1\label{eq:RepProc_(R+,E0)-Mod}
\end{equation}
by Proposition \ref{prop:FR} (c) (with $\mathbf{T}=\mathbf{R}^{+}$)
and Proposition \ref{prop:RepProc_T-Base} (a) (with $\mathbf{T}=\mathbf{R}^{+}$).
However, this does not necessarily imply $\mathbb{P}(X\in E_{0}^{\mathbf{R}^{+}})=1$
nor $\mathbb{P}(X=\widehat{X}\in E_{0}^{\mathbf{R}^{+}})=1$ since
$E_{0}$ might not be a closed subspace of $E$ or $\widehat{E}$.
\end{rem}

Moreover, we transform $\mathcal{M}^{+}(E)$-valued weakly c$\grave{\mbox{a}}$dl$\grave{\mbox{a}}$g%
\footnote{The notion of weakly c$\grave{\mbox{a}}$dl$\grave{\mbox{a}}$g was
introduced in Definition \ref{def:Weakly_Cadlag}.%
} processes into $\mathcal{P}(\widehat{E})$-valued c$\grave{\mbox{a}}$dl$\grave{\mbox{a}}$g
processes by similar construction techniques for c$\grave{\mbox{a}}$dl$\grave{\mbox{a}}$g
replica, which furthers Proposition \ref{prop:M(E)_Proc_Rep}.
\begin{lem}
\label{lem:P(E)_Proc_Rep_Cadlag}Let $E$ be a topological space,
$(E_{0},\mathcal{F};\widehat{E},\widehat{\mathcal{F}})$ be a base
over $E$, $(\Omega,\mathscr{F},\mathbb{P};X)$ be an $\mathcal{M}^{+}(E)$-valued
$(\mathbf{R}^{+},\mathfrak{mc}(\mathcal{F})^{*})$-c$\grave{\mbox{a}}$dl$\grave{\mbox{a}}$g
process satisfying (\ref{eq:Almostsure_Support_Proc}) and (\ref{eq:R+-Base_P(E)}).
Then, there exists an $\mathscr{F}_{t}^{X}$-adapted $D(\mathbf{R}^{+};\mathcal{P}(\widehat{E}))$-valued
random variable $(\Omega,\mathscr{F},\mathbb{P};Y)$ satisfying Proposition
\ref{prop:M(E)_Proc_Rep} (a, c, d).\end{lem}
\begin{rem}
\label{rem:P(Ehat)_Sko_RV}In the proof below, we let $\varphi$,
$\widehat{\varphi}$, $y_{0}$ and $\Psi$ be as in Lemma \ref{lem:M(E)_RV_Rep}
and set $S_{0}=\widehat{\varphi}[\mathcal{P}(\widehat{E})]$. Recall
that $\widehat{\varphi}$ satisfies (\ref{eq:Check_M(E)_Rep_2}).
$\mathcal{P}(\widehat{E})$ is a compact Polish space by Corollary
\ref{cor:Base_Sep_Meas} (c) (with $d=1$). Hence, $S_{0}\in\mathscr{C}(\mathbf{R}^{\infty})$
is a Polish subspace of $\mathbf{R}^{\infty}$ by Proposition \ref{prop:Compact}
(a, e) and Proposition \ref{prop:Var_Polish} (b, f). $D(\mathbf{R}^{+};\mathcal{P}(\widehat{E}))$
and $D(\mathbf{R}^{+};S_{0},\mathscr{O}_{\mathbf{R}^{\infty}}(S_{0}))$
are Polish spaces by Proposition \ref{prop:Sko_Basic_2} (d) (with
$E=\mathcal{P}(\widehat{E})$ or $(S_{0},\mathscr{O}_{\mathbf{R}^{\infty}}(S_{0}))$).
Therefore, $D(\mathbf{R}^{+};\mathcal{P}(\widehat{E}))$-valued and
$D(\mathbf{R}^{+};S_{0},\mathscr{O}_{\mathbf{R}^{\infty}}(S_{0}))$-valued
random variables are c$\grave{\mbox{a}}$dl$\grave{\mbox{a}}$g processes
by Fact \ref{fact:Sko_RV_Cadlag} (a), for which $\mathscr{F}_{t}^{X}$-adaptedness
is a proper concept.
\end{rem}
\begin{proof}
[Proof of Lemma \ref{lem:P(E)_Proc_Rep_Cadlag}]The proof of Lemma
\ref{lem:M(E)_RV_Rep} mentioned that $\mathfrak{mc}(\mathcal{F})$
is countable and $\varphi$ satisfies (\ref{eq:Check_M(E)_Rep_1}),
so $\varphi\circ X$ has an $\mathbf{R}^{\infty}$-valued c$\grave{\mbox{a}}$dl$\grave{\mbox{a}}$g
modification $\zeta$ by the converse part of Fact \ref{fact:Weakly_Cadlag}
(with $E=\mathcal{M}^{+}(E)$, $\mathcal{D}=\mathfrak{mc}(\mathcal{F})^{*}$
and $\mathbf{T}=\mathbf{R}^{+}$) and $\zeta$ is $\mathscr{F}_{t}^{X}$-adapted
by Proposition \ref{prop:Proc_Basic_2} (e).

Let $\nu_{t}^{\omega}$ be the replica of $X_{t}(\omega)$ for each
fix $\omega\in\Omega$ and $t\in\mathbf{R}^{+}$. Similar to (\ref{eq:Check_M(E)_Rep_7}),
we have that
\begin{equation}
\begin{aligned} & \inf_{t\in\mathbf{R}^{+}}\mathbb{P}\left(\zeta_{t}=\varphi\circ X_{t}\in S_{0}\right)\\
 & =\inf_{t\in\mathbf{R}^{+}}\mathbb{P}\left(\left\{ \omega\in\Omega:\varphi\circ X_{t}(\omega)=\widehat{\varphi}(\nu_{t}^{\omega})\in S_{0}\right\} \right)\\
 & =\inf_{t\in\mathbf{R}^{+}}\mathbb{P}\left(\left\{ \omega\in\Omega:X_{t}(\omega)(E)=X_{t}(\omega)(E_{0})=1\right\} \right)=1
\end{aligned}
\label{eq:R+-FR_Base_M(E)}
\end{equation}
by the countability of $\mathfrak{mc}(\mathcal{F})$, Proposition
\ref{prop:RepMeas_Basic} (a, b, e) (with $d=1$, $\mu=X_{t}(\omega)$
and $\overline{\mu}=\nu^{\omega}$), (\ref{eq:Almostsure_Support_Proc})
and (\ref{eq:R+-Base_P(E)}). It follows that
\begin{equation}
\mathbb{P}\left(\zeta\in D\left(\mathbf{R}^{+};S_{0},\mathscr{O}_{\mathbf{R}^{\infty}}(S_{0})\right)\right)=1\label{eq:M(E)_Proc_Rep_Path_Contain}
\end{equation}
by (\ref{eq:R+-FR_Base_M(E)}), the closedness of $S_{0}$ and the
c$\grave{\mbox{a}}$dl$\grave{\mbox{a}}$g property of $\zeta$. Then,
there exists a
\begin{equation}
\zeta^{\prime}\in M\left[\Omega,\mathscr{F};D\left(\mathbf{R}^{+};S_{0},\mathscr{O}_{\mathbf{R}^{\infty}}(S_{0})\right)\right]\label{eq:Construct_P(E)_Proc_Rep_Cadlag}
\end{equation}
satisfying (\ref{eq:Construct_FR_2}) by (\ref{eq:M(E)_Proc_Rep_Path_Contain}),
Proposition \ref{prop:Sko_Basic_2} (b) (with $E=\mathbf{R}^{\infty}$)
and Lemma \ref{lem:Path_Space_RV} (b) (with $E=E_{0}=\mathbf{R}^{\infty}$,
$E_{0}=S_{0}$, $S_{0}=D(\mathbf{R}^{+};S_{0},\mathscr{O}_{\mathbf{R}^{\infty}}(S_{0}))$
and $X=\zeta$). As $\zeta$ is $\mathscr{F}_{t}^{X}$-adapted, $\zeta^{\prime}$
is $\mathscr{F}_{t}^{X}$-adapted by (\ref{eq:Construct_FR_2}) and
Proposition \ref{prop:Proc_Basic_2} (e). Furthermore,
\begin{equation}
\inf_{t\in\mathbf{R}^{+}}\mathbb{P}\left(\zeta_{t}^{\prime}=\varphi\circ X_{t}\in S_{0}\right)=1.\label{eq:Check_P(E)_Proc_Rep_Cadlag_TF_Mod_1}
\end{equation}
by (\ref{eq:Construct_FR_2}) and (\ref{eq:R+-FR_Base_M(E)}).

The proof of Lemma \ref{lem:M(E)_RV_Rep} mentioned that $\widehat{\varphi}$
satisfies (\ref{eq:Check_M(E)_Rep_2}) and $\Psi$ equals $\widehat{\varphi}^{-1}$
restricted to $S_{0}$. Hence, we have that: (1)
\begin{equation}
Y\circeq\varpi(\Psi)\circ\zeta^{\prime}=\varpi(\widehat{\varphi}^{-1})\circ\zeta^{\prime}\in M\left[\Omega,\mathscr{F};D\left(\mathbf{R}^{+};\mathcal{P}(\widehat{E})\right)\right]\label{eq:P(E)_Proc_Rep_Cadlag}
\end{equation}
by (\ref{eq:Construct_P(E)_Proc_Rep_Cadlag}) and Proposition \ref{prop:Sko_Basic_1}
(d) (with $S=(S_{0},\mathscr{O}_{\mathbf{R}^{\infty}}(S_{0}))$, $E=\mathcal{P}(\widehat{E})$
and $f=\widehat{\varphi}^{-1}$), (2) $Y$ is $\mathscr{F}_{t}^{X}$-adapted
by (\ref{eq:Check_M(E)_Rep_2}), the $\mathscr{F}_{t}^{X}$-adaptedness
of $\zeta^{\prime}$ and Fact \ref{fact:Proc_Path_Mapping} (a) (with
$E=(S_{0},\mathscr{O}_{\mathbf{R}^{\infty}}(S_{0}))$, $S=\mathcal{P}(\widehat{E})$,
$f=\widehat{\varphi}^{-1}$ and $X=\zeta^{\prime}$), and (3)
\begin{equation}
\inf_{t\in\mathbf{R}^{+}}\mathbb{P}\left(\varphi\circ X_{t}=\zeta_{t}^{\prime}=\widehat{\varphi}\circ Y_{t}\right)=1\label{eq:Check_P(E)_Proc_Rep_Cadlag_TF_Mod_2}
\end{equation}
by (\ref{eq:P(E)_Proc_Rep_Cadlag}) and (\ref{eq:Check_P(E)_Proc_Rep_Cadlag_TF_Mod_1}).

Now, the result follows by Proposition \ref{prop:M(E)_Proc_Rep} (a)
- (d), (\ref{eq:Check_P(E)_Proc_Rep_Cadlag_TF_Mod_2}) and (\ref{eq:Check_M(E)_Rep_2}).\end{proof}

\section{\label{sec:RepProc_Path_Space}Weak convergence about c$\grave{\mbox{a}}$dl$\grave{\mbox{a}}$g
replica}

\subsection{\label{sub:Proc_Reg}Regularity conditions about processes}

Before discussing weak convergence of c$\grave{\mbox{a}}$dl$\grave{\mbox{a}}$g
replicas, we give some regularity conditions about stochastic processes.
\begin{defn}
\label{def:Proc_Reg}Let $E$ be a topological space and $\{(\Omega^{i},\mathscr{F}^{i},\mathbb{P}^{i};X^{i})\}_{i\in\mathbf{I}}$
be $E$-valued processes.
\begin{itemize}
\item When $(E,\mathfrak{r})$ is a metric space, $\{X^{i}\}_{i\in\mathbf{I}}$
satisfies \textbf{Mild Pointwise Containment Condition for $\mathbf{T}$}\textit{
}($\mathbf{T}$-MPCC for short) if $\mathbf{T}\subset\mathbf{R}^{+}$
is non-empty and for any $\epsilon\in(0,\infty)$ and $t\in\mathbf{T}$,
there exists a \textit{totally bounded} (see p. \pageref{Total_Bounded})
set $A_{\epsilon,t}\in\mathscr{B}(E)$ satisfying%
\footnote{The notation ``$A^{\epsilon}$'' was defind in \S \ref{sub:Topo}.%
}
\begin{equation}
\inf_{i\in\mathbf{I}}\mathbb{P}^{i}\left(X_{t}^{i}\in A_{\epsilon,t}^{\epsilon}\right)\geq1-\epsilon.\label{eq:MPCC}
\end{equation}

\item $\{X^{i}\}_{i\in\mathbf{I}}$ satisfies\textit{ }\textbf{$\mathbf{T}$-Pointwise
$\mathbf{m}$-Tightness Condition} or \textbf{$\mathbf{T}$-Pointwise
Sequential} \textbf{$\mathbf{m}$-Tightness Condition} \textbf{in
$A\subset E$ }($\mathbf{T}$-PMTC or $\mathbf{T}$-PSMTC in $A$
for short) if $\mathbf{T}\subset\mathbf{R}^{+}$ is non-empty and
$\{X_{t}^{i}\}_{i\in\mathbf{I}}$ is $\mathbf{m}$-tight or sequentially
$\mathbf{m}$-tight in $A$ for all $t\in\mathbf{T}$, respectively.
Moreover, we say $\{X^{i}\}_{i\in\mathbf{I}}$ satisfies $\mathbf{T}$-PMTC
(resp. $\mathbf{T}$-PSMTC) if it satisfies $\mathbf{T}$-PMTC (resp.
$\mathbf{T}$-PSMTC) in $E$.
\item $\{X^{i}\}_{i\in\mathbf{I}}$ satisfies \textbf{Metrizable Compact
Containment Condition in $A$} (MCCC in $A$ for short) if for each
$\epsilon,T\in(0,\infty)$, there exists a $K_{\epsilon,T}\in\mathscr{K}^{\mathbf{m}}(E)$
such that $K_{\epsilon,T}\subset A$,
\begin{equation}
\bigcap_{t\in[0,T]}(X_{t}^{i})^{-1}(K_{\epsilon,T})\in\mathscr{F}^{i},\;\forall i\in\mathbf{I}\label{eq:CCC_Measurability}
\end{equation}
and
\begin{equation}
\inf_{i\in\mathbf{I}}\mathbb{P}^{i}\left(X_{t}^{i}\in K_{\epsilon,T},\forall t\in[0,T]\right)\ge1-\epsilon.\label{eq:CCC}
\end{equation}
Moreover, by ``$\{X^{i}\}_{i\in\mathbf{I}}$ satisfies MCCC'' we
mean it satisfies MCCC in $E$.
\item When $\{X^{i}\}_{i\in\mathbf{I}}$ are measurable processes, $\{X^{i}\}_{i\in\mathbf{I}}$
is said to satisfy \textbf{Long-time-average $\mathbf{m}$-Tightness
Condition in $A$ for $\{T_{k}\}_{k\in\mathbf{N}}$} ($T_{k}$-LMTC
in $A$ for short) if $T_{k}\uparrow\infty$%
\footnote{``$T_{k}\uparrow\infty$'' as usual denotes a non-decreasing sequence
$\{T_{k}\}_{k\in\mathbf{N}}\subset\mathbf{R}$ that converges to $\infty$.%
} and
\begin{equation}
\left\{ \frac{1}{T_{k}}\int_{0}^{T_{k}}\mathbb{P}^{i}\circ(X_{\tau}^{i})^{-1}d\tau\right\} _{k\in\mathbf{N},i\in\mathbf{I}}\subset\mathcal{P}(E)\label{eq:LT_Meas}
\end{equation}
is $\mathbf{m}$-tight in $A$. Moreover, by ``$\{X^{i}\}_{i\in\mathbf{I}}$
satisfies $T_{k}$-LMTC'' we mean it satisfies $T_{k}$-LMTC in $E$.
\item $\{X_{i}\}_{i\in\mathbf{I}}$ satisfies \textbf{Modulus of Continuity
Condition for $\mathfrak{r}$} ($\mathfrak{r}$-MCC for short) if:
(1) $\mathfrak{r}$ is a pseudometric on $E$, and (2) For any $\epsilon,T\in(0,\infty)$,
there exists a $\delta_{\epsilon,T}\in(0,\infty)$ such that%
\footnote{``$w_{\left|\cdot\right|,\delta,T}^{\prime}$'' is defined by (\ref{eq:w'})
with $E=\mathbf{R}$ and $\mathfrak{r}=\left|\cdot\right|$.%
}
\begin{equation}
\left\{ \omega\in\Omega:w_{\mathfrak{r},\delta_{\epsilon,T},T}^{\prime}\circ X^{i}(\omega)\geq\epsilon\right\} \in\mathscr{F}^{i},\;\forall i\in\mathbf{I}\label{eq:w'_Measurable}
\end{equation}
and
\begin{equation}
\sup_{i\in\mathbf{I}}\mathbb{P}^{i}\left(w_{\mathfrak{r},\delta_{\epsilon,T},T}^{\prime}\circ X^{i}\geq\epsilon\right)\leq\epsilon.\label{eq:MCC}
\end{equation}

\item $\{X_{i}\}_{i\in\mathbf{I}}$ satisfies \textbf{Modulus of Continuity
Condition} (MCC for short) if there exist a family of pseudometrics
$\mathcal{R}$ that induces $\mathscr{O}(E)$%
\footnote{The meaning of $\mathcal{R}$ inducing $\mathscr{O}(E)$ was explained
in \S \ref{sub:Topo}.%
} and $\{X^{i}\}_{i\in\mathbf{I}}$ satisfies $\mathfrak{r}$-MCC for
all $\mathfrak{r}\in\mathcal{R}$.
\item $\{X_{i}\}_{i\in\mathbf{I}}$ satisfies \textbf{Functional Modulus
of Continuity Condition for $\mathcal{D}$} ($\mathcal{D}$-FMCC for
short) if: (1) $\varpi(f)\circ X^{i}$ admits a c$\grave{\mbox{a}}$dl$\grave{\mbox{a}}$g
modification $\zeta^{f,i}$ for each $f\in\mathcal{D}\subset M(E;\mathbf{R})$
and $i\in\mathbf{I}$, and (2) $\{\zeta^{f,i}\}_{i\in\mathbf{I}}$
satisfies $\left|\cdot\right|$-MCC%
\footnote{$\left|\cdot\right|$-MCC means MCC for Euclidean metric $\left|\cdot\right|$.
The notation ``$\mathfrak{ae}(\mathcal{D})$'' was defined in \S
\ref{sub:Fun}.%
} for all $f\in\mathfrak{ae}(\mathcal{D})$.
\item $\{X^{i}\}_{i\in\mathbf{I}}$ satisfies \textbf{Weak Modulus of Continuity
Condition} (WMCC for short) if: (1) There exists a $\mathcal{D}\subset C(E;\mathbf{R})$
separating points on $E$, and (2) $\{X^{i}\}_{i\in\mathbf{I}}$ satisfies
$\mathcal{D}$-FMCC.
\end{itemize}
\end{defn}
\begin{note}
\label{note:Singleton_PMTC}An $E$-valued process $X$ is said to
satisfy any of the properties above except $\mathbf{T}$-PSMTC (in
$A$)%
\footnote{Note that sequential $\mathbf{m}$-tightness is for infinite collections
of measures or random variables.%
} if the singleton $\{X\}$ does. 
\end{note}

\begin{note}
\label{note:Prog_Reg_Ind_Trans}If $\{X^{i}\}_{i\in\mathbf{I}}$ and
$\{Y^{i}\}_{i\in\mathbf{I}}$ are two bijectively indistinguishable
families of $E$-valued processes (i.e. $X^{i}$ and $Y^{i}$ are
indistinguishable for all $i\in\mathbf{I}$), then each of the conditions
above is transitive between between $\{X^{i}\}_{i\in\mathbf{I}}$
and $\{Y^{i}\}_{i\in\mathbf{I}}$. Moreover, $\mathcal{D}$-FMCC and
WMCC are transitive between $\{X^{i}\}_{i\in\mathbf{I}}$ and $\{Y^{i}\}_{i\in\mathbf{I}}$
if $Y^{i}$ is a modification of $X^{i}$ for all $i\in\mathbf{I}$.\end{note}
\begin{rem}
\label{rem:Containment}$\,$
\begin{itemize}
\item Assuming total boundedness in lieu of compactness for each $A_{\epsilon,t}$,
$\mathbf{R}^{+}$-MPCC weakens the Pointwise Containment Property
in \cite[\S 3.7, Theorem 3.7.2]{EK86} and \cite[\S 5]{K15}.
\item MCCC is a variant of the standard \textit{Compact Containment Condition}
(see \cite[\S 4, (4.8)]{J86} and \cite[\S 3.7, (7.9)]{EK86}) using
$\mathbf{m}$-tightness, which becomes standard if $E$ has metrizable
compact sets. $\mathbf{R}^{+}$-PMTC is a similar variant of the Pointwise
Tight Condition in \cite[\S 5]{K15} and \cite[\S 3.7, (7.7)]{EK86}.
\item $T_{k}$-LMTC often appears in constructing stationary distributions
(see \cite{K71} and \cite{BBK00}). The measures in (\ref{eq:LT_Meas})
are well-defined by properties of measurable process and Fubini's
Theorem.
\end{itemize}
\end{rem}

\begin{rem}
\label{rem:MCC}$\,$
\begin{itemize}
\item MCC was used in \cite{J86} and \cite{K15} (in its finite-time-horizon
form) for general Tychonoff spaces. As long as $E$ is Hausdorff,
the assumption of pseudometrics $\mathcal{R}$ inducing $\mathscr{O}(E)$
in MCC implies $E$ is Tychonoff (see Proposition \ref{prop:CR} (a,
d)).
\item WMCC is a special case of $\mathcal{D}$-FMCC. Both of them are generically
milder than MCC as $\mathcal{D}$ need not strongly separate points
on $E$.
\end{itemize}
\end{rem}
\begin{note}
\label{note:FMCC_Weakly_Cadlag}If $\{X^{i}\}_{i\in\mathbf{I}}$ satisfy
$\mathcal{D}$-FMCC, then they are apparently $(\mathbf{R}^{+},\mathcal{D})$-c$\grave{\mbox{a}}$dl$\grave{\mbox{a}}$g
processes.\end{note}
\begin{rem}
\label{rem:CCC_MCCC_Measurability}In standard texts like \cite{B68}
and \cite{EK86}, $\mathfrak{r}$-MCC and MCCC are fundamental criteria
for establishing tightness or relative compactness in Skorokhod $\mathscr{J}_{1}$-spaces.
$(E,\mathfrak{r})$ is usually a separable metric space and our measurability
conditions (\ref{eq:CCC_Measurability}) and (\ref{eq:w'_Measurable})
are automatically true. We consider general spaces so we have to specify
(\ref{eq:CCC_Measurability}) and (\ref{eq:w'_Measurable}) as part
of MCCC and $\mathfrak{r}$-MCC. Given a c$\grave{\mbox{a}}$dl$\grave{\mbox{a}}$g
$X^{i}$, Lemma \ref{lem:CCC_Measurability} justifies (\ref{eq:CCC_Measurability})
under very mild conditions about $E$ and $K_{\epsilon,T}$, and Lemma
\ref{lem:MCC_Measurability} justifies (\ref{eq:w'_Measurable}) for
the following four cases:

\renewcommand{\labelenumi}{(\arabic{enumi}) }
\begin{enumerate}
\item $(E,\mathfrak{r})$ is a metric space and $X^{i}$ is a $D(\mathbf{R}^{+};E)$-valued
random variable.
\item $(E,\mathfrak{r})$ is a separable metric space.
\item $\mathfrak{r}=\rho_{\{f\}}$%
\footnote{The pseudometric $\rho_{\{f\}}$ is defined by (\ref{eq:TF_Metric})
with $\mathcal{D}=\{f\}$.%
} with $f\in C(E;\mathbf{R})$.
\item $\mathfrak{r}=\rho_{\mathcal{D}}$ with $\mathcal{D}\subset C(E;\mathbf{R})$
being a countable point-separating collection (hence $E$ is baseable).
\end{enumerate}
Consequently, $\mathcal{D}$-FMCC never incurs measurability issue
like (\ref{eq:w'_Measurable}) by case (2) above (with $E=\mathbf{R}$),
nor does $\rho_{\{f\}}$-MCC (resp. $\rho_{\mathcal{D}}$-MCC) for
c$\grave{\mbox{a}}$dl$\grave{\mbox{a}}$g processes and $\mathcal{D}$
consistent with case (3) (resp. case (4)) above.
\end{rem}

Besides the measurability conditions (\ref{eq:CCC_Measurability})
and (\ref{eq:w'_Measurable}), \S \ref{sec:Cadlag} of Appendix \ref{chap:App1}
also provides several results about the relationship among $\mathfrak{r}$-MCC,
MCC, $\mathcal{D}$-FMCC and WMCC. The above-mentioned containment
or tightness conditions will be further discussed in \S \ref{sec:Base_Proc}.

\subsection{\label{sub:RepProc_Tight}Tightness of c$\grave{\mbox{a}}$dl$\grave{\mbox{a}}$g
replicas}

Given $E$-valued processes $\{X^{i}\}_{i\in\mathbf{I}}$, we consider
tightness of their c$\grave{\mbox{a}}$dl$\grave{\mbox{a}}$g replicas
$\{\widehat{X}^{i}\in\mathfrak{rep}_{\mathrm{c}}(X^{i};E_{0},\mathcal{F})\}_{i\in\mathbf{I}}$
in $D(\mathbf{R}^{+};\widehat{E})$.
\begin{rem}
\label{rem:RepProc_D(Ehat)_RV}C$\grave{\mbox{a}}$dl$\grave{\mbox{a}}$g
replicas are always $D(\mathbf{R}^{+};\widehat{E})$-valued random
variables (see Fact \ref{fact:Cadlag_RepProc}) and their tightness
in $D(\mathbf{R}^{+};\widehat{E})$ has the usual meaning (compared
to our general interpretation in \S \ref{sec:RV}).\end{rem}
\begin{note}
\label{note:Base_CCC}Thanks to the compactness of $\widehat{E}$
(see Lemma \ref{lem:Base} (b)), the stringent MCCC is automatic for
$\widehat{E}$-valued processes.
\end{note}
Given MCCC, tightness of $\{\widehat{X}^{i}\}_{i\in\mathbf{I}}$ in
$D(\mathbf{R}^{+};\widehat{E})$ is reduced to $\mathcal{F}$-FMCC.
\begin{prop}
\label{prop:FR_Tight}Let $E$ be a topological space, $(E_{0},\mathcal{F};\widehat{E},\widehat{\mathcal{F}})$
be a base over $E$ and $\{(\Omega^{i},\mathscr{F}^{i},\mathbb{P}^{i};X^{i})\}_{i\in\mathbf{I}}$
be $E$-valued processes satisfying
\begin{equation}
\inf_{t\in\mathbf{T},i\in\mathbf{I}}\mathbb{P}^{i}\left(\bigotimes\mathcal{F}\circ X_{t}^{i}\in\bigotimes\widehat{\mathcal{F}}(\widehat{E})\right)=1\label{eq:Common_FR-Base_I}
\end{equation}
for some dense $\mathbf{T}\subset\mathbf{R}^{+}$. Then:

\renewcommand{\labelenumi}{(\alph{enumi})}
\begin{enumerate}
\item If $\{X^{i}\}_{i\in\mathbf{I}}$ satisfies $\mathcal{F}$-FMCC, then
$\{\widehat{X}^{i}=\mathfrak{rep}_{\mathrm{c}}(X^{i};E_{0},\mathcal{F})\}_{i\in\mathbf{I}}$
exists, satisfies $\widehat{\mathcal{F}}$-FMCC, satisfies $\rho_{\widehat{\mathcal{F}}}$-MCC
and is tight in $D(\mathbf{R}^{+};\widehat{E})$.
\item The converse of (a) is true when $\mathbf{T}=\mathbf{R}^{+}$ or $\{X^{i}\}_{i\in\mathbf{I}}$
are all c$\grave{\mbox{a}}$dl$\grave{\mbox{a}}$g processes.
\item If $\mathbf{I}$ is an infinite set and any subsequence of $\{X^{i}\}_{i\in\mathbf{I}}$
has a sub-subsequence satisfying $\mathcal{F}$-FMCC, then $\{\widehat{X}^{i}=\mathfrak{rep}_{\mathrm{c}}(X^{i};E_{0},\mathcal{F})\}_{i\in\mathbf{I}\backslash\mathbf{I}_{0}}$
exists and is sequentially tight in $D(\mathbf{R}^{+};\widehat{E})$
for some $\mathbf{I}_{0}\in\mathscr{P}_{0}(\mathbf{I})$.
\end{enumerate}
\end{prop}
\begin{proof}
(a) Suppose $\zeta^{f,i}$ is a c$\grave{\mbox{a}}$dl$\grave{\mbox{a}}$g
modification of $\varpi(f)\circ X^{i}$ for each $f\in\mathfrak{ae}(\mathcal{F})$
and $i\in\mathbf{I}$ (see Note \ref{note:FMCC_Weakly_Cadlag}) and
$\{\zeta^{f,i}\}_{i\in\mathbf{I}}$ satisfies $\left|\cdot\right|$-MCC
for all $f\in\mathfrak{ae}(\mathcal{F})$ (see Note \ref{note:Prog_Reg_Ind_Trans}).
It follows by (\ref{eq:Common_FR-Base_I}) and Proposition \ref{prop:FR}
(a) (with $X=X^{i}$) that $\{\widehat{X}^{i}=\mathfrak{rep}_{\mathrm{c}}(X^{i};E_{0},\mathcal{F})\}_{i\in\mathbf{I}}$
exists and satisfies
\begin{equation}
\inf_{f\in\mathfrak{ae}(\mathcal{F}),i\in\mathbf{I}}\mathbb{P}^{i}\left(\varpi(\widehat{f})\circ\widehat{X}^{i}=\zeta^{f,i}\right)=1,\label{eq:Check_RepProc_FMCC}
\end{equation}
thus proving $\{\widehat{X}^{i}\}_{i\in\mathbf{I}}$ satisfies $\widehat{\mathcal{F}}$-FMCC.
$\{\widehat{X}^{i}\}_{i\in\mathbf{I}}$ satisfies $\rho_{\widehat{\mathcal{F}}}$-MCC
by Lemma \ref{lem:Base} (a) and Proposition \ref{prop:MCC_3} (with
$E=\widehat{E}$ and $\mathcal{D}=\widehat{\mathcal{F}}$). Now, (a)
follows by Note \ref{note:Base_CCC}, Lemma \ref{lem:Base} (c) and
Theorem \ref{thm:Sko_RV_Tight_Polish} (with $(E,\mathfrak{r})=(\widehat{E},\rho_{\widehat{\mathcal{F}}})$
and $X^{i}=\widehat{X}^{i}$).

(b) Given tightness of $\{\widehat{X}^{i}\}_{i\in\mathbf{I}}$ in
$D(\mathbf{R}^{+};\widehat{E})$, $\{\varpi(\widehat{f})\circ\widehat{X}^{i}\}_{i\in\mathbf{I}}$
is tight in $D(\mathbf{R}^{+};\mathbf{R})$ for all $\widehat{f}\in\mathfrak{ae}(\widehat{\mathcal{F}})$
by Proposition \ref{prop:Sko_Basic_1} (d) (with $E=\widehat{E}$
and $S=\mathbf{R}$) and Fact \ref{fact:Push_Forward_Tight_2} (a)
(with $E=A=D(\mathbf{R}^{+};\widehat{E})$, $S=D(\mathbf{R}^{+};\mathbf{R})$,
$f=\varpi(\widehat{f})$ and $\Gamma=\{\mathbb{P}^{i}\circ(\widehat{X}^{i})^{-1}\}_{i\in\mathbf{I}}$).
Then, $\{\varpi(\widehat{f})\circ\widehat{X}^{i}\}_{i\in\mathbf{I}}$
satisfies $\left|\cdot\right|$-MCC for all $\widehat{f}\in\mathfrak{ae}(\widehat{\mathcal{F}})$
by Theorem \ref{thm:Sko_RV_Tight_Polish} (with $(E,\mathfrak{r})=(\mathbf{R},\left|\cdot\right|)$
and $X^{i}=\varpi(\widehat{f})\circ\widehat{X}^{i}$).

We have that
\begin{equation}
\begin{aligned} & \inf_{t\in\mathbf{T},i\in\mathbf{I},f\in\mathfrak{ae}(\mathcal{F})}\mathbb{P}^{i}\left(f\circ X_{t}=\widehat{f}\circ\widehat{X}_{t}^{i}\right)\\
 & \geq\inf_{t\in\mathbf{T},i\in\mathbf{I}}\mathbb{P}^{i}\left(\bigotimes\mathcal{F}\circ X_{t}\in\bigotimes\widehat{\mathcal{F}}(\widehat{E})\right)=1
\end{aligned}
\label{eq:Check_F-FMCC}
\end{equation}
by (\ref{eq:Common_FR-Base_I}) and (\ref{eq:Define_RepProc}) (with
$X=X^{i}$ and $\widehat{X}=\widehat{X}^{i}$). If $\mathbf{T}\neq\mathbf{R}^{+}$
and $\{X^{i}\}_{i\in\mathbf{I}}$ are all c$\grave{\mbox{a}}$dl$\grave{\mbox{a}}$g,
then $\varpi(\widehat{f})\circ\widehat{X}^{i}$ is indistinguishable
from $\varpi(f)\circ X^{i}$ for all $i\in\mathbf{I}$ and $f\in\mathfrak{ae}(\mathcal{F})$
by Proposition \ref{prop:FR} (c). If $\mathbf{T}=\mathbf{R}^{+}$,
(\ref{eq:Check_F-FMCC}) shows $\varpi(\widehat{f})\circ\widehat{X}^{i}$
is a c$\grave{\mbox{a}}$dl$\grave{\mbox{a}}$g modification of $\varpi(f)\circ X^{i}$.
Therefore, $\{X^{i}\}_{i\in\mathbf{I}}$ satisfies $\mathcal{F}$-FMCC
in either case.

(c) It follows by (a) and the given condition that any infinite subset
of $\{X^{i}\}_{i\in\mathbf{I}}$ admits a sub-subsequence $\{X^{i_{n}}\}_{n\in\mathbf{N}}$
such that $\{\widehat{X}^{i_{n}}=\mathfrak{rep}_{\mathrm{c}}(X^{i_{n}};E_{0},\mathcal{F})\}_{n\in\mathbf{N}}$
exists and is tight in $D(\mathbf{R}^{+};\widehat{E})$. Therefore,
$\widehat{X}^{i}=\mathfrak{rep}_{\mathrm{c}}(X^{i};E_{0},\mathcal{F})$
must exist for all $i\in\mathbf{I}$ with only finite exceptions and
(c) follows by the definition of sequential tightness.\end{proof}

Next, we consider tightness of the indistinguishable c$\grave{\mbox{a}}$dl$\grave{\mbox{a}}$g
replicas constructed by Proposition \ref{prop:PR} in $D(\mathbf{R}^{+};E_{0},\mathscr{O}_{E}(E_{0}))$
or $D(\mathbf{R}^{+};E)$.
\begin{prop}
\label{prop:PR_Tight}Let $E$ be a topological space, $(E_{0},\mathcal{F};\widehat{E},\widehat{\mathcal{F}})$
be a base over $E$ and $\{(\Omega^{i},\mathscr{F}^{i},\mathbb{P}^{i};X^{i})\}_{i\in\mathbf{I}}$
be $E$-valued c$\grave{\mbox{a}}$dl$\grave{\mbox{a}}$g processes.
Suppose that:

\renewcommand{\labelenumi}{(\roman{enumi})}
\begin{enumerate}
\item $(E_{0},\mathscr{O}_{E}(E_{0}))$ is a Tychonoff space.
\item $\{X^{i}\}_{i\in\mathbf{I}}$ satisfies MCCC in $E_{0}$.
\end{enumerate}
Then, there exists an $S_{0}\subset\mathbb{D}_{0}\circeq D(\mathbf{R}^{+};E_{0},\mathscr{O}_{E}(E_{0}))$
such that:

\renewcommand{\labelenumi}{(\alph{enumi})}
\begin{enumerate}
\item $S_{0}$ and $\mathbb{D}_{0}$ satisfy%
\footnote{$\mathscr{B}(\mathbb{D}_{0})$ is generated by the Skorokhod $\mathscr{J}_{1}$-topology
$\mathscr{J}(E_{0},\mathscr{O}_{E}(E_{0}))$.%
}
\begin{equation}
\left.\mathscr{B}(\mathbb{D}_{0})\right|_{S_{0}}=\left.\mathscr{B}(E)^{\otimes\mathbf{R}^{+}}\right|_{S_{0}}\subset\left.\mathscr{B}(E)^{\otimes\mathbf{R}^{+}}\right|_{\mathbb{D}_{0}}\subset\mathscr{B}(\mathbb{D}_{0}).\label{eq:PR_S_Sko_Borel_Prod_Equal}
\end{equation}

\item $\{\widehat{X}^{i}=\mathfrak{rep}_{\mathrm{c}}(X^{i};E_{0},\mathcal{F})\}_{i\in\mathbf{I}}$
satisfies
\begin{equation}
\widehat{X}^{i}\in M\left[\Omega^{i},\mathscr{F}^{i};S_{0},\mathscr{O}_{\mathbb{D}_{0}}(S_{0})\right],\;\forall i\in\mathbf{I}\label{eq:RepProc_D(E0)_RV_I}
\end{equation}
and
\begin{equation}
\inf_{i\in\mathbf{I}}\mathbb{P}^{i}\left(X^{i}=\widehat{X}^{i}\in S_{0}\right)=1.\label{eq:PR_All_Indistinguishable_I}
\end{equation}

\item $\{\widehat{X}^{i}\}_{i\in\mathbf{I}}$ is $\mathbf{m}$-tight in
$S_{0}$ as $\mathbb{D}_{0}$-valued random variables if and only
if $\{X^{i}\}_{i\in\mathbf{I}}$ satisfies $\mathcal{F}$-FMCC.
\end{enumerate}
\end{prop}
\begin{rem}
\label{rem:Base_Sko_Borel}The c$\grave{\mbox{a}}$dl$\grave{\mbox{a}}$g
replicas $\{\widehat{X}^{i}\}_{i\in\mathbf{I}}$ above are $\mathbb{D}_{0}$-valued
random variables and their tightness in $S_{0}\subset\mathbb{D}_{0}$
has the usual meaning. We noted in \S \ref{sub:Sko_Meas} that $\mathbb{D}_{0}$
and $S_{0}$ always satisfy%
\footnote{$\mathscr{B}(\mathbb{D}_{0})$ is generated by the Skorokhod $\mathscr{J}_{1}$-topology
of $D(\mathbf{R}^{+};E_{0},\mathscr{O}_{E}(E_{0}))$.%
}
\begin{equation}
\left.\mathscr{B}(\mathbb{D}_{0})\right|_{S_{0}}\supset\left.\mathscr{B}(E)^{\otimes\mathbf{R}^{+}}\right|_{S_{0}}\label{eq:PR_S0_Sko_Borel_>_Prod}
\end{equation}
but not necessarily the equality in (\ref{eq:PR_S_Sko_Borel_Prod_Equal}).
So, $\mathbb{D}_{0}$-valued random variable (like $\widehat{X}^{i}$
in (\ref{eq:RepProc_D(E0)_RV_I})) is generally a stronger concept
than $(E_{0},\mathscr{O}_{E}(E))$-valued c$\grave{\mbox{a}}$dl$\grave{\mbox{a}}$g
process.
\end{rem}

In fact, the developments of Proposition \ref{prop:PR_Tight} do not
require a base. We establish the following more general result without
imposing the boundedness of the point-separating functions.
\begin{thm}
\label{thm:PR_Tight}Let $E$ be a topological space, $\{(\Omega^{i},\mathscr{F}^{i},\mathbb{P}^{i};X^{i})\}_{i\in\mathbf{I}}$
be $E$-valued c$\grave{\mbox{a}}$dl$\grave{\mbox{a}}$g processes,
$E_{0}\in\mathscr{B}(E)$ and $\mathcal{D}\subset C(E;\mathbf{R})$.
Suppose that:

\renewcommand{\labelenumi}{(\roman{enumi})}
\begin{enumerate}
\item $\mathcal{D}$ separates points on $E_{0}$.
\item $(E_{0},\mathscr{O}_{E}(E_{0}))$ is a Tychonoff space.
\item $\{X^{i}\}_{i\in\mathbf{I}}$ satisfies MCCC in $E_{0}$.
\end{enumerate}
Then, there exist $S_{0}\subset\mathbb{D}_{0}\circeq D(\mathbf{R}^{+};E_{0},\mathscr{O}_{E}(E_{0}))$
and $\{\widehat{X}^{i}\in S_{0}^{\Omega^{i}}\}_{i\in\mathbf{I}}$
such that:

\renewcommand{\labelenumi}{(\alph{enumi})}
\begin{enumerate}
\item (\ref{eq:PR_S_Sko_Borel_Prod_Equal}) holds.
\item $\{\widehat{X}^{i}\}_{i\in\mathbf{I}}$ satisfies (\ref{eq:RepProc_D(E0)_RV_I})
and (\ref{eq:PR_All_Indistinguishable_I}).
\item $\{\widehat{X}^{i}\}_{i\in\mathbf{I}}$ is $\mathbf{m}$-tight in
$S_{0}$ as $\mathbb{D}_{0}$-valued random variables if and only
if $\{X^{i}\}_{i\in\mathbf{I}}$ satisfies $\mathcal{D}$-FMCC.
\end{enumerate}
\end{thm}
\begin{proof}
We divide the proof into five steps. We equip $E_{0}$ with the subspace
topology $\mathscr{O}_{E}(E_{0})$ throughout the proof, which we
make implicit for convenience.

\textit{Step 1: Construct $S_{0}$}. By the condition (iii) above,
\begin{equation}
\inf_{i\in\mathbf{I}}\mathbb{P}^{i}\left(X_{t}^{i}\in A_{p,q},\forall t\in[0,q]\right)\geq1-2^{-p-q},\;\forall p,q\in\mathbf{N}.\label{eq:Pick_CCC_Apq}
\end{equation}
holds for some $\{A_{p,q}\}_{p,q\in\mathbf{N}}\subset\mathscr{K}^{\mathbf{m}}(E_{0})$.
It follows that
\begin{equation}
K_{p,q}\circeq\bigcup_{i=1}^{q}A_{p,i}\in\mathscr{K}^{\mathbf{m}}(E_{0})\subset\mathscr{C}(E_{0}),\;\forall p,q\in\mathbf{N}\label{eq:Pick_CCC_Kpq_E0}
\end{equation}
by the Hausdorff property of $E_{0}$, Proposition \ref{prop:Separability}
(c), Lemma \ref{lem:MC_Union} and Proposition \ref{prop:Compact}
(a). From (\ref{eq:Pick_CCC_Apq}) and (\ref{eq:Pick_CCC_Kpq_E0})
we obtain that
\begin{equation}
K_{p,q}\subset K_{p,q+1},\;\forall p,q\in\mathbf{N}\label{eq:Pick_CCC_Kpq_1}
\end{equation}
and
\begin{equation}
\inf_{i\in\mathbf{I}}\mathbb{P}^{i}\left(X_{t}^{i}\in K_{p,q},\forall t\in[0,q]\right)\geq1-2^{-p-q},\;\forall p,q\in\mathbf{N}.\label{eq:Pick_CCC_Kpq_2}
\end{equation}
Letting%
\footnote{$K_{p,q}^{[0,q)}$ in (\ref{eq:Vp_CCC}) means the Cartesian power
of $K_{p,q}$ for the index set $[0,q)$.%
}
\begin{equation}
V_{p}\circeq\bigcap_{q\in\mathbf{N}}\left\{ x\in\mathbb{D}_{0}:x|_{[0,q)}\in K_{p,q}^{[0,q)}\right\} ,\;\forall p\in\mathbf{N},\label{eq:Vp_CCC}
\end{equation}
one finds by the fact $E_{0}\in\mathscr{B}(E)$, Lemma \ref{lem:Cadlag_Measurability}
(b) (with $E=E_{0}$, $A=K_{p,q}$ and $T=q$) and and Lemma \ref{lem:Sko_Proj}
(b) (with $E=E_{0}$) that
\begin{equation}
V_{p}\in\left.\mathscr{B}_{E}(E_{0})^{\otimes\mathbf{R}^{+}}\right|_{\mathbb{D}_{0}}\subset\left.\mathscr{B}(E)^{\otimes\mathbf{R}^{+}}\right|_{\mathbb{D}_{0}}\subset\mathscr{B}(\mathbb{D}_{0}),\;\forall p\in\mathbf{N},\label{eq:Check_PR_Vp_Prod_Measurable}
\end{equation}
which immediately implies $S_{0}\in\mathscr{B}(\mathbb{D}_{0})$.
Consequently,
\begin{equation}
S_{0}\circeq\bigcup_{p\in\mathbf{N}}V_{p}\in\left.\mathscr{B}(E)^{\otimes\mathbf{R}^{+}}\right|_{S_{0}}\subset\left.\mathscr{B}(E)^{\otimes\mathbf{R}^{+}}\right|_{\mathbb{D}_{0}}\subset\mathscr{B}(\mathbb{D}_{0}).\label{eq:Check_PR_S0_Prod_Measurable}
\end{equation}

\textit{Step 2: Verify (a)}. Each of $\{K_{p,q}\}_{p,q\in\mathbf{N}}$
is a $\mathcal{D}|_{E_{0}}$-baseable subset of $E_{0}$ by (\ref{eq:Pick_CCC_Kpq_E0})
and Proposition \ref{prop:MC} (a, f) (with $E=E_{0}$, $K=K_{p,q}$
and $\mathcal{D}=\mathcal{D}|_{E_{0}}$). So, $\mathcal{D}$ has a
countable subset that separates and strongly separates points on each
of $\{K_{p,q}\}_{p,q\in\mathbf{N}}$ by Lemma \ref{lem:SP_on_Compact}.
For simplicity, we assume $\mathcal{D}$ is countable in Step 2 -
4 of the proof.

Letting $\Psi\circeq\varpi[\mathfrak{ae}(\mathcal{D})]$, we have
\begin{equation}
\Psi|_{V_{p}}\in\mathbf{imb}\left(V_{p},\mathscr{O}_{\mathbb{D}_{0}}(V_{p});D(\mathbf{R}^{+};\mathbf{R})^{\mathfrak{ae}(\mathcal{D})}\right),\;\forall p\in\mathbf{N}\label{eq:Psi|Vp_Homeo}
\end{equation}
and
\begin{equation}
\left.\mathscr{B}(\mathbb{D}_{0})\right|_{V_{p}}=\left.\mathscr{B}(E)^{\otimes\mathbf{R}^{+}}\right|_{V_{p}},\;\forall p\in\mathbf{N}\label{eq:Vp_Sko_Borel_Prod_Equal_D(E0)}
\end{equation}
by Lemma \ref{lem:Sko_Prod_Eq} (with $E=E_{0}$, $V=V_{p}$, $p=q$
and $A_{p}=K_{p,q}$). One then finds by (\ref{eq:Check_PR_Vp_Prod_Measurable}),
(\ref{eq:Vp_Sko_Borel_Prod_Equal_D(E0)}) and Fact \ref{fact:Union_Borel}
(with $E=S_{0}$, $n=p$, $A_{n}=V_{p}$, $\mathscr{U}_{1}=\mathscr{B}(\mathbb{D}_{0})|_{S_{0}}$
and $\mathscr{U}_{2}=\mathscr{B}(E)^{\otimes\mathbf{R}^{+}}|_{S_{0}}$)
that
\begin{equation}
\left.\mathscr{B}(\mathbb{D}_{0})\right|_{S_{0}}\subset\left.\mathscr{B}(E)^{\otimes\mathbf{R}^{+}}\right|_{S_{0}}.\label{eq:PR_S0_Sko_Borel_<_Prod}
\end{equation}
Now, (a) follows by (\ref{eq:Check_PR_S0_Prod_Measurable}), (\ref{eq:PR_S0_Sko_Borel_>_Prod})
and (\ref{eq:PR_S0_Sko_Borel_<_Prod}).

\textit{Step 3: Construct $\{\widehat{X}^{i}\}_{i\in\mathbf{I}}$
and verify (b)}. It follows by (\ref{eq:Pick_CCC_Kpq_1}) and (\ref{eq:Pick_CCC_Kpq_2})
that
\begin{equation}
\begin{aligned}\inf_{i\in\mathbf{I}}\mathbb{P}^{i}\left(X^{i}\in V_{p}\right) & \geq1-\sup_{i\in\mathbf{I}}\sum_{q\in\mathbf{N}}\left[1-\mathbb{P}^{i}\left(X_{t}^{i}\in K_{p,q},\forall t\in[0,q]\right)\right]\\
 & \geq1-2^{-p},\;\forall p\in\mathbf{N}.
\end{aligned}
\label{eq:Check_Vp_CCC_Large}
\end{equation}
Then, (\ref{eq:Check_PR_S0_Prod_Measurable}) and (\ref{eq:Check_Vp_CCC_Large})
imply
\begin{equation}
\inf_{i\in\mathbf{I}}\mathbb{P}^{i}\left(X^{i}\in S_{0}\subset E_{0}^{\mathbf{R}^{+}}\right)=1.\label{eq:Common_Path-Base_I}
\end{equation}
By Proposition \ref{prop:PR} (a) (with $X=X^{i}$), there exist
\begin{equation}
\widehat{X}^{i}\in M\left(\Omega^{i},\mathscr{F}^{i};S_{0},\left.\mathscr{B}(E)^{\otimes\mathbf{R}^{+}}\right|_{S_{0}}\right),\;\forall i\in\mathbf{I}\label{eq:PR_S0_Valued}
\end{equation}
satisfying (\ref{eq:PR_All_Indistinguishable_I}). Now, (b) follows
by (\ref{eq:PR_S0_Valued}) and (a).

\textit{Step 4: Verify sufficiency of (c)}. We have that%
\footnote{Here, the replica process $\widehat{X}^{i}$ is an $E_{0}$-valued
process and so $f\circ\widehat{X}_{t}^{i}$ is well-defined.%
}
\begin{equation}
\inf_{f\in\mathfrak{ae}(\mathcal{D}),i\in\mathbf{I}}\mathbb{P}^{i}\left(\varpi(f)\circ X^{i}=\varpi(f)\circ\widehat{X}^{i}\right)=1\label{eq:PR_TF_Indistinguishable}
\end{equation}
and
\begin{equation}
\varpi(f)\circ\widehat{X}^{i}\in M\left(\Omega^{i},\mathscr{F}^{i};D(\mathbf{R}^{+};\mathbf{R})\right),\;\forall f\in\mathfrak{ae}(\mathcal{D}),i\in\mathbf{I}\label{eq:PR_TF_Path_Space_RV}
\end{equation}
by (\ref{eq:PR_All_Indistinguishable_I}) and Proposition \ref{prop:Sko_Basic_1}
(d) (with $S=E_{0}$ and $E=\mathbf{R}$). Fixing $f\in\mathfrak{ae}(\mathcal{D})$,
$\varpi(f)\circ\widehat{X}^{i}$ is the unique c$\grave{\mbox{a}}$dl$\grave{\mbox{a}}$g
modification of $\varpi(f)\circ X^{i}$ up to indistinguishability
for all $i\in\mathbf{I}$ by (\ref{eq:PR_TF_Indistinguishable}) and
Proposition \ref{prop:Proc_Basic_2} (f, h). $\{\varpi(f)\circ\widehat{X}^{i}\}_{i\in\mathbf{I}}$
satisfies $\left|\cdot\right|$-MCC by (\ref{eq:PR_TF_Indistinguishable})
and $\{X^{i}\}_{i\in\mathbf{I}}$ satisfying $\mathcal{D}$-FMCC.
$\{\varpi(f)\circ\widehat{X}^{i}\}_{i\in\mathbf{I}}$ satisfies MCCC
by the boundedness of $f$. Hence, $\{\varpi(f)\circ\widehat{X}^{i}\}_{i\in\mathbf{I}}$
is tight in $D(\mathbf{R}^{+};\mathbf{R})$ by Theorem \ref{thm:Sko_RV_Tight_Polish}
(with $(E,\mathfrak{r})=(\mathbf{R},\left|\cdot\right|)$ and $X^{i}=\varpi(f)\circ\widehat{X}^{i}$).

Letting $\Psi=\varpi[\mathfrak{ae}(\mathcal{D})]$ again, $\{\Psi\circ\widehat{X}^{i}\}_{i\in\mathbf{I}}$
is tight in $D(\mathbf{R}^{+};\mathbf{R})^{\mathfrak{ae}(\mathcal{D})}$
by tightness of $\{\varpi(f)\circ\widehat{X}^{i}\}_{i\in\mathbf{I}}$
in $D(\mathbf{R}^{+};\mathbf{R})$ and Proposition \ref{lem:Sko_Prod_Tight}
(a) (with $E=E_{0}$, $\mathcal{D}=\mathfrak{ae}(\mathcal{D})$ and
$\mu^{i}=\mathbb{P}^{i}\circ(\widehat{X}^{i})^{-1}\in\mathcal{P}(\mathbb{D}_{0})$).
$\mathfrak{ae}(\mathcal{D})$ is countable by Fact \ref{fact:ac_mc_Countable}
and the fact $\mathcal{D}$ is countable in step 2 - 4, so $\mathbf{R}^{\mathfrak{ae}(\mathcal{D})}$
is a Polish space by Proposition \ref{prop:Var_Polish} (f).
\begin{equation}
\varphi\circeq\bigotimes\mathfrak{ae}(\mathcal{D})\in C\left(E;\mathbf{R}^{\mathfrak{ae}(\mathcal{D})}\right)\label{eq:Prod(ae(D))_Cont}
\end{equation}
by Fact \ref{fact:Prod_Map_2} (b). Letting $\{K_{p,q}\}_{p,q\in\mathbf{N}}$
be as in (\ref{eq:Pick_CCC_Kpq_E0}), we have that
\begin{equation}
\bigotimes\mathfrak{ae}(\mathcal{D})(K_{p,q})\in\mathscr{K}\left(\mathbf{R}^{\mathfrak{ae}(\mathcal{D})}\right)\subset\mathscr{C}\left(\mathbf{R}^{\mathfrak{ae}(\mathcal{D})}\right)\label{eq:Check_Vp_Image_Closed}
\end{equation}
by (\ref{eq:Prod(ae(D))_Cont}) and Proposition \ref{prop:Compact}
(a, e). It follows by Lemma \ref{lem:CCC_Closed} (with $E=E_{0}$,
$V=V_{p}$, $p=q$ and $A_{q}=K_{p,q}$) that
\begin{equation}
\Psi(V_{p})\in\mathscr{C}\left(D(\mathbf{R}^{+};\mathbf{R})^{\mathfrak{ae}(\mathcal{D})}\right),\;\forall p\in\mathbf{N}.\label{eq:Psi(Vp)_Closed}
\end{equation}
$\Psi\in M(S_{0},\mathscr{O}_{\mathbb{D}_{0}}(V_{p});D(\mathbf{R}^{+};\mathbf{R})^{\mathfrak{ae}(\mathcal{D})})$
by (\ref{eq:Psi|Vp_Homeo}) and the fact $S_{0}=\bigcup_{p\in\mathbf{N}}V_{p}$.
Hence, tightness of $\{\widehat{X}^{i}\}_{i\in\mathbf{I}}$ in $S_{0}$
follows by $S_{0}=\bigcup_{p\in\mathbf{N}}V_{p}$, (\ref{eq:Psi|Vp_Homeo}),
(\ref{eq:Psi(Vp)_Closed}) and Lemma \ref{lem:Pull_Back_Tight} (with
$E=\mathbb{D}_{0}$, $S=D(\mathbf{R}^{+};\mathbf{R})^{\mathfrak{ae}(\mathcal{D})}$,
$A_{p}=V_{p}$, $E_{0}=S_{0}$ and $f=\Psi$).

$A\circeq\bigcup_{p,q\in\mathbf{N}}K_{p,q}\in\mathscr{K}_{\sigma}^{\mathbf{m}}(E_{0})$
is a $\mathcal{D}|_{E_{0}}$-baseable subset of $E_{0}$ by Proposition
\ref{prop:Sigma_MC} (b, e) (with $E=E_{0}$, and $\mathcal{D}=\mathcal{D}|_{E_{0}}$).
$S_{0}$ is a baseable subspace of $\mathbb{D}_{0}$ by Proposition
\ref{prop:Sko_Baseable} (b) (with $E=E_{0}$) and Fact \ref{fact:Baseable_Space_Subset}.
Thus, $\{\widehat{X}^{i}\}_{i\in\mathbf{I}}$ is $\mathbf{m}$-tight
in $S_{0}$ by Corollary \ref{cor:Baseable_MC} (a) (with $E=S_{0}$).

\textit{Step 5: Verify necessity of (c)}. We no longer take $\mathcal{D}$
to be countable. Suppose $\{\widehat{X}^{i}\}_{i\in\mathbf{I}}$ is
$\mathbf{m}$-tight in $S_{0}$. For each fixed $f\in\mathfrak{ae}(\mathcal{D})$,
$\{\varpi(f)\circ\widehat{X}^{i}\}_{i\in\mathbf{I}}$ is tight in
$D(\mathbf{R}^{+};\mathbf{R})$ by Proposition \ref{prop:Sko_Basic_1}
(d) (with $E=E_{0}$ and $S=\mathbf{R}$) and Fact \ref{fact:Push_Forward_Tight_2}
(with $E=\mathbb{D}_{0}$, $S=D(\mathbf{R}^{+};\mathbf{R})$, $f=\varpi(f)$,
$\mu^{i}=\mathbb{P}^{i}\circ(\widehat{X}^{i})^{-1}\in\mathcal{P}(\mathbb{D}_{0})$
and $\Gamma=\{\mu^{i}\}_{i\in\mathbf{I}}$). $\{\varpi(f)\circ\widehat{X}^{i}\}_{i\in\mathbf{I}}$
satisfies $\left|\cdot\right|$-MCC by Theorem \ref{thm:Sko_RV_Tight_Polish}
(with $(E,\mathfrak{r})=(\mathbf{R},\left|\cdot\right|)$ and $X^{i}=\varpi(f)\circ\widehat{X}^{i}$).
Hence, $\{X^{i}\}_{i\in\mathbf{I}}$ satisfies $\mathcal{D}$-FMCC
by (\ref{eq:PR_TF_Indistinguishable}) and Proposition \ref{prop:Proc_Basic_2}
(f).\end{proof}

\begin{proof}
[Proof of Proposition \ref{prop:PR_Tight}]This result follows immediately
by Theorem \ref{thm:PR_Tight} (with $\mathcal{D}=\mathcal{F}$) and
Proposition \ref{prop:PR} (with $X=X^{i}$).\end{proof}

\subsection{\label{sub:RepProc_WC}Weak convergence of c$\grave{\mbox{a}}$dl$\grave{\mbox{a}}$g
replicas}

The following proposition connects the weak convergence of c$\grave{\mbox{a}}$dl$\grave{\mbox{a}}$g
replicas on path space and consequence of their finite-dimensional
convergence.
\begin{prop}
\label{prop:FR_FC_WC}Let $E$ be a topological space, $\{(\Omega^{n},\mathscr{F}^{n},\mathbb{P}^{n};X^{n})\}_{n\in\mathbf{N}}$
be $E$-valued processes, $(E_{0},\mathcal{F};\widehat{E},\widehat{\mathcal{F}})$
be a base over $E$, $\widehat{X}^{n}\in\mathfrak{rep}_{\mathrm{c}}(X^{n};E_{0},\mathcal{F})$
exists%
\footnote{Several conditions for the existence of c$\grave{\mbox{a}}$dl$\grave{\mbox{a}}$g
replicas were given in Proposition \ref{prop:FR}.%
} for each $n\in\mathbf{N}$ and $(\Omega,\mathscr{F},\mathbb{P};Y)$
be a $D(\mathbf{R}^{+};\widehat{E})$-valued random variable. Then:

\renewcommand{\labelenumi}{(\alph{enumi})}
\begin{enumerate}
\item If $\{\widehat{X}^{n}\}_{n\in\mathbf{N}}$ and $Y$ satisfy%
\footnote{The meaning of (\ref{eq:RepProc_WC_Y}) is explained in \S \ref{sec:RV}.
The notation ``$J(Y)$'' is defined in (\ref{eq:J(X)}).%
}
\begin{equation}
\widehat{X}^{n}\Longrightarrow Y\mbox{ as }n\uparrow\infty\mbox{ on }D(\mathbf{R}^{+};\widehat{E}),\label{eq:RepProc_WC_Y}
\end{equation}
then
\begin{equation}
\widehat{X}^{n}\xrightarrow{\quad\mathrm{D}(\mathbf{R}^{+}\backslash J(Y))\quad}Y\mbox{ as }n\uparrow\infty.\label{eq:RepProc_FC_along_R-J(Y)}
\end{equation}

\item If (\ref{eq:RepProc_FC_along_T_Y}) holds for some dense $\mathbf{T}\subset\mathbf{R}^{+}$,
and if $\{X^{n}\}_{n\in\mathbf{N}}$ satisfies 
\begin{equation}
\inf_{t\in\mathbf{T},n\in\mathbf{N}}\mathbb{P}^{n}\left(\bigotimes\mathcal{F}\circ X_{t}^{n}\in\bigotimes\widehat{\mathcal{F}}(\widehat{E})\right)=1\label{eq:Common_FR-Base_N}
\end{equation}
and $\mathcal{F}$-FMCC, then (\ref{eq:RepProc_WC_Y}) holds.
\item If (\ref{eq:RepProc_WC_Y}) holds, $\{X^{n}\}_{n\in\mathbf{N}}$ is
$(\mathbf{T},\mathcal{F}\backslash\{1\})$-AS%
\footnote{The notion of $(\mathbf{T},\mathcal{F}\backslash\{1\})$-AS was introduced
in Definition \ref{def:FC}.%
} and (\ref{eq:Common_FR-Base_N}) holds for some conull%
\footnote{Conull set was specified in \S \ref{sub:Prod_Space}. Conull subset
of $\mathbf{R}$ is in the Lebesgue sense.%
} $\mathbf{T}\subset\mathbf{R}^{+}$, then $Y$ is an $\widehat{E}$-valued
stationary process.
\end{enumerate}
\end{prop}
\begin{note}
\label{note:J(Y)}If $Y$ is an $\widehat{E}$-valued c$\grave{\mbox{a}}$dl$\grave{\mbox{a}}$g
process (especially a c$\grave{\mbox{a}}$dl$\grave{\mbox{a}}$g replica),
then $J(Y)\subset(0,\infty)$ is countable by Lemma \ref{lem:Base}
(c) and Proposition \ref{prop:J(Mu)_J(X)_Baseable}. In other words,
$\mathbf{R}^{+}\backslash J(Y)$ is a cocountable%
\footnote{The notion of cocountable set was defined in \S \ref{sub:Set_Num}.%
} subset of $\mathbf{R}^{+}$.
\end{note}
\begin{proof}
[Proof of Proposition \ref{prop:FR_FC_WC}](a) follows by Lemma \ref{lem:Base}
(c) and Theorem \ref{thm:Sko_RV_WC_FC_Metrizable_Separable} (a) (with
$E=\widehat{E}$, $X^{n}=\widehat{X}^{n}$ and $X=Y$).

(b) (\ref{eq:Common_FR-Base_N}) is a version of (\ref{eq:Common_FR-Base_I})
with $\mathbf{I}=\mathbf{N}$. Given dense $\mathbf{T}$ and $\mathcal{F}$-FMCC,
$\{\widehat{X}^{n}\}_{n\in\mathbf{N}}$ as the unique c$\grave{\mbox{a}}$dl$\grave{\mbox{a}}$g
replicas of $\{X^{n}\}_{n\in\mathbf{N}}$ is tight in the Polish space
$D(\mathbf{R}^{+};\widehat{E})$ by Proposition \ref{prop:FR_Tight}
(a) (with $\mathbf{I}=\mathbf{N}$). It is relatively compact in $D(\mathbf{R}^{+};\widehat{E})$%
\footnote{Relative compactness of $D(\mathbf{R}^{+};\widehat{E})$-valued random
variables $\{\widehat{X}^{n}\}_{n\in\mathbf{N}}$ follows our interpretation
in \S \ref{sec:RV}.%
} by the Prokhorov's Theorem (Theorem \ref{thm:Prokhorov} (b)). Now,
(b) follows by Theorem \ref{thm:Sko_RV_WC_FC_Metrizable_Separable}
(b) (with $E=\widehat{E}$, $X^{n}=\widehat{X}^{n}$ and $X=Y$).

(c) $\mathbf{T}\backslash J(Y)$ is a conull set, so
\begin{equation}
\mathbf{S}_{Y,\mathbf{T}_{0}}\circeq\bigcap_{t\in\mathbf{T}_{0}}\left\{ c\in(0,\infty):t+c\in\mathbf{T}\backslash J(Y)\right\} \label{eq:S_Y_T0_J(Y)}
\end{equation}
is a conull hence dense subset of $\mathbf{R}^{+}$ for any $\mathbf{T}_{0}\in\mathscr{P}_{0}(\mathbf{T}\backslash J(Y))$.
Fixing $\mathbf{T}_{0}\in\mathscr{P}_{0}(\mathbf{T}\backslash J(Y))$
and $f\in\mathfrak{mc}[\Pi^{\mathbf{T}_{0}}(\mathcal{F}\backslash\{1\})]$,
we have
\begin{equation}
\mathbb{E}\left[\widehat{f}\circ Y_{\mathbf{T}_{0}}\right]=\mathbb{E}\left[\widehat{f}\circ Y_{\mathbf{T}_{0}+c}\right]\label{eq:RepProc_Lim_Sta_Int_Test}
\end{equation}
for all $c\in\mathbf{S}_{Y,\mathbf{T}_{0}}$ by (a) and Lemma \ref{lem:RepProc_Int_Test}
(b, e) (with $\mathbf{T}=\mathbf{T}\backslash J(Y)$ and $\mathbf{S}_{\mathbf{T}_{0}}=\mathbf{S}_{Y,\mathbf{T}_{0}}$).
$\{Y_{t+c}\}_{c\geq0}$ is a c$\grave{\mbox{a}}$dl$\grave{\mbox{a}}$g
process for all $t\in\mathbf{T}_{0}$ since $Y$ is c$\grave{\mbox{a}}$dl$\grave{\mbox{a}}$g.
$\zeta\circeq\{\varpi(\widehat{f})\circ Y_{\mathbf{T}_{0}+c}\}_{c\geq0}$
is also a c$\grave{\mbox{a}}$dl$\grave{\mbox{a}}$g process by Fact
\ref{fact:Cadlag_Proc} (a, b) (with $\mathbf{I}=\mathbf{T}_{0}$,
$i=t$, $X^{i}=Y_{t+c}$, $X=\{Y_{\mathbf{T}_{0}+c}\}_{t\geq0}$ and
$f=\widehat{f}$). Then, (\ref{eq:RepProc_Lim_Sta_Int_Test}) extends
to all $c\in(0,\infty)$ by the denseness of $\mathbf{S}_{\mathbf{T}_{0},Y}$
in $\mathbf{R}^{+}$, the c$\grave{\mbox{a}}$dl$\grave{\mbox{a}}$g
property of $\zeta$ and the Dominated Convergence Theorem. Now, (c)
follows by Corollary \ref{cor:Base_Sep_Meas} (a) (with $d=\aleph(\mathbf{T}_{0})$
and $A=\widehat{E}^{d}$) and Fact \ref{fact:Sko_RV_Cadlag} (c) (with
$E=\widehat{E}$).\end{proof}

The next proposition connects weak convergence of c$\grave{\mbox{a}}$dl$\grave{\mbox{a}}$g
replicas on $D(\mathbf{R}^{+};\widehat{E})$ and that on the restricted
path space $D(\mathbf{R}^{+};E_{0},\mathscr{O}_{E}(E_{0}))$ (if well-defined).
\begin{prop}
\label{prop:PR_WC}Let $E$ be a topological space, $(E_{0},\mathcal{F};\widehat{E},\widehat{\mathcal{F}})$
be a base over $E$, $(\Omega,\mathcal{F},\mathbb{P};X)$ and $\{(\Omega^{n},\mathscr{F}^{n},\mathbb{P}^{n};X^{n})\}_{n\in\mathbf{N}}$
be $E$-valued c$\grave{\mbox{a}}$dl$\grave{\mbox{a}}$g processes,
$\widehat{X}\in\mathfrak{rep}_{\mathrm{c}}(X;E_{0},\mathcal{F})$
and $\widehat{X}^{n}\in\mathfrak{rep}_{\mathrm{c}}(X^{n};E_{0},\mathcal{F})$
exist for each $n\in\mathbf{N}$. In addition, suppose $(E_{0},\mathscr{O}_{E}(E_{0}))$
is a Tychonoff space. Then:

\renewcommand{\labelenumi}{(\alph{enumi})}
\begin{enumerate}
\item If $\widehat{X}$ and $\{\widehat{X}^{n}\}_{n\in\mathbf{N}}$ satisfy%
\footnote{As specified in \S \ref{sec:RV}, (\ref{eq:RepProc_WC_D(E0)}) abbreviates
the statement that $\{\widehat{X}^{n}\}_{n\in\mathbf{N}}$ and $\widehat{X}$
are $\mathbb{D}_{0}$-valued random variables and the distributions
of $\{\widehat{X}^{n}\}_{n\in\mathbf{N}}$ converge weakly to that
of $\widehat{X}$ in $\mathcal{P}(\mathbb{D}_{0})$.%
} 
\begin{equation}
\widehat{X}^{n}\Longrightarrow\widehat{X}\mbox{ as }n\uparrow\infty\mbox{ on }D\left(\mathbf{R}^{+};E_{0},\mathscr{O}_{E}(E_{0})\right),\label{eq:RepProc_WC_D(E0)}
\end{equation}
then they satisfy
\begin{equation}
\widehat{X}^{n}\Longrightarrow\widehat{X}\mbox{ as }n\uparrow\infty\mbox{ on }D(\mathbf{R}^{+};\widehat{E}).\label{eq:RepProc_WC}
\end{equation}

\item If there exists an $S_{0}\subset E_{0}^{\mathbf{R}^{+}}$ satisfying
(\ref{eq:Path-Base}) and
\begin{equation}
\inf_{n\in\mathbf{N}}\mathbb{P}^{n}\left(X^{n}\in S_{0}\right)=1,\label{eq:Common_Path-Base_N}
\end{equation}
and if $\mathcal{F}$ strongly separates points on $E_{0}$, then
(\ref{eq:RepProc_WC}) implies (\ref{eq:RepProc_WC_D(E0)}).
\end{enumerate}
\end{prop}
\begin{proof}
(a) For ease of notation, we let $(\Omega^{0},\mathscr{F}^{0},\mathbb{P}^{0};X^{0})\circeq(\Omega,\mathscr{F},\mathbb{P};X)$,
$\mathbb{D}_{0}\circeq D(\mathbf{R}^{+};E_{0},\mathscr{O}_{E}(E_{0}))$
and $\widehat{\mathbb{D}}\circeq D(\mathbf{R}^{+};\widehat{E})$.
$(E_{0},\mathscr{O}_{E}(E_{0}))$ is a topological refinement of $(E_{0},\mathscr{O}_{\widehat{E}}(E_{0}))$
by Lemma \ref{lem:Base} (d). $D(\mathbf{R}^{+};E_{0},\mathscr{O}_{\widehat{E}}(E_{0}))$
is a subspace of $\widehat{\mathbb{D}}$ by Corollary \ref{cor:Sko_Subspace}
(with $E=\widehat{E}$ and $A=E_{0}$). It then follows by (\ref{eq:RepProc_WC_D(E0)})
and Proposition \ref{prop:Sko_Basic_1} (e) (with $E=(E_{0},\mathscr{O}_{\widehat{E}}(E_{0}))$
and $S=(E_{0},\mathscr{O}_{E}(E_{0}))$) that $\mathbb{D}_{0}\subset\widehat{\mathbb{D}}$,
$\widehat{\mathbb{D}}_{0}\circeq(\mathbb{D}_{0},\mathscr{O}_{\widehat{\mathbb{D}}}(\mathbb{D}_{0}))$
is a topological coarsening of $\mathbb{D}_{0}$ and
\begin{equation}
\widehat{X}^{n}\in M\left(\Omega^{n},\mathscr{F}^{n};\mathbb{D}_{0}\right)\subset M\left(\Omega^{n},\mathscr{F}^{n};\widehat{\mathbb{D}}_{0}\right),\;\forall n\in\mathbf{N}_{0}.\label{eq:RepProc_D(E0,Ehat)_RV_N0}
\end{equation}

Let $\mu_{n}$, $\widehat{\nu}_{n}$ and $\nu_{n}$ denote the distribution
of $\widehat{X}^{n}$ as $\mathbb{D}_{0}$-valued, $\widehat{\mathbb{D}}$-valued
and $\widehat{\mathbb{D}}_{0}$-valued random variables for each $n\in\mathbf{N}_{0}$,
respectively. It follows by (\ref{eq:RepProc_WC_D(E0)}) and Fact
\ref{fact:Weak_Topo_Coarsen} (with $E=\mathbb{D}_{0}$, $\mathscr{U}=\mathscr{O}(\widehat{\mathbb{D}}_{0})$
and $\mu=\mu_{0}$) that
\begin{equation}
\nu_{n}\Longrightarrow\nu_{0}\mbox{ as }n\uparrow\infty\mbox{ in }\mathcal{P}(\widehat{\mathbb{D}}_{0}).\label{eq:RepProc_WC_D(E0,Ehat)}
\end{equation}
It follows by (\ref{eq:RepProc_WC_D(E0,Ehat)}) and Lemma \ref{lem:WC_Expansion}
(with $E=\widehat{\mathbb{D}}$, $A=\widehat{\mathbb{D}}_{0}$, $\mu_{n}=\nu_{n}$
and $\mu=\nu_{0}$) that
\begin{equation}
\widehat{\nu}_{n}=\left.\nu_{n}\right|^{\widehat{\mathbb{D}}}\Longrightarrow\left.\nu_{0}\right|^{\widehat{\mathbb{D}}}=\widehat{\nu}_{0}\mbox{ as }n\uparrow\infty\mbox{ in }\mathcal{P}(\widehat{\mathbb{D}}),\label{eq:RepProc_WC_D(Ehat)}
\end{equation}
which proves (\ref{eq:RepProc_WC}).

(b) The given conditions imply
\begin{equation}
\inf_{n\in\mathbf{N}_{0}}\mathbb{P}^{n}\left(X^{n}\in\mathbb{D}_{0}\right)\geq\inf_{n\in\mathbf{N}_{0}}\mathbb{P}^{n}\left(X^{n}\in S_{0}\cap\mathbb{D}_{0}\right)=1.\label{eq:Check_Concentrate_D(E0)}
\end{equation}
We have $\mathscr{O}_{E}(E_{0})=\mathscr{O}_{\widehat{E}}(E_{0})$
by (\ref{eq:F_Fhat_Coincide}), (\ref{eq:Fhat_SSP_on_Ehat}) and $\mathcal{F}$
strongly separating points on $E_{0}$, which implies $\mathbb{D}_{0}=\widehat{\mathbb{D}}_{0}$.
According to Proposition \ref{prop:PR} (c) (with $S_{0}=\mathbb{D}_{0}$
and $X=X^{n}$ or $X$), one can take
\begin{equation}
\begin{aligned}\widehat{X}^{n} & =\mathfrak{rep}_{\mathrm{c}}(X^{n};E_{0},\mathcal{F})\\
 & \in M\left(\Omega^{n},\mathscr{F}^{n};\mathbb{D}_{0}\right)=M\left(\Omega^{n},\mathscr{F}^{n};\widehat{\mathbb{D}}_{0}\right)\subset M\left(\Omega^{n},\mathscr{F}^{n};\widehat{\mathbb{D}}\right),\;\forall n\in\mathbf{N}_{0}
\end{aligned}
\label{eq:RepProc_D(E0)_RV_N0}
\end{equation}
and each $\mu_{n}$, $\widehat{\nu}_{n}$ and $\nu_{n}$ in the proof
of (a) are all well-defined with
\begin{equation}
\mu_{n}=\nu_{n}\in\mathcal{P}(\widehat{\mathbb{D}}_{0})=\mathcal{P}(\mathbb{D}_{0}),\;\forall n\in\mathbf{N}_{0}.\label{eq:RepProc_D(E0)_D(E0,Ehat)_Same_Dist}
\end{equation}

(\ref{eq:RepProc_WC}) implies (\ref{eq:RepProc_WC_D(Ehat)}). As
$D(\mathbf{R}^{+};\widehat{E})$ is a Polish space, (\ref{eq:RepProc_WC_D(Ehat)})
implies (\ref{eq:RepProc_WC_D(E0,Ehat)}) by Lemma \ref{lem:WC_Expansion}
(with $E=\widehat{\mathbb{D}}$, $A=\widehat{\mathbb{D}}_{0}$, $\mu_{n}=\nu_{n}$
and $\mu=\nu_{0}$). Now, (\ref{eq:RepProc_WC_D(E0)}) follows by
(\ref{eq:RepProc_D(E0)_RV_N0}), (\ref{eq:RepProc_WC_D(E0,Ehat)})
and (\ref{eq:RepProc_D(E0)_D(E0,Ehat)_Same_Dist}).\end{proof}

\section{\label{sec:Base_Proc}Containment in large baseable subsets}

Given $E$-valued processes $\{(\Omega^{i},\mathscr{F}^{i},\mathbb{P}^{i};X^{i})\}_{i\in\mathbf{I}}$
and a base $(E_{0},\mathcal{F};\widehat{E},\widehat{\mathcal{F}})$
over $E$, most of the developments of \S \ref{sec:RepProc}, \S
\ref{sec:RepProc_FC} and \S \ref{sec:RepProc_Path_Space} require
$E_{0}$ to have containment properties like (\ref{eq:Common_FR-Base_I}),
\begin{equation}
\inf_{t\in\mathbf{T},i\in\mathbf{I}}\mathbb{P}^{i}\left(X_{t}^{i}\in E_{0}\right)=1\label{eq:Common_T-Base_I}
\end{equation}
or (\ref{eq:Common_Path-Base_I}) for $\{X^{i}\}_{i\in\mathbf{I}}$.
In practice, one usually constructs a basesable set $E_{0}$ satisfying
the non-functional%
\footnote{(\ref{eq:Common_FR-Base_I}) by contrast is a functional condition
depending on the choice of $\mathcal{F}$.%
} conditions (\ref{eq:Common_T-Base_I}) or (\ref{eq:Common_Path-Base_I})
first, and then selects proper functions to establish the base $(E_{0},\mathcal{F};\widehat{E},\widehat{\mathcal{F}})$.
From Fact \ref{fact:T-Base_FR-Base} we immediately observe that:
\begin{fact}
\label{fact:Common_T-Base_FR-Base}(\ref{eq:Common_FR-Base_I}), (\ref{eq:Common_T-Base_I})
and (\ref{eq:Common_Path-Base_I}) are successively stronger for any
index set $\mathbf{I}$ and $\mathbf{T}\subset\mathbf{R}^{+}$.
\end{fact}
The simplest case is when $E$ itself is a baseable space. Then, one
easily obtains a base $(E,\mathcal{F};\widehat{E},\widehat{\mathcal{F}})$
by Lemma \ref{lem:Base_Construction} and the containment properties
in Fact \ref{fact:Common_T-Base_FR-Base} are automatic. When $E$
is non-baseable, one can use $\mathbf{T}$-MPCC, MCCC, $\mathbf{T}$-PMTC
or $T_{k}$-LMTC introduced in \S \ref{sub:Proc_Reg} to construct
the desired $E_{0}$ in (\ref{eq:Common_Path-Base_I}) or (\ref{eq:Common_T-Base_I}).

When $\{X^{i}\}_{i\in\mathbf{I}}$ are all c$\grave{\mbox{a}}$dl$\grave{\mbox{a}}$g,
the following proposition uses $\mathbf{T}$-MPCC and $\mathfrak{r}$-MCC
to construct an $E_{0}$ satisfying (\ref{eq:Common_Path-Base_I}).
\begin{prop}
\label{prop:Base_Proc_MPCC}Let $(E,\mathfrak{r})$ be a metric space,
$\mathcal{D}\subset C(E;\mathbf{R})$ separate points on $E$, $\mathbf{T}$
be a dense subset of $\mathbf{R}^{+}$ and $\{(\Omega^{i},\mathscr{F}^{i},\mathbb{P}^{i};X^{i})\}_{i\in\mathbf{I}}$
be $E$-valued c$\grave{\mbox{a}}$dl$\grave{\mbox{a}}$g processes
satisfying $\mathbf{T}$-MPCC and $\mathfrak{r}$-MCC. Then, there
exist $\{A_{p,q}\}_{p,q\in\mathbf{N}}\subset\mathscr{C}(E)$ satisfying
the following properties:

\renewcommand{\labelenumi}{(\alph{enumi})}
\begin{enumerate}
\item $\{A_{p,q}\}_{p,q\in\mathbf{N}}$ are totally bounded and satisfy
\begin{equation}
A_{p,q}\subset A_{p,q+1},\;\forall p,q\in\mathbf{N}\label{eq:Pick_MPCC_Apq_1}
\end{equation}
and
\begin{equation}
\inf_{i\in\mathbf{I}}\mathbb{P}^{i}\left(X_{t}^{i}\in A_{p,q},\forall t\in[0,q]\right)\geq1-2^{-p-q},\;\forall p,q\in\mathbf{N}.\label{eq:Pick_MPCC_Apq_2}
\end{equation}

\item $E_{0}\circeq\bigcup_{p,q\in\mathbf{N}}A_{p,q}$ is a second-countable
subspace and is a $\mathcal{D}$-baseable subset of $E$.
\item $E_{0}$ and $S_{0}\circeq\bigcup_{p\in\mathbf{N}}V_{p}$ satisfy
(\ref{eq:Common_Path-Base_I}), where%
\footnote{$A_{p,q}^{[0,q)}$ in (\ref{eq:Vp_MPCC}) below means the Cartesian
power of $A_{p,q}$ for the index set $[0,q)$.%
}
\begin{equation}
V_{p}\circeq\bigcap_{q\in\mathbf{N}}\left\{ x\in\mathbb{D}_{0}:x|_{[0,q)}\in A_{p,q}^{[0,q)}\right\} ,\;\forall p\in\mathbf{N}.\label{eq:Vp_MPCC}
\end{equation}

\item If $(E,\mathfrak{r})$ is complete, then $\{A_{p,q}\}_{p,q\in\mathbf{N}}\subset\mathscr{K}(E)$
and $\{X^{i}\}_{i\in\mathbf{I}}$ satisfies MCCC in $E_{0}$.
\end{enumerate}
\end{prop}
As noted in \S \ref{sub:Metrizable_Compact}, metrizable compact
subsets provided by MCCC, $\mathbf{T}$-PMTC and $T_{k}$-LMTC are
nice baseable ``blocks'' for building $E_{0}$. In metric spaces,
totally bounded subsets provided by $\mathbf{T}$-MPCC form another
category of such blocks.
\begin{fact}
\label{fact:Total_Bounded_Baseable}Let $(E,\mathfrak{r})$ be a metric
space, $\mathcal{D}\subset C(E;\mathbf{R})$ separate points on $E$
and $\{A_{n}\}_{n\in\mathbf{N}}$ be totally bounded Borel subsets
of $E$. Then, $\bigcup_{n\in\mathbf{N}}A_{n}$ is a second-countable
subspace and, in particular, is a $\mathcal{D}$-baseable subset of
$E$.
\end{fact}
\begin{proof}
$A\circeq\bigcup_{n\in\mathbf{N}}A_{n}$ is a separable subspace of
$E$ by Proposition \ref{prop:Total_Bounded_1} (a) and Proposition
\ref{prop:Countability} (e). Now, the result follows by Proposition
\ref{prop:Metrizable} (c) and Proposition \ref{prop:Hered_Lindelof_Baseable}.\end{proof}

\begin{proof}
[Proof of Proposition \ref{prop:Base_Proc_MPCC}](a) An inspection
of the proof of \cite[Theorem 17]{K15} shows that $\mathbf{T}$-MPCC
is enough for their developments. So, one follows \cite{K15} to construct
totally bounded $\{A_{p,q}\}_{p,q\in\mathbf{N}}\subset\mathscr{C}(E)$
satisfying (a). 

(b) follows by (a) and Fact \ref{fact:Total_Bounded_Baseable}.

(c) One finds by (a) that
\begin{equation}
\begin{aligned}\inf_{i\in\mathbf{I}}\mathbb{P}^{i}\left(X^{i}\in V_{p}\right) & \geq1-\sup_{i\in\mathbf{I}}\sum_{q\in\mathbf{N}}\left[1-\mathbb{P}^{i}\left(X_{t}^{i}\in A_{p,q},\forall t\in[0,q]\right)\right]\\
 & \geq1-2^{-p},\;\forall p\in\mathbf{N}.
\end{aligned}
\label{eq:Check_Vp_MPCC_Large}
\end{equation}

(d) Each $(A_{p,q},\mathfrak{r})$ is complete by the fact $A_{p,q}\subset\mathscr{C}(E)$
and Proposition \ref{prop:Completeness} (c). Then, (d) follows by
Proposition \ref{prop:Total_Bounded_2}.\end{proof}

The next proposition uses MCCC to construct an $E_{0}$ satisfying
(\ref{eq:Common_Path-Base_I}).
\begin{prop}
\label{prop:Base_Proc_CCC}Let E be a topological space, $\mathcal{D}\subset C(E;\mathbf{R})$
separate points on $E$ and $\{(\Omega^{i},\mathscr{F}^{i},\mathbb{P}^{i};X^{i})\}_{i\in\mathbf{I}}$
be $E$-valued processes satisfying MCCC in $A\subset E$. Then, there
exist $\{K_{p,q}\}_{p,q\in\mathbf{N}}\subset\mathscr{K}^{\mathbf{m}}(E)$
satisfying the following properties:

\renewcommand{\labelenumi}{(\alph{enumi})}
\begin{enumerate}
\item (\ref{eq:Pick_CCC_Kpq_1}) and (\ref{eq:Pick_CCC_Kpq_2}) hold.
\item $E_{0}\circeq\bigcup_{p,q\in\mathbf{N}}K_{p,q}\subset A$ is a $\mathcal{D}$-baseable
subset of $E$. Moreover, $\{X^{i}\}_{i\in\mathbf{I}}$ satisfies
MCCC in $E_{0}$.
\item $E_{0}$ and $S_{0}\circeq\bigcup_{p\in\mathbf{N}}V_{p}$ satisfy
(\ref{eq:Common_Path-Base_I}), where $\{V_{p}\}_{p\in\mathbf{N}}$
are defined as in (\ref{eq:Vp_CCC}).
\end{enumerate}
\end{prop}
\begin{proof}
(a) We pick $\{A_{p,q}\}_{p,q\in\mathbf{N}}\subset\mathscr{K}^{\mathbf{m}}(E)$
satisfying $A_{p,q}\subset A$ for all $p,q\in\mathbf{N}$ and (\ref{eq:Pick_CCC_Apq}).
$E$ is a Hausdorff space by Proposition \ref{prop:Fun_Sep_1} (e)
(with $A=E$). Then, 
\begin{equation}
K_{p,q}\circeq\bigcup_{i=1}^{q}A_{p,i}\in\mathscr{K}^{\mathbf{m}}(E)\subset\mathscr{C}(E)\subset\mathscr{B}(E),\;\forall p,q\in\mathbf{N}\label{eq:Pick_CCC_Kpq}
\end{equation}
by Proposition \ref{prop:Separability} (c), Lemma \ref{lem:MC_Union}
and Proposition \ref{prop:Compact} (a). Now, (\ref{eq:Pick_CCC_Kpq_1})
and (\ref{eq:Pick_CCC_Kpq_2}) follow by (\ref{eq:Pick_CCC_Apq})
and (\ref{eq:Pick_CCC_Kpq}).

(b) follows by (\ref{eq:Pick_CCC_Kpq}), (a) and Proposition \ref{prop:Sigma_MC}
(b, e).

(c) $E_{0}$ and $S$ satisfy (\ref{eq:Check_Vp_CCC_Large}) by (a),
which implies (\ref{eq:Common_Path-Base_I}) immediately.\end{proof}

The following fact gives an $E_{0}$ satisfying (\ref{eq:Common_T-Base_I})
for countable $\mathbf{T}\subset\mathbf{R}^{+}$ by $\mathbf{T}$-PMTC.
\begin{fact}
\label{fact:Base_PMTC}Let E be a topological space, $A\subset E$
and $\{(\Omega^{i},\mathscr{F}^{i},\mathbb{P}^{i};X^{i})\}_{i\in\mathbf{I}}$
be $E$-valued processes. Then:

\renewcommand{\labelenumi}{(\alph{enumi})}
\begin{enumerate}
\item If $\mathbf{I}$ is an infinite set and $\{X^{i}\}_{i\in\mathbf{I}}$
satisfies $\mathbf{T}$-PMTC in $A$, then $\{X^{i}\}_{i\in\mathbf{I}}$
satisfies $\mathbf{T}$-PSMTC in $A$.
\item $\{X^{i}\}_{i\in\mathbf{I}}$ satisfies $\mathbf{T}$-PMTC in $A$
if and only if $\{X_{\mathbf{T}_{0}}^{i}\}_{i\in\mathbf{I}}$ is $\mathbf{m}$-tight
in $A^{\mathbf{T}_{0}}$ for all $\mathbf{T}_{0}\in\mathscr{P}_{0}(\mathbf{T})$.
\item If $\{X^{i}\}_{i\in\mathbf{I}}$ satisfies $\mathbf{T}$-PSMTC in
$A$, then $\{X_{\mathbf{T}_{0}}^{i}\}_{i\in\mathbf{I}}$ is sequentially
$\mathbf{m}$-tight in $A^{\mathbf{T}_{0}}$ for all $\mathbf{T}_{0}\in\mathscr{P}_{0}(\mathbf{T})$.
\item If $\{X^{i}\}_{i\in\mathbf{I}}$ satisfies $\mathbf{T}$-PMTC in $A$
for a countable $\mathbf{T}\subset\mathbf{R}^{+}$ and $\mathcal{D}\subset C(E;\mathbf{R})$
separates points on $E$, then there exists a $\mathcal{D}$-baseable
subset $E_{0}\in\mathscr{K}_{\sigma}^{\mathbf{m}}(E)$ such that $\{X^{i}\}_{i\in\mathbf{I}}$
satisfies $\mathbf{T}$-PMTC in $E_{0}\subset A$.
\item When $(E,\mathfrak{r})$ is a metric space, $\{X^{i}\}_{i\in\mathbf{I}}$
satisfying $\mathbf{T}$-PMTC implies $\{X^{i}\}_{i\in\mathbf{I}}$
satisfying $\mathbf{T}$-MPCC and the converse is true if $(E,\mathfrak{r})$
is complete.
\item If $\{X^{i}\}_{i\in\mathbf{I}}$ satisfies MCCC in $A$, then $\{X^{i}\}_{i\in\mathbf{I}}$
satisfies $\mathbf{R}^{+}$-PMTC in $A$.
\end{enumerate}
\end{fact}
\begin{proof}
(a) and (f) are automatic by definition. (b) and (c) follow by Lemma
\ref{lem:Tightness_Prod} (with $\mathbf{I}=\mathbf{T}_{0}$, $S_{i}=E$,
$A_{i}=A$ and $\Gamma=\{\mathbb{P}^{i}\circ(X_{\mathbf{T}_{0}}^{i})^{-1}\}_{i\in\mathbf{I}}$).
(d) follows by Lemma \ref{lem:m-Tight_Base} (with $\mathbf{I}=\mathbf{T}$,
$i=t$ and $\Gamma_{i}=\{\mathbb{P}^{i}\circ(X_{t}^{i})^{-1}\}_{i\in\mathbf{I}}$).

(e) The first statement is immediate by Proposition \ref{prop:Total_Bounded_2}.
Then, we suppose $\{X^{i}\}_{i\in\mathbf{I}}$ satisfies $\mathbf{T}$-MPCC
and $(E,\mathfrak{r})$ is complete. Fixing $t\in\mathbf{T}$ and
$\epsilon\in(0,\infty)$, we pick a totally bounded $A_{\epsilon,p,t}\subset E$
for each $p\in\mathbf{N}$ such that
\begin{equation}
\inf_{i\in\mathbf{I}}\mathbb{P}^{i}\left(X_{t}^{i}\in A_{\epsilon,p,t}^{\epsilon2^{-p}}\right)\geq1-\epsilon2^{-p},\;\forall p\in\mathbf{N}.\label{eq:MPCC_p}
\end{equation}
Letting $K{}_{\epsilon,t}$ be the closure of $\bigcap_{p\in\mathbf{N}}A_{\epsilon,p,t}^{\epsilon2^{-p}}$,
one finds that 
\begin{equation}
\inf_{i\in\mathbf{I}}\mathbb{P}^{i}\left(X_{t}^{i}\in K_{\epsilon,t}\right)\geq1-\epsilon,\;\forall p\in\mathbf{N}.\label{eq:Check_PMTC}
\end{equation}
$\bigcap_{p\in\mathbf{N}}A_{\epsilon,p,t}^{\epsilon2^{-p}}$ is totally
bounded by definition and so is $K_{\epsilon,t}$ by Proposition \ref{prop:Total_Bounded_1}
(c). $(K_{\epsilon,t},\mathfrak{r})$ is complete by Proposition \ref{prop:Completeness}
(c). Hence, $K_{\epsilon,t}\in\mathscr{K}(E)$ by Proposition \ref{prop:Total_Bounded_2}
and (e) follows by (\ref{eq:Check_PMTC}).\end{proof}

Given countably many processes satisfying $T_{k}$-LMTC, the next
proposition constructs an $E_{0}$ satisfying (\ref{eq:Common_T-Base_I})
for a conull $\mathbf{T}\subset\mathbf{R}^{+}$.
\begin{prop}
\label{prop:Base_Proc_LMTC}Let E be a topological space, $\mathcal{D}\subset C(E;\mathbf{R})$
separate points on $E$ and $\{(\Omega^{n},\mathscr{F}^{n},\mathbb{P}^{n};X^{n})\}_{n\in\mathbf{N}}$
be $E$-valued measurable processes satisfying $T_{k}$-LMTC%
\footnote{Recall that the definition of $T_{k}$-LMTC includes $T_{k}\uparrow\infty$.%
} in $A\subset E$. Then, there exists a $\mathcal{D}$-baseable subset
$E_{0}\in\mathscr{K}_{\sigma}^{\mathbf{m}}(E)$ with $A\supset E_{0}$
and a conull $\mathbf{T}\subset\mathbf{R}^{+}$ such that
\begin{equation}
\inf_{t\in\mathbf{T},n\in\mathbf{N}}\mathbb{P}^{n}\left(X_{t}^{n}\in E_{0}\right)=1\label{eq:Common_T-Base_N}
\end{equation}
and $\{X^{n}\}_{n\in\mathbf{N}}$ satisfies $T_{k}$-LMTC in $E_{0}$.
\end{prop}
\begin{proof}
We take $\{K_{p}\}_{p\in\mathbf{N}}\subset\mathscr{K}^{\mathbf{m}}(A,\mathscr{O}_{E}(A))$
satisfying
\begin{equation}
\inf_{k,n\in\mathbf{N}}\frac{1}{T_{k}}\int_{0}^{T_{k}}\mathbb{P}^{n}\left(X_{\tau}^{n}\in K_{p}\right)d\tau\geq1-2^{-p},\;\forall p\in\mathbf{N}\label{eq:Pick_LT_Kp}
\end{equation}
and let $E_{0}\circeq\bigcup_{p\in\mathbf{N}}K_{p}\in\mathscr{K}_{\sigma}^{\mathbf{m}}(A,\mathscr{O}_{E}(A))$.
It follows that
\begin{equation}
\sup_{k,n\in\mathbf{N}}\int_{0}^{T_{k}}\mathbb{P}^{n}\left(X_{\tau}^{n}\notin E_{0}\right)d\tau=0\label{eq:LT_E0_Conull-1}
\end{equation}
by (\ref{eq:Pick_LT_Kp}) and continuity of measure. Hence, (\ref{eq:Common_T-Base_N})
holds for the conull set
\begin{equation}
\mathbf{T}\circeq\mathbf{R}^{+}\backslash\bigcup_{k,n\in\mathbf{N}}\left\{ t\in[0,T_{k}]:\mathbb{P}^{n}\left(X_{t}^{n}\notin E_{0}\right)>0\right\} .\label{eq:LT_E0_Conull-3}
\end{equation}
Now, the result follows by (\ref{eq:Pick_LT_Kp}) and Proposition
\ref{prop:Sigma_MC} (b, e) (with $A=E_{0}$).\end{proof}

The relationship among MCCC, $\mathbf{T}$-MPCC, $\mathbf{T}$-PMTC
and $T_{k}$-LMTC is illustrated in Figure \ref{fig:tight} below,
where green solid arrows means definite implication, blue dashed arrow
means conditional implication and red crossed arrow means false converse.

\begin{figure}[H]
\begin{centering}
\includegraphics[scale=0.95]{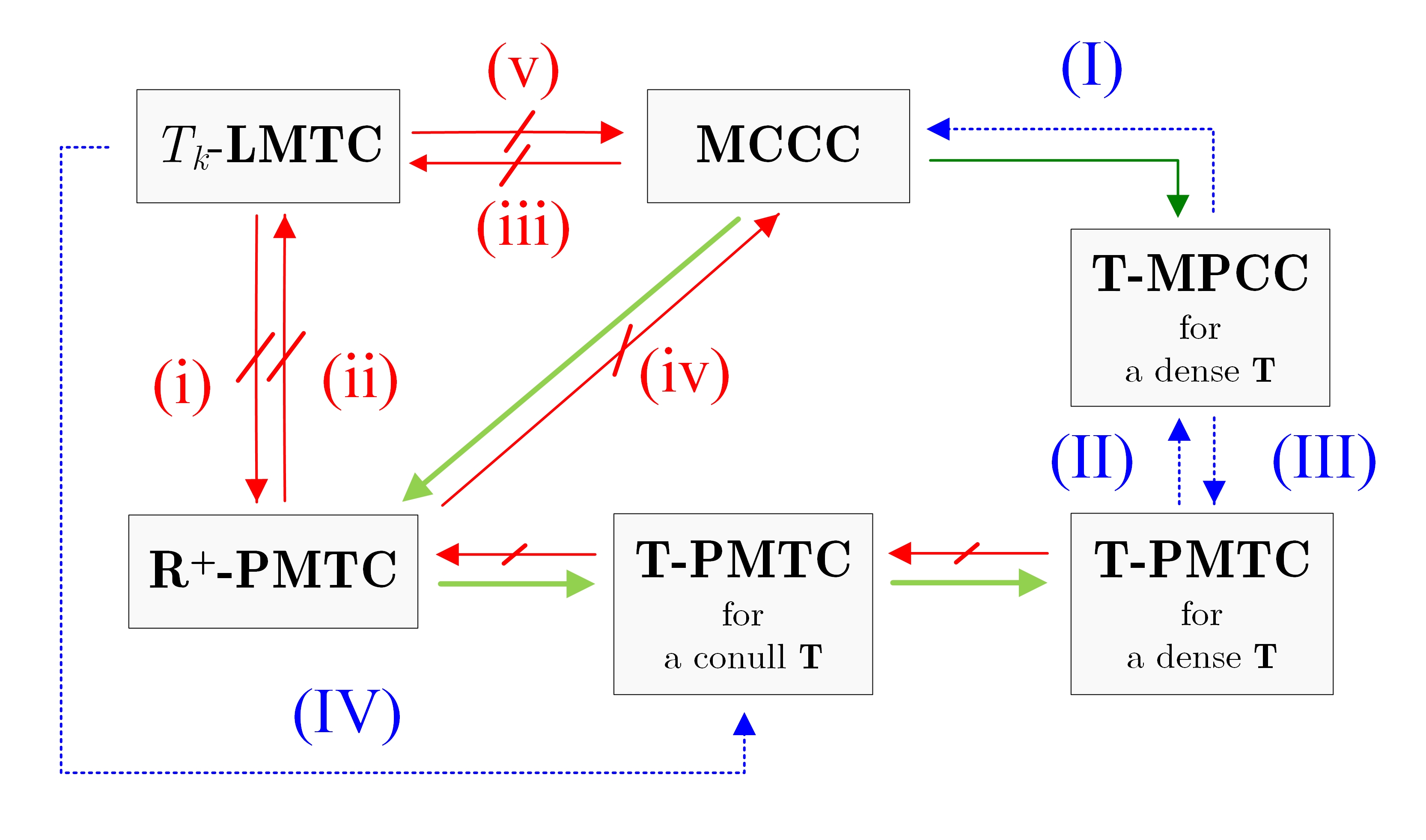}
\par\end{centering}

\caption{\label{fig:tight}\textit{The relationship among tightness/containment
conditions}}
\end{figure}

\begin{rem}
\label{rem:Containment_Figure}All the unlabelled arrows in Figure
\ref{fig:tight} are immediate. Below is some explanation for the
labelled ones:
\begin{itemize}
\item (I) was justified in Proposition \ref{prop:Base_Proc_MPCC} (a, d)
for c$\grave{\mbox{a}}$dl$\grave{\mbox{a}}$g processes living on
a complete (but not necessarily separable) metric space $(E,\mathfrak{r})$
and satisfying $\mathfrak{r}$-MCC. This is a generalization of \cite[Theorem 17]{K15}
on infinite time horizon since $\mathbf{T}$-MPCC with a dense $\mathbf{T}$
is weaker than the Pointwise Containment Property in \cite[\S 5]{K15}.
\item By Fact \ref{fact:Base_PMTC} (e), (II) is true on arbitrary metric
spaces and (III) is true on complete metric spaces.
\item (IV) was justified in Proposition \ref{prop:Base_Proc_LMTC} for a
countable collection of measurable processes.
\item (i) is not true because $T_{k}$-LMTC will not be affected by changing
the distributions of $\{X_{t}^{i}\}_{i\in\mathbf{I}}$ to a non-tight
family for each $t\in\mathbf{Q}^{+}$.
\item (ii) and (iii) are disproved by the constant process $\{t\}_{t\geq0}$.
\item (iv) and (v) are disproved by Example \ref{exp:TightnessCompare}
below, where we construct a non-stationary c$\grave{\mbox{a}}$dl$\grave{\mbox{a}}$g
process that satisfies $T_{k}$-LMTC and $\mathbf{R}^{+}$-PMTC but
violates MCCC. (iv) was also disproved by \cite[Example 2]{K15}.
\end{itemize}
\end{rem}
\begin{example}
\label{exp:TightnessCompare}Let $\mu$ be the uniform distribution
on $(0,1)$ and 
\begin{equation}
\eta_{t}(\omega)\circeq\begin{cases}
1-\omega+t, & \mbox{if }t\in[0,\omega),\\
\frac{1}{2}, & \mbox{if }t\in[\omega,\infty),
\end{cases}\;\forall\omega\in(0,1),t\in\mathbf{R}^{+}.\label{eq:eta}
\end{equation}
$\eta=\{\eta_{t}\}_{t\geq0}$ satisfies $\mathbf{R}^{+}$-PMTC since
$(0,1)$ is $\sigma$-compact. However, $\eta$ violates MCCC because
for any $a,b\in(0,1)$,
\begin{equation}
\begin{aligned} & \mu\left(\eta_{t}\in[a,b],\forall t\in[0,1]\right)\\
 & \leq1-\mu\left(\left\{ \omega\in(0,1):0\leq\omega-t<1-b,\exists t\in[0,\omega)\right\} \right)=0.
\end{aligned}
\label{eq:CCC_False}
\end{equation}
For each $\tau>0$ and $\epsilon\in(0,1/2)$,
\begin{equation}
\left\{ \eta_{\tau}\in[\epsilon,1-\epsilon]\right\} =\left((\tau\wedge1)\vee(\epsilon+\tau),1\wedge(1+\tau-\epsilon)\right)\cup(0,\tau\wedge1).\label{eq:Check_LT-1}
\end{equation}
Letting $T>1/\epsilon$, one finds by (\ref{eq:Check_LT-1}) that
\begin{equation}
\frac{1}{T}\int_{0}^{T}\mu\left(\eta_{\tau}\in[\epsilon,1-\epsilon]\right)d\tau\geq\frac{1}{T}\int_{1}^{T}1d\tau\geq1-\epsilon.\label{eq:Check_LT-2}
\end{equation}
Hence, $\eta$ satisfies $T_{k}$-LMTC for any $T_{k}\uparrow\infty$.
Moreover, $\eta$ is non-stationary since $\eta_{0}$ and $\eta_{1/2}$
have distinct expectations.
\end{example}

\chapter{\label{chap:FDDConv}Application to Finite-Dimensional Convergence}

\chaptermark{Finite-Dimensional Convergence}

The previous four chapters elaborate \ref{enu:Theme1} of this work.
With the help of replication, we have developed in \S \ref{sec:RepProc_FC}
several tool results for \ref{enu:Theme2}, the finite-dimensional
convergence of possibly non-c$\grave{\mbox{a}}$dl$\grave{\mbox{a}}$g
processes. Now, we are going to answer the target questions \ref{enu:Q_LTB},
\ref{enu:Q_FLP_Gen} and \ref{enu:Q_FLP_Prog} of \ref{enu:Theme2}
in the following three sections. \S \ref{sec:FLP_Gen}, establishing
finite-dimensional convergence to processes with general paths, answers
\ref{enu:Q_FLP_Gen}. \S \ref{sec:FLP_Prog}, establishing finite-dimensional
convergence of weakly c$\grave{\mbox{a}}$dl$\grave{\mbox{a}}$g processes
to weakly c$\grave{\mbox{a}}$dl$\grave{\mbox{a}}$g or progressive
limit processes, provides answers to both \ref{enu:Q_FLP_Gen} and
\ref{enu:Q_FLP_Prog}. In \S \ref{sec:LTB}, we answer \ref{enu:Q_LTB}
by establishing finite-dimensional convergence to long-time typical
behaviors of a given measurable process.

\section{\label{sec:FLP_Gen}Convergence of process with general paths}

Let $\{X^{i}\}_{i\in\mathbf{I}}$ be infinitely many $E$-valued processes
and $\mathbf{S}\subset\mathbf{R}^{+}$. We give in this section a
set of sufficient conditions for the unique existence of $X\in\mathfrak{flp}_{\mathbf{S}}(\{X^{i}\}_{i\in\mathbf{I}})$%
\footnote{The readers are referred to \S \ref{sec:RepProc_FC} for definitions
and notations about finite-dimensional convergence, finite-dimensional
limit point and finite-dimensional limit.%
}. The nature of establishing an $X\in\mathfrak{flp}_{\mathbf{S}}(\{X^{i}\}_{i\in\mathbf{I}})$
with general paths is establishing a Kolmogorov extension of weak
limit points of the finite-dimensional distributions of $\{X^{i}\}_{i\in\mathbf{I}}$
for each $\mathbf{T}_{0}\in\mathscr{P}_{0}(\mathbf{S})$. Hence, our
goal can be achieved by directly applying Theorem \ref{thm:WLP_Consistency}
established in \S \ref{sec:WLP_Uni}.
\begin{thm}
\label{thm:FLP_Gen}Let $E$ be a topological space, $\{(\Omega^{i},\mathscr{F}^{i},\mathbb{P}^{i};X^{i})\}_{i\in\mathbf{I}}$
be $E$-valued processes, $\mathcal{D}\subset C_{b}(E;\mathbf{R})$
separate points on $E$%
\footnote{As mentioned in Note \ref{note:Cb(E;R)_SP_Baseable}, the assumption
of $\mathcal{D}\subset C_{b}(E;\mathbf{R})$ separating points on
$E$ below does not require $E$ to be a Tychonoff or baseable space.%
} and $\mathbf{S}\subset\mathbf{R}^{+}$. Then:

\renewcommand{\labelenumi}{(\alph{enumi})}
\begin{enumerate}
\item If $\{X^{i}\}_{i\in\mathbf{I}}$ satisfies $\mathbf{S}$-PSMTC%
\footnote{$\mathbf{S}$-PMTC (in $A$) and $\mathbf{S}$-PSMTC (in $A$) were
introduced in \S \ref{sub:Proc_Reg}. As specified in Note \ref{note:Singleton_PMTC},
that $X$ satisfies $\mathbf{S}$-PMTC in $A$ means the singleton
$\{X\}$ satisfies $\mathbf{S}$-PMTC in $A$. %
} in $A\subset E$, then any $X\in\mathfrak{flp}_{\mathbf{S}}(\{X^{i}\}_{i\in\mathbf{I}})$
satisfies $\mathbf{S}$-PMTC in $A$.
\item If $\{X^{i}\}_{i\in\mathbf{I}}$ is $(\mathbf{S},\mathcal{D})$-FDC%
\footnote{The notions of $(\mathbf{S},\mathcal{D})$-FDC and $(\mathbf{T},\mathcal{D})$-AS
was introduced in \S \ref{sec:RepProc_FC}.%
} and satisfies $\mathbf{S}$-PSMTC, then there exists an $X=\mathfrak{flp}_{\mathbf{S}}(\{X^{i}\}_{i\in\mathbf{I}})$
satisfying $\mathbf{S}$-PMTC and $X=\mathfrak{fl}_{\mathbf{S}}(\{X^{i_{n}}\}_{n\in\mathbf{N}})$
for any $\{i_{n}\}_{n\in\mathbf{N}}\subset\mathbf{I}$.
\item If $\{X^{i}\}_{i\in\mathbf{I}}$ is $(\mathbf{R}^{+},\mathcal{D})$-FDC,
is $(\mathbf{R}^{+},\mathcal{D})$-AS and satisfies $\mathbf{R}^{+}$-PSMTC,
then there exists a stationary $X=\mathfrak{flp}_{\mathbf{R}^{+}}(\{X^{i}\}_{i\in\mathbf{I}})$
satisfying $\mathbf{R}^{+}$-PMTC and $X=\mathfrak{fl}_{\mathbf{R}^{+}}(\{X^{i_{n}}\}_{n\in\mathbf{N}})$
for any $\{i_{n}\}_{n\in\mathbf{N}}\subset\mathbf{I}$.
\end{enumerate}
\end{thm}
\begin{rem}
\label{rem:Lim_Proc_PMTC}We pointed out in Fact \ref{fact:Base_PMTC}
(b) (with $\mathbf{T}=\mathbf{S}$) that $X$ satisfying $\mathbf{S}$-PMTC
in $A$ is equivalent to $X_{\mathbf{T}_{0}}$ being $\mathbf{m}$-tight
in $A^{\mathbf{T}_{0}}$ for every $\mathbf{T}_{0}\in\mathscr{P}_{0}(\mathbf{S})$.
In particular, $X$ satisfying $\mathbf{R}^{+}$-PMTC is equivalent
to all finite-dimensional distributions of $X$ being $\mathbf{m}$-tight.
\end{rem}
\begin{proof}
[Proof of Theorem \ref{thm:FLP_Gen}](a) Suppose $X$ is defined on
$(\Omega,\mathscr{F},\mathbb{P})$ and
\begin{equation}
X^{i_{n}}\xrightarrow{\quad\mathrm{D}(\mathbf{S})\quad}X\mbox{ as }n\uparrow\infty.\label{eq:FC_Subseq_along_S}
\end{equation}
It follows by (\ref{eq:FC_Subseq_along_S}) and Fact \ref{fact:Proc_Basic_1}
(d) that
\begin{equation}
\mathbb{P}^{i_{n}}\circ(X_{t}^{i_{n}})^{-1}\Longrightarrow\mathbb{P}\circ X_{t}^{-1}\mbox{ as }n\uparrow\infty\mbox{ in }\mathcal{P}(E),\;\forall t\in\mathbf{S}.\label{eq:FC_Subseq_1-D_Marginal_WC}
\end{equation}
$\{X_{t}^{i_{n}}\}_{n\in\mathbf{N}}$ is sequentially $\mathbf{m}$-tight
in $A$ for all $t\in\mathbf{S}$ as $\{X^{i}\}_{i\in\mathbf{I}}$
satisfies $\mathbf{S}$-PSMTC in $A$. Hence, (a) follows by (\ref{eq:FC_Subseq_1-D_Marginal_WC})
and Lemma \ref{lem:Compact_Portmanteau} (b) (with $\Gamma=\{\mathbb{P}^{i_{n}}\circ(X_{t}^{i_{n}})^{-1}\}_{n\in\mathbf{N}}$
and $\mu=\mathbb{P}\circ X_{t}^{-1}$).

(b) For each $i\in\mathbf{I}$, we let $\mu_{i}$ be the restriction
of $\mathrm{pd}(X)$ to $\mathscr{B}(E)^{\otimes\mathbf{R}^{+}}$%
\footnote{Restriction of measure to sub-$\sigma$-algebra and $X$'s process
distribution $\mathrm{pd}(X)$ were specified in \S \ref{sub:Meas}
and \S \ref{sec:Proc} respectively.%
}. For each $\mathbf{T}_{0}\in\mathscr{P}_{0}(\mathbf{S})$, the probability
measures
\begin{equation}
\mu_{i}\circ\mathfrak{p}_{\mathbf{T}_{0}}^{-1}=\mathbb{P}^{i}\circ(X_{\mathbf{T}_{0}}^{i})^{-1}\in\mathfrak{P}\left(E^{\mathbf{T}_{0}},\mathscr{B}(E)^{\otimes\mathbf{T}_{0}}\right),\;\forall i\in\mathbf{I}\label{eq:Proc_FDD_Prob_Meas}
\end{equation}
form a sequentially $\mathbf{m}$-tight family by Fact \ref{fact:Base_PMTC}
(c) (with $A=E$). For each $\mathbf{T}_{0}\in\mathscr{P}_{0}(\mathbf{S})$
and $f\in\mathfrak{mc}[\Pi^{\mathbf{T}_{0}}(\mathcal{D})]\cup\{1\}$,
the integrals
\begin{equation}
\int_{E^{\mathbf{T}_{0}}}f(x)\mu_{i}\circ\mathfrak{p}_{\mathbf{T}_{0}}^{-1}(dx)=\mathbb{E}^{i}\left[f\circ X_{\mathbf{T}_{0}}^{i}\right],\;\forall i\in\mathbf{I}\label{eq:Proc_FDD_Int_Test}
\end{equation}
admit at most one limit point in $\mathbf{R}$ since $\{X^{i}\}_{i\in\mathbf{I}}$
is $(\mathbf{S},\mathcal{D})$-FDC. Hence, it follows by Theorem \ref{thm:WLP_Consistency}
(with $\Gamma=\{\mu_{i}\}_{i\in\mathbf{I}}$, $\mathbf{I}=\mathbf{S}$,
$\mathbf{I}_{0}=\mathbf{T}_{0}$ and $a=b=1$) that there exists a
unique $\mu\in\mathfrak{P}(E^{\mathbf{S}},\mathscr{B}(E)^{\otimes\mathbf{S}})$
and some $\{\mathbf{I}_{\mathbf{T}_{0}}\in\mathscr{P}_{0}(\mathbf{I})\}_{\mathbf{T}_{0}\in\mathscr{P}_{0}(\mathbf{S})}$
such that $\mu\circ\mathfrak{p}_{\mathbf{T}_{0}}^{-1}\in\mathcal{P}(E^{\mathbf{T}_{0}})$
is the weak limit of \textit{any subsequence of} and, hence, is the
unique weak limit point of $\{\mu_{\mathbf{T}_{0},i}=\mathfrak{be}(\mu_{i}\circ\mathfrak{p}_{\mathbf{T}_{0}}^{-1})\}_{i\in\mathbf{I}\backslash\mathbf{I}_{\mathbf{T}_{0}}}$
for all $\mathbf{T}_{0}\in\mathscr{P}_{0}(\mathbf{S})$.

We fix $t_{0}\in\mathbf{S}$ and define
\begin{equation}
X_{t}\circeq\begin{cases}
\mathfrak{p}_{t}, & \mbox{if }t\in\mathbf{S},\\
\mathfrak{p}_{t_{0}}, & \mbox{if }t\in\mathbf{R}^{+}\backslash\mathbf{S},
\end{cases}\;\forall t\in\mathbf{R}^{+}.\label{eq:Define_FLP_Gen}
\end{equation}
By Fact \ref{fact:Prod_Map_1} (a) and Fact \ref{fact:Proc_Basic_1}
(b), $X\circeq\{X_{t}\}_{t\geq0}$ well-defines an $E$-valued process
on the probability space $(E^{\mathbf{S}},\mathscr{B}(E)^{\otimes\mathbf{S}},\mu)$
that satisfies
\begin{equation}
\mu\circ X_{\mathbf{T}_{0}}^{-1}=\mu\circ\mathfrak{p}_{\mathbf{T}_{0}}^{-1},\;\forall\mathbf{T}_{0}\in\mathscr{P}_{0}(\mathbf{S}).\label{eq:FLP_Gen_FDD}
\end{equation}
Now, (b) follows by (a) and Fact \ref{fact:FLP_Uni} (with $\mathbf{T}=\mathbf{S}$).

(c) One obtains by (b) (with $\mathbf{S}=\mathbf{R}^{+}$) an $X=\mathfrak{flp}_{\mathbf{R}^{+}}(\{X^{i}\}_{i\in\mathbf{I}})$
satisfying all conclusions of (c) except for stationarity. Suppose
$X$ is defined on $(\Omega,\mathscr{F},\mathbb{P})$ and
\begin{equation}
X^{i_{n}}\xrightarrow{\quad\mathrm{D}(\mathbf{R}^{+})\quad}X\mbox{ as }n\uparrow\infty.\label{eq:FC_Subseq_along_R+}
\end{equation}
Fixing $c\in(0,\infty)$ and $\mathbf{T}_{0}\in\mathscr{P}_{0}(\mathbf{R}^{+})$,
it follows by (\ref{eq:FC_Subseq_along_R+}), Fact \ref{fact:FDC_AS}
(b) (with $n=i_{n}$) and Fact \ref{fact:FC_FDC} (with $n=i_{n}$)
that
\begin{equation}
\mathbb{E}\left[f\circ X_{\mathbf{T}_{0}}\right]-\mathbb{E}\left[f\circ X_{\mathbf{T}_{0}+c}\right]=\lim_{n\rightarrow\infty}\mathbb{E}^{n}\left[f\circ X_{\mathbf{T}_{0}}^{n}-f\circ X_{\mathbf{T}_{0}+c}^{n}\right]=0\label{eq:FLP_Gen_Sta_Int_Test}
\end{equation}
for all $f\in\mathfrak{mc}[\Pi^{\mathbf{T}_{0}}(\mathcal{D})]\cup\{1\}$.
Hence, $\mathbb{P}\circ X_{\mathbf{T}_{0}}^{-1}=\mathbb{P}\circ X_{\mathbf{T}_{0}+c}^{-1}$
by their $\mathbf{m}$-tightness and Lemma \ref{lem:Tight_Meas_Identical}
(b) (with $d=\aleph(\mathbf{T}_{0})$).\end{proof}

Theorem \ref{thm:FLP_Gen} can be used to identify a given $E$-valued
limit process as the unique finite-dimensional limit point.
\begin{cor}
\label{cor:FLP_Gen}Let $E$ be a topological space, $\mathbf{S}\subset\mathbf{R}^{+}$
and $\{(\Omega^{i},\mathscr{F}^{i},\mathbb{P}^{i};X^{i})\}_{i\in\mathbf{I}}$
and $(\Omega,\mathscr{F},\mathbb{P};X)$ be $E$-valued processes.
Suppose that:

\renewcommand{\labelenumi}{(\roman{enumi})}
\begin{enumerate}
\item $\mathcal{D}\subset C_{b}(E;\mathbf{R})$ separates points on $E$.
\item $\{X^{i}\}_{i\in\mathbf{I}}$ satisfies $\mathbf{S}$-PSMTC.
\item $X$ satisfies $\mathbf{S}$-PMTC%
\footnote{$X$ satisfying $\mathbf{S}$-PMTC means the singleton $\{X\}$ satisfies
$\mathbf{S}$-PMTC.%
}.
\item $\mathbb{E}[f\circ X_{\mathbf{T}_{0}}]$ is the unique limit point
of $\{\mathbb{E}^{i}[f\circ X_{\mathbf{T}_{0}}^{i}]\}_{i\in\mathbf{I}}$
in $\mathbf{R}$ for all $f\in\mathfrak{mc}[\Pi^{\mathbf{T}_{0}}(\mathcal{D})]$
and $\mathbf{T}_{0}\in\mathscr{P}_{0}(\mathbf{S})$.
\end{enumerate}
Then:

\renewcommand{\labelenumi}{(\alph{enumi})}
\begin{enumerate}
\item $X=\mathfrak{flp}_{\mathbf{S}}(\{X^{i}\}_{i\in\mathbf{I}})$ and $X=\mathfrak{fl}_{\mathbf{S}}(\{X^{i_{n}}\}_{n\in\mathbf{N}})$
for any $\{i_{n}\}_{n\in\mathbf{N}}\subset\mathbf{I}$.
\item If $\mathbf{S}=\mathbf{R}^{+}$ and $\{X^{i}\}_{i\in\mathbf{I}}$
is $(\mathbf{R}^{+},\mathcal{D})$-AS, then $X=\mathfrak{flp}_{\mathbf{R}^{+}}(\{X^{i}\}_{i\in\mathbf{I}})$
is a sationary process and $X=\mathfrak{fl}_{\mathbf{R}^{+}}(\{X^{i_{n}}\}_{n\in\mathbf{N}})$
for any $\{i_{n}\}_{n\in\mathbf{N}}\subset\mathbf{I}$.
\end{enumerate}
\end{cor}
\begin{proof}
(a) $\{X^{i}\}_{i\in\mathbf{I}}$ is $(\mathbf{S},\mathcal{D})$-FDC
by the condition (iv) above. By Theorem \ref{thm:FLP_Gen} (b), there
exists a $Z=\mathfrak{flp}_{\mathbf{S}}(\{X^{i}\}_{i\in\mathbf{I}})$
such that $Z$ satisfies $\mathbf{S}$-PMTC and $Z=\mathfrak{fl}_{\mathbf{S}}(\{X^{i_{n}}\}_{n\in\mathbf{N}})$
for any $\{i_{n}\}_{n\in\mathbf{N}}\subset\mathbf{I}$. We take $Z$
defined on $(\Omega,\mathscr{F},\mathbb{P})$ for simplicity, fix
$\{i_{n}\}_{n\in\mathbf{N}}\subset\mathbf{I}$ and show $\mathbb{P}\circ X_{\mathbf{T}_{0}}^{-1}=\mathbb{P}\circ Z_{\mathbf{T}_{0}}^{-1}$
for all $\mathbf{T}_{0}\in\mathscr{P}_{0}(\mathbf{S})$. Since
\begin{equation}
X^{i_{n}}\xrightarrow{\quad\mathrm{D}(\mathbf{S})\quad}Z\mbox{ as }n\uparrow\infty,\label{eq:FC_Subseq_along_S_Z}
\end{equation}
we have by (iv) and Fact \ref{fact:FC_FDC} (with $X=Z$) that
\begin{equation}
\mathbb{E}\left[f\circ X_{\mathbf{T}_{0}}\right]=\mathbb{E}\left[f\circ Z_{\mathbf{T}_{0}}\right]\label{eq:Lim_Proc_Int_Compare}
\end{equation}
for all $f\in\mathfrak{mc}[\Pi^{\mathbf{T}_{0}}(\mathcal{D})]\cup\{1\}$.
$X_{\mathbf{T}_{0}}$ and $Z_{\mathbf{T}_{0}}$ are $\mathbf{m}$-tight
by the condition (iii) above and Remark \ref{rem:Lim_Proc_PMTC} (with
$X=X$ or $Z$). Hence, $\mathbb{P}\circ X_{\mathbf{T}_{0}}^{-1}=\mathbb{P}\circ Z_{\mathbf{T}_{0}}^{-1}$
by Lemma \ref{lem:Tight_Meas_Identical} (b) (with $d=\aleph(\mathbf{T}_{0})$).

(b) follows by (a) (with $\mathbf{S}=\mathbf{R}^{+}$) and Theorem
\ref{thm:FLP_Gen} (c).\end{proof}

\begin{rem}
\label{rem:FLP_Gen_E0}A variant of Theorem \ref{thm:FLP_Gen} will
be given in Proposition \ref{prop:FLP_Gen_E0} that relies heavily
on our results herein.
\end{rem}

\section{\label{sec:FLP_Prog}Convergence of weakly c$\grave{\mbox{a}}$dl$\grave{\mbox{a}}$g
processes}

This section deals with finite-dimensional convergence of $E$-valued
processes satisfying $\mathcal{D}$-FMCC%
\footnote{The notion of $\mathcal{D}$-FMCC was introduced in Definition \ref{def:Proc_Reg}.%
}. Such processes are $(\mathbf{R}^{+},\mathcal{D})$-c$\grave{\mbox{a}}$dl$\grave{\mbox{a}}$g%
\footnote{The notion of $(\mathbf{S},\mathcal{F})$-c$\grave{\mbox{a}}$dl$\grave{\mbox{a}}$g
process was introduced in Definition \ref{def:Weakly_Cadlag}.%
} according to Note \ref{note:FMCC_Weakly_Cadlag}.

Given $\{X^{n}\}_{n\in\mathbf{N}}$ satisfying $\mathcal{D}$-DMCC,
part (a) of the next theorem establishes an $(\mathbf{S},\mathcal{D})$-c$\grave{\mbox{a}}$dl$\grave{\mbox{a}}$g
$X\in\mathfrak{flp}_{\mathbf{S}}(\{X^{n}\}_{n\in\mathbf{N}})$ and
gives an alternative answer to \ref{enu:Q_FLP_Gen} in Introduction.
Part (b) further imposes the standard Borel property and establishes
a progressive member of $\mathfrak{flp}_{\mathbf{S}}(\{X^{n}\}_{n\in\mathbf{N}})$.
In lieu of a standard Borel assumption, part (c) assumes the $(\mathbf{T},\mathcal{D})$-AS
of $\{X^{n}\}_{n\in\mathbf{N}}$ for a conull $\mathbf{T}\supset\mathbf{S}$
and establishes a stationary and progressive member of $\mathfrak{flp}_{\mathbf{S}}(\{X^{n}\}_{n\in\mathbf{N}})$.
These two parts provide answers to \ref{enu:Q_FLP_Prog} in Introduction.
\begin{thm}
\label{thm:FL_Weakly_Cadlag}Let $E$ be a topological space, $\{(\Omega^{n},\mathscr{F}^{n},\mathbb{P}^{n};X^{n})\}_{n\in\mathbf{N}}$
be $E$-valued processes and $\mathbf{S}\subset\mathbf{T}\subset\mathbf{R}^{+}$
with $\mathbf{S}$ being dense. Suppose that:

\renewcommand{\labelenumi}{(\roman{enumi})}
\begin{enumerate}
\item $C_{b}(E;\mathbf{R})$ separates points on $E$.
\item $\mathcal{D}\subset C_{b}(E;\mathbf{R})$ is countable and $E_{0}$
is a $\mathcal{D}$-baseable subset of $E$.
\item $\{X^{n}\}_{n\in\mathbf{N}}$ satisfies (\ref{eq:Common_T-Base_N}).
\item $\{X^{n}\}_{n\in\mathbf{N}}$ satisfies $\mathbf{S}$-PSMTC in $E_{0}$
and $\mathcal{D}$-FMCC.
\item $\mathfrak{flp}_{\mathbf{S}}(\{\varpi(\varphi)\circ X^{n}\}_{n\in\mathbf{N}})$
has at least one c$\grave{a}$dl$\grave{\mbox{a}}$g member with $\varphi\circeq\bigotimes\mathcal{D}$.
\end{enumerate}
Then, there exist a stochastic basis%
\footnote{The notion of stochastic basis was reviewed in \S \ref{sec:Proc}.%
} $(\Omega,\mathscr{F},\{\mathscr{G}_{t}\}_{t\geq0},\mathbb{P})$,
some $\{n_{k}\}_{k\in\mathbf{N}}\subset\mathbf{N}$ and an $X\in(E_{0}^{\mathbf{R}^{+}})^{\Omega}$
such that:

\renewcommand{\labelenumi}{(\alph{enumi})}
\begin{enumerate}
\item $X=\mathfrak{fl}_{\mathbf{S}}(\{X^{n_{k}}\}_{k\in\mathbf{N}})$ is
an $E$-valued $(\mathbf{S},\mathcal{D})$-c$\grave{\mbox{a}}$dl$\grave{\mbox{a}}$g
process and satisfies $\mathbf{S}$-PMTC in $E_{0}$.
\item If $E_{0}\in\mathscr{B}^{\mathbf{s}}(E)$, then $X=\mathfrak{fl}_{\mathbf{S}}(\{X^{n_{k}}\}_{k\in\mathbf{N}})$
can be chosen to be an $E$-valued, $(\mathbf{S},\mathcal{D})$-c$\grave{\mbox{a}}$dl$\grave{\mbox{a}}$g,
$\mathscr{G}_{t}$-progressive%
\footnote{The notion of $\mathscr{G}_{t}$-progressive processes was specified
in \S \ref{sec:Proc}.%
} process that satisfies $\mathbf{S}$-PMTC in $E_{0}$ and admits
an $(\mathbf{S},\mathcal{D})$-c$\grave{\mbox{a}}$dl$\grave{\mbox{a}}$g
progressive modification with paths in $E_{0}^{\mathbf{R}^{+}}$%
\footnote{Given $E_{0}\subset E$, an $E$-valued process with paths in $E_{0}^{\mathbf{R}^{+}}$
is equivalent to an $(E_{0},\mathscr{O}_{E}(E_{0}))$-valued process.%
}.
\item If $\mathbf{T}$ is conull and $\{X^{n_{k}}\}_{k\in\mathbf{N}}$ is
$(\mathbf{T},\mathcal{D})$-AS, then $X=\mathfrak{fl}_{\mathbf{S}}(\{X^{n_{k}}\}_{k\in\mathbf{N}})$
can be chosen to be an $E$-valued, stationary, $(\mathbf{R}^{+},\mathcal{D})$-c$\grave{\mbox{a}}$dl$\grave{\mbox{a}}$g
process that satisfies $\mathbf{R}^{+}$-PMTC in $E_{0}$ and admits
an $(\mathbf{R}^{+},\mathcal{D})$-c$\grave{\mbox{a}}$dl$\grave{\mbox{a}}$g
progressive modification with paths in $E_{0}^{\mathbf{R}^{+}}$.
\end{enumerate}
\end{thm}
\begin{rem}
\label{rem:E0_Countable_D_Baseable}If $E_{0}$ is a $\mathcal{D}$-baseable
subset with $\mathcal{D}\subset C(E;\mathbf{R})$, then $E_{0}$ is
$\mathcal{D}_{0}$-baseable for some countable $\mathcal{D}_{0}\subset\mathcal{D}$
(see Fact \ref{fact:D-Baseable} (c)) and $\mathcal{D}_{0}$-FMCC
is a weaker assumption than $\mathcal{D}$-FMCC. Hence, it is no less
general to make $\mathcal{D}$ a countable collection in the theorem
above.
\end{rem}

\begin{rem}
\label{rem:PSMTC_Baseable}Any compact subset contained in a baseable
set $E_{0}$ is metrizable by Corollary \ref{cor:Baseable_MC} (a).
Thus, the $\mathbf{m}$-tightness within $\mathbf{S}$-PSMTC in $E_{0}$
is reduced to ordinary tightness.
\end{rem}

\begin{rem}
\label{rem:E0_Valued_Prog_FLP}The proof of Theorem \ref{thm:FL_Weakly_Cadlag}
relies on Theorem \ref{thm:TransFC_1} in which the limit processs
$X$ was defined by a collection of $E_{0}$-valued mappings $\{X_{t}\}_{t\geq0}$.
This is equivalent to describing the limit process as an $E_{0}^{\mathbf{R}^{+}}$-valued
mapping (like $X$ and $X^{\prime}$ in Theorem \ref{thm:FL_Weakly_Cadlag}).
Moreover, both Theorem \ref{thm:TransFC_1} and Theorem \ref{thm:FL_Weakly_Cadlag}
consider the limit processes as $E$-valued with paths in $E_{0}^{\mathbf{R}^{+}}$
for finite-dimensional convergence.
\end{rem}

We use the next lemma to establish progressiveness in Theorem \ref{thm:FL_Weakly_Cadlag}
(b, c).
\begin{lem}
\label{lem:Prog_(T,E0)-Mod}Let $E$ be a topological space, $x_{0}\in E_{0}\subset E$,
$\mathbf{T}\subset\mathbf{R}^{+}$ and $(\Omega,\mathscr{F},\mathbb{P};Y)$
be an $E$-valued process satisfying
\begin{equation}
\inf_{t\in\mathbf{T}}\mathbb{P}(Y_{t}\in E_{0})=1.\label{eq:T-Base_Y}
\end{equation}
Then, the mapping
\begin{equation}
X\circeq\bigotimes_{t\in\mathbf{R}^{+}}\mathfrak{var}\left(Y_{t};\Omega,Y_{t}^{-1}(E_{0}),x_{0}\right)\in\left(E_{0}^{\mathbf{R}^{+}}\right)^{\Omega}\label{eq:Prog_Mod}
\end{equation}
satisfies the following statements:

\renewcommand{\labelenumi}{(\alph{enumi})}
\begin{enumerate}
\item $X_{t}\circeq\mathfrak{p}_{t}\circ X\in M(\Omega,\mathscr{F};E_{0},\mathscr{B}_{E}(E_{0}))$
for all $t\in\mathbf{T}$ and (\ref{eq:Lim_Proc_(T,E0)-Mod}) holds.
\item If $\mathbf{T}=\mathbf{R}^{+}$, then $X$ is an $E$-valued process
with paths in $E_{0}^{\mathbf{R}^{+}}$, satisfies (\ref{eq:Lim_Proc_(R+,E0)_Mod})
and is a modification of $Y$.
\item If $E_{0}\in\mathscr{B}(E)$, then $X$ is an $(E_{0},\mathscr{O}_{E}(E_{0}))$-valued
$\mathscr{F}_{t}^{Y}$-adapted process. If, in addition, $Y$ is a
measurable or progressive process, then $X$ is measurable or $\mathscr{F}_{t}^{Y}$-progressive
respectively.
\item If $E_{0}\in\mathscr{B}(E)\cap\mathscr{B}^{\mathbf{s}}(E)$, and if
$Y$ is a measurable process, then $X$ is an $(E_{0},\mathscr{O}_{E}(E_{0}))$-valued,
measurable, $\mathscr{F}_{t}^{Y}$-adapted process and admits an $\mathscr{F}_{t}^{Y}$-progressive
modification.
\end{enumerate}
\end{lem}
\begin{proof}
(a) follows by Lemma \ref{lem:var(X)} (b, c) (with $\mathscr{U}=\mathscr{B}(E)$,
$S=S_{0}=E_{0}$, $\mathscr{U}^{\prime}=\mathscr{B}_{E}(E_{0})$,
$X=Y_{t}$, $Y=X_{t}$ and $y_{0}=x_{0}$).

(b) follows by (a) (with $\mathbf{T}=\mathbf{R}^{+}$) and Fact \ref{fact:Proc_Basic_1}
(b) (with $E=(E_{0},\mathscr{O}_{E}(E_{0}))$).

(c) Let $\varphi$ denote the identity mapping on $E$. We find that
\begin{equation}
\varphi^{\prime}\circeq\mathfrak{var}(\varphi;E,E_{0},x_{0})\in M\left(E;E_{0},\mathscr{O}_{E}(E_{0})\right)\label{eq:Prog_Mod_Map}
\end{equation}
by $E_{0}\in\mathscr{B}(E)$ and Fact \ref{fact:var(f)} (b) (with
$(S,\mathscr{A})=(E,\mathscr{B}(E))$, $(E,\mathscr{U})=(E_{0},\mathscr{B}_{E}(E_{0}))$,
$A=E_{0}$, $f=\varphi$ and $y_{0}=x_{0}$). Then, (c) follows by
(\ref{eq:Prog_Mod_Map}), the fact $X=\varpi(\varphi^{\prime})\circ Y$
and Fact \ref{fact:Proc_Path_Mapping} (a) (with $S=(E_{0},\mathscr{O}_{E}(E_{0}))$,
$f=\varphi^{\prime}$, $X=Y$ and $\mathscr{G}_{t}=\mathscr{F}_{t}^{Y}$).

(d) $X$ is an $(E_{0},\mathscr{O}_{E}(E_{0}))$-valued, measurable,
$\mathscr{F}_{t}^{Y}$-adapted process by (c). Let $\varphi_{0}$
denote the identity mapping on $E_{0}$. By $E_{0}\in\mathscr{B}^{\mathbf{s}}(E)$
and Proposition \ref{prop:SB} (a, d), there exists a topology $\mathscr{U}$
on $E_{0}$ such that $(E_{0},\mathscr{U})$ is a Polish space and
$\varphi_{0}\in\mathbf{biso}(E_{0},\mathscr{O}_{E}(E_{0});E_{0},\mathscr{U})$
. Then, $(E_{0},\mathscr{O}_{E}(E_{0}))$-valued measurable (resp.
$\mathscr{F}_{t}^{Y}$-progressive) processes are equivalent to $(E_{0},\mathscr{U})$-valued
measurable (resp. $\mathscr{F}_{t}^{Y}$-progressive) processes by
Fact \ref{fact:Proc_Path_Mapping} (a) (with $E$ (or $S$) being
$(E_{0},\mathscr{O}_{E}(E_{0}))$, $S$ (or $E$) being $(E_{0},\mathscr{U})$,
$f=\varphi_{0}$ and $\mathscr{G}_{t}=\mathscr{F}_{t}^{Y}$). \cite[Theorem 0.1]{OS13}
established that every Polish-space-valued, measurable, $\mathscr{F}_{t}^{Y}$-adapted
process (like $X$) admits an $\mathscr{F}_{t}^{Y}$-progressive modification.
Thus (d) follows immediately.\end{proof}

\begin{cor}
\label{cor:Prog_Mod_Measurable}Let $E$ be a topological space, $E_{0}\in\mathscr{B}(E)\cap\mathscr{B}^{\mathbf{s}}(E)$
and $(\Omega,\mathscr{F},\mathbb{P};X)$ be an $E$-valued measurable
process satisfying (\ref{eq:R+-Base}). Then, $X$ has a progressive
modification with paths in $E_{0}^{\mathbf{R}^{+}}$.
\end{cor}
\begin{proof}
This corollary follows by Lemma \ref{lem:Prog_(T,E0)-Mod} (b, d)
(with $Y=X$ and $X=Y$) and Proposition \ref{prop:Proc_Basic_2}
(e).\end{proof}

\begin{cor}
\label{cor:Prog_Mod_Weakly_Cadlag}Let $E$ be a topological space,
$\mathcal{D}\subset C_{b}(E;\mathbf{R})$, $E_{0}$ be a $\mathcal{D}$-baseable
standard Borel subset of $E$ and $(\Omega,\mathscr{F},\mathbb{P};X)$
be an $E$-valued $(\mathbf{R}^{+},\mathcal{D})$-c$\grave{\mbox{a}}$dl$\grave{\mbox{a}}$g
process satisfying (\ref{eq:R+-Base}). Then, $X$ has a progressive
modification with paths in $E_{0}^{\mathbf{R}^{+}}$.
\end{cor}
\begin{proof}
There exists a base $(E_{0},\mathcal{F};\widehat{E},\widehat{\mathcal{F}})$
over $E$ with $\mathcal{F}\subset(\mathcal{D}\cup\{1\})$ by Lemma
\ref{lem:Base_Construction} (c). $\widehat{X}=\mathfrak{rep}_{\mathrm{c}}(X;E_{0},\mathcal{F})$
exists by the fact $(\mathcal{F}\backslash\{1\})\subset\mathcal{D}$
and Proposition \ref{prop:FR} (a). It follows by Fact \ref{fact:Cadlag_RepProc}
and Proposition \ref{prop:RepProc_T-Base} (a) (with $\mathbf{T}=\mathbf{R}^{+}$)
that $X$ and $\widehat{X}$ satisfy (\ref{eq:RepProc_(R+,E0)-Mod})
and $\widehat{X}$ is a progressive process. $\mathscr{F}^{X}=\mathscr{F}^{\widehat{X}}$%
\footnote{The notation ``$\mathscr{F}^{X}$'' as defined in \S \ref{sec:Proc}
means the augmented natural filtration of $X$.%
} by Lemma \ref{lem:Proc_Rep} (e) (with $A=E_{0}$ and $Y=\widehat{X}$),
so $\widehat{X}$ is $\mathscr{F}_{t}^{X}$-progressive. Furthermore,
we have
\begin{equation}
\mathscr{B}_{E}(E_{0})=\mathscr{B}_{\widehat{E}}(E_{0})\subset\mathscr{B}(\widehat{E})\label{eq:E0_SB_Borel_Equal}
\end{equation}
by Lemma \ref{lem:SB_Base} (a) (with $d=1$ and $A=E_{0}$). $\widehat{X}$
has an $\mathscr{F}_{t}^{X}$-progressive modification $Z$ with paths
in $E_{0}^{\mathbf{R}^{+}}$ by (\ref{eq:RepProc_(R+,E0)-Mod}), (\ref{eq:E0_SB_Borel_Equal})
and Lemma \ref{lem:Prog_(T,E0)-Mod} (b, c) (with $E=\widehat{E}$
and $Y=\widehat{X}$). $Z$ is an $(E_{0},\mathscr{O}_{E}(E_{0}))$-valued
process by (\ref{eq:E0_SB_Borel_Equal}). $Z$ is a modification of
$X$ by (\ref{eq:RepProc_(R+,E0)-Mod}). So, $Z$ is progressive by
Proposition \ref{prop:Proc_Basic_2} (e).\end{proof}

\begin{rem}
\label{rem:Prog_Mod_TF_Baseable_Space}A special case of Corollary
\ref{cor:Prog_Mod_Measurable} and Corollary \ref{cor:Prog_Mod_Weakly_Cadlag}
is when $E=E_{0}$ is a $\mathcal{D}$-baseable standard Borel space
and (\ref{eq:R+-Base}) becomes automatic.
\end{rem}

The next result shows that the condition (v) of Theorem \ref{thm:FL_Weakly_Cadlag}
is realizable.
\begin{prop}
\label{prop:FLP_Under_TF_FMCC}Let $E$ be a topological space, $\mathcal{D}\subset C_{b}(E;\mathbf{R})$
be countable, $E_{0}$ be a $\mathcal{D}$-baseable subset of $E$
and $\mathbf{I}$ be an infinite index set. If $E$-valued processes
$\{(\Omega^{i},\mathscr{F}^{i},\mathbb{P}^{i};X^{i})\}_{i\in\mathbf{I}}$
satisfy (\ref{eq:Common_FR-Base_I}) and $\mathcal{D}$-FMCC, then
there exists a countable $\mathbf{J}\subset(0,\infty)$ such that
$\mathfrak{flp}_{\mathbf{R}^{+}\backslash\mathbf{J}}(\{\varpi(\mathcal{D})\circ X^{i}\}_{i\in\mathbf{I}})$%
\footnote{Please be noted that members of $\mathfrak{flp}_{\mathbf{R}^{+}\backslash\mathbf{J}}(\{\varpi(\mathcal{D})\circ X^{i}\}_{i\in\mathbf{I}})$
are processes with time horizon $\mathbf{R}^{+}$ no matter $\mathbf{J}=\varnothing$
or not.%
} has at least one c$\grave{\mbox{a}}$dl$\grave{\mbox{a}}$g member.
\end{prop}
\begin{proof}
There exists a base $(E_{0},\mathcal{F};\widehat{E},\widehat{\mathcal{F}})$
over $E$ with $\mathcal{F}=\mathcal{D}\cup\{1\}$ by Lemma \ref{lem:Base_Construction}
(b), so $\{X^{n}\}_{n\in\mathbf{N}}$ satisfies $\mathcal{F}$-FMCC.
It follows by Proposition \ref{prop:FR_Tight} (a) and Note \ref{note:Ehat_Valued_Proc_FDD}
that $\{\widehat{X}^{i}=\mathfrak{rep}_{\mathrm{c}}(X^{i};E_{0},\mathcal{F})\}_{i\in\mathbf{I}}$
is tight in the Polish space $D(\mathbf{R}^{+};\widehat{E})$. $\{\widehat{X}^{i}\}_{i\in\mathbf{I}}$
admits at least one weak limit point $Y$ on $D(\mathbf{R}^{+};\widehat{E})$
by the Prokhorov's Theorem (Theorem \ref{thm:Prokhorov} (b)). $\mathbf{J}\circeq J(Y)\subset(0,\infty)$
is countable by Note \ref{note:J(Y)}. Now, the result follows by
Proposition \ref{prop:FR_FC_WC} (a) and Lemma \ref{lem:FLP_Under_TF}.\end{proof}

We now prove the main theorem of this section.

\begin{proof}
[Proof of Theorem \ref{thm:FL_Weakly_Cadlag}]The proof is divided
into six steps.

\textit{Step 1: Establish a base $(E_{0},\mathcal{F};\widehat{E},\widehat{\mathcal{F}})$
and c$\grave{\mbox{a}}$dl$\grave{\mbox{a}}$g replicas $\{\widehat{X}^{n}\}_{n\in\mathbf{N}}$}.
There exists a base $(E_{0},\mathcal{F};\widehat{E},\widehat{\mathcal{F}})$
over $E$ with $\mathcal{F}=\mathcal{D}\cup\{1\}$ by the condition
(ii) above and Lemma \ref{lem:Base_Construction} (b). (\ref{eq:Common_FR-Base_N})
holds by the condition (iii) above and Fact \ref{fact:Common_T-Base_FR-Base}
(with $\mathbf{I}=\mathbf{N}$). $\{X^{n}\}_{n\in\mathbf{N}}$ satisfies
$\mathcal{F}$-FMCC by the condition (iv) above and the fact $\mathcal{F}\backslash\{1\}=\mathcal{D}$,
so they are $(\mathbf{R}^{+},\mathcal{F})$-c$\grave{\mbox{a}}$dl$\grave{\mbox{a}}$g.
It then follows by Proposition \ref{prop:FR_Tight} (a) (with $\mathbf{I}=\mathbf{N}$),
Proposition \ref{prop:FR} (a) (with $X=X^{n}$) and Note \ref{note:Ehat_Valued_Proc_FDD}
that $\{\widehat{X}^{n}=\mathfrak{rep}_{\mathrm{c}}(X^{n};E_{0},\mathcal{F})\}_{n\in\mathbf{N}}$
is tight in the Polish space $D(\mathbf{R}^{+};\widehat{E})$ and
satisfies
\begin{equation}
\inf_{n\in\mathbf{N}}\mathbb{P}^{n}\left(\varphi\circ X_{t}^{n}=\widehat{\varphi}\circ\widehat{X}_{t}^{n}\right)=1\label{eq:Phi(X)_Phihat(Xhat)_R+-Mod_N}
\end{equation}
with $\widehat{\varphi}\circeq\bigotimes\widehat{\mathcal{F}}\backslash\{1\}$.

\textit{Step 2: Establish $\{n_{k}\}_{k\in\mathbf{N}}$ and a $D(\mathbf{R}^{+};\widehat{E})$-valued
random variable $Y$ such that}
\begin{equation}
\widehat{X}^{n_{k}}\xrightarrow{\quad\mathrm{D}(\mathbf{S})\quad}Y\mbox{ as }k\uparrow\infty.\label{eq:RepProc_Subseq_FC_along_S_Y}
\end{equation}
By the condition (v) above, the tightness of $\{\widehat{X}^{n}\}_{n\in\mathbf{N}}$
in $D(\mathbf{R}^{+};\widehat{E})$ and Prokhorov's Theorem (Theorem
\ref{thm:Prokhorov} (b)), there exist $\{n_{k}\}_{k\in\mathbf{N}}\subset\mathbf{N}$,
a $D(\mathbf{R}^{+};\widehat{E})$-valued random variable $Y$ (see
Remark \ref{rem:Canonical_Proc} below for an explicit construction)
and an $\mathbf{R}^{\mathcal{D}}$-valued c$\grave{\mbox{a}}$dl$\grave{\mbox{a}}$g
process $Z$ such that
\begin{equation}
\widehat{X}^{n_{k}}\Longrightarrow Y\mbox{ as }k\uparrow\infty\mbox{ on }D(\mathbf{R}^{+};\widehat{E})\label{eq:RepProc_Subseq_WC_Y}
\end{equation}
and
\begin{equation}
\varpi(\varphi)\circ X^{n_{k}}\xrightarrow{\quad\mathrm{D}(\mathbf{S})\quad}Z\mbox{ as }k\uparrow\infty.\label{eq:phi(Xn)_Subseq_FC_S}
\end{equation}
Without loss of generality, we suppose $Y$ and $Z$ are both defined
on $(\Omega,\mathscr{F},\mathbb{P})$.

Since $\mathcal{D}$ is countable, $\mathbf{R}^{\mathcal{D}}$ and
$D(\mathbf{R}^{+};\mathbf{R}^{\mathcal{D}})$ are Polish spaces as
mentioned in Note \ref{note:Ehat_Valued_Proc_FDD}. So, $Z$ can be
considered as a $D(\mathbf{R}^{+};\mathbf{R}^{\mathcal{D}})$-valued
random variable by Fact \ref{fact:Cadlag_Sko_RV} (b) (with $E=\mathbf{R}^{\mathcal{D}}$).
$\widehat{\mathcal{F}}\backslash\{1\}$ separates points on $\widehat{E}$
by Definition \ref{def:Base}. $\widehat{\mathcal{F}}\backslash\{1\}$
strongly separates points on $\widehat{E}$ and
\begin{equation}
\widehat{\varphi}\in\mathbf{imb}(\widehat{E};\mathbf{R}^{\mathcal{D}})\label{eq:phihat_Imb}
\end{equation}
by Lemma \ref{lem:Base} (a).
\begin{equation}
\varpi(\widehat{\varphi})\in C\left(D(\mathbf{R}^{+};\widehat{E});D(\mathbf{R}^{+};\mathbf{R})\right)\label{eq:phihat_PathMapping_Cont}
\end{equation}
by Proposition \ref{prop:Sko_Basic_1} (d) (with $S=\widehat{E}$,
$E=\mathbf{R}^{\mathcal{D}}$ and $f=\widehat{\varphi}$).

It follows by (\ref{eq:RepProc_Subseq_WC_Y}), (\ref{eq:phihat_PathMapping_Cont})
and Continuous Mapping Theorem (Theorem \ref{thm:ContMapTh} (a))
that
\begin{equation}
\varpi(\widehat{\varphi})\circ\widehat{X}^{n_{k}}\Longrightarrow\varpi(\widehat{\varphi})\circ Y\mbox{ as }k\uparrow\infty\mbox{ on }D\left(\mathbf{R}^{+};\mathbf{R}^{\mathcal{D}}\right).\label{eq:Phihat(Xhatn)_Subseq_WC_1}
\end{equation}
It follows by (\ref{eq:Phi(X)_Phihat(Xhat)_R+-Mod_N}) and (\ref{eq:phi(Xn)_Subseq_FC_S})
that
\begin{equation}
\varpi(\widehat{\varphi})\circ\widehat{X}^{n_{k}}\xrightarrow{\quad\mathrm{D}(\mathbf{S})\quad}Z\mbox{ as }k\uparrow\infty.\label{eq:Phihat(Xhatn)_Subseq_FC_S_1}
\end{equation}
$\{\varpi(\widehat{\varphi})\circ\widehat{X}^{n_{k}}\}_{k\in\mathbf{N}}$
is tight in $D(\mathbf{R}^{+};\mathbf{R}^{\mathcal{D}})$ by (\ref{eq:phihat_PathMapping_Cont})
and Fact \ref{fact:Push_Forward_Tight_2} (with $E=D(\mathbf{R}^{+};\widehat{E})$,
$S=D(\mathbf{R}^{+};\mathbf{R}^{\mathcal{D}})$ and $f=\varpi(\widehat{\varphi})$).
Hence, (\ref{eq:Phihat(Xhatn)_Subseq_FC_S_1}) implies
\begin{equation}
\varpi(\widehat{\varphi})\circ\widehat{X}^{n_{k}}\Longrightarrow Z\mbox{ as }k\uparrow\infty\mbox{ on }D\left(\mathbf{R}^{+};\mathbf{R}^{\mathcal{D}}\right)\label{eq:Phihat(Xhatn)_Subseq_WC_2}
\end{equation}
by the Prokhorov's Theorem (Theorem \ref{thm:Prokhorov} (b)), the
denseness of $\mathbf{S}$ and Theorem \ref{thm:Sko_RV_WC_FC_Metrizable_Separable}
(b) (with $E=\mathbf{R}^{\mathcal{D}}$, $X^{n}=\varpi(\widehat{\varphi})\circ\widehat{X}^{n_{k}}$,
$X=Z$ and $\mathbf{T}=\mathbf{S}$). $\mathcal{P}(D(\mathbf{R}^{+};\mathbf{R}^{\mathcal{D}}))$
is a Polish space by Theorem \ref{thm:P(E)_Compact_Polish} (b) (with
$E=D(\mathbf{R}^{+};\mathbf{R}^{\mathcal{D}})$) in which a sequence
converges weakly to at most one point by \cite[Theorem 17.10]{M00},
Proposition \ref{prop:Var_Polish} (c) and Proposition \ref{prop:Metrizable}
(a). So,
\begin{equation}
\mathbb{P}\circ Z^{-1}=\mathbb{P}\circ\left(\varpi(\widehat{\varphi})\circ Y\right)^{-1}\in\mathcal{P}\left(D(\mathbf{R}^{+};\mathbf{R}^{\mathcal{D}})\right)\label{eq:Y_Z_Same_Sko_Dist}
\end{equation}
by (\ref{eq:Phihat(Xhatn)_Subseq_WC_1}) and (\ref{eq:Phihat(Xhatn)_Subseq_WC_2}).
$\varpi(\widehat{\varphi})\circ Y$ and $Z$ as $\mathbf{R}^{\mathcal{D}}$-valued
processes have the same finite-dimensional distributions by (\ref{eq:Y_Z_Same_Sko_Dist})
and Fact \ref{fact:Sko_RV_Cadlag} (b) (with $E=\mathbf{R}^{\mathcal{D}}$),
so (\ref{eq:Phihat(Xhatn)_Subseq_FC_S_1}) implies
\begin{equation}
\varpi(\widehat{\varphi})\circ\widehat{X}^{n_{k}}\xrightarrow{\quad\mathrm{D}(\mathbf{S})\quad}\varpi(\widehat{\varphi})\circ Y\mbox{ as }k\uparrow\infty.\label{eq:phihat(Xhatn)_Subseq_FC_S_2}
\end{equation}

We fix $\mathbf{T}_{0}\in\mathscr{P}_{0}(\mathbf{S})$ and put $d=\aleph(\mathbf{T}_{0})$.
$\{\widehat{X}^{n_{k}}\}_{k\in\mathbf{N}}$, $\{\varpi(\widehat{\varphi})\circ\widehat{X}^{n_{k}}\}_{k\in\mathbf{N}}$,
$Y$ and $\varpi(\widehat{\varphi})\circ Y$ all have Borel finite-dimensional
distributions as mentioned in Note \ref{note:Ehat_Valued_Proc_FDD},
so (\ref{eq:phihat(Xhatn)_Subseq_FC_S_2}) implies
\begin{equation}
\left(\bigotimes_{t\in\mathbf{T}_{0}}\widehat{\varphi}\circ\mathfrak{p}_{t}\right)\circ\widehat{X}^{n_{k}}\Longrightarrow\left(\bigotimes_{t\in\mathbf{T}_{0}}\widehat{\varphi}\circ\mathfrak{p}_{t}\right)\circ Y\mbox{ as }k\uparrow\infty\mbox{ on }\mathbf{R}^{d}.\label{eq:phihat(Xhatn)_Subseq_FDD_WC}
\end{equation}
One finds that
\begin{equation}
\Psi\circeq\bigotimes_{t\in\mathbf{T}_{0}}\widehat{\varphi}^{-1}\circ\mathfrak{p}_{t}\in C\left[\widehat{\varphi}(\widehat{E})^{d},\mathscr{O}_{\mathbf{R}^{d}}\left(\widehat{\varphi}(\widehat{E})^{d}\right);\widehat{E}^{d}\right]\label{eq:FDD_WC_Inverse_Map}
\end{equation}
by (\ref{eq:phihat_Imb}) and Fact \ref{fact:Prod_Map_2} (a, b).
Hence, it follows by (\ref{eq:FDD_WC_Inverse_Map}), (\ref{eq:phihat(Xhatn)_Subseq_FDD_WC})
and Continuous Mapping Theorem (Theorem \ref{thm:ContMapTh} (a))
that
\begin{equation}
\begin{aligned}\widehat{X}_{\mathbf{T}_{0}}^{n_{k}} & =\Psi\circ\left(\bigotimes_{t\in\mathbf{T}_{0}}\widehat{\varphi}\circ\mathfrak{p}_{t}\right)\circ\widehat{X}^{n_{k}}\\
 & \Longrightarrow\Psi\circ\left(\bigotimes_{t\in\mathbf{T}_{0}}\widehat{\varphi}\circ\mathfrak{p}_{t}\right)\circ Y=Y_{\mathbf{T}_{0}}\mbox{ as }k\uparrow\infty\mbox{ on }\widehat{E}^{\mathbf{T}_{0}}.
\end{aligned}
\label{eq:Check_RepProc_Subseq_FC_S}
\end{equation}

\textit{Step 3: Construct $(\Omega,\mathscr{F},\{\mathscr{G}_{t}\}_{t\geq0},\mathbb{P})$
and }$X$. We fix an arbitrary $x_{0}\in E_{0}$ and set $\mathscr{G}_{t}=\mathscr{F}_{t}^{Y}$.
For (a), we define $X$ by (\ref{eq:Trans_FC_Lim_Proc}) with $\mathbf{T}=\mathbf{S}$.
For (b, c), we shall need an $X$ with different paths to meet the
measurability requirements. For this purpose, we define the process
$X^{\prime}$ by (\ref{eq:Prog_Mod}) with $X$ replaced by $X^{\prime}$.

\textit{Step 4: Verify (a)}. It follows by the conditions (i, iv)
above, (\ref{eq:RepProc_Subseq_FC_along_S_Y}), Lemma \ref{lem:RepProc_Int_Test}
(c, e) (with $n=n_{k}$) and Theorem \ref{thm:TransFC_1} (a, c) (with
$n=n_{k}$ and $\mathbf{T}=\mathbf{S}$) that the $X$ defined in
Step 3 is an $E$-valued process with paths in $E_{0}^{\mathbf{R}^{+}}$
that satisfies: (1) $\mathbf{S}$-PMTC in $E_{0}$, (2) 
\begin{equation}
\inf_{t\in\mathbf{S}}\mathbb{P}\left(X_{t}=Y_{t}\in E_{0}\right)=1,\label{eq:Lim_Proc_(S,E0)-Mod}
\end{equation}
and (3) $X=\mathfrak{fl}_{\mathbf{S}}(\{X^{n_{k}}\}_{k\in\mathbf{N}})$.
Hence, (a) follows by the fact $\mathcal{D}\subset\mathcal{F}$, (\ref{eq:Lim_Proc_(S,E0)-Mod})
and Lemma \ref{lem:Proc_Rep} (b) (with $\mathbf{T}=\mathbf{S}$).

\textit{Step 5: Verify (b)}. $Y$ is c$\grave{\mbox{a}}$dl$\grave{\mbox{a}}$g
hence $\mathscr{G}_{t}$-progressive by Proposition \ref{prop:Proc_Basic_2}
(a). Given $E_{0}\in\mathscr{B}^{\mathbf{s}}(E)$, (\ref{eq:E0_SB_Borel_Equal})
holds by Lemma \ref{lem:SB_Base} (a) (with $d=1$ and $A=E_{0}$).
Hence, the $X^{\prime}$ defined in Step 3 is an $(E_{0},\mathscr{O}_{E}(E_{0}))$-valued
$\mathscr{G}_{t}$-progressive process satisfying
\begin{equation}
\inf_{t\in\mathbf{S}}\mathbb{P}\left(X_{t}=Y_{t}=X_{t}^{\prime}\in E_{0}\right)=1\label{eq:Lim_Proc_Prog_(S,E0)-Mod}
\end{equation}
by (\ref{eq:Lim_Proc_(S,E0)-Mod}), (\ref{eq:E0_SB_Borel_Equal})
and Lemma \ref{lem:Prog_(T,E0)-Mod} (a, c) (with $E=\widehat{E}$
and $\mathbf{T}=\mathbf{S}$ and $X=X^{\prime}$). Then, $X^{\prime}=\mathfrak{fl}_{\mathbf{S}}(\{X^{n_{k}}\}_{k\in\mathbf{N}})$
by (a) and (\ref{eq:Lim_Proc_Prog_(S,E0)-Mod}). $X^{\prime}$ is
$(\mathbf{S},\mathcal{D})$-c$\grave{\mbox{a}}$dl$\grave{\mbox{a}}$g
by the fact $\mathcal{D}\subset\mathcal{F}$, (\ref{eq:Lim_Proc_Prog_(S,E0)-Mod})
and Lemma \ref{lem:Proc_Rep} (b) (with $\mathbf{T}=\mathbf{S}$ and
$X=X^{\prime}$). $X^{\prime}$ is a measurable process by Proposition
\ref{prop:Proc_Basic_2} (c). Now, (b) (with $X=X^{\prime}$) follows
by the fact $E_{0}\in\mathscr{B}(E)\cap\mathscr{B}^{\mathbf{s}}(E)$,
Corollary \ref{cor:Prog_Mod_Measurable} (with $X=X^{\prime}$) and
Note \ref{note:Weakly_Cadlag}.

\textit{Step 6: Verify (c)}. $\{X^{n}\}_{n\in\mathbf{N}}$ is $(\mathbf{T},\mathcal{F}\backslash\{1\})$-AS
since $\mathcal{D}=\mathcal{F}\backslash\{1\}$. Then, $Y$ is a stationary
process by (\ref{eq:RepProc_Subseq_WC_Y}) and Proposition \ref{prop:FR_FC_WC}
(c) (with $n=n_{k}$). We know from (a) that $X$ satisfies $\mathbf{S}$-PMTC
in $E_{0}$ and (\ref{eq:Lim_Proc_(S,E0)-Mod}). By $\mathbf{S}$-PMTC
in $E_{0}$ there exists an $A\in\mathscr{K}_{\sigma}^{\mathbf{m}}(E_{0},\mathscr{O}_{E}(E_{0}))$
and some $t_{0}\in\mathbf{S}$ such that $\mathbb{P}(X_{t_{0}}\in A)=1$.
(\ref{eq:Lim_Proc_(S,E0)-Mod}) implies $\mathbb{P}(Y_{t_{0}}\in A)=1$.
The stationarity of $Y$ implies
\begin{equation}
\inf_{t\in\mathbf{R}^{+}}\mathbb{P}\left(Y_{t}\in A\subset E_{0}\right)=1.\label{eq:Rep_Lim_Proc_R+-Base_A_E0}
\end{equation}
$A$ is a $\mathcal{D}$-baseable standard Borel subset of $E$ and
satisfies
\begin{equation}
\mathscr{B}_{E}(A)=\mathscr{B}_{\widehat{E}}(A)\subset\mathscr{B}(\widehat{E})\label{eq:A_SB_Borel_Equal}
\end{equation}
by Corollary \ref{cor:Base_Compact} (b) (with $d=1$), Lemma \ref{lem:SB_Base}
(a) (with $d=1$), the fact $\mathcal{F}\backslash\{1\}=\mathcal{D}$
and Fact \ref{fact:D-Baseable} (a, b). Hence, the $X^{\prime}$ defined
in Step 3 is an $(E_{0},\mathscr{O}_{\widehat{E}}(E_{0}))$-valued
process satisfying both (\ref{eq:Lim_Proc_Prog_(S,E0)-Mod}) and
\begin{equation}
\inf_{t\in\mathbf{R}^{+}}\mathbb{P}\left(X_{t}^{\prime}=Y_{t}\in A\subset E_{0}\right)=1\label{eq:Lim_Proc_(R+,A)-Mod_X^Prime}
\end{equation}
by (\ref{eq:Rep_Lim_Proc_R+-Base_A_E0}), (\ref{eq:A_SB_Borel_Equal}),
Lemma \ref{lem:Prog_(T,E0)-Mod} (b) (with $E=\widehat{E}$, $E_{0}=A$
and $X=X^{\prime}$) and (\ref{eq:Lim_Proc_(S,E0)-Mod}). $X^{\prime}$
is an $(E_{0},\mathscr{O}_{E}(E_{0}))$-valued process and satisfies
$\mathbf{R}^{+}$-PMTC in $E_{0}$ by (\ref{eq:Lim_Proc_(R+,A)-Mod_X^Prime}),
(\ref{eq:A_SB_Borel_Equal}), Lemma \ref{lem:var(X)} (b, c) (with
$E=\widehat{E}$, $\mathscr{U}=\mathscr{B}(\widehat{E})$, $S_{0}=A$,
$S=E_{0}$, $\mathscr{U}^{\prime}=\mathscr{B}_{E}(E_{0})$, $X=Y_{t}$
and $Y=X_{t}^{\prime}$) and Fact \ref{fact:Proc_Basic_1} (b). Thus,
$X^{\prime}=\mathfrak{fl}_{\mathbf{S}}(\{X^{n_{k}}\}_{k\in\mathbf{N}})$
by (a) and (\ref{eq:Lim_Proc_Prog_(S,E0)-Mod}). $X^{\prime}$ is
stationary by (\ref{eq:Lim_Proc_(R+,A)-Mod_X^Prime}) and Lemma \ref{lem:Proc_Rep}
(e) (with $X=X^{\prime}$). $X^{\prime}$ is $(\mathbf{R}^{+},\mathcal{D})$-c$\grave{\mbox{a}}$dl$\grave{\mbox{a}}$g
by the fact $\mathcal{D}\subset\mathcal{F}$, (\ref{eq:Lim_Proc_(R+,A)-Mod_X^Prime})
and Lemma \ref{lem:Proc_Rep} (b) (with $\mathbf{T}=\mathbf{R}^{+}$
and $X=X^{\prime}$). Finally, (c) (with $X=X^{\prime}$) follows
by Corollary \ref{cor:Prog_Mod_Weakly_Cadlag} (with $E_{0}=A$ and
$X=X^{\prime}$) and Note \ref{note:Weakly_Cadlag}.\end{proof}

\begin{rem}
\label{rem:Canonical_Proc}Let $\{\widehat{X}^{n_{k}}\}_{k\in\mathbf{N}}$
be as in the proof of Theorem \ref{thm:FL_Weakly_Cadlag}. One can
realize (\ref{eq:RepProc_Subseq_WC_Y}) by letting $\Omega=D(\mathbf{R}^{+};\widehat{E})$,
$\mathbb{P}$ be the weak limit of the distributions of $\{\widehat{X}^{n_{k}}\}_{k\in\mathbf{N}}$
in $\mathcal{P}(D(\mathbf{R}^{+};\widehat{E}))$, $\mathscr{F}$ be
the completion%
\footnote{Completion of measure space was specified in \S \ref{sub:Meas}.%
} of $\mathscr{B}(\Omega)$ with respect to $\mathbb{P}$ and $Y$
be the identity mapping on $D(\mathbf{R}^{+};\widehat{E})$. This
process $(\Omega,\mathscr{F},\mathbb{P};Y)$ is often called the coordinate
process or canonical process on $D(\mathbf{R}^{+};\widehat{E})$.
\end{rem}

With the help of Lemma \ref{lem:Prog_(T,E0)-Mod}, we give a variant
of Theorem \ref{thm:FLP_Gen} (b) that can be used to show uniqueness
in the settings of Theorem \ref{thm:FL_Weakly_Cadlag}.
\begin{prop}
\label{prop:FLP_Gen_E0}Let $E$ be a topological space, $\{(\Omega^{i},\mathscr{F}^{i},\mathbb{P}^{i};X^{i})\}_{i\in\mathbf{I}}$
be $E$-valued processes and $\mathbf{S}\subset\mathbf{T}\subset\mathbf{R}^{+}$.
Suppose that:

\renewcommand{\labelenumi}{(\roman{enumi})}
\begin{enumerate}
\item $C_{b}(E;\mathbf{R})$ separates points on $E$.
\item $\mathcal{D}\subset C_{b}(E;\mathbf{R})$ separates points on $E_{0}\in\mathscr{B}(E)$.
\item (\ref{eq:Common_T-Base_I}) holds.
\item $\{X^{i}\}_{i\in\mathbf{I}}$ is $(\mathbf{S},\mathcal{D})$-FDC and
satisfies $\mathbf{S}$-PSMTC in $E_{0}$.
\end{enumerate}
Then, there exists an $X=\mathfrak{flp}_{\mathbf{S}}(\{X^{i}\}_{i\in\mathbf{I}})$
with paths in $E_{0}^{\mathbf{R}^{+}}$ and satisfying $\mathbf{S}$-PMTC
in $E_{0}$. Moreover, $X=\mathfrak{fl}_{\mathbf{S}}(\{X^{i_{n}}\}_{n\in\mathbf{N}})$
for any $\{i_{n}\}_{n\in\mathbf{N}}\subset\mathbf{I}$.
\end{prop}
\begin{proof}
We let $\widetilde{f}\circeq f|_{E_{0}^{d}}$%
\footnote{Similar notations were used in Notation \ref{notation:OP}.%
} for any $f\in C_{b}(E^{d};\mathbf{R})$ and any $d\in\mathbf{N}$,
put $\widetilde{\mathcal{D}}\circeq\mathcal{D}|_{E_{0}}=\{\widetilde{f}:f\in\mathcal{D}\}$,
fix $x_{0}\in E_{0}$ and define $\{Z_{t}^{i}\}_{t\geq0}\subset E_{0}^{\Omega_{i}}$
for each $i\in\mathbf{I}$ by (\ref{eq:Trans_FC_Lim_Proc}) with $X_{t}$,
$Y_{t}$, $\Omega$ replaced by $Z_{t}^{i}$, $X_{t}^{i}$, $\Omega^{i}$
respectively.

It follows by (\ref{eq:Common_T-Base_I}) and Lemma \ref{lem:Prog_(T,E0)-Mod}
(a) (with $(\Omega,\mathscr{F},\mathbb{P};Y)=(\Omega^{i},\mathscr{F}^{i},\mathbb{P}^{i};X^{i})$
and $X_{t}=Z_{t}^{i}$) that
\begin{equation}
Z_{t}^{i}\in M\left(\Omega^{i},\mathscr{F}^{i};E_{0},\mathscr{B}_{E}(E_{0})\right),\;\forall t\in\mathbf{T},i\in\mathbf{I}\label{eq:(T,E0)-MOD_at_T}
\end{equation}
and
\begin{equation}
\inf_{t\in\mathbf{T},i\in\mathbf{I}}\mathbb{P}^{i}(X_{t}^{i}=Z_{t}^{i}\in E_{0})=1.\label{eq:Ori_Proc_(T,E0)-MOD}
\end{equation}
$E$ is a Hausdorff space by Proposition \ref{prop:Fun_Sep_1} (e)
(with $A=E$ and $\mathcal{D}=C_{b}(E;\mathbf{R})$). So, $(E_{0},\mathscr{O}_{E}(E_{0}))$
is a Hausdorff subspace and $\{x_{0}\}\in\mathscr{B}(E)$ by Proposition
\ref{prop:Separability} (a, c) and the fact $E_{0}\in\mathscr{B}(E)$.
This immediately implies
\begin{equation}
Z_{t}^{i}\in M\left(\Omega^{i},\mathscr{F}^{i};E_{0},\mathscr{O}_{E}(E_{0})\right),\;\forall t\in\mathbf{R}^{+}\backslash\mathbf{T},i\in\mathbf{I}.\label{eq:(T,E0)-MOD_Outside_T}
\end{equation}
Hence, $Z^{i}\circeq\{Z_{t}^{i}\}_{t\geq0}$ is an $(E_{0},\mathscr{O}_{E}(E_{0}))$-valued
process for all $i\in\mathbf{I}$ by Fact \ref{fact:Proc_Basic_1}
(b) (with $E=(E_{0},\mathscr{O}_{E}(E_{0}))$).

$\{Z^{i}\}_{i\in\mathbf{I}}$ satisfies $\mathbf{S}$-PSMTC by (\ref{eq:Ori_Proc_(T,E0)-MOD})
and the condition (iv) above. At the same time, we observe by (\ref{eq:Ori_Proc_(T,E0)-MOD})
that
\begin{equation}
\mathbb{E}^{i}\left[\widetilde{f}\circ Z_{\mathbf{T}_{0}}^{i}\right]=\mathbb{E}^{i}\left[f\circ X_{\mathbf{T}_{0}}^{i}\right]\label{eq:Ori_Proc_(T,E0)-Mod_Int_Compare}
\end{equation}
for all $f\in\mathfrak{mc}[\Pi^{\mathbf{T}_{0}}(C_{b}(E;\mathbf{R}))]$,
$\mathbf{T}_{0}\in\mathscr{P}_{0}(\mathbf{T})$ and $i\in\mathbf{I}$,
so $\{Z^{i}\}_{i\in\mathbf{I}}$ is $(\mathbf{S},\widetilde{\mathcal{D}})$-FDC.

Now, we apply Theorem \ref{thm:FLP_Gen} (b) (with $E=(E,\mathscr{O}_{E}(E_{0}))$,
$\mathcal{D}=\widetilde{\mathcal{D}}$, $X^{i}=Z^{i}$ and $X=Z$)
and obtain an $(E_{0},\mathscr{O}_{E}(E_{0}))$-valued process $(\Omega,\mathscr{F},\mathbb{P};Z)$
satisfying: (1) $\mathbf{S}$-PMTC (in $E_{0}$), and (2) $Z=\mathfrak{fl}_{\mathbf{S}}(\{Z^{i_{n}}\}_{n\in\mathbf{N}})$
for any $\{i_{n}\}_{n\in\mathbf{N}}\subset\mathbf{I}$.

Considering $Z$ as an $E$-valued process with paths in $E_{0}^{\mathbf{R}^{+}}$,
it follows by $Z$'s property (2) above, Fact \ref{fact:FC_FDC} (with
$X^{n}=Z^{i_{n}}$ and $X=Z$) and (\ref{eq:Ori_Proc_(T,E0)-Mod_Int_Compare})
that $\mathbb{E}[f\circ Z_{\mathbf{T}_{0}}]$ is the unique limit
point of $\{\mathbb{E}^{i}[f\circ X_{\mathbf{T}_{0}}^{i}]\}_{i\in\mathbf{I}}$
for all $f\in\mathfrak{mc}[\Pi^{\mathbf{T}_{0}}(C_{b}(E;\mathbf{R}))]$
and $\mathbf{T}_{0}\in\mathscr{P}_{0}(\mathbf{S})$. Now, the result
follows with $X=Z$ by $Z$'s property (1) and Corollary \ref{cor:FLP_Gen}
(a) (with $X=Z$ and $\mathcal{D}=C_{b}(E;\mathbf{R})$).\end{proof}

\begin{rem}
\label{rem:FLP_Weakly_Cadlag_Uni}Proposition \ref{prop:FLP_Gen_E0}
does not require $\mathcal{D}$ to separate points on the entire space
$E$ as in Theorem \ref{thm:FLP_Gen} (b).\end{rem}
\begin{cor}
\label{cor:FLP_Gen_E0}Let $E$ be a topological space, $\mathbf{S}\subset\mathbf{T}\subset\mathbf{R}^{+}$
and $(\Omega,\mathscr{F},\mathbb{P};X)$ and $\{(\Omega^{i},\mathscr{F}^{i},\mathbb{P}^{i};X^{i})\}_{i\in\mathbf{I}}$
be $E$-valued processes. Suppose that:

\renewcommand{\labelenumi}{(\roman{enumi})}
\begin{enumerate}
\item $C_{b}(E;\mathbf{R})$ separates points on $E$.
\item $\mathcal{D}\subset C_{b}(E;\mathbf{R})$ separates points on $E_{0}\in\mathscr{B}(E)$.
\item (\ref{eq:Common_T-Base_I}) holds.
\item $\{X^{i}\}_{i\in\mathbf{I}}$ satisfies $\mathbf{S}$-PSMTC in $E_{0}$.
\item $X$ satisfies $\mathbf{S}$-PMTC.
\item $\mathbb{E}[f\circ X_{\mathbf{T}_{0}}]$ is the unique limit point
of $\{\mathbb{E}^{i}[f\circ X_{\mathbf{T}_{0}}^{i}]\}_{i\in\mathbf{I}}$
in $\mathbf{R}$ for all $f\in\mathfrak{mc}[\Pi^{\mathbf{T}_{0}}(\mathcal{D})]$
and $\mathbf{T}_{0}\in\mathscr{P}_{0}(\mathbf{S})$.
\end{enumerate}
Then, $X=\mathfrak{flp}_{\mathbf{S}}(\{X^{i}\}_{i\in\mathbf{I}})$
and $X=\mathfrak{fl}_{\mathbf{S}}(\{X^{i_{n}}\}_{n\in\mathbf{N}})$
for any $\{i_{n}\}_{n\in\mathbf{N}}\subset\mathbf{I}$.
\end{cor}
\begin{proof}
This corollary follows immediately by Proposition \ref{prop:FLP_Gen_E0}
and a similar argument to the proof of Corollary \ref{cor:FLP_Gen}
(a).\end{proof}

\section{\label{sec:LTB}Stationary long-time typical behavior}

\ref{enu:Q_LTB} in Introduction motivates our interest in finite-dimensional
convergence of stochastic processes. In order to utilize our results
in \S \ref{sec:FLP_Gen} and \S \ref{sec:FLP_Prog}, we introduce
the randomly advanced processes of a given measurable process $X$
whose finite-dimensional distributions are the long-time-averaged
distributions in (\ref{eq:LTB}).
\begin{defn}
\label{def:RAP}Let $E$ be a topological space and $(\Omega,\mathscr{F},\mathbb{P};X)$
be an $E$-valued measurable process.
\begin{itemize}
\item For each $T\in(0,\infty)$, by $(\widetilde{\Omega},\widetilde{\mathscr{F}},\mathbb{P}^{T};X^{T})=\mathfrak{rap}_{T}(X)$%
\footnote{``$\mathfrak{rap}$'' is ``rap'' in fraktur font which stands
for randomly advanced process.%
} ($X^{T}=\mathfrak{rap}_{T}(X)$ for short) we denote that $\widetilde{\Omega}\circeq\mathbf{R}^{+}\times\Omega$,
$\widetilde{\mathscr{F}}\circeq\mathscr{B}(\mathbf{R}^{+})\otimes\mathscr{F}$,
\begin{equation}
\mathbb{P}^{T}(A)\circeq\frac{1}{T}\int_{0}^{T}\int_{\Omega}\mathbf{1}_{A}(\tau,\omega)\mathbb{P}(d\omega)d\tau,\;\forall A\in\widetilde{\mathscr{F}},\label{eq:RAM}
\end{equation}
and
\begin{equation}
X^{T}(\tau,\omega)(t)\circeq X_{\tau+t}(\omega),\;\forall t\in\mathbf{R}^{+},(\tau,\omega)\in\widetilde{\Omega}.\label{eq:RAP}
\end{equation}
$X^{T}\in(E^{\mathbf{R}^{+}})^{\widetilde{\Omega}}$ defined by (\ref{eq:RAP})
is called the \textbf{$T$-randomly advanced process of $X$}.
\item A \textbf{long-time typical behavior of $X$ along $\mathbf{T}$}
refers to a member of $\mathfrak{flp}_{\mathbf{T}}(\{X^{T_{k}}\}_{k\in\mathbf{N}})$
with $T_{k}\uparrow\infty$, $X^{T_{k}}=\mathfrak{rap}_{T_{k}}(X)$
for each $k\in\mathbf{N}$ and $\mathbf{R}^{+}\backslash\mathbf{T}$
being a countable subset of $(0,\infty)$%
\footnote{This means $\mathbf{T}$ is cocountable.%
}.
\end{itemize}
\end{defn}
\begin{rem}
\label{rem:RAP}As its name implies, the $T$-randomly advanced process
of $X$ is defined by advancing $X$ to start at a random time $(\tau,\omega)\mapsto\tau$
defined on $(\widetilde{\Omega},\widetilde{\mathscr{F}})$.
\end{rem}

Below is a justification of our definition of randomly advanced process. 
\begin{prop}
\label{prop:RAP}Let $E$ be a topological space, $(\Omega,\mathscr{F},\mathbb{P};X)$
be an $E$-valued measurable process and $T\in(0,\infty)$. Then,
$(\widetilde{\Omega},\widetilde{\mathscr{F}},\mathbb{P}^{T};X^{T})=\mathfrak{rap}_{T}(X)$
is an $E$-valued measurable process.
\end{prop}
\begin{proof}
(\ref{eq:RAM}) well defines $\mathbb{P}^{T}\in\mathfrak{P}(\widetilde{\Omega},\widetilde{\mathscr{F}})$
by Fubini's Theorem. Let $\xi(t,\omega)\circeq X_{t}(\omega)$, $\xi^{T}(t,(\tau,\omega))\circeq X_{\tau+t}(\omega)$
and $\varphi(t,(\tau,\omega))\circeq(\tau+t,\omega)$ for each $t\in\mathbf{R}^{+}$
and $(\tau,\omega)\in\widetilde{\Omega}$. It is well-known that
\begin{equation}
\varphi\in M\left(\mathbf{R}^{+}\times\widetilde{\Omega},\mathscr{B}(\mathbf{R}^{+})\otimes\widetilde{\mathscr{F}};\widetilde{\Omega},\widetilde{\mathscr{F}}\right).\label{eq:Check_RAP_Measurable_1}
\end{equation}
$X$ being a measurable process implies $\xi\in M(\widetilde{\Omega},\widetilde{\mathscr{F}};E)$
and
\begin{equation}
\xi^{T}=\xi\circ\varphi\in M\left(\mathbf{R}^{+}\times\widetilde{\Omega},\mathscr{B}(\mathbf{R}^{+})\otimes\widetilde{\mathscr{F}};E\right),\label{eq:Check_RAP_Measurable_2}
\end{equation}
thus proving $X^{T}$ is a measurable process.\end{proof}

We present several further properties of randomly advanced process
in \S \ref{sec:Gen_Tech} of Appendix \ref{chap:App2}. Now, we give
our answer to \ref{enu:Q_LTB}.
\begin{thm}
\label{thm:LTB}Let $E$ be a topological space, $(\Omega,\mathscr{F},\mathbb{P};X)$
be an $E$-valued measurable process satisfying $T_{k}$-LMTC in $A\subset E$%
\footnote{The terminology ``$X$ satisfying $T_{k}$-LMTC in $A$'' was defined
in Definition \ref{def:Proc_Reg} and Note \ref{note:Singleton_PMTC}.%
} and $\mathcal{D}\subset C_{b}(E;\mathbf{R})$ separate points on
$E$. Then:

\renewcommand{\labelenumi}{(\alph{enumi})}
\begin{enumerate}
\item If $\{\frac{1}{T_{k}}\int_{0}^{T_{k}}\mathbb{E}[f\circ X_{\mathbf{T}_{0}+\tau}]d\tau\}_{k\in\mathbf{N}}$
is convergent in $\mathbf{R}$ for all $f\in\mathfrak{mc}[\Pi^{\mathbf{T}_{0}}(\mathcal{D})]$
and $\mathbf{T}_{0}\in\mathscr{P}_{0}(\mathbf{R}^{+})$ with $0\in\mathbf{T}_{0}$,
then $X$ has a stationary long-time typical behavior along $\mathbf{R}^{+}$.
\item If $\mathcal{D}$ is countable and $\{X^{T_{k}}\}_{k\in\mathbf{N}}$
satisfies $\mathcal{D}$-FMCC, then there exist a cocountable $\mathbf{S}\subset\mathbf{R}^{+}$,
an $E_{0}\in\mathscr{K}_{\sigma}^{\mathbf{m}}(E)$ such that $\{\frac{1}{T_{k}}\int_{0}^{T_{k}}\mathbb{P}\circ X_{\tau}^{-1}d\tau\}_{k\in\mathbf{N}}$
is $\mathbf{m}$-tight in $E_{0}\subset A$, and a stationary long-time
typical behavior of $X$ along $\mathbf{S}$ which is an $E$-valued,
$(\mathbf{R}^{+},\mathcal{D})$-c$\grave{\mbox{a}}$dl$\grave{\mbox{a}}$g,
progressive process with paths in $E_{0}^{\mathbf{R}^{+}}$.
\end{enumerate}
\end{thm}
\begin{proof}
(a) Let $X^{T_{k}}=\mathfrak{rap}_{T_{k}}(X)$ for each $k\in\mathbf{N}$.
By Proposition \ref{prop:Base_Proc_LMTC} (with $\{X^{n}\}_{n\in\mathbf{N}}=\{X\}$),
there exists a $\mathcal{D}$-baseable subset $E_{0}\subset\mathscr{K}_{\sigma}^{\mathbf{m}}(E)$
such that $E_{0}\subset A$ and (\ref{eq:T-Base}) holds for some
conull $\mathbf{T}\subset\mathbf{R}^{+}$ and $\{X_{0}^{T_{k}}\}_{k\in\mathbf{N}}$
is $\mathbf{m}$-tight in $E_{0}$. By Lemma \ref{lem:RAP_Tight}
(b, c, d) (with $A=E_{0}$), $\{X^{T_{k}}\}_{k\in\mathbf{N}}$ satisfies
$\mathbf{R}^{+}$-PSMTC in $E_{0}$, is $(\mathbf{R}^{+},\mathcal{D})$-AS
and is $(\mathbf{R}^{+},\mathcal{D})$-FDC. Now, (a) follows by Theorem
\ref{thm:FLP_Gen} (c) (with $X_{i}=X^{T_{i}}$).

(b) Let $E_{0}$ be as above. As (\ref{eq:T-Base}) holds for the
conull set $\mathbf{T}$, we have that
\begin{equation}
\inf_{t\in\mathbf{R}^{+},k\in\mathbf{N}}\mathbb{P}^{T_{k}}\left(X_{t}^{T_{k}}\in E_{0}\right)=1.\label{eq:RAP_Common_R+-Base_N}
\end{equation}
by Lemma \ref{lem:RepProc_RAP} (a) (with $T=T_{k}$). Now, (b) follows
by (\ref{eq:RAP_Common_R+-Base_N}), Proposition \ref{prop:FLP_Under_TF_FMCC}
(with $i=T_{k}$) and Theorem \ref{thm:FL_Weakly_Cadlag} (c) (with
$n=T_{k}$, $\mathbf{T}=\mathbf{R}^{+}$ and $\mathbf{S}=\mathbf{R}^{+}\backslash\mathbf{J}$).\end{proof}

\chapter{\label{chap:Cadlag}Application to Weak Convergence on Path Space}

\chaptermark{Weak Convergence on Path Space}

The current chapter addresses the target problems \ref{enu:Q_Tight}
and \ref{enu:Q_RC} of \ref{enu:Theme3} using the replication tools
developed in \S \ref{sec:RepProc_Path_Space}. Throughout this chapter,
we consider c$\grave{\mbox{a}}$dl$\grave{\mbox{a}}$g processes taking
values in a (at least) Tychonoff space $E$, whose common path space
is the Skorokhod $\mathscr{J}_{1}$-space $D(\mathbf{R}^{+};E)$.
If necessary, the readers can look back at \S \ref{sub:Meas_Cont_Cadlag_Map},
\S \ref{sec:RV} and \S \ref{sec:Proc} for our terminologies and
notations about the Skorokhod $\mathscr{J}_{1}$-space and c$\grave{\mbox{a}}$dl$\grave{\mbox{a}}$g
process. Also, \S \ref{sec:Sko} of Appendix \ref{chap:App1} together
with \S \ref{sec:Comp_A1} of Appendix \ref{chap:App2} provide a
short, almost self-contained review of Skorokhod $\mathscr{J}_{1}$-space.

The results of this chapter are sketched in Figure \ref{fig:SkoTight}
below.

\begin{figure}[H]
\begin{centering}
\includegraphics[scale=0.9]{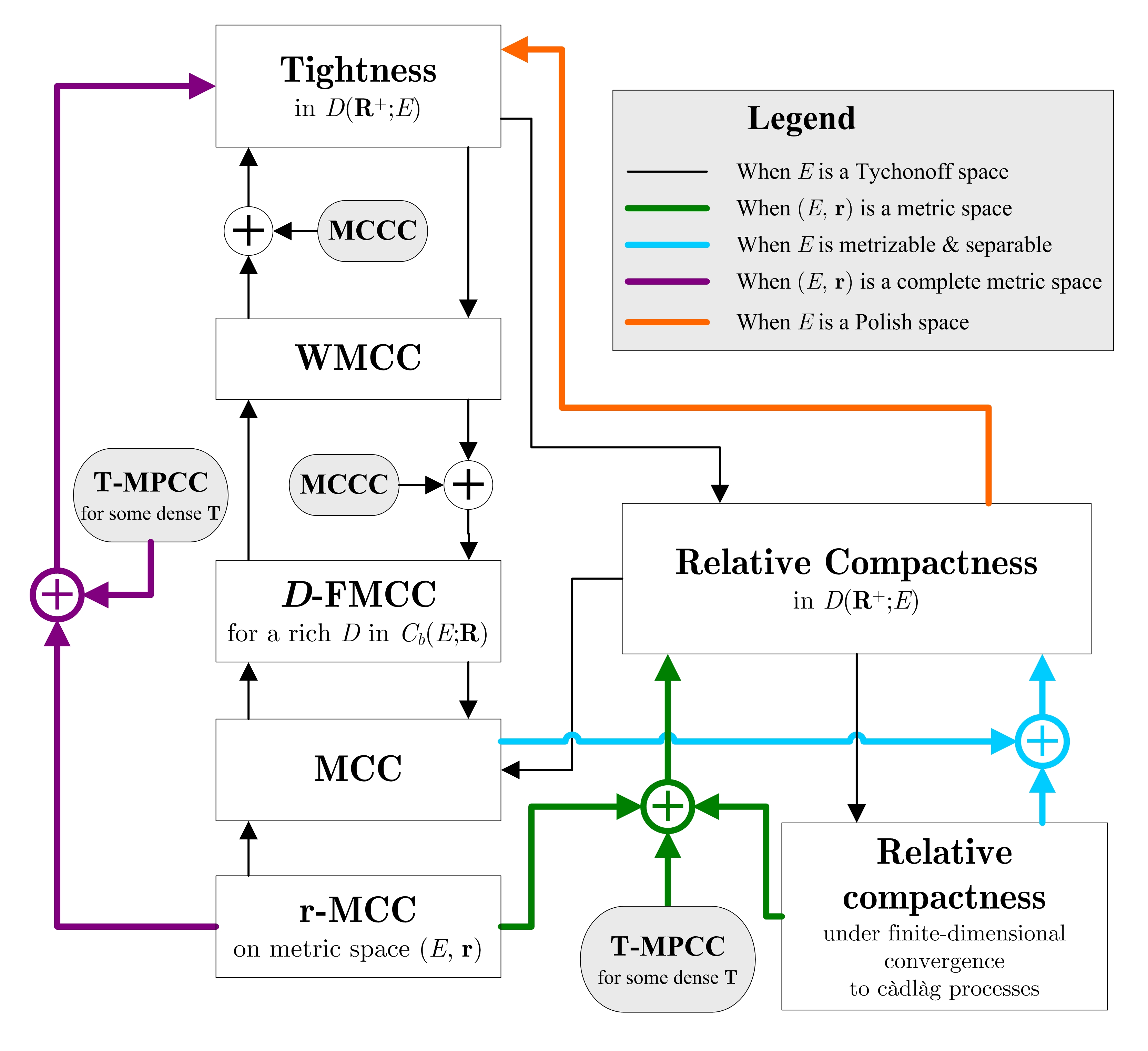}
\par\end{centering}

\caption{\textit{\label{fig:SkoTight}Tightness and relative compactness in
$D(\mathbf{R}^{+};E)$}}
\end{figure}

Given MCCC, \S \ref{sec:Sko_Tight} establishes the equivalence among
tightness in $D(\mathbf{R}^{+};E)$ and the MCC-type conditions introduced
in \S \ref{sub:RepProc_Tight}, which answers \ref{enu:Q_Tight}.
\S \ref{sec:Sko_WC_FC} looks into the relationship between weak
convergnce on $D(\mathbf{R}^{+};E)$ and finite-dimensional convergence.
\S \ref{sec:Sko_RC} establishes several results connecting finite-dimensional
convergence and relative compactness in $D(\mathbf{R}^{+};E)$, which
answers \ref{enu:Q_RC}.

Prior to the formal discussion, we recall several basic but essential
facts for this chapter. Let $E$ be a Tychonoff space, $\mu\in\mathcal{M}^{+}(D(\mathbf{R}^{+};E))$,
$\mathbf{T}_{0}\in\mathscr{P}_{0}(\mathbf{R}^{+})$ and $(\Omega,\mathscr{F},\mathbb{P};X)$
be an $E$-valued c$\grave{\mbox{a}}$dl$\grave{\mbox{a}}$g process.
Then:
\begin{itemize}
\item $D(\mathbf{R}^{+};E)$ is a Tychonoff space as mentioned in \S \ref{sub:Sko_Baseable}.
\item $\mu\circ\mathfrak{p}_{\mathbf{T}_{0}}^{-1}$%
\footnote{Herein, $\mathfrak{p}_{\mathbf{T}_{0}}$ denotes the projection on
$E^{\mathbf{R}^{+}}$ for $\mathbf{T}_{0}\subset\mathbf{R}^{+}$ restricted
to $D(\mathbf{R}^{+};E)$.%
} lies in $\mathfrak{M}^{+}(E^{\mathbf{T}_{0}},\mathscr{B}(E)^{\otimes\mathbf{T}_{0}})$
(see Corollary \ref{cor:Sko_Meas_FDD}).
\item The set $J(\mu)$%
\footnote{$J(\mu)$ was defined in (\ref{eq:J(Mu)}).%
} of fixed left-jump times of $\mu$ and the set $J(X)$%
\footnote{$J(X)$ was defined in (\ref{eq:J(X)}).%
} of fixed left-jump times of $X$ are well-defined when a countable
subset of $M(E;\mathbf{R})$ separates points on $E$. $\mathbf{R}^{+}\backslash J(\mu)$
and $\mathbf{R}^{+}\backslash J(X)$ are cocountable (hence non-empty
and dense) subsets of $\mathbf{R}^{+}$ when $E$ is baseable (see
Proposition \ref{prop:J(Mu)_J(X)_Baseable}).
\item $D(\mathbf{R}^{+};E)$-valued random variables are $E$-valued c$\grave{\mbox{a}}$dl$\grave{\mbox{a}}$g
processes (see Fact \ref{fact:Sko_RV_Cadlag} (a)) but the converse
is not necessarily true.
\item When $X$ is a $D(\mathbf{R}^{+};E)$-valued random variable, $\mathbb{P}\circ X^{-1}$
is the restriction of $\mathrm{pd}(X)|_{D(\mathbf{R}^{+};E)}$ to
$\sigma(\mathscr{J}(E))$%
\footnote{$\mathscr{J}(E)$ denotes the Skorokhod $\mathscr{J}_{1}$-topology
of $D(\mathbf{R}^{+};E)$. Restriction of measure to sub-$\sigma$-algebra
and $X$'s process distribution $\mathrm{pd}(X)$ were specified in
\S \ref{sub:Meas} and \S \ref{sec:Proc} respectively.%
} and $(\mathbb{P}\circ X^{-1})\circ\mathfrak{p}_{\mathbf{T}_{0}}^{-1}$
is the finite-dimensional distribution of $X$ for $\mathbf{T}_{0}$
as an $E$-valued process (see Fact \ref{fact:Sko_RV_Cadlag} (b,
c)).
\end{itemize}
Hereafter, we may not always make special reference for the facts
above.

\section{\label{sec:Sko_Tight}Tightness}

Our treatment of tightness in Skorokhod $\mathscr{J}_{1}$-space continues
\cite{K15} in the infinite time horizon setting. Tightness in $D(\mathbf{R}^{+};E)$
is stronger than the Compact Containment Condition in \cite{J86}
and \cite{EK86} (or MCCC if $E$ has metrizable compact subsets).
\begin{fact}
\label{fact:Tight_CCC}Let $E$ be a Tychonoff space and $\{(\Omega^{i},\mathscr{F}^{i},\mathbb{P}^{i};X^{i})\}_{i\in\mathbf{I}}$
be a tight family of $D(\mathbf{R}^{+};E)$-valued random variables.
Then, $\{X^{i}\}_{i\in\mathbf{I}}$ satisfies the Compact Containment
Condition in \cite[\S 4, (4.8)]{J86}. In particular, $\{X^{i}\}_{i\in\mathbf{I}}$
satisfies MCCC when $\mathscr{K}(E)=\mathscr{K}^{\mathbf{m}}(E)$.
\end{fact}
\begin{proof}
This fact follows by Proposition \ref{prop:Sko_Compact_Containment}
and Proposition \ref{prop:Sko_Compact}.\end{proof}

The following theorem is a version of \cite[Theorem 20]{K15} on infinite
time horizon. This result plus Fact \ref{fact:Tight_CCC}  answer
\ref{enu:Q_Tight} in Introduction.
\begin{thm}
\label{thm:Sko_Tight}Let $E$ be a Tychonoff space and $\{(\Omega^{i},\mathscr{F}^{i},\mathbb{P}^{i};X^{i})\}_{i\in\mathbf{I}}$
be $E$-valued c$\grave{\mbox{a}}$dl$\grave{\mbox{a}}$g processes.
Consider the following statements:

\renewcommand{\labelenumi}{(\alph{enumi})}
\begin{enumerate}
\item $X^{i}$ is indistinguishable from some $\widehat{X}^{i}\in M(\Omega^{i},\mathscr{F}^{i};D(\mathbf{R}^{+};E))$
for all $i\in\mathbf{I}$, and $\{\widehat{X}^{i}\}_{i\in\mathbf{I}}$
is $\mathbf{m}$-tight in $D(\mathbf{R}^{+};E)$.
\item $\{X^{i}\}_{i\in\mathbf{I}}$ satisfies $\mathcal{D}$-FMCC for some
$\mathcal{D}\subset C(E;\mathbf{R})$ and the closure of $\mathcal{D}$
under the topology of compact convergence (see \cite[\S 46, Definition, p.283]{M00})
contains $C_{b}(E;\mathbf{R})$.
\item $\{X^{i}\}_{i\in\mathbf{I}}$ satisfies MCC.
\item $\{X^{i}\}_{i\in\mathbf{I}}$ satisfies WMCC.
\end{enumerate}
Then, (a) implies (b), (c) implies (d), and (a) - (d) are all equivalent
when $\{X^{i}\}_{i\in\mathbf{I}}$ satisfies MCCC.\end{thm}
\begin{note}
\label{note:Proc_Tight}Part (a) of the theorem above addresses stronger
$\mathbf{m}$-tightness than tightness of the $D(\mathbf{R}^{+};E)$-valued
random variables $\{\widehat{X}^{i}\}_{i\in\mathbf{I}}$ above in
the usual sense%
\footnote{The $E$-valued processes $\{X^{i}\}_{i\in\mathbf{I}}$ in Theorem
\ref{thm:Sko_Tight} are $(E^{\mathbf{R}^{+}},\mathscr{B}(E)^{\otimes\mathbf{R}^{+}})$-valued
but not necessarily $D(\mathbf{R}^{+};E)$-valued random variables.
In general, $D(\mathbf{R}^{+};E)$ as a subset may not inherit the
measurability structure of $(E^{\mathbf{R}^{+}},\mathscr{B}(E)^{\otimes\mathbf{R}^{+}})$.
Hence, tightness of $\{X^{i}\}_{i\in\mathbf{I}}$ in $D(\mathbf{R}^{+};E)$
is defined in the extended way of \S \ref{sub:Tight}.%
}.\end{note}
\begin{rem}
\label{rem:Sko_Tight_1}The condition in Theorem \ref{thm:Sko_Tight}
(b) was used in \cite[p.628-629]{K75} to show tightness in $D(\mathbf{R}^{+};E)$
with $E$ being a locally compact Polish space. For general Polish
spaces, it appeared in \cite[\S 3.9]{EK86}.
\end{rem}

\begin{rem}
\label{rem:Sko_Tight_2}When $(E,\mathfrak{r})$ is a metric space,
the standard combination of $\mathfrak{r}$-MCC plus MCCC was used
as a sufficient condition for relative compactness in $D(\mathbf{R}^{+};E)$
by \cite[\S 3.7, Theorem 7.6]{EK86}. Its necessity was treated in
\cite[\S 3.7, Theorem 7.2, Remark 7.3]{EK86} with $E$ being a Polish
space. Theorem \ref{thm:Sko_Tight} refines these two results as well
as \cite[\S 3.9, Theorem 9.1]{EK86} and a few other analogues in
the following four aspects:
\begin{itemize}
\item We establish tightness which is often stronger than relative compactness.
\item The $E$ herein need not be metrizable nor separable.
\item We allow unbounded functions in $\mathcal{D}$, which can be handy
when working with algebras of polynomials for random measures as in
\cite[\S 2.1]{D91}.
\item $\{\varpi(f)\circ X^{i}\}_{i\in\mathbf{I}}$ satisfying $\left|\cdot\right|$-MCC
is milder than $\{\varpi(f)\circ X^{i}\}_{i\in\mathbf{I}}$ being
relatively compact if $f$ is not necessarily bounded. So, WMCC is
weaker than the analogous condition in \cite[Theorem 4.6, (4.9)]{J86}
which was shown very useful for establishing tightness of measure-valued
processes in \cite[\S 3.7]{D91} and \cite[\S II.4]{P02}.
\end{itemize}
\end{rem}
\begin{proof}
[Proof of Theorem \ref{thm:Sko_Tight}]((a) $\rightarrow$ (b)) For
each fixed $f\in\mathcal{D}\circeq C(E;\mathbf{R})$, $\{\varpi(f)\circ\widehat{X}^{i}\}_{i\in\mathbf{I}}$
is tight in $D(\mathbf{R}^{+};\mathbf{R})$ by Proposition \ref{prop:Sko_Basic_1}
(d) and Fact \ref{fact:Push_Forward_Tight_2} (with $E=D(\mathbf{R}^{+};E)$,
$S=D(\mathbf{R}^{+};\mathbf{R})$, $f=\varpi(f)$ and $\Gamma=\{\mathbb{P}^{i}\circ(\widehat{X}^{i})^{-1}\}_{i\in\mathbf{I}}$).
$\{\varpi(f)\circ\widehat{X}^{i}\}_{i\in\mathbf{I}}$ satisfies $\left|\cdot\right|$-MCC
by Theorem \ref{thm:Sko_RV_Tight_Polish} (a) (with $E=\mathbf{R}$).
Now, (b) follows by the bijective indistinguishability of $\{\varpi(f)\circ X^{i}\}_{i\in\mathbf{I}}$
and $\{\varpi(f)\circ\widehat{X}^{i}\}_{i\in\mathbf{I}}$ %
\footnote{Note \ref{note:Prog_Reg_Ind_Trans} mentioned the terminology ``bijective
indistinguishability'' and the transitivity of $\left|\cdot\right|$-MCC
between two bijectively indistinguishable families of processes.%
}.

((c) $\rightarrow$ (d)) is proved in Fact \ref{fact:MCC_4}.

((b) $\rightarrow$ (c) given MCCC) For each $g\in C_{b}(E;\mathbf{R})$,
$\epsilon\in(0,1/2)$ and $T\in(0,\infty)$, there exist $K_{\epsilon,T}\in\mathscr{K}^{\mathbf{m}}(E)$
and $f_{g,\epsilon,T}\in\mathcal{D}$ such that
\begin{equation}
\begin{aligned}\left\Vert f_{g,\epsilon,T}|_{K_{\epsilon,T}}-g|_{K_{\epsilon,T}}\right\Vert _{\infty} & \leq\epsilon<1-\epsilon\\
 & \leq\inf_{i\in\mathbf{I}}\mathbb{P}^{i}\left(X_{t}^{i}\in K_{\epsilon,T},\forall t\in[0,T]\right),
\end{aligned}
\label{eq:Check_MCC_1}
\end{equation}
which implies that
\begin{equation}
\begin{aligned} & \sup_{i\in\mathbf{I}}\mathbb{P}^{i}\left(\sup_{t\in[0,T]}\left|f_{g,\epsilon,T}\circ X_{t}^{i}-g\circ X_{t}^{i}\right|>\epsilon\right)\\
 & \leq1-\inf_{i\in\mathbf{I}}\mathbb{P}^{i}\left(X_{t}^{i}\in K_{\epsilon,T},\forall t\in[0,T]\right)<\epsilon.
\end{aligned}
\label{eq:Check_MCC_2}
\end{equation}
Then, $\{X^{i}\}_{i\in\mathbf{I}}$ satisfies MCC by (\ref{eq:Check_MCC_2}),
Proposition \ref{prop:CR} (a, c) and Proposition \ref{prop:MCE}
(a, b) (with $\mathcal{D}_{1}=\mathcal{D}$ and $\mathcal{D}_{2}=C_{b}(E;\mathbf{R})$).

((d) $\rightarrow$ (a) given MCCC) follows by Theorem \ref{thm:PR_Tight}
(with $E_{0}=E$).\end{proof}

When $(E,\mathfrak{r})$ is a complete (but not necessarily separable)
metric space and $\mathfrak{r}$-MCC is given, we have shown in Proposition
\ref{prop:Base_Proc_MPCC} the equivalence between MCCC and $\mathbf{T}$-MPCC
with a dense $\mathbf{T}$. This gives us one more tightness criterion.
\begin{prop}
\label{prop:SkoTight_MPCC}Let $(E,\mathfrak{r})$ be a metric space,
$\mathbf{T}$ be a dense subset of $\mathbf{R}^{+}$ and $E$-valued
c$\grave{\mbox{a}}$dl$\grave{\mbox{a}}$g processes $\{X^{i}\}_{i\in\mathbf{I}}$
satisfy $\mathfrak{r}$-MCC and $\mathbf{T}$-MPCC. Then:

\renewcommand{\labelenumi}{(\alph{enumi})}
\begin{enumerate}
\item There exists an $E_{0}\in\mathscr{B}(E)$ such that $E_{0}$ is a
separable subspace of $E$ and $X^{i}$ is indistinguishable from
a $D(\mathbf{R}^{+};E)$-valued random variable $\widehat{X}^{i}$
with paths in $E_{0}^{\mathbf{R}^{+}}$ for all $i\in\mathbf{I}$.
\item If $(E,\mathfrak{r})$ is complete, then the $\{\widehat{X}^{i}\}_{i\in\mathbf{I}}$
in (a) is tight in $D(\mathbf{R}^{+};E)$.
\end{enumerate}
\end{prop}
\begin{rem}
\label{rem:SkoTight_MPCC}$\,$
\begin{itemize}
\item Compared to \cite[\S 3.7, Theorem 7.2]{EK86}, part (b) above applies
to non-separable spaces, looses compact containment to totally bounded
containment and improves relative compactness into tightness in $D(\mathbf{R}^{+};E)$.
\item Compared to \cite[\S 3.7, Lemma 7.5]{EK86}, part (a) above replaces
MCCC by $\mathfrak{r}$-MCC plus $\mathbf{T}$-MPCC with a dense $\mathbf{T}$.
$\mathbf{T}$-MPCC is weaker than MCCC for any $\mathbf{T}\subset\mathbf{R}^{+}$
on metric spaces. In practice, $\mathfrak{r}$-MCC is usually no more
difficult than MCCC to verify.
\end{itemize}
\end{rem}
\begin{proof}
[Proof of Proposition \ref{prop:SkoTight_MPCC}](a) $C(E;\mathbf{R})$
separates points on $E$ by Proposition \ref{prop:CR_Space} (a) and
Proposition \ref{prop:CR} (a, c). By Proposition \ref{prop:Base_Proc_MPCC}
(a, b, c) (with $\mathcal{D}=C(E;\mathbf{R})$) and Proposition \ref{prop:PR}
(a) (with $X=X^{i}$ and $S_{0}\circeq D(\mathbf{R}^{+};E_{0},\mathscr{O}_{E}(E_{0}))$),
there exists an $E_{0}\in\mathscr{B}(E)$ such that $(E_{0},\mathscr{O}_{E}(E_{0}))$
is a separable subspace and $X^{i}$ is indistinguishable from an
$E$-valued process $\widehat{X}^{i}$ with paths in $D(\mathbf{R}^{+};E_{0},\mathscr{O}_{E}(E_{0}))$%
\footnote{``with paths in $D(\mathbf{R}^{+};E_{0},\mathscr{O}_{E}(E_{0}))$''
means all paths of the process lie in $D(\mathbf{R}^{+};E_{0},\mathscr{O}_{E}(E_{0}))$.%
} for all $i\in\mathbf{I}$. These $\{\widehat{X}^{i}\}_{i\in\mathbf{I}}$
are $D(\mathbf{R}^{+};E_{0},\mathscr{O}_{E}(E_{0}))$-valued variables
by Fact \ref{fact:Cadlag_Sko_RV} (a) (with $E=(E_{0},\mathscr{O}_{E}(E_{0}))$).
They are $D(\mathbf{R}^{+};E)$-valued random variables by Corollary
\ref{cor:Sko_Subspace} (with $A=E_{0}$).

(b) When $(E,\mathfrak{r})$ is complete, the $\{\widehat{X}^{i}\}_{i\in\mathbf{I}}$
above satisfies MCCC in $E_{0}$ by Proposition \ref{prop:Base_Proc_MPCC}
(d). Then, (b) follows by Corollary \ref{cor:MCC_5} (a) and Theorem
\ref{thm:Sko_Tight} (a, c).\end{proof}

\section{\label{sec:Sko_WC_FC}Weak convergence and finite-dimensional convergence}

We discuss in this section the relationship of the following properties
of $E$-valued c$\grave{\mbox{a}}$dl$\grave{\mbox{a}}$g processes
$\{(\Omega^{n},\mathscr{F}^{n},\mathbb{P}^{n};X^{n})\}_{n\in\mathbf{N}_{0}}$.
\begin{claim}
\label{claim:WC}$\:$

\begin{enumerate}
[label=\textbf{P\arabic*}, labelsep=0.5pc]
\setcounter{enumi}{7}

\item\label{enu:P8}There exists a dense subset $\mathbf{S}$ of
$\mathbf{R}^{+}$, a subset $\mathcal{D}$ of $C_{b}(E;\mathbf{R})$
such that%
\footnote{Note \ref{note:Proc_Int_Test_Integrability} argued that the expectations
in (\ref{eq:Exp_Test}) are well-defined.%
}
\begin{equation}
\lim_{n\rightarrow\infty}\mathbb{E}^{n}\left[f\circ X_{\mathbf{T}_{0}}^{n}\right]=\mathbb{E}^{0}\left[f\circ X_{\mathbf{T}_{0}}^{0}\right]\label{eq:Exp_Test_X0}
\end{equation}
for all $f\in\mathfrak{mc}[\Pi^{\mathbf{T}_{0}}(\mathcal{D})]$ and
$\mathbf{T}_{0}\in\mathscr{P}_{0}(\mathbf{S})$.

\item\label{enu:P9}\ref{enu:P8} holds with $\mathcal{D}$ separating
points on $E$.

\item\label{enu:P10}\ref{enu:P8} holds with $\mathcal{D}$ strongly
separating points on $E$.

\item\label{enu:P11}$\mathbf{S}\subset\mathbf{R}^{+}$ is dense
and $\{X^{n}\}_{n\in\mathbf{N}_{0}}$ satisfies
\begin{equation}
X^{n}\xrightarrow{\quad\mathrm{D}(\mathbf{S})\quad}X^{0}\mbox{ as }n\uparrow\infty.\label{eq:FC_along_S_X0}
\end{equation}

\item\label{enu:P12}There exist $\{\widehat{X}^{n}\in M(\Omega^{n},\mathscr{F}^{n};D(\mathbf{R}^{+};E))\}_{n\in\mathbf{N}_{0}}$
such that $X^{n}$ and $\widehat{X}^{n}$ are indistinguishable for
all $n\in\mathbf{N}_{0}$ and%
\footnote{Weak convergence, weak limit point and weak limit of random variables
were interpreted in \S \ref{sec:RV}.%
}
\begin{equation}
\widehat{X}^{n}\Longrightarrow\widehat{X}^{0}\mbox{ as }n\uparrow\infty\mbox{ on }D(\mathbf{R}^{+};E).\label{eq:RepProc_WC_D(E)_X0}
\end{equation}

\end{enumerate}

\end{claim}
\begin{rem}
\label{rem:D(E)_WL}Suppose $E$ is a Tychonoff space. Then, $\mathcal{P}(D(\mathbf{R}^{+};E))$
is a Tychonoff space by Proposition \ref{prop:E_CR_P(E)_CR}, and
so (\ref{eq:RepProc_WC_D(E)_X0}) is equivalent to that $\widehat{X}^{0}$
is the weak limit of $\{\widehat{X}^{n}\}_{n\in\mathbf{N}}$ on $D(\mathbf{R}^{+};E)$.
\end{rem}

Below are several immediate observations about \ref{enu:P8} - \ref{enu:P11}:
\begin{fact}
\label{fact:CR_Exp_Test}Let $E$ be a topological space and $\{(\Omega^{n},\mathscr{F}^{n},\mathbb{P}^{n};X^{n})\}_{n\in\mathbf{N}_{0}}$
be $E$-valued processes. Then, \ref{enu:P11} implies \ref{enu:P8}
with $\mathcal{D}=C_{b}(E;\mathbf{R})$. Moreover, if $E$ is a Hausdorff
(Tychonoff) space, then \ref{enu:P10} (resp. \ref{enu:P11}) implies
\ref{enu:P9} (resp. \ref{enu:P10} with $\mathcal{D}=C_{b}(E;\mathbf{R})$).
\end{fact}
\begin{proof}
Let $\mathcal{D}=C_{b}(E;\mathbf{R})$. \ref{enu:P11} implies \ref{enu:P8}
by Fact \ref{fact:FC_FDC} (with $X=X^{0}$ and $\mathbf{T}=\mathbf{S}$).
If $E$ is a Hausdorff space, \ref{enu:P10} implies \ref{enu:P9}
by Proposition \ref{prop:Fun_Sep_1} (a) (with $A=E$). If $E$ is
a Tychonoff space, \ref{enu:P11} implies \ref{enu:P10} since the
weaker property \ref{enu:P8} with $\mathcal{D}=C_{b}(E;\mathbf{R})$
implies \ref{enu:P10} by Proposition \ref{prop:CR} (a, c).\end{proof}

When $E$ is a metrizable and separable space, weak convergence on
$D(\mathbf{R}^{+};E)$ is well-known to imply finite-dimensional convergence
along all time points with no fixed left-jumps (see Theorem \ref{thm:Sko_RV_WC_FC_Metrizable_Separable}
(a)). Recall that every metrizable and separable space is baseable
(see Fact \ref{fact:Baseable_Metrizable_Separable} (b)), so the result
below generalizes Theorem \ref{thm:Sko_RV_WC_FC_Metrizable_Separable}
(a).
\begin{thm}
\label{thm:Sko_RV_WC_FC}Let $E$ be a Tychonoff space and $\{(\Omega^{n},\mathscr{F}^{n},\mathbb{P}^{n};X^{n})\}_{n\in\mathbf{N}_{0}}$
be $E$-valued c$\grave{\mbox{a}}$dl$\grave{\mbox{a}}$g processes.
Then:

\renewcommand{\labelenumi}{(\alph{enumi})}
\begin{enumerate}
\item If $M(E;\mathbf{R})$ has a countable subset separating points on
$E$%
\footnote{This condition ensures $J(X^{0})$ is well-defined.%
}, and if $\mathbf{S}\circeq\mathbf{R}^{+}\backslash J(X^{0})\neq\varnothing$,
then \ref{enu:P12} implies (\ref{eq:FC_along_S_X0}).
\item If $E$ is a baseable space, then \ref{enu:P12} implies \ref{enu:P11}
with $\mathbf{S}\circeq\mathbf{R}^{+}\backslash J(X^{0})$.
\end{enumerate}
\end{thm}
We prove the more general result below, and Theorem \ref{thm:Sko_RV_WC_FC}
then follows.
\begin{lem}
\label{lem:Sko_WC_FC}Let $E$ be a Tychonoff space and $\{\mu_{n}\}_{n\in\mathbf{N}_{0}}\subset\mathcal{M}^{+}(D(\mathbf{R}^{+};E))$
satisfy
\begin{equation}
\mu_{n}\Longrightarrow\mu_{0}\mbox{ as }n\uparrow\infty\mbox{ in }\mathcal{M}^{+}(D(\mathbf{R}^{+};E)).\label{eq:Mu_n_WC_Mu0_D(E)}
\end{equation}
If $M(E;\mathbf{R})$ has a countable subset separating points on
$E$%
\footnote{This condition ensures $J(\mu_{0})$ is well-defined.%
} and $\mathbf{R}^{+}\backslash J(\mu_{0})$ is non-empty, especially
if $E$ is a baseable space, then there exist $\{\nu_{\mathbf{T}_{0},n}\in\mathfrak{be}(\mu_{n}\circ\mathfrak{p}_{\mathbf{T}_{0}}^{-1})\}_{n\in\mathbf{N}_{0}}$
such that 
\begin{equation}
\nu_{\mathbf{T}_{0},n}\Longrightarrow\nu_{\mathbf{T}_{0},0}\mbox{ as }n\uparrow\infty\mbox{ in }\mathcal{M}^{+}(E^{\mathbf{T}_{0}})\label{eq:Mu_n_FDD_BExt_WC_Mu0_FDD_BExt}
\end{equation}
and%
\footnote{$\mu_{n}\circ\mathfrak{p}_{\mathbf{T}_{0}}^{-1}$ and $\mu_{0}\circ\mathfrak{p}_{\mathbf{T}_{0}}^{-1}$
are members of $\mathfrak{M}^{+}(E^{\mathbf{T}_{0}},\mathscr{B}(E)^{\mathbf{T}_{0}})$,
so the integrals in (\ref{eq:Sko_WC_FDD_Int_Test}) are well-defined
by Note \ref{note:Pi^d(D)_Mb_Cb}.%
}
\begin{equation}
\lim_{n\rightarrow\infty}\int_{E^{\mathbf{T}_{0}}}f(x)\mu_{n}\circ\mathfrak{p}_{\mathbf{T}_{0}}^{-1}(dx)=\int_{E^{\mathbf{T}_{0}}}f(x)\mu_{0}\circ\mathfrak{p}_{\mathbf{T}_{0}}^{-1}(dx)\label{eq:Sko_WC_FDD_Int_Test}
\end{equation}
for all $f\in\mathfrak{mc}[\Pi^{\mathbf{T}_{0}}(C_{b}(E;\mathbf{R}))]$
and $\mathbf{T}_{0}\in\mathscr{P}_{0}(\mathbf{R}^{+}\backslash J(\mu_{0}))$.
\end{lem}
\begin{proof}
The introductory part of this chapter noted that $\mathbf{R}^{+}\backslash J(\mu_{0})\neq\varnothing$
when $E$ is baseable. Let $\mathbb{D}\circeq D(\mathbf{R}^{+};E)$
and $(\mathbb{D},\mathscr{S}_{n},\nu_{n})$ be the completion of $(\mathbb{D},\mathscr{B}(\mathbb{D}),\mu_{n})$
for each $n\in\mathbf{N}_{0}$. It follows by Lemma \ref{lem:Sko_FDD_BExt}
(with $(\mathscr{S},\mu,\nu)=(\mathscr{S}_{n},\mu_{n},\nu_{n})$)
that
\begin{equation}
\nu_{\mathbf{T}_{0},n}\circeq\nu_{n}\circ\mathfrak{p}_{\mathbf{T}_{0}}^{-1}\in\mathfrak{be}\left(\mu_{n}\circ\mathfrak{p}_{\mathbf{T}_{0}}^{-1}\right),\;\forall n\in\mathbf{N}_{0}.\label{eq:Sko_FDD_Bext}
\end{equation}
It follows by (\ref{eq:Mu_n_WC_Mu0_D(E)}) and Fact \ref{fact:WC_Completion}
(with $E=\mathbb{D}$ and $\mathscr{U}_{n}=\mathscr{S}_{n}$) that
\begin{equation}
\nu_{n}\Longrightarrow\nu_{0}\mbox{ as }n\uparrow\infty\mbox{ in }\mathcal{M}^{+}(\mathbb{D}).\label{eq:Nu_n_WC_Nu_D(E)}
\end{equation}
The set of discontinuity points of $\mathfrak{p}_{\mathbf{T}_{0}}$
has zero measure under $\mu_{0}$ (hence under $\nu_{0}$) by the
definition of $J(\mu_{0})$ and Lemma \ref{lem:Sko_Proj} (c). Now,
(\ref{eq:Mu_n_FDD_BExt_WC_Mu0_FDD_BExt}) follows by (\ref{eq:Nu_n_WC_Nu_D(E)})
and the Continuous Mapping Theorem (Theorem \ref{thm:ContMapTh} (b)
with $E=\mathbb{D}$, $S=E^{\mathbf{T}_{0}}$, $\mu_{n}=\nu_{n}$,
$\mu=\nu_{0}$ and $f=\mathfrak{p}_{\mathbf{T}_{0}}$). (\ref{eq:Sko_WC_FDD_Int_Test})
follows by (\ref{eq:Mu_n_FDD_BExt_WC_Mu0_FDD_BExt}) and Fact \ref{fact:BExt_Same_Int}
(with $d=\aleph(\mathbf{T}_{0})$, $\mu=\mu_{n}\circ\mathfrak{p}_{\mathbf{T}_{0}}^{-1}$
and $\nu_{1}=\nu_{\mathbf{T}_{0},n}$).\end{proof}

\begin{proof}
[Proof of Theorem \ref{thm:Sko_RV_WC_FC}]Let $\{\widehat{X}^{n}\}_{n\in\mathbf{N}_{0}}$
be as in \ref{enu:P12}. Then, (a) follows by Lemma \ref{lem:Sko_WC_FC}
(with $\mu_{n}=\mathbb{P}^{n}\circ(\widehat{X}^{n})^{-1}$) and the
indistinguishability of $X^{n}$ and $\widehat{X}^{n}$. When $E$
is baseable, $\mathbf{S}$ is a dense subset of $\mathbf{R}^{+}$.
Thus (b) follows by (a).\end{proof}

The remainder of this section is about the converse of Theorem \ref{thm:Sko_RV_WC_FC}.
First of all, we establish a result about the converse of Lemma \ref{lem:Sko_WC_FC}.
\begin{lem}
\label{lem:Sko_FC_WC}Let $E$ be a baseable Tychonoff space, $\{\mu_{n}\}_{n\in\mathbf{N}_{0}}\subset\mathcal{M}^{+}(D(\mathbf{R}^{+};E))$
and $\mathcal{D}\subset C_{b}(E;\mathbf{R})$. Suppose that:

\renewcommand{\labelenumi}{(\roman{enumi})}
\begin{enumerate}
\item $\mathbf{S}$ is a dense subset of $\mathbf{R}^{+}$ and (\ref{eq:Sko_WC_FDD_Int_Test})
holds for all $f\in\mathfrak{mc}[\Pi^{\mathbf{T}_{0}}(\mathcal{D})]\cup1$
and $\mathbf{T}_{0}\in\mathscr{P}_{0}(\mathbf{S})$.
\item There exists an $S_{0}\in\mathscr{B}(D(\mathbf{R}^{+};E))$ such that
$\mu_{0}$ is supported on $S_{0}$ and
\begin{equation}
\mathscr{B}_{D(\mathbf{R}^{+};E)}(S_{0})=\left.\mathscr{B}(E)^{\otimes\mathbf{R}^{+}}\right|_{S_{0}}.\label{eq:S0_Sko_Borel_Prod_Equal}
\end{equation}

\item There exist $\{V_{p}\}_{p\in\mathbf{N}}\subset\mathscr{C}(D(\mathbf{R}^{+};E))$
such that $V_{p}\subset S_{0}$ for all $p\in\mathbf{N}$ and
\begin{equation}
\liminf_{n\rightarrow\infty}\mu_{n}\left(D(\mathbf{R}^{+};E)\backslash V_{p}\right)\leq2^{-p},\;\forall p\in\mathbf{N}.\label{eq:Vp_Sko_WLP_Support}
\end{equation}

\item $\{\mu_{n}\}_{n\in\mathbf{N}}$ is relatively compact.
\end{enumerate}
Then:

\renewcommand{\labelenumi}{(\alph{enumi})}
\begin{enumerate}
\item If $\mathcal{D}$ strongly separates points on $E$, then (\ref{eq:Mu_n_WC_Mu0_D(E)})
holds.
\item If $\mathcal{D}$ separates points on $E$, $\mu_{0}\circ\mathfrak{p}_{t}^{-1}$
is tight and $\{\mu_{n}\circ\mathfrak{p}_{t}^{-1}\}_{n\in\mathbf{N}}$
is sequentially tight for all $t$ in a conull $\mathbf{T}\subset\mathbf{R}^{+}$,
then (\ref{eq:Mu_n_WC_Mu0_D(E)}) holds.
\end{enumerate}
\end{lem}
\begin{proof}
Let $\gamma^{1}\circeq\mu_{0}$ and $\mathbb{D}\circeq D(\mathbf{R}^{+};E)$.
By the condition (iv) above and Fact \ref{fact:Uni_Seq_Lim_Conv},
(\ref{eq:Mu_n_WC_Mu0_D(E)}) follows if we show $\gamma^{1}=\gamma^{2}$
for any weak limit point $\gamma^{2}$ of $\{\mu_{n}\}_{n\in\mathbf{N}}$
in $\mathcal{M}^{+}(\mathbb{D})$. By passing to a subsequence if
necessary, we suppose
\begin{equation}
\mu_{n}\Longrightarrow\gamma^{2}\mbox{ as }n\uparrow\infty\mbox{ in }\mathcal{M}^{+}(\mathbb{D})\label{eq:Mu_n_WC_Gamma_2_D(E)}
\end{equation}
and let $\mathbf{S}^{1}\circeq\mathbf{S}$ and $\mathbf{S}^{2}\circeq J(\gamma^{2})$.
The rest of the proof is divided into three steps.

\textit{Step 1: Verify }
\begin{equation}
\lim_{n\rightarrow\infty}\int_{E^{\mathbf{T}_{0}}}f(x)\mu_{n}\circ\mathfrak{p}_{\mathbf{T}_{0}}^{-1}(dx)=\int_{E^{\mathbf{T}_{0}}}f(x)\gamma^{i}\circ\mathfrak{p}_{\mathbf{T}_{0}}^{-1}(dx)\label{eq:Sko_WLP_FDD_Uni_Int_Test}
\end{equation}
\textit{for each $f\in\mathfrak{mc}[\Pi^{\mathbf{T}_{0}}(\mathcal{D})]\cup1$,
$\mathbf{T}_{0}\in\mathscr{P}_{0}(\mathbf{S}^{i})$ and $i=1,2$}.
For $i=1$,  (\ref{eq:Sko_WLP_FDD_Uni_Int_Test}) is given by the
condition (i) above. For $i=2$, (\ref{eq:Sko_WLP_FDD_Uni_Int_Test})
follows by (\ref{eq:Mu_n_WC_Gamma_2_D(E)}) and Lemma \ref{lem:Sko_WC_FC}
(with $\mu_{0}=\gamma^{2}$).

\textit{Step 2: Verify}
\begin{equation}
\gamma^{1}\circ\mathfrak{p}_{\mathbf{T}_{0}}^{-1}=\gamma^{2}\circ\mathfrak{p}_{\mathbf{T}_{0}}^{-1}\mbox{ in }\mathfrak{M}^{+}\left(E^{\mathbf{T}_{0}},\mathscr{B}(E)^{\otimes\mathbf{T}_{0}}\right),\;\forall\mathbf{T}_{0}\in\mathscr{P}_{0}(\mathbf{R}^{+}).\label{eq:Sko_WLP_FDD_Uni}
\end{equation}
Under the conditions of (a), (\ref{eq:Sko_WLP_FDD_Uni}) follows immediately
by (\ref{eq:Sko_Meas}), Step 1 and Lemma \ref{lem:Sko_WLP_FDD_Uni}
(a). 

It takes a bit more work to show (\ref{eq:Sko_WLP_FDD_Uni}) for (b).
$E$ is baseable, so $\mathbf{S}^{2}$ is cocountable. $\mathbf{T}$
is conull by the hypothesis of (b), so $\mathbf{S}^{2}\cap\mathbf{T}$
is a conull hence dense subset of $\mathbf{R}^{+}$. Fixing $t\in\mathbf{S}^{2}\cap\mathbf{T}$,
we find from Step 1 that
\begin{equation}
\lim_{n\rightarrow\infty}\int_{E}f(x)\mu_{n}\circ\mathfrak{p}_{t}^{-1}(dx)=\int_{E}f(x)\gamma^{2}\circ\mathfrak{p}_{t}^{-1}(dx)\label{eq:Sko_WLP_Uni_1-D_Int_Test}
\end{equation}
for each $f\in\mathfrak{mc}(\mathcal{D})\cup1$. Letting $f=1$ in
(\ref{eq:Sko_WLP_Uni_1-D_Int_Test}), we find that $\{\mu_{n}\circ\mathfrak{p}_{t}^{-1}(E)\}_{n\in\mathbf{N}}$
must be contained in a compact sub-interval of $(0,\infty)$. $\gamma^{1}\circ\mathfrak{p}_{t}^{-1}$
is $\mathbf{m}$-tight and $\{\mu_{n}\circ\mathfrak{p}_{t}^{-1}\}_{n\in\mathbf{N}}$
is sequentially $\mathbf{m}$-tight by the hypothesis of (b), the
baseability of $E$ and Corollary \ref{cor:Baseable_MC} (a). It then
follows by (\ref{eq:Sko_WLP_Uni_1-D_Int_Test}) and Theorem \ref{thm:WLP_Uni}
(c) (with $d=1$ and $\Gamma=\{\mu_{n}\}_{n\in\mathbf{N}}$) that
\begin{equation}
\mathrm{w}\mbox{-}\lim_{n\rightarrow\infty}\mu_{n}\circ\mathfrak{p}_{t}^{-1}=\gamma^{2}\circ\mathfrak{p}_{t}^{-1}\label{eq:Sko_WLP_Uni_1-D_WC}
\end{equation}
and $\gamma^{2}\circ\mathfrak{p}_{t}^{-1}$ is $\mathbf{m}$-tight.

For each $\mathbf{T}_{0}\in\mathscr{P}(\mathbf{S}^{2}\cap\mathbf{T})$,
both $\gamma^{1}\circ\mathfrak{p}_{\mathbf{T}_{0}}^{-1}$ and $\gamma^{2}\circ\mathfrak{p}_{\mathbf{T}_{0}}^{-1}$
are $\mathbf{m}$-tight by Lemma \ref{lem:Tightness_Prod} (a) (with
$\mathbf{I}=\mathbf{T}_{0}$, $S_{i}=A_{i}=E$, $A=E^{\mathbf{T}_{0}}$
and $\Gamma=\{\gamma^{1}\circ\mathfrak{p}_{\mathbf{T}_{0}}^{-1}\}$
or $\{\gamma^{2}\circ\mathfrak{p}_{\mathbf{T}_{0}}^{-1}\}$). Thus,
(\ref{eq:Sko_WLP_FDD_Uni}) follows by (\ref{eq:Sko_Meas}) and Lemma
\ref{lem:Sko_WLP_FDD_Uni} (b) (with $\mathbf{S}=\mathbf{S}^{2}\cap\mathbf{T}$).

\textit{Step 3: Verify $\gamma^{1}=\gamma^{2}$ in $\mathcal{M}^{+}(\mathbb{D})$}.
As $\mathbb{D}$ is a Tychonoff space, we have that
\begin{equation}
\gamma^{2}(\mathbb{D}\backslash V_{p})\leq\liminf_{n\rightarrow\infty}\mu_{n}(\mathbb{D}\backslash V_{p})\leq2^{-p},\;\forall p\in\mathbf{N}\label{eq:Check_Sko_WLP_Supported}
\end{equation}
by $\{V_{p}\}_{p\in\mathbf{N}}\subset\mathscr{C}(\mathbb{D})$, (\ref{eq:Vp_Sko_WLP_Support})
and the Portmanteau's Theorem (Theorem \ref{thm:Portamenteau} (a,
c) with $E=\mathbb{D}$). As $S_{0}\in\mathscr{B}(\mathbb{D})$ contains
every $V_{p}$, we have by (\ref{eq:Check_Sko_WLP_Supported}) and
(ii) that
\begin{equation}
\gamma^{1}(\mathbb{D}\backslash S_{0})=\gamma^{2}(\mathbb{D}\backslash S_{0})=0.\label{eq:Sko_WLP_Support_S0}
\end{equation}
It follows by Step 2, the definition of $\mathscr{B}(E)^{\otimes\mathbf{R}^{+}}$
and (\ref{eq:Sko_WLP_Support_S0}) that
\begin{equation}
\gamma^{1}(A\cap S_{0})=\gamma^{2}(A\cap S_{0}),\;\forall A\in\left.\mathscr{B}(E)^{\otimes\mathbf{R}^{+}}\right|_{\mathbb{D}}.\label{eq:Sko_WLP_Same_Prod}
\end{equation}
It follows by (\ref{eq:Sko_WLP_Same_Prod}) and (\ref{eq:S0_Sko_Borel_Prod_Equal})
that
\begin{equation}
\gamma^{1}|_{S_{0}}=\gamma^{2}|_{S_{0}}\mbox{ in }\mathcal{M}^{+}\left(S_{0},\mathscr{O}_{\mathbb{D}}(S_{0})\right).\label{eq:Sko_WLP_Uni_on_S0}
\end{equation}
It then follows that
\begin{equation}
\gamma^{1}=\left.\gamma^{1}|_{S_{0}}\right|^{\mathbb{D}}=\left.\gamma^{2}|_{S_{0}}\right|^{\mathbb{D}}=\gamma^{2}\mbox{ in }\mathcal{M}^{+}(\mathbb{D})\label{eq:Sko_WLP_Uni}
\end{equation}
by (\ref{eq:Sko_WLP_Support_S0}), (\ref{eq:Sko_WLP_Uni_on_S0}) and
Fact \ref{fact:Meas_Concen_Expan} (c) (with $E=\mathbb{D}$, $\mathscr{U}=\mathscr{B}(E)$,
$A=S_{0}$ and $\mu=\gamma^{1}$ or $\gamma^{2}$).\end{proof}

The following proposition treats a typical case of Lemma \ref{lem:Sko_FC_WC}
where each $\mu_{n}$ is the distribution of $D(\mathbf{R}^{+};E)$-valued
random variable $X^{n}$ and the $\{V_{p}\}_{p\in\mathbf{N}}$ in
condition (iii) are compact sets provided by tightness.
\begin{prop}
\label{prop:Sko_RV_FC_WC_Tight}Let $E$ be a baseable Tychonoff space
and $X^{n}\in M(\Omega^{n},\mathscr{F}^{n},\mathbb{P}^{n};D(\mathbf{R}^{+};E))$
for each $n\in\mathbf{N}_{0}$. Suppose that:

\renewcommand{\labelenumi}{(\roman{enumi})}
\begin{enumerate}
\item There is an $S_{0}\in\mathscr{B}(D(\mathbf{R}^{+};E))$ satisfying
$\mathbb{P}^{0}(X^{0}\in S_{0})=1$ and (\ref{eq:S0_Sko_Borel_Prod_Equal}).
\item $\{X^{n}\}_{n\in\mathbf{N}}$ is tight in $S_{0}$.
\item $\mathbf{T}\subset\mathbf{R}^{+}$ is conull and $X^{0}$ satisfies
$\mathbf{T}$-PMTC.
\end{enumerate}
Then, \ref{enu:P9} implies \ref{enu:P12}.\end{prop}
\begin{rem}
\label{rem:Tight_in_S0}In Proposition \ref{prop:Sko_RV_FC_WC_Tight},
tightness in $S_{0}$ is different than tightness in $D(\mathbf{R}^{+};E)$
since $D(\mathbf{R}^{+};E)$ does not necessarily satisfy (\ref{eq:S0_Sko_Borel_Prod_Equal})
with $S_{0}=D(\mathbf{R}^{+};E)$ when $E$ is not separable.
\end{rem}
\begin{proof}
[Proof of Proposition \ref{prop:Sko_RV_FC_WC_Tight}]Let $\mu_{n}=\mathbb{P}^{n}\circ(X^{n})^{-1}\in\mathcal{P}(D(\mathbf{R}^{+};E))$
for each $n\in\mathbf{N}_{0}$. \ref{enu:P9} implies a dense subset
$\mathbf{S}$ of $\mathbf{R}^{+}$ and (\ref{eq:Sko_WC_FDD_Int_Test})
holds for $f\in\mathfrak{mc}[\Pi^{\mathbf{T}_{0}}(\mathcal{D})]\cup\{1\}$
and $\mathbf{T}_{0}\in\mathscr{P}_{0}(\mathbf{S})$. The condition
(i) above implies $\mu_{0}(S_{0})=1$. By the condition (ii) above,
there exist $\{V_{p}\}_{p\in\mathbf{N}}\subset\mathscr{K}(D(\mathbf{R}^{+};E))$
such that $V_{p}\subset S_{0}$ for all $p\in\mathbf{N}$ and $\inf_{n\in\mathbf{N}}\mu_{n}(V_{p})\geq1-2^{-p}$.
As $D(\mathbf{R}^{+};E)$ is a Tychonoff space, $\{V_{p}\}_{p\in\mathbf{N}}\subset\mathscr{C}(D(\mathbf{R}^{+};E))$
by Proposition \ref{prop:Compact} (a) and the tight sequence $\{\mu_{n}\}_{n\in\mathbf{N}}$
is relatively compact by the Prokhorov's Theorem (Theorem \ref{thm:Prokhorov}
(b)). As $E$ is a baseable space, the tight sequence $\{X^{n}\}_{n\in\mathbf{N}}$
satisfies MCCC by Fact \ref{fact:Tight_CCC} and Corollary \ref{cor:Base_Compact}
(a). $\{X^{n}\}_{n\in\mathbf{N}}$ satisfies $\mathbf{R}^{+}$-PMTC
by Fact \ref{fact:Base_PMTC} (f) (with $\mathbf{I}=\mathbf{N}$,
$i=n$ and $A=E$), so $\{\mu_{n}\circ\mathfrak{p}_{t}^{-1}\}_{n\in\mathbf{N}}$
is $\mathbf{m}$-tight for all $t\in\mathbf{R}^{+}$. Moreover, $\mu_{0}\circ\mathfrak{p}_{t}^{-1}$
is $\mathbf{m}$-tight for all $t\in\mathbf{T}$ by the condition
(iii) above. So far, we have justified all conditions of Lemma \ref{lem:Sko_FC_WC}
(b) for $\{\mu_{n}\}_{n\in\mathbf{N}_{0}}$, thus (\ref{eq:Mu_n_WC_Mu0_D(E)})
and \ref{enu:P12} hold.\end{proof}

The following proposition uses our tightness criteria established
in \S \ref{sec:Sko_Tight} to realize the ``tightness in $S_{0}$''
desired by Proposition \ref{prop:Sko_RV_FC_WC_Tight}.
\begin{prop}
\label{prop:Sko_RV_FC_WC_MCCC}Let $E$ be a Tychonoff space and $\{(\Omega^{n},\mathscr{F}^{n},\mathbb{P}^{n};X^{n})\}_{n\in\mathbf{N}_{0}}$
be $E$-valued c$\grave{\mbox{a}}$dl$\grave{\mbox{a}}$g processes.
If \ref{enu:P9} holds, $\{X^{n}\}_{n\in\mathbf{N}}$ satisfies $\mathcal{D}$-FMCC
for the $\mathcal{D}$ in \ref{enu:P9} and $\{X^{n}\}_{n\in\mathbf{N}_{0}}$
satisfies MCCC, then \ref{enu:P12} holds.
\end{prop}
\begin{proof}
The proof is divided into three steps.

\textit{Step 1: Construct a suitable base}. By Proposition \ref{prop:Base_Proc_CCC}
(b) (with $\mathbf{I}=\mathbf{N}_{0}$ and $i=n$), there exists a
$\mathcal{D}$-baseable subset $E_{0}$ of $E$ such that $\{X^{n}\}_{n\in\mathbf{N}_{0}}$
satisfies MCCC in $E_{0}$. By Lemma \ref{lem:Base_Construction}
(c) (with $\mathcal{D}_{0}=\{1\}$), there exists a base $(E_{0},\mathcal{F};\widehat{E};\widehat{\mathcal{F}})$
with $\mathcal{F}\subset(\mathcal{D}\cup\{1\})$.

\textit{Step 2: Construct $\{\widehat{X}^{n}\}_{n\in\mathbf{N}_{0}}$}.
$E_{0}$ is a Tychonoff subspace of $E$ by Proposition \ref{prop:CR_Space}
(b). $\mathbb{D}_{0}\circeq D(\mathbf{R}^{+};E_{0},\mathscr{O}_{E}(E_{0}))$
is a Tychonoff subspace of $D(\mathbf{R}^{+};E)$ by Corollary \ref{cor:Sko_Subspace}
(with $A=E_{0}$). $\{X^{n}\}_{n\in\mathbf{N}}$ satisfies $\mathcal{F}$-FMCC
since $(\mathcal{F}\backslash\{1\})\subset\mathcal{D}$. By Proposition
\ref{prop:PR_Tight} (with $\mathbf{I}=\mathbf{N}_{0}$), there exists
an $S_{0}\subset\mathbb{D}_{0}$ such that (\ref{eq:PR_S_Sko_Borel_Prod_Equal})
holds, there exist
\begin{equation}
\begin{aligned}\widehat{X}^{n} & =\mathfrak{rep}_{\mathrm{c}}(X^{n};E_{0},\mathcal{F})\in M\left(\Omega^{n},\mathscr{F}^{n};S_{0},\mathscr{O}_{\mathbb{D}_{0}}(S_{0})\right)\\
 & \subset M\left(\Omega^{n},\mathscr{F}^{n};\mathbb{D}_{0}\right)\subset M\left(\Omega^{n},\mathscr{F}^{n};D(\mathbf{R}^{+};E)\right),\;\forall n\in\mathbf{N}_{0}
\end{aligned}
\label{eq:Sko_RV_PR}
\end{equation}
satisfying
\begin{equation}
\inf_{n\in\mathbf{N}_{0}}\mathbb{P}^{n}\left(X^{n}=\widehat{X}^{n}\in S_{0}\right)=1,\label{eq:Sko_RV_PR_S0_Ind}
\end{equation}
and $\{\widehat{X}^{n}\}_{n\in\mathbf{N}}$ is $\mathbf{m}$-tight
in $S_{0}$%
\footnote{It is the subsequence $\{X^{n}\}_{n\in\mathbf{N}}$ that satisfies
$\mathcal{F}$-FMCC, so Proposition \ref{prop:PR_Tight} (c) is just
applied to $\{X^{n}\}_{n\in\mathbf{N}}$ with $X^{0}$ removed.%
}. Then, $S_{0}\in\mathscr{B}(E)^{\otimes\mathbf{R}^{+}}|_{\mathbb{D}_{0}}$
and (\ref{eq:S0_Sko_Borel_Prod_Equal}) holds by (\ref{eq:PR_S_Sko_Borel_Prod_Equal})
and Corollary \ref{cor:Sko_Subspace} (with $A=E_{0}$).

\textit{Step 3: Show} (\ref{eq:RepProc_WC_D(E)_X0}). It follows by
(\ref{eq:Sko_RV_PR_S0_Ind}), Fact \ref{fact:Common_T-Base_FR-Base}
(with $\mathbf{I}=\mathbf{N}_{0}$ and $i=n$), \ref{enu:P9}, the
fact $(\mathcal{F}\backslash\{1\})\subset\mathcal{D}$ and Lemma \ref{lem:RepProc_Int_Test}
(d) (with $X=X^{0}$) that
\begin{equation}
\lim_{n\rightarrow\infty}\mathbb{E}^{n}\left[f\circ\widehat{X}_{\mathbf{T}_{0}}^{n}\right]=\mathbb{E}^{0}\left[f\circ\widehat{X}_{\mathbf{T}_{0}}^{0}\right]\label{eq:RepProc_Exp_Test_X0}
\end{equation}
for all $f\in\mathfrak{mc}[\Pi^{\mathbf{T}_{0}}(\mathcal{F}\backslash\{1\})]$
and $\mathbf{T}_{0}\in\mathscr{P}_{0}(\mathbf{S})$. $X^{0}$ satisfies
MCCC, so it satisifes $\mathbf{R}^{+}$-PMTC by Fact \ref{fact:Base_PMTC}
(f) (with $A=E$ and $\mathbf{I}=\{0\}$). It then follows that
\begin{equation}
\widehat{X}^{n}\Longrightarrow\widehat{X}^{0}\mbox{ as }n\uparrow\infty\mbox{ on }\mathbb{D}_{0}\label{eq:RepProc_WC_D(E0)_X0}
\end{equation}
by (\ref{eq:Sko_RV_PR}) and Proposition \ref{prop:Sko_RV_FC_WC_Tight}
(with $E=E_{0}$, $\mathcal{D}=\mathcal{F}|_{E_{0}}\backslash\{1\}$,
$X^{n}=\widehat{X}^{n}$ and $\mathbf{T}=\mathbf{R}^{+}$). Now, (\ref{eq:RepProc_WC_D(E)_X0})
follows by (\ref{eq:Sko_RV_PR}), (\ref{eq:RepProc_WC_D(E0)_X0})
and Lemma \ref{lem:WC_Expansion} (with $E=D(\mathbf{R}^{+};E)$,
$A=\mathbb{D}_{0}$, $\mu_{n}=\mathbb{P}^{n}\circ(\widehat{X}^{n})^{-1}\in\mathcal{P}(\mathbb{D}_{0})$
and $\mu=\mathbb{P}^{0}\circ(\widehat{X}^{0})^{-1}\in\mathcal{P}(\mathbb{D}_{0})$).\end{proof}

Another typical case of Lemma \ref{lem:Sko_FC_WC} is when $E$ has
a metrizable and separable subspace $E_{0}$, and the $S_{0}$ in
condition (ii) and the $\{V_{p}\}_{p\in\mathbf{N}}$ in condition
(iii) of Lemma \ref{lem:Sko_FC_WC} are all taken to be $D(\mathbf{R}^{+};E_{0},\mathscr{O}_{E}(E_{0}))$.
Then, the assumption of relative compactness in Lemma \ref{lem:Sko_FC_WC}
(iv) can be loosened to $\mathcal{D}$-FMCC.
\begin{thm}
\label{thm:Sko_FC_WC_E0}Let $E$ be a Tychonoff space and $\{(\Omega^{n},\mathscr{F}^{n},\mathbb{P}^{n};X^{n})\}_{n\in\mathbf{N}_{0}}$
be $E$-valued c$\grave{\mbox{a}}$dl$\grave{\mbox{a}}$g processes.
Suppose that:

\renewcommand{\labelenumi}{(\roman{enumi})}
\begin{enumerate}
\item \ref{enu:P8} holds with $\mathcal{D}\subset C_{b}(E;\mathbf{R})$
being countable and strongly separating points on some $E_{0}\in\mathscr{B}(E)$.
\item $\{X_{n}\}_{n\in\mathbf{N}_{0}}$ satisfies
\begin{equation}
\inf_{n\in\mathbf{N}_{0}}\mathbb{P}^{n}\left(X^{n}\in E_{0}^{\mathbf{R}^{+}}\right)=1.\label{eq:Sko_RV_E0_Paths}
\end{equation}

\item $\{X^{n}\}_{n\in\mathbf{N}}$ satisfies $\mathcal{D}$-FMCC for the
$\mathcal{D}$ above.
\end{enumerate}
Then, \ref{enu:P12} holds.\end{thm}
\begin{rem}
\label{rem:RC_to_D-FMCC}The condition (i) above implies $E_{0}$
is a second-countable subspace of $E$ by Proposition \ref{prop:Fun_Sep_1}
(d) (with $A=E_{0}$). Given such $E_{0}$, the condition (i, iii)
above is weaker than relative compactness by Theorem \ref{thm:Sko_RV_RC_FDDRC}
(a) to follow, Corollary \ref{cor:MCC_2} and Proposition \ref{prop:Fun_Sep_2}
(b).
\end{rem}
\begin{proof}
[Proof of Theorem \ref{thm:Sko_FC_WC_E0}]The proof is divided into
four steps.

\textit{Step 1: Construct a suitable base by $E_{0}$ and $\mathcal{D}$}.
$E_{0}$ is a Tychonoff subspace and is a $\mathcal{D}$-baseable
subset of $E$ by Proposition \ref{prop:CR_Space} (b) and Proposition
\ref{prop:Fun_Sep_1} (a) (with $A=E_{0}$). As $\mathcal{D}$ is
countable, there exists a base $(E_{0},\mathcal{F};\widehat{E};\widehat{\mathcal{F}})$
with $\mathcal{F}=\mathcal{D}\cup\{1\}$ strongly separating points
on $E_{0}$ by Proposition \ref{lem:Base_Construction} (b).

\textit{Step 2: Construct $\{\widehat{X}^{n}\}_{n\in\mathbf{N}_{0}}$}.
Let $\mathbb{D}\circeq D(\mathbf{R}^{+};E)$, $S_{0}\circeq D(\mathbf{R}^{+};E_{0},\mathscr{O}_{E}(E_{0}))$
and $\widehat{\mathbb{D}}\circeq D(\mathbf{R}^{+};\widehat{E})$.
$S_{0}$ is a subspace of $\mathbb{D}$ by Corollary \ref{cor:Sko_Subspace}
(with $A=E_{0}$). $(E_{0},\mathscr{O}_{E}(E_{0}))$ is metrizable
and separable by Corollary \ref{cor:M_Compactification} (a, b) (with
$E=E_{0}$ and $\mathcal{D}=\mathcal{F}|_{E_{0}}$). Thus, $S_{0}$
satisfies (\ref{eq:S0_Sko_Borel_Prod_Equal}) by Proposition \ref{prop:Sko_Basic_2}
(b) (with $E=(E_{0},\mathscr{O}_{E}(E_{0}))$). Then, there exist
\begin{equation}
\begin{aligned}\widehat{X}^{n} & =\mathfrak{rep}_{\mathrm{c}}(X^{n};E_{0},\mathcal{F})\\
 & \in M\left(\Omega^{n},\mathscr{F}^{n};S_{0}\right)\cap M\left(\Omega^{n},\mathscr{F}^{n};\mathbb{D}\right)\cap M\left(\Omega^{n},\mathscr{F}^{n};\widehat{\mathbb{D}}\right),\;\forall n\in\mathbf{N}_{0}
\end{aligned}
\label{eq:Sko_RV_PR_FR}
\end{equation}
satisfying (\ref{eq:Sko_RV_PR_S0_Ind}) by (\ref{eq:Sko_RV_E0_Paths}),
Proposition \ref{prop:PR} (with $X=X^{n}$) and Fact \ref{fact:Cadlag_RepProc}.

\textit{Step 3: Show $\widehat{X}^{0}$ is the weak limit of $\{\widehat{X}^{n}\}_{n\in\mathbf{N}}$
on $\widehat{\mathbb{D}}$}. In this step, we consider $\{\widehat{X}^{n}\}_{n\in\mathbf{N}_{0}}$
as $\widehat{\mathbb{D}}$-valued random variables. As mentioned in
Note \ref{note:Ehat_Valued_Proc_FDD}, $\widehat{E}$ is a compact
Polish space, so $\{\widehat{X}^{n}\}_{n\in\mathbf{N}_{0}}$ automatically
satisfies MCCC by Note \ref{note:Base_CCC}. $\{X^{n}\}_{n\in\mathbf{N}}$
satisfies $\mathcal{F}$-FMCC since $\mathcal{F}\backslash\{1\}=\mathcal{D}$,
so $\{\widehat{X}^{n}\}_{n\in\mathbf{N}}$ satisfies $\widehat{\mathcal{F}}$-FMCC
by Proposition \ref{prop:FR_Tight} (a). (\ref{eq:RepProc_Exp_Test_X0})
holds for all $\widehat{f}\in\mathfrak{mc}[\Pi^{\mathbf{T}_{0}}(\widehat{\mathcal{F}}\backslash\{1\})]$
and $\mathbf{T}_{0}\in\mathscr{P}_{0}(\mathbf{S})$ by (\ref{eq:Sko_RV_E0_Paths}),
Fact \ref{fact:Common_T-Base_FR-Base} (with $\mathbf{I}=\mathbf{N}_{0}$
and $i=n$), \ref{enu:P8}, the fact $\mathcal{F}\backslash\{1\}=\mathcal{D}$
and Lemma \ref{lem:RepProc_Int_Test} (d) (with $X=X^{0}$). $\widehat{\mathcal{F}}\backslash\{1\}$
is a subset of $C_{b}(\widehat{E};\mathbf{R})$ by Corollary \ref{cor:Base_Fun_Dense}
(a) and separates points on $\widehat{E}$ by definition of base.
Now, the conclusion of Step 3 follows by Proposition \ref{prop:Sko_RV_FC_WC_MCCC}
(with $E=\widehat{E}$, $X^{n}=\widehat{X}^{n}$ and $\mathcal{D}=\widehat{\mathcal{F}}\backslash\{1\}$).

\textit{Step 4: Show (\ref{eq:RepProc_WC_D(E)_X0})}. $\{\widehat{X}^{n}\}_{n\in\mathbf{N}_{0}}$
as $S_{0}$-valued random variables satisfies (\ref{eq:RepProc_WC_D(E0)_X0})
with $\mathbb{D}_{0}=S_{0}$ by Step 3 and Proposition \ref{prop:PR_WC}
(b). $\{\widehat{X}^{n}\}_{n\in\mathbf{N}_{0}}$ as $\mathbb{D}$-valued
random variables satisfies (\ref{eq:RepProc_WC_D(E)_X0}) by (\ref{eq:Sko_RV_PR_FR}),
(\ref{eq:RepProc_WC_D(E0)_X0}) and Lemma \ref{lem:WC_Expansion}
(with $E=\mathbb{D}$, $A=S_{0}$, $\mu_{n}=\mathbb{P}^{n}\circ(\widehat{X}^{n})^{-1}\in\mathcal{P}(S_{0})$
and $\mu=\mathbb{P}^{0}\circ(\widehat{X}^{0})^{-1}\in\mathcal{P}(S_{0})$).\end{proof}

If $E$ itself is a metrizable and separable space, then the $E_{0}$
in Theorem \ref{thm:Sko_FC_WC_E0} can be taken to equal $E$.
\begin{cor}
\label{cor:Sko_FC_WC_Metrizable_Separable}Let $E$ be a metrizable
and separable space and $\{(\Omega^{n},\mathscr{F}^{n},\mathbb{P}^{n};X^{n})\}_{n\in\mathbf{N}_{0}}$
be $E$-valued c$\grave{\mbox{a}}$dl$\grave{\mbox{a}}$g processes.
Then, the following statements are successively weaker:

\renewcommand{\labelenumi}{(\alph{enumi})}
\begin{enumerate}
\item $\{X^{n}\}_{n\in\mathbf{N}}$ satisfies MCC and \ref{enu:P11}.
\item $\{X^{n}\}_{n\in\mathbf{N}}$ satisfies $\mathcal{D}$-FMCC and \ref{enu:P10}
for some $\mathcal{D}\subset C_{b}(E;\mathbf{R})$.
\item $\{X^{n}\}_{n\in\mathbf{N}}$ satisfies $\mathcal{D}$-FMCC and \ref{enu:P10}
for some countable $\mathcal{D}\subset C_{b}(E;\mathbf{R})$.
\item $\{X^{n}\}_{n\in\mathbf{N}}$ satisfies \ref{enu:P12}.
\end{enumerate}
\end{cor}
\begin{proof}
((a) $\rightarrow$ (b)) follows by Fact \ref{fact:CR_Exp_Test} and
Corollary \ref{cor:MCC_5} (b). ((b) $\rightarrow$ (c)) follows by
Proposition \ref{prop:Metrizable} (c) and Proposition \ref{prop:Fun_Sep_2}
(b). ((c) $\rightarrow$ (d)) follows by Theorem \ref{thm:Sko_FC_WC_E0}
(with $E_{0}=E$).\end{proof}

\begin{rem}
\label{rem:Sko_FC_WC_Metrizable_Separable}Compared to Theorem \ref{thm:Sko_RV_WC_FC_Metrizable_Separable}
(b), Corollary \ref{cor:Sko_FC_WC_Metrizable_Separable} (a, d) reduces
relative compactness to MCC which is a weaker condition by Theorem
\ref{thm:Sko_RV_RC_FDDRC} (a) to follow.
\end{rem}

When $E$ is a non-separable metric space, one can obtain the $E_{0}$
in Theorem \ref{thm:Sko_FC_WC_E0} by $\mathfrak{r}$-MCC and $\mathbf{T}$-MPCC.
\begin{prop}
\label{prop:Sko_FC_WC_Metrizable_MPCC}Let $(E,\mathfrak{r})$ be
a metric space and $\{(\Omega^{n},\mathscr{F}^{n},\mathbb{P}^{n};X^{n})\}_{n\in\mathbf{N}_{0}}$
be $E$-valued c$\grave{\mbox{a}}$dl$\grave{\mbox{a}}$g processes.
Then, the following statements are successively weaker:

\renewcommand{\labelenumi}{(\alph{enumi})}
\begin{enumerate}
\item $\{X^{n}\}_{n\in\mathbf{N}}$ satisfies $\mathfrak{r}$-MCC and $\mathbf{T}$-MPCC
with a dense $\mathbf{T}\subset\mathbf{R}^{+}$. $X^{0}$ satisfies
$\mathfrak{r}$-MCC. Moreover, \ref{enu:P11} holds.
\item $\{X^{n}\}_{n\in\mathbf{N}}$ satisfies $\mathfrak{r}$-MCC, $\mathcal{D}$-FMCC
with $\mathcal{D}\subset C_{b}(E;\mathbf{R})$ and $\mathbf{T}_{1}$-MPCC
with a dense $\mathbf{T}_{1}\subset\mathbf{R}^{+}$. $X^{0}$ satisfies
$\mathfrak{r}$-MCC and $\mathbf{T}_{2}$-MPCC with a dense $\mathbf{T}_{2}\subset\mathbf{R}^{+}$.
Moreover, \ref{enu:P10} holds.
\item $\{X^{n}\}_{n\in\mathbf{N}_{0}}$ satisfies \ref{enu:P12}.
\end{enumerate}
\end{prop}
\begin{proof}
((a) $\rightarrow$ (b)) By Corollary \ref{cor:MCC_5} (a), $\{X^{n}\}_{n\in\mathbf{N}}$
satisfies MCC. By Fact \ref{fact:CR_Exp_Test} and Corollary \ref{cor:MCC_2}
(a, b), there exists a $\mathcal{D}\subset C_{b}(E;\mathbf{R})$ such
that $\{X^{n}\}_{n\in\mathbf{N}}$ satisfies $\mathcal{D}$-FMCC and
\ref{enu:P10} holds. By Proposition \ref{prop:Base_Proc_MPCC} (a),
there exist $\{A_{p,q}\}_{p,q\in\mathbf{N}}\subset\mathscr{C}(E)$
such that each $A_{p,q}$ is a totally bounded set and
\begin{equation}
\inf_{n\in\mathbf{N}}\mathbb{P}^{n}\left(X_{t}^{n}\in A_{p,q}\right)\geq1-2^{-p},\;\forall t\in[0,q],\; p,q\in\mathbf{N}.\label{eq:Check_R+_MPCC}
\end{equation}
Hence, $\{X^{n}\}_{n\in\mathbf{N}}$ satisfies $\mathbf{R}^{+}$-MPCC.
For each $t$ in the $\mathbf{S}$ of \ref{enu:P11} and $n\in\mathbf{N}_{0}$,
$\mathbb{P}^{n}\circ(X_{t}^{n})^{-1}\in\mathcal{P}(E)$ by Fact \ref{fact:Proc_Basic_1}
(d) and so \ref{enu:P11} implies
\begin{equation}
X_{t}^{n}\Longrightarrow X_{t}^{0}\mbox{ as }n\uparrow\infty\mbox{ on }E.\label{eq:FC_1-D_Marginal_WC_X0}
\end{equation}
As $E$ is a Tychonoff space, it follows by (\ref{eq:Check_R+_MPCC}),
(\ref{eq:FC_1-D_Marginal_WC_X0}), the closedness of each $A_{p,q}$
and the Portmanteau's Theorem (Theorem \ref{thm:Portamenteau} (a,
b)) that
\begin{equation}
\begin{aligned}\mathbb{P}^{0}\left(X_{t}^{0}\in A_{p,q}\right) & \geq\inf_{n\in\mathbf{N}}\mathbb{P}^{n}\left(X_{t}^{n}\in A_{p,q}\right)\\
 & \geq1-2^{-p},\;\forall t\in\mathbf{S}\cap[0,q],p,q\in\mathbf{N},
\end{aligned}
\label{eq:Check_Sko_RV_Lim_MPCC}
\end{equation}
thus proving $X^{0}$ satisfies $\mathbf{S}$-MPCC. Now, (b) follows
by letting $\mathbf{T}_{1}=\mathbf{T}$ and $\mathbf{T}_{2}=\mathbf{S}$.

((b) $\rightarrow$ (c)) The union of two second-countable subspaces
of $E$ is still second-countable by Proposition \ref{prop:Metrizable}
(c) and Proposition \ref{prop:Countability} (b, e). So, we apply
Proposition \ref{prop:Base_Proc_MPCC} (a - c) to $\{X^{n}\}_{n\in\mathbf{N}}$
and the singleton $\{X^{0}\}$ respectively and find a second-countable
subspace $E_{0}$ of $E$ satisfying (\ref{eq:Sko_RV_E0_Paths}).
There exists a countable $\mathcal{D}_{0}\subset\mathcal{D}$ strongly
separating points on $E_{0}$ by \ref{enu:P10} and Proposition \ref{prop:Fun_Sep_2}
(b). \ref{enu:P10} implies \ref{enu:P8} so (c) follows by Theorem
\ref{thm:Sko_FC_WC_E0} (with $\mathcal{D}=\mathcal{D}_{0}$).\end{proof}

\section{\label{sec:Sko_RC}Relative compactness and finite-dimensional convergence}

When $(E,\mathfrak{r})$ is a separable metric space, relative compactness
in $D(\mathbf{R}^{+};E)$ implies $\mathfrak{r}$-MCC (see e.g. \cite[\S 3.7, Theorem 7.2]{EK86}).
For a general Tychonoff space $E$, we now justify the sufficiency
of relative compactness in $D(\mathbf{R}^{+};E)$ for MCC. If $E$
is also baseable, we leverage Theorem \ref{thm:Sko_RV_WC_FC} and
establish the sufficiency of relative compactness in $D(\mathbf{R}^{+};E)$
for ``relative compactness'' under finite-dimensional convergence.
\begin{thm}
\label{thm:Sko_RV_RC_FDDRC}Let $E$ be a Tychonoff space, $\mathbf{I}$
be an infinite index set and $\{(\Omega^{i},\mathscr{F}^{i},\mathbb{P}^{i};X^{i})\}_{i\in\mathbf{I}}$
be a relatively compact%
\footnote{Relative compactness of random variables was discussed in \S \ref{sec:Borel_Measure}.%
} family of $D(\mathbf{R}^{+};E)$-valued random variables. Then:

\renewcommand{\labelenumi}{(\alph{enumi})}
\begin{enumerate}
\item $\{X^{i}\}_{i\in\mathbf{I}}$ satisfies $C(E;\mathbf{R})$-FMCC and
MCC.
\item If $E$ is baseable, then any infinite subset of $\{X^{i}\}_{i\in\mathbf{I}}$
has a subsequence that converges finite-dimensionally to some $D(\mathbf{R}^{+};E)$-valued
random variable $X$ along $\mathbf{R}^{+}\backslash J(X)$.
\end{enumerate}
\end{thm}
\begin{proof}
By the relative compactness of $\{X^{i}\}_{i\in\mathbf{I}}$ in $D(\mathbf{R}^{+};E)$,
any infinite subset $\mathbf{J}$ of $\mathbf{I}$ contains a sequence
$\{i_{n}\}_{n\in\mathbf{N}}\subset\mathbf{J}\subset\mathbf{I}$ such
that
\begin{equation}
\mathbb{P}^{i_{n}}\circ(X^{i_{n}})^{-1}\Longrightarrow\mu\mbox{ as }n\uparrow\infty\mbox{ in }\mathcal{P}\left(D(\mathbf{R}^{+};E)\right)\label{eq:Sko_RV_WC_Subseq_P(D(E))}
\end{equation}
for some $\mu\in\mathcal{P}(D(\mathbf{R}^{+};E))$. Let $\Omega\circeq D(\mathbf{R}^{+};E)$,
$\mathscr{F}\circeq\mathscr{B}(\Omega)$, $\mathbb{P}\circeq\mu$
and $X$ be the identity mapping on $\Omega$. Then, $(\Omega,\mathscr{F},\mathbb{P};X)$%
\footnote{Such $X$ is known as the coordinate process on $D(\mathbf{R}^{+};E)$.%
} is a $D(\mathbf{R}^{+};E)$-valued random variable and satisfies
\begin{equation}
X^{i_{n}}\Longrightarrow X\mbox{ as }n\uparrow\infty\mbox{ on }D(\mathbf{R}^{+};E).\label{eq:WC_Subseq_D(E)}
\end{equation}
For each $f\in C(E;\mathbf{R})$, it follows that
\begin{equation}
\varpi(f)\circ X^{i_{n}}\Longrightarrow\varpi(f)\circ X\mbox{ as }n\uparrow\infty\mbox{ on }D(\mathbf{R}^{+};\mathbf{R})\label{eq:WC_Subseq_TF_D(R)}
\end{equation}
by (\ref{eq:WC_Subseq_D(E)}), Proposition \ref{prop:Sko_Basic_1}
(d) (with $S=\mathbf{R}$) and the Continuous Mapping Theorem (Theorem
\ref{thm:ContMapTh} (a) with $E=D(\mathbf{R}^{+};E)$, $S=D(\mathbf{R}^{+};\mathbf{R})$
and $f=\varpi(f)$). The argument above proves the relative compactness
of $\{\varpi(f)\circ X^{i}\}_{i\in\mathbf{I}}$ in $D(\mathbf{R}^{+};\mathbf{R})$.
$D(\mathbf{R}^{+};\mathbf{R})$, as mentioned in Note \ref{note:Ehat_Valued_Proc_FDD},
is a Polish space, so $\{\varpi(f)\circ X^{i}\}_{i\in\mathbf{I}}$
is tight in $D(\mathbf{R}^{+};\mathbf{R})$ by the Prokhorov's Theorem
(Theorem \ref{thm:Prokhorov} (a)) and satisfies $\left|\cdot\right|$-MCC
by Theorem \ref{thm:Sko_RV_Tight_Polish} (with $(E,\mathfrak{r})=(\mathbf{R},\left|\cdot\right|$).
$C(E;\mathbf{R})$ strongly separates points on $E$ by Proposition
\ref{prop:CR} (a, b). Now, (a) follows by Fact \ref{fact:MCC_1}
(b) (with $\mathcal{D}=C(E;\mathbf{R})$) and Corollary \ref{cor:MCC_2}
(a, d) (with $\mathcal{D}=C(E;\mathbf{R})$). If $E$ is also baseable,
then we have by Theorem \ref{thm:Sko_RV_WC_FC} (b) (with $n=i_{n}$)
that
\begin{equation}
X^{i_{n}}\xrightarrow{\quad\mathrm{D}(\mathbf{R}^{+}\backslash J(X))\quad}X\mbox{ as }n\uparrow\infty,\label{eq:FC_Subseq_along_R-J(X)}
\end{equation}
thus proving (b).\end{proof}

We then consider the converse of Theorem \ref{thm:Sko_RV_RC_FDDRC}.
\begin{thm}
\label{thm:Sko_FDDRC_RC_E0}Let $E$ be a Tychonoff space, $\mathbf{I}$
be an infinite index set and $\{(\Omega^{i},\mathscr{F}^{i},\mathbb{P}^{i};X^{i})\}_{i\in\mathbf{I}}$
be $E$-valued c$\grave{\mbox{a}}$dl$\grave{\mbox{a}}$g processes.
Suppose that for each infinite $\mathcal{I}^{*}\subset\mathbf{I}$,
there exists a subsequence $\mathcal{I}\circeq\{i_{n}\}_{n\in\mathbf{N}}\subset\mathcal{I}^{*}$,
an $E_{0,\mathcal{I}}\in\mathscr{B}(E)$, a $\mathcal{D}_{\mathcal{I}}\subset C_{b}(E;\mathbf{R})$,
an $\mathbf{S}_{\mathcal{I}}\subset\mathbf{R}^{+}$ and an $E$-valued
c$\grave{\mbox{a}}$dl$\grave{\mbox{a}}$g process $(\Omega,\mathscr{F},\mathbb{P};X^{\mathcal{I}})$
such that:

\renewcommand{\labelenumi}{(\roman{enumi})}
\begin{enumerate}
\item $\mathcal{D}_{\mathcal{I}}$ is countable and strongly separates points
on $E_{0,\mathcal{I}}$.
\item $\{X^{i_{n}}\}_{n\in\mathbf{N}}$ and $X^{\mathcal{I}}$ satisfy
\begin{equation}
\inf_{n\in\mathbf{N}}\mathbb{P}^{i_{n}}\left(X^{i_{n}}\in E_{0,\mathcal{I}}^{\mathbf{R}^{+}}\right)=\mathbb{P}\left(X^{\mathcal{I}}\in E_{0,\mathcal{I}}^{\mathbf{R}^{+}}\right)=1.\label{eq:Sko_RV_FDDRC_E0_Path}
\end{equation}

\item $\{X^{i_{n}}\}_{n\in\mathbf{N}}$ satisfies $\mathcal{D}_{\mathcal{I}}$-FMCC.
\item $\mathbf{S}_{\mathcal{I}}$ is dense in $\mathbf{R}^{+}$ and
\begin{equation}
\lim_{n\rightarrow\infty}\mathbb{E}^{i_{n}}\left[f\circ X_{\mathbf{T}_{0}}^{i_{n}}\right]=\mathbb{E}\left[f\circ X_{\mathbf{T}_{0}}^{\mathcal{I}}\right]\label{eq:Exp_Test_Subseq}
\end{equation}
for all $f\in\mathfrak{mc}[\Pi^{\mathbf{T}_{0}}(\mathcal{D}_{\mathcal{I}})]$
and $\mathbf{T}_{0}\in\mathscr{P}_{0}(\mathbf{S}_{\mathcal{I}})$.
\end{enumerate}
Then, there exist an $\mathbf{I}_{0}\in\mathscr{P}_{0}(\mathbf{I})$
and $D(\mathbf{R}^{+};E)$-valued random variables $\{\widehat{X}^{i}\}_{i\in\mathbf{I}\backslash\mathbf{I}_{0}}$
such that $\widehat{X}^{i}$ is indistinguishable from $X^{i}$ for
all $i\in\mathbf{I}\backslash\mathbf{I}_{0}$ and $\{\widehat{X}^{i}\}_{i\in\mathbf{I}\backslash\mathbf{I}_{0}}$
is relatively compact in $D(\mathbf{R}^{+};E)$.
\end{thm}
\begin{proof}
Let $\mathbb{D}\circeq D(\mathbf{R}^{+};E)$ and $\mathbb{D}_{0}\circeq D(\mathbf{R}^{+};E_{0,\mathcal{I}},\mathscr{O}_{E}(E_{0,\mathcal{I}}))$.
It follows by (\ref{eq:Sko_RV_FDDRC_E0_Path}) and the c$\grave{\mbox{a}}$dl$\grave{\mbox{a}}$g
properties of $\{X^{i}\}_{i\in\mathbf{I}}$ and $X^{\mathcal{I}}$
that
\begin{equation}
\inf_{n\in\mathbf{N}}\mathbb{P}^{i_{n}}\left(X^{i_{n}}\in\mathbb{D}_{0}\right)=\mathbb{P}\left(X^{\mathcal{I}}\in\mathbb{D}_{0}\right)=1.\label{eq:Sko_RV_D0_Paths}
\end{equation}
By Proposition \ref{prop:PR} (a) (with $S_{0}=\mathbb{D}_{0}$ and
$X=X^{i_{n}}$ or $X^{\mathcal{I}}$) and Corollary \ref{cor:Sko_Subspace}
(with $A=E_{0,\mathcal{I}}$), there exist
\begin{equation}
\left(\{\widehat{X}^{i_{n}}\}_{n\in\mathbf{N}}\cup\{\widehat{X}^{\mathcal{I}}\}\right)\subset M\left(\Omega^{i_{n}},\mathscr{F}^{i_{n}};\mathbb{D}_{0}\right)\subset M\left(\Omega^{i_{n}},\mathscr{F}^{i_{n}};\mathbb{D}\right)\label{eq:Sko_RV_FDDRC_PR}
\end{equation}
satisfying
\begin{equation}
\inf_{n\in\mathbf{N}}\mathbb{P}^{i_{n}}\left(X^{i_{n}}=\widehat{X}^{i_{n}}\in\mathbb{D}_{0}\right)=\mathbb{P}\left(X^{\mathcal{I}}=\widehat{X}^{\mathcal{I}}\in\mathbb{D}_{0}\right)=1.\label{eq:Sko_RV_FDDRC_PR_Ind}
\end{equation}
It follows by (\ref{eq:Sko_RV_FDDRC_PR_Ind}) and the condition (iv)
above that
\begin{equation}
\lim_{n\rightarrow\infty}\mathbb{E}^{i_{n}}\left[f\circ\widehat{X}_{\mathbf{T}_{0}}^{i_{n}}\right]=\mathbb{E}\left[f\circ\widehat{X}_{\mathbf{T}_{0}}^{\mathcal{I}}\right]\label{eq:Exp_Test_Subseq_PR}
\end{equation}
for each $f\in\mathfrak{mc}[\Pi^{\mathbf{T}_{0}}(\mathcal{D}_{\mathcal{I}})]$
and $\mathbf{T}_{0}\in\mathscr{P}_{0}(\mathbf{S}_{\mathcal{I}})$.
It then follows by Theorem \ref{thm:Sko_FC_WC_E0} (with $X^{n}=\widehat{X}^{i_{n}}$,
$X^{0}=\widehat{X}^{\mathcal{I}}$, $E_{0}=E_{0,\mathcal{I}}$, $\mathcal{D}=\mathcal{D}_{\mathcal{I}}$
and $\mathbf{S}=\mathbf{S}_{\mathcal{I}}$) that
\begin{equation}
\widehat{X}^{i_{n}}\Longrightarrow\widehat{X}^{\mathcal{I}}\mbox{ as }n\uparrow\infty\mbox{ on }\mathbb{D}.\label{eq:Check_Sko_RV_PR_RC}
\end{equation}

From the argument above we draw two conclusions: (1) There would be
at most finite members of $\{X^{i}\}_{i\in\mathbf{I}}$ which may
not admit an indistinguishable $\mathbb{D}$-valued copy. Let $\mathbf{I}_{0}\in\mathscr{P}_{0}(\mathbf{I})$
be the indices of these exceptions. For each $i\in\mathbf{I}\backslash\mathbf{I}_{0}$,
different $\{i_{n}\}_{n\in\mathbf{N}}\subset\mathbf{I}$ that contains
$i$ may induce different $\mathbb{D}$-valued copies of $X^{i}$.
However, such copy can be thought of as a unique one up to indistinguishability,
which we denote by $\widehat{X}^{i}$. (2) For any infinite $\mathcal{I}^{*}\subset(\mathbf{I}\backslash\mathbf{I}_{0})$,
there exist a subsequence $\{i_{n}\}_{n\in\mathbf{N}}\subset\mathcal{I}^{*}$
and a $\mathbb{D}$-valued random variable $\widehat{X}^{\mathcal{I}}$
such that (\ref{eq:Check_Sko_RV_PR_RC}) holds. In other words, $\{\widehat{X}^{i}\}_{i\in\mathbf{I}\backslash\mathbf{I}_{0}}$
is relatively compact in $\mathbb{D}$.\end{proof}

The following special cases of Theorem \ref{thm:Sko_RV_RC_FDDRC}
correspond to the settings of Corollary \ref{cor:Sko_FC_WC_Metrizable_Separable}
and Proposition \ref{prop:Sko_FC_WC_Metrizable_MPCC} respectively.
\begin{cor}
\label{cor:Sko_FDDRC_RC_Metrizable_Separable}Let $E$ be a metrizable
and separable space, $\mathbf{I}$ be an infinite index set and $\{(\Omega^{i},\mathscr{F}^{i},\mathbb{P}^{i};X^{i})\}_{i\in\mathbf{I}}$
be $E$-valued c$\grave{\mbox{a}}$dl$\grave{\mbox{a}}$g processes.
Then, the following statements are successively weaker:

\renewcommand{\labelenumi}{(\alph{enumi})}
\begin{enumerate}
\item Any infinite subset of $\{X^{i}\}_{i\in\mathbf{I}}$ has a subsequence
that satisfies MCC and converges finite-dimensionally to some $E$-valued
c$\grave{\mbox{a}}$dl$\grave{\mbox{a}}$g process along a dense subset
of $\mathbf{R}^{+}$.
\item For any infinite $\mathcal{I}^{*}\subset\mathbf{I}$, there exist
a subsequence $\mathcal{I}\circeq\{i_{n}\}_{n\in\mathbf{N}}\subset\mathcal{I}^{*}$,
a $\mathcal{D}_{\mathcal{I}}\subset C_{b}(E;\mathbf{R})$, a dense
$\mathbf{S}_{\mathcal{I}}\subset\mathbf{R}^{+}$ and an $E$-valued
c$\grave{\mbox{a}}$dl$\grave{\mbox{a}}$g process $(\Omega,\mathscr{F},\mathbb{P};X^{\mathcal{I}})$
such that: (i) $\mathcal{D}_{\mathcal{I}}$ is countable and strongly
separates points on $E$, (ii) $\{X^{i_{n}}\}_{n\in\mathbf{N}}$ satisfies
$\mathcal{D}_{\mathcal{I}}$-FMCC, and (iii) (\ref{eq:Exp_Test_Subseq})
holds for all $f\in\mathfrak{mc}[\Pi^{\mathbf{T}_{0}}(\mathcal{D}_{\mathcal{I}})]$
and $\mathbf{T}_{0}\in\mathscr{P}_{0}(\mathbf{S}_{\mathcal{I}})$.
\item There exist an $\mathbf{I}_{0}\in\mathscr{P}_{0}(\mathbf{I})$ and
$D(\mathbf{R}^{+};E)$-valued random variables $\{\widehat{X}^{i}\}_{i\in\mathbf{I}\backslash\mathbf{I}_{0}}$
such that $\widehat{X}^{i}$ is indistinguishable from $X^{i}$ for
all $i\in\mathbf{I}\backslash\mathbf{I}_{0}$ and $\{\widehat{X}^{i}\}_{i\in\mathbf{I}\backslash\mathbf{I}_{0}}$
is relatively compact in $D(\mathbf{R}^{+};E)$.
\end{enumerate}
\end{cor}
\begin{proof}
((a) $\rightarrow$ (b)) (i) and (ii) follow by Corollary \ref{cor:MCC_5}
(b). (iii) follows by (a) and Fact \ref{fact:CR_Exp_Test} (with $X^{n}=X^{i_{n}}$,
$X=X^{\mathcal{I}}$, $\mathbf{S}=\mathbf{S}_{\mathcal{I}}$ and $\mathcal{D}_{\mathcal{I}}\subset\mathcal{D}=C_{b}(E;\mathbf{R})$).

((b) $\rightarrow$ (c)) follows by Theorem \ref{thm:Sko_FDDRC_RC_E0}
(with $E_{0,\mathcal{I}}=E$).\end{proof}

\begin{prop}
\label{prop:Sko_FDDRC_RC_MPCC}Let $(E,\mathfrak{r})$ be a metric
space, $\mathbf{I}$ be an infinite index set and $\{(\Omega^{i},\mathscr{F}^{i},\mathbb{P}^{i};X^{i})\}_{i\in\mathbf{I}}$
be $E$-valued c$\grave{\mbox{a}}$dl$\grave{\mbox{a}}$g processes.
Then, the following statements are successively weaker:

\renewcommand{\labelenumi}{(\alph{enumi})}
\begin{enumerate}
\item For any infinite $\mathcal{I}^{*}\subset\mathbf{I}$, there exist
a subsequence $\mathcal{I}\circeq\{i_{n}\}_{n\in\mathbf{N}}\subset\mathcal{I}^{*}$,
an $E$-valued c$\grave{\mbox{a}}$dl$\grave{\mbox{a}}$g process
$(\Omega,\mathscr{F},\mathbb{P};X^{\mathcal{I}})$ and dense subsets
$\mathbf{T}_{\mathcal{I}}$ and $\mathbf{S}_{\mathcal{I}}$ of $\mathbf{R}^{+}$
such that: (i) $\{X^{i_{n}}\}_{n\in\mathbf{I}}$ satisfies $\mathfrak{r}$-MCC
and $\mathbf{T}_{\mathcal{I}}$-MPCC, (ii) $X^{\mathcal{I}}$ satisfies
$\mathfrak{r}$-MCC, and (iii)
\begin{equation}
X^{i_{n}}\xrightarrow{\quad\mathrm{D}(\mathbf{S}_{\mathcal{I}})\quad}X^{\mathcal{I}}\mbox{ as }n\uparrow\infty.\label{eq:FC_Subseq_along_S_I}
\end{equation}

\item For any infinite $\mathcal{I}^{*}\subset\mathbf{I}$, there exist
a sub-subsequence $\mathcal{I}\circeq\{i_{n}\}_{n\in\mathbf{N}}\subset\mathcal{I}^{*}$,
a $\mathcal{D}_{\mathcal{I}}\subset C_{b}(E;\mathbf{R})$, an $E$-valued
c$\grave{\mbox{a}}$dl$\grave{\mbox{a}}$g process $(\Omega,\mathscr{F},\mathbb{P};X^{\mathcal{I}})$
and dense subsets $\mathbf{T}_{\mathcal{I}}^{1}$, $\mathbf{T}_{\mathcal{I}}^{2}$,
$\mathbf{S}_{\mathcal{I}}$ of $\mathbf{R}^{+}$ such that: (i) $\mathcal{D}_{\mathcal{I}}$
strongly separates points on $E$, (ii) $\{X^{i_{n}}\}_{n\in\mathbf{N}}$
satisfies $\mathfrak{r}$-MCC, $\mathcal{D}_{\mathcal{I}}$-FMCC and
$\mathbf{T}_{\mathcal{I}}^{1}$-MPCC, (iii) $X^{\mathcal{I}}$ satisfies
$\mathfrak{r}$-MCC and $\mathbf{T}_{\mathcal{I}}^{2}$-MPCC, and
(iv) (\ref{eq:Exp_Test_Subseq}) holds for all $f\in\mathfrak{mc}[\Pi^{\mathbf{T}_{0}}(\mathcal{D}_{\mathcal{I}})]$
and $\mathbf{T}_{0}\in\mathscr{P}_{0}(\mathbf{S}_{\mathcal{I}})$.
\item There exist an $\mathbf{I}_{0}\in\mathscr{P}_{0}(\mathbf{I})$ and
$D(\mathbf{R}^{+};E)$-valued random variables $\{\widehat{X}^{i}\}_{i\in\mathbf{I}\backslash\mathbf{I}_{0}}$
such that $\widehat{X}^{i}$ is indistinguishable from $X^{i}$ for
all $i\in\mathbf{I}\backslash\mathbf{I}_{0}$ and $\{\widehat{X}^{i}\}_{i\in\mathbf{I}\backslash\mathbf{I}_{0}}$
is relatively compact in $D(\mathbf{R}^{+};E)$.
\end{enumerate}
\end{prop}
\begin{proof}
((a) $\rightarrow$ (b)) follows by Proposition \ref{prop:Sko_FC_WC_Metrizable_MPCC}
(a, b) (with $X^{n}=X^{i_{n}}$, $X^{0}=X^{\mathcal{I}}$, $\mathbf{T}=\mathbf{T}_{\mathcal{I}}$
and the $\mathbf{S}$ in \ref{enu:P11} being $\mathbf{S}_{\mathcal{I}}$).
((b) $\rightarrow$ (c)) follows by Proposition \ref{prop:Sko_FC_WC_Metrizable_MPCC}
(b, c) (with $X^{n}=X^{i_{n}}$, $X^{0}=X^{\mathcal{I}}$, $\mathbf{T}_{1}=\mathbf{T}_{\mathcal{I}}^{1}$
and $\mathbf{T}_{2}=\mathbf{T}_{\mathcal{I}}^{2}$) and our argument
about the finite exception set $\mathbf{I}_{0}$ in the proof of Theorem
\ref{thm:Sko_FDDRC_RC_E0}.\end{proof}

\chapter{\label{chap:App1}Background}

This appendix presents a series of background results used in Chapters
\ref{chap:Space_Change} - \ref{chap:Cadlag}. We limit our discussion
to the most necessary material.

\S \ref{sec:Point_Topo}, \S \ref{sec:SP_Fun} and \S \ref{sec:Compactification}
are about the point-set topology. More details are found in \cite[Chapter 1 - 7]{M00},
\cite[Vol. II, Chapter 6]{B07}, \cite[\S 3.4]{EK86} and \cite{BK10}.
\S \ref{sec:Meas_Sep} deals with weak topology of Borel measures
in the spirit of \cite[\S 3.1 - 3.4]{EK86}, \cite[Chapter 1]{KX95}
and \cite{BK10}. \S \ref{sec:SB} discusses standard Borel property
of topological spaces and their subsets, where we refer the readers
to \cite[\S 3.3]{S98} and \cite[Vol. II, Chapter 6]{B07} for further
materials. \S \ref{sec:Sko} gives a short review of Skorokhod $\mathscr{J}_{1}$-spaces.
Excellent treatments of this topic are available in \cite[\S 3.5 - 3.10]{EK86},
\cite{J86}, \cite{BK10} and \cite{K15}. \S \ref{sec:Cadlag} recalls
several basic properties of c$\grave{\mbox{a}}$dl$\grave{\mbox{a}}$g
processes.

This appendix compiles with all our notations, terminologies and conventions
introduced before. Several general technicalities used herein are
provided in \S \ref{sec:Gen_Tech} of Appendix \ref{chap:App2}.
A collection of miscellaneous results about the topics above are presented
in \S \ref{sec:Comp_A1} of Appendix \ref{chap:App2}.

\section{\label{sec:Point_Topo}Point-set topology}

In this appendix, $E$ denotes a topological space if not otherwise
specified.

\subsection{\label{sub:Separability}Separability}

\label{Hausdorff}$E$ is a Hausdorff space if for any distinct $x,y\in E$,
there exist disjoint $O_{x},O_{y}\in\mathscr{O}(E)$%
\footnote{$\mathscr{O}(E)$ and $\mathscr{C}(E)$ denotes the family of all
open and closed subsets of $E$ respectively.%
} such that $x\in O_{x}$ and $y\in O_{y}$. From this definition
one immediately observes that: 
\begin{fact}
\label{fact:Hausdorff_Refine}Any topological refinement%
\footnote{The terminology ``topological refinement'' was introduced in \S
\ref{sub:Topo}.%
} of a Hausdorff space is also Hausdorff.
\end{fact}
\label{T3}$E$ is a T3 space if $E$ is a Hausdorff space and for
any $x\in E$ and $F\in\mathscr{C}(E)$ excluding $x$, there exist
disjoint $O_{x},O_{F}\in\mathscr{O}(E)$ such that $x\in O_{x}$ and
$F\subset O_{F}$. \label{T4}$E$ is a T4 space if $E$ is a Hausdorff
space and for any disjoint $F_{1},F_{2}\in\mathscr{C}(E)$, there
exist disjoint $O_{1},O_{2}\in\mathscr{O}(E)$ such that $F_{1}\subset O_{1}$
and $F_{2}\subset O_{2}$. Below are several basic properties of
Hausdorff, T3 and T4 spaces.
\begin{prop}
\label{prop:Separability}The following statements are true:

\renewcommand{\labelenumi}{(\alph{enumi})}
\begin{enumerate}
\item Any finite subset of a Hausdorff space is closed.
\item The families of T4, T3 and Hausdorff spaces are successively larger.
\item Subspaces of a T3 or Hausdorff space are T3 or Hausdorff spaces respectively.
Moreover, closed subsets of a T4 space are T4 subspaces.
\item Any product space of T3 or Hausdorff spaces is a T3 or Hausdorff space
respectively.
\end{enumerate}
\end{prop}

\subsection{\label{sub:Countability}Countability}

\label{First_Countable}$E$ is first-countable if for each $x\in E$,
there exists a countable collection $\mathcal{O}_{x}\subset\mathscr{O}(E)$
such that $O\ni x$ for all $O\in\mathcal{O}_{x}$ and any $U\in\mathscr{O}(E)$
containing $x$ is the superset of some member of $\mathcal{O}_{x}$.
\label{separable}$E$ is separable if $E$ has a countable dense
subset. \label{Lindelof}$E$ is a Lindel$\ddot{\mbox{o}}$f space
if any $\{O_{i}\}_{i\in\mathbf{I}}\subset\mathscr{O}(E)$ satisfying
$E=\bigcup_{i\in\mathbf{I}}O_{i}$ admits a countable subset $\{O_{i_{n}}\}_{n\in\mathbf{N}}$
satisfying $E=\bigcup_{n\in\mathbf{N}}O_{i_{n}}$. \label{Hered_Lindelof}$E$
is a hereditary Lindel$\ddot{\mbox{o}}$f space if any subspace of
$E$ is a Lindel$\ddot{\mbox{o}}$f space. \label{Second_Countable}$E$
is a second-countable space if it admits a countable topological basis%
\footnote{The notion of topological basis was mentioned in \S \ref{sub:Topo}.%
}.
\begin{prop}
\label{prop:Countability}The following statements are true:

\renewcommand{\labelenumi}{(\alph{enumi})}
\begin{enumerate}
\item Subspaces of a first-countable, second-countable or heredetary Lindel$\ddot{\mbox{o}}$f
space are first-countable, second-countable or heredetary Lindel$\ddot{\mbox{o}}$f
spaces, respectively. Moreover, closed subsets of a Lindel$\ddot{\mbox{o}}$f
space are Lindel$\ddot{\mbox{o}}$f subspaces.
\item Every second-countable space is first-countable, separable and heriditary
Lindel$\ddot{\mbox{o}}$f simultaneously.
\item The product space of countably many first-countable, second-countable
or separable spaces is first-countable, second-countable or separable,
respectively.
\item The image of a separable, Lindel$\ddot{\mbox{o}}$f or heriditary
Lindel$\ddot{\mbox{o}}$f space under a continuous mapping is separable,
Lindel$\ddot{\mbox{o}}$f or heriditary Lindel$\ddot{\mbox{o}}$f,
respectively.
\item The union of countably many separable, Lindel$\ddot{\mbox{o}}$f or
hereditary Lindel$\ddot{\mbox{o}}$f subspaces is a separable, Lindel$\ddot{\mbox{o}}$f
or hereditary Lindel$\ddot{\mbox{o}}$f subspace, respectively.
\end{enumerate}
\end{prop}

\subsection{\label{sub:Metric}Metrizability}

\label{Metrizable}If the topology of $E$ is the same as the metric
topology induced by some metric $\mathfrak{r}$ on $E$, then $E$
is said to be \textit{metrizable}, $\mathfrak{r}$ is said to \textit{metrize}
$E$ and $(E,\mathfrak{r})$ is called a \textit{\label{Metrization}metrization}
of $E$.
\begin{prop}
\label{prop:Metrizable}The following statements are true:

\renewcommand{\labelenumi}{(\alph{enumi})}
\begin{enumerate}
\item Metrizable spaces are T4 (hence T3 and Hausdorff) spaces.
\item Subspaces of a metrizable space are metrizable.
\item Subspaces of a metrizable and separable space are metrizable and second-countable
(hence separable and hereditary Lindel$\ddot{\mbox{o}}$f).
\item Homeomprhs of metrizable spaces are metrizable.
\end{enumerate}
\end{prop}

Let $(E,\mathfrak{r})$ and $(S,\mathfrak{d})$ be metric spaces.
\label{Isometry}$f\in S^{E}$ is an \textit{isometry between $(E,\mathfrak{r})$
and $(S,\mathfrak{d})$} if $f$ is a surjective and $\mathfrak{r}(x,y)=\mathfrak{d}(f(x),f(y))$
for all $x,y\in E$. \label{Isometric}$(E,\mathfrak{r})$ and $(S,\mathfrak{d})$
are \textit{isometric} if there exists an isometry between them. \label{Completion}$(S,\mathfrak{d})$
is a completion of $(E,\mathfrak{r})$ if $(S,\mathfrak{d})$ is complete
and $E$ is isometric to a dense subspace of $S$.
\begin{note}
\label{note:Completion}Without loss of generality, a metric space
$(E,\mathfrak{r})$ can always be treated as a dense subset of its
completion $(S,\mathfrak{d})$ such that $\mathfrak{d}=\mathfrak{r}$
restricted to $E\times E$.\end{note}
\begin{prop}
\label{prop:Completeness}Let $(E,\mathfrak{r})$ and $(S,\mathfrak{d})$
be metric spaces. Then:

\renewcommand{\labelenumi}{(\alph{enumi})}
\begin{enumerate}
\item If $(E,\mathfrak{r})$ and $(S,\mathfrak{d})$ are isometric, then
they are homeomorphic. In particular, $(E,\mathfrak{r})$ is complete
precisely when $(S,\mathfrak{d})$ is complete.
\item There exists a unique completion of $(E,\mathfrak{r})$ up to isometry.
\item If $(E,\mathfrak{r})$ is complete, then the closure of $A\subset E$
equipped with (the restricted metric of) $\mathfrak{r}$ is the completion
of $A$.
\end{enumerate}
\end{prop}

The next two results are about metrizability of countable Cartesian
products.
\begin{prop}
\label{prop:Metric_Prod}Let $\{(S_{i},\mathfrak{r}_{i})\}_{i\in\mathbf{I}}$
be metric spaces and $S\circeq\prod_{i\in\mathbf{I}}S_{i}$. Then:

\renewcommand{\labelenumi}{(\alph{enumi})}
\begin{enumerate}
\item When $\mathbf{I}=\{1,...,d\}$, $S$ is metrized by%
\footnote{The notation ``$\mathfrak{p}_{i}$'' as defined in \S \ref{sub:Set_Num}
denotes one-dimensional projection on $S$ for $i\in\mathbf{I}$.%
}
\begin{equation}
\mathfrak{r}^{d}(x,y)\circeq\max_{1\leq i\leq d}\mathfrak{r}_{i}\left(\mathfrak{p}_{i}(x),\mathfrak{p}_{i}(y)\right),\;\forall x,y\in S.\label{eq:Prod_Metric}
\end{equation}
If $\{(S_{i},\mathfrak{r}_{i})\}_{1\leq i\leq d}$ are all complete,
then $(S,\mathfrak{r}^{d})$ is also.
\item When $\mathbf{I}=\mathbf{N}$, $S$ is metrized by
\begin{equation}
\mathfrak{r}_{1}^{\infty}(x,y)\circeq\sup_{i\in\mathbf{N}}i^{-1}\left[\mathfrak{r}_{i}\left(\mathfrak{p}_{i}(x),\mathfrak{p}_{i}(y)\right)\wedge1\right],\;\forall x,y\in S,\label{eq:Inf_Prod_Metric_1}
\end{equation}
or alternatively by
\begin{equation}
\mathfrak{r}_{2}^{\infty}(x,y)\circeq\sum_{i=1}^{\infty}2^{-i+1}\left[\mathfrak{r}_{i}\left(\mathfrak{p}_{i}(x),\mathfrak{p}_{i}(y)\right)\wedge1\right],\;\forall x,y\in S.\label{eq:Inf_Prod_Metric_2}
\end{equation}
If $\{(S_{i},\mathfrak{r}_{i})\}_{i\in\mathbf{N}}$ are all complete,
then $(S,\mathfrak{r}_{1}^{\infty})$ and $(S,\mathfrak{r}_{2}^{\infty})$
are also.
\end{enumerate}
\end{prop}

\begin{prop}
\label{prop:Metrizable_Prod}Let $\mathbf{I}$ be a countable index
set, $\{S_{i}\}_{i\in\mathbf{I}}$ be topological spaces and $S\circeq\prod_{i\in\mathbf{I}}S_{i}$.
Then, $A$ is a metrizable subspace of $S$ if and only if $\mathfrak{p}_{i}(A)$
is a metrizable subspace of $S_{i}$ for all $i\in\mathbf{I}$.
\end{prop}

The following property of first-countable and metrizable spaces is
indispensable.
\begin{fact}
\label{fact:First_Countable}Let $E$ be a topological space, $x\in E$
and $A\subset E$. Then:

\renewcommand{\labelenumi}{(\alph{enumi})}
\begin{enumerate}
\item If $E$ is metrizable, then $E$ is first-countable.
\item If there exist $\{x_{n}\}_{n\in\mathbf{N}}\subset A$ converging to
$x$ in $E$, then $x$ is a limit point%
\footnote{The notion of limit point was mentioned in p.\pageref{Limit_Point}.%
} of $A$. The converse is true when $E$ is first-countable space.
\end{enumerate}
\end{fact}

A subset $A$ of metric space $(E,\mathfrak{r})$ is said to be totally
bounded\label{Total_Bounded} if for any $\epsilon\in(0,\infty)$,
there exists an $A_{\epsilon}\in\mathscr{P}_{0}(E)$%
\footnote{$\mathscr{P}_{0}(E)$ denotes the family of all finite subsets of
$E$.%
} such that $E=\bigcup_{x\in A_{\epsilon}}\{y\in E:\mathfrak{r}(x,y)<\epsilon\}$.
\begin{prop}
\label{prop:Total_Bounded_1}Let $(E,\mathfrak{r})$ be a metric space.
Then:

\renewcommand{\labelenumi}{(\alph{enumi})}
\begin{enumerate}
\item If $A\subset E$ is totally bounded, $(A,\mathscr{O}_{E}(A))$ is
a second-countable space.
\item The union of finitely many totally bounded subsets of $E$ is totally
bounded.
\item If $A\subset E$ is totally bounded, then its closure is also.
\end{enumerate}
\end{prop}

\subsection{\label{sub:Var_Polish}Polish, Lusin and Souslin spaces}

Polish, Lusin and Souslin spaces are topological variations of complete
separable metric spaces. \label{Polish}Homeomorphs of complete separable
metric spaces are called Polish spaces. \label{Lusin}\label{Souslin}$E$
is a Lusin space (resp. Souslin space) if $E$ is a Hausdorff space
and there exists a bijective (resp. surjective) $f\in C(S;E)$ with
$S$ being a Polish space.
\begin{prop}
\label{prop:Var_Polish}The following statements are true:

\renewcommand{\labelenumi}{(\alph{enumi})}
\begin{enumerate}
\item Every Polish (resp. Lusin) space is a Lusin (resp. Souslin) space.
\item Open and closed subsets of a Polish, Lusin or Souslin space are Polish,
Lusin or Souslin subspaces, respectively.
\item Subspaces of a Polish space are metrizable and second-countable.
\item Subspaces of a Polish, Lusin or Souslin space are separable and heriditary
Lindel$\ddot{\mbox{o}}$f.
\item A metric space is separable if and only if its completion is a Polish
space.
\item The product space of countably many Polish, Lusin or Souslin spaces
is a Polish, Lusin or Souslin space, respectively. In particular,
$\mathbf{R}^{\infty}$ and its subspace $\mathbf{N}^{\infty}$ are
Polish spaces.
\item The intersection or union of countably many Souslin subspaces is a
Souslin subspace.
\end{enumerate}
\end{prop}

\subsection{\label{sub:Compactness}Compactness}

\label{Compact}$E$ is compact if any $\{O_{i}\}_{i\in\mathbf{I}}\subset\mathscr{O}(E)$
satisfying $E=\bigcup_{i\in\mathbf{I}}O_{i}$ admits a finite subset
$\{O_{i_{j}}\}_{1\leq j\leq n}$ satisfying $E=\bigcup_{j=1}^{n}O_{i_{j}}$.
\label{Seq_Compact}$E$ is sequentially compact (resp. \label{LP_Compact}limit
point compact) if any infinite subset of $E$ has a convergent subsequence
(resp. a limit point). Let $A\subset E$ be non-empty. \label{Compact_Set}$A$
is a compact, sequentially compact or limit point compact subset of
$E$ if $(A,\mathscr{O}_{E}(A))$ is compact, sequentially compact
or limit point compact, respectively. \label{RC}$A$ is a \textit{precompact}%
\footnote{Having a compact closure is commonly known as relative compactness.
Herein, we use the alternative terminology ``precompactness'' to
distinguish this notion from the relative compactness about finite
Borel measures.%
} subset of $E$ if $A$'s closure is compact. \label{Local_Compact}$E$
is locally compact if for any $x\in E$, there exist $K_{x}\in\mathscr{K}(E)$%
\footnote{$\mathscr{K}(E)$ denotes the family of all compact subsets of $E$.%
} and $O_{x}\in\mathscr{O}(E)$ such that $x\in O_{x}\subset K_{x}$.
\begin{prop}
\label{prop:Compact}The following statements are true:

\renewcommand{\labelenumi}{(\alph{enumi})}
\begin{enumerate}
\item Closed subsets of a compact space are compact. Moreover,  compact
subsets of a Hausdorff space are closed and hence Borel subsets.
\item The union of finitely many compact subsets is compact. Moreover, any
product space of compact spaces is a compact space.
\item Compact metric spaces are complete.
\item Hausdorff (resp. metrizable) compact spaces are T4 (resp. Polish)
spaces.
\item Compact spaces are Lindel$\ddot{\mbox{o}}$f spaces. Moreover, the
image of a compact space under any continuous mapping is a compact
space.
\item Compactness implies limit point compactness. Moreover, compactness,
sequential compactness and limit point compactness are equivalent
in metrizable spaces.
\end{enumerate}
\end{prop}

\begin{cor}
\label{cor:Compact_Prod}Let $\{S_{i}\}_{i\in\mathbf{I}}$ be topological
spaces and $(S,\mathscr{A})$ be as in (\ref{eq:(S,A)_Prod_Meas_Space}).
Then:

\renewcommand{\labelenumi}{(\alph{enumi})}
\begin{enumerate}
\item If $A_{i}\in\mathscr{K}(S_{i})$ for all $i\in\mathbf{I}$, then $\prod_{i\in\mathbf{I}}A_{i}\in\mathscr{K}(S)$.
If, in addition, $\mathbf{I}$ is countable and $\{S_{i}\}_{i\in\mathbf{I}}$
are all Hausdorff spaces, then $A_{i}\in\mathscr{B}(S_{i})$ for all
$i\in\mathbf{I}$ and $\prod_{i\in\mathbf{I}}A_{i}\in\mathscr{A}$.
\item If $A\in\mathscr{K}(S)$, then $\mathfrak{p}_{i}(A)\in\mathscr{K}(S_{i})$
for all $i\in\mathbf{I}$. If, in addition, $\{S_{i}\}_{i\in\mathbf{I}}$
are all Hausdorff spaces, then $\mathfrak{p}_{i}(A)\in\mathscr{B}(S_{i})$
for all $i\in\mathbf{I}$.
\end{enumerate}
\end{cor}

\begin{cor}
\label{cor:Sigma_Compact_Prod}Let $\mathbf{I}$ be a finite index
set, $\{S_{i}\}_{i\in\mathbf{I}}$ be topological spaces, $(S,\mathscr{A})$
be as in (\ref{eq:(S,A)_Prod_Meas_Space}) and $A_{i}\in\mathscr{K}_{\sigma}(S_{i})$
for all $i\in\mathbf{I}$. Then, $\prod_{i\in\mathbf{I}}A_{i}\in\mathscr{K}_{\sigma}(S)$.
\end{cor}

We used the following connection of total boundedness and compactness
in \S \ref{sec:Base_Proc}.
\begin{prop}
\label{prop:Total_Bounded_2}Compactness of a metric space is equivalent
to total boundedness plus completeness.
\end{prop}

\section{\label{sec:SP_Fun}Point-separation properties of functions}

\S \ref{sub:SP} introduced three functional separabilities of points:
separating points, strongly separating points and determining point
convergence. The following proposition specifies the relationship
among them.
\begin{prop}
\label{prop:Fun_Sep_1}Let $E$ be a topological space, $A\subset E$
be non-empty, $\mathcal{D}\subset\mathbf{R}^{E}$ and $d\in\mathbf{N}$.
Then:

\renewcommand{\labelenumi}{(\alph{enumi})}
\begin{enumerate}
\item If $\{\{x\}:x\in A\}\subset\mathscr{C}(E)$, especially if $A$ is
a Hausdorff subspace of $E$, then $\mathcal{D}$ strongly separating
points on $A$ implies $\mathcal{D}$ separating points on $A$.
\item If $\mathcal{D}$ strongly separates points on $A$, then \textup{$\mathcal{D}$}
determines point convergence on $A$. The converse is true when $(A,\mathscr{O}_{E}(A))$
is a Hausdorff space.
\item $\mathcal{D}$ separates points on $A$ if and only if $(A,\mathscr{O}_{\mathcal{D}}(A))$%
\footnote{The notation $\mathscr{O}_{\mathcal{D}}(A)$ was introduced in \S
\ref{sub:Topo}.%
} is a Hausdorff space.
\item $\mathscr{O}_{\mathcal{D}}(A)$ is induced by pseudometrics $\{\rho_{\{f\}}\}_{f\in\mathcal{D}}$%
\footnote{The pseudometric $\rho_{\mathcal{D}}$ was defined in \S \ref{sub:SP}
and $\rho_{\{f\}}$ refers to $\rho_{\mathcal{D}}$ with $\mathcal{D}=\{f\}$.
The meaning of $\{\rho_{\{f\}}:f\in\mathcal{D}\}$ inducing $\mathscr{O}_{\mathcal{D}}(A)$
was explained in \S \ref{sub:Topo}.%
}. If $\mathcal{D}$ is countable, then $(A,\mathscr{O}_{\mathcal{D}}(A))$
is a second-countable space pseudometrized by $\rho_{\mathcal{D}}$.
If, in addition, $\mathcal{D}$ separates points on $A$, then $\mathscr{O}_{\mathcal{D}}(A)$
is metrized by $\rho_{\mathcal{D}}$.
\item If $\mathcal{D}|_{A}\subset C(A,\mathscr{O}_{E}(A);\mathbf{R})$ separates
points (resp. strongly separates points) on $A$, then $(A,\mathscr{O}_{E}(A))$
is a Hausdorff space (resp. $\mathscr{O}_{E}(A)=\mathscr{O}_{\mathcal{D}}(A)$).
\end{enumerate}
\end{prop}

\begin{cor}
\label{cor:TF_Metric_FinDim}Let $E$ be a topological space, $\mathcal{D}\subset\mathbf{R}^{E}$
be countable and $d\in\mathbf{N}$. Then, $(E^{d},\mathscr{O}_{\mathcal{D}}(E)^{d})$
is a second-countable space pseudometrized by $\rho_{\mathcal{D}}^{d}$%
\footnote{The pseudometric $\rho_{\mathcal{D}}^{d}$ was defined in \S \ref{sub:SP}.%
}. If, in addition, $\mathcal{D}$ separates points on $E$, then $\mathscr{O}_{\mathcal{D}}(E)^{d}$
is metrized by $\rho_{\mathcal{D}}^{d}$.
\end{cor}

\begin{cor}
\label{cor:SSP_Dense}Let $E$ be a topological space and the members
of $\mathcal{D}_{0}\subset\mathbf{R}^{E}$ and $\mathcal{D}\subset\mathbf{R}^{E}$
are bounded. If $\mathcal{D}\subset\mathfrak{cl}(\mathcal{D}_{0})$%
\footnote{Recall that ``$\mathfrak{cl}(\cdot)$'' refers to closure under
supremum norm.%
}, then $\mathscr{O}_{\mathcal{D}}(E)\subset\mathscr{O}_{\mathcal{D}_{0}}(E)$.
In particular, if $\mathcal{D}$ separates points or strongly separates
points on $E$, then $\mathcal{D}_{0}$ does also.
\end{cor}

The following property of compact spaces is important for this work.
\begin{lem}
\label{lem:SP_on_Compact}Let $E$ be a compact space and $\mathcal{D}\subset C(E;\mathbf{R})$.
Then, $E$ is a Hausdorff space and $\mathcal{D}$ strongly separates
points on $E$ if and only if $\mathcal{D}$ separates points on $E$.
\end{lem}

Below are two useful properties of the function class $\Pi^{d}(\mathcal{D})$
introduced in \S \ref{sub:Fun}.
\begin{prop}
\label{prop:Pi^d_SP}Let $E$ be a topological space and $d\in\mathbf{N}$.
Then:

\renewcommand{\labelenumi}{(\alph{enumi})}
\begin{enumerate}
\item Any $\mathcal{D}\subset\mathbf{R}^{E}$ satisfies%
\footnote{The notations ``$\mathfrak{ag}(\cdot)$'' and ``$\mathfrak{ca}(\cdot)$''
were defined in \S \ref{sub:Fun}.%
}
\begin{equation}
\mathfrak{ag}\left[\Pi^{d}(\mathcal{D})\right]\subset\left[C\left(E^{d},\mathscr{O}_{\mathcal{D}}(E)^{d};\mathbf{R}\right)\cap M\left(E^{d},\sigma(\mathcal{D})^{\otimes d};\mathbf{R}\right)\right].\label{eq:ag(Pi^d(D))_C_Prod_M}
\end{equation}
If, in addition, the members of $\mathcal{D}$ are bounded, then,
\begin{equation}
\mathfrak{ca}\left[\Pi^{d}(\mathcal{D})\right]\subset\left[C_{b}\left(E^{d},\mathscr{O}_{\mathcal{D}}(E)^{d};\mathbf{R}\right)\cap M_{b}\left(E^{d},\sigma(\mathcal{D})^{\otimes d};\mathbf{R}\right)\right].\label{eq:ca(Pi^d(D))_Cb_Prod_Mb}
\end{equation}

\item If $\mathcal{D}$ separates points (resp. determines point convergence)
on $E$, then $\Pi^{d}(\mathcal{D})$%
\footnote{The definition of $\Pi^{d}(\mathcal{D})$ refers to (\ref{eq:Pi^d}).%
} separates points (resp. determines point convergence) on $E^{d}$.
\end{enumerate}
\end{prop}
\begin{rem}
\label{rem:Sigma(D)_Diff_B_D(E)}Please be reminded that the $\sigma$-algebra
$\sigma(\mathcal{D})$ induced by $\mathcal{D}$ is possibly smaller
than the Borel $\sigma$-algebra $\mathscr{B}_{\mathcal{D}}(E)$ induced
by $\mathcal{D}$ (see Fact \ref{fact:O_D_(A)_B_D_(A)}).
\end{rem}

\begin{rem}
\label{rem:Pi^d_Not_Need_1_in_D}$\Pi^{d}(\mathcal{D})$ is defined
in a way that we do not need $1\in\mathcal{D}$ in Proposition \ref{prop:Pi^d_SP}
(b).
\end{rem}

The following are two typical cases where one can select countably
many (strongly) point-separating functions.
\begin{prop}
\label{prop:Fun_Sep_2}Let $E$ be a topological space and $\mathcal{D}\subset C(E;\mathbf{R})$.
Then:

\renewcommand{\labelenumi}{(\alph{enumi})}
\begin{enumerate}
\item If $\{(x,x):x\in E\}$ is a Lindel$\ddot{\mbox{o}}$f subspace of
$E\times E$ and $\mathcal{D}$ separates points on $E$, then there
exists a countable $\mathcal{D}_{0}\subset\mathcal{D}$ that separates
points on $E$.
\item If $E$ is a second-countable space and $\mathcal{D}$ strongly separates
points on $E$, then there exists a countable $\mathcal{D}_{0}\subset\mathcal{D}$
that strongly separates points on $E$.
\end{enumerate}
\end{prop}

\section{\label{sec:Compactification}Tychonoff space and compactification}

\label{CR}$E$ is a Tychonoff (or T3$\frac{1}{2}$) space if $E$
is a Hausdorff space and for any $x\in E$ and $F\in\mathscr{C}(E)$
that excludes $x$, there exists an $f_{x,F}\in C(E;[0,1])$ such
that $f_{x,F}(x)=0$ and (the image) $f_{x,F}(F)=\{1\}$. Besides
the functional separability of points and closed sets above, Tychonoff
spaces are also defined as the spaces whose topology is induced by
some family of $\mathbf{R}$-valued functions, or alternatively by
some family of pseudometrics.
\begin{prop}
\label{prop:CR}Let $E$ be a topological space. Then, the following
statements are equivalent:

\renewcommand{\labelenumi}{(\alph{enumi})}
\begin{enumerate}
\item $E$ is a Tychonoff space.
\item $C(E;\mathbf{R})$ separates and strongly separates points on $E$.
\item $C_{b}(E;\mathbf{R})$ separates and strongly separates points on
$E$.
\item $E$ is a Hausdorff space and $\mathscr{O}(E)$ is induced by a family
of pseudometrics.
\item $E$ is a Hausdorff space and $\mathscr{O}(E)=\mathscr{O}_{\mathcal{D}}(E)$
for some $\mathcal{D}\subset\mathbf{R}^{E}$.
\end{enumerate}
\end{prop}

Below are a few more properties of Tychonoff spaces.
\begin{prop}
\label{prop:CR_Space}The following statements are true:

\renewcommand{\labelenumi}{(\alph{enumi})}
\begin{enumerate}
\item T4 spaces, especially metrizable spaces, are Tychonoff spaces. Moreover,
Tychonoff spaces are T3 spaces.
\item Subspaces of a Tychonoff space are Tychonoff spaces.
\item Any product space of Tychonoff spaces is a Tychonoff space.
\end{enumerate}
\end{prop}

Tychonoff space has close link to compactification. \label{Compactification}$S$
is called a compactification of $E$ (or $S$ compactifies $E$) if
$S$ is a compact \textit{Hausdorff}%
\footnote{As noted in \S \ref{sec:Base}, this work only considers Hausdorff
compactification.%
} space and $E$ is a dense \textit{subspace} of $S$. \label{Stone-Cech}$S$
is the Stone-$\check{\mbox{C}}$ech compactification of $E$ if $S$
compactifies $E$ and $\bigotimes C_{b}(E;\mathbf{R})$ extends to
a member of $\mathbf{imb}(S;\mathbf{R}^{C_{b}(E;\mathbf{R})})$. The
following proposition shows the equivalence of Tychonoff property,
general compactifiability and Stone-$\check{\mbox{C}}$ech compactifiability.
\begin{prop}
\label{prop:CR_Compactification}Let $E$ be a topological space.
Then, the following statements are equivalent:

\renewcommand{\labelenumi}{(\alph{enumi})}
\begin{enumerate}
\item $E$ has a compactification.
\item $E$ is a Tychonoff space.
\item $E$ has a unique Stone-$\check{\mbox{C}}$ech compactification up
to homeomorphism%
\footnote{``Unique up to homeomorphism'' means any two spaces with the relevant
property are homeomorphic.%
}.
\end{enumerate}
\end{prop}
We prove the proposition above by the following compactification result,
of which the Stone-$\check{\mbox{C}}$ech compactification is a special
case.
\begin{lem}
\label{lem:Compactification}Let $E$ be a topological space and $\mathcal{D}\subset\mathbf{R}^{E}$
be a collection of bounded functions. Then, the following statements
are equivalent:

\renewcommand{\labelenumi}{(\alph{enumi})}
\begin{enumerate}
\item $\mathcal{D}\subset C_{b}(E;\mathbf{R})$ separates and strongly separates
points on $E$.
\item $E$ admits a unique compactification $S$ up to homeomorphism such
that $\bigotimes\mathcal{D}$ extends to a homeomorphism between $S$
and the closure of $\bigotimes\mathcal{D}(E)$ in $\mathbf{R}^{\mathcal{D}}$.
\item $\bigotimes\mathcal{D}\in\mathbf{imb}(E;\mathbf{R}^{\mathcal{D}})$.
\end{enumerate}
\end{lem}

\begin{rem}
\label{rem:Metrizable_Compactification}If the $\mathcal{D}$ above
is countable, then the induced compactification has a homeomorph in
$\mathbf{R}^{\infty}$ and hence is metrizable. This is the foundation
of the replication bases.\end{rem}
\begin{cor}
\label{cor:M_Compactification}Let $E$ be a topological space. Then,
the following statements are equivalent:

\renewcommand{\labelenumi}{(\alph{enumi})}
\begin{enumerate}
\item $E$ is metrizable and separable.
\item There exists a countable subset of $\mathcal{D}\subset C_{b}(E;\mathbf{R})$
that separates and strongly separates points on $E$.
\item $E$ has a compactification that is homeomorphic to a compact subset
of $\mathbf{R}^{\infty}$.
\item $E$ admits a metrizable compactification.
\item $E$ is a dense subspace of some Polish space.
\end{enumerate}
\end{cor}

\label{One-Point}$S$ is the one-point compactification of $E$ if
$S$ compactifies $E$ and $S\backslash E$ is a singleton. Every
locally compact Hausdorff space is well-known to have a unique one-point
compactifiable up to Homeomorphism.
\begin{prop}
\label{prop:One-Point}Let $E$ be a locally compact space. Then,
the following statements are equivalent:

\renewcommand{\labelenumi}{(\alph{enumi})}
\begin{enumerate}
\item $E$ is a Hausdorff space.
\item $E$ has a unique one-point compactification up to homeomorphism.
\item $E$ is a Tychonoff space.
\end{enumerate}
\end{prop}

The next lemma is an analogue of Lemma \ref{lem:SP_on_Compact} for
locally compact spaces.
\begin{lem}
\label{lem:SP_on_LC}Let $E$ be a locally compact space and $\mathcal{D}\subset C_{0}(E;\mathbf{R})$.
Then, the following statements are equivalent:

\renewcommand{\labelenumi}{(\alph{enumi})}
\begin{enumerate}
\item $\mathcal{D}$ separates points on $E$.
\item $\mathcal{D}^{\Delta}\circeq\{f^{\Delta}\}_{f\in\mathcal{D}}\cup\{1\}$
is a subset of $C(E^{\Delta};\mathbf{R})$ that separates and strongly
separates points on $E^{\Delta}$, where $E^{\Delta}$ is a one-point
compactification of $E$ and $f^{\Delta}\circeq\mathfrak{var}(f;E^{\Delta},E,0)$%
\footnote{``$\mathfrak{var}(\cdot)$'' was defined in Notation \ref{notation:Var}.%
} for each $f\in\mathcal{D}$.
\item $E$ is a Hausdorff space and $\mathcal{D}$ strongly separates points
on $E$.
\end{enumerate}
\end{lem}

The local compactness of $E$ leads to the following point-separability
of $C_{0}(E;\mathbf{R})$.
\begin{prop}
\label{prop:LC_CR}Let $E$ be a locally compact space. Then, the
following statements are equivalent:

\renewcommand{\labelenumi}{(\alph{enumi})}
\begin{enumerate}
\item $E$ is a Hausdorff space.
\item $C_{c}(E;\mathbf{R})$ separates points on $E$.
\item $C_{c}(E;\mathbf{R})$ separates and strongly separates points on
$E$.
\item $C_{0}(E;\mathbf{R})$ separates and strongly separates points on
$E$.
\item $E$ is a Tychonoff space.
\end{enumerate}
\end{prop}

\section{\label{sec:Meas_Sep}Weak topology of finite Borel  measures}

The weak topology of $\mathcal{M}^{+}(E)$ is induced by $C_{b}(E;\mathbf{R})^{*}$%
\footnote{The notation ``$C_{b}(E;\mathbf{R})^{*}$'' was specified in \S
\ref{sec:Borel_Measure}.%
}. Hence, the Tychonoff properties of $\mathcal{M}^{+}(E)$ and $\mathcal{P}(E)$
are reduced to their Hausdorff properties.
\begin{prop}
\label{prop:P(E)_CR}Let $E$ be a topological space. Then, the following
statements are equivalent:

\renewcommand{\labelenumi}{(\alph{enumi})}
\begin{enumerate}
\item $\mathcal{M}^{+}(E)$ and $\mathcal{P}(E)$ are Tychonoff spaces.
\item $\mathcal{P}(E)$ is a Hausdorff space.
\item $C_{b}(E;\mathbf{R})^{*}$ separates points on $\mathcal{P}(E)$.
\end{enumerate}
\end{prop}

We now investigate the connection between the Tychonoff property of
$E$ and those of $\mathcal{M}^{+}(E)$ and $\mathcal{P}(E)$. On
one hand, we establish \cite[Theorem 2.1.4]{KX95} without the restriction
to Radon measures.
\begin{prop}
\label{prop:E_CR_P(E)_CR}Let $E$ be a Tychonoff space and $d\in\mathbf{N}$.
Then, $\mathcal{M}^{+}(E^{d})$ and $\mathcal{P}(E^{d})$ are Tychonoff
spaces and $\mathfrak{mc}[\Pi^{d}(C_{b}(E;\mathbf{R}))]$ is separating
on $E^{d}$%
\footnote{The terminology ``separating'' was introduced in \S \ref{sec:Borel_Measure}.%
}.
\end{prop}
To show the proposition above, we prepare a lemma that relates $E$'s
functional separabilities of points to those of $\mathcal{M}^{+}(E)$
and $\mathcal{P}(E)$.
\begin{lem}
\label{lem:Sep_Meas}Let $E$ be a topological space, $\mathcal{D}\subset M_{b}(E;\mathbf{R})$,
$d\in\mathbf{N}$ and $\mathcal{G}\circeq\mathfrak{mc}(\Pi^{d}(\mathcal{D}))$%
\footnote{As mentioned in Remark \ref{rem:Pi^d_Not_Need_1_in_D}, we need not
impose $1\in\mathcal{D}$ in Lemma \ref{lem:Sep_Meas} by the definition
of $\Pi^{d}(\mathcal{D})$.%
}. Then:

\renewcommand{\labelenumi}{(\alph{enumi})}
\begin{enumerate}
\item If $\mathcal{D}\subset C_{b}(E;\mathbf{R})$ strongly separates points
on $E$, then $\mathcal{G}^{*}$ separates points on $\mathcal{P}(E^{d})$
and $\mathcal{G}\cup\{1\}$ is separating on $E^{d}$.
\item If $\mathcal{D}$ is countable and strongly separates points on $E$,
then $\mathcal{G}^{*}$ separates and strongly separates points on
$\mathcal{P}(E^{d})$, and $\mathcal{G}\cup\{1\}$ is separating and
convergence determining on $E^{d}$.
\end{enumerate}
\end{lem}

On the other hand, we give an explicit example showing that the converse
of Proposition \ref{prop:E_CR_P(E)_CR} is not true.
\begin{example}
\label{exp:P(E)_CR}Example \ref{exp:Baseable_Space} (VII) and Example
\ref{exp:Baseable_Space_SB_2} (III) mentioned that $\mathbf{R}_{K}$
is a non-T3 (hence non-Tychonoff) topological refinement of $\mathbf{R}$
with $\mathscr{B}(\mathbf{R}_{K})=\mathscr{B}(\mathbf{R})$. $\mathcal{P}(\mathbf{R})$
is a Tychonoff space by Proposition \ref{prop:E_CR_P(E)_CR}. $\mathcal{P}(\mathbf{R}_{K})$
is a Hausdorff topological refinement of $\mathcal{P}(\mathbf{R})$
by Fact \ref{fact:Weak_Topo_Coarsen} (a) and Fact \ref{fact:Hausdorff_Refine}.
Thus, $\mathcal{P}(\mathbf{R}_{K})$ and $\mathcal{M}^{+}(\mathbf{R}_{K})$
are Tychonoff spaces by Proposition \ref{prop:P(E)_CR} (a, b).
\end{example}

The two examples below illustrate that the Hausdorff property of $E$
and those of $\mathcal{M}^{+}(E)$ and $\mathcal{P}(E)$ do not imply
each other, where the Borel measurability of singletons and the distinctiveness
of Dirac measures at distinct points play essential roles.
\begin{example}
\label{exp:P(E)_Hausdorff_1}Let $E=\{1,2,3,4\}$, $A\circeq\{1,2\}$,
$B\circeq\{3,4\}$ and equip $E$ with the topology $\mathscr{O}(E)\circeq\{\varnothing,A,B,E\}$.
Then, $\mathscr{B}(E)=\mathscr{C}(E)=\mathscr{O}(E)$ and singletons
in $E$ are neither closed nor Borel. So, $E$ is non-Hausdorff by
Proposition \ref{prop:Separability} (a). Letting%
\footnote{$\delta_{x}$ denotes the Dirac measure at $x$.%
} $\mu_{A}\circeq(\delta_{1}+\delta_{2})/2$ and $\mu_{B}=(\delta_{3}+\delta_{4})/2$,
we observe that $C_{b}(E;\mathbf{R})=\{a\mathbf{1}_{A}+b\mathbf{1}_{B}:a,b\in\mathbf{R}\}$
and $\mathcal{P}(E)=\{\mu_{a}\circeq a\mu_{A}+\left(1-a\right)\mu_{B}:a\in\left[0,1\right]\}$.
$C_{b}(E;\mathbf{R})^{*}$ separates points on $\mathcal{P}(E)$ since
for any $a_{1},a_{2}\in[0,1]$,
\begin{equation}
\int_{E}\mathbf{1}_{A}(x)\mu_{a_{1}}(dx)=a_{1}=a_{2}=\int_{E}\mathbf{1}_{A}(dx)\mu_{a_{2}}(dx)\label{eq:Check_P(E)_Hausdorff_A_B}
\end{equation}
implies $\mu_{a_{1}}=\mu_{a_{2}}$. Thus, $\mathcal{P}(E)$ and $\mathcal{M}^{+}(E)$
are Hausdorff by Proposition \ref{prop:P(E)_CR}.
\end{example}

\begin{example}
\label{exp:P(E)_Hausdorff_2}Due to the limit of space, we refer the
readers to \cite[\S 33, Exercise 11]{M00} for the non-trivial construction
of a topological space $E$ satisfying: (1) $E$ is a T3 (\textit{hence
Haudorff}) but non-Tychonoff space, and (2) there exist $a\neq b$
in $E$ such that $f^{*}(\delta_{a})=f(a)=f(b)=f^{*}(\delta_{b})$%
\footnote{The notation ``$f^{*}$'' was specified in \S \ref{sec:Borel_Measure}.%
} for all $f\in C(E;\mathbf{R})$. $\delta_{a}$ and $\delta_{b}$
are distinct measures by the Hausdorff property of $E$ and Proposition
\ref{prop:P(E)_Hausdorff} (a) below, but they can not be separated
by $C_{b}(E;\mathbf{R})^{*}$. So, neither $\mathcal{P}(E)$ nor $\mathcal{M}^{+}(E)$
is Hausdorff by Proposition \ref{prop:P(E)_CR}.
\end{example}

As long as the extreme non-Borel singletons are avoided, the Hausdorff
property of (the usually more complicated space) $\mathcal{P}(E)$
implies that of $E$.
\begin{prop}
\label{prop:P(E)_Hausdorff}Let $E$ be a topological space satisfying
$\{\{x\}:x\in E\}\subset\mathscr{B}(E)$ and $\mathcal{D}\subset M_{b}(E;\mathbf{R})$.
Then:

\renewcommand{\labelenumi}{(\alph{enumi})}
\begin{enumerate}
\item $\delta_{x}\neq\delta_{y}$ for any distinct $x,y\in E$.
\item If $\mathcal{D}^{*}$  separates points on $\mathcal{P}(E)$, then
$\mathcal{D}$ separates points on $E$.
\item If $\mathcal{P}(E)$ is a Hausdorff space, then $C_{b}(E;\mathbf{R})$
separates points on $E$ and $E$ is a Hausdorff space.
\end{enumerate}
\end{prop}

\begin{cor}
\label{cor:Metrizable_Separable_SepMeas}Let $E$ be a metrizable
and separable space. Then, there exists a countable $\mathcal{D}\subset C_{b}(E;\mathbf{R})$
satisfying the following:

\renewcommand{\labelenumi}{(\alph{enumi})}
\begin{enumerate}
\item $\mathcal{D}\ni1$ is closed under addition and multiplication.
\item $\mathcal{D}$ separates and strongly separates points on $E$.
\item $\mathcal{D}$ is separating and convergence determining on $E$.
\end{enumerate}
\end{cor}

\begin{cor}
\label{cor:Metrizable_Separable_P(E)}The following statements are
equivalent:

\renewcommand{\labelenumi}{(\alph{enumi})}
\begin{enumerate}
\item $E$ is a Tychonoff, and $\mathcal{M}^{+}(E)$ and $\mathcal{P}(E)$
are metrizable and separable spaces.
\item $E$ is a metrizable and separable space.
\end{enumerate}
\end{cor}

The properties of $\mathcal{M}^{+}$ and $\mathcal{P}(E)$ below are
vital.
\begin{thm}
\label{thm:P(E)_Compact_Polish}The following statements are true:

\renewcommand{\labelenumi}{(\alph{enumi})}
\begin{enumerate}
\item If $E$ is a compact Hausdorff space, then $\mathcal{P}(E)$ is also.
\item If $E$ is a Polish space, then $\mathcal{M}^{+}(E)$ and $\mathcal{P}(E)$
are also.
\end{enumerate}
\end{thm}

As noted in p.\pageref{Limit_Point}, the sequential concepts ``weak
limit point'' and ``relative compactness'' in $\mathcal{M}^{+}(E)$%
\footnote{Weak limit point and relative compactness of finite Borel measures
were reviewed in \S \ref{sec:Borel_Measure}.%
} may cause ambiguity in general, but one can get rid of that when
$E$ is a metriable space.
\begin{prop}
\label{prop:WLP_RC_Metrizable}Let $E$ be a topological space and
$\Gamma\subset\mathcal{M}^{+}(E)$. Then:

\renewcommand{\labelenumi}{(\alph{enumi})}
\begin{enumerate}
\item If $E$ is metrizable, then $\mathcal{M}^{+}(E)$ and $\mathcal{P}(E)$
are metrizable by the same metric.
\item If $\nu$ is a weak limit point of $\Gamma$, then $\nu$ is a limit
point of $\Gamma$ with respect to weak topology. The converse is
true when $E$ is metrizable.
\item If $E$ is metrizable, then the relative compactness of $\Gamma$
is equivalent to its precompactness with respect to weak topology.
\end{enumerate}
\end{prop}

The next lemma extends Theorem \ref{thm:Prokhorov} (b) to the non-probabilistic
case.
\begin{lem}
\label{lem:Seq_Prokhorov}Let $E$ be a Hausdorff space, $\Gamma\subset\mathcal{M}^{+}(E)$
be sequentially tight and $0<a<b$ satisfy $\{\mu(E)\}_{\mu\in\Gamma}\subset[a,b]$.
Then, $\Gamma$ is relatively compact and the total mass%
\footnote{The notion of total mass was specified in \S \ref{sub:Meas}.%
} of any weak limit point of $\Gamma$ lies in $[a,b]$.
\end{lem}

Morever, given a \textit{perfectly normal} (see e.g. \cite[\S 33, Exercise 6]{M00})
space $E$, we equate the Borel sets of $\mathcal{M}^{+}(E)$ generated
by its \textit{strong}%
\footnote{Strong topology of Borel probability measures was reviewed in Example
\ref{exp:Strong-topo}, Example \ref{exp:Baseable_Space} and Example
\ref{exp:Baseable_Space_SB_2}. It can be defined for finite Borel
measures in almost the same way.%
} and weak topologies, which generalizes \cite[Lemma 2.1]{BS89}.
\begin{lem}
\label{lem:Strong_Topo}Let $E$ be a perfectly normal (especially
metrizable or Polish) space. Then, $\mathscr{B}_{M_{b}(E;\mathbf{R})^{*}}(\mathcal{M}^{+}(E))=\mathscr{B}(\mathcal{M}^{+}(E))$.
\end{lem}

\section{\label{sec:SB}Standard Borel space}

This section reviews a few fundamental properties of standard Borel
spaces and standard Borel subsets.
\begin{fact}
\label{fact:Polish_SB}The following statements are true:

\renewcommand{\labelenumi}{(\alph{enumi})}
\begin{enumerate}
\item Borel isomorphs of standard Borel spaces are still standard Borel
spaces. In particular, Polish spaces, their Borel isomorphs and their
Borel subspaces%
\footnote{The notion of Borel subspace was introduced in Definition \ref{def:SB}.%
} are standard Borel spaces.
\item The cardinality of a standard Borel space can never exceed $\aleph(\mathbf{R})$.
\end{enumerate}
\end{fact}

Standard Borel spaces are Borel isomorphic to Borel subsets of Polish
spaces. The latter turns out to be precisely the metrizable Lusin
spaces.
\begin{prop}
\label{prop:Borel_in_Polish}Let $E$ be a metrizable space. Then,
the following statements are equivalent:

\renewcommand{\labelenumi}{(\alph{enumi})}
\begin{enumerate}
\item $E$ is a Lusin space.
\item $E$ has a Polish topological refinement $(E,\mathscr{U})$ with $\mathscr{B}(E)=\sigma(\mathscr{U})$.
\item $E$ is separable and for any metrization $(E,\mathfrak{r})$ of $E$,
$E$ is a Borel subset of the completion of $(E,\mathfrak{r})$.
\item $E$ is a Borel subspace of some Polish space.
\item There exist an $S\in\mathscr{C}(\mathbf{N}^{\infty})$ and a bijective
$f\in C(S;E)$.
\end{enumerate}
\end{prop}
The key to prove the equivalence above is the preservation of Borel
sets under bijective Borel measurable mappings. Here is a standard
result about this.
\begin{lem}
\label{lem:Polish_Lusin_Map}Let $E$ be a Lusin space, $S$ be a
Polish space, $f\in M(S;E)$ and%
\footnote{$\mathscr{U}_{f}$ is known as the ``push-forward topology of $f$''.
In any case, $f\in C(S;E,\mathscr{U}_{f})$.%
}
\begin{equation}
\mathscr{U}_{f}\circeq\left\{ O\subset E:f^{-1}(O)\in\mathscr{O}(S)\right\} .\label{eq:Push_Forward_Topo}
\end{equation}
Then:

\renewcommand{\labelenumi}{(\alph{enumi})}
\begin{enumerate}
\item If $f$ is continuous and bijective, then $f\in\mathbf{hom}(S;E,\mathscr{U}_{f})$
and $(E,\mathscr{U}_{f})$ is a Polish topological refinement of $E$.
\item If $E$ is metrizable (especially a Polish space) and $f$ is injective,
then $f(B)\in\mathscr{B}(E)$ for all $B\in\mathscr{B}(S)$.
\item If $E$ is metrizable (especially a Polish space) and $f$ is bijective,
then $f\in\mathbf{biso}(S;E)$%
\footnote{The notation ``$\mathbf{biso}$'' was defined in \S \ref{sub:Meas_Cont_Cadlag_Map}.%
}.
\end{enumerate}
\end{lem}

\begin{cor}
\label{cor:Lusin}Lusin spaces (resp. Souslin spaces) are precisely
the Hausdorff topological coarsenings%
\footnote{The terminology ``topological coarsening'' was specified in \S
\ref{sub:Topo}.%
} of Polish spaces (resp. Lusin spaces).
\end{cor}

From above we observe that a general (resp. metrizable) Lusin space
coarsens some Polish space topologically but does not necessarily
preserve (resp. does preserve) its Borel $\sigma$-algebra. By contrast,
the next proposition shows that a general standard Borel space is
a topological \textit{variant} (not necessarily a coarsening or refinement)
of some Polish space that preserves its Borel $\sigma$-algbera.
\begin{prop}
\label{prop:SB}Let $E$ be a topological space. Then, the following
statements are equivalent:

\renewcommand{\labelenumi}{(\alph{enumi})}
\begin{enumerate}
\item $E$ is a standard Borel space.
\item $E$ is Borel isomorphic to some metrizable Lusin space.
\item There exists a topology $\mathscr{U}_{1}$ on $E$ such that $(E,\mathscr{U}_{1})$
is a metrizable Lusin space and $\mathscr{B}(E)=\sigma(\mathscr{U}_{1})$.
\item There exists a topology $\mathscr{U}_{2}$ on $E$ such that $(E,\mathscr{U}_{2})$
is a Polish space and $\mathscr{B}(E)=\sigma(\mathscr{U}_{2})$.
\item $E$ is Borel isomorphic to some Polish space.
\end{enumerate}
\end{prop}

Given metrizability, standard Borel and Lusin properties become indifferent.
\begin{prop}
\label{prop:Metrizable_SB_Lusin}Let $E$ be a metrizable space. Then,
the following statements are equivalent:

\renewcommand{\labelenumi}{(\alph{enumi})}
\begin{enumerate}
\item $E$ is a standard Borel space.
\item $E$ is a Lusin space.
\item $E$ admits a Polish topological refinement $(E,\mathscr{U})$ satisfying
$\mathscr{B}(E)=\sigma(\mathscr{U})$.
\end{enumerate}
\end{prop}
Our proof is based on the following interesting result which illustrates
that a Borel measurable mapping may preserve some topological properties.
\begin{lem}
\label{lem:Separable_SB_Image}Let $E$ be a standard Borel space
and $S$ be a metrizable space. Then:

\renewcommand{\labelenumi}{(\alph{enumi})}
\begin{enumerate}
\item If there is a bijective $f\in M(E;S)$, then $S$ is separable.
\item If $E$ is metrizable, then $E$ is separable (hence second-countable).
\end{enumerate}
\end{lem}

The next proposition compares standard Borel and Borel subsets which
are likely to be different in general topological spaces.
\begin{prop}
\label{prop:SB_Borel}Let $E$ be a topological space. Then:

\renewcommand{\labelenumi}{(\alph{enumi})}
\begin{enumerate}
\item If $A\in\mathscr{B}^{\mathbf{s}}(E)$, then $\mathscr{B}_{E}(A)\subset\mathscr{B}^{\mathbf{s}}(E)$.
In particular, if $E$ is a standard Borel space, then $\mathscr{B}(E)\subset\mathscr{B}^{\mathbf{s}}(E)$.
\item If $E$ is a metrizable standard Borel space, especially if $E$ is
a Polish space, then $\mathscr{B}(E)=\mathscr{B}^{\mathbf{s}}(E)$.
\end{enumerate}
\end{prop}

If $E$ is compact and $S$ is Hausdorff, then any bijective $f\in C(E;S)$
belongs to $\mathbf{hom}(E;S)$ and $S$ is also compact (see \cite[Theorem 26.6]{M00}).
One of Kuratowski's theorems provides a similar identification for
bijective Borel measurable mappings from standard Borel spaces to
metrizable spaces. Herein, we give a short proof for integrity.
\begin{prop}
\label{prop:SB_Map}Let $E$ be a topological space, $S$ be a metrizable
space, $f\in M(E;S)$ be injective and $A\in\mathscr{B}^{\mathbf{s}}(E)$.
Then, $f(A)\in\mathscr{B}^{\mathbf{s}}(S)$ and $f|_{A}\in\mathbf{biso}(A;f(A))$.
\end{prop}

\begin{rem}
\label{rem:SB_Map}The $E$ in Lemma \ref{lem:Polish_Lusin_Map} (b,
c) is standard Borel by Proposition \ref{prop:Metrizable_SB_Lusin}
(a, b) and the Polish space $S$ in Lemma \ref{lem:Polish_Lusin_Map}
satisfies $\mathscr{B}(S)=\mathscr{B}^{\mathbf{s}}(S)$ by Proposition
\ref{prop:SB_Borel} (b). Hence, Proposition \ref{prop:SB_Map} generalizes
Lemma \ref{lem:Polish_Lusin_Map} (b, c).
\end{rem}

The next proposition is about the functional separabilities of points
and probability measures on standard Borel spaces.
\begin{prop}
\label{prop:SB_SP}Let $E$ be a standard Borel space. Then:

\renewcommand{\labelenumi}{(\alph{enumi})}
\begin{enumerate}
\item There exists a countable susbet of $M_{b}(E;\mathbf{R})$ that separates
points on $E$.
\item If $\mathcal{D}\subset\mathbf{R}^{E}$ satisfies $\mathscr{B}_{\mathcal{D}}(E)=\mathscr{B}(E)$,
then $\mathcal{D}$ separates points on $E$.
\item If $\mathcal{D}\subset M(E;\mathbf{R})$ is countable and separates
points on $E$, then $\sigma(\mathcal{D})=\mathscr{B}_{\mathcal{D}}(E)=\mathscr{B}(E)$.
\end{enumerate}
\end{prop}

\section{\label{sec:Sko}Skorokhod $\mathscr{J}_{1}$-space}

We start with several most essential properties of $D(\mathbf{R}^{+};E)$.
\begin{prop}
\label{prop:Sko_Basic_1}Let $E$ and $S$ be Tychonoff spaces. Then:

\renewcommand{\labelenumi}{(\alph{enumi})}
\begin{enumerate}
\item If $\mathcal{D}\subset C(E;\mathbf{R})$ strongly separates points
on $E$ (especially $\mathcal{D}=C(E;\mathbf{R})$), then $\{\varpi(f):f\in\mathfrak{ae}(\mathcal{D})\}\subset D(\mathbf{R}^{+};\mathbf{R})^{D(\mathbf{R}^{+};E)}$%
\footnote{``$\varpi(f)$'' and ``$\varpi(\mathcal{D})$'' were defined in
\S \ref{sub:Map}, while ``$\mathfrak{ae}(\cdot)$'' was defined
in \S \ref{sub:Fun}.%
} satisfies%
\footnote{$\mathscr{J}(E)$ denotes the Skorokhod $\mathscr{J}_{1}$-topology
of $D(\mathbf{R}^{+};E)$. %
}
\begin{equation}
\varpi[\mathfrak{ae}(\mathcal{D})]\in\mathbf{imb}\left(D(\mathbf{R}^{+};E);D(\mathbf{R}^{+};\mathbf{R})^{\mathfrak{ae}(\mathcal{D})}\right)\label{eq:Sko_Prod_Imb}
\end{equation}
and
\begin{equation}
\mathscr{J}(E)=\mathscr{O}_{\left\{ \varpi(f):f\in\mathfrak{ae}(\mathcal{D})\right\} }\left(D(\mathbf{R}^{+};E)\right).\label{eq:Sko_Topo_SSP}
\end{equation}

\item If $\mathcal{D}\subset C_{b}(E;\mathbf{R})$ separates points on $E$,
then $\{\alpha_{t,n}^{f}:f\in\mathcal{D},t\in\mathbf{Q}^{+},n\in\mathbf{N}\}$%
\footnote{$\mathbf{Q}^{+}$ denotes non-negative rational numbers.%
} is a subset of $C(D(\mathbf{R}^{+};E);\mathbf{R})$ separating points
on $D(\mathbf{R}^{+};E)$ with each $\alpha_{t,n}^{f}$ defined as
in (\ref{eq:Sko_SP_Fun}).
\item $D(\mathbf{R}^{+};E)$ is a Tychonoff space.
\item If $f\in\mathcal{D}\subset C(S;E)$, then
\begin{equation}
\varpi(f)\in C\left(D(\mathbf{R}^{+};S);D(\mathbf{R}^{+};E)\right)\label{eq:Path_Mapping_Sko_Cont}
\end{equation}
and
\begin{equation}
\varpi(\mathcal{D})\subset C\left(D(\mathbf{R}^{+};S);D(\mathbf{R}^{+};E)^{\mathcal{D}}\right).\label{eq:Path_Mapping_Class_Cont}
\end{equation}

\item If $E$ is a topological coarsening of $S$, then $D(\mathbf{R}^{+};S)\subset D(\mathbf{R}^{+};E)$
and $\mathscr{J}(S)\supset\mathscr{O}_{D(\mathbf{R}^{+};E)}(D(\mathbf{R}^{+};S))$.
\end{enumerate}
\end{prop}

\begin{rem}
\label{rem:Sko_Pseudometric}By Proposition \ref{prop:Sko_Basic_1}
(a), $\mathscr{J}(E)$ is uniquely determined by $\mathscr{O}(E)$
and does not depend on the choice of the pseudometrics in its definition.\end{rem}
\begin{cor}
\label{cor:Sko_Subspace}Let $E$ be a Tychonoff space. Then, $D(\mathbf{R}^{+};A,\mathscr{O}_{E}(A))$
is a topological subspace of $D(\mathbf{R}^{+};E)$ for any non-empty
$A\subset E$.
\end{cor}

\begin{cor}
\label{cor:D(R)^Inf_Imb}The one-dimensional projections $\mathcal{J}=\{\mathfrak{p}_{i}\}_{i\in\mathbf{I}}$
on $\mathbf{R}^{\mathbf{I}}$ satisfy
\begin{equation}
\varpi\left[\mathfrak{ae}(\mathcal{J})\right]\in\mathbf{imb}\left[D(\mathbf{R}^{+};\mathbf{R}^{\mathbf{I}});D(\mathbf{R}^{+};\mathbf{R})^{\mathfrak{ae}(\mathcal{J})}\right]\label{eq:D(R)^Inf_Imb}
\end{equation}
and%
\footnote{$\mathbf{R}^{\mathfrak{ae}(\mathcal{J})}$ is a Tychonoff space by
Proposition \ref{prop:CR_Space} (c), so $D(\mathbf{R}^{+};\mathbf{R}^{\mathfrak{ae}(\mathcal{J})})$
satisfies our definition of Skorokhod $\mathscr{J}_{1}$-spaces in
\S \ref{sub:Meas_Cont_Cadlag_Map}.%
}
\begin{equation}
\varpi\left[\bigotimes\mathfrak{ae}(\mathcal{J})\right]\in C\left[D(\mathbf{R}^{+};\mathbf{R}^{\mathbf{I}});D\left(\mathbf{R}^{+};\mathbf{R}^{\mathfrak{ae}(\mathcal{J})}\right)\right].\label{eq:D(R^Inf)_Cont}
\end{equation}

\end{cor}

The property below is well-known for compact subsets of $D(\mathbf{R}^{+};E)$.
\begin{prop}
\label{prop:Sko_Compact_Containment}Let $E$ be a Tychonoff space
and $K\in\mathscr{K}(D(\mathbf{R}^{+};E))$. Then, there exist $\{K_{n}\}_{n\in\mathbf{N}}\subset\mathscr{K}(E)$
such that
\begin{equation}
K\subset\bigcap_{n\in\mathbf{N}}\left\{ x\in D(\mathbf{R}^{+};E):x(t)\in K_{n},\forall t\in[0,n)\right\} .\label{eq:Path_CCC}
\end{equation}

\end{prop}

The next lemma is about the finite-dimensional projections on $D(\mathbf{R}^{+};E)$.
\begin{lem}
\label{lem:Sko_Proj}Let $E$ be a Tychonoff space and $\mathbf{T}_{0}\in\mathscr{P}_{0}(\mathbf{R}^{+})$.
Then:

\renewcommand{\labelenumi}{(\alph{enumi})}
\begin{enumerate}
\item $\mathfrak{p}_{\mathbf{T}_{0}}\in M(D(\mathbf{R}^{+};E);E^{\mathbf{T}_{0}},\mathscr{B}(E)^{\otimes\mathbf{T}_{0}})$%
\footnote{Herein, $\mathfrak{p}_{\mathbf{T}_{0}}$ denotes the projection on
$E^{\mathbf{R}^{+}}$ for $\mathbf{T}_{0}\subset\mathbf{R}^{+}$ restricted
to $D(\mathbf{R}^{+};E)$.%
}.
\item (\ref{eq:Sko_Borel>Prod}) and (\ref{eq:Sko_Meas}) hold.
\item $\mathfrak{p}_{\mathbf{T}_{0}}$ is continuous at $x\in D(\mathbf{R}^{+};E)$
whenever $\mathbf{T}_{0}\subset\mathbf{R}^{+}\backslash J(x)$.
\end{enumerate}
\end{lem}

\begin{cor}
\label{cor:Sko_Meas_FDD}Let $E$ be a Tychonoff space. Then, $\mu\circ\mathfrak{p}_{\mathbf{T}_{0}}^{-1}$
is a member of $\mathfrak{M}^{+}(E^{\mathbf{T}_{0}},\mathscr{B}(E)^{\otimes\mathbf{T}_{0}})$
for all $\mu\in\mathfrak{M}^{+}(D(\mathbf{R}^{+};E),\mathscr{B}(E)^{\otimes\mathbf{R}^{+}}|_{D(\mathbf{R}^{+};E)})$
(especially $\mu\in\mathcal{M}^{+}(D(\mathbf{R}^{+};E))$) and $\mathbf{T}_{0}\in\mathscr{P}_{0}(\mathbf{R}^{+})$.
\end{cor}

The following result is about the measurability of $w_{\mathfrak{r},\delta,T}^{\prime}$
on $D(\mathbf{R}^{+};E)$.
\begin{prop}
\label{prop:Sko_w'}Let $E$ be a Tychonoff space and $\delta,T\in(0,\infty)$.
Then:

\renewcommand{\labelenumi}{(\alph{enumi})}
\begin{enumerate}
\item $w_{\mathfrak{r},\delta,T}^{\prime}\in M(D(\mathbf{R}^{+};E);\mathbf{R})$%
\footnote{The notation ``$w_{\mathfrak{r},\delta,T}^{\prime}$'' was defined
in \S \ref{sub:Map}. The notation ``$w_{\rho_{\{f\}},\delta,T}^{\prime}$''
(resp. ``$w_{\left|\cdot\right|,\delta,T}^{\prime}$'') is defined
by (\ref{eq:w'}) with $\mathfrak{r}=\rho_{\{f\}}$ (resp. $E=\mathbf{R}$
and $\mathfrak{r}=\left|\cdot\right|$).%
} if $E$ allows a metrization $(E,\mathfrak{r})$.
\item $w_{\rho_{\{f\}},\delta,T}^{\prime}\in M(D(\mathbf{R}^{+};E);\mathbf{R})$
for all $f\in C(E;\mathbf{R})$.
\end{enumerate}
\end{prop}

The next fact discusses the measurability issue in (\ref{eq:J(Mu)}).
\begin{fact}
\label{fact:J(Mu)_Well_Defined}Let $E$ be a Tychonoff space. If
$M(E;\mathbf{R})$ has a countable subset separating points on $E$,
then
\begin{equation}
\left\{ x\in D(\mathbf{R}^{+};E):t\in J(x)\right\} \in\mathscr{B}(E)^{\otimes\mathbf{R}^{+}}|_{D(\mathbf{R}^{+};E)},\;\forall t\in\mathbf{R}^{+}\label{eq:J(Mu)_Well_Defined}
\end{equation}
and $J(\mu)$%
\footnote{$J(\mu)$, the set of fixed left-jump times of $\mu$, was defined
in (\ref{eq:J(Mu)}).%
} is well-defined for all $\mu\in\mathfrak{M}^{+}(D(\mathbf{R}^{+};E),\mathscr{B}(E)^{\otimes\mathbf{R}^{+}}|_{D(\mathbf{R}^{+};E)})$
(especially $\mu\in\mathcal{M}^{+}(D(\mathbf{R}^{+};E))$).
\end{fact}

The next proposition discusses the metrizability of $D(\mathbf{R}^{+};E)$.
\begin{prop}
\label{prop:Sko_Basic_2}Let $E$ be a metrizable space. Then:

\renewcommand{\labelenumi}{(\alph{enumi})}
\begin{enumerate}
\item If $(E,\mathfrak{r})$ is a metrization of $E$, then $D(\mathbf{R}^{+};E)$
is metrized by $\varrho^{\mathfrak{r}}$.
\item If $E$ is separable (espeically a Polish space), then $D(\mathbf{R}^{+};E)$
is also separable and $\mathscr{B}(D(\mathbf{R}^{+};E))=\mathscr{B}(E)^{\otimes\mathbf{R}^{+}}|_{D(\mathbf{R}^{+};E)}$.
\item If a metric $\mathfrak{r}$ completely metrizes $E$, then $\varrho^{\mathfrak{r}}$
completely metrizes $D(\mathbf{R}^{+};E)$.
\item If $E$ is a Polish space, then $D(\mathbf{R}^{+};E)$ is also.
\end{enumerate}
\end{prop}

The following three results relate convergence in $D(\mathbf{R}^{+};E)$
to that in the finite time horizon case. These results are contained
explicitly or implicitly in standard texts like \cite{JS03} and \cite{J86}.
Herein, we reestablish them for readers' convenience.
\begin{prop}
\label{prop:D_Conv_D[0,u]_Conv_Metric}Let $E$ be metrizable, $\{y_{k}\}_{k\in\mathbf{N}_{0}}\subset D(\mathbf{R}^{+};E)$%
\footnote{$\mathbf{N}_{0}$ denotes the non-negative integers.%
} and%
\footnote{Our notation $y_{k}^{u}\circeq\mathfrak{var}(x;[0,u+1],[0,u],x(u))$
represents the piecewise function defined by $y_{k}^{u}(t)\circeq y_{k}(t)$
for all $t\in[0,u]$ and $y_{k}^{u}(t)\circeq y_{k}(u)$ for all $t\in(u,u+1]$.%
}
\begin{equation}
y_{k}^{u}\circeq\mathfrak{var}\left(y_{k};[0,u+1],[0,u],y_{k}(u)\right),\;\forall k\in\mathbf{N}_{0},u\in(0,\infty).\label{eq:yk_u}
\end{equation}
Then for each $u\in\mathbf{R}^{+}\backslash J(y_{0})$,
\begin{equation}
y_{k}\longrightarrow y_{0}\mbox{\mbox{ as }k\ensuremath{\uparrow\infty}}\mbox{ in }D(\mathbf{R}^{+};E)\label{eq:yk_Conv_y0}
\end{equation}
implies%
\footnote{As mentioned in \S \ref{sub:Meas_Cont_Cadlag_Map}, $D([0,u];E)$
denotes the Skorokhod $\mathscr{J}_{1}$-space of all c$\grave{\mbox{a}}$dl$\grave{\mbox{a}}$g
mappings from $[0,u]$ to $E$.%
}
\begin{equation}
y_{k}^{u}\longrightarrow y_{0}^{u}\mbox{ as }k\uparrow\infty\mbox{ in }D([0,u+1];E).\label{eq:yk_u_Conv_y0_u}
\end{equation}

\end{prop}

\begin{lem}
\label{lem:D_Conv_D[0,u]_Conv}Let $E$ be a Tychonoff space, $\mathcal{D}\subset C(E;\mathbf{R})$
strongly separate points on $E$, $\Psi\circeq\varpi[\mathfrak{ae}(\mathcal{D})]$,
$\{y_{k}\}_{k\in\mathbf{N}_{0}}\subset D(\mathbf{R}^{+};E)$, $u\in\mathbf{R}^{+}\backslash J(y_{0})$
and $\{y_{k}^{u}\}_{k\in\mathbf{N}_{0}}$ be as in (\ref{eq:yk_u}).
Then,
\begin{equation}
\Psi(y_{k})\longrightarrow\Psi(y_{0})\mbox{ as }k\uparrow\infty\mbox{ in }D(\mathbf{R}^{+};\mathbf{R})^{\mathfrak{ae}(\mathcal{D})}\label{eq:Psi(yk)_Conv_Psi(y0)}
\end{equation}
implies (\ref{eq:yk_u_Conv_y0_u}). In particular, (\ref{eq:yk_Conv_y0})
implies (\ref{eq:yk_u_Conv_y0_u}).
\end{lem}

\begin{lem}
\label{lem:D[0,u]_Conv_D_Conv}Let $E$ be a Tychonoff space, $\{y_{k}\}_{k\in\mathbf{N}_{0}}\subset D(\mathbf{R}^{+};E)$
and $y_{k}^{u}$ be as in (\ref{eq:yk_u}) for each $k\in\mathbf{N}_{0}$
and $u\in(0,\infty)$. If $\mathbf{R}^{+}\backslash J(y_{0})$ is
dense in $\mathbf{R}^{+}$, and if (\ref{eq:yk_u_Conv_y0_u}) holds
for all $u\in\mathbf{R}^{+}\backslash J(y_{0})$, then (\ref{eq:yk_Conv_y0})
holds.
\end{lem}

\section{\label{sec:Cadlag}C$\grave{\mbox{a}}$dl$\grave{\mbox{a}}$g process}

The following two facts compare $E$-valued c$\grave{\mbox{a}}$dl$\grave{\mbox{a}}$g
processes and $D(\mathbf{R}^{+};E)$-valued random variables.
\begin{fact}
\label{fact:Sko_RV_Cadlag}Let $E$ be a Tychonoff space and $\mu\in\mathcal{M}^{+}(D(\mathbf{R}^{+};E))$
be the distribution of $D(\mathbf{R}^{+};E)$-valued random variable
$X$. Then:

\renewcommand{\labelenumi}{(\alph{enumi})}
\begin{enumerate}
\item $X$ is an $E$-valued c$\grave{\mbox{a}}$dl$\grave{\mbox{a}}$g
process
\item $\mu$ equals the restriction of $\mathrm{pd}(X)|_{D(\mathbf{R}^{+};E)}$
to $\sigma(\mathscr{J}(E))$%
\footnote{Restriction of measure to sub-$\sigma$-algebra and $X$'s process
distribution $\mathrm{pd}(X)$ were specified in \S \ref{sub:Meas}
and \S \ref{sec:Proc} respectively.%
}.
\item The finite-dimensional distribution of $X$ for each $\mathbf{T}_{0}\in\mathscr{P}_{0}(\mathbf{R}^{+})$
is $\mu\circ\mathfrak{p}_{\mathbf{T}_{0}}^{-1}$.
\end{enumerate}
\end{fact}

\begin{fact}
\label{fact:Cadlag_Sko_RV}Let $E$ be a metrizable and separable
space and $(\Omega,\mathscr{F},\mathbb{P};X)$%
\footnote{We arranged in \S \ref{sec:Convention} that $(\Omega,\mathscr{F},\mathbb{P})$
and $\{(\Omega^{i},\mathscr{F}^{i},\mathbb{P}^{i})\}_{i\in\mathbf{I}}$
are complete probability spaces. Completeness of measure space was
specified in \ref{sub:Meas}.%
} be an $E$-valued process. Then:

\renewcommand{\labelenumi}{(\alph{enumi})}
\begin{enumerate}
\item If all paths of $X$ lie in $D(\mathbf{R}^{+};E)$, then $X\in M(\Omega,\mathscr{F};D(\mathbf{R}^{+};E))$.
\item If $X$ is c$\grave{\mbox{a}}$dl$\grave{\mbox{a}}$g, then there
exists a $Y\in M(\Omega,\mathscr{F};D(\mathbf{R}^{+};E))$ that is
indistinguishable from $X$.
\end{enumerate}
In particular, the statements above are true when $E$ is a Polish
space.
\end{fact}

The next lemma solves the measurability issue in (\ref{eq:CCC_Measurability})
under mild conditions.
\begin{lem}
\label{lem:CCC_Measurability}Let $E$ be a Hausdorff space, $V$
be the family of all c$\grave{\mbox{a}}$dl$\grave{\mbox{a}}$g members
of $E^{\mathbf{R}^{+}}$%
\footnote{$E$ need not be a Tychonoff space, so we avoid the notation $D(\mathbf{R}^{+};E)$
for clarity.%
}, $(\Omega,\mathscr{F},\mathbb{P};X)$ be an $E$-valued c$\grave{\mbox{a}}$dl$\grave{\mbox{a}}$g
process and $T\in(0,\infty)$. Then, $\bigcap_{t\in[0,T]}X_{t}^{-1}(A)\in\mathscr{F}$
for all $A\in\mathscr{C}(E)$, especially for all $A\in\mathscr{K}(E)$
when $E$ is a Hausdorff space.
\end{lem}

The next lemma treats the measurability issue in (\ref{eq:w'_Measurable}).
\begin{lem}
\label{lem:MCC_Measurability}Let $E$ be a topological space, $(\Omega,\mathscr{F},\mathbb{P};X)$
be an $E$-valued c$\grave{\mbox{a}}$dl$\grave{\mbox{a}}$g process,
$\mathfrak{r}$ be a pseudometric on $E$ and $\delta,T\in(0,\infty)$.
Then, $w_{\mathfrak{r},\delta,T}^{\prime}\circ X\in M(\Omega,\mathscr{F};\mathbf{R})$
in each of the following settings:

\renewcommand{\labelenumi}{(\alph{enumi})}
\begin{enumerate}
\item $(E,\mathfrak{r})$ is a metric space and $X\in M(\Omega,\mathscr{F};D(\mathbf{R}^{+};E))$.
\item $(E,\mathfrak{r})$ is a separable metric space.
\item $\mathfrak{r}=\rho_{\{f\}}$ with $f\in C(E;\mathbf{R})$.
\item $\mathfrak{r}=\rho_{\mathcal{D}}$ with $\mathcal{D}\subset C(E;\mathbf{R})$
being countable and separating points on $E$%
\footnote{The $E$ in (d) is a $\mathcal{D}$-baseable space.%
}.
\end{enumerate}
\end{lem}

The following five results discuss the relationship among $\mathfrak{r}$-MCC%
\footnote{$\mathfrak{r}$-MCC, MCC, $\mathcal{D}$-FMCC and WMCC were introduced
in Definition \ref{def:Proc_Reg}.%
}, MCC, $\mathcal{D}$-FMCC and WMCC for c$\grave{\mbox{a}}$dl$\grave{\mbox{a}}$g
processes.
\begin{fact}
\label{fact:MCC_1}Let $E$ be a topological space, $\mathcal{D}\subset M(E;\mathbf{R})$
and $\{X^{i}\}_{i\in\mathbf{I}}$ be $E$-valued processes such that
$\{\varpi(f)\circ X^{i}\}_{f\in\mathcal{D},i\in\mathbf{I}}$ are all
c$\grave{\mbox{a}}$dl$\grave{\mbox{a}}$g. Then:

\renewcommand{\labelenumi}{(\alph{enumi})}
\begin{enumerate}
\item $\{X^{i}\}_{i\in\mathbf{I}}$ satisfies $\rho_{\{f\}}$-MCC for all
$f\in\mathcal{D}$ if and only if $\{\varpi(f)\circ X^{i}\}_{i\in\mathbf{I}}$
satisfies $\left|\cdot\right|$-MCC for all $f\in\mathcal{D}$.
\item $\{X^{i}\}_{i\in\mathbf{I}}$ satisfies $\mathcal{D}$-FMCC if and
only if $\{\varpi(f)\circ X^{i}\}_{i\in\mathbf{I}}$ satisfies $\left|\cdot\right|$-MCC
for all $f\in\mathfrak{ae}(\mathcal{D})$.
\end{enumerate}
In particular, the two statements above are true when $\mathcal{D}\subset C(E;\mathbf{R})$
and $\{X^{i}\}_{i\in\mathbf{I}}$ are all c$\grave{\mbox{a}}$dl$\grave{\mbox{a}}$g.
\end{fact}

The next result is a version of \cite[Proposition 14]{K15} for infinite
time horizon.
\begin{prop}
\label{prop:MCE}Let $E$ be a Hausdorff space and $\{(\Omega^{i},\mathscr{F}^{i},\mathbb{P}^{i};X^{i})\}_{i\in\mathbf{I}}$
be $E$-valued c$\grave{\mbox{a}}$dl$\grave{\mbox{a}}$g processes.
Then, the following statements are equivalent:

\renewcommand{\labelenumi}{(\alph{enumi})}
\begin{enumerate}
\item $\{X^{i}\}_{i\in\mathbf{I}}$ satisfies MCC.
\item There exist a $\mathcal{D}_{1}\subset C(E;\mathbf{R})$ and a $\mathcal{D}_{2}\subset C_{b}(E;\mathbf{R})$
such that: (1) $\{X^{i}\}_{i\in\mathbf{I}}$ satisfies $\mathcal{D}_{1}$-FMCC,
(2) $\mathcal{D}_{2}=\mathfrak{ac}(\mathcal{D}_{2})$%
\footnote{$\mathcal{D}_{2}=\mathfrak{ac}(\mathcal{D}_{2})$ means $\mathcal{D}_{2}$
is closed under multiplication.%
} strongly separates points on $E$, and (3) for any $g\in\mathcal{D}_{2}$
and $\epsilon,T>0$, there exists an $f_{g,\epsilon,T}\in\mathcal{D}_{1}$
satisfying
\begin{equation}
\sup_{i\in\mathbf{I}}\mathbb{P}^{i}\left(\sup_{t\in[0,T]}\left|f_{g,\epsilon,T}\circ X_{t}^{i}-g\circ X_{t}^{i}\right|\geq\epsilon\right)\leq\epsilon.\label{eq:MCE-2}
\end{equation}

\item There exist $\mathcal{D}\subset C_{b}(E;\mathbf{R})$ and $\{(\Omega^{i},\mathscr{F}^{i},\mathbb{P}^{i};\zeta^{i,f,\epsilon,T})\}_{i\in\mathbf{I},f\in\mathcal{D},\epsilon,T>0}$
such that: (1) $\mathcal{D}=\mathfrak{ac}(\mathcal{D})$ strongly
separates points on $E$, and (2) for each $f\in\mathfrak{ae}(\mathcal{D})$
and $\epsilon,T>0$, $\mathbf{R}$-valued processes $\{\zeta^{i,f,\epsilon,T}\}_{i\in\mathbf{I}}$
satisfy $\left|\cdot\right|$-MCC%
\footnote{$\left|\cdot\right|$-MCC means MCC for the Euclidean norm metric
$\left|\cdot\right|$.%
} and
\begin{equation}
\sup_{i\in\mathbf{I}}\mathbb{P}^{i}\left(\sup_{t\in[0,T]}\left|f\circ X_{t}^{i}-\zeta_{t}^{i,f,\epsilon,T}\right|\geq\epsilon\right)\leq\epsilon.\label{eq:MCE-1}
\end{equation}

\item There exists an $\mathcal{D}\subset C(E;\mathbf{R})$ such that: (1)
$\mathcal{D}$ strongly separates points on $E$, and (2) $\{\varpi(\bigotimes\mathcal{D}_{0})\circ X^{i}\}_{i\in\mathbf{I}}$
satisfies $\left|\cdot\right|$-MCC and MCCC%
\footnote{The notion of MCCC was specified in Definition \ref{def:Proc_Reg}.%
} for all $\mathcal{D}_{0}\in\mathscr{P}_{0}(\mathcal{D})$.
\item There exists a $\mathcal{D}\subset C(E;\mathbf{R})$ such that: (1)
$\mathcal{D}$ strongly separates points on $E$, and (2) $\{\varpi(g)\circ X^{i}\}_{i\in\mathbf{I}}$
satisfies $\left|\cdot\right|$-MCC and MCCC for all $g\in\mathfrak{ac}(\{af:f\in\mathcal{D},a\in\mathbf{R}\})$%
\footnote{The notation ``$\mathfrak{ac}(\cdot)$'' was defined in \S \ref{sub:Fun}.%
}.
\end{enumerate}
\end{prop}

\begin{cor}
\label{cor:MCC_2}Let $E$ be a Hausdorff space and $\{(\Omega^{i},\mathscr{F}^{i},\mathbb{P}^{i};X^{i})\}_{i\in\mathbf{I}}$
be $E$-valued c$\grave{\mbox{a}}$dl$\grave{\mbox{a}}$g processes.
Then, the following statements are equivalent:

\renewcommand{\labelenumi}{(\alph{enumi})}
\begin{enumerate}
\item $\{X^{i}\}_{i\in\mathbf{I}}$ satisfies MCC.
\item $\{X^{i}\}_{i\in\mathbf{I}}$ satisfies $\mathcal{D}$-FMCC for some
$\mathcal{D}=\mathfrak{ac}(\mathcal{D})\subset C_{b}(E;\mathbf{R})$
and $\mathcal{D}$ strongly separates points on $E$.
\item $\{X^{i}\}_{i\in\mathbf{I}}$ satisfies $\mathcal{D}$-FMCC for some
$\mathcal{D}=\mathfrak{ac}(\mathcal{D})\subset C(E;\mathbf{R})$ and
$\mathcal{D}$ strongly separates points on $E$.
\item $\{X^{i}\}_{i\in\mathbf{I}}$ satisfies $\mathcal{D}$-FMCC for some
$\mathcal{D}\subset C(E;\mathbf{R})$ and $\mathcal{D}$ strongly
separates points on $E$.
\end{enumerate}
\end{cor}

\begin{prop}
\label{prop:MCC_3}Let $E$ be a topological space, $\{X^{i}\}_{i\in\mathbf{I}}$
be $E$-valued processes and $\mathcal{D}\subset M(E;\mathbf{R})$
be countable and separate points on $E$. If $\{\varpi(f)\circ X^{i}\}_{f\in\mathcal{D},i\in\mathbf{I}}$
are all c$\grave{\mbox{a}}$dl$\grave{\mbox{a}}$g, especially if
$\{X^{i}\}_{i\in\mathbf{I}}$ are all c$\grave{\mbox{a}}$dl$\grave{\mbox{a}}$g
and $\mathcal{D}\subset C(E;\mathbf{R})$, then $\{X^{i}\}_{i\in\mathbf{I}}$
satisfying $\rho_{\mathcal{D}}$-MCC is equivalent to $\{\varpi(f)\circ X^{i}\}_{i\in\mathbf{I}}$
satisfying $\left|\cdot\right|$-MCC for all $f\in\mathcal{D}$.
\end{prop}

\begin{fact}
\label{fact:MCC_4}Let $E$ be a Hausdorff space. If $E$-valued c$\grave{\mbox{a}}$dl$\grave{\mbox{a}}$g
processes $\{X^{i}\}_{i\in\mathbf{I}}$ satisfy MCC, then $\{X^{i}\}_{i\in\mathbf{I}}$
satisfies WMCC.
\end{fact}

\begin{cor}
\label{cor:MCC_5}Let $E$ be a metrizable space and $\{X^{i}\}_{i\in\mathbf{I}}$
be $E$-valued c$\grave{\mbox{a}}$dl$\grave{\mbox{a}}$g processes.
Then:

\renewcommand{\labelenumi}{(\alph{enumi})}
\begin{enumerate}
\item If $(E,\mathfrak{r})$ is a metrization of $E$ and $\{X^{i}\}_{i\in\mathbf{I}}$
satisfies $\mathfrak{r}$-MCC, then $\{X^{i}\}_{i\in\mathbf{I}}$
satisfies MCC.
\item If $E$ is separable and $\{X^{i}\}_{i\in\mathbf{I}}$ satisfies MCC,
then $\{X^{i}\}_{i\in\mathbf{I}}$ satisfies $\mathcal{D}$-FMCC for
some countable $\mathcal{D}=\mathfrak{ac}(\mathcal{D})\subset C_{b}(E;\mathbf{R})$
and $\mathcal{D}$ strongly separates points on $E$. Moreover, $\{X^{i}\}_{i\in\mathbf{I}}$
satisfies $\mathfrak{r}$-MCC for some metrization $(E,\mathfrak{r})$
of $E$.
\end{enumerate}
\end{cor}

The classical results below are quoted from \cite{EK86} for ease
of reference.
\begin{thm}
[\textrm{\cite[\S 3.7, Theorem 7.8]{EK86}}]\label{thm:Sko_RV_WC_FC_Metrizable_Separable}Let
$E$ be a metrizable and separable space and $\{X^{n}\}_{n\in\mathbf{N}}\cup\{X\}$
be $D(\mathbf{R}^{+};E)$-valued random variables. Then:

\renewcommand{\labelenumi}{(\alph{enumi})}
\begin{enumerate}
\item (\ref{eq:WC_on_Path_Space})%
\footnote{Weak convergence and relative compactness of random variables were
specified in \S \ref{sec:RV}. %
} implies%
\footnote{$J(X)$, the set of fixed left-jump times of $X$ was defined in (\ref{eq:J(X)}).%
}
\begin{equation}
X^{n}\xrightarrow{\quad\mathrm{D}(\mathbf{R}^{+}\backslash J(X))\quad}X\mbox{ as }n\uparrow\infty.\label{eq:FC_along_R-J(X)}
\end{equation}

\item If $\{X^{n}\}_{n\in\mathbf{N}}$ is relatively compact in $D(\mathbf{R}^{+};E)$
and (\ref{eq:FC_along_T}) holds for some dense $\mathbf{T}\subset\mathbf{R}^{+}$,
then (\ref{eq:WC_on_Path_Space}) holds.
\end{enumerate}
\end{thm}

\begin{thm}
\label{thm:Sko_RV_Tight_Polish}Let $(E,\mathfrak{r})$ be a complete
separable metric space and $\{X^{i}\}_{i\in\mathbf{I}}$ be $D(\mathbf{R}^{+};E)$-valued
random variables. Then, $\{X^{i}\}_{i\in\mathbf{I}}$ is tight in
$D(\mathbf{R}^{+};E)$%
\footnote{Tightness of random variables was specified in \S \ref{sec:RV}. %
} if and only if $\{X^{i}\}_{i\in\mathbf{I}}$ satisfies MCCC and $\mathfrak{r}$-MCC.
\end{thm}

\chapter{\label{chap:App2}Miscellaneous}

This chapter consists of auxiliary results related to this work. \S
\ref{sec:Gen_Tech} contains a set of basic and general technicalities.
\S \ref{sec:Comp_A1} supplements the topics of Appendix \ref{chap:App1}.
\S \ref{sec:Aux_Rep} houses several auxiliary lemmas about replication
which are used in Chapter \ref{chap:Space_Change} - Chapter \ref{chap:Cadlag}.
All notations, terminologies and conventions introduced before apply
to this appendix. Proofs of many relatively simple results are omitted.

\section{\label{sec:Gen_Tech}General technicalities}
\begin{fact}
\label{fact:Union_Borel}Let $E$ be the union of non-empty sets $\{A_{n}\}_{n\in\mathbf{N}}$.
If $\sigma$-algebras $\mathscr{U}_{1}$ and $\mathscr{U}_{2}$ on
$E$ satisfy $A_{n}\in\mathscr{U}_{2}$ and $\mathscr{U}_{1}|_{A_{n}}=\mathscr{U}_{2}|_{A_{n}}$
for all $n\in\mathbf{N}$, then $\mathscr{U}_{1}\subset\mathscr{U}_{2}$.
\end{fact}
\begin{proof}
We observe that
\begin{equation}
\begin{aligned}\mathscr{U}_{1} & =\left\{ \bigcup_{n\in\mathbf{N}}B\cap A_{n}:B\in\mathscr{U}_{1}\right\} \subset\left\{ \bigcup_{n\in\mathbf{N}}B_{n}:B_{n}\in\mathscr{U}_{1}|_{A_{n}},\forall n\in\mathbf{N}\right\} \\
 & =\left\{ \bigcup_{n\in\mathbf{N}}B_{n}:B_{n}\in\mathscr{U}_{2}|_{A_{n}},\forall n\in\mathbf{N}\right\} \subset\sigma\left(\bigcup_{n\in\mathbf{N}}\mathscr{U}_{2}|_{A_{n}}\right)\subset\mathscr{U}_{2}.
\end{aligned}
\label{eq:Check_Union_Borel}
\end{equation}
\end{proof}

\begin{fact}
\label{fact:Indicator_Modify}Let $(E,\mathscr{U})$ be a measurable
space, $A\in\mathscr{U}$ and $k\in\mathbf{N}$. If $f\in(\mathbf{R}^{k})^{E}$
satisfies $f|_{A}\in M(A,\mathscr{U}|_{A};\mathbf{R}^{k})$, then
$f\mathbf{1}_{A}\in M(E,\mathscr{U};\mathbf{R}^{k})$%
\footnote{$\mathbf{1}_{A}$ denotes the indicator function of $A$.%
}.
\end{fact}

\begin{fact}
\label{fact:var(f)}Let $E$ and $S$ be non-empty sets, $y_{0}\in A\subset E$,
$f\in E^{S}$ and $g\circeq\mathfrak{var}(f;S,f^{-1}(A),y_{0})$.
Then:

\renewcommand{\labelenumi}{(\alph{enumi})}
\begin{enumerate}
\item $\{x\in S:f(x)=g(x)\}\supset f^{-1}(A)$.
\item If $(E,\mathscr{U})$ and $(S,\mathscr{A})$ are measurable spaces,
$f\in M(S,\mathscr{A};E,\mathscr{U})$, $A\in\mathscr{U}$ and $\{y_{0}\}\in\mathscr{U}$,
then $g\in M(S,\mathscr{A};A,\mathscr{U}|_{A})\subset M(S,\mathscr{A};E,\mathscr{U})$.
\end{enumerate}
\end{fact}

\begin{fact}
\label{fact:AS_Cont}Let $E$ and $S$ be topological space, $\mu\in\mathcal{M}^{+}(E)$
and $A$ denote the set of discontinuity points of $f\in S^{E}$.
If $(E,\mathscr{U},\nu)$ is the completion%
\footnote{Completion of measure space and the notation ``$\mathscr{N}(\mu)$''
were specified in \S \ref{sub:Meas}.%
} of $(E,\mathscr{B}(E),\mu)$ and $A\in\mathscr{N}(\mu)$, then $f\in M(E,\mathscr{U};S)$.
\end{fact}
\begin{proof}
Fixing $O\in\mathscr{O}(S)$, we have that
\begin{equation}
(f^{-1}(O)\backslash A)=(f|_{E\backslash A})^{-1}(O)\in\mathscr{O}_{E}(E\backslash A)\subset\mathscr{B}_{E}(E\backslash A)\subset\mathscr{U}.\label{eq:Check_AS_Cont_Measurable_1}
\end{equation}
Next, we recall the fact $A\in\mathscr{N}(\mu)\subset\mathscr{U}$,
so $f^{-1}(O)\in\mathscr{U}$ and get by the continuity of $f|_{E\backslash A}$
that
\begin{equation}
f^{-1}(O)\cap A\in\mathscr{N}(\mu)\subset\mathscr{U}.\label{eq:Check_AS_Cont_Measurable_2}
\end{equation}
\end{proof}

$\sigma(C(E;\mathbf{R}))$, the Baire $\sigma$-algebra on $E$ is
generally smaller than $\mathscr{B}(E)$.
\begin{fact}
\label{fact:O_D_(A)_B_D_(A)}Let $E$ be a topological space, $S$
be a non-empty set and $A\subset S$. Then, $\sigma(\mathcal{D})|_{A}\subset\mathscr{B}_{\mathcal{D}}(A)$
for any $\mathcal{D}\subset E^{S}$ and the equality holds if $E$
is a second-countable space and $\mathcal{D}$ is countable.
\end{fact}
\begin{proof}
We have that
\begin{equation}
\begin{aligned}\mathscr{B}_{\mathcal{D}}(A) & \supset\sigma\left(\bigcup_{f\in\mathcal{D}}\left\{ f^{-1}(O)\cap A:O\in\mathscr{O}(E)\right\} \right)\\
 & =\sigma\left[\bigcup_{f\in\mathcal{D}}\sigma\left(\left\{ f^{-1}(O)\cap A:O\in\mathscr{O}(E)\right\} \right)\right]\\
 & =\sigma\left(\left\{ f^{-1}(B)\cap A:B\in\sigma\left(\mathscr{O}(E)\right)=\mathscr{B}(E),f\in\mathcal{D}\right\} \right)=\left.\sigma(\mathcal{D})\right|_{A}.
\end{aligned}
\label{eq:Check_Borel_>_Baire}
\end{equation}

If $\{O_{n}\}_{n\in\mathbf{N}}$ is a countable topological basis
of $E$ and $\mathcal{D}$ is countable, then
\begin{equation}
\left\{ \bigcap_{f\in\mathcal{D}_{0}}f^{-1}(O_{n})\cap A:n\in\mathbf{N},\mathcal{D}_{0}\in\mathscr{P}_{0}(\mathcal{D})\right\} \label{eq:Check_B_D_(A)=00003DSigma_D_1}
\end{equation}
is a countable basis for $\mathscr{O}_{\mathcal{D}}(A)$ by \cite[Lemma 13.1]{M00}.
Consequently,
\begin{equation}
\begin{aligned}\mathscr{O}_{\mathcal{D}}(A) & \subset\sigma\left(\left\{ \bigcap_{f\in\mathcal{D}_{0}}f^{-1}(O_{n})\cap A:n\in\mathbf{N},\mathcal{D}_{0}\in\mathscr{P}_{0}(\mathcal{D})\right\} \right)\\
 & \subset\sigma\left(\left\{ f^{-1}(B)\cap A:B\in\mathscr{B}(E),f\in\mathcal{D}\right\} \right)=\left.\sigma(\mathcal{D})\right|_{A}.
\end{aligned}
\label{eq:Check_B_D_(A)=00003DSigma_D_2}
\end{equation}
\end{proof}

\begin{fact}
\label{fact:Uni_Seq_Lim_Conv}Let $E$ be a topological space and
$\{x_{n}\}_{n\in\mathbf{N}}\subset E$. If every convergent subsequence
of $\{x_{n}\}_{n\in\mathbf{N}}$ must converge to $x$ as $n\uparrow\infty$,
and if any infinite subset of $\{x_{n}\}_{n\in\mathbf{N}}$ has a
convergent subsequence, then $x_{n}\rightarrow x$ as $n\uparrow\infty$
in $E$.
\end{fact}

\begin{lem}
\label{lem:Imb_SSP}Let $E$ and $S$ be topological spaces, $\mathcal{D}\subset S^{E}$
and equip $V\circeq\bigotimes\mathcal{D}(E)$ with the subspace topology
$\mathscr{O}_{S^{\mathcal{D}}}(V)$. Then:

\renewcommand{\labelenumi}{(\alph{enumi})}
\begin{enumerate}
\item $(\bigotimes\mathcal{D})^{-1}\in C(V;E)$ if and only if $\mathscr{O}(E)\subset\mathscr{O}_{\mathcal{D}}(E)$
and $\bigotimes\mathcal{D}$ is injective.
\item $\bigotimes\mathcal{D}\in\mathbf{hom}(E;V)$ if and only if $\mathscr{O}(E)=\mathscr{O}_{\mathcal{D}}(E)$
and $\bigotimes\mathcal{D}$ is injective.
\end{enumerate}
\end{lem}
\begin{proof}
(a) We find by Fact \ref{fact:Prod_Map_2} (b) that $\bigotimes\mathcal{D}\in\mathbf{imb}(E,\mathscr{O}_{\mathcal{D}}(E);V)$
if and only if $\bigotimes\mathcal{D}$ is injective. Given the injectiveness
of $\bigotimes\mathcal{D}$, $(\bigotimes\mathcal{D})^{-1}\in C(V;E)$
precisely when $\mathscr{O}(E)$ is coarser than $\mathscr{O}_{\mathcal{D}}(E)$.

(b) is immediate by (a).\end{proof}

\begin{lem}
\label{lem:Biso_Hom}Let $E$ and $S$ be topological spaces, $f\in E^{S}$
and $\mathscr{U}_{f}\circeq\{O\subset E:f^{-1}(O)\in\mathscr{O}(S)\}$.
Then:

\renewcommand{\labelenumi}{(\alph{enumi})}
\begin{enumerate}
\item If $f$ is bijective, then $f\in\mathbf{hom}(S;(E,\mathscr{U}_{f}))$.
\item If $f\in\mathbf{biso}(S;E)$, then $\mathscr{B}(E)=\sigma(\mathscr{U}_{f})$.
\end{enumerate}
\end{lem}
\begin{proof}
(a) $\mathscr{U}_{f}$ is a topology, $f\in C(S;E,\mathscr{U}_{f})$
and
\begin{equation}
\mathscr{U}_{f}=\left\{ f(B):B\in\mathscr{O}(S)\right\} \label{eq:U_O(S)}
\end{equation}
by the bijectiveness of $f$, thus proving $f^{-1}\in C(E,\mathscr{U}_{f};S)$.

(b) $f\in\mathbf{biso}(S;E)$ satisfies (\ref{eq:U_O(S)}) and further
satisfies
\begin{equation}
\mathscr{B}(E)=\left\{ f(B):B\in\mathscr{B}(S)\right\} =\sigma\left(\left\{ f(B):B\in\mathscr{O}(S)\right\} \right)=\sigma(\mathscr{U}_{f}).\label{eq:Check_B(S)_U}
\end{equation}
\end{proof}

\begin{rem}
\label{rem:Biso_Hom}The lemma above shows a Borel isomorphism can
be turned into a homeomorphism by changing the topology generating
the Borel $\sigma$-algebra.\end{rem}
\begin{fact}
\label{fact:Path_Mapping}Let $\mathbf{I}$, $E$ and $S$ be non-empty
sets and $f\in S^{E}$. Then:

\renewcommand{\labelenumi}{(\alph{enumi})}
\begin{enumerate}
\item If $f$ is injective or surjective, then $\varpi_{\mathbf{I}}(f)$
is also.
\item If $(E,\mathscr{U})$ and $(S,\mathscr{A})$ are measurable spaces
and $f\in M(E,\mathscr{U};S,\mathscr{A})$, then $\varpi_{\mathbf{I}}(f)\in M(E^{\mathbf{I}},\mathscr{U}^{\otimes\mathbf{I}};S^{\mathbf{I}},\mathscr{A}^{\otimes\mathbf{I}})$.
\item If $E$ and $S$ are topological spaces and $f\in C(E;S)$, then $\varpi_{\mathbf{I}}(f)\in C(E^{\mathbf{I}};S^{\mathbf{I}})$.
\end{enumerate}
\end{fact}
\begin{proof}
(a), (b) and (c) are immediate by definition, Fact \ref{fact:Prod_Map_1}
(b) and Fact \ref{fact:Prod_Map_2} (b) respectively.\end{proof}

\begin{fact}
\label{fact:Seq_Prod_Conv}Let $\{S_{i}\}_{i\in\mathbf{I}}$ be topological
spaces. Then, $x_{k}\rightarrow x$ as $k\uparrow\infty$ in $\prod_{i\in\mathbf{I}}S_{i}$
if and only if $\mathfrak{p}_{i}(x_{k})\rightarrow\mathfrak{p}_{i}(x)$
as $k\uparrow\infty$ in $S_{i}$ for all $i\in\mathbf{I}$.
\end{fact}
\begin{proof}
This fact was justified in \cite[\S 19, Exercise 6]{M00}.\end{proof}

\begin{fact}
\label{fact:Interval_Union}Let $\mathbf{I}$ be an arbitrary index
set and $\{a_{i},b_{i}\}\subset\mathbf{R}$ satisfy $a_{i}<b_{i}$
for all $i\in\mathbf{I}$. Then, $\bigcup_{i\in\mathbf{I}}[a_{i},b_{i})\in\mathscr{B}(\mathbf{R})$.
\end{fact}
\begin{proof}
For each $\{i_{1},i_{2}\}\subset\mathbf{I}$, we define $i_{1}\sim i_{2}$
if there exist some $\mathbf{I}_{0}\subset\mathbf{I}$ such that $\bigcup_{i\in\mathbf{I}_{0}\cup\{i_{1},i_{2}\}}[a_{i},b_{i})$
is an interval. It is not difficult to see ``$\sim$'' defines an
equivalence relation on $\mathbf{I}$. Let $\{\mathbf{I}_{j}:j\in\mathbf{J}\}$
be the ``$\sim$'' equivalence classes of the members of $\mathbf{I}$
and $A_{j}\circeq\bigcup_{i\in\mathbf{I}_{j}}[a_{i},b_{i})$ for each
$j\in\mathbf{J}$. These $\{A_{j}\}_{j\in\mathbf{J}}$ are pairwisely
disjoint intervals by the definition of ``$\sim$'', so $\mathbf{J}$
is countable by \cite[\S 30, Exercise 13]{M00}. Hence, $\bigcup_{i\in\mathbf{I}}[a_{i},b_{i})=\bigcup_{j\in\mathbf{J}}A_{j}\in\mathscr{B}(\mathbf{R})$.\end{proof}

\begin{fact}
\label{fact:Cadlag_Path}Let $E$, $S$ and $\{S_{i}\}_{i\in\mathbf{I}}$
be topological spaces and $f\in C(E;S)$. Then:

\renewcommand{\labelenumi}{(\alph{enumi})}
\begin{enumerate}
\item If $x\in E^{\mathbf{R}^{+}}$ is right-continuous, then $x\in M(\mathbf{R}^{+};E)$.
\item If $x\in E^{\mathbf{R}^{+}}$ is c$\grave{\mbox{a}}$dl$\grave{\mbox{a}}$g
and $f\in C(E;S)$, then $\varpi(f)(x)\in S^{\mathbf{R}^{+}}$ is
also c$\grave{\mbox{a}}$dl$\grave{\mbox{a}}$g.
\item $\bigotimes_{i\in\mathbf{I}}f_{i}:E\rightarrow S^{i}$ is c$\grave{\mbox{a}}$dl$\grave{\mbox{a}}$g
if and only if $f_{i}:E\rightarrow S_{i}$ is c$\grave{\mbox{a}}$dl$\grave{\mbox{a}}$g
for all $i\in\mathbf{I}$.
\end{enumerate}
\end{fact}
\begin{proof}
(a) Note that
\begin{equation}
x_{n}\circeq\sum_{i=1}^{n2^{n}}x\left(\frac{i}{2^{n}}\right)\mathbf{1}_{\left[\frac{i-1}{2^{n}},\frac{i}{2^{n}}\right)}+x(n)\mathbf{1}_{[n,\infty)}\in M(\mathbf{R}^{+};E),\;\forall n\in\mathbf{N}\label{eq:D(E)_Approx_Seq}
\end{equation}
and $x_{n}\rightarrow x$ as $n\uparrow\infty$ in $E^{\mathbf{R}^{+}}$.
Then, $x\in M(\mathbf{R}^{+};E)$ as pointwise convergence preserves
measurability.

(b) and (c) are immediate by the definitions of $\varpi(f)$, $\bigotimes_{i\in\mathbf{I}}\mathscr{O}(S_{i})$
and product topology.\end{proof}

\begin{fact}
\label{fact:ac_mc_Countable}Let $E$ be a non-empty set, $d,k\in\mathbf{N}$
and $\mathcal{D}\subset(\mathbf{R}^{k})^{E}$ be a countable collection.
Then, $\mathfrak{ae}(\mathcal{D})$ and $\mathfrak{ac}(\mathcal{D})$
are countable collections. When $k=1$, $\mathfrak{mc}(\mathcal{D})$
and $\mathfrak{ag}_{\mathbf{Q}}(\mathcal{D})$ are also countable
collections.
\end{fact}
\begin{proof}
Let $\mathcal{D}=\{f_{n}\}_{n\in\mathbf{N}}$. For each $m\in\mathbf{N}$,
we observe that $\sum_{i=1}^{m}f_{n_{i}}\mapsto(n_{1},....,n_{m})$
defines an injective mapping from $\mathcal{D}_{m}\circeq\{\sum_{i=1}^{m}f_{n_{i}}:n_{1},...,n_{m}\in\mathbf{N}\}$
to the countable set $\mathbf{N}^{m}$, so $\mathcal{D}_{m}$ is countable.
As a result, $\mathfrak{ae}(\mathcal{D})=\mathcal{D}_{1}\cup\mathcal{D}_{2}$
and $\mathfrak{ac}(\mathcal{D})=\bigcup_{m\in\mathbf{N}}\mathcal{D}_{m}$
are both countable.

Next, we let $k=1$ and observe that $\prod_{i=1}^{m}f_{n_{i}}\mapsto(n_{1},....,n_{m})$
defines an injective mapping from $\mathcal{D}_{m}^{\prime}\circeq\{\prod_{i=1}^{m}f_{n_{i}}:n_{1},...,n_{m}\in\mathbf{N}\}$
to the countable set $\mathbf{N}^{m}$, so $\mathcal{D}_{m}^{\prime}$
is countable. As a result, $\mathfrak{mc}(\mathcal{D})=\bigcup_{m\in\mathbf{N}}\mathcal{D}_{m}^{\prime}$
is also countable.

Furthermore, we index $\mathfrak{mc}(\mathcal{D})$ by $\mathbf{N}$
as $\{g_{j}\}_{j\in\mathbf{N}}$ and observe that $ag_{j}\longmapsto(j,a)$
defines an injective mapping from $\mathcal{D}_{\mathbf{Q}}\circeq\{ag_{j}:j\in\mathbf{N},a\in\mathbf{Q}\}$
to the countable set $\mathbf{N}\times\mathbf{Q}$, so $\mathcal{D}_{\mathbf{Q}}$
is countable. As a result, $\mathfrak{ag}_{\mathbf{Q}}(\mathcal{D})=\mathfrak{ac}(\mathcal{D}_{\mathbf{Q}})$
is also countable.\end{proof}

\begin{fact}
\label{fact:Pi^d}Let $E$ be a non-empty set, $d\in\mathbf{N}$ and
$\mathcal{D}\subset\mathbf{R}^{E}$. Then:

\renewcommand{\labelenumi}{(\alph{enumi})}
\begin{enumerate}
\item $\Pi^{d}(\mathcal{D})$ is a countable collection whenever $\mathcal{D}$
is. Moreover,
\begin{equation}
\begin{aligned}\Pi^{d}\left(\mathfrak{ac}(\mathcal{D})\right) & \subset\mathfrak{ac}\left(\Pi^{d}(\mathcal{D})\right),\\
\Pi^{d}\left(\mathfrak{mc}(\mathcal{D})\right) & =\mathfrak{mc}\left(\Pi^{d}(\mathcal{D})\right),\\
\Pi^{d}\left(\mathfrak{ag}_{\mathbf{Q}}(\mathcal{D})\right) & \subset\mathfrak{ag}_{\mathbf{Q}}\left(\Pi^{d}(\mathcal{D})\right),\\
\Pi^{d}\left(\mathfrak{ag}(\mathcal{D})\right) & \subset\mathfrak{ag}\left(\Pi^{d}(\mathcal{D})\right).
\end{aligned}
\label{eq:Pi^d_and_Operations}
\end{equation}

\item If the members of $\mathcal{D}$ are bounded, then those of $\Pi^{d}(\mathcal{D})$
are also. Moreover,
\begin{equation}
\begin{aligned}\Pi^{d}\left(\mathfrak{cl}(\mathcal{D})\right) & \subset\mathfrak{cl}\left(\Pi^{d}(\mathcal{D})\right),\\
\Pi^{d}\left(\mathfrak{ca}(\mathcal{D})\right) & \subset\mathfrak{ca}\left(\Pi^{d}(\mathcal{D})\right).
\end{aligned}
\label{eq:Pi^d_Closure}
\end{equation}

\end{enumerate}
\end{fact}
\begin{proof}
(a) If $\mathcal{D}=\{f_{n}\}_{n\in\mathbf{N}}$, then $\prod_{i=1}^{k}f_{n_{i}}\circ\mathfrak{p}_{i}\mapsto(n_{1},...,n_{k})$
defines an injective mapping from $\mathcal{D}_{k}\circeq\{\prod_{i=1}^{k}f_{n_{i}}\circ\mathfrak{p}_{i}:n_{i}\in\mathbf{N}\}\subset\mathbf{R}^{E^{d}}$
to $\mathbf{N}^{k}$ for each $k\in\{1,...,d\}$, so $\mathcal{D}_{k}$
is countable. Hence, $\Pi^{d}(\mathcal{D})=\bigcup_{k=1}^{d}\mathcal{D}_{k}$
is countable.

Letting $k\in\{1,...,d\}$, $n_{1},....,n_{k}\in\mathbf{N}$ and $\{f_{i,j}\}_{1\leq j\leq n_{i},1\leq i\leq k}\subset\mathcal{D}$,
one notes
\begin{equation}
\prod_{i=1}^{k}\left(\sum_{j=1}^{n_{i}}f_{i,j}\right)\circ\mathfrak{p}_{i}=\sum_{j_{1}=1}^{n_{1}}...\sum_{j_{k}=1}^{n_{k}}\left(\prod_{i=1}^{k}f_{i,j_{i}}\circ\mathfrak{p}_{i}\right)\in\mathfrak{ac}\left(\Pi^{d}(\mathcal{D})\right),\label{eq:Check_Pi^d_ac}
\end{equation}
and
\begin{equation}
\prod_{i=1}^{k}\left(\prod_{j=1}^{n_{i}}f_{i,j}\right)\circ\mathfrak{p}_{i}=\prod_{1\leq j\leq n_{i},1\leq i\leq k}f_{i,j}\circ\mathfrak{p}_{i}\in\mathfrak{mc}\left(\Pi^{d}(\mathcal{D})\right).\label{eq:Check_Pi^d_mc_1}
\end{equation}
Letting $N\in\mathbf{N}$, $k_{1},...,k_{N}\in\{1,...,d\}$ and $\{f_{i,j}\}_{1\leq i\leq k_{i},1\leq j\leq N}\subset\mathcal{D}$,
we observe that
\begin{equation}
\prod_{j=1}^{N}\left(\prod_{i=1}^{k_{j}}f_{i,j}\circ\mathfrak{p}_{i}\right)=\prod_{i=1}^{k^{*}}\left(\prod_{j\in J_{i}}f_{i,j}\right)\circ\mathfrak{p}_{i}\in\Pi^{d}\left(\mathfrak{mc}(\mathcal{D})\right),\label{eq:Check_Pi^d_mc_2}
\end{equation}
where $k^{*}\circeq\max\{k_{1},....,k_{N}\}$ and $J_{i}\circeq\{j\in\{1,...,n\}:k_{j}\geq i\}$
for each $1\leq i\leq k^{*}$. Then, the first two lines of (\ref{eq:Pi^d_and_Operations})
follow by (\ref{eq:Check_Pi^d_ac}), (\ref{eq:Check_Pi^d_mc_1}) and
(\ref{eq:Check_Pi^d_mc_2}).

Using the second line of (\ref{eq:Pi^d_and_Operations}), we have
that
\begin{equation}
\begin{aligned}\Pi^{d}\left(\left\{ af:f\in\mathfrak{mc}(\mathcal{D}),a\in\mathbf{Q}\right\} \right) & =\left\{ af:f\in\Pi^{d}\left(\mathfrak{mc}(\mathcal{D})\right),a\in\mathbf{Q}\right\} \\
 & =\left\{ af:f\in\mathfrak{mc}\left(\Pi^{d}(\mathcal{D})\right),a\in\mathbf{Q}\right\} .
\end{aligned}
\label{eq:Check_Pi^d_ag_Q_1}
\end{equation}
Using (\ref{eq:Check_Pi^d_ag_Q_1}) and the first two lines of (\ref{eq:Pi^d_and_Operations}),
we have that
\begin{equation}
\begin{aligned}\Pi^{d}\left(\mathfrak{ag}_{\mathbf{Q}}(\mathcal{D})\right) & =\Pi^{d}\left[\mathfrak{ac}\left(\left\{ af:f\in\mathfrak{mc}(\mathcal{D}),a\in\mathbf{Q}\right\} \right)\right]\\
 & \subset\mathfrak{ac}\left[\Pi^{d}\left(\left\{ af:f\in\mathfrak{mc}(\mathcal{D}),a\in\mathbf{Q}\right\} \right)\right]\\
 & =\mathfrak{ac}\left(\left\{ af:f\in\mathfrak{mc}\left(\Pi^{d}(\mathcal{D})\right),a\in\mathbf{Q}\right\} \right)=\mathfrak{ag}_{\mathbf{Q}}\left(\Pi^{d}(\mathcal{D})\right),
\end{aligned}
\label{eq:Check_Pi^d_ag_Q_2}
\end{equation}
which proves the third line of (\ref{eq:Pi^d_and_Operations}). The
fourth line of (\ref{eq:Pi^d_and_Operations}) follows by a similar
argument with $\mathbf{Q}$ replaced by $\mathbf{R}$.

(b) We fix $1\leq k\leq d$. Boundedness is immediate. Suppose $\{f_{1},...,f_{k}\}\subset\mathfrak{cl}(\mathcal{D})$.
By Fact \ref{fact:First_Countable} (with $E=(\mathfrak{cl}(\mathcal{D}),\Vert\cdot\Vert_{\infty})$
and $A=\mathcal{D}$), there exist $\{f_{i,n}\}_{1\leq i\leq k,n\in\mathbf{N}}\subset\mathcal{D}$
such that $f_{i,n}\overset{u}{\rightarrow}f$ as $n\uparrow\infty$
for all $1\leq i\leq k$. We let $c\circeq(\sup_{1\leq i\leq d}\Vert f_{i}\Vert_{\infty})^{k-1}$
and find that
\begin{equation}
\lim_{n\rightarrow\infty}\left\Vert \prod_{i=1}^{k}f_{i}\circ\mathfrak{p}_{i}-\prod_{i=1}^{k}f_{i,n}\circ\mathfrak{p}_{i}\right\Vert \leq c\lim_{n\rightarrow\infty}\sum_{i=1}^{k}\Vert f_{i}-f_{i,n}\Vert_{\infty}=0,\label{eq:Check_Pi^d_cl}
\end{equation}
thus proving the first line of (\ref{eq:Pi^d_Closure}). The second
line of (\ref{eq:Pi^d_Closure}) is immediate by the first line (with
$\mathcal{D}=\mathfrak{ag}(\mathcal{D})$).\end{proof}

\begin{lem}
\label{lem:C0_Ext}Let $E$ be an open subspace of $S$ and $f\in C(E;\mathbf{R}^{k})$.
If for any $\epsilon\in(0,\infty)$, there exists an $A_{\epsilon}\subset E$
such that $A_{\epsilon}\in\mathscr{C}(S)$ and $\Vert f|_{E\backslash A_{\epsilon}}\Vert_{\infty}<\epsilon$,
then $g\circeq\mathfrak{var}(f;S,E,0)$%
\footnote{``$\mathfrak{var}(\cdot)$'' was defined in Notation \ref{notation:Var}.%
} is a continuous extension of $f$ on $S$.
\end{lem}
\begin{proof}
We need only prove the case of $k=1$ and the general result follows
by Fact \ref{fact:Prod_Map_2} (b). Let $\epsilon\in\mathbf{R}\backslash\{0\}$.
From the facts
\begin{equation}
g^{-1}\left[(-\infty,\epsilon)\right]\backslash A_{\left|\epsilon\right|}=\begin{cases}
S\backslash A_{\left|\epsilon\right|}\in\mathscr{O}(S), & \mbox{if }\epsilon>0,\\
\varnothing, & \mbox{if }\epsilon<0
\end{cases}\label{eq:Check_C0_Extend_1}
\end{equation}
and $E\in\mathscr{O}(S)$ it follows that
\begin{equation}
g^{-1}\left[(-\infty,\epsilon)\right]=\begin{cases}
f^{-1}\left[(-\infty,\epsilon)\right]\cup(S\backslash A_{\left|\epsilon\right|})\in\mathscr{O}(S), & \mbox{if }\epsilon>0,\\
f^{-1}\left[(-\infty,\epsilon)\right]\in\mathscr{O}(E)\subset\mathscr{O}(S), & \mbox{if }\epsilon<0,
\end{cases}\label{eq:Check_C0_Extend_2}
\end{equation}
thus proving the continuity of $g$.\end{proof}

\begin{fact}
\label{fact:SP_Nested_Union}Let $\{A_{n}\}_{n\in\mathbf{N}}$ be
nested%
\footnote{We explained the meaning of ``nested'' in Fact \ref{fact:D-Baseable_Union}.%
} non-empty subsets of $E$ and $\mathcal{D}_{n}\subset\mathbf{R}^{E}$
separate points on $A_{n}$ for each $n\in\mathbf{N}$. Then, $\bigcup_{n\in\mathbf{N}}\mathcal{D}_{n}$
separates points on $\bigcup_{n\in\mathbf{N}}A_{n}$.
\end{fact}

\begin{fact}
\label{fact:Path_Mapping_Injective}Let $A$ be a non-empty subset
of $E$ and $\mathcal{D}\subset\mathbf{R}^{E}$ separate points on
$A\subset E$. Then, $\varpi(\bigotimes\mathcal{D})$ and $\varpi(\mathcal{D})$
are both injective restricted to $A^{\mathbf{R}^{+}}$.
\end{fact}
\begin{proof}
The injectiveness of $\varpi(\bigotimes\mathcal{D})$ on $A^{\mathbf{R}^{+}}$
is immediate by Fact \ref{fact:Path_Mapping} (a) (with $E=A$, $\mathbf{I}=\mathbf{R}^{+}$,
$S=\mathbf{R}^{\mathcal{D}}$ and $f=\bigotimes\mathcal{D}|_{A}$).
Furthermore, we note that $\varpi(\mathcal{D})(x)=\varpi(\mathcal{D})(y)$
in $(A^{\mathbf{R}^{+}})^{\mathcal{D}}$ implies $\bigotimes\mathcal{D}[x(t)]=\bigotimes\mathcal{D}[y(t)]$
for all $t\in\mathbf{R}^{+}$. This indicates $x(t)=y(t)$ for all
$t\in\mathbf{R}^{+}$, i.e. $x=y$.\end{proof}

\begin{fact}
\label{fact:WC_M(E)_P(E)}Let $E$ be a topological space. Then:

\renewcommand{\labelenumi}{(\alph{enumi})}
\begin{enumerate}
\item $\mu_{1}=\mu_{2}$ in $\mathcal{M}^{+}(E)$ if and only if $\mu_{1}/\mu_{1}(E)=\mu_{2}/\mu_{2}(E)$%
\footnote{We mentioned in \S \ref{sec:Convention} that any measure in this
work have positive total mass.%
} in $\mathcal{P}(E)$ and $\mu_{1}(E)=\mu_{2}(E)$.
\item (\ref{eq:Mu_n_WC_Mu_M(E)}) holds if and only if $\lim_{n\rightarrow\infty}\mu_{n}(E)=\mu(E)$
and
\begin{equation}
\frac{\mu_{n}}{\mu_{n}(E)}\Longrightarrow\frac{\mu}{\mu(E)}\mbox{ as }n\uparrow\infty\mbox{ in }\mathcal{P}(E).\label{eq:Mu_n_WC_Mu_Rescaled_P(E)}
\end{equation}

\end{enumerate}
\end{fact}

\begin{fact}
\label{fact:Sep_CD}Let $E$ be a topological space and $1\in\mathcal{D}\subset M_{b}(E;\mathbf{R})$.
Then:

\renewcommand{\labelenumi}{(\alph{enumi})}
\begin{enumerate}
\item $\mathcal{D}$ is separating on $E$%
\footnote{The terminologies ``separating'' and ``convergence determining''
were introduced in \S \ref{sec:Borel_Measure}.%
} if and only if $\mathcal{D}^{*}$ separates points on $\mathcal{P}(E)$.
\item $\mathcal{D}$ is convergence determining on $E$ if and only if $\mathcal{D}^{*}$
determines point convergence on $\mathcal{P}(E)$.
\end{enumerate}
\end{fact}
\begin{proof}
This result is immediate by Fact \ref{fact:WC_M(E)_P(E)}.\end{proof}

\begin{fact}
\label{fact:P(E)_Closed_M(E)}Let $E$ be a topological space. Then,
$\mathcal{P}(E)\in\mathscr{C}[\mathcal{M}^{+}(E)]$.
\end{fact}
\begin{proof}
Let $\mu$ be a limit point%
\footnote{$\mathcal{M}^{+}(E)$ as aforementioned is not necessarily first-countable.
So, $\mu$ being a limit point of $\mathcal{P}(E)$ does not necessarily
imply a subsequence of $\mathcal{P}(E)$ converging weakly to $\mu$.%
} of $\mathcal{P}(E)$ in $\mathcal{M}^{+}(E)$. Then, there exist
$\{\mu_{p}\}_{p\in\mathbf{N}}\subset\mathcal{P}(E)$ such that $\lim_{p\rightarrow\infty}\left|\mu(E)-1\right|=\lim_{p\rightarrow\infty}\left|\mu(E)-\mu_{p}(E)\right|=0$.\end{proof}

\begin{fact}
\label{fact:Dirac_WC}If $x_{n}\rightarrow x$ as $n\uparrow\infty$
in topological space $E$, then $\delta_{x_{n}}\Rightarrow\delta_{x}$
as $n\uparrow\infty$ in $\mathcal{P}(E)$.
\end{fact}

The generalized Portmanteau's Theorem helps to establish the Continuous
Mapping Theorem on general topological spaces.
\begin{thm}
[\textbf{Continuous Mapping Theorem}]\label{thm:ContMapTh}Let $E$
and $S$ be topological spaces. Then:

\renewcommand{\labelenumi}{(\alph{enumi})}
\begin{enumerate}
\item If $f\in C(E;S)$, then (\ref{eq:Mu_n_WC_Mu_M(E)}) implies
\begin{equation}
\mu_{n}\circ f^{-1}\Longrightarrow\mu\circ f^{-1}\mbox{ as }n\uparrow\infty\mbox{ in }\mathcal{M}^{+}(S).\label{eq:ContMapTh}
\end{equation}

\item If $E$ is a Tychonoff space and the set of discontinuity points of
$f\in M(E;S)$ belongs to $\mathscr{N}(\mu)$, then (\ref{eq:Mu_n_WC_Mu_M(E)})
implies (\ref{eq:ContMapTh}).
\end{enumerate}
\end{thm}
\begin{proof}
(a) follows by the fact that
\begin{equation}
\begin{aligned}\lim_{n\rightarrow\infty}g^{*}\left(\mu_{n}\circ f^{-1}\right) & =\lim_{n\rightarrow\infty}(g\circ f)^{*}(\mu_{n})\\
 & =(g\circ f)^{*}(\mu)=g^{*}\left(\mu\circ f^{-1}\right),\;\forall g\in C_{b}(S;\mathbf{R}).
\end{aligned}
\label{eq:Check_ContMapTh_1}
\end{equation}

(b) Let $O\in\mathscr{O}(S)$, $A\subset E$ be the set of discontinuity
points of $f$, $U$ be the interior of $f^{-1}(O)$ and $\nu$ be
the completion of $\mu$. $\mu(U)=\mu\circ f^{-1}(O)$ since $(f^{-1}(O)\backslash A)\subset U\subset f^{-1}(O)$
and $A\in\mathscr{N}(\mu)$. It follows by the Tychonoff property
of $E$ and Theorem \ref{thm:Portamenteau} (a, c) that
\begin{equation}
\begin{aligned}\mu\circ f^{-1}(O) & =\mu(U)\leq\liminf_{n\rightarrow\infty}\mu_{n}(U)\leq\liminf_{n\rightarrow\infty}\mu_{n}\circ f^{-1}(O).\end{aligned}
\label{eq:Check_ContMapTh_2}
\end{equation}
Now, (b) follows by (\ref{eq:Check_ContMapTh_2}) and Theorem \ref{thm:Portamenteau}
(a, c).\end{proof}

\begin{fact}
\label{fact:Weak_Topo_Coarsen}Let $E$ be a topological space, $(E,\mathscr{U})$
be a topological coarsening of $E$ and $S\circeq\mathcal{M}^{+}(E,\mathscr{U})$.
Then:

\renewcommand{\labelenumi}{(\alph{enumi})}
\begin{enumerate}
\item $(\mathcal{M}^{+}(E),\mathscr{O}_{S}[\mathcal{M}^{+}(E)])$ and $(\mathcal{P}(E),\mathscr{O}_{S}[\mathcal{P}(E)])$
are topological coarsenings of $\mathcal{M}^{+}(E)$ and $\mathcal{P}(E)$
respectively.
\item If $\mu_{n}\Rightarrow\mu$ as $n\uparrow\infty$ in $\mathcal{M}^{+}(E)$,
then $\mu_{n}\Rightarrow\mu$ as $n\uparrow\infty$ in $S$.
\end{enumerate}
\end{fact}
\begin{proof}
(a) $\mathscr{U}\subset\mathscr{O}(E)$ implies $\mathscr{B}(E,\mathscr{U})\subset\mathscr{B}(E)$,
so every $\mu\in\mathcal{M}^{+}(E)$ is naturally a member of $S$.
$\mathscr{U}\subset\mathscr{O}(E)$ implies $C_{b}(E,\mathscr{U};\mathbf{R})\subset C_{b}(E;\mathbf{R})$.
Then, (a) follows by the fact that
\begin{equation}
\begin{aligned}\mathscr{O}\left[\mathcal{M}^{+}(E)\right] & =\mathscr{O}_{C_{b}(E;\mathbf{R})^{*}}\left[\mathcal{M}^{+}(E)\right]\\
 & \supset\mathscr{O}_{C_{b}\left(E,\mathscr{U};\mathbf{R}\right)^{*}}\left[\mathcal{M}^{+}(E)\right]=\mathscr{O}_{S}\left[\mathcal{M}^{+}(E)\right].
\end{aligned}
\label{eq:Topo_Coarsen_Weak_Topo}
\end{equation}

(b) is immediate by (a).\end{proof}

\begin{fact}
\label{fact:WC_Completion}Let $E$ be a topological space, $\mu_{n}\Rightarrow\mu_{0}$
as $n\uparrow\infty$ in $\mathcal{M}^{+}(E)$ and $(E,\mathscr{U}_{n},\nu_{n})$
be the completion of $(E,\mathscr{B}(E),\mu_{n})$ for each $n\in\mathbf{N}_{0}$.
Then, $\nu_{n}\Rightarrow\nu_{0}$ as $n\uparrow\infty$ in $\mathcal{M}^{+}(E)$%
\footnote{Recall that ``$\nu_{n}\Rightarrow\nu_{0}$ as $n\uparrow\infty$
in $\mathcal{M}^{+}(E)$'' means $\{\nu_{n}\}_{n\in\mathbf{N}}$
converges weakly to $\nu_{0}$ as members of $\mathcal{M}^{+}(E)$,
which is well-defined since $\mathscr{U}_{n}\supset\mathscr{B}(E)$
and $f^{*}(\nu_{n})=f^{*}(\mu_{n})$ for all $n\in\mathbf{N}_{0}$
and $f\in C_{b}(E;\mathbf{R})$.%
}.
\end{fact}

\begin{fact}
\label{fact:Seq_Tight_Support}Let $E$ be a topological space and
$\mathscr{U}$ be a $\sigma$-algebra on $E$. If $\Gamma\subset\mathfrak{M}^{+}(E,\mathscr{U})$
is sequentially tight in $A\subset E$%
\footnote{The notion of sequential tightness was introduced in \S \ref{sec:WLP_Uni}.%
}, then there exists a $\Gamma_{0}\in\mathscr{P}_{0}(\Gamma)$ such
that $A$ is a common support of all members of $\Gamma\backslash\Gamma_{0}$.
\end{fact}
\begin{proof}
Suppose none of $\{\mu_{n}\}_{n\in\mathbf{N}}\subset\Gamma$ is supported
on $A$. The sequential tightness of $\Gamma$ implies a subsequence
$\{\mu_{n_{k}}\}_{k\in\mathbf{N}}$ being tight in $A$. In other
words, $\{\mu_{n_{k}}\}_{k\in\mathbf{N}}$ are all supported on some
$B\in\mathscr{K}_{\sigma}(E)$ with $B\subset A$. Contradiction!\end{proof}

\begin{fact}
\label{fact:Push_Forward_Tight_1}Let $\mathscr{U}$ and $\mathscr{A}$
be $\sigma$-algebras on topological spaces $E$ and $S$ respectively.
If $\Gamma\subset\mathfrak{M}^{+}(E,\mathscr{U})$ is tight in $A\subset E$%
\footnote{The definitions of tightness and $\mathbf{m}$-tightness for possibly
non-Borel measures are in Definition \ref{def:Tight}.%
}, and if $f\in M(E,\mathscr{U};S,\mathscr{A})$ satisfies $f(K)\in\mathscr{K}(S)\cap\mathscr{A}$
for all $K\in\mathscr{K}(E)\cap\mathscr{U}$, then $\{\mu\circ f^{-1}\}_{\mu\in\Gamma}$
is tight (resp. $\mathbf{m}$-tight) in $f(A)$. This implication
is also true if tightness, $\mathscr{K}(S)$ and $\mathscr{K}(S)$
are replaced by $\mathbf{m}$-tightness, $\mathscr{K}^{\mathbf{m}}(E)$
and $\mathscr{K}^{\mathbf{m}}(S)$, respectively.
\end{fact}

\begin{lem}
\label{lem:var(X)}Let $(E,\mathscr{U})$ be a measurable space, $S_{0}\subset S\subset E$,
$y_{0}\in S$, $X$ be a mapping from $(\Omega,\mathscr{F},\mathbb{P})$%
\footnote{We mentioned in \S \ref{sec:Convention} that $(\Omega,\mathscr{F},\mathbb{P})$
denotes a complete probability space. Completeness of measure space
was specified in \ref{sub:Meas}.%
} to $E$ and $Y\circeq\mathfrak{var}(X;\Omega,X^{-1}(S),y_{0})$.
Then:

\renewcommand{\labelenumi}{(\alph{enumi})}
\begin{enumerate}
\item If $\mathbb{P}(X=Z)=1$ for some $Z\in M(\Omega,\mathscr{F};S,\mathscr{U}|_{S})$,
then $X\in M(\Omega,\mathscr{F};E,\mathscr{U})$.
\item If $X\in M(\Omega,\mathscr{F};E,\mathscr{U})$ satisfies $\mathbb{P}(X\in S_{0})=1$,
then $\mathbb{P}(X=Y\in S)=1$.
\item If, in addition to the condition of (b), $(S,\mathscr{U}^{\prime})$
is a measurable space satisfying $\mathscr{U}^{\prime}|_{S_{0}}=\mathscr{U}|_{S_{0}}$,
then $Y\in M(\Omega,\mathscr{F};S,\mathscr{U}^{\prime})$.
\end{enumerate}
\end{lem}
\begin{proof}
(a) Let $\Omega_{0}\circeq\{\omega\in\Omega:X(\omega)=Z(\omega)\}$.
It follows by $\mathbb{P}(X=Z)=1$ and the completeness of $(\Omega,\mathscr{F},\mathbb{P})$
that $\Omega_{0}\in\mathscr{F}$ and $X^{-1}(A)\backslash\Omega_{0}\in\mathscr{N}(\mathbb{P})\subset\mathscr{F}$
for all $A\in\mathscr{U}$. Hence, we have that
\begin{equation}
X^{-1}(A)=\left[Z^{-1}(A\cap S)\cap\Omega_{0}\right]\cup\left(X^{-1}(A)\backslash\Omega_{0}\right)\in\mathscr{F},\;\forall A\in\mathscr{U}.\label{eq:Check_RV_Ind_Measurable}
\end{equation}

(b) We find $\mathbb{P}(X=Y\in S)\geq\mathbb{P}(X\in S_{0})=1$ by
Fact \ref{fact:var(f)} (a) (with $(S,\mathscr{A})=(\Omega,\mathscr{F})$
and $A=S$) and the completeness of $(\Omega,\mathscr{F},\mathbb{P})$.

(c) We fix $A\in\mathscr{U}^{\prime}$ and find $A\cap S_{0}\in\mathscr{U}|_{S_{0}}$
by $\mathscr{U}^{\prime}|_{S_{0}}=\mathscr{U}|_{S_{0}}$. It follows
by $\mathbb{P}(X\in S_{0})=1$ and the completeness of $(\Omega,\mathscr{F},\mathbb{P})$
that $X^{-1}(S_{0})\in\mathscr{F}$ and $Y^{-1}(A\cap S)\backslash X^{-1}(S_{0})\in\mathscr{N}(\mathbb{P})\subset\mathscr{F}$.
Hence, we have that
\begin{equation}
Y^{-1}(A\cap S)=\left[X^{-1}(A)\cap X^{-1}(S_{0})\right]\cup\left[Y^{-1}(A\cap S)\backslash X^{-1}(S_{0})\right]\in\mathscr{F}.\label{eq:Check_(S,y0)_Ver_Measurable}
\end{equation}
\end{proof}

\begin{fact}
\label{fact:Proc_Path_Mapping}Let $E$ and $S$ be topological spaces,
$\{X^{i}\}_{i\in\mathbf{I}}$ and $X$ be $E$-valued processes defined
on stochastic basis $(\Omega,\mathscr{F},\{\mathscr{G}_{t}\}_{t\geq0},\mathbb{P})$
and $f\in M(E;S)$. Then:

\renewcommand{\labelenumi}{(\alph{enumi})}
\begin{enumerate}
\item If $X$ is a general, $\mathscr{G}_{t}$-adapted, measurable or $\mathscr{G}_{t}$-progressive
process, then $\varpi(f)\circ X$ is an $S$-valued process with the
corresponding measurability.
\item If $\{X^{i}\}_{i\in\mathbf{I}}$ are general, $\mathscr{G}_{t}$-adapted,
measurable or $\mathscr{G}_{t}$-progressive processes, then $\{\bigotimes_{i\in\mathbf{I}}X_{t}^{i}\}_{t\geq0}$
is an $E^{\mathbf{I}}$-valued process with the corresponding measurability.
\end{enumerate}
\end{fact}
\begin{proof}
This result follows by Fact \ref{fact:Path_Mapping} (b), Fact \ref{fact:Prod_Map_1}
(b), Fact \ref{fact:Proc_Basic_1} (b) and the definitions of measurable
processes, $\mathscr{G}_{t}$-progressive processes and product topology.\end{proof}

\begin{prop}
\label{prop:Proc_Basic_2}Let $E$ be a topological space and $X$
and $Y$ be $E$-valued processes defined on stochastic basis $(\Omega,\mathscr{F},\{\mathscr{G}_{t}\}_{t\geq0},\mathbb{P})$%
\footnote{The notion of stochastic basis was specified in \S \ref{sec:Proc}.%
}. Then:

\renewcommand{\labelenumi}{(\alph{enumi})}
\begin{enumerate}
\item If $X$ has c$\grave{\mbox{a}}$dl$\grave{\mbox{a}}$g paths (resp.
is a c$\grave{\mbox{a}}$dl$\grave{\mbox{a}}$g process), then it
is (resp. is indistinguishable from) an $E$-valued progressive process.
\item If $X$ is progressive and $\mathscr{G}_{t}$-adapted, then it is
$\mathscr{G}_{t}$-progressive.
\item If $X$ is $\mathscr{G}_{t}$-progressive, then it is $\mathscr{G}_{t}$-adapted
and measurable.
\item If $X$ is measurable, then $X(\omega)\in M(\mathbf{R}^{+};E)$ for
all $\omega\in\Omega$.
\item If $X$ and $Y$ are modifications of each other, then $\mathscr{F}^{X}=\mathscr{F}^{Y}$%
\footnote{The filtrations $\mathscr{F}^{X}$ and $\mathscr{F}^{Y}$ were defined
in \S \ref{sec:Proc}.%
}.
\item If $X$ and $Y$ are indistinguishable, then they are modifications
of each other. If, in addition, $X$ is a measurable, $\mathscr{G}_{t}$-progressive,
progressive or c$\grave{\mbox{a}}$dl$\grave{\mbox{a}}$g process,
then $Y$ is also.
\item If $\inf_{t\in\mathbf{T}}\mathbb{P}(X_{t}=Y_{t})=1$ for some dense
$\mathbf{T}\subset\mathbf{R}^{+}$, and if $X$ and $Y$ are c$\grave{\mbox{a}}$dl$\grave{\mbox{a}}$g,
then $X$ and $Y$ are indistinguishable.
\item If $X$ is c$\grave{\mbox{a}}$dl$\grave{\mbox{a}}$g, then it is
indistinguishable from any of its c$\grave{\mbox{a}}$dl$\grave{\mbox{a}}$g
modifications and such modification is at most unique up to indistinguishability.
\end{enumerate}
\end{prop}
\begin{proof}
The well-known facts above are treated in standard texts like \cite[Chapter 2]{EK86},
\cite[Chapter 1]{P90} and \cite{N06} for $E$ being a Euclidean
or metric space. An inspection into their proofs shows that there
is no problem to make $E$ a general topological space.\end{proof}

\begin{fact}
\label{fact:Cadlag_Proc}Let $E$ and $S$ be topological spaces,
$\{X^{i}\}_{i\in\mathbf{I}}$ and $X$ be $E$-valued c$\grave{\mbox{a}}$dl$\grave{\mbox{a}}$g
processes defined on $(\Omega,\mathscr{F},\mathbb{P})$ and $f\in C(E;S)$.
Then:

\renewcommand{\labelenumi}{(\alph{enumi})}
\begin{enumerate}
\item $\varpi(f)\circ X$ is an $S$-valued c$\grave{\mbox{a}}$dl$\grave{\mbox{a}}$g
process.
\item If $\mathbf{I}$ is countable, then $\{\bigotimes_{i\in\mathbf{I}}X_{t}^{i}\}_{t\geq0}$
is an $E^{\mathbf{I}}$-valued c$\grave{\mbox{a}}$dl$\grave{\mbox{a}}$g
process.
\item If $S$ is a topological coarsening of $E$, then $X$ is an $S$-valued
c$\grave{\mbox{a}}$dl$\grave{\mbox{a}}$g process.
\end{enumerate}
\end{fact}
\begin{proof}
This result is immediate by Fact \ref{fact:Proc_Path_Mapping} and
Fact \ref{fact:Cadlag_Path} (b, c).\end{proof}

\begin{fact}
\label{fact:FDC_AS}Let $E$ be a topological space, $\mathbf{T}\subset\mathbf{R}^{+}$
and $\{(\Omega^{n},\mathscr{F}^{n},\mathbb{P}^{n};X^{n})\}_{n\in\mathbf{N}}$
and $(\Omega,\mathscr{F},\mathbb{P};X)$ be $E$-valued processes.
Then:

\renewcommand{\labelenumi}{(\alph{enumi})}
\begin{enumerate}
\item If $\{X^{n}\}_{n\in\mathbf{N}}$ is $(\mathbf{T},\mathcal{D})$-FDC%
\footnote{The notions of $(\mathbf{T},\mathcal{D})$-FDC and $(\mathbf{T},\mathcal{D})$-AS
were introduced in Definition \ref{def:FC}.%
}, then $\{\mathbb{E}^{n}[f\circ X_{\mathbf{T}_{0}}^{n}]\}_{n\in\mathbf{N}}$
is a convergent sequence in $\mathbf{R}$ for all $f\in\mathfrak{mc}[\Pi^{\mathbf{T}_{0}}(\mathcal{D})]$
and $\mathbf{T}_{0}\in\mathscr{P}_{0}(\mathbf{T})$.
\item If $\{X^{n}\}_{n\in\mathbf{N}}$ is $(\mathbf{T},\mathcal{D})$-AS,
then $\lim_{n\rightarrow\infty}\mathbb{E}^{n}[f\circ X_{\mathbf{T}_{0}}^{n}-f\circ X_{\mathbf{T}_{0}+c}^{n}]=0$
for all $c\in(0,\infty)$, $f\in\mathfrak{mc}[\Pi^{\mathbf{T}_{0}}(\mathcal{D})]$
and $\mathbf{T}_{0}\in\mathscr{P}_{0}(\mathbf{T})$.
\end{enumerate}
\end{fact}
\begin{proof}
This result is immediate by the Bolzano-Weierstrass Theorem and Fact
\ref{fact:Uni_Seq_Lim_Conv} (with $E=\mathbf{R}$ and $x_{n}=\mathbb{E}^{n}[f\circ X_{\mathbf{T}_{0}}^{n}]$
or $\mathbb{E}^{n}[f\circ X_{\mathbf{T}_{0}}^{n}-f\circ X_{\mathbf{T}_{0}+c}^{n}]$).\end{proof}

\begin{fact}
\label{fact:FLP_Uni}Let $E$ be a topological space, $\mathbf{T}\subset\mathbf{R}^{+}$
and $\{(\Omega^{i},\mathscr{F}^{i},\mathbb{P}^{i};X^{i})\}_{i\in\mathbf{I}}$
be $E$-valued processes. If for each $\mathbf{T}_{0}\in\mathscr{P}_{0}(\mathbf{T})$,
there exists some $\mathbf{I}_{\mathbf{T}_{0}}\in\mathscr{P}_{0}(\mathbf{I})$
such that $\{\mu_{\mathbf{T}_{0},i}=\mathfrak{be}(\mathbb{P}\circ(X_{\mathbf{T}_{0}}^{i})^{-1})\}_{i\in\mathbf{I}\backslash\mathbf{I}_{\mathbf{T}_{0}}}$
has at most one weak limit point, then $\mathfrak{flp}_{\mathbf{T}}(\{X^{i}\}_{i\in\mathbf{I}})$%
\footnote{The notation ``$\mathfrak{flp}_{\mathbf{T}}(\{X^{n}\}_{n\in\mathbf{N}})$''
was introduced in \S \ref{sec:RepProc_FC} and stands for the family
of all equivalence classes of finite-dimensional limit points of $\{X^{n}\}_{n\in\mathbf{N}}$
along $\mathbf{T}$.%
} is at most a singleton.
\end{fact}
\begin{proof}
Suppose $(\Omega,\mathscr{F},\mathbb{P};Y^{j})\in\mathfrak{flp}_{\mathbf{T}}(\{X^{i}\}_{i\in\mathbf{I}})$
for each $j=1,2$. Fixing $\mathbf{T}_{0}\in\mathscr{P}_{0}(\mathbf{T})$,
there exist $\{\nu_{j}\in\mathfrak{be}(\mathbb{P}\circ(Y_{\mathbf{T}_{0}}^{j})^{-1})\}_{j=1,2}$
such that $\nu_{1}$ and $\nu_{2}$ are both weak limit points of
$\{\mu_{\mathbf{T}_{0},i}\}_{i\in\mathbf{I}\backslash\mathbf{I}_{\mathbf{T}_{0}}}$.
So, $\nu_{1}=\nu_{2}$ and hence $\mathbb{P}\circ(Y_{\mathbf{T}_{0}}^{1})^{-1}=\mathbb{P}\circ(Y_{\mathbf{T}_{0}}^{2})^{-1}$.\end{proof}

\begin{fact}
\label{fact:Weakly_Cadlag}Let $E$ be a topological space, $(\Omega,\mathscr{F},\mathbb{P};X)$
be an $E$-valued process and $\mathbf{T}\subset\mathbf{R}^{+}$.
If there exists an $\mathbf{R}^{\mathcal{D}}$-valued c$\grave{\mbox{a}}$dl$\grave{\mbox{a}}$g
process $(\Omega,\mathscr{F},\mathbb{P};\zeta)$ such that
\begin{equation}
\inf_{t\in\mathbf{T}}\mathbb{P}\left(\bigotimes\mathcal{D}\circ X_{t}=\zeta_{t}\right)=1,\label{eq:Prod(D)_Weakly_Cadlag}
\end{equation}
then $X$ is $(\mathbf{T},\mathcal{D})$-c$\grave{\mbox{a}}$dl$\grave{\mbox{a}}$g.
The converse is true when $\mathcal{D}$ is a countable collection.
\end{fact}
\begin{proof}
$\{\zeta^{f}\circeq\varpi(\mathfrak{p}_{f})\circ\zeta\}_{f\in\mathcal{D}}$%
\footnote{Recall that $\mathfrak{p}_{f}$ denotes the projection on $\mathbf{R}^{\mathcal{D}}$
for $f\in\mathcal{D}$.%
} are $\mathbf{R}$-valued c$\grave{\mbox{a}}$dl$\grave{\mbox{a}}$g
processes satisfying
\begin{equation}
\begin{aligned} & \inf_{t\in\mathbf{T}}\mathbb{P}\left(\mathfrak{p}_{f}\circ\zeta_{t}=\zeta_{t}^{f}=f\circ X_{t}=\mathfrak{p}_{f}\circ\bigotimes\mathcal{D}\circ X_{t},\forall f\in\mathcal{D}\right)\\
 & =\inf_{t\in\mathbf{T}}\mathbb{P}\left(\zeta_{t}=\bigotimes\mathcal{D}\circ X_{t}\right)=1
\end{aligned}
\label{eq:Check_Prod(D)(X)_Mod}
\end{equation}
by Fact \ref{fact:Prod_Map_2} (a) and Fact \ref{fact:Cadlag_Proc}
(a) (with $E=\mathbf{R}^{\mathcal{D}}$ and $f=\mathfrak{p}_{f}$).

Conversely, we suppose $\mathcal{D}$ is countable and $\mathbf{R}$-valued
c$\grave{\mbox{a}}$dl$\grave{\mbox{a}}$g processes $\{\zeta^{f}\}_{f\in\mathcal{D}}$
satisfy (\ref{eq:(T,D)_Weakly_Cadlag}). Letting $\zeta_{t}\circeq\bigotimes_{f\in\mathcal{D}}\zeta_{t}^{f}$
for each $t\in\mathbf{R}^{+}$, we find that $\{\zeta_{t}\}_{t\geq0}$
is an $\mathbf{R}^{\infty}$-valued c$\grave{\mbox{a}}$dl$\grave{\mbox{a}}$g
process satisfying (\ref{eq:Check_Prod(D)(X)_Mod}) by the countability
of $\mathcal{D}$, Fact \ref{fact:Cadlag_Proc} (b) (with $\mathbf{I}=\mathcal{D}$,
$i=f$ and $X^{i}=\zeta^{f}$) and (\ref{eq:(T,D)_Weakly_Cadlag}).\end{proof}

\begin{lem}
\label{lem:RAP_Seq_Tight}Let $E$ be a topological space, $\{(\Omega^{i},\mathscr{F}^{i},\mathbb{P}^{i};X^{i})\}_{i\in\mathbf{I}}$
be $E$-valued measurable processes, $T_{k}\uparrow\infty$, $\{A_{p}\}_{p\in\mathbf{N}}\subset\mathscr{B}(E)$
and $(\widetilde{\Omega}^{i},\widetilde{\mathscr{F}}^{i},\mathbb{P}^{i,T_{k}};X^{i,T_{k}})=\mathfrak{rap}_{T_{k}}(X^{i})$%
\footnote{Randomly advanced process and related notations were introduced in
\S \ref{sec:LTB}.%
} for each $i\in\mathbf{I}$ and $k\in\mathbf{N}$. If\textup{
\begin{equation}
\inf_{i\in\mathbf{I},k\in\mathbf{N}}\mathbb{P}^{i,T_{k}}\left(X_{0}^{i,T_{k}}\in A_{p}\right)\geq1-2^{-p},\;\forall p\in\mathbf{N},\label{eq:Xi_RAP_Ap}
\end{equation}
}then for each $\mathbf{T}_{0}\in\mathscr{P}_{0}(\mathbf{R}^{+})$
with $d\circeq\aleph(\mathbf{T}_{0})$, there are $\{N_{\mathbf{T}_{0},p}\}_{p\in\mathbf{N}}\subset\mathbf{N}$
such that 
\begin{equation}
\inf_{i\in\mathbf{I},k>N_{\mathbf{T}_{0},p}}\mathbb{P}^{i,T_{k}}\left(X_{\mathbf{T}_{0}}^{i,T_{k}}\in A_{p}^{d}\right)\geq1-(d+1)2^{-p},\;\forall p\in\mathbf{N}.\label{eq:RAP_Seq_Tight}
\end{equation}

\end{lem}
\begin{proof}
Let $\mathbf{T}_{0}=\{t_{1},...,t_{d})$, $t\circeq\sum_{j=1}^{d}t_{j}$,
$N_{\mathbf{T}_{0},0}\circeq0$ and $T_{N_{\mathbf{T}_{0},0}}\circeq0$.
Define $\{N_{\mathbf{T}_{0},p}\}_{p\in\mathbf{N}}$ inductively by
\begin{equation}
N_{\mathbf{T}_{0},p}\circeq\min\left\{ k\in\mathbf{N}:T_{k}>(2^{p+1}td)\vee T_{N_{\mathbf{T}_{0},p-1}}\right\} ,\;\forall p\in\mathbf{N}.\label{eq:Choose_N_T0_p}
\end{equation}
For each $p\in\mathbf{N}$, it follows by (\ref{eq:Xi_RAP_Ap}) and
(\ref{eq:Choose_N_T0_p}) that
\begin{equation}
\begin{aligned} & \inf_{i\in\mathbf{I},k>N_{\mathbf{T}_{0},p}}\mathbb{P}^{i,T_{k}}\left(X_{\mathbf{T}_{0}}^{i,T_{k}}\in A_{p}^{d}\right)\\
 & \geq1-\sup_{i\in\mathbf{I},k>N_{\mathbf{T}_{0},p}}\sum_{j=1}^{d}\mathbb{P}^{i,T_{k}}\left(X_{t_{j}}^{i,T_{k}}\notin A_{p}\right)\\
 & \ge1-d\sup_{i\in\mathbf{I},k>N_{\mathbf{T}_{0},p}}\mathbb{P}^{i,T_{k}}\left(X_{0}^{i,T_{k}}\notin A_{p}\right)\\
 & -\sup_{i\in\mathbf{I},k>N_{\mathbf{T}_{0},p}}\frac{d}{T_{k}}\int_{[0,t]\cup[T_{k},T_{k}+t]}\mathbb{P}^{i}\left(X_{\tau}^{i}\notin A_{p}\right)d\tau\\
 & \geq1-d2^{-p}-\frac{2td}{2^{p+1}td}=1-(d+1)2^{-p}.
\end{aligned}
\label{eq:Check_RAP_ST}
\end{equation}
\end{proof}

\begin{lem}
\label{lem:RAP_Cadlag}Let $E$ be a topological space, $(\Omega,\mathscr{F},\mathbb{P};X)$
be an $E$-valued c$\grave{\mbox{a}}$dl$\grave{\mbox{a}}$g process,
$\epsilon,\delta,T,c\in(0,\infty)$ and $(\widetilde{\Omega},\widetilde{\mathscr{F}},\mathbb{P}^{T};X^{T})=\mathfrak{rap}_{T}(X)$.
Then:

\renewcommand{\labelenumi}{(\alph{enumi})}
\begin{enumerate}
\item $\xi^{\tau}\circeq\{X_{\tau+t}\}_{t\geq0}$ well defines an $E$-valued
c$\grave{\mbox{a}}$dl$\grave{\mbox{a}}$g process for all $\tau\in\mathbf{R}^{+}$.
\item $X^{T}$ is an $E$-valued c$\grave{\mbox{a}}$dl$\grave{\mbox{a}}$g
process.
\item If $(E,\mathfrak{r})$ is a separable metric space, then
\begin{equation}
\frac{1}{T}\int_{0}^{T}\mathbb{P}\left(w_{\mathfrak{r},\delta,c}^{\prime}\circ\xi^{\tau}\geq\epsilon\right)d\tau=\mathbb{P}^{T}\left(w_{\mathfrak{r},\delta,c}^{\prime}\circ X^{T}\geq\epsilon\right).\label{eq:RAP_MCC_Measurability}
\end{equation}

\end{enumerate}
\end{lem}
\begin{proof}
$\{\xi^{\tau}\}_{\tau\in\mathbf{R}^{+}}$ are $E$-valued processes
by Fact \ref{fact:Proc_Basic_1} (b). Letting $\Omega_{0}\circeq\{\omega\in\Omega:X(\omega)\mbox{ is c}\grave{\mbox{a}}\mbox{dl}\grave{\mbox{a}}\mbox{g}\}$,
we find that
\begin{equation}
\left\{ \omega\in\Omega:\xi^{\tau}(\omega)\mbox{ is c}\grave{\mbox{a}}\mbox{dl}\grave{\mbox{a}}\mbox{g}\right\} \supset\Omega_{0},\;\forall\tau\in\mathbf{R}^{+}\label{eq:Proc_Shifted_Cadlag}
\end{equation}
and
\begin{equation}
\left\{ (\tau,\omega)\in\widetilde{\Omega}:X^{T}(\tau,\omega)\mbox{ is c}\grave{\mbox{a}}\mbox{dl}\grave{\mbox{a}}\mbox{g}\right\} \supset\mathbf{R}^{+}\times\Omega_{0}.\label{eq:Check_RAP_Cadlag}
\end{equation}
Then, (a, b) follows by (\ref{eq:Proc_Shifted_Cadlag}), (\ref{eq:Check_RAP_Cadlag})
and the fact $\mathbb{P}^{T}(\mathbf{R}^{+}\times\Omega_{0})=\mathbb{P}(\Omega_{0})=1$.

(c) It follows by (a), (b) and Lemma \ref{lem:MCC_Measurability}
(b) (with $X=\xi^{\tau}$ or $X^{T}$) that $w_{\mathfrak{r},\delta,c}^{\prime}\circ\xi^{\tau}\in M(\Omega,\mathscr{F};\mathbf{R})$
for all $\tau\in\mathbf{R}^{+}$ and $w_{\mathfrak{r},\delta,c}^{\prime}\circ X^{T}\in M(\widetilde{\Omega},\widetilde{\mathscr{F}};\mathbf{R})$.
Hence, both sides of (\ref{eq:RAP_MCC_Measurability}) are all well-defined
and (\ref{eq:RAP_MCC_Measurability}) is true since $\xi^{\tau}(\omega)=X^{T}(\tau,\omega)$
for all $(\tau,\omega)\in\widetilde{\Omega}$.\end{proof}

\begin{lem}
\label{lem:RAP_Tight}Let $E$ be a topological space, $A\subset E$,
$\mathcal{D}\subset M_{b}(E;\mathbf{R})$ and $(\Omega,\mathscr{F},\mathbb{P};X)$
be an $E$-valued measurable process, $T_{k}\uparrow\infty$ and $(\widetilde{\Omega},\widetilde{\mathscr{F}},\mathbb{P}^{T_{k}};X^{T_{k}})=\mathfrak{rap}_{T_{k}}(X)$
for each $k\in\mathbf{N}$. Then:

\renewcommand{\labelenumi}{(\alph{enumi})}
\begin{enumerate}
\item $X$ satisfies $T_{k}$-LMTC in $A$%
\footnote{The terminology ``$X$ satisfying $T_{k}$-LMTC in $A$'' was set
in Definition \ref{def:Proc_Reg} and Note \ref{note:Singleton_PMTC}.%
} if and only if $\{X_{0}^{T_{k}}\}_{n\in\mathbf{N}}$ is $\mathbf{m}$-tight
in $A$.
\item If $\{X_{0}^{T_{k}}\}_{n\in\mathbf{N}}$ is $\mathbf{m}$-tight in
$A$, then $\{X^{T_{k}}\}_{k\in\mathbf{N}}$ satisfies $\mathbf{R}^{+}$-PSMTC
in $A$%
\footnote{The notion of $\mathbf{R}^{+}$-PSMTC was introduced in Definition
\ref{def:Proc_Reg}.%
}.
\item $\{X^{T_{k}}\}_{k\in\mathbf{N}}$ is $(\mathbf{R}^{+},M_{b}(E;\mathbf{R}))$-AS.
\item If $\{X^{T_{k}}\}_{k\in\mathbf{N}}$is $(\mathbf{T},\mathcal{D})$-FDC,
then it is $(\mathbf{T}+c,\mathcal{D})$-FDC for all $c\in(0,\infty)$.
\end{enumerate}
\end{lem}
\begin{proof}
(a) is automatic by the definition of $\{X^{T_{k}}\}_{k\in\mathbf{N}}$.

(b) follows by (a) and Lemma \ref{lem:RAP_Seq_Tight} (with $\{X^{i}\}_{i\in\mathbf{I}}=\{X\}$).

(c) and (d) follow immediately by the fact that
\begin{equation}
\begin{aligned} & \lim_{k\rightarrow\infty}\left|\mathbb{E}^{T_{k}}\left[f\circ X_{\mathbf{T}_{0}}^{T_{k}}-f\circ X_{\mathbf{T}_{0}+c}^{T_{k}}\right]\right|\\
 & \leq\lim_{k\rightarrow\infty}\frac{1}{T_{k}}\int_{0}^{c}\mathbb{E}\left[\left|f\circ X_{\mathbf{T}_{0}}-f\circ X_{\mathbf{T}_{0}+\tau+T_{k}}\right|\right]d\tau\leq\lim_{k\rightarrow\infty}\frac{2c\Vert f\Vert_{\infty}}{T_{k}}=0
\end{aligned}
\label{eq:RAP_AS}
\end{equation}
for all $c\in(0,\infty)$, $f\in\Pi^{\mathbf{T}_{0}}(M_{b}(E;\mathbf{R}))$
and $\mathbf{T}_{0}\in\mathscr{P}_{0}(\mathbf{R}^{+})$, where $\mathbb{E}^{T_{k}}$
denotes the expectation operator of $(\widetilde{\Omega},\widetilde{\mathscr{F}},\mathbb{P}^{T_{k}})$
for each $k\in\mathbf{N}$.\end{proof}

\section{\label{sec:Comp_A1}Supplementary results for Appendix \ref{chap:App1}}
\begin{fact}
\label{fact:Cc_C0_Cb}Let $E$ be a topological space and $k\in\mathbf{N}$.
Then, $C_{c}(E;\mathbf{R}^{k})\subset C_{0}(E;\mathbf{R}^{k})\subset C_{b}(E;\mathbf{R}^{k})$
and they are indifferent if $E$ is compact.
\end{fact}
\begin{proof}
This result follows by \cite[Theorem 27.4]{M00}.\end{proof}

\begin{prop}
\label{prop:Cc_Lattice}Let $E$ be a Hausdorff space. Then:

\renewcommand{\labelenumi}{(\alph{enumi})}
\begin{enumerate}
\item $C_{c}(E;\mathbf{R})$ is a subalgebra of $C_{b}(E;\mathbf{R})$ and
is a function lattice.
\item $C_{c}(E;\mathbf{R})\subset C_{0}(E;\mathbf{R})\subset\mathfrak{cl}(C_{c}(E;\mathbf{R}))$.
\end{enumerate}
\end{prop}
\begin{proof}
It is straightforward to show $C_{c}(E;\mathbf{R})$ is an algebra.
Observing that
\begin{equation}
\left|f\right|\in C_{c}(E;\mathbf{R}),\;\forall f\in C_{c}(E;\mathbf{R}),\label{eq:Check_Cc_Lattice_2}
\end{equation}
we have that
\begin{equation}
f\vee g=\frac{1}{2}\left(f+g\right)+\frac{1}{2}\left|f-g\right|\in C_{c}(E;\mathbf{R}),\;\forall f,g\in C_{c}(E;\mathbf{R}),\label{eq:Check_Cc_Lattice_3}
\end{equation}
and that
\begin{equation}
f\wedge g=\frac{1}{2}\left(f+g\right)-\frac{1}{2}\left|f-g\right|\in C_{c}(E;\mathbf{R}),\;\forall f,g\in C_{c}(E;\mathbf{R}),\label{eq:Check_Cc_Lattice_4}
\end{equation}
thus proving $C_{c}(E;\mathbf{R})$ is a function lattice.

(b) The first inclusion is immediate. We prove the second one. We
fix $f\in C_{0}(E;\mathbf{R})$, $p\in\mathbf{N}$ and a $K_{p}\in\mathscr{K}(E)$
such that $\Vert f|_{E\backslash K_{p}}\Vert_{\infty}<2^{-p}$. The
case where $K_{p}=E$ is trivial. Otherwise, we define $A\circeq f^{-1}((2^{-p},\infty))$,
$B\circeq f^{-1}[(-\infty,-2^{-p})]$ and 
\begin{equation}
f_{p}\circeq\left(f^{+}(x)-2^{-p}\right)^{+}-\left(f^{-}(x)-2^{-p}\right)^{+}.\label{eq:Check_Cc_Dense_1}
\end{equation}
$A$ and $B$ are disjoint subsets of $K_{p}$ such that $A\cup B=E\backslash f^{-1}(\{0\})$.
Letting $F$ be the closure of $A\cup B$ in $E$, we have by Proposition
\ref{prop:Compact} (a) that $K_{p}\in\mathscr{C}(E)$, $F\subset K_{p}$
and $F\in\mathscr{K}(E)$, thus proving $f_{p}\in C_{c}(E;\mathbf{R})$.
Furthermore, from the fact
\begin{equation}
f_{p}(x)-f(x)=\begin{cases}
-2^{-p}, & \mbox{if }x\in A,\\
2^{-p}, & \mbox{if }x\in B,\\
-f(x)\in(-2^{-p},2^{-p},), & \mbox{if }x\in E\backslash(A\cup B)
\end{cases}\label{eq:Check_Cc_Dense_4}
\end{equation}
it follows that $\Vert f_{p}-f\Vert_{\infty}\leq2^{-p}$.\end{proof}

\begin{fact}
\label{fact:Cc_Ext}Let $E$ be a topological space, $f\in C(E;\mathbf{R})$
and $A$ be a dense subset of $E$ with $E\backslash A\neq\varnothing$.
Then:

\renewcommand{\labelenumi}{(\alph{enumi})}
\begin{enumerate}
\item If $E$ is a first-countable space and $(A\backslash B)\subset f^{-1}(\{0\})$
for some $B\in\mathscr{C}(E)$ with $B\subset A$, then $f|_{E\backslash A}=0$.
\item If $E$ is a metrizable space and $f\in C_{c}(A,\mathscr{O}_{E}(A))$,
then $f|_{E\backslash A}=0$.
\end{enumerate}
\end{fact}
\begin{proof}
(a) For each $x\in E\backslash A$, the first-countability of $E$
implies a sequence $\{x_{n}\}_{n\in\mathbf{N}}\subset A$ converging
to $x$ as $n\uparrow\infty$. As $x\in E\backslash B$ and $E\backslash B\in\mathscr{O}(E)$,
there exists an $N\in\mathbf{N}$ such that $x_{n}\in A\backslash B$
and $f(x_{n})=0$ for all $n>N$. Hence, the continuity of $f$ implies
$f(x)=\lim_{n\rightarrow\infty}f(x_{n})=0$.

(b) If $E$ is a metrizable space, then it is first-countable by Fact
\ref{fact:First_Countable}. If, in addition, $f\in C_{c}(A,\mathscr{O}_{E}(A))$,
then we let $B$ be the closure of $A\backslash f^{-1}(\{0\})$ and
$B\in\mathscr{C}(E)$ by Proposition \ref{prop:Metrizable} (a) and
Proposition \ref{prop:Compact} (a).\end{proof}

\begin{prop}
\label{prop:Prod_Space}Let $\{S_{i}\}_{i\in\mathbf{I}}$ be topological
spaces and $(S,\mathscr{A})$ be defined as in (\ref{eq:(S,A)_Prod_Meas_Space}).
Then:

\renewcommand{\labelenumi}{(\alph{enumi})}
\begin{enumerate}
\item $\mathscr{B}(S)\supset\mathscr{A}$.
\item If $\mathbf{I}$ is countable and $S$ is hereditary Lindel$\ddot{\mbox{o}}$f,
then $\mathscr{B}(S)=\mathscr{A}$.
\item If $\mathbf{I}$ is countable and $\{S_{i}\}_{i\in\mathbf{I}}$ are
all second-countable, then $\mathscr{B}(S)=\mathscr{A}$.
\item If $\mathbf{I}$ is countable and $\{S_{i}\}_{i\in\mathbf{I}}$ are
all metrizable and separable spaces (especially Polish spaces), then
$\mathscr{B}(S)=\mathscr{A}$.
\end{enumerate}
\end{prop}
\begin{note}
\label{note:Prod_Space}As arranged in \S \ref{sec:Convention},
the Cartesian product $S\circeq\prod_{i\in\mathbf{I}}S_{i}$ above
is equipped with the product topology $\mathscr{O}(S)\circeq\bigotimes_{i\in\mathbf{I}}\mathscr{O}(S_{i})$
and its Borel $\sigma$-algebra is $\mathscr{B}(S)\circeq\sigma[\mathscr{O}(S)]$.
\end{note}
\begin{proof}
[Proof of Proposition \ref{prop:Prod_Space}](a) follows by the argument
establishing \cite[Vol. II, Lemma 6.4.1]{B07}.

(b) and (c) were proved in \cite[Vol. II, Lemma 6.4.2 (ii)]{B07}.

(d) follows by (c), Proposition \ref{prop:Metrizable} (c) and Proposition
\ref{prop:Var_Polish} (c).\end{proof}

\begin{lem}
\label{lem:Union_Borel_Prod_Equal}Let $\{S_{i}\}_{i\in\mathbf{I}}$
be topological spaces, $(S,\mathscr{A})$ be as in (\ref{eq:(S,A)_Prod_Meas_Space}),
$A\in\mathscr{A}$, $\{A_{n}\}_{n\in\mathbf{N}}\subset\mathscr{A}$
and $\mu\in\mathfrak{M}^{+}(S,\mathscr{A})$. Then:

\renewcommand{\labelenumi}{(\alph{enumi})}
\begin{enumerate}
\item If $S=\bigcup_{n\in\mathbf{N}}A_{n}$ and $\mathscr{A}|_{A_{n}}=\mathscr{B}_{S}(A_{n})$
for all $n\in\mathbf{N}$, then $\mathscr{B}(S)=\mathscr{A}$.
\item If $\mu$ is supported on $A$ and $\mathfrak{be}(\mu|_{A})$%
\footnote{``$\mu|_{A}$'' and $\nu|^{E}$ denote the concentration of $\mu$
on $A$ and the expansion of $\nu$ onto $E$. ``$\mathfrak{be}(\mu|_{A})$''
denotes the Borel extension(s) of $\mu$.%
} is a singleton, then $\mathfrak{be}(\mu)$ is a singleton.
\item If $\mu$ is supported on $A$ and $\mathscr{B}_{S}(A)=\mathscr{A}|_{A}$,
then $\mathfrak{be}(\mu)$ is a singleton.
\end{enumerate}
\end{lem}
\begin{proof}
(a) follows by Proposition \ref{prop:Prod_Space} (a) and Fact \ref{fact:Union_Borel}
(with $E=S$, $\mathscr{U}_{1}=\mathscr{B}(S)$ and $\mathscr{U}_{2}=\mathscr{A}$).

(b) Let $\nu=\mathfrak{be}(\mu|_{A})$%
\footnote{``$\nu=\mathfrak{be}(\mu|_{A})$'' means $\nu$ is the unique Borel
extension of $\mu|_{A}$.%
}. $\mu_{1}\circeq\nu|^{S}\in\mathcal{M}^{+}(S)$ by Fact \ref{fact:Meas_Concen_Expan}
(b) (with $E=S$ and $\mathscr{U}=\mathscr{B}(E)$). Since $A\in\mathscr{A}$,
we have that
\begin{equation}
B\cap A\in\mathscr{A}|_{A}\subset\mathscr{B}_{S}(A),\;\forall B\in\mathscr{A}.\label{eq:Check_Expansion_BExt_1}
\end{equation}
$\nu$ and $\mu|_{A}$ are identical restricted to $\mathscr{A}|_{A}$.
It then follows by (\ref{eq:Check_Expansion_BExt_1}) that
\begin{equation}
\mu|_{A}(B\cap A)=\nu(B\cap A)=\mu_{1}(B),\;\forall B\in\mathscr{A}.\label{eq:Check_Expansion_BExt_2}
\end{equation}
It follows by the fact $\mu(A)=1$, the fact $A\in\mathscr{A}$ and
(\ref{eq:Check_Expansion_BExt_2}) that
\begin{equation}
\mu(B)=\mu(B\cap A)=\mu|_{A}(B\cap A)=\mu_{1}(B),\;\forall B\in\mathscr{A},\label{eq:Check_Expansion_BExt_3}
\end{equation}
thus proving $\mu_{1}\in\mathfrak{be}(\mu)$. If $\mu_{2}\in\mathfrak{be}(\mu)$,
then we have that
\begin{equation}
\mu_{2}|_{A}=\mathfrak{be}(\mu|_{A})=\nu\in\mathscr{M}^{+}\left(A,\mathscr{O}_{S}(A)\right).\label{eq:Check_Expansion_BExt_Uni_1}
\end{equation}
It follows that
\begin{equation}
\mu_{2}=(\mu_{2}|_{A})|^{S}=\nu|^{S}=\mu_{1}\label{eq:Check_Expansion_BExt_Uni_2}
\end{equation}
by (\ref{eq:Check_Expansion_BExt_Uni_1}) and Fact \ref{fact:Meas_Concen_Expan}
(a, c) (with $E=S$, $\mathscr{U}=\mathscr{B}(E)$ and $\mu=\mu_{2}$).

(c) $\mathscr{B}_{S}(A)=\mathscr{A}|_{A}$ implies $\mu|_{A}=\mathfrak{be}(\mu|_{A})$.
Then, (c) follows by (b).\end{proof}

\begin{lem}
\label{lem:f=00003Dg_Measurable_Set}Let $E$ and $S$ be measurable
spaces and $f,g\in M(S;E)$. If there exists a countable subset of
$M(E;\mathbf{R})$ separating points on $E$, then $\{x\in S:f(x)=g(x)\}$
is a measurable subset of $S$. In particular, this is true when $E$
is baseable.
\end{lem}
\begin{proof}
Let $\{h_{n}\}_{n\in\mathbf{N}}\subset M(E;\mathbf{R})$ separate
points on $E$. Then,
\begin{equation}
\begin{aligned}\left\{ x\in S:f(x)=g(x)\right\}  & =\left\{ x\in S:h_{n}\circ f(x)=h_{n}\circ g(x)\right\} \\
 & =\bigcap_{n\in\mathbf{N}}\left(h_{n}\circ f-h_{n}\circ g\right)^{-1}(\{0\})
\end{aligned}
\label{eq:Check_f=00003Dg_Measurable_2}
\end{equation}
is a measurable subset of $S$.\end{proof}

\begin{lem}
\label{lem:Cadlag_Measurability}Let $E$ be a topological space,
$V$ be the family of all c$\grave{\mbox{a}}$dl$\grave{\mbox{a}}$g
members of $E^{\mathbf{R}^{+}}$%
\footnote{$E$ need not be a Tychonoff space, so we avoid the notation $D(\mathbf{R}^{+};E)$
for clarity.%
}, $t\in\mathbf{R}^{+}$ and $T\in(0,\infty)$. Then:

\renewcommand{\labelenumi}{(\alph{enumi})}
\begin{enumerate}
\item If $M(E;\mathbf{R})$ has a countable subset separating points on
$E$, then $\{x\in V:t\in J(x)\}\in\mathscr{B}(E)^{\otimes\mathbf{R}^{+}}|_{V}$.
\item $\{x\in V:x|_{[0,T)}\in A^{[0,T)}\}\in\mathscr{B}(E)^{\otimes\mathbf{R}^{+}}|_{V}$
for all $A\in\mathscr{C}(E)$, especially for all $A\in\mathscr{K}(E)$
when $E$ is a Hausdorff space.
\end{enumerate}
\end{lem}
\begin{proof}
(a) We fix $t\in\mathbf{R}^{+}$, let $\mathfrak{p}_{t-}$ denote
the mapping associating each $x\in V$ to its left limit at $t$ and
find by Fact \ref{fact:Prod_Map_1} (a) that
\begin{equation}
\mathfrak{p}_{t-}^{-1}(A)=\bigcap_{p\in\mathbf{Q}^{+}\cap[0,t)}\bigcup_{q\in\mathbf{Q}^{+}\cap(p,t)}\mathfrak{p}_{q}^{-1}(A)\in\mathscr{B}(E)^{\otimes\mathbf{R}^{+}}|_{V},\;\forall A\in\mathscr{O}(E),\label{eq:Check_Left_Limit_Map_Measurable}
\end{equation}
so $\mathfrak{p}_{t-}\in M(V,\mathscr{B}(E)^{\otimes\mathbf{R}^{+}}|_{V};E)$.
It then follows by Lemma \ref{lem:f=00003Dg_Measurable_Set} (with
$S=V$, $f=\mathfrak{p}_{t-}$ and $g=\mathfrak{p}_{t}$) that
\begin{equation}
\left\{ x\in V:t\in J(x)\right\} =\left\{ x\in V:\mathfrak{p}_{t-}(x)=\mathfrak{p}_{t}(x)\right\} \in\mathscr{B}(E)^{\otimes\mathbf{R}^{+}}|_{V}.\label{eq:Check_J(Mu)_Well_Defined}
\end{equation}

(b) When $E$ is Hausdorff, $\mathscr{K}(E)\subset\mathscr{C}(E)$
by Proposition \ref{prop:Compact} (a). It follows by the closedness
of $A$ and the right-continuity of each $x$ that
\begin{equation}
\left\{ x\in V:x|_{[0,T)}\in A^{[0,T)}\right\} =\bigcap_{t\in\mathbf{Q}\cap[0,T)}V\cap\mathfrak{p}_{t}^{-1}(A).\label{eq:Check_CCC_Set_Measurability_1}
\end{equation}
It follows by the fact $A\in\mathscr{C}(E)\subset\mathscr{B}(E)$
and Fact \ref{fact:Prod_Map_1} (a) that
\begin{equation}
V\cap\mathfrak{p}_{t}^{-1}(A)\in\mathscr{B}(E)^{\otimes\mathbf{R}^{+}}|_{V},\;\forall t\in\mathbf{R}^{+}.\label{eq:Check_CCC_Set_Measurability_2}
\end{equation}
Now, (b) follows by (\ref{eq:Check_CCC_Set_Measurability_1}), (\ref{eq:Check_CCC_Set_Measurability_2})
and the countability of $\mathbf{Q}\cap[0,T)$.\end{proof}

\begin{fact}
\label{fact:Compact_Topo_Coarsen}Let $E$ is a topological space
and $K\in\mathscr{K}(E)$. Then:

\renewcommand{\labelenumi}{(\alph{enumi})}
\begin{enumerate}
\item $K\in\mathscr{K}(E,\mathscr{U})$ for any topological coarsening $(E,\mathscr{U})$
of $E$.
\item If $\mathcal{D}\subset C(E;\mathbf{R})$ separate points on $K$,
then $\mathscr{O}_{E}(K)=\mathscr{O}_{\mathcal{D}}(K)$ and $K\in\mathscr{K}(E,\mathscr{O}_{\mathcal{D}}(E))\subset\mathscr{C}(E,\mathscr{O}_{\mathcal{D}}(E))$.
\end{enumerate}
\end{fact}
\begin{proof}
(a) is immediate by the definition of compactness.

(b) $\mathscr{O}_{E}(K)=\mathscr{O}_{\mathcal{D}}(K)$ is a Hausdorff
topology by Lemma \ref{lem:SP_on_Compact} (with $E=K$ and $\mathcal{D}=\mathcal{D}|_{K}$)
and Proposition \ref{prop:Fun_Sep_1} (c) (with $A=E$). Now, (b)
follows by Proposition \ref{prop:Compact} (a).\end{proof}

\begin{lem}
\label{lem:SP_Bounded}Let $E$ be a non-empty set, $\mathcal{G}\subset\mathbf{R}^{E}$
and $\mathcal{H}\subset\mathbf{R}^{E}$. Suppose that for any $g\in\mathcal{G}$
and $n\in\mathbf{N}$, there exists a bounded function $f_{g,n}\in\mathcal{H}$
such that
\begin{equation}
A_{g,n}\circeq\left\{ x\in E:\left|g(x)\right|<n\right\} =\left\{ x\in E:\left|f_{g,n}(x)\right|<n\right\} \label{eq:A_gn}
\end{equation}
and
\begin{equation}
g\mathbf{1}_{A_{g,n}}=f_{g,n}\mathbf{1}_{A_{g,n}}.\label{eq:f_gn}
\end{equation}
Then, there exists a subset $\mathcal{F}\subset\mathcal{H}$ such
that:

\renewcommand{\labelenumi}{(\alph{enumi})}
\begin{enumerate}
\item The members of $\mathcal{F}$ are all bounded and include all the
bounded members of $\mathcal{G}$. In particular, $\mathcal{F}=\mathcal{G}$
when the members of $\mathcal{G}$ are all bounded.
\item $\mathcal{F}$ is countable if $\mathcal{G}$ is.
\item $\mathscr{O}_{\mathcal{G}}(E)\subset\mathscr{O}_{\mathcal{F}}(E)$.
Moreover, if $\mathcal{G}$ separates points on $E$, or if $E$ is
a topological space and $\mathcal{G}$ strongly separates points on
$E$, then $\mathcal{F}$ has the same property.
\end{enumerate}
\end{lem}
\begin{proof}
(a, b) are immediate.

(c) It follows by (\ref{eq:A_gn}) and (\ref{eq:f_gn}) that 
\begin{equation}
\begin{aligned} & \left\{ x\in E:g(x)<a\right\} =\bigcup_{n>a}\left\{ x\in E:g(x)<a,\left|g(x)\right|<n\right\} \\
 & =\bigcup_{n>a}\left\{ x\in E:f_{g,n}(x)<a,\left|f_{g,n}(x)\right|<n\right\} \\
 & =\bigcup_{n>a}\left\{ x\in E:-n<f_{g,n}(x)<a\right\} \in\mathscr{O}_{\mathcal{F}}(E),\;\forall a\in\mathbf{R},g\in\mathcal{G},
\end{aligned}
\label{eq:C_Cb_Same_Topo}
\end{equation}
thus proving $\mathscr{O}_{\mathcal{G}}(E)\subset\mathscr{O}_{\mathcal{F}}(E)$.

The Hausdorff property of $\mathscr{O}_{\mathcal{G}}(E)$ implies
that of $\mathscr{O}_{\mathcal{F}}(E)$ by Fact \ref{fact:Hausdorff_Refine}.
So, $\mathcal{G}$ separating points on $E$ implies $\mathcal{F}$
separating points on $E$ by Proposition \ref{prop:Fun_Sep_1} (c).

If $\mathcal{G}$ strongly separates points on topological space $E$,
then $\mathscr{O}(E)\subset\mathscr{O}_{\mathcal{G}}(E)\subset\mathscr{O}_{\mathcal{F}}(E)$
and so $\mathcal{F}$ strongly separates points on $E$.\end{proof}

\begin{cor}
\label{cor:C(E;R)_Cb(E;R)_SP}Let $E$ be a topological space. Then,
$C(E;\mathbf{R})$ separates points (resp. strongly separates points)
on $E$ if and only if $C_{b}(E;\mathbf{R})$ does.
\end{cor}

\begin{fact}
\label{fact:BExt_Same_Int}Let $E$ be a topological space, $\mu\in\mathfrak{M}^{+}(E^{d},\mathscr{B}(E)^{\otimes d})$,
$\nu_{1}\in\mathfrak{be}(\mu)$, $X\in M(\Omega,\mathscr{F},\mathbb{P};E^{d},\mathscr{B}(E)^{\otimes d})$
and $\nu_{2}\in\mathfrak{be}(\mathbb{P}\circ X^{-1})$. Then, $\int_{E^{d}}f(x)\mu(dx)=f^{*}(\nu_{1})$
and $\mathbb{E}[f\circ X]=f^{*}(\nu_{2})$ for all $f\in\mathfrak{ca}[\Pi^{d}(M_{b}(E;\mathbf{R}))]$.
\end{fact}
\begin{proof}
This result follows by Proposition \ref{prop:Pi^d_SP} (a) (with $\mathcal{D}=M_{b}(E;\mathbf{R})$)
and the fact that $\nu_{1}$ (resp. $\nu_{2}$) and $\mu$ (resp.
$\mathbb{P}\circ X^{-1}$) are the same measures on $(E^{d},\mathscr{B}(E)^{\otimes d})$.\end{proof}

\begin{lem}
\label{lem:WC_Expansion}Let $E$ be a topological space and $A\subset E$.
Then,
\begin{equation}
\mu_{n}\Longrightarrow\mu\mbox{ as }n\uparrow\infty\mbox{ in }\mathcal{M}^{+}\left(A,\mathscr{O}_{E}(A)\right)\label{eq:WC_on_Subspace}
\end{equation}
implies
\begin{equation}
\mu_{n}|^{E}\Longrightarrow\mu|^{E}\mbox{ as }n\uparrow\infty\mbox{ in }\mathcal{M}^{+}(E).\label{eq:WC_Expansion}
\end{equation}
The converse is true when $E$ is a Tychonoff space.
\end{lem}
\begin{proof}
It follows by (\ref{eq:WC_on_Subspace}) and $C_{b}(E;\mathbf{R})|_{A}\subset C_{b}(A,\mathscr{O}_{E}(A);\mathbf{R})$
that
\begin{equation}
\begin{aligned}\lim_{n\rightarrow\infty}\int_{E}f(x)\mu_{n}|^{E}(dx) & =\lim_{n\rightarrow\infty}\int_{A}f|_{A}(x)\mu_{n}(dx)=\int_{A}f|_{A}(x)\mu(dx)\\
 & =\int_{E}f(x)\mu|^{E}(dx),\;\forall f\in C_{b}(E;\mathbf{R}),
\end{aligned}
\label{eq:Func_Test_WC_on_A}
\end{equation}
proving (\ref{eq:WC_Expansion}). Conversely, if $E$ is Tychonoff,
then (\ref{eq:WC_Expansion}) implies
\begin{equation}
\limsup_{n\rightarrow\infty}\mu_{n}(F\cap A)=\limsup_{n\rightarrow\infty}\mu_{n}|^{E}(F)\leq\mu|^{E}(F)=\mu(F\cap A),\;\forall F\in\mathscr{C}(E)\label{eq:Portamental_Test_on_A}
\end{equation}
by Theorem \ref{thm:Portamenteau} (a, b) (with $\mu_{n}=\mu_{n}|^{E}$
and $\mu=\mu|^{E}$). $(A,\mathscr{O}_{E}(A))$ is a Hausdorff subspace
by Proposition \ref{prop:CR_Space} (b). Now, (\ref{eq:WC_on_Subspace})
follows by (\ref{eq:Portamental_Test_on_A}) and Theorem \ref{thm:Portamenteau}
(a, b).\end{proof}

\begin{cor}
\label{cor:Concent_Expan_RC}Let $E$ be a topological space and $\Gamma\subset\mathcal{M}^{+}(A,\mathscr{O}_{E}(A))$
with $A\subset E$. Then, $\{\mu|^{E}\}_{\mu\in\Gamma}$ is relatively
compact in $\mathcal{M}^{+}(E)$ whenever $\Gamma$ is.
\end{cor}

\begin{cor}
\label{cor:Identify_Seq_WC}Let $E$ be a topological space, $A\in\mathscr{B}(E)$
be a Hausdorff subspace, $\{\mu_{n}\}_{n\in\mathbf{N}}\subset\mathcal{M}^{+}(E)$
be sequentially tight in $A$ and $\{\mu_{n}(E)\}_{n\in\mathbf{N}}\subset[a,b]$
for some $0<a<b$. If $\mu$ is the unique weak limit point of $\{\mu_{n}\}_{n\in\mathbf{N}}$
in $\mathcal{M}^{+}(E)$, then (\ref{eq:WL}) holds.
\end{cor}
\begin{proof}
$\{\mu_{n}\}_{n\in\mathbf{N}}$ are all supported on $A$ with finite
exception by Fact \ref{fact:Seq_Tight_Support} (with $\mathscr{U}=\mathscr{B}(E)$
and $\Gamma=\{\mu_{n}\}_{n\in\mathbf{N}}$) and $\{\mu_{n}|_{A}\}_{n\in\mathbf{N}}$
is relatively compact in $\mathcal{M}^{+}(A,\mathscr{O}_{E}(A))$
by Lemma \ref{lem:Seq_Prokhorov} (with $E=(A,\mathscr{O}_{E}(A))$
and $\Gamma=\{\mu_{n}|_{A}\}_{n\in\mathbf{N}}$). Then, $\{\mu_{n}\}_{n\in\mathbf{N}}$
is relatively compact in $\mathcal{M}^{+}(E)$ by Fact \ref{fact:Meas_Concen_Expan}
(c) (with $\mathscr{U}=\mathscr{B}(E)$ and $\nu=\mu|_{A}$) and Corollary
\ref{cor:Concent_Expan_RC}. Now, the corollary follows by Fact \ref{fact:Uni_Seq_Lim_Conv}
(with $(E,x_{n},x)=(\mathcal{M}^{+}(E),\mu_{n},\mu)$).\end{proof}

\begin{lem}
\label{lem:CR_Dirac_Meas}Let $E$ be a Hausdorff space. Then, the
following statements are equivalent:

\renewcommand{\labelenumi}{(\alph{enumi})}
\begin{enumerate}
\item $E$ is a Tychonoff space.
\item $\delta_{x_{n}}\Rightarrow\delta_{x}$ as $n\uparrow\infty$ in $\mathcal{P}(E)$
implies $x_{n}\rightarrow x$ as $n\uparrow\infty$ in $E$.
\item Convergence determining implies determining point convergence on $E$.
\end{enumerate}
\end{lem}
\begin{proof}
((a) $\rightarrow$ (b)) $\delta_{x_{n}}\Rightarrow\delta_{x}$ as
$n\uparrow\infty$ in $\mathcal{P}(E)$ implies
\begin{equation}
\lim_{n\rightarrow\infty}f(x_{n})=\lim_{n\rightarrow\infty}f^{*}(\delta_{x_{n}})=f(\delta_{x})=f(x)\label{eq:CD_means_SSP}
\end{equation}
for all $f\in C_{b}(E;\mathbf{R})$. $C_{b}(E;\mathbf{R})$ determines
point convergence on $E$ by Proposition \ref{prop:CR} (a, c) and
Proposition \ref{prop:Fun_Sep_1} (b). Hence, (\ref{eq:CD_means_SSP})
implies $x_{n}\rightarrow x$ as $n\uparrow\infty$.

((b) $\rightarrow$ (c)) If $\mathcal{D}\subset M_{b}(E;\mathbf{R})$
satisfies $\bigotimes\mathcal{D}(x_{n})\rightarrow\bigotimes\mathcal{D}(x)$
as $n\uparrow\infty$, then (\ref{eq:CD_means_SSP}) holds for all
$f\in\mathcal{D}$. This implies $\delta_{x_{n}}\Rightarrow\delta_{x}$
as $n\uparrow\infty$ since $\mathcal{D}$ is convergence determining
on $E$. Now, we have $x_{n}\rightarrow x$ as $n\uparrow\infty$
by (b).

((c) $\rightarrow$ (a)) $C_{b}(E;\mathbf{R})$ determines point convergence
on $E$ by (c). It strongly separates points and separates points
on $E$ by Proposition \ref{prop:Fun_Sep_1} (a, b). Now, (a) follows
by Proposition \ref{prop:CR} (a, c).\end{proof}

\begin{lem}
\label{lem:Tight_Meas_Identical}Let $E$ be a topological space,
$\mathcal{D}\subset C_{b}(E;\mathbf{R})$ separate points on $E$
and $d\in\mathbf{N}$. Then:

\renewcommand{\labelenumi}{(\alph{enumi})}
\begin{enumerate}
\item If each of $\mu_{1},\mu_{2}\in\mathcal{M}^{+}(E)$ is tight and $\mathcal{D}$
is closed under multiplication, then $f^{*}(\mu_{1})=f^{*}(\mu_{2})$
for all $f\in\mathcal{D}\cup\{1\}$ implies $\mu_{1}=\mu_{2}$.
\item If each of $\mu_{1},\mu_{2}\in\mathfrak{M}^{+}(E^{d},\mathscr{B}(E)^{\otimes d})$
is $\mathbf{m}$-tight, then $f^{*}(\mu_{1})=f^{*}(\mu_{2})$ for
all \textup{$f\in\mathfrak{mc}[\Pi^{d}(\mathcal{D})]\cup\{1\}$ }implies
$\mu_{1}=\mu_{2}$.
\end{enumerate}
\end{lem}
\begin{proof}
(a) Let $a\circeq\mu_{1}(E)=\mu_{2}(E)>0$%
\footnote{We mentioned in \S \ref{sec:Convention} that any measure in this
work has positive total mass.%
} and $\nu_{i}\circeq\mu_{i}/a$ for each $i\in\{1,2\}$. Each of $\nu_{1}$
and $\nu_{2}$ is a tight member of $\mathcal{P}(E)$ and they satisfy
$f^{*}(\nu_{1})=f^{*}(\nu_{2})$ for all $f\in\mathcal{D}$. Then,
$\nu_{1}=\nu_{2}$ by \cite[Theorem 11 (d)]{BK10} and so $\mu_{1}=\mu_{2}$.

(b) $E$ is a Hausdorff space by Proposition \ref{prop:Fun_Sep_1}
(e) (with $A=E$). For each $j=1,2$, there exists an $\mathbf{m}$-tight
$\mu_{j}^{\prime}=\mathfrak{be}(\mu_{j})$ by Proposition \ref{prop:m-Tight_BExt}
(with $\mathbf{I}=\{1,...,d\}$, $S_{i}=E$, $A=E^{d}$ and $\Gamma=\{\mu_{j}\}$).
$\mathfrak{mc}[\Pi^{d}(\mathcal{D})]$ separates points on $E^{d}$
by Proposition \ref{prop:Pi^d_SP} (b). Now, (b) follows by (a) (with
$E=E^{d}$, $\mathcal{D}=\mathfrak{mc}[\Pi^{d}(\mathcal{D})]$ and
$\mu_{j}=\mu_{j}^{\prime}$).\end{proof}

\begin{fact}
\label{fact:Push_Forward_Tight_2}Let $E$ be a topological space
and $S$ a Hausdorff space. If $\Gamma\subset\mathcal{M}^{+}(E)$
is tight in $A\subset E$ and $f\in C(E;S)$, then $\{\mu\circ f^{-1}:\mu\in\Gamma\}$
is tight in $f(A)$.
\end{fact}
\begin{proof}
$f(K)\in\mathscr{K}(S)\subset\mathscr{B}(S)$ for all $K\in\mathscr{K}(E)$
by Proposition \ref{prop:Compact} (a, e). Now, the result follows
by Fact \ref{fact:Push_Forward_Tight_1} (with $\mathscr{U}=\mathscr{B}(E)$
and  $\mathscr{A}=\mathscr{B}(S)$).\end{proof}

\begin{lem}
\label{lem:Tightness_Prod}Let $\mathbf{I}$ be a countable index
set, $\{S_{i}\}_{i\in\mathbf{I}}$ be topological spaces, $(S,\mathscr{A})$
be as in (\ref{eq:(S,A)_Prod_Meas_Space}), $\Gamma\subset\mathfrak{M}^{+}(S,\mathscr{A})$,
$A_{i}\subset S_{i}$ for each $i\in\mathbf{I}$ and $A\circeq\prod_{i\in\mathbf{I}}A_{i}$.
Then:

\renewcommand{\labelenumi}{(\alph{enumi})}
\begin{enumerate}
\item If $\{\mu\circ\mathfrak{p}_{i}^{-1}\}_{\mu\in\Gamma}$ is tight (resp.
$\mathbf{m}$-tight) in $A_{i}$ for all $i\in\mathbf{I}$, then $\Gamma$
is tight (resp. $\mathbf{m}$-tight) in $A$. The converse is true
when $(A_{i},\mathscr{O}_{S_{i}}(A_{i}))$ is a Hausdorff subspace
of $S_{i}$ and $A_{i}\in\mathscr{B}(S_{i})$.
\item If $\{\mu\circ\mathfrak{p}_{i}^{-1}\}_{\mu\in\Gamma}$ is sequentially
tight (resp. $\mathbf{m}$-tight) in $A_{i}$ for all $i\in\mathbf{I}$,
then $\Gamma$ is sequentially tight (resp. $\mathbf{m}$-tight) in
$A$. The converse is true when $(A_{i},\mathscr{O}_{S_{i}}(A_{i}))$
is a Hausdorff subspace of $S_{i}$ and $A_{i}\in\mathscr{B}(S_{i})$.
\end{enumerate}
\end{lem}
\begin{proof}
(a) Without loss of generality, we suppose $\mathbf{I}=\mathbf{N}$.
Each $A_{i}$ is equipped with the subspace topology $\mathscr{O}_{S_{i}}(A_{i})$
throughout the proof. If $\{\mu\circ\mathfrak{p}_{i}^{-1}\}_{\mu\in\Gamma}$
is tight in $A_{i}$ for all $i\in\mathbf{I}$, then there exist
\begin{equation}
A_{i}\supset K_{p,i}\in\mathscr{K}(S_{i})\cap\mathscr{B}(S_{i}),\;\forall i,p\in\mathbf{N}\label{eq:Check_Joint_Tight_1}
\end{equation}
such that
\begin{equation}
\sup_{\mu\in\Gamma}\mu\circ\mathfrak{p}_{i}^{-1}\left(S_{i}\backslash\mathfrak{p}_{i}(K_{p})\right)\leq2^{-p-i},\;\forall i,p\in\mathbf{N}.\label{eq:Check_Joint_Tight_2}
\end{equation}
It follows that
\begin{equation}
A\supset\prod_{i\in\mathbf{N}}K_{p,i}\in\mathscr{K}(S)\cap\mathscr{A},\;\forall p\in\mathbf{N}\label{eq:Check_Joint_Tight_3}
\end{equation}
by Proposition \ref{prop:Compact} (b), Fact \ref{fact:Prod_Map_2}
(a) and the fact $\prod_{i\in\mathbf{I}}K_{p,i}=\bigcap_{i\in\mathbf{I}}\mathfrak{p}_{i}^{-1}(K_{p,i})$.
Now, we conclude the tightness of $\Gamma$ in $A$ by observing that
\begin{equation}
\sup_{\mu\in\Gamma}\mu\left(S\backslash\prod_{i\in\mathbf{N}}K_{p,i}\right)\leq\sum_{i=1}^{\infty}\sup_{\mu\in\Gamma}\mu\left(S_{i}\backslash K_{p,i}\right)\leq2^{-p},\;\forall p\in\mathbf{N}.\label{eq:Check_Joint_Tight_4}
\end{equation}
If $\{\mu\circ\mathfrak{p}_{i}^{-1}\}_{\mu\in\Gamma}$ is $\mathbf{m}$-tight
in $A_{i}$ for all $i\in\mathbf{I}$, then we retake each $K_{p,i}$
above from $\mathscr{K}^{\mathbf{m}}(S_{i})\cap\mathscr{B}(S_{i})$,
find $\prod_{i\in\mathbf{I}}K_{p,i}\in\mathscr{K}^{\mathbf{m}}(S)$
by Lemma \ref{lem:MC_Prod} (a) (with $A_{i}=K_{p,i}$) and verify
the $\mathbf{m}$-tightness of $\Gamma$ by a similar argument.

Next, we suppose $A_{i}\in\mathscr{B}(S_{i})$ is a Hausdorff subspace
for all $i\in\mathbf{I}$ and justify the converse statement. We have
that
\begin{equation}
\mathfrak{p}_{i}(K)\in\mathscr{K}(A_{i})\subset\mathscr{B}(A_{i})\subset\mathscr{B}(S_{i}),\;\forall K\in\mathscr{K}(S),i\in\mathbf{I}\label{eq:Check_p_i(K)_Borel_1}
\end{equation}
and
\begin{equation}
\mathfrak{p}_{i}(K)\in\mathscr{K}^{\mathbf{m}}(A_{i})\subset\mathscr{B}(A_{i})\subset\mathscr{B}(S_{i}),\;\forall K\in\mathscr{K}^{\mathbf{m}}(S),i\in\mathbf{I}\label{eq:Check_p_i(K)_Borel_2}
\end{equation}
by Corollary \ref{cor:Compact_Prod} (b) (with $A=K$ and $S_{i}=A_{i}$),
Lemma \ref{lem:MC_Prod} (b) (with $A=K$ and $S_{i}=A_{i}$) and
the fact $A_{i}\in\mathscr{B}(S_{i})$. If $\Gamma$ is tight (resp.
$\mathbf{m}$-tight) in $A$, then for each $i\in\mathbf{I}$, the
tightness (resp. $\mathbf{m}$-tightness) of $\{\mu\circ\mathfrak{p}_{i}^{-1}\}_{\mu\in\Gamma}$
in $A_{i}$ follows by (\ref{eq:Check_p_i(K)_Borel_1}), (\ref{eq:Check_p_i(K)_Borel_2}),
the fact $A_{i}=\mathfrak{p}_{i}(A)$, Fact \ref{fact:Prod_Map_1}
(a) and Fact \ref{fact:Push_Forward_Tight_1} (with $(E,\mathscr{U})=(S,\mathscr{A})$,
$(S,\mathscr{A})=(S_{i},\mathscr{B}(S_{i}))$ and $f=\mathfrak{p}_{i}$).

(b) follows immediately by (a) and a triangular array argument.\end{proof}

\begin{lem}
\label{lem:Sko_FDD_BExt}Let $E$ be a Tychonoff space and $(D(\mathbf{R}^{+};E),\mathscr{S},\nu)$
be the completion of $(D(\mathbf{R}^{+};E),\sigma(\mathscr{J}(E)),\mu)$%
\footnote{$\mathscr{B}(D(\mathbf{R}^{+};E))\circeq\sigma[\mathscr{J}(E)]$,
so the measure space notation ``$(D(\mathbf{R}^{+};E),\sigma(\mathscr{J}(E)),\mu)$''
implies $\mu\in\mathcal{M}^{+}(D(\mathbf{R}^{+};E))$.%
}. If $M(E;\mathbf{R})$ has a countable subset separating points on
$E$, especially if $E$ is baseable, then $\nu\circ\mathfrak{p}_{\mathbf{T}_{0}}^{-1}\in\mathcal{M}^{+}(E^{\mathbf{T}_{0}})$%
\footnote{Herein, we show the domain of $\nu\circ\mathfrak{p}_{\mathbf{T}_{0}}^{-1}$
contains $\mathscr{B}(E^{\mathbf{T}_{0}})$ so $\nu\circ\mathfrak{p}_{\mathbf{T}_{0}}^{-1}$
can be viewed as a member of $\mathcal{M}^{+}(E^{\mathbf{T}_{0}})$.%
} is a Borel extension of $\mu\circ\mathfrak{p}_{\mathbf{T}_{0}}^{-1}\in\mathfrak{M}^{+}(E^{\mathbf{T}_{0}},\mathscr{B}(E)^{\otimes\mathbf{T}_{0}})$
for all non-empty $\mathbf{T}_{0}\in\mathscr{P}_{0}(\mathbf{R}^{+}\backslash J(\mu))$%
\footnote{The $J(\mu)$ herein is well-defined by Lemma \ref{lem:Sko_Proj}
(b) and Fact \ref{fact:J(Mu)_Well_Defined}.%
}.
\end{lem}
\begin{proof}
$\mu\circ\mathfrak{p}_{\mathbf{T}_{0}}^{-1}$ is a member of $\mathfrak{M}^{+}(E^{\mathbf{T}_{0}},\mathscr{B}(E)^{\otimes\mathbf{T}_{0}})$
by Lemma \ref{lem:Sko_Proj} (a). $\mathfrak{p}_{\mathbf{T}_{0}}\in M[D(\mathbf{R}^{+};E),\mathscr{S};E^{\mathbf{T}_{0}},\mathscr{B}(E^{\mathbf{T}_{0}})]$
by Lemma \ref{lem:Sko_Proj} (c), the definition of $J(\mu)$ and
Fact \ref{fact:AS_Cont} (with $E=D(\mathbf{R}^{+};E)$, $\mathscr{U}=\mathscr{S}$,
$S=E^{\mathbf{T}_{0}}$ and $f=\mathfrak{p}_{\mathbf{T}_{0}}$). Hence,
$\nu\circ\mathfrak{p}_{\mathbf{T}_{0}}^{-1}\in\mathfrak{be}(\mu\circ\mathfrak{p}_{\mathbf{T}_{0}}^{-1})$
as $\nu$ is an extension%
\footnote{Extension of measure was specified in \S \ref{sub:Meas}.%
} of $\mu$ to $\mathscr{S}$.\end{proof}

\begin{lem}
\label{lem:Sko_Prod_Tight}Let $E$ be a Tychonoff space, $\mathcal{D}\subset C(E;\mathbf{R})$
be countable and $\{\mu_{i}\}_{i\in\mathbf{I}}\subset\mathcal{M}^{+}(D(\mathbf{R}^{+};E))$.
Then:

\renewcommand{\labelenumi}{(\alph{enumi})}
\begin{enumerate}
\item If $\{\mu_{i}\circ\varpi(f)^{-1}\}_{i\in\mathbf{I}}$ is tight in
$D(\mathbf{R}^{+};\mathbf{R})$ for all $f\in\mathcal{D}$, then $\{\mu_{i}\circ\varpi(\mathcal{D})^{-1}\}_{i\in\mathbf{I}}$
is tight in $D(\mathbf{R}^{+};\mathbf{R})^{\mathcal{D}}$.
\item If $\mathcal{D}$ strongly separates points on $E$ and $\{\mu_{i}\circ\varpi(f)^{-1}\}_{i\in\mathbf{I}}$
is tight in $D(\mathbf{R}^{+};\mathbf{R})$ for all $f\in\mathfrak{ae}(\mathcal{D})$,
then $\{\mu_{i}\circ\varpi[\bigotimes\mathfrak{ae}(\mathcal{D})]^{-1}\}_{i\in\mathbf{I}}$
is tight in $D(\mathbf{R}^{+};\mathbf{R}^{\mathfrak{ae}(\mathcal{D})})$.
\end{enumerate}
\end{lem}
\begin{note}
\label{note:D(R^Inf)_Polish}When $\mathcal{D}$ is countable, $\mathfrak{ae}(\mathcal{D})$
is also by Fact \ref{fact:ac_mc_Countable}. Then, $\mathbf{R}^{\mathfrak{ae}(\mathcal{D})}$,
$D(\mathbf{R}^{+};\mathbf{R})$, $D(\mathbf{R}^{+};\mathbf{R}^{\mathfrak{ae}(\mathcal{D})})$
and $D(\mathbf{R}^{+};\mathbf{R})^{\mathfrak{ae}(\mathcal{D})}$ are
all Polish spaces by Proposition \ref{prop:Var_Polish} (f) and Proposition
\ref{prop:Sko_Basic_2} (d).
\end{note}
\begin{proof}
[Proof of Lemma \ref{lem:Sko_Prod_Tight}](a) follows by the fact
\begin{equation}
\left(\mu_{i}\circ\varpi(\mathcal{D})^{-1}\right)\circ\mathfrak{p}_{f}^{-1}=\mu_{i}\circ\varpi(f)^{-1}\in\mathcal{M}^{+}\left(D(\mathbf{R}^{+};\mathbf{R})\right),\;\forall f\in\mathcal{D}\label{eq:Path_Mapping_Push_Forward_Meas}
\end{equation}
and Lemma \ref{lem:Tightness_Prod} (a) (with $\mathbf{I}=\mathcal{D}$,
$S_{i}=A_{i}=D(\mathbf{R}^{+};\mathbf{R})$ and $\Gamma=\{\mu_{i}\circ\varpi(\mathcal{D})^{-1}\}_{i\in\mathbf{I}}$).
Here, $\mathfrak{p}_{f}$ denotes the one-dimensional projection on
$D(\mathbf{R}^{+};\mathbf{R})^{\mathcal{D}}$ for $f\in\mathcal{D}$.

(b) Letting $\mathcal{J}=\{\mathfrak{p}_{f}\}_{f\in\mathcal{D}}$
be the one-dimensional projections on $\mathbf{R}^{\mathcal{D}}$,
we have by Corollary \ref{cor:D(R)^Inf_Imb} (with $\mathbf{I}=\mathcal{D}$)
that%
\footnote{Remark \ref{rem:Path_Mapping} noted that $\phi_{1}$ and $\phi_{2}$
are different.%
}
\begin{equation}
\phi_{1}\circeq\varpi\left[\bigotimes\mathfrak{ae}(\mathcal{J})\right]\in C\left[D(\mathbf{R}^{+};\mathbf{R}^{\mathcal{D}});D\left(\mathbf{R}^{+};\mathbf{R}^{\mathfrak{ae}(\mathcal{D})}\right)\right]\label{eq:Check_D(R)^Inf_D(R^Inf)_Cont_1}
\end{equation}
and
\begin{equation}
\phi_{2}\circeq\varpi[\mathfrak{ae}(\mathcal{J})]\in\mathbf{imb}\left[D(\mathbf{R}^{+};\mathbf{R}^{\mathcal{D}});D(\mathbf{R}^{+};\mathbf{R})^{\mathfrak{ae}(\mathcal{D})}\right],\label{eq:Check_D(R)^Inf_D(R^Inf)_Cont_2}
\end{equation}
which implies
\begin{equation}
\phi_{1}\circ\phi_{2}^{-1}\in C\left[D(\mathbf{R}^{+};\mathbf{R})^{\mathfrak{ae}(\mathcal{D})};D\left(\mathbf{R}^{+};\mathbf{R}^{\mathfrak{ae}(\mathcal{D})}\right)\right].\label{eq:D(R)^Inf_D(R^Inf)_Cont}
\end{equation}
$\{\mu_{i}\circ\varpi(\mathfrak{ae}(\mathcal{D}))^{-1}\}_{i\in\mathbf{I}}$
is tight in $D(\mathbf{R}^{+};\mathbf{R})^{\mathfrak{ae}(\mathcal{D})}$
by (a) (with $\mathcal{D}=\mathfrak{ae}(\mathcal{D})$). Observing
that
\begin{equation}
\mu_{i}\circ\phi_{1}^{-1}=\left[\mu_{i}\circ\varpi\left(\mathfrak{ae}(\mathcal{D})\right)^{-1}\right]\circ\left[\phi_{1}\circ\phi_{2}^{-1}\right]^{-1},\;\forall i\in\mathbf{I},\label{eq:Path_Mapping_Prod_Dist_Proj}
\end{equation}
we have the desired tightness of $\{\mu_{i}\circ\phi_{1}^{-1}\}_{i\in\mathbf{I}}$
by (\ref{eq:D(R)^Inf_D(R^Inf)_Cont}), (\ref{eq:Path_Mapping_Prod_Dist_Proj})
and Fact \ref{fact:Push_Forward_Tight_2} (with $E=D(\mathbf{R}^{+};\mathbf{R})^{\mathfrak{ae}(\mathcal{D})}$,
$S=D(\mathbf{R}^{+};\mathbf{R}^{\mathfrak{ae}(\mathcal{D})})$ and
$f=\phi_{1}\circ\phi_{2}^{-1}$).\end{proof}

\begin{lem}
\label{lem:Sko_Prod_Eq}Let $E$ be a Tychonoff space, $\mathcal{D}\subset C(E;\mathbf{R})$
be countable, $\Psi\circeq\varpi[\mathfrak{ae}(\mathcal{D})]$, $V\subset D(\mathbf{R}^{+};E)$
and $\{A_{p}\}_{p\in\mathbf{N}}\subset\mathscr{B}(E)$. If $A_{p}\subset A_{p+1}$,
$\mathcal{D}$ strongly separates points on $A_{p}$ and $x|_{[0,p)}\in A_{p}^{[0,p)}$
for all $x\in V$ and $p\in\mathbf{N}$, then
\begin{equation}
\Psi|_{V}\in\mathbf{imb}\left(V,\mathscr{O}_{D(\mathbf{R}^{+};E)}(V);D(\mathbf{R}^{+}\mathbf{R})^{\mathfrak{ae}(\mathcal{D})}\right)\label{eq:Psi_Homeo_on_V}
\end{equation}
and
\begin{equation}
\mathscr{B}_{D(\mathbf{R}^{+};E)}(V)=\left.\mathscr{B}(E)^{\otimes\mathbf{R}^{+}}\right|_{V}.\label{eq:V_Sko_Borel_Prod_Equal}
\end{equation}

\end{lem}
\begin{proof}
\textit{Step 1: Show $\Psi|_{V}$ is injective}. $E_{0}\circeq\bigcup_{p\in\mathbf{N}}A_{p}\in\mathscr{B}(E)$
satisfies
\begin{equation}
V\subset D\left(\mathbf{R}^{+};E_{0},\mathscr{O}_{E}(E_{0})\right).\label{eq:Path_in_E0}
\end{equation}
$\mathfrak{ae}(\mathcal{D})$ separates points on $E_{0}$ by the
Hausdorff property of $E$, Proposition \ref{prop:Fun_Sep_1} (a)
and Fact \ref{fact:SP_Nested_Union} (with $n=p$ and $\mathcal{D}_{n}=\mathfrak{ae}(\mathcal{D})$).
So, $\Psi|_{V}$ is injective by (\ref{eq:Path_in_E0}) and Fact \ref{fact:Path_Mapping_Injective}
(with $A=E_{0}$ and $\mathcal{D}=\mathfrak{ae}(\mathcal{D})$).

\textit{Step 2: Show the continuity of $\Psi|_{V}$}. We have by Fact
\ref{fact:ac_mc_Countable} and Proposition \ref{prop:Fun_Sep_1}
(d) (with $A=E_{0}$) that $(E_{0},\rho_{\mathcal{D}})$ is a separable
metric space, $\mathfrak{ae}(\mathcal{D})$ is a countable subset
of $C(E_{0},\rho_{\mathcal{D}};\mathbf{R})$ and $\mathfrak{ae}(\mathcal{D})$
strongly separates points on $(E_{0},\rho_{\mathcal{D}})$. It follows
that
\begin{equation}
\Psi|_{D\left(\mathbf{R}^{+};E_{0},\rho_{\mathcal{D}}\right)}\in\mathbf{imb}\left(D\left(\mathbf{R}^{+};E_{0},\rho_{\mathcal{D}}\right);D(\mathbf{R}^{+};\mathbf{R})^{\mathfrak{ae}(\mathcal{D})}\right)\label{eq:Psi_Cont_rho_D}
\end{equation}
by Proposition \ref{prop:Sko_Basic_1} (a) (with $\mathcal{D}=\mathfrak{ae}(\mathcal{D})$).
$(E_{0},\rho_{\mathcal{D}})=(E_{0},\mathscr{O}_{\mathcal{D}}(E_{0}))$
is a topological coarsening of $(E_{0},\mathscr{O}_{E}(E_{0}))$ since
$\mathcal{D}\subset C(E;\mathbf{R}^{+})$. So,
\begin{equation}
\Psi|_{D\left(\mathbf{R}^{+};E_{0},\mathscr{O}_{E}(E_{0})\right)}\in C\left(D\left(\mathbf{R}^{+};E_{0},\mathscr{O}_{E}(E_{0})\right);D(\mathbf{R}^{+};\mathbf{R})^{\mathfrak{ae}(\mathcal{D})}\right)\label{eq:Psi_Cont}
\end{equation}
by (\ref{eq:Psi_Cont_rho_D}) and Proposition \ref{prop:Sko_Basic_1}
(e) (with $E=(E_{0},\rho_{\mathcal{D}})$ and $S=(E_{0},\mathscr{O}_{E}(E_{0}))$).
Hence, we have by (\ref{eq:Path_in_E0}) and (\ref{eq:Psi_Cont})
that
\begin{equation}
\Psi|_{V}\in C\left(V,\mathscr{O}_{D(\mathbf{R}^{+};E)}(V);D(\mathbf{R}^{+};\mathbf{R})^{\mathfrak{ae}(\mathcal{D})}\right).\label{eq:Psi_V_Cont}
\end{equation}

\textit{Step 3: Show the continuity of $(\Psi|_{V})^{-1}$}. $D(\mathbf{R}^{+};\mathbf{R})^{\mathfrak{ae}(\mathcal{D})}$
is a Polish space, so its subspace $\Psi(V)$ is metrizable by Proposition
\ref{prop:Metrizable} (b). According to \cite[Theorem 21.3]{M00},
showing the continuity of $(\Psi|_{V})^{-1}$ is equivalent to showing
that (\ref{eq:Psi(yk)_Conv_Psi(y0)}) implies (\ref{eq:yk_Conv_y0})
for all $\{y_{k}\}_{k\in\mathbf{N}_{0}}\subset V$.

We suppose (\ref{eq:Psi(yk)_Conv_Psi(y0)}) holds, fix $u\in\mathbf{R}^{+}\backslash J(y_{0})$,
define $p_{u}\circeq\min\{p\in\mathbf{N}:p>u+1\}$ and let $\{y_{k}^{u}\}_{k\in\mathbf{N}_{0}}$
be as in (\ref{eq:yk_u}). Observing that
\begin{equation}
\left\{ y_{k}|_{[0,u+1]}\right\} _{k\in\mathbf{N}}\cup\left\{ y|_{[0,u+1]}\right\} \subset D\left([0,u+1];A_{p_{u}},\mathscr{O}_{E}(A_{p_{u}})\right),\label{eq:Check_Sko_Conv-2}
\end{equation}
we have
\begin{equation}
y_{k}^{u}\longrightarrow y_{0}^{u}\mbox{ as }k\uparrow\infty\mbox{ in }D\left([0,u+1];A_{p_{u}},\mathscr{O}_{E}(A_{p_{u}})\right)\label{eq:yku_Conv_y0u_D(A_pu)}
\end{equation}
by Lemma \ref{lem:D_Conv_D[0,u]_Conv} (with $E=(A_{p_{u}},\mathscr{O}_{E}(A_{p_{u}}))$
and $\mathcal{D}=\mathcal{D}|_{A_{p_{u}}}$). This implies
\begin{equation}
y_{k}^{u}\longrightarrow y_{0}^{u}\mbox{ as }k\uparrow\infty\mbox{ in }D\left([0,u+1];E_{0},\mathscr{O}_{E}(E_{0})\right)\label{eq:yku_Conv_y0u_D(E0)}
\end{equation}
by Corollary \ref{cor:Sko_Subspace} (with $E=(E_{0},\mathscr{O}_{E}(E_{0}))$
and $A=A_{p_{u}}$). The countable collection $\mathcal{D}\subset C(E;\mathbf{R})$
separates points on $E_{0}$, $(E_{0},\mathscr{O}_{E}(E_{0}))$ is
a baseable space. Hence, it follows by (\ref{eq:yku_Conv_y0u_D(E0)})
and Lemma \ref{lem:D[0,u]_Conv_D_Conv} (with $E=(E_{0},\mathscr{O}_{E}(E_{0}))$)
that
\begin{equation}
y_{k}\longrightarrow y_{0}\mbox{ as }k\uparrow\infty\mbox{ in }D\left(\mathbf{R}^{+};E_{0},\mathscr{O}_{E}(E_{0})\right).\label{eq:yk_Conv_y0_D(E0)}
\end{equation}
Now, the desired (\ref{eq:yk_Conv_y0}) follows by Corollary \ref{cor:Sko_Subspace}
(with $A=E_{0}$).

\textit{Step 4: Show (\ref{eq:V_Sko_Borel_Prod_Equal})}. The three
steps above established (\ref{eq:Psi_Homeo_on_V}).
\begin{equation}
\varpi(f)\in M\left(E^{\mathbf{R}^{+}},\mathscr{B}(E)^{\otimes\mathbf{R}^{+}};\mathbf{R}^{\mathbf{R}^{+}},\mathscr{B}(\mathbf{R})^{\otimes\mathbf{R}^{+}}\right),\;\forall f\in\mathfrak{ae}(\mathcal{D})\label{eq:psi_Prod_Measurable}
\end{equation}
by $\mathcal{D}\subset C(E;\mathbf{R})$ and Fact \ref{fact:Path_Mapping}
(b) (with $\mathbf{I}=\mathbf{R}^{+}$).
\begin{equation}
\left[\left.\mathscr{B}(\mathbf{R})^{\otimes\mathbf{R}^{+}}\right|_{D(\mathbf{R}^{+};\mathbf{R})}\right]^{\otimes\mathfrak{ae}(\mathcal{D})}=\sigma\left[\mathscr{J}(\mathbf{R})\right]^{\otimes\mathfrak{ae}(\mathcal{D})}=\sigma\left[\mathscr{J}(\mathbf{R})^{\mathfrak{ae}(\mathcal{D})}\right]\label{eq:Psi(V)_Sko_Borel_Prod_Equal}
\end{equation}
by Proposition \ref{prop:Sko_Basic_2} (b) (with $E=\mathbf{R}$),
the fact that $D(\mathbf{R}^{+};\mathbf{R})$ is a Polish space and
Proposition \ref{prop:Prod_Space} (d) (with $S_{i}=D(\mathbf{R}^{+};\mathbf{R})$).
\begin{equation}
\Psi|_{V}\in M\left(V,\left.\mathscr{B}(E)^{\otimes\mathbf{R}^{+}}\right|_{V};\Psi(V),\left.\sigma\left[\mathscr{J}(\mathbf{R})^{\mathfrak{ae}(\mathcal{D})}\right]\right|_{\Psi(V)}\right)\label{eq:Check_Sko_Borel_1}
\end{equation}
by (\ref{eq:psi_Prod_Measurable}), Fact \ref{fact:Prod_Map_1} (b)
and (\ref{eq:Psi(V)_Sko_Borel_Prod_Equal}). 
\begin{equation}
\begin{aligned}\mathscr{B}_{D(\mathbf{R}^{+};E)}(V) & =\sigma\left(\left\{ \left(\Psi|_{V}\right)^{-1}(O):O\in\mathscr{J}(\mathbf{R})^{\mathfrak{ae}(\mathcal{D})}\right\} \right)\\
 & =\left\{ \left(\Psi|_{V}\right)^{-1}(B):B\in\sigma\left(\mathscr{J}(\mathbf{R})^{\mathfrak{ae}(\mathcal{D})}\right)\right\} \subset\left.\mathscr{B}(E)^{\otimes\mathbf{R}^{+}}\right|_{V}
\end{aligned}
\label{eq:Check_Sko_Borel_2}
\end{equation}
by (\ref{eq:Psi_Homeo_on_V}) and (\ref{eq:Check_Sko_Borel_1}). Now,
(\ref{eq:V_Sko_Borel_Prod_Equal}) follows by Lemma \ref{lem:Sko_Proj}
(b).\end{proof}

\begin{lem}
\label{lem:CCC_Closed}Let $E$ be a Tychonoff space, $\mathcal{D}\subset C(E;\mathbf{R})$
be countable, $\Psi\circeq\varpi[\mathfrak{ae}(\mathcal{D})]$, $\varphi\circeq\bigotimes\mathfrak{ae}(\mathcal{D})$,
$\{A_{p}\}_{p\in\mathbf{N}}\subset\mathscr{B}(E)$ and
\begin{equation}
V\circeq\bigcap_{p\in\mathbf{N}}\left\{ x\in D(\mathbf{R}^{+};E):x|_{[0,p)}\in A_{p}^{[0,p)}\right\} .\label{eq:V_CCC}
\end{equation}
If $A_{p}\subset A_{p+1}$, $\mathcal{D}$ strongly separates points
on $A_{p}$ and $\varphi(A_{p})\in\mathscr{C}(\mathbf{R}^{\mathfrak{ae}(\mathcal{D})})$
for all $p\in\mathbf{N}$, then $\Psi(V)\in\mathscr{C}(D(\mathbf{R}^{+};\mathbf{R})^{\mathfrak{ae}(\mathcal{D})})$.
\end{lem}
\begin{proof}
$D(\mathbf{R}^{+};\mathbf{R}^{\mathfrak{ae}(\mathcal{D})})$ is a
Polish space, so $\Psi(V)$ as a subspace is metrizable by Proposition
\ref{prop:Metrizable} (b). Hence, showing the closeness of $\Psi(V)$
is reduced by Fact \ref{fact:First_Countable} (with $E=D(\mathbf{R}^{+};\mathbf{R})^{\mathfrak{ae}(\mathcal{D})}$
and $A=\Psi(V)$) to showing that
\begin{equation}
\Psi(y_{k})\longrightarrow z\mbox{ as }k\uparrow\infty\mbox{ in }D(\mathbf{R}^{+};\mathbf{R})^{\mathfrak{ae}(\mathcal{D})}\label{eq:Check_Psi(V)_Closed}
\end{equation}
imply $z\in\Psi(V)$ for any $\{y_{k}\}_{k\in\mathbf{N}}\subset V$.

Let $\{\mathfrak{p}_{f}\}_{f\in\mathcal{D}}$ be the one-dimensional
projections on $D(\mathbf{R}^{+};\mathbf{R})^{\mathfrak{ae}(\mathcal{D})}$.
\begin{equation}
z^{\prime}(t)\circeq\bigotimes_{f\in\mathfrak{ae}(\mathcal{D})}\mathfrak{p}_{f}(z)(t),\;\forall t\in\mathbf{R}^{+}\label{eq:z^prime}
\end{equation}
defines a member of $D(\mathbf{R}^{+};\mathbf{R}^{\mathfrak{ae}(\mathcal{D})})$
by Fact \ref{fact:Cadlag_Path} (c).
\begin{equation}
\mathbf{T}\circeq\bigcap_{f\in\mathfrak{ae}(\mathcal{D})}\mathbf{R}^{+}\backslash J(\mathfrak{p}_{f}(z))\label{eq:Common_Cont_Times}
\end{equation}
is cocountable by Proposition \ref{prop:Sko_Baseable} (b) (with $E=\mathbf{R}$
and $x=\mathfrak{p}_{f}(z)$).
\begin{equation}
\varphi\circ y_{k}(t)\longrightarrow z^{\prime}(t)\mbox{ as }k\uparrow\infty\mbox{ in }\mathbf{R}^{\mathfrak{ae}(\mathcal{D})},\;\forall t\in\mathbf{T}\label{eq:phi(yk)_Conv_Cont_Times}
\end{equation}
by (\ref{eq:Check_Psi(V)_Closed}), Fact \ref{fact:Seq_Prod_Conv}
and Lemma \ref{lem:Sko_Proj} (c). It then follows that
\begin{equation}
z^{\prime}(t)\in\varphi(A_{p}),\;\forall t\in[0,p),p\in\mathbf{N}\label{eq:Limit_in_phi(Ap)}
\end{equation}
by (\ref{eq:phi(yk)_Conv_Cont_Times}), the closedness of each $\varphi(A_{p})$
in $\mathbf{R}^{\mathfrak{ae}(\mathcal{D})}$, the denseness of $\mathbf{T}$
in $\mathbf{R}^{+}$ and the right-continuity of $z^{\prime}$. 
\begin{equation}
\varphi|_{A_{p}}\in\mathbf{imb}\left(A_{p},\mathscr{O}_{E}(A_{p});\mathbf{R}^{\mathfrak{ae}(\mathcal{D})}\right),\;\forall p\in\mathbf{N}\label{eq:phi_Imb_Ap}
\end{equation}
by Lemma \ref{lem:Compactification} (a, c). So, $(\varphi|_{A_{p}})^{-1}\circ z^{\prime}|_{[0,p)}$
is a c$\grave{\mbox{a}}$dl$\grave{\mbox{a}}$g mapping from $[0,p)$
to $(A_{p},\mathscr{O}_{E}(A_{p}))$ for all $p\in\mathbf{N}$ and,
hence,
\begin{equation}
y(t)\circeq(\varphi|_{A_{p}})^{-1}(z^{\prime}(t)),\;\forall t\in[0,p),p\in\mathbf{N}\label{eq:Construct_Psi(yk)_Lim}
\end{equation}
well defines a member of $V$. Now, one observes $\Psi(y)=z$ from
(\ref{eq:z^prime}), (\ref{eq:Construct_Psi(yk)_Lim}) and the definitions
of $\Psi$ and $\varphi$.\end{proof}

The next lemma is adapted from \cite[Lemma 24]{K15}, restated befittingly.
\begin{lem}
\label{lem:Pull_Back_Tight}Let $E$ and $S$ be topological spaces,
$\{A_{p}\}_{p\in\mathbf{N}}\subset\mathscr{B}(E)$ and $f\in S^{E}$
satisfy $\{f(A_{p})\}_{p\in\mathbf{N}}\subset\mathscr{C}(S)$ and
\begin{equation}
f|_{A_{p}}\in\mathbf{hom}\left[A_{p},\mathscr{O}_{E}(A_{p});f(A_{p}),\mathscr{O}_{S}\left(f(A_{p})\right)\right],\;\forall p\in\mathbf{N},\label{eq:Pull_Back_Tight_Homeo}
\end{equation}
$E_{0}\circeq\bigcup_{p\in\mathbf{N}}A_{p}$ satisfy $f\in M(E_{0},\mathscr{O}_{E}(E_{0});S)$
and $\{\mu_{i}\}_{i\in\mathbf{I}}\subset\mathcal{P}(E)$ satisfy
\begin{equation}
\inf_{i\in\mathbf{I}}\mu_{i}(A_{p})\geq1-2^{-p},\;\forall p\in\mathbf{N}.\label{eq:Pull_Back_Tight_Ap}
\end{equation}
Then, tightness of $\{\mu_{i}\circ f^{-1}\}_{i\in\mathbf{I}}$ implies
that of $\{\mu_{i}\}_{i\in\mathbf{I}}$ in $E_{0}$.
\end{lem}

\begin{lem}
\label{lem:Sko_Meas_FDD_Int_Test}Let $E$ be a Tychonoff space, $\mathbf{T}\subset\mathbf{R}^{+}$
be dense, $d\in\mathbf{N}$ and $f\in C_{b}(E^{d};\mathbf{R})$. Then:

\renewcommand{\labelenumi}{(\alph{enumi})}
\begin{enumerate}
\item For each $\mathbf{T}_{0}=\{t_{1},...,t_{d}\}\in\mathscr{P}_{0}(\mathbf{R}^{+})$,
there exists a $\mathbf{T}_{p}=\{t_{p,1},...,t_{p,d}\}\in\mathscr{P}_{0}(\mathbf{T})$
for each $p\in\mathbf{N}$ such that%
\footnote{Recall that Corollary \ref{cor:Sko_Meas_FDD} verified $\mu\circ\mathfrak{p}_{\mathbf{T}_{0}}^{-1}\in\mathfrak{M}^{+}(E^{\mathbf{T}_{0}},\mathscr{B}(E)^{\otimes\mathbf{T}_{0}})$
for all $\mathbf{T}_{0}\in\mathscr{P}_{0}(\mathbf{R}^{+})$ and $\mu\in\mathfrak{M}^{+}(D(\mathbf{R}^{+};E),\mathscr{B}(E)^{\otimes\mathbf{R}^{+}}|_{D(\mathbf{R}^{+};E)})$.%
}
\begin{equation}
\lim_{p\rightarrow\infty}\int_{E^{\mathbf{T}_{0}}}f(x)\mu\circ\mathfrak{p}_{\mathbf{T}_{p}}^{-1}(dx)=\int_{E^{\mathbf{T}_{0}}}f(x)\mu\circ\mathfrak{p}_{\mathbf{T}_{0}}^{-1}(dx).\label{eq:Sko_Meas_FDD_Int_Time_Approx}
\end{equation}
for all $\mu\in\mathfrak{M}^{+}(D(\mathbf{R}^{+};E),\mathscr{B}(E)^{\otimes\mathbf{R}^{+}}|_{D(\mathbf{R}^{+};E)})$.
\item If $\gamma^{1},\gamma^{2}\in\mathfrak{M}^{+}(D(\mathbf{R}^{+};E),\mathscr{B}(E)^{\otimes\mathbf{R}^{+}}|_{D(\mathbf{R}^{+};E)})$
satisfy (\ref{eq:Sko_WLP_FDD_Uni_Int_Test}) for all $\mathbf{T}_{0}=\{t_{1},...,t_{d}\}\in\mathscr{P}_{0}(\mathbf{T})$,
then they satisfy
\begin{equation}
\int_{E^{\mathbf{T}_{0}}}f(x)\gamma^{1}\circ\mathfrak{p}_{\mathbf{T}_{0}}^{-1}(dx)=\int_{E^{\mathbf{T}_{0}}}f(x)\gamma^{2}\circ\mathfrak{p}_{\mathbf{T}_{0}}^{-1}(dx)\label{eq:Sko_WLP_FDD_Int_Same}
\end{equation}
for all $\mathbf{T}_{0}=\{t_{1},...,t_{d}\}\in\mathscr{P}_{0}(\mathbf{R}^{+})$.
\end{enumerate}
\end{lem}
\begin{proof}
(a) For ease of notation, we define $\phi_{\mathbf{S}}\circeq f\circ\mathfrak{p}_{\mathbf{S}}$
for each $\mathbf{S}=\{s_{1},...,s_{d}\}\in\mathscr{P}_{0}(\mathbf{R}^{+})$.
The denseness of $\mathbf{T}$ in $\mathbf{R}^{+}$ allows us to take
$\mathbf{T}_{p}=\{t_{p,1},...,t_{p,d}\}\in\mathscr{P}_{0}(\mathbf{T})$
for each $p\in\mathbf{N}$ such that
\begin{equation}
\lim_{p\rightarrow\infty}\sup_{1\leq i\leq d}\left|t_{i}-t_{p,i}\right|=0.\label{eq:Time_Set_Approx}
\end{equation}
It follows that
\begin{equation}
\lim_{p\rightarrow\infty}\left|\phi_{\mathbf{T}_{0}}(x)-\phi_{\mathbf{T}_{p}}(x)\right|=0,\;\forall x\in D(\mathbf{R}^{+};E)\label{eq:Sko_Meas_FDD_TF_Approx}
\end{equation}
by (\ref{eq:Time_Set_Approx}), the right-continuity of $x\in S$,
Fact \ref{fact:Seq_Prod_Conv} and the continuity of $f$. The boundeness
of $f$ implies
\begin{equation}
\sup_{p\in\mathbf{N}}\left\Vert \phi_{\mathbf{T}_{p}}\right\Vert _{\infty}\leq\Vert f\Vert_{\infty}<\infty.\label{eq:Sko_Meas_FDD_TF_Bounded}
\end{equation}
Now, we have by (\ref{eq:Sko_Meas_FDD_TF_Approx}), (\ref{eq:Sko_Meas_FDD_TF_Bounded})
and the Dominated Convergence Theorem that
\begin{equation}
\begin{aligned}\lim_{n\rightarrow\infty}\int_{E^{\mathbf{T}_{0}}}f(x)\mu\circ\mathfrak{p}_{\mathbf{T}_{p}}^{-1}(dx) & =\lim_{n\rightarrow\infty}\int_{S}\phi_{\mathbf{T}_{p}}(y)\mu(dy)\\
 & =\int_{S}\phi_{\mathbf{T}_{0}}(y)\mu(dy)=\int_{E^{\mathbf{T}_{0}}}f(x)\mu\circ\mathfrak{p}_{\mathbf{T}_{0}}^{-1}(dx).
\end{aligned}
\label{eq:Check_Sko_Int_Time_Approx}
\end{equation}

(b) follows immediately by (a) (with $\mu=\gamma^{1}$ or $\gamma^{2}$).\end{proof}

\begin{lem}
\label{lem:Sko_WLP_FDD_Uni}Let $E$ be a Tychonoff space, $\mathbf{S}^{1}$,
$\mathbf{S}^{2}$ and $\mathbf{S}$ be dense subsets of $\mathbf{R}^{+}$,
$\mathcal{D}\subset C_{b}(E;\mathbf{R})$ and $\{\mu_{n}\}_{n\in\mathbf{N}}\cup\{\gamma^{1},\gamma^{2}\}\subset\mathfrak{M}^{+}(D(\mathbf{R}^{+};E),\mathscr{B}(E)^{\otimes\mathbf{R}^{+}}|_{D(\mathbf{R}^{+};E)})$
satisfy (\ref{eq:Sko_WLP_FDD_Uni_Int_Test}) for each $f\in\mathfrak{mc}[\Pi^{\mathbf{T}_{0}}(\mathcal{D})]\cup\{1\}$,
$\mathbf{T}_{0}\in\mathscr{P}_{0}(\mathbf{S}^{i})$ and $i=1,2$.
Then:

\renewcommand{\labelenumi}{(\alph{enumi})}
\begin{enumerate}
\item If $\mathcal{D}$ strongly separates points on $E$, then (\ref{eq:Sko_WLP_FDD_Uni})
holds.
\item If $\mathcal{D}$ separates points on $E$, and if each of $\gamma^{1}\circ\mathfrak{p}_{\mathbf{T}_{0}}^{-1}$
and $\gamma^{2}\circ\mathfrak{p}_{\mathbf{T}_{0}}^{-1}$ is $\mathbf{m}$-tight
for all $\mathbf{T}_{0}\in\mathscr{P}_{0}(\mathbf{S})$, then (\ref{eq:Sko_WLP_FDD_Uni})
holds.
\end{enumerate}
\end{lem}
\begin{proof}
For ease of notation, we define $\phi_{\mathbf{T}}\circeq f\circ\mathfrak{p}_{\mathbf{T}}$
for each $\mathbf{T}\in\mathscr{P}_{0}(\mathbf{R}^{+})$ and $f\in M_{b}(E^{\mathbf{T}};\mathbf{R})$.
The proof is divided into four steps.

\textit{Step 1: Show }(\ref{eq:Sko_WLP_FDD_Int_Same})\textit{ for
each $f\in\mathfrak{mc}[\Pi^{\mathbf{T}_{0}}(\mathcal{D})]\cup\{1\}$
and $\mathbf{T}_{0}=\{t_{1},...,t_{d}\}\in\mathscr{P}_{0}(\mathbf{S}^{1})$}.
Note \ref{note:Pi^d(D)_Mb_Cb} argued that $f\in C_{b}(E^{d};\mathbf{R})$.
By Lemma \ref{lem:Sko_Meas_FDD_Int_Test} (a) (with $\mathbf{T}=\mathbf{S}^{2}$),
there exists a $\mathbf{T}_{p}=\{t_{p,1},...,t_{p,d}\}\in\mathscr{P}_{0}(\mathbf{S}^{2})$
for each $p\in\mathbf{N}$ such that
\begin{equation}
\begin{aligned} & \lim_{p\rightarrow\infty}\int_{D(\mathbf{R}^{+};E)}\left(\phi_{\mathbf{T}_{p}}-\phi_{\mathbf{T}_{0}}\right)(x)\gamma^{2}(dx)\\
 & =\lim_{p\rightarrow\infty}\int_{D(\mathbf{R}^{+};E)}\left(\phi_{\mathbf{T}_{p}}-\phi_{\mathbf{T}_{0}}\right)(x)\mu_{n}(dx)=0,\;\forall n\in\mathbf{N}.
\end{aligned}
\label{eq:Check_Sko_WLP_FDD_Uni_1}
\end{equation}
From \textit{(\ref{eq:Sko_WLP_FDD_Uni_Int_Test})} we get
\begin{equation}
\lim_{n\rightarrow\infty}\int_{D(\mathbf{R}^{+};E)}\phi_{\mathbf{T}_{0}}(x)\mu_{n}(dx)=\int_{D(\mathbf{R}^{+};E)}\phi_{\mathbf{T}_{0}}(x)\gamma^{1}(dx)\label{eq:Check_Sko_WLP_FDD_Uni_2}
\end{equation}
and
\begin{equation}
\lim_{n\rightarrow\infty}\int_{D(\mathbf{R}^{+};E)}\phi_{\mathbf{T}_{p}}(x)\mu_{n}(dx)=\int_{D(\mathbf{R}^{+};E)}\phi_{\mathbf{T}_{p}}(x)\gamma^{2}(dx),\;\forall p\in\mathbf{N}.\label{eq:Check_Sko_WLP_FDD_Uni_3}
\end{equation}

Let $\epsilon\in(0,\infty)$ be arbitrary and $n_{0}\circeq1$. By
(\ref{eq:Check_Sko_WLP_FDD_Uni_2}) and (\ref{eq:Check_Sko_WLP_FDD_Uni_3}),
we inductively choose an $n_{p}\in\mathbf{N}\cap(n_{p-1},\infty)$
for each $p\in\mathbf{N}$ such that
\begin{equation}
\begin{aligned} & \left|\int_{D(\mathbf{R}^{+};E)}\phi_{\mathbf{T}_{0}}(x)\mu_{n_{p}}(dx)-\int_{D(\mathbf{R}^{+};E)}\phi_{\mathbf{T}_{0}}(x)\gamma^{1}(dx)\right|\vee\\
 & \left|\int_{D(\mathbf{R}^{+};E)}\phi_{\mathbf{T}_{p}}(x)\mu_{n_{p}}(dx)-\int_{D(\mathbf{R}^{+};E)}\phi_{\mathbf{T}_{p}}(x)\gamma^{2}(dx)\right|<\frac{\epsilon}{2}.
\end{aligned}
\label{eq:Check_Sko_WLP_FDD_Uni_4}
\end{equation}
From Triangle Inequality and (\ref{eq:Check_Sko_WLP_FDD_Uni_4}) it
follows that
\begin{equation}
\begin{aligned} & \left|\int_{E^{\mathbf{T}_{0}}}f(x)\gamma^{1}\circ\mathfrak{p}_{\mathbf{T}_{0}}^{-1}(dx)-\int_{E^{\mathbf{T}_{0}}}f(x)\gamma^{2}\circ\mathfrak{p}_{\mathbf{T}_{0}}^{-1}(dx)\right|\\
 & =\left|\int_{D(\mathbf{R}^{+};E)}\phi_{\mathbf{T}_{0}}(x)\gamma^{1}(dx)-\int_{D(\mathbf{R}^{+};E)}\phi_{\mathbf{T}_{0}}(x)\gamma^{2}(dx)\right|\\
 & \leq\left|\int_{D(\mathbf{R}^{+};E)}\left(\phi_{\mathbf{T}_{0}}-\phi_{\mathbf{T}_{p}}\right)(x)\gamma^{2}(dx)\right|\\
 & +\left|\int_{D(\mathbf{R}^{+};E)}\left(\phi_{\mathbf{T}_{p}}-\phi_{\mathbf{T}_{0}}\right)(x)\mu_{n_{p}}(dx)\right|+\epsilon,\;\forall p\in\mathbf{N}.
\end{aligned}
\label{eq:Check_Sko_WLP_FDD_Uni_5}
\end{equation}
Now, (\ref{eq:Sko_WLP_FDD_Int_Same}) follows by (\ref{eq:Check_Sko_WLP_FDD_Uni_1}),
letting $p\uparrow\infty$ in (\ref{eq:Check_Sko_WLP_FDD_Uni_5})
and then letting $\epsilon\downarrow0$.

\textit{Step 2: Show (\ref{eq:Sko_WLP_FDD_Int_Same}) for each $f\in\mathfrak{mc}[\Pi^{\mathbf{T}_{0}}(\mathcal{D})]\cup\{1\}$
and $\mathbf{T}_{0}=\{t_{1},...,t_{d}\}\in\mathscr{P}_{0}(\mathbf{R}^{+})$}.
This step follows by Step 1, the denseness of $\mathbf{S}^{1}$ in
$\mathbf{R}^{+}$ and Lemma \ref{lem:Sko_Meas_FDD_Int_Test} (b) (with
$\mathbf{T}=\mathbf{S}^{1}$).

\textit{Step 3: Verify $\gamma^{1}\circ\mathfrak{p}_{\mathbf{T}_{0}}^{-1}=\gamma^{2}\circ\mathfrak{p}_{\mathbf{T}_{0}}^{-1}$
for each $\mathbf{T}_{0}\in\mathscr{P}_{0}(\mathbf{R}^{+})$ in (a)}.
For each $i=1,2$, we let $(D(\mathbf{R}^{+};E),\mathscr{A}_{n_{0}},\nu^{i})$
be the completion of $(D(\mathbf{R}^{+};E),\sigma[\mathscr{J}(E)],\gamma^{i})$
and find by Lemma \ref{lem:Sko_FDD_BExt} (with $\mu=\gamma^{i}$
and $\nu=\nu^{i}$) that $\nu^{i}\circ\mathfrak{p}_{\mathbf{T}_{0}}^{-1}$%
\footnote{Herein, the domain of $\nu^{i}\circ\mathfrak{p}_{\mathbf{T}_{0}}^{-1}$
is larger than $\mathscr{B}(E^{\mathbf{T}_{0}})$, so we consider
$\nu^{i}\circ\mathfrak{p}_{\mathbf{T}_{0}}^{-1}$ as a member of $\mathcal{M}^{+}(E^{\mathbf{T}_{0}})$.%
} is a Borel extension of $\gamma^{i}\circ\mathfrak{p}_{\mathbf{T}_{0}}^{-1}$.
It follows that
\begin{equation}
\int_{E^{\mathbf{T}_{0}}}f(x)\nu^{1}\circ\mathfrak{p}_{\mathbf{T}_{0}}^{-1}(dx)=\int_{E^{\mathbf{T}_{0}}}f(x)\nu^{2}\circ\mathfrak{p}_{\mathbf{T}_{0}}^{-1}(dx)\label{eq:Sko_WLP_FDD_BExt_Int_Same}
\end{equation}
for all $f\in\mathfrak{mc}[\Pi^{\mathbf{T}_{0}}(\mathcal{D})]\cup\{1\}$
by Step 2 and Fact \ref{fact:BExt_Same_Int} (with $d=\aleph(\mathbf{T}_{0})$,
$\mu=\gamma^{i}\circ\mathfrak{p}_{\mathbf{T}_{0}}^{-1}$ and $\nu_{1}=\nu^{i}\circ\mathfrak{p}_{\mathbf{T}_{0}}^{-1}$).
Then, $\nu^{1}\circ\mathfrak{p}_{\mathbf{T}_{0}}^{-1}=\nu^{2}\circ\mathfrak{p}_{\mathbf{T}_{0}}^{-1}$
by Lemma \ref{lem:Sep_Meas} (a) (with $d=\aleph(\mathbf{T}_{0})$)
and Fact \ref{fact:Sep_CD} (a) (with $E=E^{\mathbf{T}_{0}}$ and
$\mathcal{D}=\mathfrak{mc}[\Pi^{\mathbf{T}_{0}}(\mathcal{D})]\cup\{1\}$),
which of course implies \textit{$\gamma^{1}\circ\mathfrak{p}_{\mathbf{T}_{0}}^{-1}=\gamma^{2}\circ\mathfrak{p}_{\mathbf{T}_{0}}^{-1}$}.

\textit{Step 4: Verify $\gamma^{1}\circ\mathfrak{p}_{\mathbf{T}_{0}}^{-1}=\gamma^{2}\circ\mathfrak{p}_{\mathbf{T}_{0}}^{-1}$
for each $\mathbf{T}_{0}\in\mathscr{P}_{0}(\mathbf{R}^{+})$ in (b)}.
When $\mathbf{T}_{0}\in\mathscr{P}_{0}(\mathbf{S})$, $\gamma^{1}\circ\mathfrak{p}_{\mathbf{T}_{0}}^{-1}=\gamma^{1}\circ\mathfrak{p}_{\mathbf{T}_{0}}^{-1}$
by Step 2 and Lemma \ref{lem:Tight_Meas_Identical} (b) (with $d=\aleph(\mathbf{T}_{0})$)
and so (\ref{eq:Sko_WLP_FDD_Int_Same}) holds for all $f\in\mathfrak{mc}[\Pi^{\mathbf{T}_{0}}(C_{b}(E;\mathbf{R}))]$.
For general $\mathbf{T}_{0}\in\mathscr{P}_{0}(\mathbf{R}^{+})$, the
key equality (\ref{eq:Sko_WLP_FDD_Int_Same}) holds for all $f\in\mathfrak{mc}[\Pi^{\mathbf{T}_{0}}(C_{b}(E;\mathbf{R}))]$
by (\ref{eq:ca(Pi^d(D))_Cb}) (with $\mathcal{D}=C_{b}(E;\mathbf{R})$
and $d=\aleph(\mathbf{T}_{0})$), the denseness of $\mathbf{S}$ in
$\mathbf{R}^{+}$ and Lemma \ref{lem:Sko_Meas_FDD_Int_Test} (b) (with
$\mathbf{T}=\mathbf{S}$). $C_{b}(E;\mathbf{R})$ strongly separates
points on Proposition \ref{prop:CR} (a, c). Then, one follows the
argument of Step 3 (with $\mathcal{D}=C_{b}(E;\mathbf{R})$) to show
\textit{$\gamma^{1}\circ\mathfrak{p}_{\mathbf{T}_{0}}^{-1}=\gamma^{2}\circ\mathfrak{p}_{\mathbf{T}_{0}}^{-1}$.}\end{proof}

\begin{lem}
\label{lem:Path_Space_RV}Let $E$ be a topological space, $(E_{0},\mathscr{O}_{E}(E_{0}))$
be a Tychonoff subspace of $E$, $y_{0}\in S_{0}\subset\mathbb{D}_{0}\circeq D(\mathbf{R}^{+};E_{0},\mathscr{O}_{E}(E_{0}))$,
$\mathscr{U}\circeq\mathscr{B}(E)^{\otimes\mathbf{R}^{+}}$ and $X$
be a mapping from $(\Omega,\mathscr{F},\mathbb{P})$ to $E^{\mathbf{R}^{+}}$.
Then:

\renewcommand{\labelenumi}{(\alph{enumi})}
\begin{enumerate}
\item If $\Omega\backslash\{X=Z\}\in\mathscr{N}(\mathbb{P})$ for some $Z\in M(\Omega,\mathscr{F};\mathbb{D}_{0})$,
then $X$ is an $E$-valued c$\grave{\mbox{a}}$dl$\grave{\mbox{a}}$g
process.
\item If $X$ is an $E$-valued process, $\mathbb{P}(X\in S_{0})=1$ and
$S_{0}$ satisfies $\mathscr{B}(\mathbb{D}_{0})|_{S_{0}}=\mathscr{U}|_{S_{0}}$,
then
\begin{equation}
Y\circeq\mathfrak{var}\left(X;\Omega,X^{-1}(S_{0}),y_{0}\right)\in M\left(\Omega,\mathscr{F};S_{0},\mathscr{O}_{\mathbb{D}_{0}}(S_{0})\right)\label{eq:Var_Path_Space_RV}
\end{equation}
and $\mathbb{P}(X=Y\in S_{0})=1$.
\end{enumerate}
\end{lem}
\begin{proof}
(a) follows by Lemma \ref{lem:Sko_Proj} (b) (with $E=(E_{0},\mathscr{O}_{E}(E_{0}))$)
and Lemma \ref{lem:var(X)} (a) (with $E=E^{\mathbf{R}^{+}}$ and
$S=\mathbb{D}_{0}$).

(b) follows by Lemma \ref{lem:var(X)} (b, c) (with $(E,S,\mathscr{U}^{\prime})=(E^{\mathbf{R}^{+}},S_{0},\mathscr{B}_{\mathbb{D}_{0}}(S_{0}))$).\end{proof}

\begin{prop}
\label{prop:M(E)_RV_Proc_TF}Let $E$ be a topological space, $(\Omega,\mathscr{F},\{\mathscr{G}_{t}\}_{t\geq0},\mathbb{P})$
be a stochastic basis, $k\in\mathbf{N}$, $\xi\in M(\Omega,\mathscr{F};\mathcal{M}^{+}(E))$
and $(\Omega,\mathscr{F},\mathbb{P};X)$ be an $\mathcal{M}^{+}(E)$-valued
process. In addition, suppose either of the following hypotheses is
true:

\renewcommand{\labelenumi}{(\roman{enumi})}
\begin{enumerate}
\item $f\in C_{b}(E;\mathbf{R}^{k})$.
\item $E$ is a perfectly normal%
\footnote{The notion of perfectly normal space was mentioned in \S \ref{lem:Strong_Topo}.%
} (especially metrizable or Polish) space and $f\in M_{b}(E;\mathbf{R}^{k})$.
\end{enumerate}
Then, $f^{*}\circ\xi\in M(\Omega,\mathscr{F};\mathbf{R}^{k})$ and
$\varpi(f^{*})\circ X$ is an $\mathbf{R}^{k}$-valued process. If,
in addition, $X$ is a $\mathscr{G}_{t}$-adapted, measurable or $\mathscr{G}_{t}$-progressive
process, then $\varpi(f^{*})\circ X$ also has the corresponding measurability.
\end{prop}
\begin{proof}
Under (i), $f^{*}\in C_{b}(\mathcal{M}^{+}(E);\mathbf{R}^{k})$ by
the definition of weak topology and Fact \ref{fact:Prod_Map_2} (b).
Under (ii), $f^{*}\in M_{b}(\mathcal{M}^{+}(E);\mathbf{R}^{k})$ by
Lemma \ref{lem:Strong_Topo} and Fact \ref{fact:Prod_Map_2} (b).
The result now follows by Fact \ref{fact:Proc_Path_Mapping} (a) (with
$E=\mathcal{M}^{+}(E)$, $S=\mathbf{R}^{k}$ and $f=f^{*}$).\end{proof}

\section{\label{sec:Aux_Rep}Auxiliary results about replication}
\begin{fact}
\label{fact:f+_Rep}Let $E$ be a topological space, $(E_{0},\mathcal{F};\widehat{E},\widehat{\mathcal{F}})$
be a base over $E$ and $d,k\in\mathbf{N}$. If $f\in C(E^{d};\mathbf{R}^{k})$
has a replica $\widehat{f}$, then:

\renewcommand{\labelenumi}{(\alph{enumi})}
\begin{enumerate}
\item $\Vert\widehat{f}\Vert_{\infty}=\Vert f|_{E_{0}^{d}}\Vert_{\infty}\leq\Vert f\Vert_{\infty}$.
\item $\widehat{f^{+}}=\widehat{f}^{+}$ and $\widehat{f^{-}}=\widehat{f}^{-}$.
\end{enumerate}
\end{fact}
\begin{proof}
(a) follows by the fact $f|_{E_{0}^{d}}=\widehat{f}|_{E_{0}^{d}}$,
the denseness of $E_{0}$ in $\widehat{E}$ and the continuities of
$f$ and $\widehat{f}$.

(b) follows by the facts $\widetilde{f^{+}}=\widehat{f}^{+}|_{E_{0}^{d}}$,
$\widetilde{f^{-}}=\widehat{f}^{-}|_{E_{0}^{d}}$ and the continuities
of $\widehat{f}^{+}$, $\widehat{f}^{-}$.\end{proof}

\begin{lem}
\label{lem:m-Tight_Base}Let $E$ be a topological space, $\mathcal{D}\subset C(E;\mathbf{R})$
separate points on $E$, $d\in\mathbf{N}$, $\mathbf{I}$ be a countable
index set and $\Gamma_{i}\subset\mathfrak{M}^{+}(E^{d},\mathscr{B}(E)^{\otimes d})$
be $\mathbf{m}$-tight for each $i\in\mathbf{I}$. Then, there exists
a base $(E_{0},\mathcal{F};\widehat{E},\widehat{\mathcal{F}})$ over
$E$ such that $E_{0}\in\mathscr{K}_{\sigma}^{\mathbf{m}}(E)$ and
$\Gamma_{i}$ is tight in $(E_{0}^{d},\mathscr{O}_{E}(E_{0})^{d})$
for all $i\in\mathbf{I}$. In particular, $\mathcal{F}$ can be taken
within $\mathcal{D}\cup\{1\}$ when $\mathcal{D}\subset C_{b}(E;\mathbf{R})$.
\end{lem}
\begin{proof}
Without loss of generality, we let $\mathbf{I}=\mathbf{N}$. By the
$\mathbf{m}$-tightness of each $\Gamma_{i}$, there exist $\{K_{p,i}\}_{p,i\in\mathbf{N}}\subset\mathscr{K}^{\mathbf{m}}(E^{d})$
satisfying
\begin{equation}
\sup_{\mu\in\Gamma_{i}}\mu(E^{d}\backslash K_{p,i})\geq1-2^{-p},\;\forall p,i\in\mathbf{N}.\label{eq:Pick_Kpi_m-Tight_Base}
\end{equation}
$E$ is a Hausdorff space by Proposition \ref{prop:Fun_Sep_1} (e)
(with $A=E$).
\begin{equation}
\left\{ K_{p,i,j}\circeq\mathfrak{p}_{j}(K_{p,i}):1\leq j\leq d,p,i\in\mathbf{N}\right\} \subset\mathscr{K}^{\mathbf{m}}(E)\subset\mathscr{B}(E)\label{eq:K_p_i_j_Metrizable_Compact}
\end{equation}
by Proposition \ref{prop:Separability} (c) and Lemma \ref{lem:MC_Prod}
(b) (with $A=K_{p,i}$). So, 
\begin{equation}
E_{0}\circeq\bigcup_{i\in\mathbf{N}}\bigcup_{p\in\mathbf{N}}\bigcup_{j=1}^{d}K_{p,i,j}\in\mathscr{K}_{\sigma}^{\mathbf{m}}(E).\label{eq:E0=00003DUnion_K_p_i_j}
\end{equation}
We have by Corollary \ref{cor:Compact_Prod} (a) and (\ref{eq:Pick_Kpi_m-Tight_Base})
that
\begin{equation}
\prod_{j=1}^{d}K_{p,i,j}\subset\mathscr{K}\left(E_{0}^{d},\mathscr{O}_{E}(E_{0})^{d}\right),\;\forall p,i\in\mathbf{N}\label{eq:K_p_i_j_Prod_Compact_E0d}
\end{equation}
and
\begin{equation}
\mu(E_{0}^{d})\geq\mu\left(\prod_{j=1}^{d}K_{p,i,j}\right)\geq\mu(K_{p,i})\geq1-2^{p},\;\forall\mu\in\Gamma_{i},p,i\in\mathbf{N}.\label{eq:Check_Mu_Tight_E0D}
\end{equation}
thus proving the tightness of each $\Gamma_{i}$ in $(E_{0}^{d},\mathscr{O}_{E}(E_{0})^{d})$.
$E_{0}$ is a $\mathcal{D}$-baseable subset of $E$ by (\ref{eq:E0=00003DUnion_K_p_i_j})
and Proposition \ref{prop:Sigma_MC} (b, e) (with $A=E_{0}$). Now,
the result follows by Lemma \ref{lem:Base_Construction} (a, c) (with
$\mathcal{D}_{0}=\varnothing$).\end{proof}

\begin{fact}
\label{fact:TF_(T,E0)-Mod}Let $E$ be a topological space, $(E_{0},\mathcal{F};\widehat{E},\widehat{\mathcal{F}})$
be base over $E$, $(\Omega,\mathscr{F},\mathbb{P};X)$ be an $E$-valued
process, $(\Omega,\mathscr{F},\mathbb{P};Y)$ be an $\widehat{E}$-valued
process and $\mathbf{T}\subset\mathbf{R}^{+}$. If $X$ satisfies
(\ref{eq:T-Base}), and if
\begin{equation}
\mathbb{P}\left(\bigotimes\mathcal{F}\circ X_{t}=\bigotimes\widehat{\mathcal{F}}\circ Y_{t}\right)\geq\mathbb{P}\left(\bigotimes\mathcal{F}\circ X_{t}\in\bigotimes\widehat{\mathcal{F}}(\widehat{E})\right),\;\forall t\in\mathbf{T},\label{eq:RepProc_T}
\end{equation}
then $X$ and $Y$ satisfy (\ref{eq:Lim_Proc_(T,E0)-Mod}).
\end{fact}
\begin{proof}
 It follows by (\ref{eq:RepProc_T}), (\ref{eq:F_Fhat_Coincide})
and (\ref{eq:Base_Imb}) that
\begin{equation}
\begin{aligned}1 & =\mathbb{P}\left(X_{t}\in E_{0}\right)\\
 & =\mathbb{P}\left(\bigotimes\widehat{\mathcal{F}}\circ Y_{t}=\bigotimes\mathcal{F}\circ X_{t}\in\bigotimes\widehat{\mathcal{F}}(\widehat{E}),X_{t}\in E_{0}\right)\\
 & \leq\mathbb{P}\left(Y_{t}=\left(\bigotimes\widehat{\mathcal{F}}\right)^{-1}\circ\bigotimes\mathcal{F}\circ X_{t}=X_{t}\in E_{0}\right),\;\forall t\in\mathbf{T}.
\end{aligned}
\label{eq:Check_(T,E0)-Mod}
\end{equation}
\end{proof}

\begin{fact}
\label{fact:Proc_FDD_BExt}Let $E$ be a topological space, $(E_{0},\mathcal{F};\widehat{E},\widehat{\mathcal{F}})$
be a base over $E$ and $\{(\Omega^{i},\mathscr{F}^{i},\mathbb{P}^{i};X^{i})\}_{i\in\mathbf{I}}$
be $E$-valued processes satisfying $\mathbf{T}$-PSMTC in $E_{0}$%
\footnote{The notion of $\mathbf{T}$-PSMTC was introduced in Definition \ref{def:Proc_Reg}.%
}. Then for each $\mathbf{T}_{0}\in\mathscr{P}_{0}(\mathbf{T})$:

\renewcommand{\labelenumi}{(\alph{enumi})}
\begin{enumerate}
\item $\{X_{\mathbf{T}_{0}}^{i}\}_{i\in\mathbf{I}}$ is sequentially $\mathbf{m}$-tight
in $E_{0}^{\mathbf{T}_{0}}$ and $\mathbb{P}^{i}(X_{\mathbf{T}_{0}}^{i}\in E_{0}^{\mathbf{T}_{0}})=1$.
\item $\{\mu^{i}=\mathfrak{be}(\mathbb{P}^{i}\circ X_{\mathbf{T}_{0}}^{-1})\}_{i\in\mathbf{I}\backslash\mathbf{I}_{\mathbf{T}_{0}}}$
exists for some $\mathbf{I}_{\mathbf{T}_{0}}\in\mathscr{P}_{0}(\mathbf{I})$.
\end{enumerate}
\end{fact}
\begin{proof}
Let $\Gamma\circeq\{\mathbb{P}^{i}\circ X_{\mathbf{T}_{0}}^{-1}\}_{i\in\mathbf{I}}$
and $A\circeq E_{0}^{\mathbf{T}_{0}}$. Then, (a) follows by Lemma
\ref{lem:Tightness_Prod} (b) (with $\mathbf{I}=\mathbf{T}_{0}$,
$S_{i}=E$ and $A_{i}=E_{0}$) and Fact \ref{fact:Seq_Tight_Support}
(with $E=E^{\mathbf{T}_{0}}$ and $\mathscr{U}=\mathscr{B}(E)^{\otimes\mathbf{T}_{0}}$).
(b) follows by Lemma \ref{lem:Base} (e) (with $A=E_{0}$) and Proposition
\ref{prop:Seq_m-Tight_BExt} (with $\mathbf{I}=\mathbf{T}_{0}$ and
$S_{i}=E$).\end{proof}

\begin{lem}
\label{lem:Proc_Rep}Let $E$ be a topological space, $(E_{0},\mathcal{F};\widehat{E},\widehat{\mathcal{F}})$
be a base over $E$, $\mathbf{T}\subset\mathbf{R}^{+}$, $(\Omega,\mathcal{F},\mathbb{P};X)$
be $E$-valued process and $(\Omega,\mathcal{F},\mathbb{P};Y)$ be
an $\widehat{E}$-valued process. Then:

\renewcommand{\labelenumi}{(\alph{enumi})}
\begin{enumerate}
\item If
\begin{equation}
\inf_{f\in\mathcal{F},t\in\mathbf{T}}\mathbb{P}\left(f\circ X_{t}=\widehat{f}\circ Y_{t}\right)=1,\label{eq:Lim_Proc_TF_T-Mod}
\end{equation}
then
\begin{equation}
\mathbb{P}\left(f\circ X_{\mathbf{T}_{0}}=\widehat{f}\circ Y_{\mathbf{T}_{0}}\right)=1,\;\forall f\in\mathfrak{ca}\left[\Pi^{\mathbf{T}_{0}}(\mathcal{F})\right],\mathbf{T}_{0}\in\mathscr{P}_{0}(\mathbf{T}).\label{eq:Lim_Proc_TF_T-Mod_FD}
\end{equation}
Moreover, (\ref{eq:Lim_Proc_(T,E0)-Mod}) implies (\ref{eq:Lim_Proc_TF_T-Mod}).
\item If (\ref{eq:Lim_Proc_TF_T-Mod}) holds (especially (\ref{eq:Lim_Proc_(T,E0)-Mod})
holds) and $Y$ is c$\grave{\mbox{a}}$dl$\grave{\mbox{a}}$g, then
$X$ is $(\mathbf{T},\mathcal{F})$-c$\grave{\mbox{a}}$dl$\grave{\mbox{a}}$g.
\item If $\mathbb{E}[f\circ X_{\mathbf{T}_{0}}]=\mathbb{E}[\widehat{f}\circ Y_{\mathbf{T}_{0}}]$
for all $f\in\mathfrak{ca}[\Pi^{\mathbf{T}_{0}}(\mathcal{F})]$ and
$\mathbf{T}_{0}\in\mathscr{P}_{0}(\mathbf{R}^{+})$ (especially (\ref{eq:Lim_Proc_TF_T-Mod})
or (\ref{eq:Lim_Proc_(T,E0)-Mod}) holds) and $X$ is stationary,
then $Y$ is stationary.
\item If $\mathbb{E}[f\circ X_{\mathbf{T}_{0}}]=\mathbb{E}[\widehat{f}\circ Y_{\mathbf{T}_{0}}]$
for all $f\in\mathfrak{ca}[\Pi^{\mathbf{T}_{0}}(\mathcal{F})]$ and
$\mathbf{T}_{0}\in\mathscr{P}_{0}(\mathbf{T})$ (especially (\ref{eq:Lim_Proc_TF_T-Mod})
or (\ref{eq:Lim_Proc_(T,E0)-Mod}) holds), $\mathbf{T}$ is conull,
$X$ is stationary and $Y$ is c$\grave{\mbox{a}}$dl$\grave{\mbox{a}}$g,
then $Y$ is stationary.
\item If $A\in\mathscr{B}^{\mathbf{s}}(E)$ (especially $A\in\mathscr{K}_{\sigma}(E)$)
satisfies $A\subset E_{0}$ and (\ref{eq:Lim_Proc_(R+,A)-Mod}) holds,
then $\mathscr{F}^{X}=\mathscr{F}^{Y}$. If, in addition, $Y$ is
stationary, then $X$ is stationary.
\item If (\ref{eq:Lim_Proc_(T,E0)-Mod}) holds, and if $f\in M_{b}(E^{\mathbf{T}_{0}};\mathbf{R})$
and $\mathbf{T}_{0}\in\mathscr{P}_{0}(\mathbf{T})$ satisfy $\overline{f}\in M_{b}(\widehat{E}^{\mathbf{T}_{0}};\mathbf{R})$
(especially if $E_{0}^{\mathbf{T}_{0}}\in\mathscr{B}^{\mathbf{s}}(E^{\mathbf{T}_{0}})$),
then
\begin{equation}
\mathbb{P}\left(f\circ X_{\mathbf{T}_{0}}=\overline{f}\circ Y_{\mathbf{T}_{0}}\right)=1.\label{eq:Lim_Proc_TF_T-Mod_Bar}
\end{equation}

\end{enumerate}
\end{lem}
\begin{proof}
(a) (\ref{eq:Lim_Proc_TF_T-Mod}) implies (\ref{eq:Lim_Proc_TF_T-Mod_FD})
by properties of uniform convergence. (\ref{eq:Lim_Proc_(T,E0)-Mod})
implies (\ref{eq:Lim_Proc_TF_T-Mod}) by (\ref{eq:F_Fhat_Coincide}).

(b) $\{\varpi(\widehat{f})\circ Y\}_{f\in\mathcal{F}}$ are all c$\grave{\mbox{a}}$dl$\grave{\mbox{a}}$g
processes by Fact \ref{fact:Cadlag_Proc} (a) (with $E=\widehat{E}$,
$S=\mathbf{R}$ and $X=Y$). Then, (b) follows by (a).

(c) One finds by (a) (with $\mathbf{T}=\mathbf{R}^{+}$) and the stationarity
of $X$ that
\begin{equation}
\mathbb{E}\left[\widehat{f}\circ Y_{\mathbf{T}_{0}}-\widehat{f}\circ Y_{\mathbf{T}_{0}+c}\right]=\mathbb{E}\left[f\circ X_{\mathbf{T}_{0}}-f\circ X_{\mathbf{T}_{0}+c}\right]=0\label{eq:Check_R+_Mod_Sta}
\end{equation}
for all $c\in(0,\infty)$ and $\mathbf{T}_{0}\in\mathscr{P}_{0}(\mathbf{R}^{+})$.
Then, (c) follows by Corollary \ref{cor:Base_Sep_Meas} (a) (with
$d=\aleph(\mathbf{T}_{0})$ and $A=\widehat{E}^{d}$).

(d) Fixing $\mathbf{T}_{0}\in\mathscr{P}_{0}(\mathbf{R}^{+})$, one
finds by (a) and the stationarity of $X$ that (\ref{eq:Check_R+_Mod_Sta})
holds for all $c$ in the conull set
\begin{equation}
\mathbf{S}_{\mathbf{T}_{0}}\circeq\bigcap_{t\in\mathbf{T}_{0}}\left\{ c\in(0,\infty):t+c\in\mathbf{T}\right\} .\label{eq:S_T0}
\end{equation}
Then, (d) follows by a similar argument to the proof of Proposition
\ref{prop:FR_FC_WC} (c).

(e) Any $A\in\mathscr{K}_{\sigma}(E)$ satisfying $A\subset E_{0}$
belongs to $\mathscr{B}^{\mathbf{s}}(E)$ by Corollary \ref{cor:Base_Compact}
(b) (with $d=1$). For each fixed $t\in\mathbf{R}^{+}$, we let $\Omega_{t}^{0}\circeq\{\omega\in\Omega:X_{t}(\omega)=Y_{t}(\omega)\in A\}$
and find by (\ref{eq:Lim_Proc_(R+,A)-Mod}), the $\mathbb{P}$-completeness
of $\mathscr{F}$%
\footnote{Completeness of measure space was specified in \ref{sub:Meas}. Completeness
of filtration was specified in \S \ref{sec:Proc}.%
} and Lemma \ref{lem:SB_Base} (a) (with $d=1$) that $\Omega\backslash\Omega_{t}^{0}\in\mathscr{N}(\mathbb{P})\subset\mathscr{F}$,
\begin{equation}
\begin{aligned}X_{t}^{-1}(B)\cap\Omega_{t}^{0} & =X_{t}^{-1}(B\cap A)\cap\Omega_{t}^{0}\\
 & \in\left\{ Y_{t}^{-1}(V)\cap\Omega_{t}^{0}:V\in\mathscr{B}(\widehat{E})\right\} ,\;\forall B\in\mathscr{B}(E)
\end{aligned}
\label{eq:Check_Same_Filtration_1}
\end{equation}
and
\begin{equation}
\begin{aligned}Y_{t}^{-1}(V)\cap\Omega_{t}^{0} & =Y_{t}^{-1}(V\cap A)\cap\Omega_{t}^{0}\\
 & \in\left\{ X_{t}^{-1}(B)\cap\Omega_{t}^{0}:B\in\mathscr{B}(E)\right\} ,\;\forall V\in\mathscr{B}(\widehat{E}).
\end{aligned}
\label{eq:Check_Same_Filtration_2}
\end{equation}
Thus, $\mathscr{F}^{X}=\mathscr{F}^{Y}$ by their $\mathbb{P}$-completeness.
When $Y$ is stationary, we fix $\mathbf{T}_{0}\in\mathscr{P}_{0}(\mathbf{R}^{+})$
and find by Lemma \ref{lem:SB_Base} (d) (with $d=\aleph(\mathbf{T}_{0})$)
and (\ref{eq:Lim_Proc_(R+,A)-Mod}) that
\begin{equation}
\begin{aligned}\mathbb{P}\left(X_{\mathbf{T}_{0}}\in B\right) & =\mathbb{P}\left(Y_{\mathbf{T}_{0}}\in B\cap A^{\mathbf{T}_{0}}\right)\\
 & =\mathbb{P}\left(Y_{\mathbf{T}_{0}+c}\in B\cap A^{\mathbf{T}_{0}}\right)=\mathbb{P}\left(X_{\mathbf{T}_{0}+c}\in B\right)
\end{aligned}
\label{eq:Check_Lim_Proc_Sta}
\end{equation}
for all $B\in\mathscr{B}(E)^{\otimes\mathbf{T}_{0}}$ and $c\in(0,\infty)$,
which gives the stationarity of $X$.

(f) follows by the definition of $\overline{f}$ and Proposition \ref{prop:RepFun_Basic}
(b) (with $d=\aleph(\mathbf{T}_{0})$).\end{proof}

\begin{lem}
\label{lem:RepProc_Int_Test}Let $E$ be a topological space, $(E_{0},\mathcal{F};\widehat{E},\widehat{\mathcal{F}})$
be a base over $E$, $\mathbf{T}\subset\mathbf{R}^{+}$, $\mathcal{G}_{\mathbf{T}_{0}}\circeq\mathfrak{mc}[\Pi^{\mathbf{T}_{0}}(\mathcal{F}\backslash\{1\})]$
for each $\mathbf{T}_{0}\in\mathscr{P}_{0}(\mathbf{T})$, $\{(\Omega^{n},\mathscr{F}^{n},\mathbb{P}^{n};X^{n})\}_{n\in\mathbf{N}}$
be $E$-valued processes satisfying (\ref{eq:Common_FR-Base_N}),
$\widehat{X}^{n}\in\mathfrak{rep}(X^{n};E_{0},\mathcal{F})$ for each
$n\in\mathbf{N}$, $(\Omega,\mathscr{F},\mathbb{P};Y)$ be an $\widehat{E}$-valued
process, $(\Omega,\mathscr{F},\mathbb{P};X)$ be an $E$-valued process
satisfying (\ref{eq:FR-Base}) and $\widehat{X}\in\mathfrak{rep}(X;E_{0},\mathcal{F})$.
Then:

\renewcommand{\labelenumi}{(\alph{enumi})}
\begin{enumerate}
\item $\{X^{n}\}_{n\in\mathbf{N}}$ is $(\mathbf{T},\mathcal{F}\backslash\{1\})$-FDC
if and only if $\{\widehat{X}^{n}\}_{n\in\mathbf{N}}$ is $(\mathbf{T},\widehat{\mathcal{F}}\backslash\{1\})$-FDC.
\item If $\{X^{n}\}_{n\in\mathbf{N}}$ is $(\mathbf{T},\mathcal{F}\backslash\{1\})$-AS,
then
\begin{equation}
\lim_{n\rightarrow\infty}\mathbb{E}^{n}\left[\widehat{f}\circ\widehat{X}_{\mathbf{T}_{0}}^{n}-\widehat{f}\circ\widehat{X}_{\mathbf{T}_{0}+c}^{n}\right]=0\label{eq:RepProc_AS}
\end{equation}
for all $f\in\mathcal{G}_{\mathbf{T}_{0}}$, $\mathbf{T}_{0}\in\mathscr{P}_{0}(\mathbf{T})$
and $c$ (if any) in the set $\mathbf{S}_{\mathbf{T}_{0}}$ defined
in (\ref{eq:S_T0}).
\item (\ref{eq:Exp_Test_Y}) is equivalent to
\begin{equation}
\lim_{n\rightarrow\infty}\mathbb{E}^{n}\left[\widehat{f}\left(\widehat{X}_{\mathbf{T}_{0}}^{n}\right)\right]=\mathbb{E}\left[\widehat{f}\left(Y_{\mathbf{T}_{0}}\right)\right]\label{eq:RepProc_Exp_Test_Y}
\end{equation}
for all $f\in\mathcal{G}_{\mathbf{T}_{0}}$ and $\mathbf{T}_{0}\in\mathscr{P}_{0}(\mathbf{T})$.
\item If (\ref{eq:Lim_Proc_TF_T-Mod_FD}) holds (especially (\ref{eq:Lim_Proc_(T,E0)-Mod})
holds), then (\ref{eq:Exp_Test}) is equivalent to
\begin{equation}
\lim_{n\rightarrow\infty}\mathbb{E}^{n}\left[\widehat{f}\left(\widehat{X}_{\mathbf{T}_{0}}^{n}\right)\right]=\mathbb{E}\left[\widehat{f}\left(\widehat{X}_{\mathbf{T}_{0}}\right)\right]\label{eq:RepProc_Exp_Test}
\end{equation}
for all $f\in\mathcal{G}_{\mathbf{T}_{0}}$ and $\mathbf{T}_{0}\in\mathscr{P}_{0}(\mathbf{T})$.
\item (\ref{eq:RepProc_FC_along_T_Y}) holds if and only if (\ref{eq:RepProc_Exp_Test_Y})
holds for all $f\in\mathcal{G}_{\mathbf{T}_{0}}$ and $\mathbf{T}_{0}\in\mathscr{P}_{0}(\mathbf{T})$.
\item (\ref{eq:Exp_Test}) holds for all $f\in\mathcal{G}_{\mathbf{T}_{0}}$
and $\mathbf{T}_{0}\in\mathscr{P}_{0}(\mathbf{T})$ if and only if
\begin{equation}
\widehat{X}^{n}\xrightarrow{\quad\mathrm{D}(\mathbf{T})\quad}\widehat{X}\mbox{ as }n\uparrow\infty.\label{eq:RepProc_FC_along_T}
\end{equation}

\item (\ref{eq:FC_along_T}) implies (\ref{eq:RepProc_FC_along_T}).
\end{enumerate}
In particular, the conclusions above are when if $\{X^{n}\}_{n\in\mathbf{N}}$
(resp. $X$) satisfies the stronger condition%
\footnote{We compared these conditions in Fact \ref{fact:T-Base_FR-Base} and
Fact \ref{fact:Common_T-Base_FR-Base}.%
} (\ref{eq:Common_T-Base_N}) (resp. (\ref{eq:T-Base})) than (\ref{eq:Common_FR-Base_N})
(resp. (\ref{eq:FR-Base})).
\end{lem}
\begin{proof}
(a) - (c) follow by Fact \ref{fact:FDC_AS} (b) and Proposition \ref{prop:RepProc_FR-Base}
(a) (with $X=X^{n}$).

(d) (\ref{eq:Lim_Proc_(T,E0)-Mod}) implies (\ref{eq:Lim_Proc_TF_T-Mod_FD})
by Lemma \ref{lem:Proc_Rep} (a). (\ref{eq:Lim_Proc_TF_T-Mod_FD})
implies (\ref{eq:FR-Base}). Hence, (d) follows by (c) and Proposition
\ref{prop:RepProc_FR-Base} (a).

(e) follows by (\ref{eq:Q-Algebra_RepTF}) and Corollary \ref{cor:Base_Sep_Meas}
(a) (with $(d,A)=(\aleph(\mathbf{T}_{0}),\widehat{E}^{d})$).

(f) follows by Proposition \ref{prop:RepProc_FR-Base} (a) and (d,
e) (with $Y=\widehat{X}$).

(g) follows by $\mathcal{F}\subset C_{b}(E;\mathbf{R})$, Fact \ref{fact:FC_FDC}
and (f).\end{proof}

\begin{lem}
\label{lem:RepProc_RAP}Let $E$ be a topological space, $(\Omega,\mathscr{F},\mathbb{P};X)$
be an $E$-valued measurable process, $(E_{0},\mathcal{F};\widehat{E},\widehat{\mathcal{F}})$
be a base over $E$, $T\in(0,\infty)$ and $(\widetilde{\Omega},\widetilde{\mathcal{F}},\mathbb{P}^{T};X^{T})=\mathfrak{rap}_{T}(X)$.
Then:

\renewcommand{\labelenumi}{(\alph{enumi})}
\begin{enumerate}
\item If (\ref{eq:FR-Base}) or (\ref{eq:T-Base}) holds for some conull
$\mathbf{T}\subset\mathbf{R}^{+}$, then $X^{T}$ satisfies
\begin{equation}
\inf_{t\in\mathbf{R}^{+}}\mathbb{P}^{T}\left(\bigotimes\mathcal{F}\circ X_{t}^{T}\in\bigotimes\widehat{\mathcal{F}}(\widehat{E})\right)=1\label{eq:RAP_FR-Base}
\end{equation}
or
\begin{equation}
\inf_{t\in\mathbf{R}^{+}}\mathbb{P}^{T}\left(X_{t}^{T}\in E_{0}\right)=1.\label{eq:RAP_R+-Base}
\end{equation}
respectively.
\item If (\ref{eq:Path-Base}) holds for $S_{0}\subset E_{0}^{\mathbf{R}^{+}}$,
then
\begin{equation}
\mathbb{P}^{T}\left(X^{T}\in S_{0}\right)=1.\label{eq:RAP_Path_Base}
\end{equation}

\item If $\widehat{X}\in\mathfrak{rep}_{\mathrm{m}}(X;E_{0},\mathcal{F})$,
then $\widehat{X}^{T}=\mathfrak{rap}_{T}(X)\in\mathfrak{rep}_{\mathrm{m}}(X^{T};E_{0},\mathcal{F})$%
\footnote{The notations ``$\mathfrak{rep}_{\mathrm{m}}(\cdot;\cdot)$'' and
``$\mathfrak{rep}_{\mathrm{c}}(\cdot;\cdot)$'' were specified in
Notation \ref{notation:RepProc}.%
}.
\item If $\widehat{X}\in\mathfrak{rep}_{\mathrm{c}}(X;E_{0},\mathcal{F})$,
then $\widehat{X}^{T}=\mathfrak{rap}_{T}(X)\in\mathfrak{rep}_{\mathrm{c}}(X^{T};E_{0},\mathcal{F})$.
\end{enumerate}
\end{lem}
\begin{proof}
(a) One finds by the conullity of $\mathbf{T}$ that
\begin{equation}
\begin{aligned} & \inf_{t\in\mathbf{R}^{+}}\mathbb{P}^{T}\left(\bigotimes\mathcal{F}\circ X_{t}^{T}\in\bigotimes\widehat{\mathcal{F}}(\widehat{E})\right)\\
 & \geq\frac{1}{T}\int_{[t,T+t]\cap\mathbf{T}}\mathbb{P}\left(\bigotimes\mathcal{F}\circ X_{\tau}\in\bigotimes\widehat{\mathcal{F}}(\widehat{E})\right)d\tau=1
\end{aligned}
\label{eq:Check_RAP_FR-Base}
\end{equation}
and
\begin{equation}
\inf_{t\in\mathbf{R}^{+}}\mathbb{P}^{T}\left(X_{t}^{T}\in E_{0}\right)\geq\frac{1}{T}\int_{[t,T+t]\cap\mathbf{T}}\mathbb{P}\left(X_{\tau}\in E_{0}\right)d\tau=1.\label{eq:Check_RAP_R+-Basis}
\end{equation}

(b) follows by the fact that $(X^{T})^{-1}(S_{0})\supset\mathbf{R}^{+}\times X^{-1}(S_{0})$.

(c) follows by Proposition \ref{prop:RAP} (with $X=\widehat{X}$),
(\ref{eq:Define_RepProc}) and the fact that
\begin{equation}
\begin{aligned} & \frac{1}{T}\int_{0}^{T}\mathbb{P}\left(\bigotimes\mathcal{F}\circ X_{\tau+t}=\bigotimes\widehat{\mathcal{F}}\circ\widehat{X}_{\tau+t}\right)d\tau\\
 & \geq\frac{1}{T}\int_{0}^{T}\mathbb{P}\left(\bigotimes\mathcal{F}\circ X_{\tau+t}\in\bigotimes\widehat{\mathcal{F}}(\widehat{E})\right)d\tau,\;\forall t\in\mathbf{R}^{+}.
\end{aligned}
\label{eq:Check_RepProc_RAP}
\end{equation}

(d) follows by (c), Fact \ref{fact:Cadlag_RepProc} and Lemma \ref{lem:RAP_Cadlag}
(b) (with $X=\widehat{X}$).\end{proof}

\begin{lem}
\label{lem:FLP_Under_TF}Let $E$ be a topological space, $(E_{0},\mathcal{F};\widehat{E},\widehat{\mathcal{F}})$
be a base over $E$, $\mathbf{T}\subset\mathbf{R}^{+}$, $\{(\Omega^{n},\mathscr{F}^{n},\mathbb{P}^{n};X^{n})\}_{n\in\mathbf{N}}$
be $E$-valued processes satisfying (\ref{eq:Common_FR-Base_N}) (especially
(\ref{eq:Common_T-Base_N})), $\widehat{X}^{n}\in\mathfrak{rep}(X^{n};E_{0},\mathcal{F})$
for each $n\in\mathbf{N}$ and $(\Omega,\mathscr{F},\mathbb{P};Y)$
be an $\widehat{E}$-valued process. Then, (\ref{eq:RepProc_FC_along_T_Y})
implies $\varpi(\bigotimes\widehat{\mathcal{F}})\circ Y=\mathfrak{fl}_{\mathbf{T}}(\{\varpi(\bigotimes\mathcal{F})\circ X^{n}\}_{n\in\mathbf{N}})$.
\end{lem}
\begin{proof}
(\ref{eq:Common_FR-Base_N}) is weaker than (\ref{eq:Common_T-Base_N})
by Fact \ref{fact:Common_T-Base_FR-Base} (with $\mathbf{I}=\mathbf{N}$).
Define $Z\circeq\bigotimes\widehat{\mathcal{F}}\circ Y$, $Z^{n}\circeq\bigotimes\widehat{\mathcal{F}}\circ\widehat{X}^{n}$
and $\zeta^{n}\circeq\bigotimes\mathcal{F}\circ X^{n}$ for each $n\in\mathbf{N}$.
One finds that
\begin{equation}
\inf_{t\in\mathbf{T},n\in\mathbf{N}}\mathbb{P}^{n}\left(\zeta_{t}^{n}=Z_{t}^{n}\right)=1\label{eq:Phi(X)_Phihat(Xhat)_T-Mod_N}
\end{equation}
by Proposition \ref{prop:RepProc_FR-Base} (a). We fix $\mathbf{T}_{0}\in\mathscr{P}_{0}(\mathbf{T})$
and put $d\circeq\aleph(\mathbf{T}_{0})$. As mentioned in Note \ref{note:Ehat_Valued_Proc_FDD},
$\{\widehat{X}\}_{n\in\mathbf{N}}$, $\{\zeta^{n}\}_{n\in\mathbf{N}}$,
$\{Z^{n}\}_{n\in\mathbf{N}}$, $Y$ and $Z$ all have Borel finite-dimensional
distributions, so (\ref{eq:RepProc_FC_along_T_Y}) implies
\begin{equation}
\widehat{X}_{\mathbf{T}_{0}}^{n}\Longrightarrow Y_{\mathbf{T}_{0}}\mbox{ as }k\uparrow\infty\mbox{ on }\mathbf{R}^{d}.\label{eq:RepProc_FDD_WC_Y}
\end{equation}
One finds that
\begin{equation}
\varphi\circeq\bigotimes_{t\in\mathbf{T}_{0}}\left(\bigotimes\widehat{\mathcal{F}}\right)\circ\mathfrak{p}_{t}\in C\left[\widehat{E}^{d};(\mathbf{R}^{\infty})^{d}\right]\label{eq:FDD_Under_TF_Map}
\end{equation}
by (\ref{eq:Base_Imb}) and Fact \ref{fact:Prod_Map_2} (a, b). Hence,
it follows by (\ref{eq:FDD_Under_TF_Map}), (\ref{eq:RepProc_FDD_WC_Y})
and Continuous Mapping Theorem (Theorem \ref{thm:ContMapTh} (a))
that
\begin{equation}
Z_{\mathbf{T}_{0}}^{n}=\varphi\circ\widehat{X}_{\mathbf{T}_{0}}^{n}\Longrightarrow\varphi\circ Y_{\mathbf{T}_{0}}=Z_{\mathbf{T}_{0}}\mbox{ as }k\uparrow\infty\mbox{ on }(\mathbf{R}^{\infty})^{d}.\label{eq:Check_RepProc_Under_TF_Subseq_FC_T}
\end{equation}
$(\mathbf{R}^{\infty})^{d}$ and $\mathcal{P}((\mathbf{R}^{\infty})^{d})$
are Polish spaces by Proposition \ref{prop:Var_Polish} (f) and Theorem
\ref{thm:P(E)_Compact_Polish} (b) (with $E=(\mathbf{R}^{\infty})^{d}$).
Now, the result follows by (\ref{eq:Check_RepProc_Under_TF_Subseq_FC_T}),
(\ref{eq:Phi(X)_Phihat(Xhat)_T-Mod_N}) and Fact \ref{fact:FLP_Uni}
(with $E=\mathbf{R}^{\infty}$ and $X^{i}=\zeta^{n}$).\end{proof}

\addtocontents{toc}{\protect\setcounter{tocdepth}{0}}
\cleardoublepage
\newpage{}
\phantomsection
\addcontentsline{toc}{chapter}{Bibliography}
\renewcommand{\addcontentsline}[3]{}

\bibliographystyle{amsalpha}
\bibliography{nonlinearfiltering,probability,mathematics,appliedprobability}

\end{document}